\documentclass[twoside]{article}

\usepackage{tocbibind}
\usepackage{fancyhdr,bm,mathtools,graphicx}
\usepackage{amssymb,amsmath}
\usepackage{marvosym}
\usepackage{mathstyle}
\usepackage{hieroglf}
\usepackage[T1]{fontenc}
\usepackage{eufrak}
\usepackage[T1]{fontenc}
\usepackage{textcomp}
\usepackage{romanbar}
\usepackage[dvipsnames]{xcolor}
\usepackage{ifthen}
\usepackage{tikz-cd}
\tikzcdset{
arrow style=tikz,
diagrams={>={Stealth[scale=1]}}
}
\usepackage[T1]{fontenc}
\usepackage{textcomp}

\usepackage{stmaryrd}
\usepackage[mathscr]{euscript}
\usepackage{chngcntr}
\usepackage[style]{fncychap}
\usepackage{latexsym}

\def\Ref#1{[{equation~\ref{#1}}]}

\newcounter{chapterr}
\counterwithin{section}{chapterr}

\newcounter{definition}
\counterwithin{definition}{chapterr}
\def\definition
{\refstepcounter{definition}\vskip0.4\baselineskip\noindent{{\textbf{Definition~}}}{$\bf\arabic{chapterr}.$}{$\bf\arabic{definition}$}:\hskip 0.4\parindent}
\newcounter{theorem}
\counterwithin{theorem}{chapterr}
\def\theorem
{\refstepcounter{theorem}\vskip0.4\baselineskip\noindent{{\textbf{Theorem~}}}{$\bf\arabic{chapterr}.$}{$\bf\arabic{theorem}$}:\hskip 0.4\parindent}
\newcounter{symbol}
\counterwithin{symbol}{chapterr}

\newcounter{QMaxiom}
\counterwithin{QMaxiom}{chapterr}

\newcounter{corollary}
\counterwithin{corollary}{chapterr}
\def\corollary
{\refstepcounter{corollary}\vskip0.4\baselineskip\noindent{{\textbf{Corollary~}}}{$\bf\arabic{chapterr}.$}{$\bf\arabic{corollary}$}:\hskip 0.4\parindent}
\newcounter{lemma}
\counterwithin{lemma}{chapterr}
\def\lemma
{\refstepcounter{lemma}\vskip0.4\baselineskip\noindent{{\textbf{Lemma~}}}{$\bf\arabic{chapterr}.$}{$\bf\arabic{lemma}$}:\hskip 0.4\parindent}
\def\proof
{\vskip0.3\baselineskip\noindent{\textbf{proof}}\hskip.08\baselineskip:\hskip 0.4\parindent}

\newcounter{fixed}
\counterwithin{fixed}{chapterr}
\def\fixed
{\refstepcounter{fixed}\vskip0.4\baselineskip\noindent{{\textbf{Fixed Objects~}}}{$\bf\arabic{chapterr}.$}{$\bf\arabic{fixed}$}:\hskip 0.4\parindent}
\def\caution{\textasteriskcentered\hskip0.25\baselineskip}
\newcommand{\pr}[1]{{\textbf{\textsf{p{#1}}}}}

\def\V{V}
\def\C{\mathbb{C}}

\def\c{c}

\def\R{\mathbb{R}}

\def\X{{\rm X}}

\def\Z{\mathbb{Z}}

\newcommand{\p}[1]{{#1^{\prime}}}

\def\F{\mathscr{F}}
\def\C{\mathbb{C}}
\def\R{\mathbb{R}}

\def\f{{\rm F}}
\def\ff{{\rm F}}

\def\card{{\textsf{card}}}
\def\endef{{\flushright{\noindent$\blacksquare$\\}}\noindent}
\def\endthm{{\flushright{\noindent$\square$\\}}\noindent}
\def\endcor{{\flushright{\noindent$\square$\\}}\noindent}

\def\endlem{{\flushright{\noindent$\square$\\}}\noindent}
\def\endp{}

\def\endlem{{\flushright{\noindent$\square$\\}}\noindent}
\def\endfixed{{\flushright{\noindent$\Diamond$\\}}\noindent}

\def\then{\Rightarrow}
\def\thenn{\Leftrightarrow}

\def\({\left(}
\def\){\right)}
\def\[{\left[}
\def\]{\right]}

\newcommand{\CSs}[1]{{\mathcal{P}}\(#1\)}

\newcommand{\Card}[1]{\left| #1\right|}
\newcommand{\CarD}[1]{\card\(#1\)}
\def\cardeq{\overset{\underset{\mathrm{card}}{}}{=}}

\newcommand{\union}[1]{\bigcup#1}
\newcommand{\intersection}[1]{\bigcap#1}

\newcommand{\Union}[3]{\bigcup_{#1\in#2}#3}
\newcommand{\Unionn}[4]{\bigcup_{#1=#2}^{#3}#4}

\newcommand{\Dproduct}[3]{\prod_{{#1}\in{#2}}{#3}}
\newcommand{\cmp}[2]{#1\circ#2}
\newcommand{\Func}[2]{{\textsf{F}}\bpair{#1}{#2}}
\newcommand{\surFunc}[2]{{{\mathfrak{s}}\textsf{F}}\opair{#1}{#2}}
\newcommand{\IF}[2]{{\mathsf{B\hphantom{}F}}\opair{#1}{#2}}
\newcommand{\InF}[2]{{\mathfrak{i}}\hphantom{}{\mathsf{F}}\opair{#1}{#2}}

\newcommand{\domain}[1]{{\mathsf{dom}}\bbsingle{#1}}
\newcommand{\codomain}[1]{{\mathsf{codom}}\bbsingle{#1}}
\newcommand{\funcimage}[1]{{\mathsf{img}}\bbsingle{#1}}
\newcommand{\EqR}[1]{{\mathsf{EqR}}\bbsingle{#1}}
\newcommand{\EqClass}[2]{{#1}\left/{#2}\right.}
\newcommand{\pEqclass}[2]{\[{#1}\]_{#2}}
\newcommand{\PEqclass}[2]{{\mathsf{EqC}}\bbpair{#1}{#2}}
\newcommand{\Cprod}[2]{{#1}\times{#2}}

\newcommand{\OR}[2]{#1\thinspace\lor\thinspace#2}
\newcommand{\AND}[2]{#1\thinspace\land\thinspace#2}
\newcommand*{\suchthat}{\;\ifnum\currentgrouptype=16 \middle\fi|\;}
\newcommand{\Foreach}[2]{\forall\thinspace#1\in#2\negthinspace:\thickspace}
\newcommand{\Exists}[2]{\exists\thinspace#1\in#2\negthinspace:\thickspace}

\newcommand{\Existsu}[2]{\exists!\thinspace#1\in#2\negthinspace:\thickspace}
\newcommand{\Existsis}[2]{\(\exists\)\thinspace#1\in#2\negthinspace:\thickspace}

\newcommand{\defset}[3]{\left\{#1\in#2\thickspace:\thickspace#3\right\}}
\newcommand{\defsets}[3]{\left\{#1\subseteq#2\thickspace:\thickspace#3\right\}}
\newcommand{\defSet}[2]{\left\{#1\thickspace|\thickspace#2\right\}}
\newcommand{\negation}[1]{\neg{#1}}
\def\eqdef{\overset{\underset{\mathrm{def}}{}}{=}}
\def\indef{\overset{\underset{\mathrm{def}}{}}{\in}}

\newcommand{\resd}[1]{{\mathfrak{res}}{\mathsf{D}}_{#1}}
\newcommand{\rescd}[1]{{\mathfrak{res}}{\mathsf{C\negthinspace D}}_{#1}}
\newcommand{\res}[1]{{\mathfrak{res}}_{#1}}
\newcommand{\finv}[1]{{#1}^{-1}}

\def\hh{h}

\def\K{{\rm K}}

\newcommand{\Gen}[2]{{\langle#1\rangle}_{#2}}

\newcommand{\IG}[1]{e_{#1}}

\newcommand{\invG}[2]{{#1}^{-1}_{#2}}
\newcommand{\invg}[2]{{#1}^{-1}_{\tiny{#2}}}

\def\index{{\mathscr{I}}}
\def\Zp{\Z^{+}}
\def\Zn{\Z^{-}}

\def\Rp{\R^{+}}

\def\empty{\varnothing}
\def\Xt{{\mathbb{X}}}
\def\Yt{{\mathbb{Y}}}
\newcommand{\opair}[2]{\(#1,\thinspace #2\)}

\newcommand{\bpair}[2]{\negthinspace\left(#1,\thinspace #2\right)}
\newcommand{\bbpair}[2]{\negthinspace\left(#1;\thinspace #2\right)}
\newcommand{\bbsingle}[1]{\negthinspace\left({#1}\right)}

\newcommand{\topology}[1]{{\mathsf{T}}_{#1}}
\newcommand{\alltopologies}[1]{{\mathbf{Top}}\bbsingle{#1}}

\newcommand{\closedsets}[1]{{\mathsf{Closed}}\bbsingle{#1}\negthinspace}
\newcommand{\compl}[2]{#1\setminus#2}

\def\y{y}
\newcommand{\seta}[1]{\left\{#1\right\}}
\newcommand{\func}[2]{#1\(#2\)}

\def\x{x}

\def\U{U}

\def\ff{{\rm{f}}}
\def\a{a}

\newcommand{\stopology}[2]{\left. #1\right|_{#2}}

\def\asubset{A}

\def\csubset{C}

\def\point{p}

\newcommand{\nei}[1]{{\mathsf{Nei}}_{#1}}

\newcommand{\Cl}[1]{{\mathsf{Cl}}_{#1}}

\def\cf{f}
\def\cg{g}

\def\hf{h}

\newcommand{\image}[1]{#1^{\rightarrow}}
\newcommand{\pimage}[1]{#1^{\leftarrow}}
\newcommand{\CF}[2]{{\mathsf{CF}}\bpair{#1}{#2}}
\newcommand{\HOF}[2]{{\mathsf{HF}}\bpair{#1}{#2}}

\newcommand{\function}[3]{#1:\thinspace#2\to#3}

\newcommand{\atopology}[1]{{\boldsymbol{\tau}}_{#1}}
\newcommand{\binary}[2]{#1,\thinspace#2}

\def\aconnectedset{\csubset}

\newcommand{\connecteds}[1]{{\mathsf{Cnd}}\bbsingle{#1}}
\newcommand{\maxcon}[1]{{\mathsf{MaxCnd}}\bbsingle{#1}}

\newcommand{\Ointerval}[2]{\opair{#1}{#2}}
\newcommand{\COinterval}[2]{\[#1,\thinspace #2\)}
\newcommand{\OCinterval}[2]{\(#1,\thinspace #2\]}

\newcommand{\compacts}[1]{{\mathsf{CM\negthinspace P}}\bbsingle{#1}}

\newcommand{\identity}[1]{{\rm{Id}}_{#1}}

\newcommand{\collection}[1]{{\mathcal{S}}_{#1}}

\def\VSS{{\mathrm{V}}}
\newcommand{\VV}[1]{\VSS_{#1}}
\newcommand{\VVS}[1]{\VS_{#1}}
\newcommand{\W}[1]{{\mathrm{W}}_{#1}}
\newcommand{\WW}[1]{W_{#1}}

\newcommand{\subvec}[2]{{\mathsf{SubVec}}\bbpair{#1}{#2}}
\newcommand{\Vspan}[1]{{\mathsf{span}}_{#1}}
\def\VS{{\mathbb{V}}}
\newcommand{\NVS}[1]{{\mathscr{V}}_{#1}}

\newcommand{\vsum}[1]{+_{#1}}
\newcommand{\vsprod}[1]{\times_{#1}}
\newcommand{\spro}[1]{\cdot_{#1}}
\newcommand{\Lin}[2]{{\mathsf{L}}\bpair{#1}{#2}}
\newcommand{\Linisom}[2]{{\mathsf{LIsom}}\bpair{#1}{#2}}

\newcommand{\VLin}[2]{{\mathbb{L}}\bpair{#1}{#2}}

\newcommand{\ovecbasis}[1]{{\mathsf{O\negthinspace VBases}}\bbsingle{#1}}

\newcommand{\Det}[1]{{\mathrm{det}}_{#1}}
\newcommand{\triple}[3]{\opair{#1}{\binary{#2}{#3}}}
\newcommand{\tuple}[4]{\opair{#1}{\binary{#2}{\binary{#3}{#4}}}}
\newcommand{\mtuple}[2]{\(\suc{#1}{#2}\)}

\def\FF{{\mathbb{F}}}

\def\u{u}

\newcommand{\vv}[1]{v_{#1}}
\newcommand{\ww}[1]{w_{#1}}

\newcommand{\suc}[2]{{#1},\ldots,{#2}}

\newcommand{\Times}[2]{{#1}\times\ldots\times{#2}}

\newcommand{\Succ}[4]{\underbrace{{#1}{#3}\ldots{#3}{#2}}_{#4}}
\newcommand{\Succc}[3]{{#1}{#3},\ldots,{#3}{#2}}

\def\RR{\tilde{\R}}
\newcommand{\Ropenman}[2]{{\tilde{#1}}_{\RR^{#2}}}
\newcommand{\topR}[1]{\widehat{#1}}
\newcommand{\Rtanidentity}[2]{{\tilde{1}}^{\(#1\)}_{#2}}
\def\ointR{{\mathsf{OInterval}}_{0}\bbsingle{\R}}
\newcommand{\aninterval}[1]{I_{#1}}
\newcommand{\maxinterval}[1]{{\mathfrak{I}}_{#1}}
\newcommand{\maxintcurve}[1]{{\boldsymbol{\gamma}}_{#1}}
\newcommand{\funcprod}[2]{{#1}{\overset{\underset{\mathrm{\mathfrak{f}}}{}}{\times}}{#2}}
\newcommand{\Injection}[2]{{\mathrm{Inj}}_{{#1}\to{#2}}}
\newcommand{\mantop}[1]{{{\mathsf{Top}}^{\mathfrak{m}}}\bbsingle{#1}}
\newcommand{\mantops}[1]{\widehat{#1}}
\newcommand{\NVLin}[2]{{\overline{\mathbb{L}}}\bpair{#1}{#2}}
\newcommand{\Mat}[3]{{\mathsf{MAT}}\bbpair{#1}{\binary{#2}{#3}}}
\newcommand{\BMat}[3]{{\mathbb{MAT}}\bbpair{#1}{\binary{#2}{#3}}}
\newcommand{\TMat}[3]{{\widehat{\mathbb{MAT}}}\bbpair{#1}{\binary{#2}{#3}}}
\newcommand{\matelement}[3]{\[#1\]_{#2#3}}
\newcommand{\Man}[1]{\mathscr{M}_{#1}}
\newcommand{\Emsubman}[1]{{\mathsf{EmSubMan}}\bbsingle{#1}}
\newcommand{\emsubman}[2]{{\mathsf{emsubman}}\bbpair{#1}{#2}}
\newcommand{\asmooth}[1]{\xi_{#1}}
\newcommand{\charttransfer}[3]{{\mathsf{chart}}^{\star}_{\triple{#1}{#2}{#3}}}
\newcommand{\manprod}[2]{{#1}{\overset{\underset{\mathrm{\mathfrak{m}}}{}}{\times}}{#2}}
\newcommand{\topprod}[2]{{#1}{\overset{\underset{\mathrm{\mathfrak{t}}}{}}{\times}}{#2}}
\newcommand{\topologyofspace}[1]{{\mathcal{T}}\bbsingle{#1}}
\newcommand{\maxatlases}[3]{{\mathsf{maxAtl}}^{\(#1\)}\bpair{#2}{#3}}
\newcommand{\atlases}[3]{{\mathsf{Atl}}^{\(#1\)}\bpair{#2}{#3}}
\newcommand{\rectatlas}[3]{{\mathsf{RectAtl}}^{\opair{#1}{#2}}\bbsingle{#3}}
\newcommand{\maxatlasgen}[3]{{\mathfrak{maxAtl}}^{\(#1\)}_{\opair{#2}{#3}}}
\newcommand{\M}[1]{\mathrm{M}_{#1}}
\newcommand{\maxatlas}[1]{{\boldsymbol{\mathscr{A}}}_{#1}}
\newcommand{\atlas}[1]{\mathscr{A}_{#1}}
\newcommand{\difclass}[1]{\mathsf{C}^{#1}}
\newcommand{\vecf}[2]{{\mathsf{V\negthinspace F}}^{#2}\bbsingle{#1}}
\newcommand{\Vecf}[2]{{\mathbb{V\negthinspace F}}^{#2}\bbsingle{#1}}
\newcommand{\Lievecf}[2]{{\boldsymbol{\mathfrak{vf}}}^{#2}\bbsingle{#1}}
\newcommand{\avecf}[1]{{\mathscr{X}}_{#1}}
\newcommand{\avecff}[1]{{\mathscr{Y}}_{#1}}
\newcommand{\tanspace}[2]{{\mathsf{T}}_{#1}{#2}}
\newcommand{\Tanspace}[2]{{\mathbb{T}}_{#1}{#2}}
\newcommand{\avec}[1]{{\boldsymbol{v}}_{#1}}
\newcommand{\tanbun}[1]{{\mathsf{T}}#1}
\newcommand{\Tanbun}[1]{{\mathbf{T}}#1}
\newcommand{\basep}[1]{\pi_{#1}}
\newcommand{\mapdifclass}[3]{{\mathsf{C}}^{#1}\bpair{#2}{#3}}
\newcommand{\Lmapdifclass}[3]{{\mathbf{C}}^{#1}\bpair{#2}{#3}}
\newcommand{\Diffeo}[3]{{\mathsf{Diff}}^{#1}\bpair{#2}{#3}}
\newcommand{\Diff}[2]{{\mathsf{Diff}}^{#1}\bbsingle{#2}}
\newcommand{\GDiff}[2]{{\mathbf{Diff}}^{#1}\bbsingle{#2}}
\newcommand{\banachmapdifclass}[5]{{\mathsf{C}}^{#1}_{\opair{#2}{#3}}\bpair{#4}{#5}}
\newcommand{\Derivation}[2]{{\mathsf{Der}}^{#2}\bbsingle{#1}}
\newcommand{\LDerivation}[2]{{\mathbb{D}}^{#2}\bbsingle{#1}}
\newcommand{\aderivation}[1]{\Delta_{#1}}

\def\rdot{\cdot}
\newcommand{\tanchart}[2]{{\mathscr T}^{#1}_{#2}}
\newcommand{\tanatlas}[1]{{\mathfrak{TA}}_{#1}}
\newcommand{\tanspaceiso}[3]{{\boldsymbol{\theta}}_{#1}^{\opair{#2}{#3}}}

\newcommand{\der}[3]{{\mathsf{D}}^{\opair{#2}{#3}}{#1}}
\newcommand{\derop}[2]{{\mathsf{D}}^{\opair{#1}{#2}}}
\newcommand{\derr}[1]{{\mathsf{D}}{#1}}
\newcommand{\banachder}[3]{{\mathscr{D}}^{\opair{#2}{#3}}{#1}}

\newcommand{\Rder}[2]{{\mathsf{d}}_{#2}{#1}}
\newcommand{\Rderop}[1]{{\mathsf{d}}_{#1}}
\newcommand{\LieDer}[1]{\bar{\mathcal{L}}^{#1}}
\newcommand{\Lieder}[1]{{\mathcal{L}}^{#1}}
\newcommand{\LieDerinv}[1]{\bar{\chi}_{#1}}
\newcommand{\Liederinv}[1]{{\chi}_{#1}}
\newcommand{\projection}[2]{{\boldsymbol{\mathscr{P}}}^{\(#1\)}_{#2}}

\newcommand{\Eucbase}[2]{{\boldsymbol{e}}^{\(#1\)}_{#2}}
\newcommand{\deltaf}[2]{\delta_{{#1},{#2}}}
\newcommand{\support}[1]{{\mathsf{supp}}\bbsingle{#1}}
\newcommand{\zerovec}[1]{{\boldsymbol{0}}_{#1}}
\newcommand{\fextension}[4]{{\mathsf{Ex}}_{#1}\bbpair{#2}{\binary{#3}{#4}}}
\newcommand{\convex}[1]{{\mathsf{Cnx}}\bbsingle{#1}}
\newcommand{\Liebracket}[1]{{\mathfrak{L\negthinspace B}}_{#1}}
\newcommand{\liebracket}[3]{\left\llbracket#1,~#2\right\rrbracket_{#3}}
\newcommand{\aliealg}[1]{{\boldsymbol{\ell}}_{#1}}
\newcommand{\aliealgebra}[1]{{\boldsymbol{\mathfrak{l}}}_{#1}}
\newcommand{\subspace}[2]{\left. #1\right|_{#2}}

\newcommand{\sublie}[1]{{\mathsf{SubLieAl}}\bbsingle{#1}}
\newcommand{\subliedim}[2]{{\mathsf{SubLieAl}}^{\(#1\)}\bbsingle{#2}}
\newcommand{\Lker}[2]{{\mathsf{LieKer}}_{\opair{#1}{#2}}}
\newcommand{\liedim}[1]{{\mathrm{dim}}_{\mathrm{lie}}\bbsingle{#1}}
\newcommand{\mandim}[1]{{\mathrm{dim}}_{\mathrm{man}}\bbsingle{#1}}
\newcommand{\LS}[1]{L_{#1}}
\newcommand{\oliebasis}[1]{{\mathsf{O\negthinspace LBases}}\bbsingle{#1}}
\newcommand{\cinterval}[2]{\[#1,\thinspace#2\]}
\newcommand{\partialder}[1]{\partial_{#1}}
\newcommand{\derived}[1]{\p{#1}}
\newcommand{\immersion}[3]{{\mathsf{Immersion}}^{#1}\bpair{#2}{#3}}
\newcommand{\submersion}[3]{{\mathsf{Submersion}}^{#1}\bpair{#2}{#3}}
\newcommand{\embedding}[3]{{\mathsf{Embedding}}^{#1}\bpair{#2}{#3}}
\newcommand{\foliation}[1]{\boldsymbol{\mathscr{F}}_{#1}}
\newcommand{\leaf}[1]{{\mathscr{L}}_{#1}}

\newcommand{\Foliations}[2]{{\mathsf{Foliation}}^{\(\infty\)}\bbpair{#2}{#1}}
\newcommand{\Leaves}[2]{{\mathsf{Leaf}}_{\opair{#1}{#2}}}
\newcommand{\folatlas}[3]{{\mathsf{FolAtl}}{\bbpair{\binary{#1}{#3}}{#2}}}
\newcommand{\leafatlas}[3]{{\mathsf{LeafAtl}}_{\bbpair{\binary{#1}{#3}}{#2}}}
\newcommand{\leafman}[3]{\boldsymbol{\mathsf{LeafMan}}_{\bbpair{\binary{#1}{#3}}{#2}}}
\newcommand{\Leafman}[1]{\widetilde{#1}}
\newcommand{\Projection}[3]{{\boldsymbol{\mathscr{P}}}^{\(#1\)}_{\suc{#2}{#3}}}
\newcommand{\foltanbun}[3]{{\mathsf{FT}}_{\bbpair{\binary{#1}{#3}}{#2}}}
\newcommand{\distribution}[1]{\boldsymbol{\mathscr{D}}_{#1}}
\newcommand{\Distributions}[2]{{\mathsf{Distribution}}^{\(\infty\)}\bbpair{#2}{#1}}
\newcommand{\subtanbun}[2]{{\mathsf{T}}^{#2}#1}
\newcommand{\integrabledist}[2]{{\mathsf{IntDist}}^{\(\infty\)}\bbpair{#2}{#1}}
\newcommand{\involutivedist}[2]{{\mathsf{InvolDist}}^{\(\infty\)}\bbpair{#2}{#1}}
\newcommand{\distvecf}[3]{{\mathsf{DistVF}}{\bbpair{\binary{#1}{#3}}{#2}}}
\newcommand{\intdistsolution}[2]{\int_{\bbpair{#1}{#2}}}
\newcommand{\curve}[1]{\alpha_{#1}}
\newcommand{\integralcurves}[1]{{\mathsf{IntCurves}}_{#1}}
\newcommand{\vfflow}[2]{{\boldsymbol{\Phi}}_{\bbpair{#1}{#2}}}
\newcommand{\vfFlow}[2]{{\boldsymbol{\Psi}}_{\bbpair{#1}{#2}}}
\newcommand{\leftparinj}[3]{{\mathfrak{LI}}^{\opair{#1}{#2}}_{#3}}
\newcommand{\rightparinj}[3]{{\mathfrak{RI}}^{\opair{#1}{#2}}_{#3}}
\newcommand{\prodmantan}[4]{{\boldsymbol{\Theta}}^{\opair{#1}{#2}}_{\opair{#3}{#4}}}
\newcommand{\aliegroup}[1]{{\mathscr{J}}_{#1}}
\newcommand{\bu}[1]{{\boldsymbol{u}}_{#1}}
\newcommand{\bw}[1]{{\boldsymbol{w}}_{#1}}
\newcommand{\G}[1]{{\mathrm{G}}_{#1}}
\newcommand{\g}[1]{g_{#1}}
\newcommand{\gop}[1]{\bullet_{#1}}
\newcommand{\gpower}[3]{{#1}^{{#2}{#3}}}
\newcommand{\gsetprod}[3]{\({#1}{#2}\)_{#3}}
\newcommand{\Group}[1]{{\mathbb{G}}_{#1}}
\newcommand{\ginv}[1]{{\mathfrak{inv}}_{#1}}
\newcommand{\gopr}[2]{\bar{#1}_{#2}}
\newcommand{\constmap}[2]{{\mathfrak{C}}_{\opair{#1}{#2}}}
\newcommand{\diagmap}[1]{{\mathfrak{diag}}_{#1}}
\newcommand{\gltrans}[2]{{\mathfrak{Lt}}^{\(#1\)}_{#2}}
\newcommand{\grtrans}[2]{{\mathfrak{Rt}}^{\(#1\)}_{#2}}
\newcommand{\gconj}[2]{{\mathfrak{Cnj}}^{\(#1\)}_{#2}}
\newcommand{\nuclei}[1]{{\mathsf{Nuc}}\bbsingle{#1}}
\newcommand{\aliebra}[1]{{\boldsymbol{\ell}}_{#1}}
\newcommand{\aliealgmor}[1]{{\mathcal{H}}_{#1}}
\newcommand{\LiealgMor}[2]{{\mathsf{LieAlgMor}}\bpair{#1}{#2}}
\newcommand{\LiealgIsom}[2]{{\mathsf{LieAlgIsom}}\bpair{#1}{#2}}
\newcommand{\LiealgMon}[2]{{\mathsf{LieAlgMonmor}}\bpair{#1}{#2}}
\newcommand{\LiealgEpi}[2]{{\mathsf{LieAlgEpimor}}\bpair{#1}{#2}}
\newcommand{\liealgisomorphic}[2]{{#1}\overset{\underset{\mathfrak{l}}{}}{\bumpeq}{#2}}
\newcommand{\Topgroup}[1]{{\mathcal{G}}_{#1}}
\newcommand{\Liegroup}[1]{{\mathscr{G}}_{#1}}
\newcommand{\lietop}[1]{{{\mathsf{Top}}^{\mathfrak{l}}}\bbsingle{#1}}
\newcommand{\lietops}[1]{\widehat{#1}}
\newcommand{\topgtops}[1]{\widehat{#1}}
\newcommand{\LieG}[1]{\overline{#1}}
\newcommand{\TopgG}[1]{\overline{#1}}
\newcommand{\Lieman}[1]{\widetilde{#1}}
\newcommand{\Lietopg}[1]{\overline{\widehat{#1}}}
\newcommand{\Topsubgroups}[1]{{\mathsf{TopSubGr}}\bbsingle{#1}}
\newcommand{\IndTop}[1]{{\mathsf{IndTop}}_{#1}}
\newcommand{\Liepoints}[1]{\overset{\underset{\mathrm{\therefore}}{}}{#1}}
\newcommand{\Leftinvvf}[1]{{\mathsf{L\negthinspace IV\negthinspace F}}\bbsingle{#1}}
\newcommand{\VLeftinvvf}[1]{{\mathbb{L\negthinspace IV\negthinspace F}}\bbsingle{#1}}
\newcommand{\LiegroupLiealgebra}[1]{{\mathfrak{lg}}\bbsingle{#1}}
\newcommand{\LiegroupLiealgebratan}[1]{{\bar{\mathfrak{lg}}}\bbsingle{#1}}
\newcommand{\liegvftan}[1]{\Xi_{#1}}
\newcommand{\tanliebracket}[1]{\aliealg{#1}}
\newcommand{\Lliebracket}[1]{{\mathfrak{lb}}_{#1}}
\newcommand{\GHom}[2]{{\mathsf{GHom}}\bpair{#1}{#2}}
\newcommand{\GIsom}[2]{{\mathsf{GIsom}}\bpair{#1}{#2}}
\newcommand{\Subgroups}[1]{{\mathsf{SubGr}}\bbsingle{#1}}
\newcommand{\NSubgroups}[1]{{\mathsf{NSubGr}}\bbsingle{#1}}
\newcommand{\asubgroup}[1]{{\mathrm{H}}_{#1}}
\newcommand{\LCoset}[1]{{\mathsf{Lcoset}}_{#1}}
\newcommand{\aliemor}[1]{{\mathscr{H}}_{#1}}
\newcommand{\LieMor}[2]{{\mathsf{LieGrMor}}\bpair{#1}{#2}}
\newcommand{\LieIsom}[2]{{\mathsf{LieGrIsom}}\bpair{#1}{#2}}
\newcommand{\LieAut}[1]{{\mathsf{LieGrAut}}\bbsingle{#1}}
\newcommand{\GLieAut}[1]{{\mathbf{LieGrAut}}\bbsingle{#1}}
\newcommand{\indliemor}[2]{\Upsilon_{\opair{#1}{#2}}}
\newcommand{\lieisomorphic}[2]{{#1}\overset{\underset{\mathfrak{L}}{}}{\bumpeq}{#2}}

\newcommand{\imsubgroup}[1]{{\mathsf{ImSubGr}}\bbsingle{#1}}
\newcommand{\connectedimsubgroup}[1]{{\mathsf{CndImSubGr}}\bbsingle{#1}}
\newcommand{\emsubgroup}[1]{{\mathsf{EmSubGr}}\bbsingle{#1}}
\newcommand{\subliealgdist}[2]{{\boldsymbol{\mathscr{D}}}_{\opair{#1}{#2}}}
\newcommand{\liesublie}[2]{{\boldsymbol{\mathfrak{G}}}_{\bpair{#1}{#2}}}
\newcommand{\liesublieset}[2]{{\mathfrak{G}}_{\bpair{#1}{#2}}}
\newcommand{\liesubgalcor}[1]{{\boldsymbol{\Omega}}_{#1}}
\newcommand{\expLie}[1]{{\mathrm{exp}}_{#1}}
\newcommand{\ring}[1]{{\mathrm{R}}_{#1}}
\newcommand{\Ring}[1]{{\mathbf{R}}_{#1}}
\newcommand{\ringid}[1]{\overset{\underset{\mathrm{r}}{}}{{\boldsymbol{1}}}_{#1}}
\newcommand{\ringzero}[1]{\overset{\underset{\mathrm{r}}{}}{{\boldsymbol{0}}}_{#1}}
\newcommand{\rr}[1]{r_{#1}}
\newcommand{\module}[1]{\mu_{#1}}
\newcommand{\Module}[1]{{\boldsymbol{\mu}}_{#1}}
\newcommand{\modulezero}[1]{\overset{\underset{\mathrm{m}}{}}{{\boldsymbol{0}}}_{#1}}
\newcommand{\mm}[1]{{\mathrm{m}}_{#1}}
\newcommand{\Submodules}[1]{{\mathsf{SubMod}}\bbsingle{#1}}
\newcommand{\modulegen}[1]{{\mathsf{ModGen}}_{#1}}
\newcommand{\Modulebases}[1]{{\mathsf{ModBases}}\bbsingle{#1}}
\newcommand{\amodmor}[1]{{\mathcal{H}}_{#1}}
\newcommand{\ModuleMors}[2]{{\mathsf{ModMor}}\bpair{#1}{#2}}
\newcommand{\ModuleIsoms}[2]{{\mathsf{ModIsom}}\bpair{#1}{#2}}
\newcommand{\smoothfprod}[1]{\cdot_{#1}}
\newcommand{\smoothfsum}[1]{\overset{\underset{\mathrm{r}}{}}{+}_{#1}}
\newcommand{\smoothfsub}[1]{\overset{\underset{\mathrm{r}}{}}{-}_{#1}}
\newcommand{\smoothvfprod}[1]{\overset{\underset{\mathrm{vf}}{}}{\odot}_{#1}}
\newcommand{\smoothvfsum}[1]{\overset{\underset{\mathrm{vf}}{}}{+}_{#1}}
\newcommand{\smoothvfsub}[1]{\overset{\underset{\mathrm{vf}}{}}{-}_{#1}}
\newcommand{\smoothderprod}[1]{\overset{\underset{\mathrm{der}}{}}{\odot}_{#1}}
\newcommand{\smoothdersum}[1]{\overset{\underset{\mathrm{der}}{}}{+}_{#1}}
\newcommand{\smoothdersub}[1]{\overset{\underset{\mathrm{der}}{}}{-}_{#1}}
\newcommand{\smoothring}[2]{{\mathscr{R}}\difclass{#1}\bbsingle{#2}}
\newcommand{\smoothvfmodule}[2]{{\boldsymbol{\mathscr{M}}}\vecf{#2}{#1}}
\newcommand{\smoothdermodule}[2]{{\boldsymbol{\mathscr{M}}}\Derivation{#2}{#1}}
\newcommand{\Injec}[1]{{\mathrm{I}}_{#1}}

\newcommand{\myitem}[1]{${\fontsize{6.65}{7}\selectfont{\textbf{#1}}}$}
\def\varfill{\dotfill}

\def\toclevel@section{1}\def\toclevel@subection{2}
\addtocontents{toc}{\string\let\string\l@section\string\l@subsection}
\def\toclevel@subsection{2}\def\toclevel@subsubection{3}
\addtocontents{toc}{\string\let\string\l@subsection\string\l@subsubsection}

\usepackage{tocloft}
\newcommand{\chapteR}[1]{\cleardoublepage
{\refstepcounter{chapterr}\vskip\baselineskip\centering{\fontsize{21}{21}\selectfont${\bf\Roman{chapterr}}$
\vskip0.6\baselineskip{\fontsize{21}{21}\selectfont{\textbf{#1}}}}
\vskip 5.9\baselineskip}
\addcontentsline{toc}{0}
{\protect\vskip0.5\baselineskip\noindent\bf{\Roman{chapterr}\hskip0.5\baselineskip#1\hspace{\fill}}}\par
\fancyhead[LO]{\ifthenelse{\value{chapterr}=0}{#1}{$\bf\Roman{chapterr}$\hskip0.7\baselineskip#1}}
}
\newcommand{\Bibliography}[1]{\vskip0.5\baselineskip\centering{\huge\bf{References}}\vskip \baselineskip
\addcontentsline{toc}{0}
{\protect\vskip0.5\baselineskip\noindent\bf{References}\hspace{\fill}}\par
\fancyhead[LO,RE]{\bf{References}#1}
}

\def\refthm#1{[{theorem~\ref{#1}}]}
\def\reflem#1{[{lemma~\ref{#1}}]}
\def\refdef#1{[{definition~\ref{#1}}]}
\def\refcor#1{[{corollary~\ref{#1}}]}

\newcommand{\mathleft}{\@fleqntrue\@mathmargin0pt}

\renewcommand{\sectionmark}[1]{\ifthenelse{\value{section}=0}{\markright{#1}{}}
{\markright{${\arabic{section}}$ #1}{}}}
\makeatletter
\def\cleardoublepage{\clearpage\if@twoside \ifodd\c@page\else
\hbox{}
\vspace*{\fill}
\vspace{\fill}
\thispagestyle{empty}
\newpage
\if@twocolumn\hbox{}\newpage\fi\fi\fi}
\makeatother
\fancypagestyle{plain}{
\fancyhf{}
\fancyhead{}
\fancyhead[R]{\thepage}
\fancyhead[L]{\leftmark}
}
\newcommand{\newsymp}[1]{{#1}\equiv}
\newcommand{\newsymb}[1]{{#1}\equiv}

\newcommand{\SET}[1]{{\mathscr{S}}_{#1}}
\newcommand{\FUNCTION}[1]{{{f}}_{#1}}

\newcommand{\prop}[1]{{\mathfrak{p}}\llparenthesis{#1}\rrparenthesis}
\def\Prop{\mathfrak{p}}
\newcommand{\propos}[1]{{\mathbf{\mathfrak{p}}}_{#1}}
\def\dummy{\centerdot}

\usepackage[hang,flushmargin]{footmisc}
\renewcommand{\footnoterule}{%
  \kern 20pt
  \hrule width \textwidth height 0.5pt
  \kern 5pt
}
\def\quotl{``}
\def\quotr{"}

\usepackage[inner=3.8cm,outer=3cm,bottom=3.9cm,twoside]{geometry}
\begin{document}
\thispagestyle{empty}
\noindent
{\\ \\\textbf{\fontsize{40}{40}\selectfont
{\textsf{Real Lie Groups}}}}
\\[0.8\baselineskip]{\textbf{\fontsize{20}{20}\selectfont
{\textsf{of Finite Dimension}}}}
\\[8\baselineskip]
\noindent
{\fontsize{21}{21}\selectfont
{\textsf{Farzad Shahi}}
}
\\[11\baselineskip]
\noindent
{\fontsize{11}{11}\selectfont
{\textrm{Version:} 1.00}
}
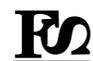
\vfill\hfill
$\underline{\Huge{\textsf{\bf F}}\negthickspace\negthickspace\negthinspace{{\rotatebox{90}{\textsf{\bf S}}}}}$
\newpage
\thispagestyle{empty}
\noindent
{\fontsize{9.4}{9.4}\selectfont
{\underline{{\bf\textsf{Author:}} Farzad Shahi}}}\\
{{\bf\textsf{email}}}:~{\texttt{shahi.farzad@gmail.com}}
\vskip 4\baselineskip
\noindent
{\fontsize{9.4}{9.4}\selectfont
{\underline{{\bf\textsf{Version:}} 1.00}}
}\\
\noindent
{\fontsize{9.4}{9.4}\selectfont
{\bf\textsf{2021}}}
\vskip 4\baselineskip
\noindent
{\fontsize{9.4}{9.4}\selectfont
{\bf\textsf{Typesetting:}} By the author, using \TeX}
\vskip 5\baselineskip
\noindent
{\fontsize{9.4}{9.4}\selectfont
{\bf\textsf{Abstract:}} The current version offers an introduction to the basics of the theory of finite-dimensional Lie groups over the field of real numbers.
The notion of the tangent-space of a manifold at a point is considered to be defined via the well known chart-vector formalism, here, a formalism equivalent
to other commonly used ones (namely, the curve and derivation methods).
The proofs of all assertions about Lie groups are in alignment with this formalism, here.}
\newpage
\thispagestyle{empty}
\noindent
{\fontsize{20}{20}\selectfont
{\bf Preface}}
\vskip1.5\baselineskip
\noindent
The set of all invertible matrices of an arbitrary fixed order $n$ over the field of complex numbers endowed
with the binary operation of matrix multiplication, possesses the structure of a group, which goes under the name of
$\quotl$general linear group of order $n$$\quotr$.
Any subgroup of a general linear group
is called a $\quotl$linear group$\quotr$.
A natural topology can be adapted to the set of all square matrices of an arbitrary order $n$ over $\C$,
transferred from that of $\C^{n^2}$ via a canonical one-to-one correspondence between them.
Therefore, any linear group, and in particular any general linear group, can be equipped with a natural topology
inherited from that of the set of all matrices of the corresponding order over $\C$. Each linear group is
considered as a topological space with this topology.
Since the spaces $\C^{n^2}$ and $\R^{2n^2}$ with their canonical topologies are homeomorphic, this topological space
is also homeomorphic to $\R^{2n^{2}}$.
It can be verified that the group operation and the
inverse mapping of a linear group are continuous with respect to its canonical topology.
This behaviour can be generalized to the abstract notion of a topological group. Any set
endowed with compatible structures of a group and a topological space is considered to be a mathematical
structure of this type. Thereby, linear groups, having their canonical topologies attached,
are obvious instances of this structure.\\
Additionally, general linear groups can be viewed inherently as geometric objects. More precisely, given a positive integer $n$,
the canonical differentiable structure of  $\R^{2n^2}$ can be carried to the set of all square matrices of order $n$
through their canonical isomorphism, and thus turning it into a $2n^2$-dimensional smooth manifold modeled on real Euclidean-space of dimension $2n^2$.
As an straightforward feature of this differentiable structure, its underlying topological-space is in alignment
with the natural topology of the set of $n$-matrices up to homeomorphism. It is worthwhile to mention that there is an alternative natural treatment, that is to
carry the canonical complex differentiable structure of $\C^{n^2}$ to the set of all $n$-matrices, making it as an
$n$-dimensional complex manifold in alignment with its initially given topology. This treatment is not the concern of this survey.
\\According to the continuity of the determinant function on the set of all $n$-matrices, and
making use of the fact that the general linear group of order $n$ consists of the $n$-matrices
with non-zero determinant, it becomes clear that the general linear group of order $n$ is an open subset of
the set of all $n$-matrices. Therefore the general linear group of order $n$ canonically inherits the differentiable
structure of the set of all $n$-matrices, and hence becomes a $2n^2$-dimensional smooth and real manifold.
This differentiable structure is compatible with the group structure of general linear group of order $n$, in the sense that
the binary operation and inverse mapping of this group are smooth (infinitely differentiable) with respect to it.
The abstract setting that unites such objects goes under the name of $\quotl$real smooth group$\quotr$ or
$\quotl$real Lie group$\quotr$. Precisely, a real smooth group refers to any set endowed simultaneously with structures of a
group and a smooth real manifold in a way that these structures are compatible with each other in the above sense.
The dimension of a real smooth group is defined to be the dimension of its underlying real manifold. Thus, for any positive integer $n$,
the general linear group of order $n$ is canonically a $2n^2$-dimensional real smooth group.\\
Since any smooth mapping between any pair of smooth real manifolds is continuous with respect
to the topological spaces corresponded to them, it becomes evident that any real smooth group holds interinsically
the structure of a corresponding topological group. Hence any behaviour of the abstract topological groups is essentially
a behaviour of the abstract real smooth groups.\\
Any subgroup of a real smooth group is referred to as a closed subgroup of it if the set of elements of
that subgroup is a closed set within the underlying topological space of that smooth group.
It is a well-known fact in the theory of real smooth groups that a differentiable structure can be adapted to
any closed subgroup of a given real smooth group in a way that makes it a real smooth group that can be embedded
in the initially given smooth group. Hence, as a consequence of a well-known feature of embedded (regular) submanifolds of
a smooth real manifold, any closed subgroup of a real smooth group canonically inherits a unique such differentiable structure from that of
the original smooth group. Any closed subgroup of a smooth group
endowed with this differentiable structure is called a $\quotl$closed smooth subgroup of that smooth group$\quotr$.
The definition of an embedded submanifold of
a real manifold implies that the underlying topological space of a closed smooth subgroup of a real smooth group
coincides with the induced topology from that of the real smooth group.
Therefore, for a given positive
integer $n$, each closed subgroup of the general linear group of order $n$ is canonically a real smooth group
whose underlying topology is in alignment with the subgroup's initially given canonical topology.
So every linear group that is closed in its corresponded general linear group can be regarded naturally
as a real smooth group. Such real smooth groups are simply called $\quotl$closed linear groups$\quotr$.\\
Closed linear groups constitute an important class of real smooth groups. Aside from them, there are
other crucial instances of real smooth groups that arise in different theories. For example, the Poincaré group
that consists of all isometries of space-time can be regarded as a real smooth group that is a natural setting
for symmetries of space-time. The Lorentz groups is a subgroup of the Poincaré group that consists of those isometries of space-time that
leave the origin of space-time intact, which is also canonically a real smooth group. These smooth groups are among the
important real smooth groups that arise in mathematical physics.\\
The study of the notion of abstract real smooth group is an essential tool to reveal many commonly shared characteristics of
every instance of this mathematical notion, and in particular those important ones discussed above.
Any map between (the set of points of) a pair of real smooth groups
is considered to be structure-preserving if it is both a group homomorphism and a smooth map
when the intrinsic group strucures and the intrinsic smooth manifold structures of the smooth groups are taken into account, respectively.
Real smooth groups along with the structure-preserving maps between them form the objects and morphisms of a category. This category
is called the $\quotl$smooth group category$\quotr$. A pair of smooth groups is regarded to be isomorphic when there is a
smooth group isomorphism between them.\\
The most fundamental characteristic of a real smooth group is an algebraic structure called the Lie-algebra of that
smooth group. The Lie-algebra of a smooth group is a unique Lie-algebraic structure on a specific subspace of the
vector-space of all smooth vector-fields on the underlying manifold of the smooth group, called the $\quotl$left-invariant
vector-fields on the smooth group$\quotr$. The Lie-algebras of any pair of isomorphic smooth groups are also isomorphic Lie-algebraically.
The converse is not true generally, but there is a sub-category of smooth group category where this behavior becomes two-sided.\\
The set of all loops of a topological-space based at a fixed point is partitioned via the equivalence relation
of path-homotopy on this set, the quotient of which possesses the structure of a group when endowed with an approporiate
operation, going under the name of the $\quotl$fundamental-group of the topological-space at that point$\quotr$.
The fundamental-groups of a path-connected topological-space at all points are the same up to group isomorphism.
So, the fundamental-group of a path-connected topological-space is a well-defined notion which can be identified
with the fundamental-group at any arbitrary point (up to group isomorphism). A simply-connected topological-space
is defined to be a path-connected topological-space with trivial fundamental-group (the group with one element).
When the underlying topological-space of a smooth group is simply-connected, the smooth group is also called a $\quotl$simply-connected
smooth group$\quotr$.
Since path-connectedness and connectedness are equivalent properties of the underlying topological-space of a smooth manifold
(and hence that of a smooth group), a simply-connected smooth group is actually a connected smooth-group,
all loops at any point of which are path-homotopic to the stationary path of that point (with respect to its intrinsic topology).\\
Any smooth group isomorphic to a simply-connected smooth group must be also a simply-connected one. The reason is, the intrinsic
topological-spaces of a pair of isomorphic smooth groups are homeomorphic, connectedness (or path-connectedness) is invariant under homeomorphisms,
and fundamental-groups are preserved under homeomorphisms (up to group-isomorphism).
When the Lie-algebras of a pair of smooth groups are isomorphic and one of them is known to be simply-connected, then the smooth groups themselves
must be isomorphic, and consequently the other one is also revealed to be simply-connected.
Therefore, by restricting the attention to simply-connected smooth groups, the isomorphism of a pair of smooth groups
is equivalent to the isomorphism of the Lie-algebras of them. Hence, it is anticipated that in this sub-category of smooth groups,
the classification problem of smooth groups can be carried to the task of classifying Lie-algebras to a considerable extent.\\
The Lie-algebra structure of smooth groups plays further crucial roles in a variety of other natural and important problems
that arise in the theory of smooth groups. Some direct or indirect applications include
identifying all possible immersed subgroups of a smooth group and
providing a proof of the famous closed subgroup theorem, investigations of the actions of smooth groups on manifolds,
and analysis on smooth groups.\\
As the final consideration, it is important to be aware of a problem which is not discussed in the main body of this text.
The problem is, any set with a topological group structure that simultaneously possesses the structure of a topological manifold in accordance with
it (meaning that the underlying topological-spaces of both structures agree, and the operation and inverse mapping of
the intrinsic group of the topological group are both continuous with regard to this topology), inherently
holds the structure of a unique smooth group. This problem is a famous for the $\quotl$Hilbert's fifth problem$\quotr$
which is solved by \textit{D. Montgomery}. \textit{L. Zippin}, and \textit{A. M. Gleason} \cite{Gleason, MontgomeryZippin}.
In addition, the converse of this problem is known to be clearly true. That is, every smooth group
holds uniquely the structure of a topological group and a consistent topological manifold.
So, any topological group that is simultaneously a consistent topological manifold, ultimately can be seen as another face of a
corresponded smooth group.
It is also worthwhile to notice that in this problem, when the assumption of existence of a topological manifold consistent
with the topological group is dropped, there is no gurantee for the existence of a smooth group with that topological group.
As an instance, consider the case of the group $\opair{\R}{+}$ that is endowed with the discrete topology on $\R$.
There can be found no smooth group in accordance with this topological group, because in that case
every singleton $\seta{\x}$ ($\x\in\R$) would be homeomorphic to an open set of a Euclidean-space, which is impossible.\\
\vskip0.5\baselineskip
\noindent
{\bf{A brief explanation about chapters:}}
The first and fourth chapters are devoted to the review of basic definitions and problems of differential geometry necessary to the
development of the theory of smooth groups.
The second chapter studies some fundamental problems of the more general theory of topological groups. Since any
smooth group is also a topological group, this chapter actually gives an amount of information about the structure
of smooth groups. Additionally, the theorems and corollaries of this chapter may be of considerable importance in dealing
with some problems in the smooth groups theory.
The third chapter is a very brief introduction to the abstract theory of Lie-algebras, just to the extent that
is crucial for investigating the smooth groups.
The fifth chapter is considered to be the core one which studies several aspects of smooth groups.
\vskip0.5\baselineskip
\noindent
{\bf{Prerequisites:}}
A considerable knowledge of linear algebra, general topology, and differential geometry
(the theory of differentiable manifolds), along with
a thorough comprehension of the differential calculus of Cartan over Banach spaces \cite{Cartan}
is presumed. A basic knowledge of the algebraic stuctures of a group, ring, and module is needed as well,
even though a brief review of essentials of each is provided. There is no need for an extensive
awareness of the theories of Lie algebras and topological groups in order to grasp the core part of this text,
and the the amount of information provided in the related chapters will suffice. Actually, the theory of Lie algebras
extends far beyond the concerns of this text.
\vskip\baselineskip
\hfill
{\textsf{Farzad Shahi}}
\newpage
\thispagestyle{empty}
\section*{\fontsize{21}{21}\selectfont\bf{Contents}}
\addtocontents{toc}{\protect\setcounter{tocdepth}{-1}}
\tableofcontents
\addtocontents{toc}{\protect\setcounter{tocdepth}{3}}
\newpage
\newpage
\chapteR{
Mathematical Notations
}
\thispagestyle{fancy}
\section*{
Set-theory
}
\sectionmark{Mathematical Notations}
$\newsymb{\empty}$
empty-set
\varfill
$\empty$\\
$\newsymp{\SET{1}=\SET{2}}$
$\SET{1}$
equals
$\SET{2}$.
\varfill
$=$\\
$\newsymp{\SET{1}\in\SET{2}}$
$\SET{1}$
is an element of
$\SET{2}$.
\varfill
$\in$\\
$\newsymp{\SET{1}\ni\SET{2}}$
$\SET{1}$
contains
$\SET{2}$.
\varfill
$\ni$\\
$\newsymb{\seta{\binary{\SET{1}}{\SET{2}}}}$
the set composed of
$\SET{1}$
and
$\SET{2}$
\varfill
$\seta{\binary{\dummy}{\dummy}}$\\
$\newsymb{\defset{\SET{1}}{\SET{2}}{\prop{\SET{1}}}}$
all elements of
$\SET{2}$
having the property
$\Prop$
\varfill
$\defset{\dummy}{\dummy}{\dummy}$\\
$\newsymb{\union{\SET{}}}$
union of all elements of
$\SET{}$
\varfill
$\bigcup$\\
$\newsymb{\intersection{\SET{}}}$
intersection of all elements of
$\SET{}$
\varfill
$\bigcap$\\
$\newsymb{\CSs{\SET{}}}$
power-set of
$\SET{}$
\varfill
$\CSs{\dummy}$\\
$\newsymb{\Dproduct{\alpha}{\index}{\SET{\alpha}}}$
Cartesian-product of the collection of indexed sets
${\seta{\SET{\alpha}}}_{\alpha\in\index}$
\varfill
$\prod$\\
$\newsymp{\SET{1}\subseteq\SET{2}}$
$\SET{1}$
is a subset of
$\SET{2}$.
\varfill
$\subseteq$\\
$\newsymp{\SET{1}\supseteq\SET{2}}$
$\SET{1}$
includes
$\SET{2}$.
\varfill
$\supseteq$\\
$\newsymp{\SET{1}\subset\SET{2}}$
$\SET{1}$
is a proper subset of
$\SET{2}$.
\varfill
$\subset$\\
$\newsymp{\SET{1}\supset\SET{2}}$
$\SET{1}$
properly includes
$\SET{2}$.
\varfill
$\supset$\\
$\newsymb{\SET{1}\cup\SET{2}}$
union of
$\SET{1}$
and
$\SET{2}$
\varfill
$\cup$\\
$\newsymb{\SET{1}\cap\SET{2}}$
intersection of
$\SET{1}$
and
$\SET{2}$
\varfill
$\cap$\\
$\newsymb{\SET{1}\times\SET{2}}$
Cartesian-product of
$\SET{1}$
and
$\SET{2}$
\varfill
$\times$\\
$\newsymb{\compl{\SET{1}}{\SET{2}}}$
the relative complement of
$\SET{2}$
with respect to
$\SET{1}$
\varfill
$\setminus$\\
$\newsymb{\func{\FUNCTION{}}{\SET{}}}$
value of the function
$\FUNCTION{}$
at
$\SET{}$
\varfill
$\func{\dummy}{\dummy}$\\
$\newsymb{\domain{\FUNCTION{}}}$
domain of the function
$\FUNCTION{}$
\varfill
$\domain{\dummy}$\\
$\newsymb{\codomain{\FUNCTION{}}}$
codomain of the function
$\FUNCTION{}$
\varfill
$\codomain{\dummy}$\\
$\newsymb{\funcimage{\FUNCTION{}}}$
image of the function
$\FUNCTION{}$
\varfill
$\funcimage{\dummy}$\\
$\newsymb{\image{\FUNCTION{}}}$
image-map of the function
$\FUNCTION{}$
\varfill
$\image{\dummy}$\\
$\newsymb{\pimage{\FUNCTION{}}}$
inverse-image-map of the function
$\FUNCTION{}$
\varfill
$\pimage{\dummy}$\\
$\newsymb{\resd{\FUNCTION{}}}$
domain-restriction-map of the function
$\FUNCTION{}$
\varfill
$\resd{\dummy}$\\
$\newsymb{\rescd{\FUNCTION{}}}$
codomain-restriction-map of the function
$\FUNCTION{}$
\varfill
$\rescd{\dummy}$\\
$\newsymb{\func{\res{\FUNCTION{}}}{\SET{}}}$
domain-restriction and codomain-restriction of\\ the function
$\FUNCTION{}$ to $\SET{}$ and $\func{\image{\FUNCTION{}}}{\SET{}}$, respectively
\varfill
$\res{\dummy}$\\
$\newsymb{\Func{\SET{1}}{\SET{2}}}$
the set of all maps from
$\SET{1}$
to
$\SET{2}$
\varfill
$\Func{\dummy}{\dummy}$\\
$\newsymb{\IF{\SET{1}}{\SET{2}}}$
the set of all bijective functions from
$\SET{1}$
to
$\SET{2}$
\varfill
$\IF{\dummy}{\dummy}$\\
$\newsymb{\finv{\FUNCTION{}}}$
the inverse mapping of the bijective function $\FUNCTION{}$
\varfill
$\finv{\dummy}$\\
$\newsymb{\surFunc{\SET{1}}{\SET{2}}}$
the set of all surjective functions from
$\SET{1}$
to
$\SET{2}$
\varfill
$\surFunc{\dummy}{\dummy}$\\
$\newsymb{\cmp{\FUNCTION{1}}{\FUNCTION{2}}}$
composition of the function
$\FUNCTION{1}$
with the function
$\FUNCTION{2}$
\varfill
$\cmp{}{}$\\
$\newsymb{\Injection{\SET{1}}{\SET{2}}}$
the injection-mapping of the set $\SET{1}$ into the set $\SET{2}$
\varfill
$\Injection{\dummy}{\dummy}$\\
$\newsymb{\funcprod{\cf_1}{\cf_2}}$
the function-product of the function $\cf_1$ and $\cf_2$
\varfill
$\Cprod{\dummy}{\dummy}$\\
$\newsymb{\EqR{\SET{}}}$
the set of all equivalence relations on the set
$\SET{}$
\varfill
$\EqR{\dummy}$\\
$\newsymb{\EqClass{\SET{1}}{\SET{2}}}$
quotient-set of
$\SET{1}$
by the equivalence-relation
$\SET{1}$
\varfill
$\EqClass{\dummy}{\dummy}$\\
$\newsymb{\pEqclass{\SET{1}}{\SET{2}}}$
equivalence-class of
$\SET{1}$
by the equivalence-relation
$\SET{2}$
\varfill
$\pEqclass{\dummy}{\dummy}$\\
$\newsymb{\PEqclass{\SET{1}}{\SET{2}}}$
equivalence-class of
$\SET{2}$
by the equivalence-relation
$\SET{1}$
\varfill
$\PEqclass{\dummy}{\dummy}$\\
$\newsymp{\Card{\SET{1}}\cardeq\Card{\SET{2}}}$
$\IF{\SET{1}}{\SET{2}}$
is non-empty.
\varfill
$\Card{\dummy}\cardeq\Card{\dummy}$\\
$\newsymb{\CarD{\SET{}}}$
cardinality of
$\SET{}$
\varfill
$\CarD{\dummy}$
\section*{
Logic
}
$\newsymp{\AND{\propos{1}}{\propos{2}}}$
$\propos{1}$
and
$\propos{2}$.
\varfill
$\AND{}{}$\\
$\newsymp{\OR{\propos{1}}{\propos{2}}}$
$\propos{1}$
or
$\propos{2}$.
\varfill
$\OR{}{}$\\
$\newsymp{{\propos{1}}\then{\propos{2}}}$
if
$\propos{1}$,
then
$\propos{2}$.
\varfill
$\then$\\
$\newsymp{{\propos{1}}\thenn{\propos{2}}}$
$\propos{1}$,
if-and-only-if
$\propos{2}$.
\varfill
$\thenn$\\
$\newsymp{\negation{\propos{}}}$
$\propos{1}$,
negation of
$\propos{}$.
\varfill
$\negation{}$\\
$\newsymp{\Foreach{\SET{1}}{\SET{2}}{\prop{\SET{1}}}}$
for every
$\SET{2}$
in
$\SET{1}$,
$\prop{\SET{1}}$.
\varfill
$\Foreach{\dummy}{\dummy}\dummy$\\
$\newsymp{\Exists{\SET{1}}{\SET{2}}{\prop{\SET{1}}}}$
exists
$\SET{2}$
in
$\SET{1}$ such that
$\prop{\SET{1}}$.
\varfill
$\Foreach{\dummy}{\dummy}\dummy$\\
\section*{
Group Theory
}
$\newsymp{\GHom{\Group{1}}{\Group{2}}}$
the set of all group-homomorphisms from the group $\Group{1}$ targeted to the
group $\Group{2}$
\varfill
$\GHom{\dummy}{\dummy}$\\
$\newsymp{\GIsom{\Group{1}}{\Group{2}}}$
the set of all group-isomorphisms from the group $\Group{1}$ targeted to the
group $\Group{2}$
\varfill
$\GIsom{\dummy}{\dummy}$\\
$\newsymp{\Subgroups{\Group{}}}$
the set of all subgroups of the group $\Group{}$
\varfill
$\Subgroups{\dummy}$\\
$\newsymp{\func{\LCoset{\Group{}}}{\asubgroup{}}}$
the set of all left-cosets of the subgroup $\asubgroup{}$ of the group $\Group{}$
\varfill
$\func{\LCoset{\dummy}}{\dummy}$
\section*{
Differential Calculus
}
$\newsymp{\banachmapdifclass{r}{\NVS{1}}{\NVS{2}}{\U_1}{\U_2}}$
the set of all $r$-times differentiable maps from the Banach-space $\NVS{1}$
to the Banach-space $\NVS{2}$ with domain $\U_1$ and codomain $\U_2$
\varfill
$\banachmapdifclass{\dummy}{\dummy}{\dummy}{\dummy}{\dummy}$\\
$\newsymp{\banachder{\FUNCTION{}}{\NVS{1}}{\NVS{2}}}$
derived map of $\FUNCTION{}$, $\FUNCTION{}$ being a $\difclass{r}$ map
from the Banach-space $\NVS{1}$ to the Banach-space $\NVS{2}$
(with open domain and codomain
in $\NVS{1}$ and $\NVS{2}$, respectively)
\varfill
$\banachder{\dummy}{\dummy}{\dummy}$
\section*{
Linear Algebra
}
$\newsymp{\subvec{\VVS{}}{m}}$
the set of all sets of vectors of all $m$-dimensional vector-subspaces of
the vector-space $\VVS{}$
\varfill
$\subvec{\dummy}{\dummy}$\\
$\newsymp{\func{\Vspan{\VVS{}}}{\asubset}}$
the vector-subspace of the vector-space $\VVS{}$ spanned by the subset $\asubset$
of the set of all vectors of $\VVS{}$
\varfill
$\func{\Vspan{\dummy}}{\dummy}$\\
$\newsymp{\ovecbasis{\VVS{}}}$
the set of all ordered-bases of the vector-space $\VVS{}$
\varfill
$\ovecbasis{\dummy}$\\
$\newsymp{\Lin{\VVS{1}}{\VVS{2}}}$
the set of all linear maps from the vector-space $\VVS{1}$ to the vector-space $\VVS{2}$
\varfill
$\Lin{\dummy}{\dummy}$\\
$\newsymp{\Linisom{\VVS{1}}{\VVS{2}}}$
the set of all linear isomorphisms from the vector-space $\VVS{1}$ to the vector-space
$\VVS{2}$
\varfill
$\Linisom{\dummy}{\dummy}$\\
$\newsymp{\VLin{\VVS{1}}{\VVS{2}}}$
the canonical vector-space of all linear maps from the vector-space
$\VVS{1}$ to the vector-space $\VVS{2}$
\varfill
$\VLin{\dummy}{\dummy}$\\
$\newsymp{\NVLin{\NVS{1}}{\NVS{2}}}$
the canonical Banach-space of all linear maps from the finite-dimensional Banach-space
$\NVS{1}$ to the finite-dimensional Banach-space $\VVS{2}$
\varfill
$\NVLin{\dummy}{\dummy}$\\
$\newsymp{\Mat{\F}{m}{n}}$
the set of all $m\times n$ matrices over the field $\F$
\varfill
$\Mat{\dummy}{\dummy}{\dummy}$\\
$\newsymp{\BMat{\F}{m}{n}}$
the canonical Banach-space of all $m\times n$ matrices over the field $\F$
(endowed with its natural norm)
\varfill
$\BMat{\dummy}{\dummy}{\dummy}$\\
$\newsymp{\Det{n}}$
the determinant function on the set of all square $\R$-matrices of degree $n$
\varfill
$\Det{\dummy}$
\section*{Topology}
$\newsymp{\alltopologies{\SET{}}}$
the set of all topologies on the set $\SET{}$
\dotfill
$\alltopologies{\dummy}$\\
$\newsymp{\topologyofspace{\Xt}}$
topology of the topological-space $\Xt$
\dotfill
$\topologyofspace{\dummy}$\\
$\newsymp{\closedsets{\Xt}}$
the set of all closed sets of the topological-space $\Xt$
\dotfill
$\closedsets{\dummy}$\\
$\newsymp{\func{\Cl{\Xt}}{\SET{}}}$
the closure of $\SET{}$ in the topological-space $\Xt$
\dotfill
$\func{\Cl{\dummy}}{\dummy}$\\
$\newsymp{\func{\nei{\Xt}}{\asubset}}$
the set of all open sets of the topological-space $\Xt$ including the subset $\asubset$
of $\Xt$ (the set of all neighborhoods of $\asubset$ in $\Xt$)
\dotfill
$\func{\nei{\dummy}}{\dummy}$\\
$\newsymp{\func{\IndTop{\Xt}}{\SET{}}}$
the topology on the subset $\SET{}$ of the topological-space $\Xt$ induced (inherited)
from that of $\Xt$
\dotfill
$\func{\IndTop{\dummy}}{\dummy}$\\
$\newsymp{\CF{\Xt}{\Yt}}$
the set of all continuous maps from the topological-space $\Xt$ to the topological-space
$\Yt$
\varfill
$\CF{\dummy}{\dummy}$\\
$\newsymp{\HOF{\Xt}{\Yt}}$
the set of all homeomorphisms from the topological-space $\Xt$ to the topological-space
$\Yt$
\varfill
$\HOF{\dummy}{\dummy}$\\
$\newsymp{\topprod{\Xt}{\Yt}}$
the topological-product of the topological-spaces $\Xt$ and $\Yt$
\varfill
$\topprod{\dummy}{\dummy}$\\
$\newsymp{\connecteds{\Xt}}$
the set of all connected sets of the topological-space $\Xt$
\varfill
$\connecteds{\dummy}$\\
$\newsymp{\maxcon{\Xt}}$
the set of all connected components of the topological-space $\Xt$
\varfill
$\maxcon{\dummy}$
$\newsymp{\TMat{\F}{m}{n}}$
the canonical topological-space of all $m\times n$ matrices over the field $\F$
(with the topology induced by its natural norm)
\varfill
$\TMat{\dummy}{\dummy}{\dummy}$\\
$\newsymp{\topR{\R^n}}$
the canonical topological-space of the Euclidean-space $\R^n$
(with the topology induced by its natural norm)
\varfill
$\topR{}$
\section*{Differential Geometry}
$\newsymp{\atlases{r}{\M{}}{\NVS{}}}$
the set of all atlases of differentiablity class $\difclass{r}$
on the set $\M{}$
modeled on the Banach-space $\NVS{}$
\varfill
$\atlases{\dummy}{\dummy}{\dummy}$\\
$\newsymp{\maxatlases{r}{\M{}}{\NVS{}}}$
the set of all maximal-atlases of differentiablity class $\difclass{r}$
on the set $\M{}$
constructed upon the Banach-space $\NVS{}$
\varfill
$\maxatlases{\dummy}{\dummy}{\dummy}$\\
$\newsymp{\func{\maxatlasgen{r}{\M{}}{\NVS{}}}{\atlas{}}}$
the maximal-atlas of differentiablity class $\difclass{r}$
on the set $\M{}$ modeled on
the Banach-space $\NVS{}$, generated by the atlas $\atlas{}$ in
$\atlases{r}{\M{}}{\NVS{}}$
\varfill
$\func{\maxatlasgen{\dummy}{\dummy}{\dummy}}{\dummy}$
\section*{Mathematical Environments}
$\newsymp{\blacksquare}$
end of definition
\varfill
$\blacksquare$\\
$\newsymp{\square}$
end of theorem, lemma, or corollary
\varfill
$\square$\\
$\newsymp{\Diamond}$
end of the introduction of new fixed objects
\varfill
$\Diamond$
\chapteR{Review of Differential Geometry}
\thispagestyle{fancy}
\section{Basic Concepts of smooth manifolds}
\textit{Here, the category of $\difclass{\infty}$ manifolds is decided to be the category that consists
$\difclass{\infty}$ differentiable-structures without boundary modeled on a Banach-space of any finite non-zero dimension with the
Hausdorff and second-countable underlying topological-space as the objects, and the smooth
mappings between such $\difclass{\infty}$ differentiable-structures as the morphisms.}\\
\fixed
$\Man{}=\opair{\M{}}{\maxatlas{}}$ is fixed as an $n$-dimensional and $\difclass{\infty}$ manifold
modeled the Banach-space $\R^n$. So, $\maxatlas{}$ is a $\difclass{\infty}$ maximal-atlas on
$\M{}$ modeled on the Banach-space $\R^{n}$, and the topology on $\M{}$ induced by this maximal-atlas
is Hausdorff and second-countable (i.e. possessing a countable base). The dimension of $\Man{}$
is denoted by $\mandim{\Man{}}$.\\
Also, $\defSet{\Man{i}=\opair{\M{i}}{\maxatlas{i}}}{i\in\Zp}$ is fixed as a collection of manifolds such that
for each positive integer $i$, $\Man{i}$ is an
$m_{i}$-dimensional and $\difclass{\infty}$ manifold modeled on $\R^{m_i}$, where each $m_i$
is a positive integer.
\endfixed
\begin{itemize}
\item[\myitem{DG~1.}]
For every point $\point$ of $\Man{}$, any chart $\phi$ of $\Man{}$ whose domain contains
$\point$ is called a $\quotl$chart of $\Man{}$ around $\point$$\quotr$ or
a $\quotl$neighbourhood chart of $\point$ in the manifold $\Man{}$$\quotr$.
For every point $\point$ of $\Man{}$, any chart $\phi$
of $\Man{}$ around $\point$ such that $\func{\phi}{\point}=\zerovec{}$ is called a
$\quotl$chart of $\Man{}$ centered at $\point$$\quotr$, and such a chart exists at every point of $\Man{}$.
\item[\myitem{DG~2.}]
For every positive integer $n$,
$\RR^n$ denotes the $n$-dimensional Euclidean-space endowed with its canonical differentiable-structure
constructed upon the Banach-space $\R^n$ that arises
from the trivial atlas consisting merely of the identity map on $\R^{n}$.
So $\identity{\R^n}$ is a chart of the smooth manifold $\RR^n$.
$\RR^n$ is called the $\quotl$canonical differentiable structure of $\R^n$$\quotr$.
\item[\myitem{DG~3.}]
The topology on $\M{}$ induced by the maximal-atlas $\maxatlas{}$ is denoted by $\mantop{\Man{}}$,
and the corresponding topological-space is denoted by $\mantops{\Man{}}$, that is,
$\mantops{\Man{}}:=\opair{\M{}}{\mantop{\Man{}}}$. When $\mantop{\Man{}}$ coincides with
an initially given topology on $\M{}$, it is customary to say that $\quotl$the manifold $\Man{}$
is compatible with that topology$\quotr$.
\item[\myitem{DG~4.}]
For every positive integer $n$ and every integer $k$ in $\seta{\suc{1}{n}}$, $\projection{n}{k}$
denotes the projection mapping of $\R^{n}$ onto its $k$-th factor. That is, $\projection{n}{k}$
is the element of $\Func{\R^n}{\R}$ such that for every $\mtuple{\x_1}{\x_n}$ in $\R^n$,
$\func{\projection{n}{k}}{\suc{\x_1}{\x_n}}=\x_k$. Furthermore, for every positive integer $n$,
$\Eucbase{n}{}$ denotes the standard ordered-base of $\R^n$ endowed with its canonical linear-structure.
That is, $\Eucbase{n}{}$ is the element of $\Func{\seta{\suc{1}{n}}}{\R^n}$ such that for every $k$
in $\seta{\suc{1}{n}}$, $\func{\projection{n}{j}}{\Eucbase{n}{k}}=\deltaf{k}{j}$.
\item[\myitem{DG~5.}]
$\mapdifclass{\infty}{\Man{}}{\Man{1}}$ denotes the set of all mappings
$\cf$ in $\Func{\M{}}{\M{1}}$ such that for every point $\point$ of $\Man{}$, there exists a chart $\phi$ in $\maxatlas{}$
and a chart $\psi$ in $\maxatlas{1}$ such that $\point\in\domain{\phi}$,
$\func{\image{\cf}}{\domain{\phi}}\subseteq\domain{\psi}$, and,
\begin{equation}\label{eqdesmoothnessdefinition}
\cmp{\psi}{\cmp{\cf}{\finv{\phi}}}\in\banachmapdifclass{\infty}{\R^{n}}{\R^{m}}{\funcimage{\phi}}{\funcimage{\psi}}.
\end{equation}
Additionally, given a mapping $\cf$ in $\mapdifclass{\infty}{\Man{}}{\Man{1}}$, for any chart $\phi$ of $\Man{}$
and any chart $\psi$ of $\Man{1}$ such that $\func{\pimage{\cf}}{\domain{\psi}}\cap\domain{\phi}\neq\empty$,
\Ref{eqdesmoothnessdefinition} is satisfied.
Each element of $\mapdifclass{\infty}{\Man{}}{\Man{1}}$ is called a $\quotl$smooth map from $\Man{}$ to $\Man{1}$$\quotr$.
Every smooth map from the manifold $\Man{}$ to the manifold $\Man{1}$ is a continuous map from the underlying topological-space of
$\Man{}$ to the underlying topological-space of $\Man{1}$. That is,
$\mapdifclass{\infty}{\Man{}}{\Man{1}}\subseteq\CF{\mantops{\Man{}}}{\mantops{\Man{1}}}$.
\item[\myitem{DG~6.}]
$\Diffeo{\infty}{\Man{}}{\Man{1}}$ denotes the set of all bijective maps $\function{\cf}{\M{}}{\M{1}}$ such that
$\cf$ is a smooth map from $\Man{}$ to $\Man{1}$ and $\finv{\cf}$ is a smooth map from $\Man{1}$ to $\Man{}$. That is,
\begin{equation}
\Diffeo{\infty}{\Man{}}{\Man{1}}:=\defset{\cf}{\IF{\M{}}{\M{1}}}{\cf\in\mapdifclass{\infty}{\Man{}}{\Man{1}},~
{\finv{\cf}}\in\mapdifclass{\infty}{\Man{1}}{\Man{}}}.
\end{equation}
Each element of $\Diffeo{\infty}{\Man{}}{\Man{1}}$ is called an $\quotl$$\infty$-diffeomorphism from $\Man{}$to $\Man{1}$$\quotr$.
It is said that $\quotl$$\Man{}$ is diffeomorphic to $\Man{1}$$\quotr$ iff there exists at least one $\infty$-diffeomorphism
from $\Man{}$ to $\Man{1}$. The existence of diffeomorphism between manifolds clearly induces an equivalence-relation
on a collection of manifolds. Furthermore, evidently
every $\infty$-diffeomorphism from the manifold $\Man{}$ to $\Man{1}$ is a homeomorphism from the underlying topological-space of
$\Man{}$ to the underlying topological-space of $\Man{1}$. That is,
$\Diffeo{\infty}{\Man{}}{\Man{1}}\subseteq\HOF{\mantops{\Man{}}}{\mantops{\Man{1}}}$. So, diffeomorphic manifolds
are homeomorphic, but the converse can not be true in general.\\
The set of all $\infty$-diffeomorphisms from $\Man{}$ to $\Man{}$ is simply denoted by $\Diff{\infty}{\Man{}}$, which
together with the binary operation of function-composition
on them, that is the pair $\opair{\Diff{\infty}{\Man{}}}{\cmp{}{}}$, is obviously a group which is denoted by
$\GDiff{\infty}{\Man{}}$. This group is referred to as the $\quotl$$\infty$-diffeomorphism group of the manifold $\Man{}$$\quotr$, or
the $\quotl$group of $\infty$-automorphisms of $\Man{}$$\quotr$. Each element of the set $\Diff{\infty}{\Man{}}$ is also
referred to as an $\quotl$$\infty$-automorphism of the manifold $\Man{}$$\quotr$.
The composition of any pair of $\infty$-automorphisms of the manifold $\Man{}$
is an $\infty$-automorphism of $\Man{}$.\\
Every $\infty$-automorphism of $\Man{}$ brings about a transference of charts of $\Man{}$. Precisely,
given an $\infty$-automorphism $\cf$ of $\Man{}$, for every point $\point$ of $\Man{}$ and every chart $\phi$ of
$\Man{}$ around $\point$, $\cmp{\phi}{\cf}$ is a chart of $\Man{}$ around $\func{\finv{\cf}}{\point}$.
So such transference of charts exhibits a pullback-like behavior when considered pointwise.
Thus for every $\infty$-automorphism $\cf$ of $\Man{}$ and every point $\point$ of $\Man{}$,
it can be defined a mapping $\charttransfer{\Man{}}{\cf}{\point}$ such that,
\begin{align}\label{eqcharttransfer}
&\function{\charttransfer{\Man{}}{\cf}{\point}}{\defset{\psi}{\maxatlas{}}{\point\in\domain{\psi}}}
{\defset{\psi}{\maxatlas{}}{\func{\finv{\cf}}{\point}\in\domain{\psi}}},\cr
&\Foreach{\phi}{\defset{\psi}{\maxatlas{}}{\point\in\domain{\psi}}}
\func{\charttransfer{\Man{}}{\cf}{\point}}{\phi}\eqdef\cmp{\phi}{\cf}.
\end{align}
\item[\myitem{DG~7.}]
$\mapdifclass{\infty}{\Man{}}{\RR}$ endowed with its canonical linear-structure is denoted by $\Lmapdifclass{\infty}{\Man{}}{\RR}$
whose addition and scalar-multiplication operations are defined pointwise.
$\rdot$ denotes the usual maltiplication of real-valued functions on $\M{}$, which together with the addition operation
in the linear-structure of $\mapdifclass{\infty}{\Man{}}{\RR}$ induces the structure of a ring on
$\mapdifclass{\infty}{\Man{}}{\RR}$. Also, $\Lmapdifclass{\infty}{\Man{}}{\RR}$
together with the binary operation $\rdot$ is a commutative and associative $\R$-algebra.
Each element of $\mapdifclass{\infty}{\Man{}}{\RR}$ is referred to as a $\quotl$real-valued smooth map on $\Man{}$$\quotr$.
\item[\myitem{DG~8.}]
For every $\point$ in $\M{}$, the set of all tangent vectors to $\Man{}$ at $\point$ is denoted by $\tanspace{\point}{\Man{}}$,
and correspondingly the tangent-space of $\Man{}$ at $\point$ is denoted by $\Tanspace{\point}{\Man{}}$ which is a real vector-space,
whose addition and scalar-product operations are defined via any chart around $\point$ and simultaneously independent of the choice of such a chart.
For a given point $\point$ of $\Man{}$, the addition and scalar-product operations of $\Tanspace{\point}{\Man{}}$ are denoted by
$\vsum{\opair{\Man{}}{\point}}$ and $\vsprod{\opair{\Man{}}{\point}}$, respectively. When there is no chance of confusion,
$\vv{1}\vsum{}\vv{2}$ and $\c\vv{}$ can replace $\vv{1}\vsum{\opair{\Man{}}{\point}}\vv{2}$ and $\c\vsprod{\opair{\Man{}}{\point}}\vv{}$, respectively.
For every $\point$ in $\Man{}$ and every chart $\phi$ in $\maxatlas{}$ whose domain contains $\point$, there exists a canonical linear-isomorphism
from $\Tanspace{\point}{\Man{}}$ to $\R^n$ denoted by $\tanspaceiso{\point}{\Man{}}{\phi}$.
The neutral element of addition in the linear-structure of $\Tanspace{\point}{\Man{}}$, that is $\zerovec{\Tanspace{\point}{\Man{}}}$,
is simply denoted by $\zerovec{\point}$.
For every point $\point$ of $\Man{}$, every chart $\phi$ of $\Man{}$, and every $\vv{}$ in $\R^n$,
the element $\func{\finv{\[\tanspaceiso{\point}{\Man{}}{\phi}\]}}{\vv{}}$ of $\tanspace{\point}{\Man{}}$ is called the
$\quotl$vector of the tangent-space of $\Man{}$ at the point $\point$ corresponded to $\vv{}$ with respect to the chart $\phi$$\quotr$.
\item[\myitem{DG~9.}]
The set of all tangent vectors to $\Man{}$ is denoted by $\tanbun{\Man{}}$.
The base-point-identifier of $\tanbun{\Man{}}$ is denoted by $\basep{\Man{}}$. That is,
$\basep{\Man{}}$ is defined to be a map from $\tanbun{\Man{}}$ to $\M{}$ such that for every $\point$ in $\M{}$
and every vector $\avec{}$ in $\tanspace{\point}{\Man{}}$, $\func{\basep{\Man{}}}{\avec{}}=\point$.
The tangent-bundle of $\Man{}$ is denoted by $\Tanbun{\Man{}}$, which is a $2n$-dimensional and $\difclass{\infty}$ manifold
built on the set $\tanbun{\Man{}}$ with its canonical differentiable-structure constructed upon the Banach-space $\Cprod{\R^n}{\R^n}$,
obtained from that of $\Man{}$.
Precisely, the set $\defSet{\tanchart{\Man{}}{\phi}}{\phi\in\maxatlas{}}$ is an atlas on $\tanbun{\Man{}}$, where,
\begin{align}\label{eqtangentbundlemaps}
\Foreach{\phi}{\maxatlas{}}
\left\{
\begin{aligned}
&\tanchart{\Man{}}{\phi}\in\Func{\func{\pimage{\basep{\Man{}}}}{\domain{\phi}}}{\Cprod{\funcimage{\phi}}{\R^n}},\cr
&\Foreach{\avec{}}{\func{\pimage{\basep{\Man{}}}}{\domain{\phi}}}
\func{\tanchart{\Man{}}{\phi}}{\avec{}}\eqdef
\opair{\func{\(\cmp{\phi}{\basep{\Man{}}}\)}{\avec{}}}{\func{\(\tanspaceiso{\func{\basep{\Man{}}}{\avec{}}}{\Man{}}{\phi}\)}{\avec{}}}.
\end{aligned}\right.
\end{align}
The $\difclass{\infty}$ maximal-atlas on $\tanbun{\Man{}}$ corresponded to the atlas $\defSet{\tanchart{\Man{}}{\phi}}{\phi\in\maxatlas{}}$
is denoted by $\tanatlas{\Man{}}$ and $\Tanbun{\Man{}}=\opair{\tanbun{\Man{}}}{\tanatlas{\Man{}}}$. For every $\phi$ in $\maxatlas{}$,
$\tanchart{\Man{}}{\phi}$ is referred to as the $\quotl$tangent-bundle chart of $\Man{}$ associated with $\phi$$\quotr$.
\item[\myitem{DG~10.}]
The set of all smooth ($\difclass{\infty}$) vector-fields on $\Man{}$ is denoted by $\vecf{\Man{}}{\infty}$.
In other words, $\vecf{\Man{}}{\infty}$ is the set of all maps $\avecf{}$ in $\mapdifclass{\infty}{\Man{}}{\Tanbun{\Man{}}}$
such that $\cmp{\basep{\Man{}}}{\avecf{}}=\identity{\M{}}$.
$\Vecf{\Man{}}{\infty}$ denotes the vector-space obtained from the set $\vecf{\Man{}}{\infty}$ endowed with its canonical linear-structure.
\item[\myitem{DG~11.}]
The set of all $\infty$-derivations on $\Man{}$ is denoted by $\Derivation{\Man{}}{\infty}$. In other words,
$\Derivation{\Man{}}{\infty}$ is the set of all linear maps
$\aderivation{}$ in $\Lin{\Lmapdifclass{\infty}{\Man{}}{\RR}}{\Lmapdifclass{\infty}{\Man{}}{\RR}}$ such that
for every pair $\cf$ and $\cg$ of elements of $\mapdifclass{\infty}{\Man{}}{\RR}$,
$\func{\aderivation{}}{\cf\rdot\cg}=\[\func{\aderivation{}}{\cf}\]\rdot\cg+\cf\rdot\[\func{\aderivation{}}{\cg}\]$.
$\LDerivation{\Man{}}{\infty}$ denotes the vector-space obtained from the set $\Derivation{\Man{}}{\infty}$ endowed with its
canonical linear-structure which is inherited from that of $\Lin{\Lmapdifclass{\infty}{\Man{}}{\RR}}{\Lmapdifclass{\infty}{\Man{}}{\RR}}$.
\item[\myitem{DG~12.}]
$\derop{\Man{}}{\Man{1}}$ denotes the differential operator on the set $\mapdifclass{\infty}{\Man{}}{\Man{1}}$. That is,
$\derop{\Man{}}{\Man{1}}$ is defined to be the element of
$\Func{\mapdifclass{\infty}{\Man{}}{\Man{1}}}{\mapdifclass{\infty}{\Tanbun{\Man{}}}{\Tanbun{\Man{1}}}}$ such that
for every map $\cf$ in $\mapdifclass{\infty}{\Man{}}{\Man{1}}$, every $\avec{}$ in $\tanbun{\Man{}}$,
and for each chart $\phi$ in $\maxatlas{}$ whose domain contains $\func{\basep{\Man{}}}{\avec{}}$ and
each chart $\psi$ in $\maxatlas{1}$ whose domain contains $\func{\cf}{\func{\basep{\Man{}}}{\avec{}}}$,
\begin{equation}\label{eqdefinitionofdifferentialofamap}
\func{\[\der{\cf}{\Man{}}{\Man{1}}\]}{\avec{}}=
\func{\bigg[\cmp{\finv{\(\tanspaceiso{\func{\cf}{\func{\basep{\Man{}}}{\avec{}}}}{\Man{1}}{\psi}\)}}{\cmp{\bigg(\func{\[\banachder{\(\cmp{\psi}{\cmp{f}{\finv{\phi}}}\)}{\R^n}{\R^m}\]}
{\func{\phi}{\func{\basep{\Man{}}}{\avec{}}}}\bigg)}{\tanspaceiso{\func{\basep{\Man{}}}{\avec{}}}{\Man{}}{\phi}}}\bigg]}{\avec{}}.
\end{equation}
It can be easily seen that for a smooth map $\cf$ in $\mapdifclass{\infty}{\Man{}}{\Man{1}}$,
$\der{\cf}{\Man{}}{\Man{1}}$ maps $\tanspace{\point}{\Man{}}$ into $\tanspace{\func{\cf}{\point}}{\Man{1}}$
for every point $\point$ of $\Man{}$. In other words,
\begin{equation}\label{eqtangentmapbasepoint}
\Foreach{\avec{}}{\tanbun{\Man{}}}
\func{\basep{\Man{}}}{\func{\[\der{\cf}{\Man{}}{\Man{1}}\]}{\avec{}}}=
\func{\[\cmp{\cf}{\basep{\Man{}}}\]}{\avec{}}.
\end{equation}
Furthermore, for every $\cf$ in $\mapdifclass{\infty}{\Man{}}{\Man{1}}$ and each point $\point$ of $\Man{}$,
the restriction of $\der{\cf}{\Man{}}{\Man{1}}$ to $\tanspace{\point}{\Man{}}$ is a linear map from
$\Tanspace{\point}{\Man{}}$ to $\Tanspace{\func{\cf}{\point}}{\Man{1}}$, that is,
\begin{equation}
\Foreach{\point}{\M{}}
\func{\res{\der{\cf}{\Man{}}{\Man{1}}}}{\tanspace{\point}{\Man{}}}\in\Lin{\Tanspace{\point}{\Man{}}}{\Tanspace{\func{\cf}{\point}}{\Man{1}}}.
\end{equation}
The differential operator has local behavior. This means, for any $\cf$ and $\cg$ in $\mapdifclass{\infty}{\Man{}}{\RR}$
that coincide in an open set $\U$ of $\mantops{\Man{}}$ (an element of $\mantop{\Man{}}$),
$\der{\cf}{\Man{}}{\Man{1}}$ and $\der{\cg}{\Man{}}{\Man{1}}$ coincide in $\func{\pimage{\basep{\Man{}}}}{\U}$.\\
The chain rule of differentiation (in the category of smooth manifolds) asserts that for every $\cf$ in
$\mapdifclass{\infty}{\Man{}}{\Man{1}}$ and every $\cg$ in $\mapdifclass{\infty}{\Man{1}}{\Man{2}}$,
\begin{equation}\label{eqchainrule}
\der{\(\cmp{\cg}{\cf}\)}{\Man{}}{\Man{2}}=\cmp{\(\der{\cg}{\Man{1}}{\Man{2}}\)}{\(\der{\cf}{\Man{}}{\Man{1}}\)}.
\end{equation}
So since for every $\cf$ in $\Diffeo{\infty}{\Man{}}{\Man{1}}$, $\cmp{\finv{\cf}}{\cf}=\identity{\M{}}$, and
$\der{\identity{\M{}}}{\Man{}}{\Man{}}=\identity{\tanbun{\Man{}}}$, clearly,
\begin{equation}\label{eqdiffeomorphismdifrule}
\Foreach{\cf}{\Diffeo{\infty}{\Man{}}{\Man{1}}}
\finv{\(\der{\cf}{\Man{}}{\Man{1}}\)}=\der{\finv{\cf}}{\Man{1}}{\Man{}}.
\end{equation}
\\$\der{\cf}{\Man{}}{\Man{1}}$ can simply be denoted by $\derr{\cf}{}{}$ when there is no ambiguity about the underlying
source and target manifolds $\Man{}$ and $\Man{1}$.
\item[\myitem{DG~13.}]
$\Rderop{\Man{}}$ denotes the derivative operator on $\mapdifclass{\infty}{\Man{}}{\RR}$. That is,
$\Rderop{\Man{}}$ is defined to be the element of
$\Func{\mapdifclass{\infty}{\Man{}}{\RR}}{\mapdifclass{\infty}{\Tanbun{\Man{}}}{\RR}}$ such that
for every smooth map $\cf$ in $\mapdifclass{\infty}{\Man{}}{\RR}$, every $\avec{}$ in $\tanbun{\Man{}}$,
\begin{align}\label{eqdefderivativeoperator}
\func{\[\Rder{\cf}{\Man{}}\]}{\avec{}}:=&\func{\tanspaceiso{\func{\cf}{\func{\basep{\Man{}}}{\avec{}}}}{\RR}
{\identity{\R}}}{\func{\[\der{\cf}{\Man{}}{\RR}\]}{\avec{}}}\cr
=&\func{\[\cmp{\bigg(\func{\[\banachder{\(\cmp{f}{\finv{\phi}}\)}{\R^n}{\R}\]}
{\func{\phi}{\func{\basep{\Man{}}}{\avec{}}}}\bigg)}{\tanspaceiso{\func{\basep{\Man{}}}{\avec{}}}{\Man{}}{\phi}}\]}{\avec{}},
\end{align}
where $\phi$ can be any of charts in $\maxatlas{}$ whose domain contains $\func{\basep{\Man{}}}{\point}$.
For every $\cf$ in $\mapdifclass{\infty}{\Man{}}{\RR}$, $\Rderop{\Man{}}$
can be regarded as a replacement for the differential of $\cf$.
Equipping $\mapdifclass{\infty}{\Man{}}{\RR}$ and $\mapdifclass{\infty}{\Tanbun{\Man{}}}{\RR}$ with
their natural linear-structure, $\Rderop{\Man{}}$ is a linear map as a direct result of a similar property for
the derivative operator of real-valued smooth maps on Banach-spaces. That is,
$\Rderop{\Man{}}\in\Lin{\mapdifclass{\infty}{\Man{}}{\RR}}{\mapdifclass{\infty}{\Tanbun{\Man{}}}{\RR}}$.
As another crucial property of $\Rderop{\Man{}}$, for every $\cf$ and $\cg$ in $\mapdifclass{\infty}{\Man{}}{\RR}$,
and every $\avec{}$ in $\tanbun{\Man{}}$,
\begin{equation}\label{Jacobipropertyofderivative}
\func{\[\func{\Rderop{\Man{}}}{\cf\rdot\cg}\]}{\avec{}}=\[\func{\(\cmp{\cg}{\basep{\Man{}}}\)}{\avec{}}\]\[\func{\Rder{\cf}{\Man{}}}{\avec{}}\]+
\[\func{\(\cmp{\cf}{\basep{\Man{}}}\)}{\avec{}}\]\[\func{\Rder{\cg}{\Man{}}}{\avec{}}\].
\end{equation}
As a trivial consequence of the general case,
for every $\cf$ in $\mapdifclass{\infty}{\Man{}}{\RR}$ and each point $\point$ of $\Man{}$,
the restriction of $\Rder{\cf}{\Man{}}$ to $\tanspace{\point}{\Man{}}$ is a linear map from
$\Tanspace{\point}{\Man{}}$ to $\R$, that is,
$\func{\res{\Rder{\cf}{\Man{}}}}{\tanspace{\point}{\Man{}}}\in\Lin{\Tanspace{\point}{\Man{}}}{\R}$.
\\According to \Ref{eqtangentmapbasepoint}, \Ref{eqchainrule}, and \Ref{eqdefderivativeoperator},
it is evident that,
\begin{equation}\label{eqspecialchainrule1}
\Foreach{\cf}{\mapdifclass{\infty}{\Man{1}}{\RR}}
\Foreach{\cg}{\mapdifclass{\infty}{\Man{}}{\Man{1}}}
\Rder{\(\cmp{\cf}{\cg}\)}{\Man{}}=\cmp{\(\Rder{\cf}{\Man{1}}\)}{\(\der{\cg}{\Man{}}{\Man{1}}\)}.
\end{equation}
This is called the $\quotl$special chain-rule of differentiation$\quotr$, here.
$\Rder{\cf}{\Man{}}$ can simply be denoted by $\Rder{\cf}{}$ when there is no ambiguity about the underlying
manifold $\Man{}$, and it can be simply referred to as the $\quotl$derivative of $\cf$$\quotr$. $\func{\Rder{\cf}{\Man{}}}{\avec{}}$
can simply be referred to as the $\quotl$derivative of $\cf$ in the direction $\avec{}$$\quotr$.
\item[\myitem{DG~14.}]
The set $\atlas{}:=\defSet{\funcprod{\phi}{\psi}}{\phi\in\maxatlas{},~\psi\in\maxatlas{1}}$
is an element of $\atlases{\infty}{\Cprod{\M{}}{\M{1}}}{\R^{n+m}}$, that is, a $\difclass{\infty}$ atlas on $\Cprod{\M{}}{\M{1}}$
constructed upon $\R^{n+m}$. The topology on $\Cprod{\M{}}{\M{1}}$ induced by the $\difclass{\infty}$ maximal-atlas on
$\Cprod{\M{}}{\M{1}}$ constructed upon $\R^{n+m}$ that is generated by the atlas $\atlas{}$
(that is, $\func{\maxatlasgen{\infty}{\Cprod{\M{}}{\M{1}}}{\R^{n+m}}}{\atlas{}}$),
is the topological-product of the underlying topological-spaces of the manifolds $\Man{}$ and $\Man{1}$, that is,
$\topprod{\mantops{\Man{}}}{\mantops{\Man{1}}}$. Since $\mantops{\Man{}}$ and $\mantops{\Man{1}}$
are Hausdorff and second-countable spaces, so is the topological space $\topprod{\mantops{\Man{}}}{\mantops{\Man{1}}}$.
Consequently, the set $\Cprod{\M{}}{\M{1}}$ endowed with the $\difclass{\infty}$ maximal-atlas
$\func{\maxatlasgen{\infty}{\Cprod{\M{}}{\M{1}}}{\R^{n+m}}}{\atlas{}}$ is a differentiable-structure
whose underlying topological-space is Hausdorff and second-countable. Thus,
$\opair{\Cprod{\M{}}{\M{1}}}{\func{\maxatlasgen{\infty}{\Cprod{\M{}}{\M{1}}}{\R^{n+m}}}{\atlas{}}}$ is a manifold.
$\manprod{\Man{}}{\Man{1}}$ denotes this manifold and is referred to as
the $\quotl$manifold-product of the manifolds $\Man{}$ and $\Man{1}$$\quotr$.
\item[\myitem{DG~15.}]
For every $\cf$ in $\mapdifclass{\infty}{\Man{}}{\Man{2}}$ and every $\cg$ in $\mapdifclass{\infty}{\Man{1}}{\Man{3}}$,
the product mapping $\funcprod{\cf}{\cg}$ is a smooth map from the manifold-product of $\Man{}$ and $\Man{1}$ to the
manifold-product of $\Man{2}$ and $\Man{3}$, that is,
$\(\funcprod{\cf}{\cg}\)\in\mapdifclass{\infty}{\manprod{\Man{}}{\Man{1}}}{\manprod{\Man{2}}{\Man{3}}}$.
\item[\myitem{DG~16.}]
Every smooth map $\cf$ from $\Man{}$ to $\Man{1}$ whose differential becomes a a linear monomorphism from
$\tanspace{\point}{\Man{}}$ to $\tanspace{\func{\cf}{\point}}{\Man{1}}$ when restricted to $\tanspace{\point}{\Man{}}$
for every point $\point$ of $\Man{}$,
is referred to as a $\quotl$$\infty$-immersion from the manifold $\Man{}$ to the manifold $\Man{1}$$\quotr$.
The set of all $\infty$-immersions from $\Man{}$ to $\Man{1}$ is denoted by $\immersion{\infty}{\Man{}}{\Man{1}}$.\\
Every $\infty$-immersion from $\Man{}$ to $\Man{1}$ which is additionally a topological embedding of
$\mantops{\Man{}}$ (the underlying topological-space of $\Man{}$) into $\mantops{\Man{1}}$ (the underlying
topological-space of $\Man{1}$), is referred to as an $\quotl$$\infty$-embedding of $\Man{}$ into $\Man{1}$$\quotr$.
The set of all $\infty$-embedding from $\Man{}$ to $\Man{1}$ is denoted by $\embedding{\infty}{\Man{}}{\Man{1}}$.\\
Every smooth map $\cf$ from $\Man{}$ to $\Man{1}$ whose differential becomes a a linear epimorphism from
$\tanspace{\point}{\Man{}}$ and $\tanspace{\func{\cf}{\point}}{\Man{1}}$ when restricted to $\tanspace{\point}{\Man{}}$
for every point $\point$ of $\Man{}$,
is referred to as a $\quotl$$\infty$-submersion from the manifold $\Man{}$ to the manifold $\Man{1}$$\quotr$.
The set of all $\infty$-submersions from $\Man{}$ to $\Man{1}$ is denoted by $\submersion{\infty}{\Man{}}{\Man{1}}$.
\item[\myitem{DG~17.}]
$\Man{1}$ is said to be an $\quotl$($m_{1}$-dimensional) $\infty$-immersed submanifold of $\Man{}$$\quotr$ or
simply an $\quotl$immersed submanifold of $\Man{}$$\quotr$ iff $\M{1}\subseteq\M{}$ and
the injection of $\M{1}$ into $\M{}$ is an $\infty$-immersion from $\Man{1}$ to $\Man{}$, that is
$\Injection{\M{1}}{\M{}}\in\immersion{\infty}{\Man{1}}{\Man{}}$.\\
$\Man{1}$ is said to be an $\quotl$($m_{1}$-dimensional) $\infty$-embedded submanifold of $\Man{}$$\quotr$ or
simply an $\quotl$embedded submanifold of $\Man{}$$\quotr$ iff $\M{1}\subseteq\M{}$ and
the injection of $\M{1}$ into $\M{}$ is an $\infty$-embedding from $\Man{1}$ to $\Man{}$, that is
$\Injection{\M{1}}{\M{}}\in\embedding{\infty}{\Man{1}}{\Man{}}$. It is clear that every
embedded submanifold of $\Man{}$ is necessarily an immersed submanifold of $\Man{}$.\\
given a subset $\SET{}$ of $\M{}$, there is at most one corresponded differentiable structure modeled on $\R^{q}$ for some $q\in\Zp$
that makes it a $\difclass{\infty}$ manifold which is an embedded submanifold of $\Man{}$;
If this unique differentiable structure exists, it is inherited from the maximal atlas of $\Man{}$ in a
canonical way. The set of all subsets of $\M{}$ accepting a differentiable structure modeled on a Euclidean-space
that makes it an embedded submanifold of $\Man{}$ is denoted by $\Emsubman{\Man{}}$.
For every $\SET{}$ in $\Emsubman{\Man{}}$, the unique $\difclass{\infty}$
manifold with the set of points $\SET{}$ that is an embedded submanifold of $\Man{}$ is denoted by $\emsubman{\Man{}}{\SET{}}$.
Each element of $\Emsubman{\Man{}}$ is called a $\quotl$embedded set of the manifold $\Man{}$$\quotr$.
Every open set of $\Man{}$ is an embedded set of $\Man{}$, that is $\mantop{\Man{}}\subseteq\Emsubman{\Man{}}$.
Furthermore, every embedded set of $\Man{}$ is a locally-closed set of the underlying topological-space of $\Man{1}$,
which means every embedded set of $\Man{}$ is open in its closure (with respect to its topology inherited from the topology of $\Man{}$).\\
If $\Man{2}$ is an immersed submanifold of $\Man{1}$ and $\Man{1}$ is an immersed submanifold of $\Man{}$, then
$\Man{2}$ is an immersed submanifold of $\Man{}$. Additionally, given subsets $\SET{1}$ and $\SET{2}$ of $\M{}$
such that $\SET{2}\subseteq\SET{1}$, if $\SET{1}\in\Emsubman{\Man{}}$ and $\SET{2}\in\Emsubman{\emsubman{\Man{}}{\SET{1}}}$,
then $\SET{2}\in\Emsubman{\Man{}}$.
\item[\myitem{DG~18.}]
For every smooth map $\cf$ from $\Man{}$ to $\Man{1}$, and every element $\SET{}$ of $\Emsubman{\Man{}}$ (every embedded set of $\Man{}$),
the domain-restriction of $\cf$ to $\SET{}$ is also a smooth map from the embedded submanifold $\emsubman{\Man{}}{\SET{}}$ of
$\Man{}$ to $\Man{1}$.\\
For every map $\cf$ from $\Man{}$ to $\Man{1}$, and every element $\SET{}$ of $\Emsubman{\Man{1}}$ (every embedded set of $\Man{1}$)
that includes the image of $\cf$,
$\cf$ is a smooth map from $\Man{}$ to $\Man{1}$ if and only if
the codomain-restriction of $\cf$ to $\SET{}$ is a smooth map from the $\Man{}$ to the embedded submanifold
$\emsubman{\Man{1}}{\SET{}}$ of $\Man{1}$.\\
The embedded sets of the manifold $\Man{}$ is in a one-to-one correspondence with the embedded sets of
any manifold diffeomorphic to $\Man{}$. Actually,
Every $\infty$-diffeomorphism $\cf$ from $\Man{}$ to $\Man{1}$ maps an embedded set of $\Man{}$ to an embedded set of
$\Man{1}$. Therefeore the restriction (that is simultaneously domain-restriction and codomain-restriction) of
a smooth map $\cf$ from $\Man{}$ to $\Man{1}$ to an embedded set $\SET{}$ of $\Man{}$ is again a smooth map
from the embedded submanifold $\emsubman{\Man{}}{\SET{}}$ of $\Man{}$ to the embedded submanifold
$\emsubman{\Man{1}}{\func{\image{\cf}}{\SET{}}}$ of $\Man{1}$.
\item[\myitem{DG~19.}]
Given an $\infty$-immersion $\cf$ from $\Man{}$ to $\Man{1}$, for every point $\point$ of $\Man{}$,
there exists an open set $\U$ of the manifold $\Man{}$ such that $\func{\image{\cf}}{\U}$ is an embedded
set of the manifold $\Man{1}$, and the restriction of $\cf$ to $\U$ is an $\infty$-diffeomorphism from
the embedded submanifold $\emsubman{\Man{}}{\U}$ to the embedded submanifold $\emsubman{\Man{1}}{\func{\image{\cf}}{\U}}$
of $\Man{1}$.
\item[\myitem{DG~20.}]
Given a positive integer $n$, and a non-empty open set $\U$ of $\topR{\R^n}$, trivially $\U$ is an embedded
set of the canonical smooth manifold of $\R^n$, and the manifold $\emsubman{\RR^n}{\U}$ is simply denoted by $\Ropenman{\U}{n}$,
and is called the $\quotl$canonical differentiable structure of $\U$ inherited from $\RR^n$$\quotr$.
It is also a trivial fact that $\Ropenman{\U}{n}$ is a $\difclass{\infty}$ differentiable structure modeled on $\R^n$, the singleton
$\seta{\identity{\U}}$ being an atlas and consequently $\identity{\U}$ a chart of it.
\item[\myitem{DG~21.}]
The notion of partial derivative of differentiable mappings from a product of Banach spaces to a Banach space
as introduced in \cite{Cartan}, can be generalized to the case of smooth mappings from a product of manifolds
to a manifold. Here, only the simplest case is considered, where the domain of a smooth mapping is the product of
a pair of manifolds.\\
$\point_1$ and $\point_2$ are taken as points of the manifolds $\Man{1}$ and $\Man{2}$, respectively.
The mappings $\function{\leftparinj{\M{1}}{\M{2}}{\point_2}}{\M{1}}{\Cprod{\M{1}}{\M{2}}}$ and
$\function{\rightparinj{\M{1}}{\M{2}}{\point_2}}{\M{2}}{\Cprod{\M{1}}{\M{2}}}$ are defined as,
\begin{align}
&\Foreach{\x}{\M{1}}
\func{\leftparinj{\M{1}}{\M{2}}{\point_2}}{\x}\eqdef\opair{\x}{\point_2},\\
&\Foreach{\x}{\M{2}}
\func{\rightparinj{\M{1}}{\M{2}}{\point_1}}{\x}\eqdef\opair{\point_1}{\x}.
\end{align}
It is a trivial fact that $\leftparinj{\M{1}}{\M{2}}{\point_2}$ is a smooth map from $\Man{1}$ to
$\manprod{\Man{1}}{\Man{2}}$, and $\rightparinj{\M{1}}{\M{2}}{\point_1}$ a smooth map
from $\Man{2}$ to $\manprod{\Man{1}}{\Man{2}}$. The mapping $\prodmantan{\Man{1}}{\Man{2}}{\point_1}{\point_2}$
is defined as,
\begin{align}
&\prodmantan{\Man{1}}{\Man{2}}{\point_1}{\point_2}\indef\Func{\Cprod{\tanspace{\point_1}{\Man{1}}}{\tanspace{\point_2}{\Man{2}}}}
{\tanspace{\opair{\point_1}{\point_2}}{\manprod{\Man{1}}{\Man{2}}}},\cr
&\begin{aligned}
&\Foreach{\opair{\vv{1}}{\vv{2}}}{\Cprod{\tanspace{\point_1}{\Man{1}}}{\tanspace{\point_2}{\Man{2}}}}\cr
&\func{\prodmantan{\Man{1}}{\Man{2}}{\point_1}{\point_2}}{\binary{\vv{1}}{\vv{2}}}\eqdef
\func{\[\der{\leftparinj{\M{1}}{\M{2}}{\point_2}}{\Man{1}}{\manprod{\Man{1}}{\Man{2}}}\]}{\vv{1}}+
\func{\[\der{\rightparinj{\M{1}}{\M{2}}{\point_1}}{\Man{2}}{\manprod{\Man{1}}{\Man{2}}}\]}{\vv{2}}.
\end{aligned}\cr
&{}
\end{align}
Actually, $\prodmantan{\Man{1}}{\Man{2}}{\point_1}{\point_2}$ is a bijection and hence the tangent-space of $\manprod{\Man{1}}{\Man{2}}$
at the point $\opair{\point_1}{\point_2}$ can be identified with $\Cprod{\tanspace{\point_1}{\Man{1}}}{\tanspace{\point_2}{\Man{2}}}$
via $\prodmantan{\Man{1}}{\Man{2}}{\point_1}{\point_2}$.
\\$\cf$ is taken as a smooth map from the manifold-product of $\Man{1}$ and $\Man{2}$ to the manifold $\Man{}$, that is an element of
$\mapdifclass{\infty}{\manprod{\Man{1}}{\Man{2}}}{\Man{}}$.
\begin{align}\label{EQpartialsLemma}
&\Foreach{\opair{\vv{1}}{\vv{2}}}{\Cprod{\tanspace{\point_1}{\Man{1}}}{\tanspace{\point_2}{\Man{2}}}}\cr
&\begin{aligned}
\func{\[\der{\cf}{\manprod{\Man{1}}{\Man{2}}}{\Man{}}\]}{\func{\prodmantan{\Man{1}}{\Man{2}}{\point_1}{\point_2}}{\binary{\vv{1}}{\vv{2}}}}=
&\hskip0.5\baselineskip\func{\[\der{\(\cmp{\cf}{\leftparinj{\M{1}}{\M{2}}{\point_2}}\)}{\Man{1}}{\Man{}}\]}{\vv{1}}\cr
&+\func{\[\der{\(\cmp{\cf}{\rightparinj{\M{1}}{\M{2}}{\point_1}}\)}{\Man{2}}{\Man{}}\]}{\vv{2}},
\end{aligned}
\end{align}
where the addition takes place in the tangent-space of $\manprod{\Man{1}}{\Man{2}}$ at $\func{\cf}{\binary{\point_1}{\point_2}}$.
\end{itemize}
\section{Flow of a Smooth Vector-Field}
\textit{The theorems of this section are not provided with proofs,
and considered to be standard elementary facts of differential geometry.}
\definition
The set of all connected and open sets of $\topR{\R}$ (the canonical topological-space of $\R$) containing $0$
is denoted by $\ointR$, which is known to be the set of all open intervals of $\R$ containing $0$.
\begin{align}
\ointR:=
&\func{\nei{\topR{\R}}}{\seta{0}}\cap\connecteds{\topR{\R}}\cr
&\begin{aligned}
&=\hskip0.5\baselineskip\defsets{\aninterval{}}{\R}{0\in\aninterval{},~\bigg(\Exists{a}{\R}\Exists{b}{\R}a<b,~I=\Ointerval{a}{b}\bigg)}\cr
&\hskip0.5\baselineskip\cup\defsets{\aninterval{}}{\R}{0\in\aninterval{},~\bigg(\Exists{a}{\R}I=\Ointerval{a}{+\infty}\bigg)}\cr
&\hskip0.5\baselineskip\cup\defsets{\aninterval{}}{\R}{0\in\aninterval{},~\bigg(\Exists{a}{\R}I=\Ointerval{-\infty}{a}\bigg)}\cr
&\hskip0.5\baselineskip\cup\seta{\R}.
\end{aligned}
\end{align}
\endef
\definition\label{defidentityvectorofthetangentspaceofreals}
For every open set $\aninterval{}$ of $\topR{\R}$, and every $\x$ in $\aninterval{}$,
the vector of the tangent-space of $\Ropenman{\aninterval{}}{}$ at the point $\x$ corresponded to $1$ with respect to the chart
$\identity{\aninterval{}}$ is denoted by $\Rtanidentity{\aninterval{}}{\x}$.
\begin{equation}
\Foreach{\aninterval{}}{\topologyofspace{\topR{\R}}}
\Foreach{\x}{\aninterval{}}
\Rtanidentity{\aninterval{}}{\x}:=\func{\finv{\[\tanspaceiso{\x}{\Ropenman{\aninterval{}}{}}{\identity{\aninterval{}}}\]}}{1}.
\end{equation}
\endef
\definition\label{defintegralcurvesofasmoothvectorfield}
$\avecf{}$ is taken as an element of $\vecf{\Man{}}{\infty}$ (a smooth vector-field on the manifold $\Man{}$),
and $\point$ as a point of the manifold $\Man{}$.
Given an open interval $\aninterval{}$ of $\R$ containing $0$, a mapping $\function{\curve{}}{\aninterval{}}{\M{}}$ is
called an $\quotl$integral-curve of the smooth vector-field $\avecf{}$ on $\Man{}$ with the initial condition $\point$$\quotr$ iff
$\curve{}$ is a smooth map from the manifold $\Ropenman{\aninterval{}}{}$ to the manifold $\Man{}$, and,
\begin{align}
\left\{
\begin{aligned}
&\func{\curve{}}{0}=\point,\\
&\Foreach{t}{\aninterval{}}
\func{\[\der{\curve{}}{\Ropenman{\aninterval{}}{}}{\Man{}}\]}{\Rtanidentity{\aninterval{}}{t}}=\func{\(\cmp{\avecf{}}{\curve{}}\)}{t}.
\end{aligned}\right.
\end{align}
When there is no ambiguity in the context, $\func{\[\der{\curve{}}{\Ropenman{\aninterval{}}{}}{\Man{}}\]}{\Rtanidentity{\aninterval{}}{t}}$
can simply be denoted by $\func{\dot{\curve{}}}{t}$ and called the $\quotl$velocity of the curve $\curve{}$ at time $t$$\quotr$.\\
The set of all integral-curves of the smooth vector-field $\avecf{}$ of $\Man{}$ with the initial condition $\point$
is denoted bt $\func{\integralcurves{\Man{}}}{\binary{\avecf{}}{\point}}$. That is,
\begin{align}
&\hskip0.8\baselineskip\func{\integralcurves{\Man{}}}{\binary{\avecf{}}{\point}}\cr
&:=\Union{\aninterval{}}{\ointR}
{\defset{\curve{}}{\mapdifclass{\infty}{\Ropenman{\aninterval{}}{}}{\Man{}}}{\func{\curve{}}{0}=\point,~
\bigg(\Foreach{t}{\aninterval{}}
\func{\[\der{\curve{}}{\Ropenman{\aninterval{}}{}}{\Man{}}\]}{\Rtanidentity{\aninterval{}}{t}}=
\func{\(\cmp{\avecf{}}{\curve{}}\)}{t}\bigg)}}.\cr
&{}
\end{align}
\endef
\theorem
$\avecf{}$ is taken as an element of $\vecf{\Man{}}{\infty}$ (a smooth vector-field on the manifold $\Man{}$).
For every point $\point$ of $\Man{}$, there exists at least one integral-curve of $\avecf{}$
with the initial condition $\point$ (and hence an uncountably infinite collection of such integral-curves). That is,
\begin{equation}
\Foreach{\point}{\M{}}
\func{\integralcurves{\Man{}}}{\binary{\avecf{}}{\point}}\neq\empty.
\end{equation}
\endthm
\theorem\label{thmintegralcurvescoincideinthedomainintersection}
$\point$ is taken as a point of $\Man{}$. Every pair of integral-curves of $\avecf{}$ with the initial condition $\point$,
coincide in the intersection of their domains. That is,
\begin{equation}
\Foreach{\opair{\curve{1}}{\curve{2}}}{\func{\integralcurves{\Man{}}}{\binary{\avecf{}}{\point}}^{\times 2}}
\Foreach{t}{\domain{\curve{1}}\cap\domain{\curve{2}}}\func{\curve{1}}{t}=\func{\curve{2}}{t}.
\end{equation}
\endthm
\definition\label{deftimeintervalofmaximalintegralcurve}
$\avecf{}$ is taken as an element of $\vecf{\Man{}}{\infty}$ (a smooth vector-field on the manifold $\Man{}$),
and $\point$ as a point of the manifold $\Man{}$.
The union of the domains of all integral-curves of $\avecf{}$ with the initial condition $\point$
is denoted by $\func{\maxinterval{\Man{}}}{\binary{\avecf{}}{\point}}$.
\begin{equation}
\func{\maxinterval{\Man{}}}{\binary{\avecf{}}{\point}}:=
\union{\defSet{\domain{\curve{}}}{\curve{}\in\func{\integralcurves{\Man{}}}{\binary{\avecf{}}{\point}}}}.
\end{equation}
Furthermore, the mapping $\func{\maxintcurve{\Man{}}}{\binary{\avecf{}}{\point}}$ is defined as,
\begin{align}
&\func{\maxintcurve{\Man{}}}{\binary{\avecf{}}{\point}}\indef
\Func{\func{\maxinterval{\Man{}}}{\binary{\avecf{}}{\point}}}{\M{}},\cr
&\Foreach{\curve{}}{\func{\integralcurves{\Man{}}}{\binary{\avecf{}}{\point}}}
\Foreach{t}{\domain{\curve{}}}
\func{\[\func{\maxintcurve{\Man{}}}{\binary{\avecf{}}{\point}}\]}{t}\eqdef\func{\curve{}}{t}.
\end{align}
which is referred to as the $\quotl$maximal integral-curve of the smooth vector-field $\avecf{}$ on $\Man{}$ with the
initial condition $\point$$\quotr$.
\endef
\corollary\label{corsubcurvesofmaximalintegralcurve}
$\avecf{}$ is taken as an element of $\vecf{\Man{}}{\infty}$ (a smooth vector-field on the manifold $\Man{}$),
and $\point$ as a point of the manifold $\Man{}$.
\begin{itemize}
\item[$\centerdot$]
Given an integral-curve $\curve{}$ of $\avecf{}$ with the initial condition $\point$, the domain-restriction of
$\curve{}$ to every open interval of $\topR{\R}$ containing $0$ and included in $\domain{\curve{}}$ is again
an integral-curve $\curve{}$ of $\avecf{}$ with the initial condition $\point$.
\begin{align}
&\Foreach{\curve{}}{\func{\integralcurves{\Man{}}}{\binary{\avecf{}}{\point}}}
\Foreach{\aninterval{}}{\ointR\cap\CSs{\domain{\curve{}}}}\cr
&\func{\resd{\curve{}}}{\aninterval{}}\in\func{\integralcurves{\Man{}}}{\binary{\avecf{}}{\point}}.
\end{align}
Consequently, by considering the fact that there exists at least one integral-curve of $\avecf{}$
with the initial condition $\point$, it becomes evident that there exists a positive real number $\varepsilon$
and an integral-curve of $\avecf{}$ with the initial condition $\point$ whose domain is $\Ointerval{-\varepsilon}{+\varepsilon}$.
\begin{equation}
\Exists{\varepsilon}{\Rp}\Exists{\curve{}}{\func{\integralcurves{\Man{}}}{\binary{\avecf{}}{\point}}}
\domain{\curve{}}=\Ointerval{-\varepsilon}{+\varepsilon}.
\end{equation}
\item[$\centerdot$]
The domain of every integral-curve of $\avecf{}$ with the initial condition
$\point$ is a subset of the domain of maximal integral-curve of $\avecf{}$ with the initial condition $\point$.
\begin{equation}
\Foreach{\curve{}}{\func{\integralcurves{\Man{}}}{\binary{\avecf{}}{\point}}}
\domain{\curve{}}\subseteq\func{\maxinterval{\Man{}}}{\binary{\avecf{}}{\point}}.
\end{equation}
\item[$\centerdot$]
The domain-restriction of $\func{\maxintcurve{\Man{}}}{\binary{\avecf{}}{\point}}$ to the domain of every
integral-curve $\curve{}$ of $\avecf{}$ with the initial condition $\point$, coincides with $\curve{}$. That is,
\begin{equation}
\Foreach{\curve{}}{\func{\integralcurves{\Man{}}}{\binary{\avecf{}}{\point}}}
\func{\resd{\func{\maxintcurve{\Man{}}}{\binary{\avecf{}}{\point}}}}{\domain{\curve{}}}=\curve{}.
\end{equation}
\end{itemize}
\endcor
\theorem\label{thmmaximalintegralcurveisanintegralcurve}
$\avecf{}$ is taken as an element of $\vecf{\Man{}}{\infty}$ (a smooth vector-field on the manifold $\Man{}$).
For every point $\point$ of $\Man{}$,
$\func{\maxintcurve{\Man{}}}{\binary{\avecf{}}{\point}}$ is an integral-curve of $\avecf{}$ with the initial condition $\point$.
\begin{equation}
\Foreach{\point}{\M{}}
\func{\maxintcurve{\Man{}}}{\binary{\avecf{}}{\point}}\in
\func{\integralcurves{\Man{}}}{\binary{\avecf{}}{\point}}.
\end{equation}
\endthm
\definition
$\avecf{}$ is taken as an element of $\vecf{\Man{}}{\infty}$ (a smooth vector-field on the manifold $\Man{}$).
$\avecf{}$ is called a $\quotl$complete (smooth) vector-field on $\Man{}$$\quotr$ iff,
\begin{equation}
\Foreach{\point}{\M{}}
\func{\maxinterval{\Man{}}}{\binary{\avecf{}}{\point}}=\R.
\end{equation}
\endef
\theorem
If the underlying topological-space $\mantops{\Man{}}$ of the manifold $\Man{}$ is compact, then every
smooth vector-field of $\Man{}$ is complete.
\endthm
\definition\label{defsmoothvectorfieldflow}
$\avecf{}$ is taken as an element of $\vecf{\Man{}}{\infty}$ (a smooth vector-field on the manifold $\Man{}$).
The mapping $\vfflow{\Man{}}{\avecf{}}$ is defined as,
\begin{align}
&\vfflow{\Man{}}{\avecf{}}\indef\Func{\Union{\point}{\M{}}
\[\Cprod{\func{\maxinterval{\Man{}}}{\binary{\avecf{}}{\point}}}{\seta{\point}}\]}{\M{}},\cr
&\Foreach{\point}{\M{}}
\Foreach{t}{\func{\maxinterval{\Man{}}}{\binary{\avecf{}}{\point}}}
\func{\vfflow{\Man{}}{\avecf{}}}{\binary{t}{\point}}\eqdef
\func{\[\func{\maxintcurve{\Man{}}}{\binary{\avecf{}}{\point}}\]}{t},
\end{align}
which is referred to as the $\quotl$flow of the smooth vector-field $\avecf{}$ on the manifold $\Man{}$$\quotr$.
\endef
\theorem
$\avecf{}$ is taken as an element of $\vecf{\Man{}}{\infty}$ (a smooth vector-field on the manifold $\Man{}$),
and $\point$ as a point of the manifold $\Man{}$.
\begin{align}
&\Foreach{t}{\func{\maxinterval{\Man{}}}{\binary{\avecf{}}{\point}}}\cr
&\left\{\begin{aligned}
&\func{\maxinterval{\Man{}}}{\binary{\avecf{}}{\func{\[\func{\maxintcurve{\Man{}}}{\binary{\avecf{}}{\point}}\]}{t}}}=
\defSet{s-t}{s\in\func{\maxinterval{\Man{}}}{\binary{\avecf{}}{\point}}},\cr
&\Foreach{s}{\func{\maxinterval{\Man{}}}{\binary{\avecf{}}{\func{\[\func{\maxintcurve{\Man{}}}{\binary{\avecf{}}{\point}}\]}{t}}}}
\func{\[\func{\maxintcurve{\Man{}}}{\binary{\avecf{}}{\func{\[\func{\maxintcurve{\Man{}}}{\binary{\avecf{}}{\point}}\]}{t}}}\]}{s}=
\func{\[\func{\maxintcurve{\Man{}}}{\binary{\avecf{}}{\point}}\]}{s+t}.
\end{aligned}\right.
\end{align}
Consequently,
\begin{align}
&\Foreach{t}{\func{\maxinterval{\Man{}}}{\binary{\avecf{}}{\point}}}\cr
&\Foreach{s}{\defSet{s-t}{s\in\func{\maxinterval{\Man{}}}{\binary{\avecf{}}{\point}}}}
\func{\vfflow{\Man{}}{\avecf{}}}{\binary{s}{\func{\vfflow{\Man{}}{\avecf{}}}{\binary{t}{\point}}}}=
\func{\vfflow{\Man{}}{\avecf{}}}{\binary{s+t}{\point}}.
\end{align}
\endthm
\definition\label{defcompletevectorfieldFlow}
$\avecf{}$ is taken as a complete vector-field on $\Man{}$. The mapping $\vfFlow{\Man{}}{\avecf{}}$ is defined as,
\begin{align}
&\vfFlow{\Man{}}{\avecf{}}\indef{\R}{\Func{\M{}}{\M{}}},\cr
&\Foreach{t}{\R}
\Foreach{\point}{\M{}}
\func{\[\func{\vfFlow{\Man{}}{\avecf{}}}{t}\]}{\point}\eqdef
\func{\vfflow{\Man{}}{\avecf{}}}{\binary{t}{\point}}.
\end{align}
\endef
\theorem\label{thmcompletevectorfieldFlowaction}
$\avecf{}$ is taken as a complete vector-field on $\Man{}$. For every $t$ in $\R$, $\func{\vfFlow{\Man{}}{\avecf{}}}{t}$
is an $\infty$-automorphism of the manifold $\Man{}$. Furthermore, $\vfFlow{\Man{}}{\avecf{}}$ is a group-homomorphism
from the additive group of $\R$ to the group of $\infty$-automorphisms of $\Man{}$. That is,
\begin{equation}
\vfFlow{\Man{}}{\avecf{}}\in\GHom{\opair{\R}{+}}{\GDiff{\infty}{\Man{}}}.
\end{equation}
So $\vfFlow{\Man{}}{\avecf{}}$ is in particular an action of the group $\opair{\R}{+}$ on the set of points $\M{}$ of the manifold $\Man{}$.
It is straightforward to check that the set of orbits of this action coincides with the set of images of
all integral-curves of $\avecf{}$ on $\Man{}$. Therefore, considering that the set of orbits of an action of a group on a set is a
partition of that set, it becomes evident that the set of images of all integral-curves of $\avecf{}$ on $\Man{}$ is a partition of
$\M{}$. As an immediate consequence, two distinct integral-curves of $\avecf{}$ do not intersect each other (which is also true
when $\avecf{}$ is a non-complete smooth vector-field).
\endthm
\theorem\label{thmintegralcurvesofvectorfieldswithvariablescale}
$\avecf{}$ is taken as a smooth vector-field on $\Man{}$, and $s$ as as a non-zero real number.
Changing the scale of $\avecf{}$ by the factor $s$, directly leads to a re-scaling of the maximal integral-curves of $\avecf{}$
by the factor $s$, that is re-sizing the domain of them by the factor $s^{-1}$ and altering their speed by the factor $s$.
Precisely, given a point $\point$ of $\Man{}$,
\begin{align}
&\func{\maxinterval{\Man{}}}{\binary{s\avecf{}}{\point}}=\defSet{\frac{t}{s}}
{t\in\func{\maxinterval{\Man{}}}{\binary{\avecf{}}{\point}}},\cr
&\Foreach{t}{\func{\maxinterval{\Man{}}}{\binary{s\avecf{}}{\point}}}
\func{\[\func{\maxintcurve{\Man{}}}{\binary{s\avecf{}}{\point}}\]}{t}=
\func{\[\func{\maxintcurve{\Man{}}}{\binary{\avecf{}}{\point}}\]}{st}.
\end{align}
\caution
Here, $s\avecf{}$ denotes the scalar-product of $\avecf{}$ with the scalar $s$ within the
canonical linear structure of the smooth vector-fields of $\Man{}$, $\Vecf{\Man{}}{\infty}$.
\endthm
\section{Foliations of a Manifold}
\textit{The theorems of this section are not provided with proofs,
and considered to be standard elementary facts of the theory of foliations.}
\fixed
$r$ is fixed as a positive integer less than or equal to the dimension of $\Man{}$ ($r\leq n$).
Also $q$ is defined to be $n-r$.
\endfixed
\definition\label{defrectangularcharts}
\begin{itemize}
\item[$\cdot$]
$\phi$ is taken as an element of $\maxatlas{}$ (a chart of the manifold $\Man{}$).
$\phi$ is referred to as an $\quotl$$\opair{r}{q}$-rectangular chart of the manifold $\Man{}$$\quotr$
iff there exists a rectangular open set $\U_1$ of $\R^r$ and a rectangular open set $\U_2$ of $\R^q$
such that $\funcimage{\phi}=\Cprod{\U_1}{\U_2}$ (by implicitly considering the natural isomorphism
between $\Cprod{\R^r}{\R^{q}}$ and $\R^n$).
\item[$\cdot$]
The set of all $\opair{r}{q}$-rectangular charts of $\Man{}$ is denoted by $\rectatlas{r}{q}{\Man{}}$, and
is called the $\quotl$$\opair{r}{q}$-rectangular atlas of $\Man{}$$\quotr$.
\end{itemize}
\endef
\definition\label{deffoliation}
$\foliation{}$ is taken as an equivalence relation on $\M{}$. $\foliation{}$ is referred to as an
$\quotl$$r$-dimensional $\infty$-foliation (or smooth foliation) of the manifold $\Man{}$$\quotr$ iff
each element of $\EqClass{\M{}}{\foliation{}}$ is a connected set of the topological-space $\mantops{\Man{}}$, and
for every point $\point$ of $\Man{}$,
there exists an $\opair{r}{q}$-rectangular chart $\phi$ of $\Man{}$ with $\funcimage{\phi}=\Cprod{\U_1}{\U_2}$,
$\U_1$ and $\U_2$ being rectanular open sets of $\R^r$ and $\R^q$ respectively,
such that for every $\leaf{}$ in $\EqClass{\M{}}{\foliation{}}$,
and for each connected component $\aconnectedset$ of $\domain{\phi}\cap\leaf{}$
(with respect to its inherited topology from that of $\mantops{\Man{}}$),
there exists a point $y$ of $\R^{q}$ in a way that
$\aconnectedset=\func{\pimage{\phi}}{\Cprod{\U_1}{\seta{y}}}$ (by implicitly considering the natural isomorphism
between $\Cprod{\R^r}{\R^{q}}$ and $\R^n$).
\endef
\definition
The set of all $r$-dimensional $\infty$-foliations of the manifold $\Man{}$ is denoted by $\Foliations{r}{\Man{}}$. That is,
\begin{align}
&\Foliations{r}{\Man{}}:=\cr
&
\defset{\foliation{}}{\EqR{\M{}}}
{\EqClass{\M{}}{\foliation{}}\subseteq\connecteds{\mantops{\Man{}}},~\(\begin{aligned}
&\Foreach{\point}{\M{}}
\Exists{\phi}{\rectatlas{r}{q}{\Man{}}}
\Foreach{\leaf{}}{\EqClass{\M{}}{\foliation{}}}\cr
&\Foreach{\aconnectedset}{\maxcon{\opair{\domain{\phi}\cap\leaf{}}{\stopology{\mantop{\Man{}}}{\domain{\phi}\cap\leaf{}}}}}\cr
&\Exists{y}{\R^q}
\aconnectedset=\func{\pimage{\phi}}{\Cprod{\R^r}{\seta{y}}}
\end{aligned}\)}.\cr
&{}
\end{align}
For each $\foliation{}$ in $\Foliations{r}{\Man{}}$,
the pair $\opair{\Man{}}{\foliation{}}$ is referred to as a $\quotl$smooth $\opair{n}{r}$-foliated space$\quotr$,
or simply as a $\quotl$foliated space$\quotr$.
\endef
\definition\label{deffoliatedatlas}
$\foliation{}$ is taken as an element of $\Foliations{r}{\Man{}}$ (an $r$-dimensional smooth foliation of $\Man{}$).
The set $\folatlas{\Man{}}{\foliation{}}{r}$ is defined as,
\begin{align}
&\folatlas{\Man{}}{\foliation{}}{r}:=\cr
&\defset{\phi}{\rectatlas{r}{q}{\Man{}}}
{\(\begin{aligned}
&\Foreach{\leaf{}}{\EqClass{\M{}}{\foliation{}}}\cr
&\Foreach{\aconnectedset}{\maxcon{\opair{\domain{\phi}\cap\leaf{}}{\stopology{\mantop{\Man{}}}{\domain{\phi}\cap\leaf{}}}}}\cr
&\Exists{y}{\R^q}
\aconnectedset=\func{\pimage{\phi}}{\Cprod{\R^r}{\seta{y}}}
\end{aligned}\)}.
\end{align}
$\folatlas{\Man{}}{\foliation{}}{r}$ is referred to as the $\quotl$foliated atlas of the $r$-dimensional smooth foliation
$\foliation{}$ of $\Man{}$$\quotr$.
\endef
\definition\label{defleafatlas}
$\foliation{}$ is taken as an element of $\Foliations{r}{\Man{}}$. For every element $\leaf{}$ of $\EqClass{\M{}}{\foliation{}}$,
the foliated atlas of $\foliation{}$ induces a set of charts on $\leaf{}$ which is denoted by
$\func{\leafatlas{\Man{}}{\foliation{}}{r}}{\leaf{}}$ and defined as,
\begin{equation}
\func{\leafatlas{\Man{}}{\foliation{}}{r}}{\leaf{}}:=
\defSet{\cmp{\Projection{n}{1}{r}}{\func{\res{\phi}}{\domain{\phi}\cap\leaf{}}}}{\phi\in\folatlas{\Man{}}{\foliation{}}{r},~
\domain{\phi}\cap\leaf{}\neq\empty},
\end{equation}
where $\Projection{n}{1}{r}$ denotes the projection of $\R^n$ onto the subspace spanned by its first $r$ standard
basis vectors.
\endef
\theorem
$\foliation{}$ is taken as an element of $\Foliations{r}{\Man{}}$. For every element $\leaf{}$ of $\EqClass{\M{}}{\foliation{}}$,
$\func{\leafatlas{\Man{}}{\foliation{}}{r}}{\leaf{}}$ is a $\difclass{\infty}$ atlas on $\leaf{}$ modeled on the Banach space $\R^r$.
That is,
\begin{equation}
\Foreach{\leaf{}}{\EqClass{\M{}}{\foliation{}}}
\func{\leafatlas{\Man{}}{\foliation{}}{r}}{\leaf{}}\in
\atlases{\infty}{\leaf{}}{\R^r}.
\end{equation}
\endthm
\theorem
$\foliation{}$ is taken as an element of $\Foliations{r}{\Man{}}$. Every element $\leaf{}$ of $\EqClass{\M{}}{\foliation{}}$ endowed with
the $\difclass{\infty}$ maximal-atlas on it modeled on $\R^r$ and generated by $\func{\leafatlas{\Man{}}{\foliation{}}{r}}{\leaf{}}$
is an $r$-dimensional and $\difclass{\infty}$ manifold modeled on $\R^r$.
That is, for every element $\leaf{}$ of $\EqClass{\M{}}{\foliation{}}$,
$\opair{\leaf{}}{\func{\maxatlasgen{\infty}{\leaf{}}{\R^r}}{\func{\leafatlas{\Man{}}{\foliation{}}{r}}{\leaf{}}}}$
is an $r$-dimensional and $\difclass{\infty}$ differentiable structure with the underlying topological-space being Hausdorff
and second-countable, and hence an $r$-dimensional and $\difclass{\infty}$ manifold. In addition, for every $\leaf{}$
in $\EqClass{\M{}}{\foliation{}}$, the underlying topology of
$\opair{\leaf{}}{\func{\maxatlasgen{\infty}{\leaf{}}{\R^r}}{\func{\leafatlas{\Man{}}{\foliation{}}{r}}{\leaf{}}}}$
is finer than that of the manifold $\Man{}$.
\endthm
\definition\label{defleaves}
$\foliation{}$ is taken as an element of $\Foliations{r}{\Man{}}$. The set of all $r$-dimensional $\difclass{\infty}$
manifolds constructed by endowing an element of $\EqClass{\M{}}{\foliation{}}$ by its canonically induced maximal atlas
by the smooth foliation $\foliation{}$ is denoted by $\func{\Leaves{\Man{}}{r}}{\foliation{}}$. That is,
\begin{align}
\func{\Leaves{\Man{}}{r}}{\foliation{}}:=
\defSet{\opair{\leaf{}}{\func{\maxatlasgen{\infty}{\leaf{}}{\R^r}}{\func{\leafatlas{\Man{}}{\foliation{}}{r}}{\leaf{}}}}}
{\leaf{}\in\EqClass{\M{}}{\foliation{}}}.
\end{align}
Each element of $\func{\Leaves{\Man{}}{r}}{\foliation{}}$ is called a
$\quotl$leaf of the $\infty$-foliation $\foliation{}$ of $\Man{}$$\quotr$.\\
For every $\leaf{}$ in $\EqClass{\M{}}{\foliation{}}$, the corresponding canonical differentiable structure induced by $\foliation{}$
is denoted by $\func{\leafman{\Man{}}{\foliation{}}{r}}{\leaf{}}$. That is,
\begin{equation}
\Foreach{\leaf{}}{\EqClass{\M{}}{\foliation{}}}
\func{\leafman{\Man{}}{\foliation{}}{r}}{\leaf{}}:=
\opair{\leaf{}}{\func{\maxatlasgen{\infty}{\leaf{}}{\R^r}}{\func{\leafatlas{\Man{}}{\foliation{}}{r}}{\leaf{}}}}.
\end{equation}
Also, for every $\leaf{}$ in $\EqClass{\M{}}{\foliation{}}$, this corresponded manifold can simply be denoted by
$\Leafman{\leaf{}}$ when everything is clear about the basic manifold and the considered foliation.
\endef
\theorem\label{thmleavesareimmersedsubmanifolds}
$\foliation{}$ is taken as an element of $\Foliations{r}{\Man{}}$.
Every leaf of the foliation $\foliation{}$ is an $r$-dimensional immersed submanifold of $\Man{}$.
That is, for every $\leaf{}$ in $\EqClass{\M{}}{\foliation{}}$, the injection of $\leaf{}$
into $\M{}$ is an $\infty$-immersion from $\Leafman{\leaf{}}$ to $\Man{}$.
\begin{align}
\Foreach{\leaf{}}{\EqClass{\M{}}{\foliation{}}}
\left\{\begin{aligned}
&\Injection{\leaf{}}{\M{}}\in\mapdifclass{\infty}{\Leafman{\leaf{}}}{\Man{}},\\
&\Foreach{\point}{\leaf{}}\func{\(\der{\Injection{\leaf{}}{\M{}}}{\Leafman{\leaf{}}}{\Man{}}\)}{\point}\in
\Linisom{\tanspace{\point}{\Leafman{\leaf{}}}}
{\tanspace{\point}{\Man{}}}.
\end{aligned}\right.
\end{align}
\endthm
\corollary
$\foliation{}$ is taken as an element of $\Foliations{r}{\Man{}}$.
The union of the underlying topologies of all leaves of the $\infty$-foliation $\foliation{}$ of $\Man{}$
is a topology on $\M{}$ which is finer than the underlying topology of $\Man{}$. That is,
\begin{equation}
\mantop{\Man{}}\subseteq
\Union{\leaf{}}{\EqClass{\M{}}{\foliation{}}}
{\mantop{\func{\leafman{\Man{}}{\foliation{}}{r}}{\leaf{}}}}
\in\alltopologies{\M{}}.
\end{equation}
\endcor
\definition
The union of all collections of $r$-dimensional vector-subspaces of $\Tanspace{\point}{\Man{}}$,
$\point$ ranging over all points of the manifold $\Man{}$, is denoted by $\subtanbun{\Man{}}{r}$.
\begin{equation}
\subtanbun{\Man{}}{r}:=\Union{\point}{\M{}}{\subvec{\Tanspace{\point}{\Man{}}}{r}}.
\end{equation}
\endef
\definition\label{deftangentbundleoffoliation}
$\foliation{}$ is taken as an element of $\Foliations{r}{\Man{}}$. $\foltanbun{\Man{}}{\foliation{}}{r}$
is defined to be the mapping from $\M{}$ to $\subtanbun{\Man{}}{r}$ that assigns to every
point $\point$ of $\Man{}$ the $r$-dimensionl vector-subspace of $\Tanspace{\point}{\Man{}}$ that is the image of
$\tanspace{\point}{\Leafman{\pEqclass{\point}{\foliation{}}}}$ under
the differential of the injection of the leaf $\Leafman{\pEqclass{\point}{\foliation{}}}$ (equivalence class of $\point$ by $\foliation{}$
considered with its canonical differentiable structure) into $\Man{}$. Precisely,
\begin{align}
&\foltanbun{\Man{}}{\foliation{}}{r}\indef\Func{\M{}}{\subtanbun{\Man{}}{r}},\cr
&\Foreach{\point}{\M{}}
\func{\foltanbun{\Man{}}{\foliation{}}{r}}{\point}\eqdef
\func{\image{\[\func{\(\der{\Injection{\pEqclass{\point}{\foliation{}}}{\M{}}}
{\Leafman{\pEqclass{\point}{\foliation{}}}}{\Man{}}\)}{\point}\]}}
{\tanspace{\point}{\Leafman{\pEqclass{\point}{\foliation{}}}}}.
\end{align}
$\foltanbun{\Man{}}{\foliation{}}{r}$ is called the $\quotl$tangent-bundle of the foliation $\foliation{}$ of the manifold $\Man{}$$\quotr$.
\endef
\section{Distributions of a Manifold}
\definition\label{defdistribution}
$\distribution{}$ is taken as an element of $\Func{\M{}}{\subtanbun{\Man{}}{r}}$, that is a mapping from the set of all points of $\Man{}$
to $\subtanbun{\Man{}}{r}$. $\distribution{}$ is referred to as an
$\quotl$$r$-dimensional $\infty$-distribution (or smooth distribution) of the manifold $\Man{}$$\quotr$ iff
for each point $\point$ of $\Man{}$, $\func{\distribution{}}{\point}$ lies in $\subvec{\Tanspace{\point}{\Man{}}}{r}$ and
there exists an open neighborhood $\U$ of $\seta{\point}$ in the
topological-space $\mantops{\Man{}}$ and an $r$-tuple $\mtuple{\avecf{1}}{\avecf{r}}$ of smooth vector-fields on $\Man{}$
such that for every point $\p{\point}$ in $\U$, the set $\seta{\suc{\func{\avecf{1}}{\p{\point}}}{\func{\avecf{r}}{\p{\point}}}}$
spans the vector-space $\func{\distribution{}}{\p{\point}}$, that is,
\begin{align}
\Foreach{\point}{\M{}}
\func{\distribution{}}{\point}\in\subvec{\Tanspace{\point}{\Man{}}}{r},~
\(\begin{aligned}
&\Exists{\mtuple{\avecf{1}}{\avecf{r}}}{\[\vecf{\Man{}}{\infty}\]^{\times n}}
\Exists{\U}{\func{\nei{\mantops{\Man{}}}}{\seta{\point}}}\cr
&\Foreach{\p{\point}}{\U}\func{\Vspan{\Tanspace{\p{\point}}{\Man{}}}}{\seta{\suc{\func{\avecf{1}}{\p{\point}}}{\func{\avecf{r}}{\p{\point}}}}}=
\func{\distribution{}}{\p{\point}}.
\end{aligned}\)
\end{align}
\endef
\definition
The set of all $r$-dimensional $\infty$-distributions of the manifold $\Man{}$ is denoted by $\Distributions{r}{\Man{}}$. That is,
\begin{align}
&\Distributions{r}{\Man{}}:=\cr
&\defset{\distribution{}}{\Func{\M{}}{\subtanbun{\Man{}}{r}}}
{\(\begin{aligned}
&\Foreach{\point}{\M{}}
\func{\distribution{}}{\point}\in\subvec{\Tanspace{\point}{\Man{}}}{r},\cr
&\Foreach{\point}{\M{}}
\Exists{\mtuple{\avecf{1}}{\avecf{r}}}{\[\vecf{\Man{}}{\infty}\]^{\times n}}\cr
&\Exists{\U}{\func{\nei{\mantops{\Man{}}}}{\seta{\point}}}\cr
&\Foreach{\p{\point}}{\U}\func{\Vspan{\Tanspace{\p{\point}}{\Man{}}}}{\seta{\suc{\func{\avecf{1}}{\p{\point}}}{\func{\avecf{r}}{\p{\point}}}}}=
\func{\distribution{}}{\p{\point}}.
\end{aligned}\)}.\cr
&{}
\end{align}
Each element of $\Distributions{1}{\Man{}}$ is alternatively referred to as a $\quotl$smooth line-field of $\Man{}$$\quotr$. Also,
each element of $\Distributions{2}{\Man{}}$ is alternatively referred to as a $\quotl$smooth plane-field of $\Man{}$$\quotr$.
\endef
\definition
$\[Integrability~of~a~Distribution\]$
\begin{itemize}
\item[$\centerdot$]
$\distribution{}$ is taken as an element of $\Distributions{r}{\Man{}}$.
$\distribution{}$ is referred to as an $\quotl$$r$-dimensional integrable smooth distribution of $\Man{}$$\quotr$ iff
there exists an $r$-dimensional smooth foliation $\foliation{}$ of $\Man{}$ such that $\foltanbun{\Man{}}{\foliation{}}{r}$
coincides with $\distribution{}$.
\item[$\centerdot$]
The set of all $r$-dimensional integrable smooth distributions of the manifold $\Man{}$ is denoted by $\integrabledist{r}{\Man{}}$.
That is,
\begin{align}
&~\integrabledist{r}{\Man{}}\cr
:=&\defset{\distribution{}}{\Distributions{r}{\Man{}}}
{\(\Exists{\foliation{}}{\Foliations{r}{\Man{}}}\foltanbun{\Man{}}{\foliation{}}{r}=\distribution{}\)}.\cr
&{}
\end{align}
\end{itemize}
\endef
\theorem\label{thmintegrabledistributionhasuniquesolution}
For every $r$-dimensional integrable smooth distribution $\distribution{}$ of $\Man{}$,
there exists a unique $r$-dimensional smooth foliation $\foliation{}$ of $\Man{}$ such that $\foltanbun{\Man{}}{\foliation{}}{r}$
coincides with $\distribution{}$.
\begin{equation}
\Foreach{\distribution{}}{\integrabledist{r}{\Man{}}}
\Existsu{\foliation{}}{\Foliations{r}{\Man{}}}
\foltanbun{\Man{}}{\foliation{}}{r}=\distribution{}.
\end{equation}
\endthm
\definition\label{defsolutionofanintegrabledistributionasafoliation}
For every $r$-dimensional integrable smooth distribution $\distribution{}$ of $\Man{}$,
the unique $r$-dimensional smooth foliation $\foliation{}$ of $\Man{}$ for which $\foltanbun{\Man{}}{\foliation{}}{r}$
coincides with $\distribution{}$ is referred to as the $\quotl$unique solution (or integral)
of the integrable smooth distribution $\distribution{}$
of $\Man{}$$\quotr$ and is denoted by $\intdistsolution{\Man{}}{r}\distribution{}$.
\endef
\definition\label{defdistributionpreservingmap}
\begin{itemize}
\item[$\centerdot$]
For every $\distribution{}$ in $\Distributions{r}{\Man{}}$, the ordered pair $\opair{\Man{}}{\distribution{}}$ is called a
$\quotl$$\difclass{\infty}$ distributed-space$\quotr$.
\item[$\centerdot$]
$\distribution{}$ and $\distribution{1}$ are taken as elements of $\Distributions{r}{\Man{}}$ and $\Distributions{r}{\Man{1}}$,
respectively. $\cf$ is taken as an element of $\mapdifclass{\infty}{\Man{}}{\Man{1}}$. $\cf$ is called a $\quotl$distribution-preserving
smooth map from the $\difclass{\infty}$ distributed-space $\opair{\Man{}}{\distribution{}}$ to the $\difclass{\infty}$ distributed-space
iff
\begin{equation}
\Foreach{\point}{\M{}}
\func{\image{\[\der{\cf}{\Man{}}{\Man{1}}\]}}{\func{\distribution{}}{\point}}=
\func{\distribution{1}}{\func{\cf}{\point}}.
\end{equation}
More specifically, if $\cf$ is an element of $\Diffeo{\infty}{\Man{}}{\Man{1}}$ with this property, it is
also called a $\quotl$distribution-preserving $\infty$-diffeomorphism from $\opair{\Man{}}{\distribution{}}$ to
$\opair{\Man{1}}{\distribution{}}$$\quotr$.
\end{itemize}
\endef
\theorem\label{thmdistpreservingdiffinterchangestheleavesofsolutions}
$\distribution{}$ and $\distribution{1}$ are taken as elements of $\integrabledist{r}{\Man{}}$ and $\integrabledist{r}{\Man{1}}$
($r$-dimensional integrable smooth distributions of $\Man{}$ and $\Man{1}$), respectively. $\foliation{}$ and $\foliation{1}$
are defined to be the $r$-dimensional smooth foliations of $\Man{}$ that are unique solutions of $\distribution{}$ and
$\distribution{1}$, respectively. For every distribution-preserving $\infty$-diffeomorphism $\cf$ from
the distributed space $\opair{\Man{}}{\distribution{}}$ to the distributed-space $\opair{\Man{1}}{\distribution{1}}$,
\begin{equation}
\EqClass{\M{1}}{\foliation{1}}=\defSet{\func{\image{\cf}}{\leaf{}}}{\leaf{}\in\EqClass{\M{}}{\foliation{}}},
\end{equation}
which equivalently means that the image-map $\image{\cf}$ forms a one-to-one correspondence between the
leaves of $\foliation{}$ and those of $\foliation{1}$.
\endthm
\theorem\label{thmintegralmanifoldofadistributionisanopensetofaleafogthesolutionfoliation}
$\distribution{}$ is taken as an element of $\integrabledist{r}{\Man{}}$ (an $r$-dimensional integrable smooth distribution of $\Man{}$),
and $\point_0$ a point of $\Man{}$. If $\Man{1}$ is an $\infty$-immersed submanifold of $\Man{}$ such that $\point_0\in\M{1}$,
the underlying topological-space $\mantops{\Man{1}}$ of $\Man{1}$ is connected, and,
\begin{equation}
\Foreach{\point}{\M{1}}
\func{\image{\[\func{\(\der{\Injection{\M{1}}{\M{}}}
{\Man{1}}{\Man{}}\)}{\point}\]}}{\tanspace{\point}{\Man{1}}}=
\func{\distribution{}}{\point},
\end{equation}
then the set of points $\M{1}$ of the manifold $\Man{1}$ is an open set of the unique leaf of the foliation
$\intdistsolution{\Man{}}{r}\distribution{}$ (the unique solution of the integrable smooth distribution $\distribution{}$ of $\Man{}$)
that contains the point $\point_0$ (with respect to the leaf's own intrinsic topology). That is, by defining
$\foliation{}:=\intdistsolution{\Man{}}{r}\distribution{}$,
\begin{equation}
\M{1}\in
\mantop{\Leafman{\pEqclass{\point}{\foliation{}}}}.
\end{equation}
Furthermore, $\Man{1}$ is an embedded submanifold of $\Leafman{\pEqclass{\point}{\foliation{}}}$.\\
\caution
Any manifold such $\Man{1}$ is called a $\quotl$connected integral manifold of the integrable smooth distribution $\distribution{}$ of $\Man{}$
passing through the point $\point_0$$\quotr$.
\endthm
\textit{A detailed proof of the assertion above is addressed in \cite[page~94,~Theorem~2]{Chevalley},
not exactly in the language of foliations, but an equivalent one.
It is important to be aware that \cite{Chevalley} takes
the underlying topological-spaces of manifolds to be connected, by definition.}
\chapteR{Abstract Topological Group}
\thispagestyle{fancy}
\section{Basic Structure of a Topological Group}
\definition\label{defgrouptranslations}
$\Group{}=\opair{\G{}}{\gop{}}$ is taken to be a group, where $\gop{}$ is called the $\quotl$group-operation of the
group $\Group{}$$\quotr$.
\begin{itemize}
\item[$\centerdot$]
The mappings $\function{\ginv{\Group{}}}{\G{}}{\G{}}$
is defined as,
\begin{align}
\Foreach{\g{}}{\G{}}\func{\ginv{\Group{}}}{\g{}}\eqdef\invg{\g{}}{\Group{}},
\end{align}
which is called the $\quotl$inverse-mapping of the group $\Group{}$$\quotr$.
\item[$\centerdot$]
The mapping $\function{\gopr{\gop{}}{\Group{}}}{\Cprod{\G{}}{\G{}}}{\G{}}$ is defined as,
\begin{equation}
\Foreach{\opair{\g{1}}{\g{2}}}{\Cprod{\G{}}{\G{}}}
\func{\gopr{\gop{}}{\Group{}}}{\binary{\g{1}}{\g{2}}}\eqdef\g{1}\gop{}\invg{\(\g{2}\)}{\Group{}},
\end{equation}
Which is called the $\quotl$secondary group-operation of the group $\Group{}$$\quotr$.
When there is no ambiguity about the underlying group $\Group{}$, $\gopr{\gop{}}{}$ can replace
$\gopr{\gop{}}{\Group{}}$.
\item[$\centerdot$]
For every $\g{0}$ in $\G{}$, the mapping $\function{\gltrans{\Group{}}{\g{0}}}{\G{}}{\G{}}$ is defined as,
\begin{equation}
\Foreach{\g{}}{\G{}}
\func{\gltrans{\Group{}}{\g{0}}}{\g{}}\eqdef
\g{0}\gop{}\g{},
\end{equation}
or equivalently,
\begin{align}
\gltrans{\Group{}}{\g{0}}=\cmp{\cmp{\gop{}}{\(\constmap{\G{}}{\g{0}}\times\identity{\G{}}\)}}{\diagmap{\G{}}},
\end{align}
which is called the $\quotl$left-translation of $\Group{}$ by $\g{0}$$\quotr$.
\item[$\centerdot$]
For every $\g{0}$ in $\G{}$, the mapping $\function{\grtrans{\Group{}}{\g{0}}}{\G{}}{\G{}}$ is defined as,
\begin{equation}
\Foreach{\g{}}{\G{}}
\func{\grtrans{\Group{}}{\g{0}}}{\g{}}\eqdef
\g{}\gop{}\g{0},
\end{equation}
or equivalently,
\begin{align}
\grtrans{\Group{}}{\g{0}}=\cmp{\cmp{\gop{}}{\(\identity{\G{}}\times\constmap{\G{}}{\g{0}}\)}}{\diagmap{\G{}}},
\end{align}
which is called the $\quotl$right-translation of $\Group{}$ by $\g{0}$$\quotr$.
\item[$\centerdot$]
For every $\g{0}$ in $\G{}$, the mapping $\function{\gconj{\Group{}}{\g{0}}}{\G{}}{\G{}}$ is defined as,
\begin{equation}
\Foreach{\g{}}{\G{}}
\func{\gconj{\Group{}}{\g{0}}}{\g{}}\eqdef
\g{0}\gop{}\g{}\gop{}\invg{\(\g{0}\)}{\Group{}},
\end{equation}
or equivalently,
\begin{equation}
\gconj{\Group{}}{\g{0}}:=\cmp{\gltrans{\Group{}}{\g{0}}}{\grtrans{\Group{}}{\invG{\g{0}}{}}},
\end{equation}
which is called the $\quotl$$\g{0}$-conjugation of the group $\Group{}$$\quotr$.
\item[$\centerdot$]
For every $\g{}$ in $\G{}$,
\begin{align}
&\gpower{\g{}}{\gop{}}{0}\eqdef\IG{},\cr
&\Foreach{n}{\Zp}\gpower{\g{}}{\gop{}}{n}\eqdef\g{}\gop{}\gpower{\g{}}{\gop{}}{n-1}
=\Succ{\g{}}{\g{}}{\gop{}}{n},\cr
&\Foreach{n}{\Zn}\gpower{\g{}}{\gop{}}{n}\eqdef\gpower{\(\invG{\g{}}{}\)}{\gop{}}{\(-n\)}.
\end{align}
For every $\g{}$ in $\G{}$ and every $n$ in $\Z$, $\gpower{\g{}}{\gop{}}{n}$
is called the $\quotl$$n$-th power of $\g{}$ in group $\Group{}$$\quotr$.
\item[$\centerdot$]
For every pair of subsets $\V_1$ and $\V_2$ of $\G{}$,
\begin{equation}
\gsetprod{\V_1}{\V_2}{\Group{}}:=\defSet{\g{1}\gop{}\g{2}}{\g{1}\in\V_1,~\g{2}\in\V_2},
\end{equation}
which is referred to as the $\quotl$product of $\V_1$ and $\V_2$ in the group $\Group{}$$\quotr$.
\item[$\centerdot$]
For every subset $\V$ of $\G{}$,
\begin{align}
&\gpower{\V}{\gop{}}{0}\eqdef\seta{\IG{}},\cr
&\Foreach{n}{\Zp}
\gpower{\V}{\gop{}}{n}\eqdef
\gsetprod{\V}{\gpower{\V}{\gop{}}{n-1}}{\Group{}}
=\defSet{\Succc{\g{1}}{\g{n}}{\gop{}}}{\mtuple{\g{1}}{\g{n}}\in\G{}^{\times n}},\cr
&\Foreach{n}{\Zn}\gpower{\V}{\gop{}}{n}\eqdef\gpower{\(\func{\image{\ginv{\TopgG{\Topgroup{}}}}}{\V}\)}{\gop{}}{\(-n\)}.
\end{align}
For every subset $\V$ of $\G{}$ and every $n$ in $\Z$, $\gpower{\V}{\gop{}}{n}$
is called the $\quotl$$n$-th power of the set $\V$ in group $\Group{}$$\quotr$.
\item[$\centerdot$]
Any subset $\V$ of $\G{}$ is called a $\quotr$symmetric set of the group $\Group{}$ iff,
\begin{equation}
\func{\image{\ginv{\Group{}}}}{\V}=\V.
\end{equation}
\item[$\centerdot$]
For every non-empty subset $\V$ of $\G{}$, $\Gen{\V}{\Group{}}$ is defined to be the smallest subgroup
of $\Group{}$ including $\V$, that is,
\begin{equation}
\Gen{\V}{\Group{}}:=\intersection{\defset{\asubgroup{}}{\Subgroups{\Topgroup{}}}{\asubgroup{}\supseteq\V}}=
\Union{n}{\Z}{\gpower{\V}{\gop{}}{n}}.
\end{equation}
$\Gen{\V}{\Group{}}$ is referred to as the $\quotl$subgroup of $\Group{}$ generated by $\V$$\quotr$.
Also, $\V$ is called a $\quotl$generator of $\Gen{\V}{\Group{}}$$\quotr$.
\item[$\centerdot$]
Given a subgroup $\asubgroup{}$ of $\Group{}$, $\V$ is referred to as a $\quotl$normal subgroup of $\Group{}$$\quotr$
iff,
\begin{equation}
\Foreach{\g{}}{\G{}}
\Foreach{\hh}{\asubgroup{}}
\func{\gconj{\Group{}}{\g{}}}{\hh}=
\g{}\gop{}\hh\invG{\g{}}{}\in\asubgroup{}.
\end{equation}
It is straightforward to check that for a normal subgroup $\asubgroup{}$ of $\Group{}$,
\begin{equation}
\func{\image{\[\gconj{\Group{}}{\g{0}}\]}}{\asubgroup{}}=\asubgroup{}.
\end{equation}
The set of all normal subgroups of $\Group{}$ is denoted by $\NSubgroups{\Group{}}$.
\end{itemize}
\endef
\definition\label{deftopologicalgroup}
$\G{}$ is taken as a set, $\gop{}$ as an element of $\Func{\Cprod{\G{}}{\G{}}}{\G{}}$ (a binary operation on $\G{}$),
and $\atopology{}$ as an element of $\CSs{\CSs{\G{}}}$ (a collection of subsets of $\G{}$).
The triple $\triple{\G{}}{\gop{}}{\atopology{}}$ is referred to as a $\quotl$topological group$\quotr$ iff these properties are satisfied.
\begin{itemize}
\item[\myitem{TG~1.}]
$\Group{}:=\opair{\G{}}{\gop{}}$ is a group.
\item[\myitem{TG~2.}]
$\opair{\G{}}{\atopology{}}$ is a topological-space. Equivalently, $\atopology{}$ is a topology on the set $\G{}$. That is,
\begin{equation}
\atopology{}\in\alltopologies{\G{}}.
\end{equation}
\item[\myitem{TG~3.}]
$\gop{}$ is a continuous map from the topological-product
$\topprod{\opair{\G{}}{\atopology{}}}{\opair{\G{}}{\atopology{}}}$ to the topological-space $\opair{\G{}}{\atopology{}}$. That is,
\begin{equation}
\gop{}\in\CF{\topprod{\opair{\G{}}{\atopology{}}}{\opair{\G{}}{\atopology{}}}}{\opair{\G{}}{\atopology{}}}.
\end{equation}
\item[\myitem{TG~4.}]
$\ginv{\Group{}}$ is a continuous map from $\opair{\G{}}{\atopology{}}$ to $\opair{\G{}}{\atopology{}}$. That is,
\begin{equation}
\ginv{\Group{}}\in\CF{\opair{\G{}}{\atopology{}}}{\opair{\G{}}{\atopology{}}}.
\end{equation}
\end{itemize}
\endef
\theorem\label{thmtopologicalgroupequiv0}
$\G{}$ is taken as a set, $\gop{}$ as an element of $\Func{\Cprod{\G{}}{\G{}}}{\G{}}$,
and $\atopology{}$ as an element of $\alltopologies{\G{}}$.
$\triple{\G{}}{\gop{}}{\atopology{}}$ is a topological group if and only if $\gopr{\gop{}}{\Group{}}$
is a continuous map, that is,
\begin{equation}
\gopr{\gop{}}{\Group{}}\in\CF{\topprod{\opair{\G{}}{\atopology{}}}{\opair{\G{}}{\atopology{}}}}{\opair{\G{}}{\atopology{}}}.
\end{equation}
\proof
Let $\IG{}$ denote the identity element of the group $\opair{\G{}}{\gop{}}$.
It can be easily verified that,
\begin{align}
\gopr{\gop{}}{\Group{}}&=\cmp{\gop{}}{\(\identity{\G{}}\times\ginv{\Group{}}\)},
\label{thmtopologicalgroupequiv0peq1}\\
\ginv{\Group{}}&=\cmp{\cmp{\gopr{\gop{}}{\Group{}}}{\(\constmap{\G{}}{\IG{}}\times\identity{\G{}}\)}}{\diagmap{\G{}}},
\label{thmtopologicalgroupequiv0peq2}\\
\gop{}&=\cmp{\gopr{\gop{}}{\Group{}}}{\(\identity{\G{}}\times\ginv{\Group{}}\)},
\label{thmtopologicalgroupequiv0peq3}
\end{align}
where $\function{\constmap{\G{}}{\IG{}}}{\G{}}{\G{}}$ denotes the constnant map on $\G{}$ with single value $\IG{}$,
$\function{\diagmap{\G{}}}{\G{}}{\Cprod{\G{}}{\G{}}}$ denotes the diagonal map on $\G{}$ sending
each $\g{}$ in $\G{}$ to $\opair{\g{}}{\g{}}$, and $\function{\identity{\G{}}}{\G{}}{\G{}}$ denotes the
identity map on $\G{}$ sending each $\g{}$ in $\G{}$ to $\g{}$. All mappings $\constmap{\G{}}{\IG{}}$,
$\identity{\G{}}$, and $\diagmap{\G{}}$ are continuous. That is,
\begin{align}
\constmap{\G{}}{\IG{}}\in\CF{\opair{\G{}}{\atopology{}}}{\opair{\G{}}{\atopology{}}},
\label{thmtopologicalgroupequiv0peq4}\\
\identity{\G{}}\in\CF{\opair{\G{}}{\atopology{}}}{\opair{\G{}}{\atopology{}}},
\label{thmtopologicalgroupequiv0peq5}\\
\diagmap{\G{}}\in\CF{\opair{\G{}}{\atopology{}}}{\manprod{\opair{\G{}}{\atopology{}}}{\opair{\G{}}{\atopology{}}}}.
\label{thmtopologicalgroupequiv0peq6}
\end{align}
\begin{itemize}
\item[$\pr{1}$]
It is assumed that $\triple{\G{}}{\gop{}}{\atopology{}}$ is a topological group. So according to \refdef{deftopologicalgroup},
\begin{align}
\gop{}&\in\CF{\manprod{\opair{\G{}}{\atopology{}}}{\opair{\G{}}{\atopology{}}}}{\opair{\G{}}{\atopology{}}},\\
\ginv{\Group{}}&\in\CF{\opair{\G{}}{\atopology{}}}{\opair{\G{}}{\atopology{}}}.
\end{align}
So since $\identity{\G{}}\in\CF{\opair{\G{}}{\atopology{}}}{\opair{\G{}}{\atopology{}}}$,
obviously,
\begin{equation}
\(\identity{\G{}}\times\ginv{\Group{}}\)\in
\CF{\topprod{\opair{\G{}}{\atopology{}}}{\opair{\G{}}{\atopology{}}}}
{\topprod{\opair{\G{}}{\atopology{}}}{\opair{\G{}}{\atopology{}}}},
\end{equation}
and hence,
\begin{equation}
\cmp{\gop{}}{\(\identity{\G{}}\times\ginv{\Group{}}\)}
\in\CF{\topprod{\opair{\G{}}{\atopology{}}}{\opair{\G{}}{\atopology{}}}}{\opair{\G{}}{\atopology{}}},
\end{equation}
which according to \Ref{thmtopologicalgroupequiv0peq1} means,
\begin{equation}
\gopr{\gop{}}{\Group{}}
\in\CF{\manprod{\opair{\G{}}{\atopology{}}}{\opair{\G{}}{\atopology{}}}}{\opair{\G{}}{\atopology{}}}.
\end{equation}
\endp
\end{itemize}
\begin{itemize}
\item[$\pr{2}$]
It is now assumed that $\gopr{\gop{}}{\Group{}}
\in\CF{\topprod{\opair{\G{}}{\atopology{}}}{\opair{\G{}}{\atopology{}}}}{\opair{\G{}}{\atopology{}}}$.
Then according to \Ref{thmtopologicalgroupequiv0peq2}, \Ref{thmtopologicalgroupequiv0peq4}, \Ref{thmtopologicalgroupequiv0peq5}, and
\Ref{thmtopologicalgroupequiv0peq6},
\begin{equation}
\ginv{\Group{}}\in\CF{\opair{\G{}}{\atopology{}}}{\opair{\G{}}{\atopology{}}},
\end{equation}
and accordingly, \Ref{thmtopologicalgroupequiv0peq3} and \Ref{thmtopologicalgroupequiv0peq5} imply,
\begin{equation}
\gop{}\in\CF{\topprod{\opair{\G{}}{\atopology{}}}{\opair{\G{}}{\atopology{}}}}{\opair{\G{}}{\atopology{}}}.
\end{equation}
Therefore $\triple{\G{}}{\gop{}}{\atopology{}}$ is a topological group.
\endp
\end{itemize}
\endthm
\definition\label{deftopologyoftopologicalgroup}
$\Topgroup{}=\triple{\G{}}{\gop{}}{\atopology{}}$ is taken to be a topological group.
\begin{itemize}
\item[$\centerdot$]
The intrinsic (or underlying) group structure of $\Topgroup{}$ is denoted by $\TopgG{\Topgroup{}}$. That is,
$\TopgG{\Topgroup{}}:=\opair{\G{}}{\gop{}}$.
\item[$\centerdot$]
The underlying topological-space of $\Topgroup{}$ is denoted by $\topgtops{\Topgroup{}}$. That is,
$\topgtops{\Topgroup{}}:=\opair{\G{}}{\atopology{}}$. This is called the $\quotl$intrinsic topological-structure
of the topological group $\Topgroup{}$$\quotr$.
\item[$\centerdot$]
Each element of $\atopology{}$, that is each open set of the topological-space $\topgtops{\Topgroup{}}$,
is called an $\quotl$open set of the topological group $\Topgroup{}$$\quotr$.
\item[$\centerdot$]
Each closed set of the topological-space $\topgtops{\Topgroup{}}$,
is called a $\quotl$closed set of the topological group $\Topgroup{}$$\quotr$.
\item[$\centerdot$]
Each element of $\G{}$ is called a $\quotl$point of the topological group $\Topgroup{}$$\quotr$.
\item[$\centerdot$]
Any symmetric set of $\TopgG{\Topgroup{}}$ is called a $\quotl$symmetric set of the topological group $\Topgroup{}$$\quotr$.
\end{itemize}
\endef
\fixed
\begin{itemize}
\item[$\centerdot$]
$\Topgroup{}=\triple{\G{}}{\gop{}}{\atopology{}}$ is fixed as a topological group.
$\IG{}$ is defined to be the identity element of the group $\TopgG{\Topgroup{}}$.
\item[$\centerdot$]
$\Topgroup{1}=\triple{\G{1}}{\gop{1}}{\atopology{1}}$ is fixed as a topological group.
$\IG{1}$ is defined to be the identity element of the group $\TopgG{\Topgroup{1}}$.
\item[$\centerdot$]
$\Topgroup{2}=\triple{\G{2}}{\gop{2}}{\atopology{2}}$ is fixed as a topological group.
$\IG{2}$ is defined to be the identity element of the group $\TopgG{\Topgroup{2}}$.
\item[$\centerdot$]
$\Xt=\opair{\X}{\atopology{0}}$ is fixed as a topological-space.
\end{itemize}
\endfixed
\theorem\label{thminversemappingishomeomorphism}
The inverse-mapping of the group $\TopgG{\Topgroup{}}$ is a homeomorphism from
the underlying topological-space $\topgtops{\Topgroup{}}$ of $\Topgroup{}$ to itself. That is,
\begin{equation}
\ginv{\LieG{\Liegroup{}}}\in\HOF{\topgtops{\Topgroup{}}}{\topgtops{\Topgroup{}}}.
\end{equation}
\proof
According to \refdef{deftopologicalgroup}, $\ginv{\LieG{\Liegroup{}}}$ is a continuous map from $\topgtops{\Topgroup{}}$
to $\topgtops{\Topgroup{}}$.
\begin{equation}\label{thminversemappingishomeomorphismpeq1}
\ginv{\LieG{\Liegroup{}}}\in\CF{\topgtops{\Topgroup{}}}{\topgtops{\Topgroup{}}}.
\end{equation}
Moreover, it is well-known that $\ginv{\LieG{\Liegroup{}}}$ is a bijective mapping
from $\G{}$ to $\G{}$, and coincides with its own inverse, that is,
\begin{equation}
\finv{\ginv{\LieG{\Liegroup{}}}}=\ginv{\LieG{\Liegroup{}}}.
\end{equation}
Thus, in addition to \Ref{thminversemappingishomeomorphismpeq1},
\begin{equation}
\ginv{\LieG{\Liegroup{}}}\in\CF{\topgtops{\Topgroup{}}}{\topgtops{\Topgroup{}}},
\end{equation}
and consequently it becomes evident that $\ginv{\LieG{\Liegroup{}}}$ is a homeomorphism from the topological-space
$\topgtops{\Topgroup{}}$ to itself.
\endthm
\theorem\label{thmtranslationsarehomeomorphisms}
For every point $\g{}$ of $\Liegroup{}$, the left-translation and right-translation of $\TopgG{\Topgroup{}}$ by $\g{}$, and
the $\g{}$-conjugation of $\TopgG{\Topgroup{}}$
are homeomorphisms from the underlying manifold $\topgtops{\Topgroup{}}$ of the topological group $\Liegroup{}$ to itself. That is,
\begin{align}
\Foreach{\g{}}{\G{}}
\begin{cases}
\gltrans{\TopgG{\Topgroup{}}}{\g{}}\in\HOF{\topgtops{\Topgroup{}}}{\topgtops{\Topgroup{}}},\cr
\grtrans{\TopgG{\Topgroup{}}}{\g{}}\in\HOF{\topgtops{\Topgroup{}}}{\topgtops{\Topgroup{}}},\cr
\gconj{\TopgG{\Topgroup{}}}{\g{}}\in\HOF{\topgtops{\Topgroup{}}}{\topgtops{\Topgroup{}}}.
\end{cases}
\end{align}
\proof
According to \refdef{defgrouptranslations} and considering the continuity of the mappings $\constmap{\G{}}{\g{}}$ for all $\g{}$ in $\G{}$,
$\identity{\G{}}$, and $\diagmap{\G{}}$, it is clear that,
\begin{equation}\label{thmtranslationsarehomeomorphismspeq1}
\Foreach{\g{}}{\G{}}
\gltrans{\TopgG{\Topgroup{}}}{\g{}}\in\CF
{\topgtops{\Topgroup{}}}{\topgtops{\Topgroup{}}}.
\end{equation}
In addition, it is known that for every $\g{}$ in $\G{}$, $\gltrans{\TopgG{\Topgroup{}}}{\g{}}$ is a bijective mapping
and,
\begin{equation}\label{thmtranslationsarehomeomorphismspeq2}
\Foreach{\g{}}{\G{}}
\finv{\(\gltrans{\TopgG{\Topgroup{}}}{\g{}}\)}=\gltrans{\TopgG{\Topgroup{}}}{\invg{\g{}}{}}.
\end{equation}
\Ref{thmtranslationsarehomeomorphismspeq1} and \Ref{thmtranslationsarehomeomorphismspeq2} imply that,
\begin{align}
\Foreach{\g{}}{\G{}}
\gltrans{\TopgG{\Topgroup{}}}{\g{}}\in\CF
{\topgtops{\Topgroup{}}}{\topgtops{\Topgroup{}}},~
\finv{\(\gltrans{\TopgG{\Topgroup{}}}{\g{}}\)}\in\CF
{\topgtops{\Topgroup{}}}{\topgtops{\Topgroup{}}}.
\end{align}
This means $\gltrans{\TopgG{\Topgroup{}}}{\g{}}$ is ahomeomorphism from the topological-space $\topgtops{\Topgroup{}}$ to itself,
for all $\g{}$ in $\G{}$. That is,
\begin{equation}
\Foreach{\g{}}{\G{}}
\gltrans{\TopgG{\Topgroup{}}}{\g{}}\in\HOF{\topgtops{\Topgroup{}}}{\topgtops{\Topgroup{}}}.
\end{equation}
In an obviously similar way it can be seen that,
\begin{equation}
\Foreach{\g{}}{\G{}}
\grtrans{\TopgG{\Topgroup{}}}{\g{}}\in\HOF{\topgtops{\Topgroup{}}}{\topgtops{\Topgroup{}}}.
\end{equation}
Furthermore, since $\gltrans{\TopgG{\Topgroup{}}}{\g{}}$ and $\grtrans{\TopgG{\Topgroup{}}}{\g{}}$
are homeomorphisms from $\topgtops{\Topgroup{}}$ to itself and
$\gconj{\TopgG{\Topgroup{}}}{\g{}}=\cmp{\gltrans{\TopgG{\Topgroup{}}}{\g{}}}{\grtrans{\Group{}}{\g{}}}$, for all
$\g{}$ is in $\G{}$, and considering that the composition of any pair of homeomorphisms from $\topgtops{\Topgroup{}}$ to itself
is again a homeomorphism from $\topgtops{\Topgroup{}}$ to itself, it is evident that,
\begin{equation}
\Foreach{\g{}}{\G{}}
\gconj{\TopgG{\Topgroup{}}}{\g{}}\in\HOF{\topgtops{\Topgroup{}}}{\topgtops{\Topgroup{}}}.
\end{equation}
\endthm
\lemma\label{lemtranslationsofanopensetisopeninatopologicalgroup}
$\g{}$ is taken as an element of $\G{}$ (a point of the topological group $\Topgroup{}$).
\begin{itemize}
\item[$\centerdot$]
The image of every open (resp.~closed) set of the topological group $\Topgroup{}$
under the inverse-mapping of $\TopgG{\Topgroup{}}$ is also an open (resp. closed) set of $\Liegroup{}$. That is,
\begin{align}
&\Foreach{\U}{\atopology{}}
\func{\image{\[\ginv{\TopgG{\Topgroup{}}}\]}}{\U}\in\atopology{},\\
&\Foreach{\V}{\(\compl{\G{}}{\atopology{}}\)}
\func{\image{\[\ginv{\TopgG{\Topgroup{}}}\]}}{\V}\in\(\compl{\G{}}{\atopology{}}\).
\end{align}
\item[$\centerdot$]
The image of every open (resp.~closed) set of the smooth group $\Liegroup{}$
under the left-translation of $\TopgG{\Topgroup{}}$ by $\g{}$ is also an open (resp. closed) set of $\Liegroup{}$. That is,
\begin{align}
&\Foreach{\U}{\atopology{}}
\func{\image{\[\gltrans{\TopgG{\Topgroup{}}}{\g{}}\]}}{\U}\in\atopology{},\\
&\Foreach{\V}{\(\compl{\G{}}{\atopology{}}\)}
\func{\image{\[\gltrans{\TopgG{\Topgroup{}}}{\g{}}\]}}{\V}\in\(\compl{\G{}}{\atopology{}}\).
\end{align}
\item[$\centerdot$]
The image of every open (resp.~closed) set of the smooth group $\Liegroup{}$
under the right-translation of $\TopgG{\Topgroup{}}$ by $\g{}$ is also an open (resp. closed) set of $\Liegroup{}$. That is,
\begin{align}
&\Foreach{\U}{\atopology{}}
\func{\image{\[\grtrans{\TopgG{\Topgroup{}}}{\g{}}\]}}{\U}\in\atopology{},\\
&\Foreach{\V}{\(\compl{\G{}}{\atopology{}}\)}
\func{\image{\[\grtrans{\TopgG{\Topgroup{}}}{\g{}}\]}}{\V}\in\(\compl{\G{}}{\atopology{}}\).
\end{align}
\item[$\centerdot$]
The image of every open (resp.~closed) set of the smooth group $\Liegroup{}$
under the $\g{}$-conjugation of $\TopgG{\Topgroup{}}$ is also an open (resp. closed) set of $\Liegroup{}$. That is,
\begin{align}
&\Foreach{\U}{\atopology{}}
\func{\image{\[\gconj{\TopgG{\Topgroup{}}}{\g{}}\]}}{\U}\in\atopology{},\\
&\Foreach{\V}{\(\compl{\G{}}{\atopology{}}\)}
\func{\image{\[\gconj{\TopgG{\Topgroup{}}}{\g{}}\]}}{\V}\in\(\compl{\G{}}{\atopology{}}\).
\end{align}
\end{itemize}
\proof
By considering that $\ginv{\TopgG{\Topgroup{}}}$, $\gltrans{\TopgG{\Topgroup{}}}{\g{}}$, $\grtrans{\TopgG{\Topgroup{}}}{\g{}}$,
and $\gconj{\TopgG{\Topgroup{}}}{\g{}}$ are homeomorphisms from the topological-space $\topgtops{\Topgroup{}}$
to $\topgtops{\Topgroup{}}$, and every homeomorphism between topological-spaces
maps each open (resp.~closed) set of the source space onto an open (resp.~closed) set of the target space,
and $\atopology{}$ is the corresponded topology of the topological-space $\topgtops{\Topgroup{}}$,
the verification of the regarded assertions is straightforward.
\endlem
\lemma\label{lemproductfopensetsisopeninatopologicalgroup}
The product of every pair of a subset $\U_1$ and an open set $\U_2$ of the topological group $\Topgroup{}$ is an open set of $\Topgroup{}$.
\begin{equation}
\Foreach{\opair{\U_1}{\U_2}}{\Cprod{\CSs{\G{}}}{\atopology{}}}
\gsetprod{\U_1}{\U_2}{\Group{}}\in\atopology{}.
\end{equation}
\proof
$\U_1$ is taken as an arbitrary subset of $\G{}$ and and $\U_2$ is taken as an arbitrary element of $\atopology{}$.
If $\U_1$ is empty, trivially $\gsetprod{\U_1}{\U_2}{\Group{}}$ is the empty-set
and hence an open set of $\Topgroup{}$. If both $\U_1$ is non-empty, then clearly,
\begin{equation}
\gsetprod{\U_1}{\U_2}{\Group{}}=\Union{\g{}}{\U_1}{\func{\image{\[\gltrans{\TopgG{\Topgroup{}}}{\g{}}\]}}{\U_2}},
\end{equation}
and thus according to \reflem{lemtranslationsofanopensetisopeninatopologicalgroup},
$\gsetprod{\U_1}{\U_2}{\Group{}}$ is the union of a collection of open sets and hence an open set of $\Topgroup{}$.
\endlem
\lemma\label{lempowersofanopensetisopen}
For every open set $\V$ of $\Topgroup{}$ and every positive integer $n$, $\gpower{\V}{\gop{}}{n}$ is an open set of $\Topgroup{}$.
\begin{equation}
\Foreach{\V}{\atopology{}}
\Foreach{n}{\Zp}
\gpower{\V}{\gop{}}{n}\in\atopology{}.
\end{equation}
\proof
According to \reflem{lemproductfopensetsisopeninatopologicalgroup} and
\refdef{defgrouptranslations}, it is trivial.
\endlem
\section{Nuclei of a Topological Group}
\definition\label{defnucleioftopologicalgroup}
The set of all open sets of $\topgtops{\Topgroup{}}$ containing the identity element of the group
$\TopgG{\Topgroup{}}$, that is, the set of all neighborhoods of $\IG{}$ in the topological-space
$\topgtops{\Topgroup{}}$ is denoted by $\nuclei{\Topgroup{}}$. That is,
\begin{align}
\nuclei{\Topgroup{}}:=\func{\nei{\topgtops{\Topgroup{}}}}{\seta{\IG{}}}
=\defset{\U}{\atopology{}}{\IG{}\in\U}.
\end{align}
$\nuclei{\Topgroup{}}$ is referred to as the $\quotl$nuclei of the topological group $\Topgroup{}$$\quotr$,
and each element of $\nuclei{\Topgroup{}}$ is called a $\quotl$nucleus of the topological group $\Liegroup{}$$\quotr$.
Also a subset of $\G{}$ is called a $\quotl$symmetric nucleus of $\Topgroup{}$$\quotr$ iff it is simultaneously a nucleus of $\Topgroup{}$
and a symmetric set of the group $\TopgG{\Topgroup{}}$.
\endef
\lemma\label{lemeachnucleusincludesarestrictednucleus}
Every nucleus $\U$ of $\Liegroup{}$ includes a nucleus $\V$ of $\Topgroup{}$ such that
for every $\g{1}$ and $\g{2}$ in $\V$, $\g{1}\gop{}\g{2}$ belongs to $\U$.
\begin{equation}
\Foreach{\U}{\nuclei{\Topgroup{}}}
\Exists{\V}{\nuclei{\Topgroup{}}}
\bigg(\V\subseteq\U,~\func{\image{\gop{}}}{\Cprod{\V}{\V}}\subseteq\U\bigg),
\end{equation}
or equivalently,
\begin{equation}
\Foreach{\U}{\nuclei{\Topgroup{}}}
\Exists{\V}{\nuclei{\Topgroup{}}}
\bigg(\V\subseteq\U,~\gpower{\V}{\gop{}}{2}\subseteq\U\bigg).
\end{equation}
\proof
$\U$ is taken as an element of $\nuclei{\Topgroup{}}$. Thus according to \refdef{defnucleioftopologicalgroup},
\begin{equation}
\IG{}\in\U\in\atopology{}.
\end{equation}
So since $\U$ is an open set of $\topgtops{\Topgroup{}}$ and the group-operation of $\TopgG{\Topgroup{}}$, that is $\gop{}$,
is a continuous map from $\topprod{\topgtops{\Topgroup{}}}{\topgtops{\Topgroup{}}}$ to
$\topgtops{\Topgroup{}}$ (\refdef{deftopologicalgroup}), evidently the definition of continuity implies that
the inverse-image of $\U$ under $\gop{}$ is an open set of $\topprod{\topgtops{\Topgroup{}}}{\topgtops{\Topgroup{}}}$.
Also, since $\U$ contains $\IG{}$ and $\IG{}\gop{}\IG{}=\IG{}$, $\opair{\IG{}}{\IG{}}$ belongs to the
inverse-image of $\U$ under $\gop{}$. So briefly,
\begin{equation}
\opair{\IG{}}{\IG{}}\in\func{\pimage{\gop{}}}{\U}\in
\topologyofspace{\topprod{\topgtops{\Topgroup{}}}{\topgtops{\Topgroup{}}}}.
\end{equation}
Consequently, since the set $\defSet{\Cprod{\U_1}{\U_2}}{\U_1,~\U_2\in\atopology{}}$ is a base for
the topological-space $\topprod{\topgtops{\Topgroup{}}}{\topgtops{\Topgroup{}}}$,
\begin{equation}
\Existsis{\opair{\U_1}{\U_2}}{\Cprod{\atopology{}}{\atopology{}}}
\opair{\IG{}}{\IG{}}\in\Cprod{\U_1}{\U_2}\subseteq
\func{\pimage{\gop{}}}{\U}.
\end{equation}
Thus, each $\U_1$, $\U_2$, and $\U$ is an open set of $\topgtops{\Topgroup{}}$ containing $\IG{}$, and hence
so is their intersection. Thus the intersection of $\U_1$, $\U_2$, and $\U$ is a nucleus of $\Topgroup{}$, that is,
\begin{equation}
\U_1\cap\U_2\cap\U\in\nuclei{\Topgroup{}}.
\end{equation} 
Additionally, since $\Cprod{\U_1}{\U_2}\subseteq\func{\pimage{\gop{}}}{\U}$, evidently,
\begin{equation}
\func{\image{\gop{}}}{\Cprod{\U_1}{\U_2}}\subseteq\U,
\end{equation}
and hence as a trivial consequence,
\begin{equation}
\func{\image{\gop{}}}{\Cprod{\(\U_1\cap\U_2\cap\U\)}{\(\U_1\cap\U_2\cap\U\)}}\subseteq\U.
\end{equation}
\endthm
\lemma\label{lemeachnucleusincludesasymmetricnucleus}
Every nucleus $\U$ of $\Liegroup{}$ includes a symmetric nucleus $\V$ of $\Topgroup{}$. That is,
\begin{equation}
\Foreach{\U}{\nuclei{\Topgroup{}}}
\Exists{\V}{\nuclei{\Topgroup{}}}
\bigg(\V\subseteq\U,~\func{\image{\ginv{\TopgG{\Topgroup{}}}}}{\V}=\V\bigg).
\end{equation}
\proof
$\U$ is taken as an element of $\nuclei{\Topgroup{}}$.
According to \refthm{thminversemappingishomeomorphism}, $\ginv{\TopgG{\Topgroup{}}}$ is a homeomorphism from $\topgtops{\Topgroup{}}$
to itself, and thus maps open sets onto open sets. Therefore, since $\U$ is an open set of $\topgtops{\Topgroup{}}$,
so is $\func{\image{\[\ginv{\TopgG{\Topgroup{}}}\]}}{\U}$ and accordingly their intersection. That is,
\begin{equation}\label{lemeachnucleusincludesasymmetricnucleuspeq1}
\V:=\U\cap\func{\image{\[\ginv{\TopgG{\Topgroup{}}}\]}}{\U}\in\atopology{}.
\end{equation}
In addition, since $\IG{}\in\U$, obviously $\IG{}\in\func{\image{\[\ginv{\TopgG{\Topgroup{}}}\]}}{\U}$,
and thus $\IG{}\in\V$, which together with \Ref{lemeachnucleusincludesasymmetricnucleuspeq1} implies that
$\V$ is also e nucleus of $\Topgroup{}$ included in $\U$.
\begin{equation}\label{lemeachnucleusincludesasymmetricnucleuspeq2}
\U\supseteq\V\in\nuclei{\Topgroup{}}.
\end{equation}
Moreover, since $\ginv{\TopgG{\Topgroup{}}}$ is a bijection and
$\finv{\[\ginv{\TopgG{\Topgroup{}}}\]}=\ginv{\TopgG{\Topgroup{}}}$,
\begin{align}
\func{\image{\ginv{\TopgG{\Topgroup{}}}}}{\V}&=
\func{\image{\ginv{\TopgG{\Topgroup{}}}}}{\U}\cap
\func{\image{\ginv{\TopgG{\Topgroup{}}}}}{\func{\image{\[\ginv{\TopgG{\Topgroup{}}}\]}}{\U}}\cr
&=\func{\image{\ginv{\TopgG{\Topgroup{}}}}}{\U}\cap\U\cr
&=\V,
\end{align}
which means $\V$ is a symmetric set of the group $\TopgG{\Topgroup{}}$.
\endlem
\lemma\label{lemeverynucleushasasymmetricnucleuspower}
For every nucleus $\U$ of $\Topgroup{}$ and every positive integer $n$, there exists a symmetric nucleus of $\Topgroup{}$
whose first and $n$-th power is included in $\U$.
\begin{equation}
\Foreach{\U}{\nuclei{\Topgroup{}}}
\Foreach{n}{\Zp}
\Exists{\V}{\nuclei{\Topgroup{}}}
\bigg(
\V\subseteq\U,~\gpower{\V}{\gop{}}{n}\subseteq\U,~\func{\image{\ginv{\TopgG{\Topgroup{}}}}}{\V}=\V
\bigg).
\end{equation}
\proof
$\U$ is taken as an arbitrary element of $\nuclei{\Topgroup{}}$.
According to \reflem{lemeachnucleusincludesasymmetricnucleus},
\begin{equation}\label{lemeverynucleushasasymmetricnucleuspowerpeq1}
\Exists{\V}{\nuclei{\Topgroup{}}}
\bigg(\V\subseteq\U,~\func{\image{\ginv{\TopgG{\Topgroup{}}}}}{\V}=\V\bigg).
\end{equation}
\begin{itemize}
\item[$\pr{1}$]
$n$ is taken as  an arbitrary positive integer, and it is assumed that,
\begin{equation}
\Existsis{\V}{\nuclei{\Topgroup{}}}
\bigg(
\V\subseteq\U,~\gpower{\V}{\gop{}}{n}\subseteq\U,~\func{\image{\ginv{\TopgG{\Topgroup{}}}}}{\V}=\V
\bigg).
\end{equation}
Then according to \reflem{lemeachnucleusincludesarestrictednucleus},
\begin{equation}
\Existsis{\V_1}{\nuclei{\Topgroup{}}}
\bigg(\V_1\subseteq\V,~\gpower{\V_1}{\gop{}}{2}\subseteq\V\bigg).
\end{equation}
Hence,
\begin{equation}
\gpower{\V_1}{\gop{}}{n+1}=\gsetprod{\gpower{\V_1}{\gop{}}{2}}{\gpower{\V_1}{\gop{}}{n-1}}{\Group{}}\subseteq
\gsetprod{\V}{\gpower{\V}{\gop{}}{n-1}}{\Group{}}=
\gpower{\V}{\gop{}}{n}\subseteq\U.
\end{equation}
In addition, according to \reflem{lemeachnucleusincludesasymmetricnucleus},
\begin{equation}
\Existsis{\V_2}{\nuclei{\Topgroup{}}}
\bigg(\V_2\subseteq\V_1,~\func{\image{\ginv{\TopgG{\Topgroup{}}}}}{\V_2}=\V_2\bigg).
\end{equation}
Therefore,
\begin{align}
&\V_2\subseteq\V_1\subseteq\V\subseteq\U,\cr
&\gpower{\V_2}{\gop{}}{n+1}\subseteq\gpower{\V_1}{\gop{}}{n+1}\subseteq\U.
\end{align}
\endp
\end{itemize}
Therefore,
\begin{align}\label{lemeverynucleushasasymmetricnucleuspowerpeq2}
&\[\Exists{\V}{\nuclei{\Topgroup{}}}
\bigg(
\V\subseteq\U,~\gpower{\V}{\gop{}}{n}\subseteq\U,~\func{\image{\ginv{\TopgG{\Topgroup{}}}}}{\V}=\V
\bigg)\]\cr
\then&\[\Exists{\V}{\nuclei{\Topgroup{}}}
\bigg(
\V\subseteq\U,~\gpower{\V}{\gop{}}{n+1}\subseteq\U,~\func{\image{\ginv{\TopgG{\Topgroup{}}}}}{\V}=\V
\bigg)\].
\end{align}
\Ref{lemeverynucleushasasymmetricnucleuspowerpeq1} and \Ref{lemeverynucleushasasymmetricnucleuspowerpeq2}
imply by induction the intended assertion.
\endlem
\section{Topological Subgroups}
\definition\label{deftopologicalsubgroup}
$\asubgroup{}$ is taken as a subset of $\G{}$. $\asubgroup{}$ is referred to as a $\quotl$topological subgroup of the topological group
$\Topgroup{}$$\quotr$ iff these properties are satisfied.
\begin{itemize}
\item[\myitem{TS~1.}]
$\asubgroup{}$ is a subgroup of the intrinsic group $\TopgG{\Topgroup{}}$ of the topological group $\Topgroup{}$.
\begin{equation}
\asubgroup{}\in\Subgroups{\TopgG{\Topgroup{}}}.
\end{equation}
\item[\myitem{TS~2.}]
The set $\asubgroup{}$ endowed with the group operation inherited from that of $\TopgG{\Topgroup{}}$ and the topology inherited from that of
$\topgtops{\Topgroup{}}$ (the underlying topological-space of $\Topgroup{}$) on it, is a topological group. That is, the triple
$\triple{\asubgroup{}}{\func{\res{\gop{}}}{\Cprod{\asubgroup{}}{\asubgroup{}}}}{\func{\IndTop{\topgtops{\Topgroup{}}}}{\asubgroup{}}}$
is a topological group.
\end{itemize}
The set of all topological subgroups of the topological group $\Topgroup{}$ is denoted by $\Topsubgroups{\Topgroup{}}$.
\endef
\theorem\label{thmsubgroupsaretopologicalsubgroups}
Every subgroup $\asubgroup{}$ of the intrinsic group structure $\TopgG{\Topgroup{}}$ of the topological group $\Topgroup{}$,
is a topological subfroup of $\Topgroup{}$. This means the set of subgroups of the intrinsic group of $\Topgroup{}$
coincides with the set of all topological subgroups of $\Topgroup{}$.
\begin{equation}
\Topsubgroups{\Topgroup{}}=\Subgroups{\TopgG{\Topgroup{}}}.
\end{equation}
\proof
According to \refdef{deftopologicalsubgroup}, it is obvious that,
\begin{equation}
\Topsubgroups{\Topgroup{}}\subseteq\Subgroups{\TopgG{\Topgroup{}}}.
\end{equation}
\begin{itemize}
\item[$\pr{1}$]
$\asubgroup{}$ is taken as an arbitrary element of $\Subgroups{\TopgG{\Topgroup{}}}$ (a subgroup of
$\TopgG{\Topgroup{}}$). Then, clearly $\opair{\asubgroup{}}{\func{\res{\gop{}}}{\Cprod{\asubgroup{}}{\asubgroup{}}}}$
is a group.\\
Moreover, since $\Topgroup{}=\triple{\G{}}{\gop{}}{\atopology{}}$ is a topological group, according to \refdef{deftopologicalgroup},
$\gop{}$ is a continuous map from $\topprod{\opair{\G{}}{\atopology{}}}{\opair{\G{}}{\atopology{}}}$ to
$\opair{\G{}}{\atopology{}}$, and $\ginv{\TopgG{\Topgroup{}}}$ is a continuous map from $\opair{\G{}}{\atopology{}}$
to $\opair{\G{}}{\atopology{}}$. Therefore, $\func{\res{\gop{}}}{\Cprod{\asubgroup{}}{\asubgroup{}}}$ is a continuous map
from the topological-subspace
$\topprod{\opair{\asubgroup{}}{\func{\IndTop{\topgtops{\Topgroup{}}}}{\asubgroup{}}}}
{\opair{\asubgroup{}}{\func{\IndTop{\topgtops{\Topgroup{}}}}{\asubgroup{}}}}$ of the
topological-space $\topprod{\opair{\G{}}{\atopology{}}}{\opair{\G{}}{\atopology{}}}$ to the topological-subspace
$\opair{\asubgroup{}}{\func{\IndTop{\topgtops{\Topgroup{}}}}{\asubgroup{}}}$ of $\topgtops{\Topgroup{}}$.
Also $\func{\res{\ginv{\TopgG{\Topgroup{}}}}}{\asubgroup{}}$ becomes a continuous map from
$\opair{\asubgroup{}}{\func{\IndTop{\topgtops{\Topgroup{}}}}{\asubgroup{}}}$ to
$\opair{\asubgroup{}}{\func{\IndTop{\topgtops{\Topgroup{}}}}{\asubgroup{}}}$. But,
\begin{equation}
\func{\res{\ginv{\TopgG{\Topgroup{}}}}}{\asubgroup{}}
=\ginv{\opair{\asubgroup{}}{\func{\res{\gop{}}}{\Cprod{\asubgroup{}}{\asubgroup{}}}}},
\end{equation}
and hence,
\begin{align}
\func{\res{\gop{}}}{\Cprod{\asubgroup{}}{\asubgroup{}}}&\in
\CF{\topprod{\opair{\asubgroup{}}{\func{\IndTop{\topgtops{\Topgroup{}}}}{\asubgroup{}}}}
{\opair{\asubgroup{}}{\func{\IndTop{\topgtops{\Topgroup{}}}}{\asubgroup{}}}}}
{\opair{\asubgroup{}}{\func{\IndTop{\topgtops{\Topgroup{}}}}{\asubgroup{}}}},\\
\ginv{\opair{\asubgroup{}}{\func{\res{\gop{}}}{\Cprod{\asubgroup{}}{\asubgroup{}}}}}&\in
\CF{\opair{\asubgroup{}}{\func{\IndTop{\topgtops{\Topgroup{}}}}{\asubgroup{}}}}
{\opair{\asubgroup{}}{\func{\IndTop{\topgtops{\Topgroup{}}}}{\asubgroup{}}}}.
\end{align}
Therefore, according to \refdef{deftopologicalgroup},
$\triple{\asubgroup{}}{\func{\res{\gop{}}}{\Cprod{\asubgroup{}}{\asubgroup{}}}}{\func{\IndTop{\topgtops{\Topgroup{}}}}{\asubgroup{}}}$
is a topological group and hence,
\begin{equation}
\asubgroup{}\in\Topsubgroups{\Topgroup{}}.
\end{equation}
\endp
\end{itemize}
\endthm
\theorem\label{thmopentopologicalsubgroupsareclosed}
$\asubgroup{}$ is taken as an element of $\Subgroups{\TopgG{\Topgroup{}}}$ (a subgroup of $\TopgG{\Topgroup{}}$).
If $\asubgroup{}$ is an open set of the topological group $\Topgroup{}$,
then $\asubgroup{}$ is also a closed set of $\topgtops{\Topgroup{}}$.
\begin{equation}
\asubgroup{}\in\atopology{}\then
\asubgroup{}\in\closedsets{\topgtops{\Topgroup{}}}.
\end{equation}
\proof
It is assumed that $\asubgroup{}\in\atopology{}$.
\begin{itemize}
\item[$\pr{1}$]
$\SET{}$ is taken as an arbitrary element of $\func{\LCoset{\TopgG{\Topgroup{}}}}{\asubgroup{}}$, that is an arbitrary left-coset of
the subgroup $\asubgroup{}$ of $\TopgG{\Topgroup{}}$. Then there exists an element $\g{}$ of $\G{}$ such that,
\begin{equation}
\SET{}=\g{}\asubgroup{}=\func{\image{\[\gltrans{\TopgG{\Topgroup{}}}{\g{}}\]}}{\asubgroup{}},
\end{equation}
and thus according to \reflem{lemtranslationsofanopensetisopeninatopologicalgroup}, $\SET{}$ is an open set of $\topgtops{\Topgroup{}}$.
\begin{equation}
\SET{}\in\atopology{}.
\end{equation}
\endp
\end{itemize}
Therefore, all left-cosets of the subgroup $\asubgroup{}$ of $\TopgG{\Topgroup{}}$ are open sets of
$\topgtops{\Topgroup{}}$.
\begin{equation}
\func{\LCoset{\TopgG{\Topgroup{}}}}{\asubgroup{}}\subseteq\atopology{}.
\end{equation}
Moreover, since $\func{\LCoset{\TopgG{\Topgroup{}}}}{\asubgroup{}}$ is a partition of $\G{}$,
and $\asubgroup{}\in\func{\LCoset{\TopgG{\Topgroup{}}}}{\asubgroup{}}$, evidently,
the complement of $\asubgroup{}$ in $\G{}$ is the union of all left-cosets of $\asubgroup{}$ in $\TopgG{\Topgroup{}}$
distinct from $\asubgroup{}$ itself.
\begin{equation}
\compl{\G{}}{\asubgroup{}}=\union{\bigg(\compl{\func{\LCoset{\TopgG{\Topgroup{}}}}{\asubgroup{}}}{\seta{\asubgroup{}}}\bigg)}.
\end{equation}
Since all elements of $\compl{\func{\LCoset{\TopgG{\Topgroup{}}}}{\asubgroup{}}}{\seta{\asubgroup{}}}$ are open sets of $\topgtops{\Topgroup{}}$,
so is $\union{\bigg(\compl{\func{\LCoset{\TopgG{\Topgroup{}}}}{\asubgroup{}}}{\seta{\asubgroup{}}}\bigg)}$. Therefore $\compl{\G{}}{\asubgroup{}}$
is an open set of $\topgtops{\Topgroup{}}$, and hence $\asubgroup{}$ is a closed set of $\topgtops{\Topgroup{}}$.
\endthm
\theorem\label{thmclosureofatopologicalsubgroupisatopologicalsubgroup}
For every subgroup $\asubgroup{}$ of $\TopgG{\Topgroup{}}$), the closure of $\asubgroup{}$ in the topological-space $\topgtops{\Topgroup{}}$
is also a subgroup of $\TopgG{\Topgroup{}}$.
\begin{equation}
\Foreach{\asubgroup{}}{\Subgroups{\TopgG{\Topgroup{}}}}
\func{\Cl{\topgtops{\Topgroup{}}}}{\asubgroup{}}\in
\Subgroups{\TopgG{\Topgroup{}}}.
\end{equation}
\proof
$\asubgroup{}$ is taken as an element of $\Subgroups{\TopgG{\Topgroup{}}}$ (a subgroup of $\TopgG{\Topgroup{}}$).
Then, clearly $\asubgroup{}$ is a non-empty set such that,
\begin{equation}\label{thmclosureofatopologicalsubgroupisatopologicalsubgrouppeq1}
\Foreach{\opair{\g{}}{\hh}}{\Cprod{\asubgroup{}}{\asubgroup{}}}
\func{\gopr{\gop{}}{\TopgG{\Topgroup{}}}}{\binary{\g{}}{\hh}}=
\g{}\gop{}\invG{\hh}{}\in\asubgroup{}.
\end{equation}
Since the closure of $\asubgroup{}$ in the topological-space $\topgtops{\Topgroup{}}$ includes $\asubgroup{}$,
it is also a non-empty set.
\begin{equation}\label{thmclosureofatopologicalsubgroupisatopologicalsubgrouppeq2}
\func{\Cl{\topgtops{\Topgroup{}}}}{\asubgroup{}}\neq\empty.
\end{equation}
\begin{itemize}
\item[$\pr{1}$]
Each $\g{}$ and $\hh$ is taken as an arbitrary element of $\func{\Cl{\topgtops{\Topgroup{}}}}{\asubgroup{}}$.
Then every neighborhood of each one of them in the topological-space $\topgtops{\Topgroup{}}$, intersects $\asubgroup{}$. That is,
\begin{align}
&\Foreach{\U}{\func{\nei{\topgtops{\Topgroup{}}}}{\seta{\g{}}}}\U\cap\asubgroup{}\neq\empty,
\label{thmclosureofatopologicalsubgroupisatopologicalsubgroupp1eq1}\\
&\Foreach{\U}{\func{\nei{\topgtops{\Topgroup{}}}}{\seta{\hh}}}\U\cap\asubgroup{}\neq\empty.
\label{thmclosureofatopologicalsubgroupisatopologicalsubgroupp1eq2}
\end{align}
\begin{itemize}
\item[$\pr{1-1}$]
$\V$ is taken as an arbitrary element of
$\func{\nei{\topgtops{\Topgroup{}}}}{\seta{\func{\gopr{\gop{}}{\TopgG{\Topgroup{}}}}{\binary{\g{}}{\hh}}}}$,
that is an open set of $\topgtops{\Topgroup{}}$ containing $\g{}\gop{}\invG{\hh}{}$.
Since $\gopr{\gop{}}{\TopgG{\Topgroup{}}}$ is a continuous map from $\topgtops{\Topgroup{}}$ to
$\topgtops{\Topgroup{}}$ (\refthm{thmtopologicalgroupequiv0}) and $\V$ is an open set of $\topgtops{\Topgroup{}}$, according to the definition of
continuity, the inverse-image of $\V$ under $\gopr{\gop{}}{\TopgG{\Topgroup{}}}$ is an open set of
$\topprod{\topgtops{\Topgroup{}}}{\topgtops{\Topgroup{}}}$. That is,
\begin{equation}\label{thmclosureofatopologicalsubgroupisatopologicalsubgroupp11eq1}
\func{\pimage{\gopr{\gop{}}{\TopgG{\Topgroup{}}}}}{\V}\in\topologyofspace{\topprod{\topgtops{\Topgroup{}}}{\topgtops{\Topgroup{}}}},
\end{equation}
where
$\topprod{\topgtops{\Topgroup{}}}{\topgtops{\Topgroup{}}}=
\opair{\Cprod{\G{}}{\G{}}}{\topologyofspace{\topprod{\topgtops{\Topgroup{}}}{\topgtops{\Topgroup{}}}}}$.
Also, since $\func{\gopr{\gop{}}{\TopgG{\Topgroup{}}}}{\binary{\g{}}{\hh}}\in\V$, evidently,
\begin{equation}\label{thmclosureofatopologicalsubgroupisatopologicalsubgroupp11eq2}
\opair{\g{}}{\hh}\in
\func{\pimage{\gopr{\gop{}}{\TopgG{\Topgroup{}}}}}{\V}.
\end{equation}
Therefore, according to \Ref{thmclosureofatopologicalsubgroupisatopologicalsubgroupp11eq1} and
\Ref{thmclosureofatopologicalsubgroupisatopologicalsubgroupp11eq2}, and considering that
$\defSet{\Cprod{\U_1}{\U_2}}{\U_1,~\U_2\in\topology{}}$ is a topological-basis for the
topological-space $\topprod{\topgtops{\Topgroup{}}}{\topgtops{\Topgroup{}}}$, it is clear that
there exists a neighborhood $\U_1$ of $\g{}$ and a neighborhood $\U_2$ of $\hh$ (in $\topgtops{\Topgroup{}}$)
such that $\Cprod{\U_1}{\U_2}$ is included in $\func{\pimage{\gopr{\gop{}}{\TopgG{\Topgroup{}}}}}{\V}$.
\begin{align}\label{thmclosureofatopologicalsubgroupisatopologicalsubgroupp11eq3}
\Existsis{\opair{\U_1}{\U_2}}{\Cprod{\func{\nei{\topgtops{\Topgroup{}}}}{\seta{\g{}}}}
{\func{\nei{\topgtops{\Topgroup{}}}}{\seta{\hh}}}}\Cprod{\U_1}{\U_2}\subseteq
\func{\pimage{\gopr{\gop{}}{\TopgG{\Topgroup{}}}}}{\V}.
\end{align}
Thus,
\begin{equation}\label{thmclosureofatopologicalsubgroupisatopologicalsubgroupp11eq4}
\func{\image{\gopr{\gop{}}{\TopgG{\Topgroup{}}}}}{\Cprod{\U_1}{\U_2}}\subseteq\V,
\end{equation}
and according to \Ref{thmclosureofatopologicalsubgroupisatopologicalsubgroupp1eq1} and
\Ref{thmclosureofatopologicalsubgroupisatopologicalsubgroupp1eq2},
\begin{align}
\U_1\cap\asubgroup{}&\neq\empty,
\label{thmclosureofatopologicalsubgroupisatopologicalsubgroupp11eq5}\\
\U_2\cap\asubgroup{}&\neq\empty.
\label{thmclosureofatopologicalsubgroupisatopologicalsubgroupp11eq6}
\end{align}
Hence,
\begin{equation}\label{thmclosureofatopologicalsubgroupisatopologicalsubgroupp11eq7}
\Existsis{\opair{\g{1}}{\g{2}}}{\Cprod{\asubgroup{}}{\asubgroup{}}}
\g{1}\in\U_1,~\g{2}\in\U_2.
\end{equation}
Therefore, according to \Ref{thmclosureofatopologicalsubgroupisatopologicalsubgrouppeq1},
\begin{equation}
\func{\gopr{\gop{}}{\TopgG{\Topgroup{}}}}{\binary{\g{1}}{\g{2}}}\in\asubgroup{},
\end{equation}
and according to \Ref{thmclosureofatopologicalsubgroupisatopologicalsubgroupp11eq4},
\begin{equation}
\func{\gopr{\gop{}}{\TopgG{\Topgroup{}}}}{\binary{\g{1}}{\g{2}}}\in\V,
\end{equation}
and therefore,
\begin{equation}
\V\cap\asubgroup{}\neq\empty.
\end{equation}
\endp
\end{itemize}
Therefore,
\begin{equation}
\Foreach{\V}
{\func{\nei{\topgtops{\Topgroup{}}}}{\seta{\func{\gopr{\gop{}}{\TopgG{\Topgroup{}}}}{\binary{\g{}}{\hh}}}}}
\V\cap\asubgroup{}\neq\empty,
\end{equation}
which means $\func{\gopr{\gop{}}{\TopgG{\Topgroup{}}}}{\binary{\g{}}{\hh}}$ is in the closure of $\asubgroup{}$.
\endp
\end{itemize}
Therefore,
\begin{equation}
\Foreach{\opair{\g{}}{\hh}}{\Cprod{\func{\Cl{\topgtops{\Topgroup{}}}}{\asubgroup{}}}{\func{\Cl{\topgtops{\Topgroup{}}}}{\asubgroup{}}}}
\func{\gopr{\gop{}}{\TopgG{\Topgroup{}}}}{\binary{\g{}}{\hh}}\in
\func{\Cl{\topgtops{\Topgroup{}}}}{\asubgroup{}},
\end{equation}
which, by considering that $\func{\gopr{\gop{}}{\TopgG{\Topgroup{}}}}{\binary{\g{}}{\hh}}$ is non-empty,
implies that $\func{\gopr{\gop{}}{\TopgG{\Topgroup{}}}}{\binary{\g{}}{\hh}}$ is a subgroup of the underlying topological structure
$\TopgG{\Topgroup{}}$ of $\Topgroup{}$. That is,
\begin{equation}
\func{\gopr{\gop{}}{\TopgG{\Topgroup{}}}}{\binary{\g{}}{\hh}}\in\Subgroups{\TopgG{\Topgroup{}}}.
\end{equation}
\endthm
\theorem\label{thmlocallyclosedtopologicalsubgroupisclosed}
$\asubgroup{}$ is taken as an element of $\Subgroups{\TopgG{\Topgroup{}}}$ (a subgroup of $\TopgG{\Topgroup{}}$).
If $\asubgroup{}$ is a locally-closed set of the underlying topological-space $\topgtops{\Topgroup{}}$ of the smooth group $\Topgroup{}$,
that is if $\asubgroup{}$ is an open set of its closure (endowed with the inherited topology from $\topgtops{\Topgroup{}}$),
then $\asubgroup{}$ is also a closed set of $\topgtops{\Topgroup{}}$.
\begin{equation}
\asubgroup{}\in\func{\IndTop{\topgtops{\Topgroup{}}}}{\func{\Cl{\topgtops{\Topgroup{}}}}{\asubgroup{}}}
\then
\asubgroup{}\in\closedsets{\topgtops{\Topgroup{}}}.
\end{equation}
\proof
It is assumed that
$\asubgroup{}$ is an open set of its closure endowed with its topology inherited from $\topgtops{\Topgroup{}}$, that is,
\begin{equation}
\asubgroup{}\in\func{\IndTop{\topgtops{\Topgroup{}}}}{\func{\Cl{\topgtops{\Topgroup{}}}}{\asubgroup{}}}.
\end{equation}
Since $\asubgroup{}$ is a subgroup of $\TopgG{\Topgroup{}}$, according to \refthm{thmclosureofatopologicalsubgroupisatopologicalsubgroup},
$\func{\Cl{\topgtops{\Topgroup{}}}}{\asubgroup{}}$ is also a subgroup of $\TopgG{\Topgroup{}}$, and hence according to
\refthm{thmsubgroupsaretopologicalsubgroups}, $\func{\Cl{\topgtops{\Topgroup{}}}}{\asubgroup{}}$ is also a topological subgroup of
$\Topgroup{}$. Thus, according to \refdef{deftopologicalsubgroup}, the triple
$\triple{\asubgroup{}}{\func{\res{\gop{}}}{\Cprod{\func{\Cl{\topgtops{\Topgroup{}}}}{\asubgroup{}}}{\func{\Cl{\topgtops{\Topgroup{}}}}{\asubgroup{}}}}}
{\func{\IndTop{\topgtops{\Topgroup{}}}}{\func{\Cl{\topgtops{\Topgroup{}}}}{\asubgroup{}}}}$
is a topological group. Furthermore, it is trivial that since $\asubgroup{}$ is a subgroup of the group
$\opair{\G{}}{\gop{}}$, it is also a subgroup of the intrinsic group of this topological group,
which is the group
$\opair{\func{\Cl{\topgtops{\Topgroup{}}}}{\asubgroup{}}}
{\func{\res{\gop{}}}{\Cprod{\func{\Cl{\topgtops{\Topgroup{}}}}{\asubgroup{}}}{\func{\Cl{\topgtops{\Topgroup{}}}}{\asubgroup{}}}}}$.
In addition, since $\asubgroup{}$
is an open set of this topological group, \refthm{thmopentopologicalsubgroupsareclosed} implies that it is also a closed set of
this topological group. That is,
\begin{equation}
\asubgroup{}\in\closedsets{\opair{\func{\Cl{\topgtops{\Topgroup{}}}}{\asubgroup{}}}
{\func{\IndTop{\topgtops{\Topgroup{}}}}{\func{\Cl{\topgtops{\Topgroup{}}}}{\asubgroup{}}}}}.
\end{equation}
Accordingly, considering that $\func{\Cl{\topgtops{\Topgroup{}}}}{\asubgroup{}}$ is itself a closed set of the topological-space
$\topgtops{\Topgroup{}}$, it becomes evident that $\asubgroup{}$ is also a closed set of $\topgtops{\Topgroup{}}$,
via the transitivity property of closed (or open) sets of a topological-space.
\endthm
\lemma\label{lemsubgroupgeneratedbyanonemptysymmetricopenset}
For every non-empty, symmetric, and open set $\V$ of the topological group $\Topgroup{}$,
$\displaystyle\Unionn{n}{1}{\infty}{\gpower{\V}{\gop{}}{n}}$ is simultaneously an open-and-closed set and
a topological subgroup of $\Topgroup{}$.
\begin{equation}
\Foreach{\V}{\defset{\U}{\compl{\atopology{}}{\seta{\empty}}}{\func{\image{\ginv{\TopgG{\Topgroup{}}}}}{\U}=\U}}
\Unionn{n}{1}{\infty}{\gpower{\V}{\gop{}}{n}}\in\atopology{}\cap\closedsets{\topgtops{\Topgroup{}}}\cap\Subgroups{\TopgG{\Topgroup{}}}.
\end{equation}
\proof
$\V$ is taken as an arbitray non-empty element of $\atopology{}$ such that,
\begin{equation}\label{lemsubgroupgeneratedbyanonemptysymmetricopensetpeq1}
\func{\image{\ginv{\TopgG{\Topgroup{}}}}}{\V}=\V.
\end{equation}
So clearly,
\begin{equation}\label{lemsubgroupgeneratedbyanonemptysymmetricopensetpeq2}
\Foreach{\g{}}{\V}\invG{\g{}}{}\in\V.
\end{equation}
Since $\V$ is non-emty, obviously $\displaystyle\Unionn{n}{1}{\infty}{\gpower{\V}{\gop{}}{n}}$ is also non-empty.
\begin{itemize}
\item[$\pr{1}$]
Each $\g{}$ and $\hh$ is taken as an arbitrary element of $\displaystyle\Unionn{n}{1}{\infty}{\gpower{\V}{\gop{}}{n}}$.
Hence,
\begin{align}
&\Existsis{m}{\Zp}
\Existsis{\mtuple{\g{1}}{\g{m}}}{\V^{\times m}}
\g{}=\Succc{\g{1}}{\g{m}}{\gop{}},\\
&\Existsis{l}{\Zp}
\Existsis{\mtuple{\hh_1}{\hh_l}}{\V^{\times l}}
\hh=\Succc{\hh_1}{\hh_l}{\gop{}}.
\end{align}
Hence,
\begin{equation}
\invG{\g{}}{}=\Succc{\invG{\g{m}}{}}{\invG{\g{1}}{}}{\gop{}},
\end{equation}
and according to \Ref{lemsubgroupgeneratedbyanonemptysymmetricopensetpeq2},
\begin{equation}
\invG{\g{1}}{}\in\V,\ldots,\invG{\g{m}}{}\in\V.
\end{equation}
Therefore,
\begin{equation}
\hh\gop{}{}\invG{\g{}}{}\in\V^{\times(l+m)},
\end{equation}
and hence,
\begin{equation}
\hh\gop{}{}\invG{\g{}}{}\in\Unionn{n}{1}{\infty}{\gpower{\V}{\gop{}}{n}}.
\end{equation}
\endp
\end{itemize}
Therefore,
\begin{equation}
\Foreach{\opair{\g{}}{\hh}}{\Cprod{\Unionn{n}{1}{\infty}{\gpower{\V}{\gop{}}{n}}}{\Unionn{n}{1}{\infty}{\gpower{\V}{\gop{}}{n}}}}
\hh\gop{}{}\invG{\g{}}{}\in\Unionn{n}{1}{\infty}{\gpower{\V}{\gop{}}{n}},
\end{equation}
which together with the fact that $\displaystyle\Unionn{n}{1}{\infty}{\gpower{\V}{\gop{}}{n}}$ is non-empty,
implies that $\displaystyle\in\Unionn{n}{1}{\infty}{\gpower{\V}{\gop{}}{n}}$ is a subgroup of $\TopgG{\Topgroup{}}$.
\begin{equation}
\Unionn{n}{1}{\infty}{\gpower{\V}{\gop{}}{n}}\in\Subgroups{\TopgG{\Topgroup{}}}.
\end{equation}
Moreover, considering that $\V$ is an open set of $\Topgroup{}$, so is $\gpower{\V}{\gop{}}{n}$ and accordingly
$\displaystyle\Unionn{n}{1}{\infty}{\gpower{\V}{\gop{}}{n}}$.
\begin{equation}
\Unionn{n}{1}{\infty}{\gpower{\V}{\gop{}}{n}}\in\atopology.
\end{equation}
Since $\displaystyle\Unionn{n}{1}{\infty}{\gpower{\V}{\gop{}}{n}}$ is simultaneously a subgroup
and an open set of of $\Topgroup{}$, based on \refthm{thmopentopologicalsubgroupsareclosed}, it is evident that
$\displaystyle\Unionn{n}{1}{\infty}{\gpower{\V}{\gop{}}{n}}$ is also a closed set of $\Topgroup{}$. That is,
\begin{equation}
\Unionn{n}{1}{\infty}{\gpower{\V}{\gop{}}{n}}\in\closedsets{\topgtops{\Topgroup{}}}.
\end{equation}
\endlem
\theorem\label{thmeverynucleusgeneratestheconnectedtopologicalgroup}
If the underlying topological-space $\topgtops{\Topgroup{}}$ of the topological group $\Topgroup{}$ is connected,
then every nucleus of $\Topgroup{}$ is generates the the intrinsic group $\TopgG{\Topgroup{}}$ of the topological group $\Topgroup{}$.
That is,
\begin{equation}
\G{}\in\connecteds{\topgtops{\Topgroup{}}}\then
\(\Foreach{\U}{\nuclei{\Topgroup{}}}\Union{n}{\Z}{\gpower{\U}{\gop{}}{n}}=\G{}\).
\end{equation}
\proof
It is assumed that $\topgtops{\Topgroup{}}$ is a connected topological-space. Thus, $\empty$ and $\G{}$ are the only
open-and-closed sets of $\topgtops{\Topgroup{}}$. That is,
\begin{equation}\label{thmeverynucleusgeneratestheconnectedtopologicalgrouppeq1}
\atopology{}\cap\closedsets{\topgtops{\Topgroup{}}}=\seta{\binary{\empty}{\G{}}}.
\end{equation}
\begin{itemize}
\item[$\pr{1}$]
$\U$ is taken as an arbitrary element of $\nuclei{\Topgroup{}}$. According to \reflem{lemeachnucleusincludesasymmetricnucleus},
\begin{equation}
\Exists{\V}{\nuclei{\Topgroup{}}}
\bigg(\V\subseteq\U,~\func{\image{\ginv{\TopgG{\Topgroup{}}}}}{\V}=\V\bigg),
\end{equation}
and therefore according to \reflem{lemsubgroupgeneratedbyanonemptysymmetricopenset} and
\refdef{defnucleioftopologicalgroup},
\begin{equation}
\empty\neq\Unionn{n}{1}{\infty}{\gpower{\V}{\gop{}}{n}}\in
\atopology{}\cap\closedsets{\topgtops{\Topgroup{}}},
\end{equation}
and thus according to \Ref{thmeverynucleusgeneratestheconnectedtopologicalgrouppeq1},
\begin{equation}
\Unionn{n}{1}{\infty}{\gpower{\V}{\gop{}}{n}}=\G{}.
\end{equation}
This implies trivially that,
\begin{equation}
\Union{n}{\Z}{\gpower{\V}{\gop{}}{n}}=\G{}.
\end{equation}
So, since $\V\subseteq\U$, trivially $\displaystyle\Union{n}{\Z}{\gpower{\V}{\gop{}}{n}}\subseteq\Union{n}{\Z}{\gpower{\U}{\gop{}}{n}}$, and thus,
\begin{equation}
\Union{n}{\Z}{\gpower{\U}{\gop{}}{n}}=\G{}.
\end{equation}
\endp
\end{itemize}
\endthm
\theorem
The connected component of the underlying topological-space $\topgtops{\Topgroup{}}$ of the topological group $\Topgroup{}$ containing
the identity element $\IG{}$ is
a normanl subgroup of $\TopgG{\Topgroup{}}$, and simultaneously a closed set of $\topgtops{\Topgroup{}}$.
\proof
Let $\asubgroup{}$ denote the connected component of $\topgtops{\Topgroup{}}$ that contains $\IG{}$.
It is a well-known fact of point-set topology that every maximally-connected set of a topological-space is a closed set of it.
Therefore, $\asubgroup{}$ is a closed set of $\topgtops{\Topgroup{}}$.\\
According to \refthm{thmtopologicalgroupequiv0}, $\gopr{\gop{}}{\TopgG{\Topgroup{}}}$ is a continuous map from
$\topprod{\topgtops{\Topgroup{}}}{\topgtops{\Topgroup{}}}$ to $\topgtops{\Topgroup{}}$, and since a continuous map
maps any connected set of the source topological-space to a connected set of the target topological-space,
it becomes evident that $\func{\image{\gopr{\gop{}}{\TopgG{\Topgroup{}}}}}{\Cprod{\asubgroup{}}{\asubgroup{}}}$
is a connected set of $\topgtops{\Topgroup{}}$ that contains $\IG{}$. Therefore, considering that $\asubgroup{}$ is a maxmally-connected
set of $\topgtops{\Topgroup{}}$ that contains $\IG{}$, it is clear that
$\func{\image{\gopr{\gop{}}{\TopgG{\Topgroup{}}}}}{\Cprod{\asubgroup{}}{\asubgroup{}}}$ must be included in $\asubgroup{}$,
which means $\asubgroup{}$ is a subgroup of $\TopgG{\Topgroup{}}$.\\
Moreover, according to \refthm{thmtranslationsarehomeomorphisms}, for every $\g{}$ in $\G{}$, the
$\g{}$-conjugation of $\TopgG{\Topgroup{}}$ is a homeomorphism from $\topgtops{\Topgroup{}}$ to itself.
Hence, considering that a homeomorphism maps a maximally-connected set of the source topological-space
onto a maximally-connected set of the target topological-space, it is evident that for every $\g{}$ in $\G{}$,
$\func{\image{\[\gconj{\TopgG{\Topgroup{}}}{\g{}}\]}}{\asubgroup{}}$ is a maximally connected set of $\topgtops{\Topgroup{}}$
containing $\IG{}$ and therefore must be equal to $\asubgroup{}$, because there exists only one maximally-connected set of
$\topgtops{\Topgroup{}}$ that contains $\IG{}$. Thus, $\asubgroup{}$ is a normal subgroup of $\TopgG{\Topgroup{}}$.
\endthm
\chapteR{Lie-Algebras}
\thispagestyle{fancy}
\section{Basic structure of a Lie-Algebra}
\fixed
$\F=\triple{\FF}{\vsum{\f}}{\spro{\ff}}$ is fixed as a field of characteristic $0$.
\endfixed
\definition\label{defliealgebra}
$\VS=\tuple{\VV{}}{+}{\times}{\F}$ is taken to be a vector-space (over the field $\F$),
and the neutral element of addition operation of the vector-space $\VS$ is denoted by $\zerovec{}$.
$\aliebra{}$ is taken as an
element of $\Func{\Cprod{\VV{}}{\VV{}}}{\VV{}}$ (a binary operation on $\VV{}$).
The pair $\opair{\VS}{\aliebra{}}$ is called a
$\quotl$Lie-algebra (over the field $\F$)$\quotr$ iff these axioms are satisfied in respect of $\aliebra{}$.
\begin{itemize}
\item[\myitem{LA~1.}]
Bilinearity:
\begin{align}
\Foreach{\triple{\vv{1}}{\vv{2}}{\vv{3}}}{\VV{}^{\times 3}}
\Foreach{\c}{\F}
\begin{cases}
\func{\aliebra{}}{\binary{\c\vv{1}+\vv{2}}{\vv{3}}}=\c\func{\aliebra{}}{\binary{\vv{1}}{\vv{3}}}+
\func{\aliebra{}}{\binary{\vv{2}}{\vv{3}}},\cr
\func{\aliebra{}}{\binary{\vv{3}}{\c\vv{1}+\vv{2}}}=\c\func{\aliebra{}}{\binary{\vv{3}}{\vv{1}}}+
\func{\aliebra{}}{\binary{\vv{3}}{\vv{2}}}.
\end{cases}
\end{align}
\item[\myitem{LA~2.}]
Alternativity:
\begin{equation}
\Foreach{\vv{}}{\VV{}}\func{\aliebra{}}{\binary{\vv{}}{\vv{}}}=\zerovec.
\end{equation}
\item[\myitem{LA~3.}]
Jacobi identity:
\begin{equation}
\Foreach{\triple{\vv{1}}{\vv{2}}{\vv{3}}}{\VV{}^{\times 3}}
\func{\aliebra{}}{\binary{\vv{1}}{\func{\aliebra{}}{\binary{\vv{2}}{\vv{3}}}}}+
\func{\aliebra{}}{\binary{\vv{2}}{\func{\aliebra{}}{\binary{\vv{3}}{\vv{1}}}}}+
\func{\aliebra{}}{\binary{\vv{3}}{\func{\aliebra{}}{\binary{\vv{1}}{\vv{2}}}}}=\zerovec{}.
\end{equation}
\end{itemize}
When $\opair{\VVS{}}{\aliebra{}}$ is a Lie-algebra,
$\VVS{}$ is called the $\quotl$underlying vector-space of the Lie-algebra $\opair{\VVS{}}{\aliebra{}}$, and
the dimension of $\opair{\VVS{}}{\aliebra{}}$
is defined to be the same as that of the underlying vector-space $\VVS{}$, and is denoted by
$\liedim{\opair{\VVS{}}{\aliebra{}}}$. Also any basis (resp. ordered-basis) of the vector-space $\VVS{}$ is also considered to be
a basis (resp. ordered-basis) of the Lie-algebra $\opair{\VVS{}}{\aliebra{}}$, and $\oliebasis{\opair{\VVS{}}{\aliebra{}}})$
is defined to be identical to $\ovecbasis{\VVS{}}$.
Furthermore, the binary operation $\aliebra{}$ is called the $\quotl$Lie-operation of the Lie-algebra
$\opair{\VVS{}}{\aliebra{}}$$\quotr$.\\
Furthermore, when $\opair{\VVS{}}{\aliebra{}}$ is a Lie-algebra such that,
\begin{equation}
\Foreach{\opair{\vv{1}}{\vv{2}}}{\VV{}^{\times 2}}
\func{\aliebra{}}{\binary{\vv{1}}{\vv{2}}}=\zerovec{\VS},
\end{equation}
then $\opair{\VVS{}}{\aliebra{}}$ is called an $\quotl$abelian Lie-algebra$\quotr$.
\endef
\fixed
\begin{itemize}
\item[$\centerdot$]
$\VVS{}=\tuple{\VV{}}{+}{\times}{\F}$ is taken to be a vector-space, and $\aliebra{}$
an element of $\Func{\Cprod{\VV{}}{\VV{}}}{\VV{}}$ (a binary operation on $\VV{}$), such that
$\aliealgebra{}:=\opair{\VVS{}}{\aliebra{}}$ is a Lie-algebra.
\item[$\centerdot$]
$\defSet{\VVS{i}=\tuple{\VV{i}}{\vsum{i}}{\vsprod{i}}{\F}}{i\in\Zp}$ is taken to be a collection of vector-spaces, and $\aliealg{i}$
an element of $\Func{\Cprod{\VV{i}}{\VV{i}}}{\VV{i}}$ (a binary operation on $\VV{i}$), such that
$\aliealgebra{i}:=\opair{\VVS{i}}{\aliealg{i}}$ is a Lie-algebra, for every $i$ in $\Zp$.
\end{itemize}
\endfixed
\definition\label{defliesubalgebra}
$\W{}$ is taken as a (non-empty) subset of $\VV{}$.
\begin{itemize}
\item[$\centerdot$]
$\W{}$ is called a $\quotl$Lie-subalgebra of $\aliealgebra{}$$\quotr$
iff $\W{}$ is both a vector-subspace of $\VVS{}$ and closed under the binary operation $\aliebra{}$, that is,
\begin{align}
\Foreach{\opair{\ww{1}}{\ww{2}}}{\Cprod{\W{}}{\W{}}}
\Foreach{\c}{\f}
\begin{cases}
\c\ww{1}+\ww{2}\in\W{},\cr
\func{\aliebra{}}{\binary{\ww{1}}{\ww{2}}}\in\W{}.
\end{cases}
\end{align}
\item[$\centerdot$]
When $\W{}$ is a Lie-subalgebra of $\aliealgebra{}$,
\begin{equation}
\subspace{\aliealgebra{}}{\W{}}:=
\opair{\tuple{\W{}}{\func{\res{+}}{\Cprod{\W{}}{\W{}}}}{\func{\res{\times}}{\Cprod{\W{}}{\W{}}}}{\F}}
{\func{\res{\aliebra{}}}{\Cprod{\W{}}{\W{}}}}.
\end{equation}
\item[$\centerdot$]
The set of all Lie-subalgebras of $\aliealgebra{}$ is denoted by $\sublie{\aliealgebra{}}$.
\item[$\centerdot$]
Given a non-negative integer $m$,
the set of all $m$-dimensional Lie-subalgebras of $\aliealgebra{}$ is denoted by $\subliedim{m}{\aliealgebra{}}$.
\end{itemize}
\endef
\corollary\label{corliesubalgebra}
For every $\W{}$ in $\sublie{\aliealgebra{}}$ (every Lie-subalgebra of $\aliealgebra{}$),
$\subspace{\aliealgebra{}}{\W{}}$ is a Lie-algebra over the field $\F$.
\endcor
\definition\label{defliealgebramorphism}
\begin{itemize}
\item[$\centerdot$]
$\LiealgMor{\aliealgebra{}}{\aliealgebra{1}}$ is defined as the set of all mappings from $\VV{}$ to $\VV{1}$
that preserve the Lie-algebra structure,
that is the set of all linear maps from the underlying vector-space of
the Lie-algebra $\aliealgebra{}$ to the underlying vector-space of $\aliealgebra{1}$ that
preserve the Lie-operation as well. Precisely,
\begin{align}
&\hskip 0.8\baselineskip\LiealgMor{\aliealgebra{}}{\aliealgebra{1}}\cr
&:=\defset{\aliealgmor{}}{\Lin{\VVS{}}{\VVS{1}}}{\Foreach{\opair{\vv{1}}{\vv{2}}}{\Cprod{\VV{}}{\VV{}}}
\func{\aliealgmor{}}{\func{\aliebra{}}{\binary{\vv{1}}{\vv{2}}}}=
\func{\aliebra{1}}{\binary{\func{\aliealgmor{}}{\vv{1}}}{\func{\aliealgmor{}}{\vv{2}}}}}.\cr
&{}
\end{align}
Each element of $\LiealgMor{\aliealgebra{}}{\aliealgebra{1}}$ is referred to as a $\quotl$Lie-algebra-morphism from
the Lie-algebra $\aliealgebra{}$ to the Lie-algebra $\aliealgebra{1}$$\quotr$.
\item[$\centerdot$]
$\LiealgMon{\aliealgebra{}}{\aliealgebra{1}}$ is defined as the set of all injective maps $\function{\aliealgmor{}}{\VV{}}{\VV{1}}$
such that $\aliealgmor{}$ is a Lie-algebra-morphism from $\aliealgebra{}$ to $\aliealgebra{1}$. That is,
\begin{align}
\LiealgMon{\aliealgebra{}}{\aliealgebra{1}}:=\InF{\VV{}}{\VV{1}}\cap\LiealgMor{\aliealgebra{}}{\aliealgebra{1}}.
\end{align}
Each element of $\LiealgMon{\aliealgebra{}}{\aliealgebra{1}}$ is referred to as a $\quotl$Lie-algebra-monomorphism from
the Lie-algebra $\aliealgebra{}$ to the Lie-algebra $\aliealgebra{1}$$\quotr$.
\item[$\centerdot$]
$\LiealgEpi{\aliealgebra{}}{\aliealgebra{1}}$ is defined as the set of all surjective maps $\function{\aliealgmor{}}{\VV{}}{\VV{1}}$
such that $\aliealgmor{}$ is a Lie-algebra-morphism from $\aliealgebra{}$ to $\aliealgebra{1}$. That is,
\begin{align}
\LiealgEpi{\aliealgebra{}}{\aliealgebra{1}}:=\surFunc{\VV{}}{\VV{1}}\cap\LiealgMor{\aliealgebra{}}{\aliealgebra{1}}.
\end{align}
Each element of $\LiealgEpi{\aliealgebra{}}{\aliealgebra{1}}$ is referred to as a $\quotl$Lie-algebra-epimorphism from
the Lie-algebra $\aliealgebra{}$ to the Lie-algebra $\aliealgebra{1}$$\quotr$.
\item[$\centerdot$]
$\LiealgIsom{\aliealgebra{}}{\aliealgebra{1}}$ is defined as the set of all bijective maps $\function{\aliealgmor{}}{\VV{}}{\VV{1}}$
such that each $\aliealgmor{}$ and $\finv{\aliealgmor{}}$ is a Lie-algebra-morphism from its source Lie-algebra to its
target Lie-algebra. That is,
\begin{align}
\LiealgIsom{\aliealgebra{}}{\aliealgebra{1}}:=\defset{\aliealgmor{}}{\IF{\VV{}}{\VV{1}}}
{\aliealgmor{}\in\LiealgMor{\aliealgebra{}}{\aliealgebra{1}},~
\finv{\aliealgmor{}}\in\LiealgMor{\aliealgebra{1}}{\aliealgebra{}}},
\end{align}
or equivalently,
\begin{equation}
\LiealgIsom{\aliealgebra{}}{\aliealgebra{1}}=
\LiealgMon{\aliealgebra{}}{\aliealgebra{1}}\cap
\LiealgEpi{\aliealgebra{}}{\aliealgebra{1}}.
\end{equation}
Each element of $\LiealgIsom{\aliealgebra{}}{\aliealgebra{1}}$ is referred to as a $\quotl$Lie-algebra-isomorphism from
the Lie-algebra $\aliealgebra{}$ to the Lie-algebra $\aliealgebra{1}$$\quotr$.
\item[$\centerdot$]
By definition,
\begin{equation}
\liealgisomorphic{\aliealgebra{}}{\aliealgebra{1}}
:\thenn
\LiealgIsom{\aliealgebra{}}{\aliealgebra{1}}\neq\empty.
\end{equation}
It is said that $\quotl$$\aliealgebra{}$ is Lie-algebraically-isomorphic (or simply isomorphic) to $\aliealgebra{1}$$\quotr$
iff $\liealgisomorphic{\aliealgebra{}}{\aliealgebra{1}}$.
\end{itemize}
\endef
\definition\label{defliealgebramorphismkernel}
$\aliealgmor{}$ is taken as an element of $\LiealgMor{\aliealgebra{}}{\aliealgebra{1}}$.
The set of all elements of $\VV{}$ mapped by $\aliealgmor{}$ to the neutral element of the vector-space $\VVS{1}$
is referred to as the $\quotl$kernel of the Lie-algebra-morphism $\aliealgmor{}$ from $\aliealgebra{}$ to $\aliealgebra{1}$$\quotr$, and
denoted by $\func{\Lker{\aliealgebra{}}{\aliealgebra{1}}}{\aliealgmor{}}$.
\begin{equation}
\func{\Lker{\aliealgebra{}}{\aliealgebra{1}}}{\aliealgmor{}}:=
\defset{\vv{}}{\VV{}}{\func{\aliealgmor{}}{\vv{}}=\zerovec{\VVS{1}}}.
\end{equation}
$\func{\Lker{\aliealgebra{}}{\aliealgebra{1}}}{\aliealgmor{}}$ is actually the same as the kernel of $\aliealgmor{}$
when regarded as a linear map from the vector-space $\VVS{}$ to the vector-space $\VVS{1}$.
\endef
\theorem\label{thmliealgebramorphismimageandkernel}
$\aliealgmor{}$ is taken as an element of $\LiealgMor{\aliealgebra{}}{\aliealgebra{1}}$.
The kernel of $\aliealgmor{}$ is a Lie-subalgebra of $\aliealgebra{}$. That is,
\begin{equation}
\func{\Lker{\aliealgebra{}}{\aliealgebra{1}}}{\aliealgmor{}}\in\sublie{\aliealgebra{}}.
\end{equation}
The image of $\aliealgmor{}$ is a Lie-subalgebra of $\aliealgebra{1}$. That is,
\begin{equation}
\func{\image{\aliealgmor{}}}{\VV{}}\in\sublie{\aliealgebra{1}}.
\end{equation}
\proof
It is trivial.
\endthm
\theorem
$\aliealgmor{}$ is taken as an element of $\LiealgMor{\aliealgebra{}}{\aliealgebra{1}}$
\begin{itemize}
\item[$\centerdot$]
For every Lie-subalgebra $\WW{}$ of $\aliealgebra{}$,
the domain-restriction of $\aliealgmor{}$ to $\WW{}$ is a Lie-algebra-morphism from
$\subspace{\aliealgebra{}}{\WW{}}$ to $\aliealgebra{1}$.
\begin{equation}
\Foreach{\WW{}}{\sublie{\aliealgebra{}}}
\func{\resd{\aliealgmor{}}}{\WW{}}\in\LiealgMor{\subspace{\aliealgebra{}}{\WW{}}}{\aliealgebra{1}}.
\end{equation}
\item[$\centerdot$]
For every Lie-subalgebra $\WW{1}$ of $\aliealgebra{1}$ including the image of $\aliealgmor{}$,
the codomain-restriction of $\aliealgmor{}$ to $\WW{1}$ is a Lie-algebra-morphism from
$\aliealgebra{}$ to $\subspace{\aliealgebra{1}}{\WW{1}}$.
\begin{equation}
\Foreach{\WW{}}{\defset{\p{\WW{}}}{\sublie{\aliealgebra{1}}}{\p{\WW{}}\supseteq\funcimage{\aliealgmor{}}}}
\func{\rescd{\aliealgmor{}}}{\WW{}}\in\LiealgMor{\aliealgebra{}}{\subspace{\aliealgebra{1}}{\WW{1}}}.
\end{equation}
\end{itemize}
\proof
It is trivial.
\endthm
\corollary\label{corimageofliealgebramonomorphism}
For every Lie-algebra-monomorphism $\aliealgmor{}$ from $\aliealgebra{}$ to $\aliealgebra{1}$, the codomain-restriction
of $\aliealgmor{}$ to its image is a Lie-algebra-isomorphism from $\aliealgebra{}$ to
$\subspace{\aliealgebra{1}}{\funcimage{\aliealgmor{}}}$. That is,
\begin{equation}
\Foreach{\aliealgmor{}}{\LiealgMon{\aliealgebra{}}{\aliealgebra{1}}}
\func{\rescd{\aliealgmor{}}}{\funcimage{\aliealgmor{}}}\in
\LiealgIsom{\aliealgebra{}}{\subspace{\aliealgebra{1}}{\funcimage{\aliealgmor{}}}}.
\end{equation}
\endcor
\chapteR{
Lie Algebra of Smooth Vector-Fields}
\thispagestyle{fancy}
\fixed
\begin{itemize}
\item[$\centerdot$]
$\Man{}=\opair{\M{}}{\maxatlas{}}$ is fixed as an $n$-dimensional and $\difclass{\infty}$ manifold
modeled on the Banach-space $\R^n$.
\item[$\centerdot$]
$\Man{1}=\opair{\M{1}}{\maxatlas{1}}$ is fixed as an $n_1$-dimensional and $\difclass{\infty}$ manifold
modeled on the Banach-space $\R^{n_1}$.
\end{itemize}
\endfixed
\section{Lie-Derivative}
\lemma
For every smooth vector-field $\avecf{}$ on $\Man{}$, and every smooth real-valued map $\cf$ on $\Man{}$,
$\cmp{\Rder{\cf}{\Man{}}}{\avecf{}}$ is a smooth real-valued map on $\Man{}$. That is,
\begin{equation}
\Foreach{\avecf{}}{\vecf{\Man{}}{\infty}}\Foreach{\cf}{\mapdifclass{\infty}{\Man{}}{\RR}}
\cmp{\(\Rder{\cf}{\Man{}}\)}{\avecf{}}\in\mapdifclass{\infty}{\Man{}}{\RR}.
\end{equation}
\proof
$\avecf{}$ is taken as an arbitrary element of $\vecf{\Man{}}{\infty}$, and $\cf$
as an arbitrary element of $\mapdifclass{\infty}{\Man{}}{\RR}$. Then,
\begin{align}
\avecf{}&\in\mapdifclass{\infty}{\Man{}}{\Tanbun{\Man{}}},\\
\Rder{\cf}{\Man{}}&\in\mapdifclass{\infty}{\Tanbun{\Man{}}}{\RR},
\end{align}
and thus,
\begin{equation}
\cmp{\(\Rder{\cf}{\Man{}}\)}{\avecf{}}\in\mapdifclass{\infty}{\Man{}}{\RR}.
\end{equation}
\endlem
\definition\label{defliederivative}
The mapping $\LieDer{\Man{}}$ is defined as,
\begin{align}
&\LieDer{\Man{}}\indef\Func{\vecf{\Man{}}{\infty}}{\Func{\mapdifclass{\infty}{\Man{}}{\RR}}{\mapdifclass{\infty}{\Man{}}{\RR}}},\cr
&\Foreach{\avecf{}}{\vecf{\Man{}}{\infty}}
\Foreach{\cf}{\mapdifclass{\infty}{\Man{}}{\RR}}
\func{\[\func{\LieDer{\Man{}}}{\avecf{}}\]}{\cf}=\cmp{\(\Rder{\cf}{\Man{}}\)}{\avecf{}}.
\end{align}
Equivalenty,
\begin{equation}
\Foreach{\avecf{}}{\vecf{\Man{}}{\infty}}
\Foreach{\cf}{\mapdifclass{\infty}{\Man{}}{\RR}}\Foreach{\point}{\M{}}
\func{\(\func{\[\func{\LieDer{\Man{}}}{\avecf{}}\]}{\cf}\)}{\point}\eqdef
\func{\Rder{\cf}{\Man{}}}{\func{\avecf{}}{\point}}.
\end{equation}
For every $\avecf{}$ in $\vecf{\Man{}}{\infty}$, every $\cf$ in $\mapdifclass{\infty}{\Man{}}{\RR}$,
and every $\point$ in $\M{}$, $\func{\(\func{\[\func{\LieDer{\Man{}}}{\avecf{}}\]}{\cf}\)}{\point}$
is actually defined to be the differential of $\cf$ at the point $\point$ in the direction
$\func{\avecf{}}{\point}$.
\endef
\lemma\label{lemmaliederivativeofavectorfieldisaderivation}
For every smooth vector-field $\avecf{}$ on $\Man{}$, $\func{\LieDer{\Man{}}}{\avecf{}}$ is an $\infty$-derivation
on $\Man{}$. That is,
\begin{equation}
\Foreach{\avecf{}}{\vecf{\Man{}}{\infty}}
\func{\LieDer{\Man{}}}{\avecf{}}\in\Derivation{\Man{}}{\infty},
\end{equation}
or equivalently,
\begin{equation}
\func{\image{\[\LieDer{\Man{}}\]}}{\vecf{\Man{}}{\infty}}\subseteq\Derivation{\Man{}}{\infty}.
\end{equation}
\proof
$\avecf{}$ is taken as an arbitrary element of $\vecf{\Man{}}{\infty}$. So,
\begin{equation}\label{lemmaliederivativeofavectorfieldisaderivationpeq1}
\cmp{\basep{\Man{}}}{\avecf{}}=\identity{\M{}}.
\end{equation}
\begin{itemize}
\item[$\pr{1}$]
Each $\cf$ and $\cg$ is taken as an arbitrary element of $\mapdifclass{\infty}{\Man{}}{\RR}$,
and $\c$ as an arbitrary element of $\R$.
Since $\Rderop{\Man{}}\in\Lin{\mapdifclass{\infty}{\Man{}}{\RR}}{\mapdifclass{\infty}{\Tanbun{\Man{}}}{\RR}}$,
\begin{equation}
\func{\Rderop{\Man{}}}{\c\cf+\cg}=\c\(\Rder{\cf}{\Man{}}\)+\Rder{\cg}{\Man{}},
\end{equation}
and thus according to \refdef{defliederivative}, and invoking the operations of addition
and scalar-multiplication in the canonical linear structures of $\mapdifclass{\infty}{\Man{}}{\RR}$
and $\mapdifclass{\infty}{\Tanbun{\Man{}}}{\RR}$,
\begin{align}
\Foreach{\point}{\M{}}
\func{\(\func{\[\func{\LieDer{\Man{}}}{\avecf{}}\]}{\c\cf+\cg}\)}{\point}&=
\func{\[\Rder{\(\c\cf+\cg\)}{\Man{}}\]}{\func{\avecf{}}{\point}}\cr
&=\func{\[\c\(\Rder{\cf}{\Man{}}\)+\Rder{\cg}{\Man{}}\]}{\func{\avecf{}}{\point}}\cr
&=\c\[\func{\Rder{\cf}{\Man{}}}{\func{\avecf{}}{\point}}\]+\func{\Rder{\cg}{\Man{}}}{\func{\avecf{}}{\point}}\cr
&=\c\[\func{\(\cmp{\Rder{\cf}{\Man{}}}{\avecf{}}\)}{\point}\]+\func{\(\cmp{\Rder{\cg}{\Man{}}}{\avecf{}}\)}{\point}\cr
&=\func{\(\c\[\cmp{\Rder{\cf}{\Man{}}}{\avecf{}}\]+\cmp{\Rder{\cg}{\Man{}}}{\avecf{}}\)}{\point},
\end{align}
and hence,
\begin{equation}
\func{\[\func{\LieDer{\Man{}}}{\avecf{}}\]}{\c\cf+\cg}=
\c\[\cmp{\Rder{\cf}{\Man{}}}{\avecf{}}\]+\cmp{\Rder{\cg}{\Man{}}}{\avecf{}}.
\end{equation}
\endp
\end{itemize}
Therefore,
\begin{align}\label{lemmaliederivativeofavectorfieldisaderivationpeq2}
&\Foreach{\opair{\cf}{\cg}}{\Cprod{\mapdifclass{\infty}{\Man{}}{\RR}}{\mapdifclass{\infty}{\Man{}}{\RR}}}
\Foreach{\c}{\R}\cr
&\func{\[\func{\LieDer{\Man{}}}{\avecf{}}\]}{\c\cf+\cg}=
\c\[\cmp{\Rder{\cf}{\Man{}}}{\avecf{}}\]+\cmp{\Rder{\cg}{\Man{}}}{\avecf{}},
\end{align}
which means $\func{\LieDer{\Man{}}}{\avecf{}}$ is a linear map from $\func{\LieDer{\Man{}}}{\avecf{}}$
to $\func{\LieDer{\Man{}}}{\avecf{}}$. That is,
\begin{equation}
\func{\LieDer{\Man{}}}{\avecf{}}\in
\Lin{\mapdifclass{\infty}{\Man{}}{\RR}}{\mapdifclass{\infty}{\Man{}}{\RR}}.
\end{equation}
\begin{itemize}
\item[$\pr{2}$]
Each $\cf$ and $\cg$ is taken as an arbitrary element of $\mapdifclass{\infty}{\Man{}}{\RR}$.
Considering that,
\begin{equation}
\Foreach{\avec{}}{\tanbun{\Man{}}}
\func{\[\func{\Rderop{\Man{}}}{\cf\rdot\cg}\]}{\avec{}}=
\[\func{\(\cmp{\cg}{\basep{\Man{}}}\)}{\avec{}}\]\[\func{\Rder{\cf}{\Man{}}}{\avec{}}\]+
\[\func{\(\cmp{\cf}{\basep{\Man{}}}\)}{\avec{}}\]\[\func{\Rder{\cg}{\Man{}}}{\avec{}}\],
\end{equation}
and according to \refdef{defliederivative} and \Ref{lemmaliederivativeofavectorfieldisaderivationpeq1},
\begin{align}
\Foreach{\point}{\M{}}
&~~~\func{\(\func{\[\func{\LieDer{\Man{}}}{\avecf{}}\]}{\cf\rdot\cg}\)}{\point}\cr&=
\func{\[\Rder{\(\cf\rdot\cg\)}{\Man{}}\]}{\func{\avecf{}}{\point}}\cr
&=\[\func{\(\cmp{\cg}{\basep{\Man{}}}\)}{\func{\avecf{}}{\point}}\]\[\func{\Rder{\cf}{\Man{}}}{\func{\avecf{}}{\point}}\]+
\[\func{\(\cmp{\cf}{\basep{\Man{}}}\)}{\func{\avecf{}}{\point}}\]\[\func{\Rder{\cg}{\Man{}}}{\func{\avecf{}}{\point}}\]\cr
&=\[\func{\cg}{\point}\]\bigg(\func{\(\func{\[\func{\LieDer{\Man{}}}{\avecf{}}\]}{\cf}\)}{\point}\bigg)+
\[\func{\cf}{\point}\]\bigg(\func{\(\func{\[\func{\LieDer{\Man{}}}{\avecf{}}\]}{\cg}\)}{\point}\bigg)\cr
&=\func{\bigg(\func{\[\func{\LieDer{\Man{}}}{\avecf{}}\]}{\cf}\rdot\cg+
\cf.\func{\[\func{\LieDer{\Man{}}}{\avecf{}}\]}{\cg}\bigg)}{\point},
\end{align}
and hence,
\begin{equation}
\func{\[\func{\LieDer{\Man{}}}{\avecf{}}\]}{\cf\rdot\cg}=
\func{\[\func{\LieDer{\Man{}}}{\avecf{}}\]}{\cf}\rdot\cg+
\cf\rdot\func{\[\func{\LieDer{\Man{}}}{\avecf{}}\]}{\cg}.
\end{equation}
\endp
\end{itemize}
Therefore,
\begin{align}\label{lemmaliederivativeofavectorfieldisaderivationpeq3}
&\Foreach{\opair{\cf}{\cg}}{\Cprod{\mapdifclass{\infty}{\Man{}}{\RR}}{\mapdifclass{\infty}{\Man{}}{\RR}}}\cr
&\func{\[\func{\LieDer{\Man{}}}{\avecf{}}\]}{\cf\rdot\cg}=
\func{\[\func{\LieDer{\Man{}}}{\avecf{}}\]}{\cf}\rdot\cg+
\cf\rdot\func{\[\func{\LieDer{\Man{}}}{\avecf{}}\]}{\cg}.
\end{align}
\Ref{lemmaliederivativeofavectorfieldisaderivationpeq2} and \Ref{lemmaliederivativeofavectorfieldisaderivationpeq3} imply
that $\func{\LieDer{\Man{}}}{\avecf{}}$ is a $\infty$-derivation on $\Man{}$. That is,
\begin{equation}
\func{\LieDer{\Man{}}}{\avecf{}}\in\Derivation{\Man{}}{\infty}.
\end{equation}
\endlem
\definition\label{defliederivative00}
$\Lieder{\Man{}}$ is defined to be the codomain-restriction of $\LieDer{\Man{}}$ to
$\Derivation{\Man{}}{\infty}$.
\begin{equation}
\Lieder{\Man{}}:=\func{\rescd{\LieDer{\Man{}}}}{\Derivation{\Man{}}{\infty}}.
\end{equation}
When there is no ambiguity about the manifold $\Man{}$, $\Lieder{}$ can replace $\Lieder{\Man{}}$.
$\Lieder{\Man{}}$ is referred to as the $\quotl$Lie-derivative operator on $\Man{}$$\quotr$, and for each
smooth vector-field $\avecf{}$ on $\Man{}$, $\func{\Lieder{\Man{}}}{\avecf{}}$ is simply referred to as the
$\quotl$Lie-derivative of $\avecf{}$$\quotr$.
\endef
\lemma\label{lemmacutofffunctions}
$\K$ is taken as a compact set of $\mantops{\Man{}}$, and $\U$ as an open set of $\mantops{\Man{}}$
such that $\K\subseteq\U$. There exists a real-valued smooth map on $\Man{}$ with compact support included in $\U$ and
constant value $1$ on $\K$. That is,
\begin{equation}
\Exists{\cf}{\mapdifclass{\infty}{\Man{}}{\RR}}
\[\compacts{\mantops{\Man{}}}\ni\support{\cf}\subseteq\U,~
\func{\image{\cf}}{\K}=\seta{1}\].
\end{equation}
\proof
A detailed proof is given in \cite[page~28,~Lemma~1.69]{JLee}.
\endthm
\lemma\label{lemexistanceofcompactsetsofdomainsofchartscontainingopensets}
$\point$ is taken as a point of $\Man{}$ and $\phi$ as a chart of $\Man{}$ whose domain contains $\point$.
There exists a subset of $\U$ that is a compact set of $\mantops{\Man{}}$
including an open set of $\mantops{\Man{}}$ which contains $\point$. That is,
\begin{equation}
\Existsis{\K}{\compacts{\mantops{\Man{}}}}
\Existsis{\V}{\mantop{\Man{}}}
\point\in\V\subseteq\K\subseteq\U.
\end{equation}
\proof
Set $\U:=\domain{\phi}$ and $\p{\U}:=\funcimage{\phi}$. $\U$ is a sub-space of the topological-space $\mantops{\Man{}}$ and $\p{\U}$
a sub-space of $\RR^{n}$ that are homeomorphic via the map $\phi$. Thus they possess the same topological properties.
So, since $\p{\U}$ includes a compact set of $\p{\U}$ including an open set of $\p{\U}$ which contains $\func{\phi}{\point}$
(there is no need to emphasize the
ambient space $\p{\U}$, since the set of all compact (or open) sets of $\p{\U}$ coincides with the set of all compact
(or open) sets of $\RR^{n}$), there also exists a compact set of $\U$ including an open set of $\U$ which contains $\point$.
Furthermore, every compact set of $\U$ is a compact set of $\mantops{\Man{}}$, and since $\U$ is an open set of $\mantops{\Man{}}$
every open set of $\U$ is an open set of $\mantops{\Man{}}$.
\endlem
\lemma\label{lemmadifferentiablemapexpansion}
$\point$ is taken as an element of $\M{}$ (a point of the manifold $\Man{}$), and
$\phi$ as an element of $\maxatlas{}$ (a chart of the manifold $\Man{}$) such that $\point\in\domain{\phi}$.
$\lambda$ is taken as an element of $\banachmapdifclass{\infty}{\R^n}{\R}{\funcimage{\phi}}{\R}$
(a smooth map from $\funcimage{\phi}$ to $\R$).\\
There exists a real-valued smooth map on $\Man{}$ that coincides with $\cmp{\lambda}{\phi}$ in a neighbourhood
of $\point$ included in $\domain{\phi}$. That is,
\begin{equation}
\Exists{\cf}{\mapdifclass{\infty}{\Man{}}{\RR}}
\Exists{\V}{\mantop{\Man{}}}
\[\point\in\V\subseteq\domain{\phi},~
\func{\resd{\cf}}{\V}=\func{\resd{\cmp{\lambda}{\phi}}}{\V}\].
\end{equation}
\proof
According to \reflem{lemexistanceofcompactsetsofdomainsofchartscontainingopensets},
\begin{equation}\label{lemmadifferentiablemapexpansionpeq1}
\Existsis{\K}{\compacts{\mantops{\Man{}}}}
\Existsis{\V}{\mantop{\Man{}}}
\point\in\V\subseteq\K\subseteq\U.
\end{equation}
So according to \reflem{lemmacutofffunctions},
\begin{equation}\label{lemmadifferentiablemapexpansionpeq2}
\Existsis{\rho}{\mapdifclass{\infty}{\Man{}}{\RR}}
\[\compacts{\mantops{\Man{}}}\ni\support{\rho}\subseteq\U,~
\func{\image{\rho}}{\K}=\seta{1}\].
\end{equation}
$\cf$ is defined as the element of $\Func{\M{}}{\R}$ such that,
\begin{align}\label{lemmadifferentiablemapexpansionpeq3}
\Foreach{\x}{\M{}}
\func{\cf}{\x}:=
\begin{cases}
\func{\rho}{\x}\[\func{\(\cmp{\lambda}{\phi}\)}{\x}\],~&\x\in\U,\cr
0,~&\x\in\compl{\M{}}{\U}.
\end{cases}
\end{align}
Since $\rho$ is a real-valued smooth map on $\Man{}$,
\begin{equation}\label{lemmadifferentiablemapexpansionpeq4}
\cmp{\rho}{\finv{\phi}}\in\banachmapdifclass{\infty}{\R^n}{\R}{\p{\U}}{\R}.
\end{equation}
Thus since $\lambda\in\banachmapdifclass{\infty}{\R^n}{\R}{\p{\U}}{\R}$, clearly,
\begin{equation}\label{lemmadifferentiablemapexpansionpeq5}
\(\cmp{\rho}{\finv{\phi}}\)\rdot\lambda
\in\banachmapdifclass{\infty}{\R^n}{\R}{\p{\U}}{\R}.
\end{equation}
In addition, it is clear that,
\begin{equation}\label{lemmadifferentiablemapexpansionpeq6}
\cmp{\cf}{\finv{\phi}}=\(\cmp{\rho}{\finv{\phi}}\)\rdot\lambda.
\end{equation}
Therefore,
\begin{align}\label{lemmadifferentiablemapexpansionpeq7}
\cmp{\identity{\R}}{\cmp{\cf}{\finv{\phi}}}=
\cmp{\cf}{\finv{\phi}}
\in\banachmapdifclass{\infty}{\R^n}{\R}{\p{\U}}{\R},
\end{align}
and hence considering that $\identity{\R}$ is a chart of the manifold $\RR$,
and $\func{\image{\cf}}{\U}\subseteq\R=\domain{\identity{\R}}$,
\begin{align}\label{lemmadifferentiablemapexpansionpeq8}
&\Foreach{\x}{\U}\Existsis{\phi}{\maxatlas{}}\cr
&\[\x\in\domain{\phi},~\func{\image{\cf}}{\domain{\phi}}\subseteq\domain{\identity{\R}},~
\cmp{\identity{\R}}{\cmp{\cf}{\finv{\phi}}}\in\banachmapdifclass{\infty}{\R^n}{\R}{\p{\U}}{\R}\].
\end{align}
This means $\cf$ is infinitely defferentiable at every $\x$ in $\U$.\\
It is a trivial fact that,
\begin{equation}\label{lemmadifferentiablemapexpansionpeq9}
\Foreach{\x}{\(\compl{\M{}}{\support{\cf}}\)}
\func{\cf}{\x}=0.
\end{equation}
Also since $\support{\cf}\subseteq\U$,
\begin{equation}\label{lemmadifferentiablemapexpansionpeq10}
\(\compl{\M{}}{\U}\)\subseteq
\(\compl{\M{}}{\support{\cf}}\),
\end{equation}
and since $\support{\cf}$ is a closed set of $\mantops{\Man{}}$, $\(\compl{\M{}}{\support{\cf}}\)$
is an open set of $\mantops{\Man{}}$, that is,
\begin{equation}\label{lemmadifferentiablemapexpansionpeq11}
\(\compl{\M{}}{\support{\cf}}\)\in\mantop{\Man{}}.
\end{equation}
Thus considering that $\defSet{\domain{\psi}}{\psi\in\maxatlas{}}$ is a base for the topological-space
$\mantops{\Man{}}$, for every $\x$ in $\(\compl{\M{}}{\U}\)$ there exists a chart $\psi_{\x}$ of $\Man{}$
whose domain contains $\x$ and is included in $\(\compl{\M{}}{\support{\cf}}\)$, and hence
$\func{\(\cmp{\cf}{\finv{\psi_{\x}}}\)}{\y}=0$ for every $\y$ in $\domain{\psi_{\x}}$. That is,
\begin{align}\label{lemmadifferentiablemapexpansionpeq12}
&\Foreach{\x}{\(\compl{\M{}}{\U}\)}\Existsis{\psi_{\x}}{\maxatlas{}}\cr
&\[\x\in\domain{\psi_{\x}},~\func{\image{\cf}}{\domain{\psi_{\x}}}\subseteq\domain{\identity{\R}},~
\cmp{\identity{\R}}{\cmp{\cf}{\finv{\psi_{\x}}}}\in\banachmapdifclass{\infty}{\R^n}{\R}{\p{\U}}{\R}\],
\end{align}
which means $\cf$ is infinitely defferentiable at every $\x$ in $\compl{\M{}}{\U}$.\\
\Ref{lemmadifferentiablemapexpansionpeq8} and \Ref{lemmadifferentiablemapexpansionpeq12} imply that
$\cf$ is infinitely differentiable at every point of $\Man{}$, and hence a real-valued amooth map on $\Man{}$. That is,
\begin{equation}
\cf\in\mapdifclass{\infty}{\Man{}}{\RR}.
\end{equation}
In addition, according to \Ref{lemmadifferentiablemapexpansionpeq1}, \Ref{lemmadifferentiablemapexpansionpeq2}, and
\Ref{lemmadifferentiablemapexpansionpeq3}, $\cf$ coincides with $\cmp{\lambda}{\phi}$ in the neighbourhood $\V$ of $\point$. That is,
\begin{equation}
\Existsis{\V}{\mantop{\Man{}}}
\[\point\in\V\subseteq\domain{\phi},~
\func{\resd{\cf}}{\V}=\func{\resd{\cmp{\lambda}{\phi}}}{\V}\].
\end{equation}
\endlem
\lemma\label{lemmanonvanishingderivative}
For every point $\point$ of $\Man{}$ and every vector $\avec{}$ in $\compl{\tanspace{\point}{\Man{}}}{\seta{\zerovec{\point}}}$,
there exists a smooth real-valued map on $\Man{}$ with
non-vanishing derivative in the direction $\avec{}$. That is,
\begin{equation}
\Foreach{\point}{\M{}}
\Foreach{\avec{}}{\compl{\tanspace{\point}{\Man{}}}{\seta{\zerovec{\point}}}}
\Exists{\cf}{\mapdifclass{\infty}{\Man{}}{\RR}}
\func{\Rder{\cf}{\Man{}}}{\avec{}}\neq 0.
\end{equation}
\proof
$\point$ is taken as an arbitrary element of $\M{}$ and
$\avec{}$ as an arbitrary element of $\compl{\tanspace{\point}{\Man{}}}{\seta{\zerovec{\point}}}$. So,
$\point=\func{\basep{\Man{}}}{\avec{}}$. A chart $\phi$ of $\maxatlas{}$ whose domain contains $\point$ is chosen, and
$\U:=\domain{\phi}$. By definition,
\begin{align}\label{lemmanonvanishingderivativepeq1}
\nu:=\func{\bigg(\tanspaceiso{\point}{\Man{}}{\phi}\bigg)}{\avec{}}.
\end{align}
$\nu$ is actually the canonically corresponded $\R^n$-vector to $\avec{}$, where the correspondence takes place in the
tangent-space of $\Man{}$ at the base-point of $\avec{}$. Furtheremore, for each $k$ in $\seta{\suc{1}{n}}$,
$\nu_k$ is defined to be the $k$-th component of $\nu$, that is,
\begin{equation}\label{lemmanonvanishingderivativepeq2}
\Foreach{k}{\seta{\suc{1}{n}}}
\nu_k:=\func{\projection{n}{k}}{\nu}.
\end{equation}
Now, the mapping $\lambda$ is defined as,
\begin{equation}\label{lemmanonvanishingderivativepeq3}
\lambda:=\sum_{k=1}^{n}\nu_{k}\projection{n}{k},
\end{equation}
which is evidently a linear map from $\R^n$ to $\R$, and hence a smooth map from $\R^{n}$ to $\R$.
According to \reflem{lemmadifferentiablemapexpansion}, there exists a real-valued smooth map
on $\Man{}$ that coincides with $\cmp{\lambda}{\phi}$ in a neighbourhood of $\point$. That is,
\begin{align}\label{lemmanonvanishingderivativepeq4}
\Existsis{\cf}{\mapdifclass{\infty}{\Man{}}{\RR}}
\Existsis{\V}{\mantop{\Man{}}}
\[\point\in\V\subseteq\U,~
\func{\resd{\cf}}{\V}=\func{\resd{\cmp{\lambda}{\phi}}}{\V}\].
\end{align}
Since $\V$ is an open subset of $\U$, clearly the domain-restriction of $\phi$ to $\V$ is also an element of $\maxatlas{}$
(a chart of $\Man{}$), which will be denoted by $\varphi$. That is,
\begin{equation}\label{lemmanonvanishingderivativepeq5}
\varphi:=\func{\resd{\phi}}{\V}\in\maxatlas{}.
\end{equation}
Therefore,
\begin{align}\label{lemmanonvanishingderivativepeq6}
\func{\(\Rder{\cf}{\Man{}}\)}{\avec{}}&=\func{\[\cmp{\bigg(\func{\[\banachder{\(\cmp{f}{\finv{\varphi}}\)}{\R^n}{\R}\]}
{\func{\varphi}{\point}}\bigg)}{\tanspaceiso{\point}{\Man{}}{\varphi}}\]}{\avec{}}\cr
&=\func{\[\cmp{\bigg(\func{\[\banachder{\bar{\lambda}}{\R^n}{\R}\]}
{\func{\phi}{\point}}\bigg)}{\tanspaceiso{\point}{\Man{}}{\varphi}}\]}{\avec{}},
\end{align}
where,
\begin{align}\label{lemmanonvanishingderivativepeq7}
\bar{\lambda}:=&\cmp{\cmp{\lambda}{\varphi}}{\finv{\varphi}}\cr
=&\func{\resd{\lambda}}{\func{\image{\phi}}{\V}}.
\end{align}
So since $\bar{\lambda}$ is the domain-restriction of a linear map to an open set of $\R^n$, clearly,
\begin{equation}\label{lemmanonvanishingderivativepeq8}
\bigg(\func{\[\banachder{\bar{\lambda}}{\R^n}{\R}\]}
{\func{\phi}{\point}}\bigg)=\lambda.
\end{equation}
\Ref{lemmanonvanishingderivativepeq1}, \Ref{lemmanonvanishingderivativepeq2}, \Ref{lemmanonvanishingderivativepeq3},
\Ref{lemmanonvanishingderivativepeq6}, and \Ref{lemmanonvanishingderivativepeq8} imply,
\begin{align}\label{lemmanonvanishingderivativepeq9}
\func{\(\Rder{\cf}{\Man{}}\)}{\avec{}}&=\func{\lambda}{\nu}\cr
&=\sum_{k=1}^{n}{\nu_{k}}\func{\projection{n}{k}}{\nu}\cr
&=\sum_{k=1}^{n}\nu_{k}^{2}.
\end{align}
Since $\tanspaceiso{\point}{\Man{}}{\phi}$ is a linear isomorphism from $\Tanspace{\point}{\Man{}}$ to $\R^n$,
it sends non-zero elements of $\tanspace{\point}{\Man{}}$ to non-zero elements of $\R^n$. Thus $\nu\neq\zerovec{}$
since $\avec{}\neq\zerovec{\point}$, and hence,
\begin{equation}\label{lemmanonvanishingderivativepeq10}
\sum_{k=1}^{n}\nu_{k}^{2}\neq 0,
\end{equation}
which together with \Ref{lemmanonvanishingderivativepeq9} implies,
\begin{equation}
\func{\(\Rder{\cf}{\Man{}}\)}{\avec{}}\neq 0.
\end{equation}
\endlem
\lemma\label{lemmaliederivativeislinear}
The Lie-derivative operator on $\Man{}$ is a linear map from the vector-space of smooth vector-fields on $\Man{}$
to the vector-space of $\infty$-derivations on $\Man{}$. That is,
\begin{equation}
\Lieder{\Man{}}\in\Lin{\Vecf{\Man{}}{\infty}}{\LDerivation{\Man{}}{\infty}}.
\end{equation}
\proof
Each $\avecf{1}$ and $\avecf{2}$ is taken as an arbitrary  element of $\vecf{\Man{}}{\infty}$, and $\c$ as an arbitrary element of $\R$.
According to \refdef{defliederivative}, and considering that every smooth vector-field sends every point of $\Man{}$
to a vector in the tangent-space of $\Man{}$ at that point,
and invoking the operations of addition
and scalar-multiplication of the vector-spaces $\Vecf{\Man{}}{\infty}$, $\LDerivation{\Man{}}{\infty}$, and
$\Tanspace{\point}{\Man{}}$ for each $\point$ in $\M{}$, and also considering that for every real-valued smooth map $\cf$ on $\Man{}$,
$\Rder{\cf}{\Man{}}$ operates linearly on the tangent-space of $\Man{}$ at each point of $\Man{}$,
\begin{align}
&\begin{aligned}
\Foreach{\cf}{\mapdifclass{\infty}{\Man{}}{\RR}}\Foreach{\point}{\M{}}
\end{aligned}\cr
&\begin{aligned}
\func{\(\func{\[\func{\Lieder{\Man{}}}{\c\avecf{1}+\avecf{2}}\]}{\cf}\)}{\point}&=
\func{\Rder{\cf}{\Man{}}}{\func{\[\c\avecf{1}+\avecf{2}\]}{\point}}\cr
&=\func{\Rder{\cf}{\Man{}}}{\c\func{\avecf{1}}{\point}+\func{\avecf{2}}{\point}}\cr
&=\c\[\func{\Rder{\cf}{\Man{}}}{\func{\avecf{1}}{\point}}\]+\func{\Rder{\cf}{\Man{}}}{\func{\avecf{2}}{\point}}\cr
&=\c\func{\(\func{\[\func{\Lieder{\Man{}}}{\avecf{1}}\]}{\cf}\)}{\point}+
\func{\(\func{\[\func{\Lieder{\Man{}}}{\avecf{2}}\]}{\cf}\)}{\point}\cr
&=\func{\(\c\func{\[\func{\Lieder{\Man{}}}{\avecf{1}}\]}{\cf}+\func{\[\func{\Lieder{\Man{}}}{\avecf{2}}\]}{\cf}\)}{\point}\cr
&=\func{\(\func{\[\c\func{\Lieder{\Man{}}}{\avecf{1}}+\func{\Lieder{\Man{}}}{\avecf{2}}\]}{\cf}\)}{\point}.
\end{aligned}\cr
&{}
\end{align}
Therefore,
\begin{equation}
\func{\Lieder{\Man{}}}{\c\avecf{1}+\avecf{2}}=
\c\func{\Lieder{\Man{}}}{\avecf{1}}+\func{\Lieder{\Man{}}}{\avecf{2}}.
\end{equation}
\endlem
\theorem\label{thmliederivativeisinjective}
$\Lieder{\Man{}}$ is an injective map. That is,
\begin{equation}
\Lieder{\Man{}}\in\InF{\vecf{\Man{}}{\infty}}{\Derivation{\Man{}}{\infty}}.
\end{equation}
\proof
Each $\avecf{1}$ and $\avecf{2}$ is taken as an element of $\vecf{\Man{}}{\infty}$ such that,
\begin{equation}
\func{\Lieder{\Man{}}}{\avecf{1}}=\func{\Lieder{\Man{}}}{\avecf{2}}.
\end{equation}
According to the linearity of $\Lieder{\Man{}}$,
\begin{equation}
\func{\Lieder{\Man{}}}{\avecf{1}-\avecf{2}}=\zerovec{},
\end{equation}
which means,
\begin{equation}
\Foreach{\cf}{\mapdifclass{\infty}{\Man{}}{\RR}}\Foreach{\point}{\M{}}
\func{\(\func{\[\func{\Lieder{\Man{}}}{\avecf{1}-\avecf{2}}\]}{\cf}\)}{\point}=0.
\end{equation}
Thus according to \refdef{defliederivative},
\begin{equation}
\Foreach{\cf}{\mapdifclass{\infty}{\Man{}}{\RR}}\Foreach{\point}{\M{}}
\func{\Rder{\cf}{\Man{}}}{\func{\[\avecf{1}-\avecf{2}\]}{\point}}=0,
\end{equation}
and hence according to \reflem{lemmanonvanishingderivative},
\begin{equation}
\Foreach{\point}{\M{}}
\func{\(\avecf{1}-\avecf{2}\)}{\point}=0,
\end{equation}
and hence,
\begin{equation}
\avecf{1}=\avecf{2}.
\end{equation}
\endthm
\lemma\label{lemmaderivationislocal}
$\point$ is taken as an element of $\Man{}$, and each $\cf$ and $\cg$ is taken as an element of $\mapdifclass{\infty}{\Man{}}{\RR}$.
If there exists a neighbourhood $O$ of $\point$ (in the topological-space $\mantops{\Man{}}$) on which $\cf$ and $\cg$ coincide,
then for every $\infty$-derivation $\aderivation{}$ on $\Man{}$, $\func{\aderivation{}}{\cf}$ and $\func{\aderivation{}}{\cg}$ coincide at
$\point$. That is,
\begin{align}
&~~~~~\bigg(\Exists{O}{\mantop{\Man{}}}\[\point\in\U,~\Foreach{\x}{O}\func{\cf}{\x}=\func{\cg}{\x}\]\bigg)\cr
&\then
\bigg(\Foreach{\aderivation{}}{\Derivation{\Man{}}{\infty}}
\func{\[\func{\aderivation{}}{\cf}\]}{\point}=\func{\[\func{\aderivation{}}{\cg}\]}{\point}\bigg).
\end{align}
\proof
$\aderivation{}$ is taken as an arbitrary element of $\Derivation{\Man{}}{\infty}$.
It is assumed that there exists an open set $O$ of $\mantops{\Man{}}$ containing $\point$
(a neighbourhood of $\point$ in $\mantops{\Man{}}$) that,
\begin{equation}\label{lemmaderivationislocalpeq1}
\Foreach{\x}{O}\func{\cf}{\x}=\func{\cg}{\x}.
\end{equation}
Since $\defSet{\domain{\phi}}{\phi\in\maxatlas{}}$ is a base for the topology $\mantop{\Man{}}$,
there exists a chart $\phi$ of $\Man{}$ such that,
\begin{equation}\label{lemmaderivationislocalpeq2}
\point\in\domain{\phi}\subseteq O.
\end{equation}
By defining $\U:=\domain{\phi}$ and,
\begin{equation}\label{lemmaderivationislocalpeq3}
\hf:=\cf-\cg,
\end{equation}
it is clear that,
\begin{equation}\label{lemmaderivationislocalpeq4}
\hf\in\mapdifclass{\infty}{\Man{}}{\RR},~\Existsis{\U}{\func{\nei{\mantops{\Man{}}}}{\seta{\point}}}\Foreach{\x}{\U}\func{\hf}{\x}=0.
\end{equation}
According to \reflem{lemexistanceofcompactsetsofdomainsofchartscontainingopensets},
there exists a compact set $\K$ of $\mantops{\Man{}}$ included in $\U$ and containing $\point$. That is,
\begin{equation}\label{lemmaderivationislocalpeq5}
\Existsis{\K}{\compacts{\mantops{\Man{}}}}
\point\in\K\subseteq\U,
\end{equation}
and hence according to \reflem{lemmacutofffunctions}, there exists a real-valued smooth map $\rho$ on $\Man{}$
with a compact support included in $\U$ (so vanishing outside $\U$), possessing the value $1$ all over the set $\K$.
\begin{equation}\label{lemmaderivationislocalpeq6}
\Existsis{\rho}{\mapdifclass{\infty}{\Man{}}{\RR}}
\[\compacts{\mantops{\Man{}}}\ni\support{\rho}\subseteq\U,~
\func{\image{\rho}}{\K}=\seta{1}\].
\end{equation}
It is evident that,
\begin{equation}\label{lemmaderivationislocalpeq7}
\Foreach{\x}{\M{}}
\func{\(\rho\rdot\hf\)}{\x}=0,
\end{equation}
which means
\begin{equation}\label{lemmaderivationislocalpeq8}
\rho\rdot\hf=\zerovec{},
\end{equation}
where $\zerovec{}$ denotes the neutral element of addition operation of the vector-space $\Lmapdifclass{\infty}{\Man{}}{\RR}$.
Thus, since $\aderivation{}$ is a linear map from $\Lmapdifclass{\infty}{\Man{}}{\RR}$ to $\Lmapdifclass{\infty}{\Man{}}{\RR}$,
clearly,
\begin{equation}\label{lemmaderivationislocalpeq9}
\func{\aderivation{}}{\rho\rdot\hf}=\zerovec{}.
\end{equation}
In addition, as a property of $\infty$-derivations on $\Man{}$,
and considering that $\func{\hf}{\point}=0$ and $\func{\rho}{\point}=1$,
\begin{align}\label{lemmaderivationislocalpeq10}
\func{\[\func{\aderivation{}}{\rho\rdot\hf}\]}{\point}&=
\bigg(\func{\[\func{\aderivation{}}{\rho}\]}{\point}\bigg)\bigg(\func{\hf}{\point}\bigg)+
\bigg(\func{\rho}{\point}\bigg)\bigg(\func{\[\func{\aderivation{}}{\hf}\]}{\point}\bigg)\cr
&=\func{\[\func{\aderivation{}}{\hf}\]}{\point}.
\end{align}
\Ref{lemmaderivationislocalpeq9} and \Ref{lemmaderivationislocalpeq10} imply,
\begin{equation}\label{lemmaderivationislocalpeq11}
\func{\[\func{\aderivation{}}{\hf}\]}{\point}=0.
\end{equation}
Using the linearity of $\aderivation{}$, \Ref{lemmaderivationislocalpeq3} implies,
\begin{equation}\label{lemmaderivationislocalpeq12}
\func{\[\func{\aderivation{}}{\hf}\]}{\point}=
\func{\[\func{\aderivation{}}{\cf}\]}{\point}-
\func{\[\func{\aderivation{}}{\cg}\]}{\point},
\end{equation}
and hence according to \Ref{lemmaderivationislocalpeq11},
\begin{equation}
\func{\[\func{\aderivation{}}{\cf}\]}{\point}=
\func{\[\func{\aderivation{}}{\cg}\]}{\point}.
\end{equation}
\endlem
\definition\label{defsmoothextensionoflocalmaps}
$\point$ is taken as a point of $\Man{}$, $\phi$ as an element of $\maxatlas{}$ (a chart of $\Man{}$) such that
$\point\in\domain{\phi}$. $\lambda$ is taken as an element of $\banachmapdifclass{\infty}{\R^n}{\R}{\U}{\R}$ for some
open $\U$ in $\R^{n}$ such that $\funcimage{\phi}\subseteq\U$.
\begin{align}
&\fextension{\Man{}}{\point}{\phi}{\lambda}:=\cr
&\defset{\cf}{\mapdifclass{\infty}{\Man{}}{\RR}}
{\Exists{\V}{\mantop{\Man{}}}
\[\point\in\V\subseteq\domain{\phi},~
\func{\resd{\cf}}{\V}=\func{\resd{\cmp{\lambda}{\phi}}}{\V}\]}.
\end{align}
According to \reflem{lemmadifferentiablemapexpansion}, this set is non-empty.
Each element of $\fextension{\Man{}}{\point}{\phi}{\lambda}$ is referred to as a
$\quotl$smooth extension of $\cmp{\lambda}{\phi}$ on $\Man{}$ fixed at $\point$$\quotr$.
\endef
\lemma\label{lemlocalpowersum}
$\U$ is taken as an open and convex subset of $\R^n$, and $\cf$ as an element of $\banachmapdifclass{\infty}{\R^n}{\R}{\U}{\R}$
(a smooth map from $\U$ to $\R$). 
\begin{align}
&\Foreach{\a}{\U}\Exists{\mtuple{\cg_1}{\cg_n}}{{\banachmapdifclass{\infty}{\R^n}{\R}{\U}{\R}}^{\times n}}\cr
&\Foreach{k}{\seta{\suc{1}{n}}}\func{\cg_k}{\a}=\func{\[\func{\(\banachder{\cf}{\R^n}{\R}\)}{\a}\]}{\Eucbase{n}{k}},\cr
&\Foreach{\x}{\U}\func{\cf}{\x}=\func{\cf}{\a}+\sum_{k=1}^{n}\[\func{\projection{n}{k}}{\x}-\func{\projection{n}{k}}{\a}\]\func{\cg_k}{\x}.
\end{align}
\proof
For each $\x$ in $\U$, $\function{s_\x}{\cinterval{0}{1}}{\R^n}$ is defined as,
\begin{equation}\label{lemlocalpowersumpeq1}
\Foreach{t}{\cinterval{0}{1}}
\func{s_\x}{t}\eqdef\(1-t\)\y+t\x,
\end{equation}
and $\function{\xi_\x}{\cinterval{0}{1}}{\R}$ is defined as,
\begin{equation}\label{lemlocalpowersumpeq2}
\xi_\x:=\cmp{\cf}{s_\x}.
\end{equation}
For every $\x$ in $\U$, $\xi_\x$ is smooth since both $\cf$ and $s_\x$ are smooth.
It should be noticed that $\cinterval{0}{1}$
is not an open set of $\R$, and since objects of differentiablity are mappings defined on an open subset
of a Banach-space, the differentiablity of $s_\x$ merits reconsideration. In fact $s_\x$ is smooth in the sense that
there exists a unique mapping defined on an open subset of $\R$ including $\Cprod{\U}{\cinterval{0}{1}}$ that
is infinitely differentiable and coincides with $s_\x$ in the interval $\cinterval{0}{1}$.
Furthermore according to the chain rule
for differentiation,
\begin{align}\label{lemlocalpowersumpeq3}
\Foreach{\x}{\U}
\Foreach{t}{\cinterval{0}{1}}
\func{\derived{\xi_\x}}{t}&=
\func{\[\func{\(\banachder{\xi_\x}{\R}{\R}\)}{t}\]}{1}\cr
&=\func{\bigg(\cmp{\[\func{\(\banachder{\cf}{\R^n}{\R}\)}{\func{s_\x}{t}}\]}{\[\func{\(\banachder{s_\x}{\R}{\R^n}\)}{t}\]}\bigg)}{1}\cr
&=\func{\[\func{\(\banachder{\cf}{\R^n}{\R}\)}{\func{s_\x}{t}}\]}{\func{\derived{s_\x}}{t}}\cr
&=\func{\[\func{\(\banachder{\cf}{\R^n}{\R}\)}{\func{s_\x}{t}}\]}{\x-\y}.
\end{align}
where $\derived{\xi_\x}$ and $\derived{s_\x}$ are the usual
derivatives of $\xi_\x$ and $s_\x$ as considered in calculus.
For every $t$ in $\cinterval{0}{1}$, since the differential of $\cf$ at $\func{s_\x}{t}$,
that is $\banachder{\cf}{\R^n}{\R}$, is a linear map from $\R^n$ to $\R$,
\begin{align}\label{lemlocalpowersumpeq4}
&\Foreach{\x}{\U}
\Foreach{t}{\cinterval{0}{1}}\cr
&\func{\[\func{\(\banachder{\cf}{\R^n}{\R}\)}{\func{s_\x}{t}}\]}{\x-\y}=
\sum_{k=1}^{n}\[\func{\partialder{k}\cf}{\(1-t\)\y+t\x}\]\[\func{\projection{n}{k}}{\x}-\func{\projection{n}{k}}{\y}\],
\end{align}
where,
\begin{align}\label{lemlocalpowersumpeq5}
\Foreach{k}{\seta{\suc{1}{n}}}\Foreach{z}{\U}
\func{\partialder{k}\cf}{z}:=\func{\[\func{\(\banachder{\cf}{\R^n}{\R}\)}{z}\]}{\Eucbase{n}{k}}.
\end{align}
For every $\x$ in $\U$, $\derived{\xi_\x}$ is continuous because it is differentiable,
and hence it is integrable. Furthermore, based on the the fundamental theorem of calculus,
\begin{equation}\label{lemlocalpowersumpeq6}
\Foreach{\x}{\U}
\func{\xi_\x}{1}-\func{\xi_\x}{0}=\int_{0}^{1}\func{\derived{\xi_\x}}{\tau}{\rm{d}\tau}.
\end{equation}
In addition,
\begin{align}\label{lemlocalpowersumpeq7}
\Foreach{\x}{\U}
\func{\xi_\x}{1}-\func{\xi_\x}{0}=\func{\cf}{\x}-\func{\cf}{\y}.
\end{align}
\Ref{lemlocalpowersumpeq1}, \Ref{lemlocalpowersumpeq2}, \Ref{lemlocalpowersumpeq6}, and \Ref{lemlocalpowersumpeq7} imply,
\begin{align}\label{lemlocalpowersumpeq8}
\Foreach{\x}{\U}
\func{\cf}{\x}=\func{\cf}{\y}+
\sum_{k=1}^{n}\[\func{\projection{n}{k}}{\x}-\func{\projection{n}{k}}{\y}\]
\int_{0}^{1}\func{\alpha_\x^{\(k\)}}{\tau}{\rm{d}\tau},
\end{align}
where, the smooth (and hence continuous) mappings $\function{\alpha_\x^{\(k\)}}{\cinterval{0}{1}}{\R}$ are defined as,
\begin{equation}\label{lemlocalpowersumpeq9}
\Foreach{k}{\seta{\suc{1}{n}}}
\Foreach{\x}{\U}
\Foreach{t}{\cinterval{0}{1}}
\func{\alpha_\x^{\(k\)}}{t}\eqdef\func{\partialder{k}\cf}{\(1-t\)\y+t\x}.
\end{equation}
By defining the mappings $\function{\cg_{k}}{\U}{\R}$ as,
\begin{equation}\label{lemlocalpowersumpeq10}
\Foreach{k}{\seta{\suc{1}{n}}}
\Foreach{\x}{\U}
\func{\cg_k}{\x}:=\int_{0}^{1}\func{\alpha_\x^{\(k\)}}{\tau}{\rm{d}\tau},
\end{equation}
\Ref{lemlocalpowersumpeq8} becomes,
\begin{equation}\label{lemlocalpowersumpeq11}
\Foreach{\x}{\U}
\func{\cf}{\x}=\func{\cf}{\y}+
\sum_{k=1}^{n}\[\func{\projection{n}{k}}{\x}-\func{\projection{n}{k}}{\y}\]\func{\cg_{k}}{\x}.
\end{equation}
Furthermore, \Ref{lemlocalpowersumpeq5}, \Ref{lemlocalpowersumpeq9},
and \Ref{lemlocalpowersumpeq10} imply,
\begin{equation}\label{lemlocalpowersumpeq12}
\Foreach{k}{\seta{\suc{1}{n}}}\func{\cg_k}{\y}=\func{\[\func{\(\banachder{\cf}{\R^n}{\R}\)}{\y}\]}{\Eucbase{n}{k}}.
\end{equation}
By defining the mappings $\function{\phi_{k}}{\Cprod{\U}{\cinterval{0}{1}}}{\R}$ as,
\begin{equation}\label{lemlocalpowersumpeq13}
\Foreach{k}{\seta{\suc{1}{n}}}
\Foreach{\x}{\U}
\Foreach{t}{\cinterval{0}{1}}
\func{\phi_{k}}{\binary{\x}{t}}\eqdef
\func{\partialder{k}\cf}{\(1-t\)\y+t\x},
\end{equation}
\Ref{lemlocalpowersumpeq10} becomes,
\begin{equation}
\Foreach{k}{\seta{\suc{1}{n}}}
\Foreach{\x}{\U}
\func{\cg_k}{\x}:=\int_{0}^{1}\func{\phi_{k}}{\binary{\x}{\tau}}{\rm{d}\tau}.
\end{equation}
Also it is clear that each $\phi_{k}$ is infinitely differentiable with respect to its first factor.
So according to \cite[page~35,~COROLLARY~2.12.3]{Cartandifforms}, each $\cg_{k}$ is infinitely differentiable. That is,
\begin{equation}\label{lemlocalpowersumpeq14}
\Foreach{k}{\seta{\suc{1}{n}}}
\cg_{k}\in\banachmapdifclass{\infty}{\R^n}{\R}{\U}{\R}.
\end{equation}
\endlem
\corollary\label{corneighborhoodofderivativeoftransitionmaps}
$\point$ is taken as a point of $\Man{}$, and
each $\phi$ and $\psi$ as an element of $\maxatlas{}$ (a chart of $\Man{}$) such that $\point\in\domain{\phi}\cap\domain{\psi}$
and $\func{\phi}{\point}=\func{\psi}{\point}=\zerovec{}$. There exists a convex and open subset $\U$ of $\R^{n}$ such that
$\func{\phi}{\point}\in\U\subseteq\func{\image{\phi}}{\domain{\phi}\cap\domain{\psi}}$,
\begin{align}
&\Foreach{j}{\seta{\suc{1}{n}}}\Exists{\mtuple{\lambda_1^{j}}{\lambda_n^{j}}}{{\banachmapdifclass{\infty}{\R^n}{\R}{\U}{\R}}^{\times n}}\cr
&\Foreach{j,~k}{\seta{\suc{1}{n}}}\func{\(\cmp{\lambda_{k}^{j}}{\bar{\phi}}\)}{\point}=
\func{\[\func{\(\banachder{\(\cmp{\projection{n}{j}}{\cmp{\psi}{\finv{\phi}}}\)}{\R^n}{\R}\)}{\func{\phi}{\point}}\]}{\Eucbase{n}{k}},\cr
&\Foreach{\x}{\func{\pimage{\phi}}{\U}}\func{\(\cmp{\projection{n}{j}}{\psi}\)}{\x}=
\sum_{k=1}^{n}\[\func{\(\cmp{\projection{n}{k}}{\phi}\)}{\x}\]\[\func{\(\cmp{\lambda_{k}^{j}}{\bar{\phi}}\)}{\x}\],
\end{align}
where,
\begin{equation}
\bar{\phi}:=\func{\res{\phi}}{\func{\pimage{\phi}}{\U}}\in\defset{\p{\phi}}{\maxatlas{}}{\point\in\domain{\p{\phi}}}.
\end{equation}
\endcor
\corollary\label{corneighborhoodofderivativeofrealvaluedmaps}
$\point$ is taken as a point of $\Man{}$.
$\phi$ is taken as an element of $\maxatlas{}$ (a chart of $\Man{}$) such that $\point\in\domain{\phi}$,
$\funcimage{\phi}$ is a covex set of $\R^n$,
and $\func{\phi}{\point}=\zerovec{}$. $\cf$ is taken as an element of $\mapdifclass{\infty}{\Man{}}{\RR}$.
\begin{align}
&\Exists{\mtuple{\lambda_1}{\lambda_n}}{{\banachmapdifclass{\infty}{\R^n}{\R}{\funcimage{\phi}}{\R}}^{\times n}}\cr
&\Foreach{k}{\seta{\suc{1}{n}}}\func{\(\cmp{\lambda_{k}}{\phi}\)}{\point}=
\func{\[\func{\(\banachder{\(\cmp{\cf}{\finv{\phi}}\)}{\R^n}{\R}\)}{\func{\phi}{\point}}\]}{\Eucbase{n}{k}},\cr
&\Foreach{\x}{\domain{\phi}}\func{\cf}{\x}=
\sum_{k=1}^{n}\[\func{\(\cmp{\projection{n}{k}}{\phi}\)}{\x}\]\[\func{\(\cmp{\lambda_{k}}{\phi}\)}{\x}\].
\end{align}
\endcor
\lemma\label{lemmaprederivationinducedvectorfield1}
$\aderivation{}$ is taken as an arbitrary element of $\Derivation{\Man{}}{\infty}$, and $\point$ as a point of $\Man{}$.
The value of $\displaystyle\sum_{k=1}^{n}\bigg(\func{\[\func{\aderivation{}}{\cf_k}\]}{\point}\bigg)
\bigg(\func{\finv{\[\tanspaceiso{\point}{\Man{}}{\phi}\]}}{\Eucbase{n}{k}}\bigg)$ is independent of the
choice of a chart $\phi$ of $\Man{}$ centered at $\point$ and choices of smooth extension of
$\suc{\cmp{\projection{n}{1}}{\phi}}{\cmp{\projection{n}{n}}{\phi}}$ on $\Man{}$ fixed at $\point$.
That is,
\begin{align}
&\Foreach{\opair{\phi}{\psi}}
{{\defset{\p{\phi}}{\maxatlas{}}{\point\in\domain{\p{\phi}},~\func{\p{\phi}}{\point}=\zerovec{}}}^{\times 2}}\cr
&\Foreach{\mtuple{\cf_1}{\cf_n}}
{\Times{\fextension{\Man{}}{\point}{\phi}{\projection{n}{1}}}{\fextension{\Man{}}{\point}{\phi}{\projection{n}{n}}}}\cr
&\Foreach{\mtuple{\cg_1}{\cg_n}}
{\Times{\fextension{\Man{}}{\point}{\psi}{\projection{n}{1}}}{\fextension{\Man{}}{\point}{\psi}{\projection{n}{n}}}}\cr
&\sum_{k=1}^{n}\bigg(\func{\[\func{\aderivation{}}{\cf_k}\]}{\point}\bigg)
\bigg(\func{\finv{\[\tanspaceiso{\point}{\Man{}}{\phi}\]}}{\Eucbase{n}{k}}\bigg)=
\sum_{k=1}^{n}\bigg(\func{\[\func{\aderivation{}}{\cg_k}\]}{\point}\bigg)
\bigg(\func{\finv{\[\tanspaceiso{\point}{\Man{}}{\psi}\]}}{\Eucbase{n}{k}}\bigg).
\end{align}
\proof
Each $\phi$ and $\psi$ is taken as an arbitrary element of
$\defset{\p{\phi}}{\maxatlas{}}{\point\in\domain{\p{\phi}},~\func{\p{\phi}}{\point}=\zerovec{}}$.
$\suc{\cf_1}{\cf_n}$ are taken as arbitrary elements of
$\suc{\fextension{\Man{}}{\point}{\phi}{\projection{n}{1}}}{\fextension{\Man{}}{\point}{\phi}{\projection{n}{n}}}$, respectively, and
$\V$ is taken as such an open set of $\mantops{\Man{}}$ such that,
\begin{equation}\label{lemmaprederivationinducedvectorfield1peq1}
\Foreach{k}{\seta{\suc{1}{n}}}
\[\point\in\V\subseteq\domain{\phi},~
\func{\resd{\cf_k}}{\V}=\func{\resd{\cmp{\projection{n}{k}}{\phi}}}{\V}\].
\end{equation}
$\suc{\cg_1}{\cg_n}$ are taken as arbitrary elements of
$\suc{\fextension{\Man{}}{\point}{\psi}{\projection{n}{1}}}{\fextension{\Man{}}{\point}{\psi}{\projection{n}{n}}}$, respectively, and
$\W{}$ is taken as such an open set of $\mantops{\Man{}}$ such that,
\begin{equation}\label{lemmaprederivationinducedvectorfield1peq2}
\Foreach{k}{\seta{\suc{1}{n}}}
\[\point\in\W{}\subseteq\domain{\psi},~
\func{\resd{\cg_k}}{\W{}}=\func{\resd{\cmp{\projection{n}{k}}{\psi}}}{\W{}}\].
\end{equation}
According to \refcor{corneighborhoodofderivativeoftransitionmaps},
there exists an open set $\p{\U}$ of $\mantops{\Man{}}$ such that
$\point\in{\p{\U}}\subseteq\domain{\phi}\cap\domain{\psi}$,
\begin{align}\label{lemmaprederivationinducedvectorfield1peq3}
&\Foreach{j}{\seta{\suc{1}{n}}}\Existsis{\mtuple{\lambda_1^{j}}{\lambda_n^{j}}}{{\banachmapdifclass{\infty}{\R^n}{\R}{\func{\image{\phi}}{\p{\U}}}{\R}}^{\times n}}\cr
&\Foreach{j,~k}{\seta{\suc{1}{n}}}\func{\(\cmp{\lambda_{k}^{j}}{\bar{\phi}}\)}{\point}=
\func{\[\func{\(\banachder{\(\cmp{\projection{n}{j}}{\cmp{\psi}{\finv{\phi}}}\)}{\R^n}{\R}\)}{\func{\phi}{\point}}\]}{\Eucbase{n}{k}},\cr
&\Foreach{\x}{\p{\U}}\func{\(\cmp{\projection{n}{j}}{\psi}\)}{\x}=
\sum_{k=1}^{n}\[\func{\(\cmp{\projection{n}{k}}{\phi}\)}{\x}\]\[\func{\(\cmp{\lambda_{k}^{j}}{\bar{\phi}}\)}{\x}\],
\end{align}
where,
\begin{equation}
\bar{\phi}:=\func{\res{\phi}}{\p{\U}}\in\defset{\p{\phi}}{\maxatlas{}}{\point\in\domain{\p{\phi}}}.
\end{equation}
For every $j$ and $k$ in $\seta{\suc{1}{n}}$, $\suc{F_1^1}{F_n^n}$ are taken as representative elements of\\
$\suc{\fextension{\Man{}}{\point}{\bar{\phi}}{\lambda_1^1}}{\fextension{\Man{}}{\point}{\bar{\phi}}{\lambda_n^n}}$, respectively, and
$\U$ is taken as such an element of $\mantops{\Man{}}$ such that,
\begin{equation}\label{lemmaprederivationinducedvectorfield1peq4}
\Foreach{j,~k}{\seta{\suc{1}{n}}}
\[\point\in\U\subseteq\p{\U},~
\func{\resd{F_k^j}}{\U}=\func{\resd{\cmp{\lambda_k^j}{\phi}}}{\U}\].
\end{equation}
Therefore by defining,
$O:=\V\cap\W{}\cap\U\cap\p{\U}$,
evidently $\point\in O\in\mantop{\Man{}}$ and,
\begin{align}\label{lemmaprederivationinducedvectorfield1peq5}
\Foreach{j}{\seta{\suc{1}{n}}}
\Foreach{\x}{O}
\func{\cg_j}{\x}&=\sum_{k=1}^{n}\func{\cf_k}{\x}\func{F_k^j}{\x}\cr
&=\func{\(\sum_{k=1}^{n}\cf_{k}\rdot F_{k}^{j}\)}{\x},
\end{align}
and hence according to \reflem{lemmaderivationislocal}, \Ref{lemmaprederivationinducedvectorfield1peq1},
\Ref{lemmaprederivationinducedvectorfield1peq3}, \Ref{lemmaprederivationinducedvectorfield1peq4},
properties of an $\infty$-derivation on $\Man{}$, chain-rule of differentiation in Banach-spaces
and the linearity of the derivative of a differentiable map between Banach-spaces at any point of its domain,
considering that $\func{\phi}{\point}=\zerovec{}$,
and also considering the relation
$\[\func{\(\banachder{\(\cmp{\psi}{\finv{\phi}}\)}{\R^n}{\R^n}\)}{\func{\phi}{\point}}\]=
\cmp{\[\tanspaceiso{\point}{\Man{}}{\psi}\]}{\finv{\[\tanspaceiso{\point}{\Man{}}{\phi}\]}}$ and
linearity of $\finv{\[\tanspaceiso{\point}{\Man{}}{\phi}\]}$,
\begin{align}\label{lemmaprederivationinducedvectorfield1peq6}
&\begin{aligned}
\Foreach{j}{\seta{\suc{1}{n}}}
\end{aligned}\cr
&\begin{aligned}
\func{\[\func{\aderivation{}}{\cg_j}\]}{\point}&=
\func{\[\func{\aderivation{}}{\sum_{k=1}^{n}\cf_{k}\rdot F_{k}^{j}}\]}{\point}\cr
&=\sum_{k=1}^{n}\func{\[\func{\aderivation{}}{\cf_k}\]}{\point}\func{F_{k}^{j}}{\point}+
\func{\cf_k}{\point}\func{\[\func{\aderivation{}}{F_k^j}\]}{\point}\cr
&=\sum_{k=1}^{n}\func{\[\func{\aderivation{}}{\cf_k}\]}{\point}
\bigg(\func{\[\func{\(\banachder{\(\cmp{\projection{n}{j}}
{\cmp{\psi}{\finv{\phi}}}\)}{\R^n}{\R}\)}{\func{\phi}{\point}}\]}{\Eucbase{n}{k}}\bigg)\cr
&=\[\func{\(\banachder{\(\cmp{\projection{n}{j}}
{\cmp{\psi}{\finv{\phi}}}\)}{\R^n}{\R}\)}{\func{\phi}{\point}}\]
\bigg(\sum_{k=1}^{n}\func{\[\func{\aderivation{}}{\cf_k}\]}{\point}\Eucbase{n}{k}\bigg)\cr
&=\bigg(\cmp{\projection{n}{j}}{\[\func{\(\banachder{\(\cmp{\psi}{\finv{\phi}}\)}{\R^n}{\R^n}\)}{\func{\phi}{\point}}\]}\bigg)
\bigg(\sum_{k=1}^{n}\func{\[\func{\aderivation{}}{\cf_k}\]}{\point}\Eucbase{n}{k}\bigg)\cr
&=\bigg(\cmp{\projection{n}{j}}{\cmp{\[\tanspaceiso{\point}{\Man{}}{\psi}\]}{\finv{\[\tanspaceiso{\point}{\Man{}}{\phi}\]}}}\bigg)
\bigg(\sum_{k=1}^{n}\func{\[\func{\aderivation{}}{\cf_k}\]}{\point}\Eucbase{n}{k}\bigg)\cr
&=\bigg(\cmp{\projection{n}{j}}{\[\tanspaceiso{\point}{\Man{}}{\psi}\]}\bigg)
\[\sum_{k=1}^{n}\bigg(\func{\[\func{\aderivation{}}{\cf_k}\]}{\point}\bigg)
\bigg(\func{\finv{\[\tanspaceiso{\point}{\Man{}}{\phi}\]}}{\Eucbase{n}{k}}\bigg)\].
\end{aligned}
\end{align}
Thus by defining,
\begin{equation}
\tanspace{\point}{\Man{}}\ni
\alpha:=\sum_{k=1}^{n}\bigg(\func{\[\func{\aderivation{}}{\cf_k}\]}{\point}\bigg)
\bigg(\func{\finv{\[\tanspaceiso{\point}{\Man{}}{\phi}\]}}{\Eucbase{n}{k}}\bigg),
\end{equation}
and using the linearity of $\finv{\[\tanspaceiso{\point}{\Man{}}{\phi}\]}$,
\begin{align}
\sum_{j=1}^{n}\func{\[\func{\aderivation{}}{\cg_j}\]}{\point}\bigg(\func{\finv{\[\tanspaceiso{\point}{\Man{}}{\psi}\]}}{\Eucbase{n}{j}}\bigg)&=
\sum_{j=1}^{n}\[\func{\bigg(\cmp{\projection{n}{j}}{\[\tanspaceiso{\point}{\Man{}}{\psi}\]}\bigg)}{\alpha}\]
\bigg(\func{\finv{\[\tanspaceiso{\point}{\Man{}}{\psi}\]}}{\Eucbase{n}{j}}\bigg)\cr
&=\func{\finv{\[\tanspaceiso{\point}{\Man{}}{\psi}\]}}{\sum_{j=1}^{n}
{\[\func{\bigg(\cmp{\projection{n}{j}}{\[\tanspaceiso{\point}{\Man{}}{\psi}\]}\bigg)}{\alpha}\]\Eucbase{n}{j}}}\cr
&=\func{\finv{\[\tanspaceiso{\point}{\Man{}}{\psi}\]}}{\func{\[\tanspaceiso{\point}{\Man{}}{\psi}\]}{\alpha}}\cr
&=\alpha.
\end{align}
So,
\begin{equation}
\sum_{k=1}^{n}\bigg(\func{\[\func{\aderivation{}}{\cf_k}\]}{\point}\bigg)
\bigg(\func{\finv{\[\tanspaceiso{\point}{\Man{}}{\phi}\]}}{\Eucbase{n}{k}}\bigg)=
\sum_{k=1}^{n}\bigg(\func{\[\func{\aderivation{}}{\cg_k}\]}{\point}\bigg)
\bigg(\func{\finv{\[\tanspaceiso{\point}{\Man{}}{\psi}\]}}{\Eucbase{n}{k}}\bigg).
\end{equation}
\endlem
\definition\label{definvliederivative}
The mapping $\LieDerinv{\Man{}}$ is defined as,
\begin{align}
&\LieDerinv{\Man{}}\indef\Func{\Derivation{\Man{}}{\infty}}{\Func{\M{}}{\tanbun{\Man{}}}},\cr
&\Foreach{\aderivation{}}{\Derivation{\Man{}}{\infty}}\Foreach{\point}{\M{}}
\func{\[\func{\LieDerinv{\Man{}}}{\aderivation{}}\]}{\point}\eqdef
\sum_{k=1}^{n}\bigg(\func{\[\func{\aderivation{}}{\cf_k}\]}{\point}\bigg)
\bigg(\func{\finv{\[\tanspaceiso{\point}{\Man{}}{\phi}\]}}{\Eucbase{n}{k}}\bigg),~~~~
\end{align}
where $\phi$ is an element of $\defset{\p{\phi}}{\maxatlas{}}{\point\in\domain{\p{\phi}},~\func{\p{\phi}}{\point}=\zerovec{}}$,
and $\suc{\cf_1}{\cf_n}$ are elements of
$\suc{\fextension{\Man{}}{\point}{\phi}{\projection{n}{1}}}{\fextension{\Man{}}{\point}{\phi}{\projection{n}{n}}}$, respectively.
\endef
\lemma\label{lemtotalsmoothextension}
$\lambda$ is taken as an element of $\mapdifclass{\infty}{\RR^n}{\RR}$ (a smooth map from $\R^n$ to $\R$).
For every point $\point$ of $\Man{}$ and every chart $\psi$ of $\Man{}$ centered at $\point$,
there exists a chart $\phi$ of $\Man{}$ centered at $\point$ which is the restriction of $\psi$ to an open set of $\mantops{\Man{}}$,
and there exists a smooth extension
$\cf$ of $\cmp{\lambda}{\phi}$ on $\Man{}$ fixed at $\point$ such that $\cf$ coincides with $\cmp{\lambda}{\phi}$
totally on the domain of $\cmp{\lambda}{\phi}$ (which is the same as the domain of $\phi$). That is,
\begin{align}
&\Foreach{\point}{\M{}}\Foreach{\psi}{\defset{\p{\phi}}{\maxatlas{}}{\point\in\domain{\phi},~\func{\p{\phi}}{\point}=\zerovec{}}}\cr
&\Exists{\phi}{\defset{\p{\phi}}{\maxatlas{}}{\point\in\domain{\phi},~\func{\p{\phi}}{\point}=\zerovec{}}}
\Exists{\cf}{\fextension{\Man{}}{\point}{\phi}{\lambda}}\cr
&\bigg(\func{\resd{\cf}}{\domain{\phi}}=\cmp{\lambda}{\phi},~\domain{\phi}\subseteq\domain{\psi},~\phi=\func{\res{\psi}}{\domain{\phi}}\bigg).
\end{align}
\proof
$\point$ is taken as an arbitrary point of $\Man{}$.
The set $\defset{\p{\phi}}{\maxatlas{}}{\point\in\domain{\phi},~\func{\p{\phi}}{\point}=\zerovec{}}$ is non-empty
and $\psi$ is taken as an element of it. According to \refdef{defsmoothextensionoflocalmaps} and
\reflem{lemmadifferentiablemapexpansion},
\begin{equation}\label{lemtotalsmoothextensionpeq1}
\fextension{\Man{}}{\point}{\psi}{\lambda}\neq\empty.
\end{equation}
$\cf$ is taken as an element of $\fextension{\Man{}}{\point}{\psi}{\lambda}$. So according to
\refdef{defsmoothextensionoflocalmaps},
\begin{equation}\label{lemtotalsmoothextensionpeq2}
\Existsis{\V}{\mantop{\Man{}}}
\[\point\in\V\subseteq\domain{\psi},~
\func{\resd{\cf}}{\V}=\func{\resd{\cmp{\lambda}{\psi}}}{\V}\].
\end{equation}
By defining,
\begin{equation}\label{lemtotalsmoothextensionpeq3}
\phi:=\func{\res{\psi}}{\V},
\end{equation}
it is obvious that $\phi$ is also a chart of $\Man{}$ centered at $\point$, that is,
\begin{equation}\label{lemtotalsmoothextensionpeq4}
\phi\in\defset{\p{\phi}}{\maxatlas{}}{\point\in\domain{\phi},~\func{\p{\phi}}{\point}=\zerovec{}}.
\end{equation}
In addition, \Ref{lemtotalsmoothextensionpeq3} implies,
\begin{equation}\label{lemtotalsmoothextensionpeq5}
\cmp{\lambda}{\phi}=\func{\resd{\cmp{\lambda}{\psi}}}{\V},
\end{equation}
and hence according to \Ref{lemtotalsmoothextensionpeq2},
\begin{equation}
\func{\resd{\cf}}{\V}=\cmp{\lambda}{\phi}.
\end{equation}
Also according to \refdef{defsmoothextensionoflocalmaps}, it is clear that,
\begin{equation}
\cf\in\fextension{\Man{}}{\point}{\phi}{\lambda}.
\end{equation}
\endlem
\lemma\label{lemmainvliederivativeissmooth}
For every $\infty$-derivation $\aderivation{}$ on $\Man{}$, $\func{\LieDerinv{\Man{}}}{\aderivation{}}$ is a
smooth vector-field on $\Man{}$. That is,
\begin{equation}
\Foreach{\aderivation{}}{\Derivation{\Man{}}{\infty}}
\func{\LieDerinv{\Man{}}}{\aderivation{}}\in\vecf{\Man{}}{\infty}.
\end{equation}
\proof
$\aderivation{}$ is taken as an arbitrary element of $\Derivation{\Man{}}{\infty}$.
According to \refdef{definvliederivative}, obviously for every point $\point$ of $\Man{}$,
$\func{\[\func{\LieDerinv{\Man{}}}{\aderivation{}}\]}{\point}\in\tanspace{\point}{\Man{}}$, and thus,
\begin{equation}\label{lemmainvliederivativeissmoothpeq1}
\cmp{\basep{\Man{}}}{\[\func{\LieDerinv{\Man{}}}{\aderivation{}}\]}=\identity{\M{}}.
\end{equation}
\begin{itemize}
\item[$\pr{1}$]
$\point$ is taken as an arbitrary point of $\Man{}$, and a $\psi$ is chosen from the non-empty set\\
$\defset{\p{\phi}}{\maxatlas{}}{\point\in\domain{\phi},~\func{\p{\phi}}{\point}=\zerovec{}}$.
According to \reflem{lemtotalsmoothextension},
\begin{align}\label{lemmainvliederivativeissmoothp1eq1}
&\Foreach{k}{\seta{\suc{1}{n}}}\cr
&\Existsis{\phi_{k}}{\defset{\p{\phi}}{\maxatlas{}}{\point\in\domain{\phi},~\func{\p{\phi}}{\point}=\zerovec{}}}
\Existsis{\cf_{k}}{\fextension{\Man{}}{\point}{\phi_{k}}{\projection{n}{k}}}\cr
&\bigg(\func{\resd{\cf_{k}}}{\domain{\phi_{k}}}=\cmp{\projection{n}{k}}{\phi_{k}},~\domain{\phi_{k}}\subseteq\domain{\psi},~
\phi_{k}=\func{\res{\psi}}{\domain{\phi_{k}}}\bigg).
\end{align}
By defining $\displaystyle\U:=\bigcap_{k=1}^{n}\domain{\phi_{k}}$ which is obviously an open set of $\mantops{\Man{}}$,
since each $\phi_{k}$ is the restriction of $\psi$ to $\V_{k}$ and $\U$ is a subset of each $\V_{k}$, it is clear that,
\begin{equation}\label{lemmainvliederivativeissmoothp1eq2}
\Foreach{k}{\seta{\suc{1}{n}}}
\func{\res{\phi_{k}}}{\U}=\phi,
\end{equation}
where,
\begin{equation}\label{lemmainvliederivativeissmoothp1eq3}
\phi:=\func{\res{\psi}}{\U}.
\end{equation}
Since $\U$ is an open set of $\mantops{\Man{}}$ containing $\point$ and $\psi$ is a chert of $\Man{}$
centered at $\point$, it is evident that $\phi$ is also a chart of $\Man{}$ centered at $\point$, that is,
\begin{equation}\label{lemmainvliederivativeissmoothp1eq4}
\phi\in\defset{\p{\phi}}{\maxatlas{}}{\point\in\domain{\phi},~\func{\p{\phi}}{\point}=\zerovec{}}.
\end{equation}
In addition, since for each $k$ in $\seta{\suc{1}{n}}$, $\func{\resd{\cf_{k}}}{\domain{\phi_{k}}}=\cmp{\projection{n}{k}}{\phi_{k}}$
and $\U$ is a subset of $\domain{\phi_{k}}$, it is evident that,
\begin{equation}\label{lemmainvliederivativeissmoothp1eq5}
\Foreach{k}{\seta{\suc{1}{n}}}
\func{\resd{\cf_{k}}}{\U}=\cmp{\projection{n}{k}}{\phi},
\end{equation}
and hence according to \refdef{defsmoothextensionoflocalmaps}
it can be easily seen that each $\cf_{k}$ is a smooth extension of $\cmp{\projection{n}{k}}{\phi}$ on $\Man{}$
fixed at any point in $\U$. That is,
\begin{equation}\label{lemmainvliederivativeissmoothp1eq6}
\Foreach{\p{\point}}{\U}
\Foreach{k}{\seta{\suc{1}{n}}}
\cf_{k}\in\fextension{\Man{}}{\p{\point}}{\phi}{\projection{n}{k}}.
\end{equation}
Therefore according to \refdef{definvliederivative},
\begin{equation}\label{lemmainvliederivativeissmoothp1eq7}
\Foreach{\p{\point}}{\U}
\func{\[\func{\LieDerinv{\Man{}}}{\aderivation{}}\]}{\p{\point}}\eqdef
\sum_{k=1}^{n}\bigg(\func{\[\func{\aderivation{}}{\cf_k}\]}{\p{\point}}\bigg)
\bigg(\func{\finv{\[\tanspaceiso{\p{\point}}{\Man{}}{\phi}\]}}{\Eucbase{n}{k}}\bigg).
\end{equation}
On the other hand, since $\func{\LieDerinv{\Man{}}}{\aderivation{}}$ sends every point of $\Man{}$ to the tangent-space of
$\Man{}$ associated with that point, it is clear that,
\begin{equation}\label{lemmainvliederivativeissmoothp1eq8}
\func{\image{\[\func{\LieDerinv{\Man{}}}{\aderivation{}}\]}}{U}\subseteq
\func{\pimage{\basep{\Man{}}}}{\U}.
\end{equation}
and hence considering that $\domain{\tanchart{\Man{}}{\phi}}=\func{\pimage{\basep{\Man{}}}}{\U}$,
\begin{equation}\label{lemmainvliederivativeissmoothp1eq9}
\func{\image{\[\func{\LieDerinv{\Man{}}}{\aderivation{}}\]}}{\domain{\phi}}\subseteq
\domain{\tanchart{\Man{}}{\phi}}.
\end{equation}
It is known that,
\begin{equation}\label{lemmainvliederivativeissmoothp1eq10}
\funcimage{\tanchart{\Man{}}{\phi}}=\Cprod{\funcimage{\phi}}{\R^n}=\Cprod{\func{\image{\phi}}{\U}}{\R^n},
\end{equation}
and hence
\begin{equation}\label{lemmainvliederivativeissmoothp1eq11}
\cmp{\tanchart{\Man{}}{\phi}}{\cmp{\[\func{\LieDerinv{\Man{}}}{\aderivation{}}\]}{\finv{\phi}}}
\in\Func{\func{\image{\phi}}{\U}}{\Cprod{\func{\image{\phi}}{\U}}{\R^n}}.
\end{equation}
By invoking the definition of the tangent-bundle chart of $\Man{}$ associated with $\phi$ and using
\Ref{lemmainvliederivativeissmoothpeq1},
\begin{align}\label{lemmainvliederivativeissmoothp1eq12}
&\begin{aligned}
\Foreach{\x}{\func{\image{\phi}}{\U}}
\end{aligned}\cr
&\begin{aligned}
&~~~\func{\(\cmp{\tanchart{\Man{}}{\phi}}{\cmp{\[\func{\LieDerinv{\Man{}}}{\aderivation{}}\]}{\finv{\phi}}}\)}{\x}\cr
&=\opair{\func{\(\cmp{\phi}{\basep{\Man{}}}\)}
{\func{\[\func{\LieDerinv{\Man{}}}{\aderivation{}}\]}
{\func{\finv{\phi}}{\x}}}}{\func{\(\tanspaceiso{\func{\basep{\Man{}}}
{\func{\[\func{\LieDerinv{\Man{}}}{\aderivation{}}\]}{\func{\finv{\phi}}{\x}}}}{\Man{}}{\phi}\)}
{\func{\[\func{\LieDerinv{\Man{}}}{\aderivation{}}\]}{\func{\finv{\phi}}{\x}}}}\cr
&=\opair{\x}{\func{\(\tanspaceiso{\func{\finv{\phi}}{\x}}{\Man{}}{\phi}\)}
{\func{\[\func{\LieDerinv{\Man{}}}{\aderivation{}}\]}{\func{\finv{\phi}}{\x}}}}.
\end{aligned}\cr
&{}
\end{align}
So by defining the map $\alpha$ as,
\begin{align}\label{lemmainvliederivativeissmoothp1eq13}
&\alpha\indef\Func{\func{\image{\phi}}{\U}}{\R^n},\cr
&\Foreach{\x}{\func{\image{\phi}}{\U}}
\func{\alpha}{\x}\eqdef\func{\(\tanspaceiso{\func{\finv{\phi}}{\x}}{\Man{}}{\phi}\)}
{\func{\[\func{\LieDerinv{\Man{}}}{\aderivation{}}\]}{\func{\finv{\phi}}{\x}}},
\end{align}
\Ref{lemmainvliederivativeissmoothp1eq12} becomes,
\begin{equation}\label{lemmainvliederivativeissmoothp1eq14}
\Foreach{\x}{\func{\image{\phi}}{\U}}
\func{\(\cmp{\tanchart{\Man{}}{\phi}}{\cmp{\[\func{\LieDerinv{\Man{}}}{\aderivation{}}\]}{\finv{\phi}}}\)}{\x}=
\opair{\func{\identity{\func{\image{\phi}}{\U}}}{\x}}{\func{\alpha}{\x}}.
\end{equation}
Therefore,
since $\identity{\func{\image{\phi}}{\U}}\in\banachmapdifclass{\infty}{\R^n}{\R^n}{\func{\image{\phi}}{\U}}{\func{\image{\phi}}{\U}}$,
as a result of a theorem in differential calculus,
\begin{equation}\label{lemmainvliederivativeissmoothp1eq15}
\(\cmp{\tanchart{\Man{}}{\phi}}{\cmp{\[\func{\LieDerinv{\Man{}}}{\aderivation{}}\]}{\finv{\phi}}}\)\in
\banachmapdifclass{\infty}{\R^n}{\R^n}{\func{\image{\phi}}{\U}}{\Cprod{\func{\image{\phi}}{\U}}{\R^n}}
\thenn
\alpha\in\banachmapdifclass{\infty}{\R^n}{\R^n}{\func{\image{\phi}}{\U}}{\R^n},
\end{equation}
which simply means $\(\cmp{\tanchart{\Man{}}{\phi}}{\cmp{\[\func{\LieDerinv{\Man{}}}{\aderivation{}}\]}{\finv{\phi}}}\)$
is smooth if and only if so is $\alpha$.
By exploiting the linearity of $\finv{\[\tanspaceiso{\func{\finv{\phi}}{\x}}{\Man{}}{\phi}\]}$,
\Ref{lemmainvliederivativeissmoothp1eq7} and \Ref{lemmainvliederivativeissmoothp1eq13} imply,
\begin{equation}\label{lemmainvliederivativeissmoothp1eq16}
\Foreach{\x}{\func{\image{\phi}}{\U}}
\func{\alpha}{\x}=
\sum_{k=1}^{n}\bigg(\func{\[\func{\aderivation{}}{\cf_k}\]}{\func{\finv{\phi}}{\x}}\bigg)
{\Eucbase{n}{k}},
\end{equation}
and hence,
\begin{equation}\label{lemmainvliederivativeissmoothp1eq17}
\Foreach{k}{\seta{\suc{1}{n}}}
\cmp{\projection{n}{k}}{\alpha}=\cmp{\[\func{\aderivation{}}{\cf_k}\]}{\finv{\phi}}.
\end{equation}
Also since for each $k$ in $\seta{\suc{1}{n}}$, $\[\func{\aderivation{}}{\cf_k}\]\in\mapdifclass{\infty}{\Man{}}{\RR}$
and $\phi$ is a chart of $\Man{}$,
\begin{equation}\label{lemmainvliederivativeissmoothp1eq18}
\Foreach{k}{\seta{\suc{1}{n}}}
\cmp{\[\func{\aderivation{}}{\cf_k}\]}{\finv{\phi}}\in
\banachmapdifclass{\infty}{\R^n}{\R}{\func{\image{\phi}}{\U}}{\R},
\end{equation}
Therefore,
\begin{equation}\label{lemmainvliederivativeissmoothp1eq19}
\Foreach{k}{\seta{\suc{1}{n}}}
\cmp{\projection{n}{k}}{\alpha}\in
\banachmapdifclass{\infty}{\R^n}{\R}{\func{\image{\phi}}{\U}}{\R},
\end{equation}
and hence obviously,
\begin{equation}\label{lemmainvliederivativeissmoothp1eq20}
\alpha\in\banachmapdifclass{\infty}{\R^n}{\R^n}{\func{\image{\phi}}{\U}}{\R^n},
\end{equation}
which according to \Ref{lemmainvliederivativeissmoothp1eq15} implies,
\begin{equation}\label{lemmainvliederivativeissmoothp1eq21}
\(\cmp{\tanchart{\Man{}}{\phi}}{\cmp{\[\func{\LieDerinv{\Man{}}}{\aderivation{}}\]}{\finv{\phi}}}\)\in
\banachmapdifclass{\infty}{\R^n}{\R^n}{\func{\image{\phi}}{\U}}{\Cprod{\func{\image{\phi}}{\U}}{\R^n}}.
\end{equation}
\endp
\end{itemize}
So,
\begin{align}\label{lemmainvliederivativeissmoothpeq2}
&\Foreach{\point}{\M{}}\cr
&~~~\Existsis{\phi}{\defset{\p{\phi}}{\maxatlas{}}{\point\in\domain{\p{\phi}}}}
\Existsis{\tanchart{\Man{}}{\phi}}{\tanatlas{\Man{}}}\cr
&~~~\func{\image{\[\func{\LieDerinv{\Man{}}}{\aderivation{}}\]}}{\domain{\phi}}\subseteq
\domain{\tanchart{\Man{}}{\phi}},\cr
&~~~\(\cmp{\tanchart{\Man{}}{\phi}}{\cmp{\[\func{\LieDerinv{\Man{}}}{\aderivation{}}\]}{\finv{\phi}}}\)\in
\banachmapdifclass{\infty}{\R^n}{\R^n}{\funcimage{\phi}}{\funcimage{\tanchart{\Man{}}{\phi}}},
\end{align}
which means,
\begin{equation}\label{lemmainvliederivativeissmoothpeq3}
\[\func{\LieDerinv{\Man{}}}{\aderivation{}}\]\in\mapdifclass{\infty}{\Man{}}{\Tanbun{\Man{}}}.
\end{equation}
\Ref{lemmainvliederivativeissmoothpeq1} and \Ref{lemmainvliederivativeissmoothpeq3} imply that
$\[\func{\LieDerinv{\Man{}}}{\aderivation{}}\]$ is a smooth vector-field on $\Man{}$.
\begin{equation}
\[\func{\LieDerinv{\Man{}}}{\aderivation{}}\]\in\vecf{\Man{}}{\infty}.
\end{equation}
\endlem
\definition\label{defliederinv}
$\Liederinv{\Man{}}$ is defined to be the codomain-restriction of $\LieDerinv{\Man{}}$ to
$\vecf{\Man{}}{\infty}$.
\begin{equation}
\Liederinv{\Man{}}:=\func{\rescd{\LieDerinv{\Man{}}}}{\vecf{\Man{}}{\infty}}.
\end{equation}
So $\Liederinv{\Man{}}$ is an element of $\Func{\Derivation{\Man{}}{\infty}}{\vecf{\Man{}}{\infty}}$.
When there is no ambiguity about the manifold $\Man{}$, $\Liederinv{}$ can replace $\Liederinv{\Man{}}$.
\endef
\theorem\label{thmliederivativeissurjective}
$\Liederinv{\Man{}}$ is a right-inverse to $\Lieder{\Man{}}$. That is,
\begin{equation}
\cmp{\Lieder{\Man{}}}{\Liederinv{\Man{}}}=\identity{\Derivation{\Man{}}{\infty}},
\end{equation}
or equivalently,
\begin{equation}
\Foreach{\aderivation{}}{\Derivation{\Man{}}{\infty}}
\func{\Lieder{\Man{}}}{\func{\Liederinv{\Man{}}}{\aderivation{}}}=\aderivation{}.
\end{equation}
\proof
$\aderivation{}$ is taken as an arbitrary element of $\Derivation{\Man{}}{\infty}$.
\begin{itemize}
\item[$\pr{1}$]
$\cf$ is taken as an arbitrary element of $\mapdifclass{\infty}{\Man{}}{\RR}$ and $\point$ as an arbitrary point of $\Man{}$.
$\phi$ is chosen as an element of $\defset{\p{\phi}}{\maxatlas{}}{\point\in\domain{\p{\phi}},~\func{\p{\phi}}{\point}=\zerovec{},~\funcimage{\phi}\in\convex{\R^n}}$, and
$\suc{\cf_1}{\cf_n}$ are chosen as elements of
$\suc{\fextension{\Man{}}{\point}{\phi}{\projection{n}{1}}}{\fextension{\Man{}}{\point}{\phi}{\projection{n}{n}}}$, respectively.
$\V$ is taken as such an open set of $\mantops{\Man{}}$ such that,
\begin{equation}\label{thmliederivativeissurjectivep1eq0}
\Foreach{k}{\seta{\suc{1}{n}}}
\[\point\in\V\subseteq\domain{\phi},~
\func{\resd{\cf_k}}{\V}=\func{\resd{\cmp{\projection{n}{k}}{\phi}}}{\V}\].
\end{equation}
According to \refdef{definvliederivative},
\begin{align}\label{thmliederivativeissurjectivep1eq1}
\func{\[\func{\Liederinv{\Man{}}}{\aderivation{}}\]}{\point}=
\sum_{k=1}^{n}\bigg(\func{\[\func{\aderivation{}}{\cf_k}\]}{\point}\bigg)
\bigg(\func{\finv{\[\tanspaceiso{\point}{\Man{}}{\phi}\]}}{\Eucbase{n}{k}}\bigg),
\end{align}
and according to \refdef{defliederivative},
\begin{align}\label{thmliederivativeissurjectivep1eq2}
\func{\(\func{\[\func{\Lieder{\Man{}}}{\func{\Liederinv{\Man{}}}{\aderivation{}}}\]}{\cf}\)}{\point}=
\func{\Rder{\cf}{\Man{}}}{\func{\[\func{\Liederinv{\Man{}}}{\aderivation{}}\]}{\point}}.
\end{align}
In addition, since $\func{\Liederinv{\Man{}}}{\aderivation{}}$ is a smooth vector-field on $\Man{}$,
$\func{\basep{\Man{}}}{\func{\[\func{\Liederinv{\Man{}}}{\aderivation{}}\]}{\point}}=\point$ and hence,
\begin{align}\label{thmliederivativeissurjectivep1eq3}
\func{\Rder{\cf}{\Man{}}}{\func{\[\func{\Liederinv{\Man{}}}{\aderivation{}}\]}{\point}}&=
\func{\[\cmp{\bigg(\func{\[\banachder{\(\cmp{f}{\finv{\phi}}\)}{\R^n}{\R}\]}
{\func{\phi}{\point}}\bigg)}{\tanspaceiso{\point}{\Man{}}{\phi}}\]}
{\func{\[\func{\Liederinv{\Man{}}}{\aderivation{}}\]}{\point}}.
\end{align}
\Ref{thmliederivativeissurjectivep1eq1}, \Ref{thmliederivativeissurjectivep1eq2}, and
\Ref{thmliederivativeissurjectivep1eq3} together with linearity of\\
$\cmp{\bigg(\func{\[\banachder{\(\cmp{f}{\finv{\phi}}\)}{\R^n}{\R}\]}
{\func{\phi}{\point}}\bigg)}{\tanspaceiso{\point}{\Man{}}{\phi}}$ imply,
\begin{equation}\label{thmliederivativeissurjectivep1eq4}
\func{\(\func{\[\func{\Lieder{\Man{}}}{\func{\Liederinv{\Man{}}}{\aderivation{}}}\]}{\cf}\)}{\point}=
\sum_{k=1}^{n}\bigg(\func{\[\func{\aderivation{}}{\cf_k}\]}{\point}\bigg)
\[\func{\bigg(\func{\[\banachder{\(\cmp{f}{\finv{\phi}}\)}{\R^n}{\R}\]}
{\func{\phi}{\point}}\bigg)}{\Eucbase{n}{k}}\].
\end{equation}
According to \refcor{corneighborhoodofderivativeofrealvaluedmaps},
\begin{align}\label{thmliederivativeissurjectivep1eq5}
&\Existsis{\mtuple{\lambda_1}{\lambda_n}}{{\banachmapdifclass{\infty}{\R^n}{\R}{\funcimage{\phi}}{\R}}^{\times n}}\cr
&\Foreach{k}{\seta{\suc{1}{n}}}\func{\(\cmp{\lambda_{k}}{\phi}\)}{\point}=
\func{\[\func{\(\banachder{\(\cmp{\cf}{\finv{\phi}}\)}{\R^n}{\R}\)}{\func{\phi}{\point}}\]}{\Eucbase{n}{k}},\cr
&\Foreach{\x}{\domain{\phi}}\func{\cf}{\x}=
\sum_{k=1}^{n}\[\func{\(\cmp{\projection{n}{k}}{\phi}\)}{\x}\]\[\func{\(\cmp{\lambda_{k}}{\phi}\)}{\x}\].
\end{align}
For every $k$ in $\seta{\suc{1}{n}}$, $\suc{F_1}{F_n}$ are taken as representative elements of\\
$\suc{\fextension{\Man{}}{\point}{\phi}{\lambda_1}}{\fextension{\Man{}}{\point}{\phi}{\lambda_n}}$, respectively, and
$\U$ is taken as such an element of $\mantops{\Man{}}$ such that,
\begin{equation}\label{thmliederivativeissurjectivep1eq6}
\Foreach{k}{\seta{\suc{1}{n}}}
\[\point\in\U\subseteq\domain{\phi},~
\func{\resd{F_k}}{\U}=\func{\resd{\cmp{\lambda_k}{\phi}}}{\U}\].
\end{equation}
Therefore by defining,
$O:=\V\cap\U$,
evidently $\point\in O\in\mantop{\Man{}}$ and,
\begin{align}\label{thmliederivativeissurjectivep1eq7}
\Foreach{j}{\seta{\suc{1}{n}}}
\Foreach{\x}{O}
\func{\cf}{\x}&=\sum_{k=1}^{n}\func{\cf_k}{\x}\func{F_k}{\x}\cr
&=\func{\(\sum_{k=1}^{n}\cf_{k}\rdot F_{k}\)}{\x},
\end{align}
and hence according to \reflem{lemmaderivationislocal}, \Ref{thmliederivativeissurjectivep1eq0},
\Ref{thmliederivativeissurjectivep1eq4}, \Ref{thmliederivativeissurjectivep1eq5},
\Ref{thmliederivativeissurjectivep1eq6},
properties of an $\infty$-derivation on $\Man{}$, considering that $\func{\phi}{\point}=\zerovec{}$,
\begin{align}\label{thmliederivativeissurjectivep1eq8}
\func{\[\func{\aderivation{}}{\cf}\]}{\point}&=
\func{\[\func{\aderivation{}}{\sum_{k=1}^{n}\cf_{k}\rdot F_{k}}\]}{\point}\cr
&=\sum_{k=1}^{n}\func{\[\func{\aderivation{}}{\cf_k}\]}{\point}\func{F_{k}}{\point}+
\func{\cf_k}{\point}\func{\[\func{\aderivation{}}{F_k}\]}{\point}\cr
&=\sum_{k=1}^{n}\bigg(\func{\[\func{\aderivation{}}{\cf_k}\]}{\point}\bigg)
\[\func{\bigg(\func{\[\banachder{\(\cmp{f}{\finv{\phi}}\)}{\R^n}{\R}\]}
{\func{\phi}{\point}}\bigg)}{\Eucbase{n}{k}}\]\cr
&=\func{\(\func{\[\func{\Lieder{\Man{}}}{\func{\Liederinv{\Man{}}}{\aderivation{}}}\]}{\cf}\)}{\point}.
\end{align}
\endp
\end{itemize}
Therefore,
\begin{align}
\Foreach{\cf}{\mapdifclass{\infty}{\Man{}}{\RR}}
\Foreach{\point}{\Man{}}
\func{\(\func{\[\func{\Lieder{\Man{}}}{\func{\Liederinv{\Man{}}}{\aderivation{}}}\]}{\cf}\)}{\point}=
\func{\[\func{\aderivation{}}{\cf}\]}{\point},
\end{align}
and hence,
\begin{equation}
\func{\Lieder{\Man{}}}{\func{\Liederinv{\Man{}}}{\aderivation{}}}=\aderivation{}.
\end{equation}
\endthm
\corollary\label{corliederivativeissurjective1}
$\Lieder{\Man{}}$ is an surjective map. That is,
\begin{equation}
\Lieder{\Man{}}\in\surFunc{\vecf{\Man{}}{\infty}}{\Derivation{\Man{}}{\infty}}.
\end{equation}
\endcor
\corollary\label{corliederivativeisalinearisomorphism}
$\Lieder{\Man{}}$ is a is a linear-isomorphism from $\Vecf{\Man{}}{\infty}$ to $\LDerivation{\Man{}}{\infty}$. That is,
it is a bijective and linear map from $\vecf{\Man{}}{\infty}$ to $\Derivation{\Man{}}{\infty}$ endowd with their canonical linear-structures.
\begin{equation}
\Lieder{\Man{}}\in\IF{\vecf{\Man{}}{\infty}}{\Derivation{\Man{}}{\infty}}
\cap\Lin{\Vecf{\Man{}}{\infty}}{\LDerivation{\Man{}}{\infty}}.
\end{equation}
Additionally,
\begin{equation}
\finv{\(\Lieder{\Man{}}\)}=\Liederinv{\Man{}}.
\end{equation}
\proof
It is obvious according to \reflem{lemmaliederivativeislinear}, \refthm{thmliederivativeisinjective},
\refcor{corliederivativeissurjective1}, and \refthm{thmliederivativeissurjective}.
\endcor
\section{Lie-Bracket}
\theorem\label{thmpreliebracket}
For every pair of smooth vector-fields $\avecf{1}$ and $\avecf{2}$ on $\Man{}$,
$\cmp{\func{\Lieder{\Man{}}}{\avecf{1}}}{\func{\Lieder{\Man{}}}{\avecf{2}}}-
\cmp{\func{\Lieder{\Man{}}}{\avecf{2}}}{\func{\Lieder{\Man{}}}{\avecf{1}}}$ is an $\infty$-derivation on $\Man{}$. That is,
\begin{equation}
\Foreach{\opair{\avecf{1}}{\avecf{2}}}{\Cprod{\vecf{\Man{}}{\infty}}{\vecf{\Man{}}{\infty}}}
\[\cmp{\func{\Lieder{\Man{}}}{\avecf{1}}}{\func{\Lieder{\Man{}}}{\avecf{2}}}-
\cmp{\func{\Lieder{\Man{}}}{\avecf{2}}}{\func{\Lieder{\Man{}}}{\avecf{1}}}\]\in
\Derivation{\Man{}}{\infty}.
\end{equation}
Here, the canonical linear-structure of the space of all linear maps from $\Lmapdifclass{\infty}{\Man{}}{\RR}$ to
$\Lmapdifclass{\infty}{\Man{}}{\RR}$ is invoked when subtracting 
$\cmp{\func{\Lieder{\Man{}}}{\avecf{2}}}{\func{\Lieder{\Man{}}}{\avecf{1}}}$ from
$\cmp{\func{\Lieder{\Man{}}}{\avecf{1}}}{\func{\Lieder{\Man{}}}{\avecf{2}}}$.
\proof
Each $\avecf{1}$ and $\avecf{2}$ is taken as an arbitrary element of $\vecf{\Man{}}{\infty}$. So,
\begin{equation}\label{thmpreliebracketpeq0}
\cmp{\basep{\Man{}}}{\avecf{j}}=\identity{\M{}},~j\in\seta{\binary{1}{2}}.
\end{equation}
Since each $\func{\Lieder{\Man{}}}{\avecf{1}}$ and $\func{\Lieder{\Man{}}}{\avecf{2}}$ is an $\infty$-derivation on $\Man{}$, it is a
linear map from $\Lmapdifclass{\infty}{\Man{}}{\RR}$ to $\Lmapdifclass{\infty}{\Man{}}{\RR}$, and hence
so is $\[\cmp{\func{\Lieder{\Man{}}}{\avecf{1}}}{\func{\Lieder{\Man{}}}{\avecf{2}}}-
\cmp{\func{\Lieder{\Man{}}}{\avecf{2}}}{\func{\Lieder{\Man{}}}{\avecf{1}}}\]$. That is,
\begin{equation}\label{thmpreliebracketpeq1}
\[\cmp{\func{\Lieder{\Man{}}}{\avecf{1}}}{\func{\Lieder{\Man{}}}{\avecf{2}}}-
\cmp{\func{\Lieder{\Man{}}}{\avecf{2}}}{\func{\Lieder{\Man{}}}{\avecf{1}}}\]\in
\Lin{\Lmapdifclass{\infty}{\Man{}}{\RR}}{\Lmapdifclass{\infty}{\Man{}}{\RR}}.
\end{equation}
\begin{itemize}
\item[$\pr{1}$]
Each $\cf$ and $\cg$ is taken as an arbitrary element of $\mapdifclass{\infty}{\Man{}}{\RR}$ and $\point$ as an arbitrary
point of $\Man{}$.
Considering the addition operation in the linear-structure of the space of all linear maps from $\Lmapdifclass{\infty}{\Man{}}{\RR}$ to
$\Lmapdifclass{\infty}{\Man{}}{\RR}$,
\begin{align}\label{thmpreliebracketp1eq1}
&~~~\func{\[\cmp{\func{\Lieder{\Man{}}}{\avecf{1}}}{\func{\Lieder{\Man{}}}{\avecf{2}}}-
\cmp{\func{\Lieder{\Man{}}}{\avecf{2}}}{\func{\Lieder{\Man{}}}{\avecf{1}}}\]}{\cf\rdot\cg}\cr
&=\func{\[\cmp{\func{\Lieder{\Man{}}}{\avecf{1}}}{\func{\Lieder{\Man{}}}{\avecf{2}}}\]}{\cf\rdot\cg}-
\func{\[\cmp{\func{\Lieder{\Man{}}}{\avecf{2}}}{\func{\Lieder{\Man{}}}{\avecf{1}}}\]}{\cf\rdot\cg}\cr
&=\func{\[\func{\Lieder{\Man{}}}{\avecf{1}}\]}{\func{\[\func{\Lieder{\Man{}}}{\avecf{2}}\]}{\cf\rdot\cg}}-
\func{\[\func{\Lieder{\Man{}}}{\avecf{2}}\]}{\func{\[\func{\Lieder{\Man{}}}{\avecf{1}}\]}{\cf\rdot\cg}}.
\end{align}
Since each $\func{\Lieder{\Man{}}}{\avecf{1}}$ and $\func{\Lieder{\Man{}}}{\avecf{2}}$ is an $\infty$-derivation on $\Man{}$,
\begin{align}\label{thmpreliebracketp1eq2}
&~~~\func{\[\func{\Lieder{\Man{}}}{\avecf{1}}\]}{\func{\[\func{\Lieder{\Man{}}}{\avecf{2}}\]}{\cf\rdot\cg}}\cr
&=\func{\[\func{\Lieder{\Man{}}}{\avecf{1}}\]}
{\bigg[\func{\(\func{\Lieder{\Man{}}}{\avecf{2}}\)}{\cf}\bigg]\rdot\cg+
\cf\rdot\bigg[\func{\(\func{\Lieder{\Man{}}}{\avecf{2}}\)}{\cg}\bigg]}\cr
&=\func{\[\func{\Lieder{\Man{}}}{\avecf{1}}\]}
{\bigg[\func{\(\func{\Lieder{\Man{}}}{\avecf{2}}\)}{\cf}\bigg]\rdot\cg}+
\func{\[\func{\Lieder{\Man{}}}{\avecf{1}}\]}
{\cf\rdot\bigg[\func{\(\func{\Lieder{\Man{}}}{\avecf{2}}\)}{\cg}\bigg]}\cr
&=~~\bigg(\func{\[\func{\Lieder{\Man{}}}{\avecf{1}}\]}{\func{\[\func{\Lieder{\Man{}}}{\avecf{2}}\]}{\cf}}\bigg)\rdot\cg+
\bigg[\func{\(\func{\Lieder{\Man{}}}{\avecf{2}}\)}{\cf}\bigg]\rdot\bigg[\func{\(\func{\Lieder{\Man{}}}{\avecf{1}}\)}{\cg}\bigg]\cr
&~~~+\bigg[\func{\(\func{\Lieder{\Man{}}}{\avecf{1}}\)}{\cf}\bigg]\rdot\bigg[\func{\(\func{\Lieder{\Man{}}}{\avecf{2}}\)}{\cg}\bigg]+
\cf\rdot\bigg(\func{\[\func{\Lieder{\Man{}}}{\avecf{1}}\]}{\func{\[\func{\Lieder{\Man{}}}{\avecf{2}}\]}{\cg}}\bigg).~~~~~~
\end{align}
Similarly,
\begin{align}\label{thmpreliebracketp1eq3}
&~~~\func{\[\func{\Lieder{\Man{}}}{\avecf{2}}\]}{\func{\[\func{\Lieder{\Man{}}}{\avecf{1}}\]}{\cf\rdot\cg}}\cr
&=~~\bigg(\func{\[\func{\Lieder{\Man{}}}{\avecf{2}}\]}{\func{\[\func{\Lieder{\Man{}}}{\avecf{1}}\]}{\cf}}\bigg)\rdot\cg+
\bigg[\func{\(\func{\Lieder{\Man{}}}{\avecf{1}}\)}{\cf}\bigg]\rdot\bigg[\func{\(\func{\Lieder{\Man{}}}{\avecf{2}}\)}{\cg}\bigg]\cr
&~~~+\bigg[\func{\(\func{\Lieder{\Man{}}}{\avecf{2}}\)}{\cf}\bigg]\rdot\bigg[\func{\(\func{\Lieder{\Man{}}}{\avecf{1}}\)}{\cg}\bigg]+
\cf\rdot\bigg(\func{\[\func{\Lieder{\Man{}}}{\avecf{2}}\]}{\func{\[\func{\Lieder{\Man{}}}{\avecf{1}}\]}{\cg}}\bigg).~~~~~~
\end{align}
Considering the addition operation in the linear-structure of the space of all linear maps from $\Lmapdifclass{\infty}{\Man{}}{\RR}$ to
$\Lmapdifclass{\infty}{\Man{}}{\RR}$ and invoking the $\R$-algebra structure of $\Lmapdifclass{\infty}{\Man{}}{\RR}$,
\Ref{thmpreliebracketp1eq1}, \Ref{thmpreliebracketp1eq2}, and \Ref{thmpreliebracketp1eq3} imply,
\begin{align}\label{thmpreliebracketp1eq4}
&~~~\func{\[\cmp{\func{\Lieder{\Man{}}}{\avecf{1}}}{\func{\Lieder{\Man{}}}{\avecf{2}}}-
\cmp{\func{\Lieder{\Man{}}}{\avecf{2}}}{\func{\Lieder{\Man{}}}{\avecf{1}}}\]}{\cf\rdot\cg}\cr
&=~~\bigg(\func{\[\cmp{\func{\Lieder{\Man{}}}{\avecf{1}}}{\func{\Lieder{\Man{}}}{\avecf{2}}}-
\cmp{\func{\Lieder{\Man{}}}{\avecf{2}}}{\func{\Lieder{\Man{}}}{\avecf{1}}}\]}{\cf}\bigg)\rdot\cg\cr
&~~~+
\cf\rdot\bigg(\func{\[\cmp{\func{\Lieder{\Man{}}}{\avecf{1}}}{\func{\Lieder{\Man{}}}{\avecf{2}}}-
\cmp{\func{\Lieder{\Man{}}}{\avecf{2}}}{\func{\Lieder{\Man{}}}{\avecf{1}}}\]}{\cg}\bigg).
\end{align}
\endp
\end{itemize}
Therefore,
\begin{align}\label{thmpreliebracketpeq2}
\Foreach{\opair{\cf}{\cg}}{{\mapdifclass{\infty}{\Man{}}{\RR}}^{\times 2}}
&~~~\func{\[\cmp{\func{\Lieder{\Man{}}}{\avecf{1}}}{\func{\Lieder{\Man{}}}{\avecf{2}}}-
\cmp{\func{\Lieder{\Man{}}}{\avecf{2}}}{\func{\Lieder{\Man{}}}{\avecf{1}}}\]}{\cf\rdot\cg}\cr
&=~~\bigg(\func{\[\cmp{\func{\Lieder{\Man{}}}{\avecf{1}}}{\func{\Lieder{\Man{}}}{\avecf{2}}}-
\cmp{\func{\Lieder{\Man{}}}{\avecf{2}}}{\func{\Lieder{\Man{}}}{\avecf{1}}}\]}{\cf}\bigg)\rdot\cg\cr
&~~~+
\cf\rdot\bigg(\func{\[\cmp{\func{\Lieder{\Man{}}}{\avecf{1}}}{\func{\Lieder{\Man{}}}{\avecf{2}}}-
\cmp{\func{\Lieder{\Man{}}}{\avecf{2}}}{\func{\Lieder{\Man{}}}{\avecf{1}}}\]}{\cg}\bigg).~~~~~~~
\end{align}
According to the definition of an $\infty$-derivation on $\Man{}$, \Ref{thmpreliebracketpeq1} and \Ref{thmpreliebracketpeq2}
equivalently mean that,
\begin{equation}
\[\cmp{\func{\Lieder{\Man{}}}{\avecf{1}}}{\func{\Lieder{\Man{}}}{\avecf{2}}}-
\cmp{\func{\Lieder{\Man{}}}{\avecf{2}}}{\func{\Lieder{\Man{}}}{\avecf{1}}}\]\in
\Derivation{\Man{}}{\infty}.
\end{equation}
\endthm
\definition\label{defliebracket}
The binary operation $\Liebracket{\Man{}}$ on $\vecf{\Man{}}{\infty}$ is defined as,
\begin{align}
&\Liebracket{\Man{}}\indef\Func{\Cprod{\vecf{\Man{}}{\infty}}{\vecf{\Man{}}{\infty}}}{\vecf{\Man{}}{\infty}},\cr
&\Foreach{\opair{\avecf{1}}{\avecf{2}}}{{\vecf{\Man{}}{\infty}}^{\times 2}}\cr
&\func{\Liebracket{\Man{}}}{\binary{\avecf{1}}{\avecf{2}}}\eqdef\func{\finv{\(\Lieder{\Man{}}\)}}
{\cmp{\func{\Lieder{\Man{}}}{\avecf{1}}}{\func{\Lieder{\Man{}}}{\avecf{2}}}-
\cmp{\func{\Lieder{\Man{}}}{\avecf{2}}}{\func{\Lieder{\Man{}}}{\avecf{1}}}}.\cr
&{}
\end{align}
$\Liebracket{\Man{}}$ is called the $\quotl$Lie-bracket operation on $\vecf{\Man{}}{\infty}$$\quotr$, and for every pair
$\avecf{1}$ and $\avecf{2}$ of smooth vector-fields on $\Man{}$, $\func{\Liebracket{\Man{}}}{\binary{\avecf{1}}{\avecf{2}}}$
is referred to as the $\quotl$Lie-bracket (with respect to the manifold $\Man{}$) of $\avecf{1}$ and $\avecf{2}$$\quotr$.
$\liebracket{\avecf{1}}{\avecf{2}}{\Man{}}$ conventionally denotes $\func{\Liebracket{\Man{}}}{\binary{\avecf{1}}{\avecf{2}}}$.
Also, when there is no ambiguity about the underlying manifold $\Man{}$, $\liebracket{\avecf{1}}{\avecf{2}}{}$
can replace $\liebracket{\avecf{1}}{\avecf{2}}{\Man{}}$. When $\liebracket{\avecf{1}}{\avecf{2}}{\Man{}}$ vanishes, that is
when it equals the zero element $\zerovec{\Vecf{\Man{}}{\infty}}$ of $\Vecf{\Man{}}{\infty}$, then it is said that
$\quotl$the pair of smooth vector-fields $\avecf{1}$ and $\avecf{2}$ on $\Man{}$ commute$\quotr$.
\endef
\theorem\label{thmliealgebraofvectorfields0}
The set of all smooth vector-fields on $\Man{}$ endowed with its canonical linear-structure together with the
Lie-bracket operation on $\vecf{\Man{}}{\infty}$ is a Lie-algebra. That is, the pair
$\opair{\Vecf{\Man{}}{\infty}}{\Liebracket{\Man{}}}$ is a Lie-algebra.
\proof
\begin{itemize}
\item[$\pr{1}$]
Each $\avecf{1}$, $\avecf{2}$, and $\avecf{3}$ is taken as an arbitrary element of $\vecf{\Man{}}{\infty}$, and $\c$ as
an arbitrary element of $\R$. According to \refdef{defliebracket} and \refcor{corliederivativeisalinearisomorphism}
($\Lieder{\Man{}}$ is a linear isomorphism),
\begin{align}
&~~~\liebracket{\c\avecf{1}+\avecf{2}}{\avecf{3}}{\Man{}}\cr
&=\func{\finv{\(\Lieder{\Man{}}\)}}
{\cmp{\func{\Lieder{\Man{}}}{\c\avecf{1}+\avecf{2}}}{\func{\Lieder{\Man{}}}{\avecf{3}}}-
\cmp{\func{\Lieder{\Man{}}}{\avecf{3}}}{\func{\Lieder{\Man{}}}{\c\avecf{1}+\avecf{2}}}}\cr
&=\func{\finv{\(\Lieder{\Man{}}\)}}
{\cmp{\[\c\func{\Lieder{\Man{}}}{\avecf{1}}+\func{\Lieder{\Man{}}}{\avecf{2}}\]}{\func{\Lieder{\Man{}}}{\avecf{3}}}-
\cmp{\func{\Lieder{\Man{}}}{\avecf{3}}}{\[\c\func{\Lieder{\Man{}}}{\avecf{1}}+\func{\Lieder{\Man{}}}{\avecf{2}}\]}}\cr
&=~~\c\bigg[\func{\finv{\(\Lieder{\Man{}}\)}}
{\[\cmp{\func{\Lieder{\Man{}}}{\avecf{1}}}{\func{\Lieder{\Man{}}}{\avecf{3}}}-
\cmp{\func{\Lieder{\Man{}}}{\avecf{3}}}{\func{\Lieder{\Man{}}}{\avecf{1}}}\]}\bigg]\cr
&~~~+\bigg[\func{\finv{\(\Lieder{\Man{}}\)}}
{\[\cmp{\func{\Lieder{\Man{}}}{\avecf{2}}}{\func{\Lieder{\Man{}}}{\avecf{3}}}-
\cmp{\func{\Lieder{\Man{}}}{\avecf{3}}}{\func{\Lieder{\Man{}}}{\avecf{2}}}\]}\bigg]\cr
&=\c\liebracket{\avecf{1}}{\avecf{3}}{\Man{}}+\liebracket{\avecf{2}}{\avecf{3}}{\Man{}}.
\end{align}
Similarly,
\begin{equation}
\liebracket{\avecf{3}}{\c\avecf{1}+\avecf{2}}{\Man{}}=
\c\liebracket{\avecf{3}}{\avecf{1}}{\Man{}}+\liebracket{\avecf{3}}{\avecf{2}}{\Man{}}.
\end{equation}
\endp
\end{itemize}
\begin{itemize}
\item[$\pr{2}$]
$\avecf{}$ is taken as an arbitrary element of $\vecf{\Man{}}{\infty}$.
According to \refdef{defliebracket} and \refcor{corliederivativeisalinearisomorphism}
($\Lieder{\Man{}}$ is a linear isomorphism),
\begin{align}
\liebracket{\avecf{}}{\avecf{}}{\Man{}}&=
\func{\finv{\(\Lieder{\Man{}}\)}}
{\cmp{\func{\Lieder{\Man{}}}{\avecf{}}}{\func{\Lieder{\Man{}}}{\avecf{}}}-
\cmp{\func{\Lieder{\Man{}}}{\avecf{}}}{\func{\Lieder{\Man{}}}{\avecf{}}}}\cr
&=\func{\finv{\(\Lieder{\Man{}}\)}}{\zerovec{\LDerivation{\Man{}}{\infty}}}\cr
&=\zerovec{\Vecf{\Man{}}{\infty}},
\end{align}
where $\zerovec{\LDerivation{\Man{}}{\infty}}$ and $\zerovec{\Vecf{\Man{}}{\infty}}$ denote the neutral elements
of addition operations of the vector-spaces $\LDerivation{\Man{}}{\infty}$ and $\Vecf{\Man{}}{\infty}$, respectively.
\endp
\end{itemize}
\begin{itemize}
\item[$\pr{3}$]
Each $\avecf{1}$, $\avecf{2}$, and $\avecf{3}$ is taken as an arbitrary element of $\vecf{\Man{}}{\infty}$.
According to \refdef{defliebracket},
\begin{align}
\liebracket{\avecf{1}}{\liebracket{\avecf{2}}{\avecf{3}}{}}{}=\func{\finv{\(\Lieder{}\)}}
{\cmp{\func{\Lieder{}}{\avecf{1}}}{\func{\Lieder{}}{\liebracket{\avecf{2}}{\avecf{3}}{}}}-
\cmp{\func{\Lieder{}}{\liebracket{\avecf{2}}{\avecf{3}}{}}}{\func{\Lieder{}}{\avecf{1}}}},
\end{align}
and,
\begin{align}
\func{\Lieder{}}{\liebracket{\avecf{2}}{\avecf{3}}{}}&=
\func{\Lieder{}}{\func{\finv{\(\Lieder{}\)}}
{\cmp{\func{\Lieder{}}{\avecf{2}}}{\func{\Lieder{}}{\avecf{3}}}-
\cmp{\func{\Lieder{}}{\avecf{3}}}{\func{\Lieder{}}{\avecf{2}}}}}\cr
&=\cmp{\func{\Lieder{}}{\avecf{2}}}{\func{\Lieder{}}{\avecf{3}}}-
\cmp{\func{\Lieder{}}{\avecf{3}}}{\func{\Lieder{}}{\avecf{2}}}.
\end{align}
Thus,
\begin{align}
&~~~\func{\Lieder{}}{\liebracket{\avecf{1}}{\liebracket{\avecf{2}}{\avecf{3}}{}}{}}\cr
&=~~\cmp{\func{\Lieder{}}{\avecf{1}}}{\cmp{\func{\Lieder{}}{\avecf{2}}}{\func{\Lieder{}}{\avecf{3}}}}-
\cmp{\func{\Lieder{}}{\avecf{1}}}{\cmp{\func{\Lieder{}}{\avecf{3}}}{\func{\Lieder{}}{\avecf{2}}}}\cr
&~~~-\cmp{\func{\Lieder{}}{\avecf{2}}}{\cmp{\func{\Lieder{}}{\avecf{3}}}{\func{\Lieder{}}{\avecf{1}}}}+
\cmp{\func{\Lieder{}}{\avecf{3}}}{\cmp{\func{\Lieder{}}{\avecf{2}}}{\func{\Lieder{}}{\avecf{1}}}}.
\end{align}
Similarly,
\begin{align}
&~~~\func{\Lieder{}}{\liebracket{\avecf{2}}{\liebracket{\avecf{3}}{\avecf{1}}{}}{}}\cr
&=~~\cmp{\func{\Lieder{}}{\avecf{2}}}{\cmp{\func{\Lieder{}}{\avecf{3}}}{\func{\Lieder{}}{\avecf{1}}}}-
\cmp{\func{\Lieder{}}{\avecf{2}}}{\cmp{\func{\Lieder{}}{\avecf{1}}}{\func{\Lieder{}}{\avecf{3}}}}\cr
&~~~-\cmp{\func{\Lieder{}}{\avecf{3}}}{\cmp{\func{\Lieder{}}{\avecf{1}}}{\func{\Lieder{}}{\avecf{2}}}}+
\cmp{\func{\Lieder{}}{\avecf{1}}}{\cmp{\func{\Lieder{}}{\avecf{3}}}{\func{\Lieder{}}{\avecf{2}}}},
\end{align}
and,
\begin{align}
&~~~\func{\Lieder{}}{\liebracket{\avecf{3}}{\liebracket{\avecf{1}}{\avecf{2}}{}}{}}\cr
&=~~\cmp{\func{\Lieder{}}{\avecf{3}}}{\cmp{\func{\Lieder{}}{\avecf{1}}}{\func{\Lieder{}}{\avecf{2}}}}-
\cmp{\func{\Lieder{}}{\avecf{3}}}{\cmp{\func{\Lieder{}}{\avecf{2}}}{\func{\Lieder{}}{\avecf{1}}}}\cr
&~~~-\cmp{\func{\Lieder{}}{\avecf{1}}}{\cmp{\func{\Lieder{}}{\avecf{2}}}{\func{\Lieder{}}{\avecf{3}}}}+
\cmp{\func{\Lieder{}}{\avecf{2}}}{\cmp{\func{\Lieder{}}{\avecf{1}}}{\func{\Lieder{}}{\avecf{3}}}},
\end{align}
and hence,
\begin{align}
&~~~\func{\Lieder{}}{\liebracket{\avecf{1}}{\liebracket{\avecf{2}}{\avecf{3}}{}}{}+
\liebracket{\avecf{2}}{\liebracket{\avecf{3}}{\avecf{1}}{}}{}+
\liebracket{\avecf{3}}{\liebracket{\avecf{1}}{\avecf{2}}{}}{}}\cr
&=\func{\Lieder{}}{\liebracket{\avecf{1}}{\liebracket{\avecf{2}}{\avecf{3}}{}}{}}+
\func{\Lieder{}}{\liebracket{\avecf{2}}{\liebracket{\avecf{3}}{\avecf{1}}{}}{}}+
\func{\Lieder{}}{\liebracket{\avecf{3}}{\liebracket{\avecf{1}}{\avecf{2}}{}}{}}\cr
&=\zerovec{\LDerivation{\Man{}}{\infty}}.
\end{align}
Thus according to \refcor{corliederivativeisalinearisomorphism},
\begin{align}
\liebracket{\avecf{1}}{\liebracket{\avecf{2}}{\avecf{3}}{}}{}+
\liebracket{\avecf{2}}{\liebracket{\avecf{3}}{\avecf{1}}{}}{}+
\liebracket{\avecf{3}}{\liebracket{\avecf{1}}{\avecf{2}}{}}{}=
\zerovec{\Vecf{\Man{}}{\infty}}.
\end{align}
\endp
\end{itemize}
\endthm
\definition\label{defliealgebraofvectorfields}
The Lie-algebra $\opair{\Vecf{\Man{}}{\infty}}{\Liebracket{\Man{}}}$ is referred to as the
$\quotl$Lie-algebra of smooth vector-fields on the manifold $\Man{}$$\quotr$, and is denoted by
$\Lievecf{\Man{}}{\infty}$.
\endef
\theorem\label{thmliebracketofrelatedvectorfields}
$\asmooth{}$ is taken as an element of $\mapdifclass{\infty}{\Man{}}{\Man{1}}$ (a smooth map from $\Man{}$ to $\Man{1}$).
Each $\avecf{1}$ and $\avecf{2}$ is taken as an element of $\vecf{\Man{}}{\infty}$, and each
$\avecff{1}$ and $\avecff{2}$ as an element of $\vecf{\Man{1}}{\infty}$. If
\begin{align}
\cmp{\avecff{1}}{\asmooth{}}=\cmp{\[\der{\asmooth{}}{\Man{}}{\Man{1}}\]}{\avecf{1}},\cr
\cmp{\avecff{2}}{\asmooth{}}=\cmp{\[\der{\asmooth{}}{\Man{}}{\Man{1}}\]}{\avecf{2}},
\end{align}
then,
\begin{align}
\cmp{\liebracket{\avecff{1}}{\avecff{2}}{\Man{1}}}{\asmooth{}}=
\cmp{\[\der{\asmooth{}}{\Man{}}{\Man{1}}\]}{\liebracket{\avecf{1}}{\avecf{2}}{\Man{}}}.
\end{align}
\proof
It is assumed that,
\begin{align}
\cmp{\avecff{1}}{\asmooth{}}=\cmp{\[\der{\asmooth{}}{\Man{}}{\Man{1}}\]}{\avecf{1}},
\label{thmliebracketofrelatedvectorfieldspeq1}\\
\cmp{\avecff{2}}{\asmooth{}}=\cmp{\[\der{\asmooth{}}{\Man{}}{\Man{1}}\]}{\avecf{2}}.
\label{thmliebracketofrelatedvectorfieldspeq2}
\end{align}
$\point$ is taken as an arbitrary point of the manifold $\Man{}$.
According to \refdef{defliederivative} and \refdef{defliederivative00},
\begin{equation}\label{thmliebracketofrelatedvectorfieldspeq3}
\Foreach{\cf}{\mapdifclass{\infty}{\Man{1}}{\RR}}
\func{\[\func{\Lieder{\Man{1}}}{\liebracket{\avecff{1}}{\avecff{2}}{\Man{1}}}\]}{\cf}=
\cmp{\(\Rder{\cf}{\Man{1}}\)}{\liebracket{\avecff{1}}{\avecff{2}}{\Man{1}}}.
\end{equation}
In addition, according to \refdef{defliederivative}, \refdef{defliederivative00}, and \refdef{defliebracket},
\begin{align}\label{thmliebracketofrelatedvectorfieldspeq4}
\Foreach{\cf}{\mapdifclass{\infty}{\Man{1}}{\RR}}
&\hskip 0.5\baselineskip\func{\[\func{\Lieder{\Man{1}}}
{\liebracket{\avecff{1}}{\avecff{2}}{\Man{1}}}\]}{\cf}\cr
&=\func{\[\cmp{\func{\Lieder{\Man{1}}}{\avecff{1}}}{\func{\Lieder{\Man{1}}}{\avecff{2}}}-
\cmp{\func{\Lieder{\Man{1}}}{\avecff{2}}}{\func{\Lieder{\Man{1}}}{\avecff{1}}}\]}{\cf}\cr
&=\func{\[\func{\Lieder{\Man{1}}}{\avecff{1}}\]}{\func{\[\func{\Lieder{\Man{1}}}{\avecff{2}}\]}{\cf}}-
\func{\[\func{\Lieder{\Man{1}}}{\avecff{2}}\]}{\func{\[\func{\Lieder{\Man{1}}}{\avecff{2}}\]}{\cf}}\cr
&=\func{\[\func{\Lieder{\Man{1}}}{\avecff{1}}\]}{\cmp{\Rder{\cf}{\Man{1}}}{\avecff{2}}}-
\func{\[\func{\Lieder{\Man{1}}}{\avecff{2}}\]}{\cmp{\Rder{\cf}{\Man{1}}}{\avecff{1}}}\cr
&=\cmp{\Rder{\(\cmp{\Rder{\cf}{\Man{1}}}{\avecff{2}}\)}{\Man{1}}}{\avecff{1}}-
\cmp{\Rder{\(\cmp{\Rder{\cf}{\Man{1}}}{\avecff{1}}\)}{\Man{1}}}{\avecff{2}}.
\end{align}
Additionally, it is clear that, according to the composition rule of smooth maps,
\begin{equation}\label{thmliebracketofrelatedvectorfieldspeq5}
\Foreach{\cf}{\mapdifclass{\infty}{\Man{1}}{\RR}}
\cmp{\cf}{\asmooth{}}\in\mapdifclass{\infty}{\Man{}}{\RR},
\end{equation}
and hence according to \refdef{defliederivative} and \refdef{defliederivative00},
and \Ref{eqspecialchainrule1} (the special chain-rule of differentiation),
\begin{align}\label{thmliebracketofrelatedvectorfieldspeq6}
\Foreach{\cf}{\mapdifclass{\infty}{\Man{1}}{\RR}}
\func{\[\func{\Lieder{\Man{}}}{\liebracket{\avecf{1}}{\avecf{2}}{\Man{}}}\]}{\cmp{\cf}{\asmooth{}}}&=
\cmp{\[\Rder{\(\cmp{\cf}{\asmooth{}}\)}{\Man{}}\]}{\liebracket{\avecf{1}}{\avecf{2}}{\Man{}}}\cr
&=\cmp{\(\Rder{\cf}{\Man{1}}\)}{\cmp{\(\der{\asmooth{}}{\Man{}}{\Man{1}}\)}{\liebracket{\avecf{1}}{\avecf{2}}{\Man{}}}}.
\end{align}
Also, according to \refdef{defliederivative}, \refdef{defliederivative00}, and \refdef{defliebracket},
passing through a similar reasoning as that of \Ref{thmliebracketofrelatedvectorfieldspeq4},
\begin{align}\label{thmliebracketofrelatedvectorfieldspeq7}
\Foreach{\cf}{\mapdifclass{\infty}{\Man{1}}{\RR}}
&\hskip 0.5\baselineskip\func{\[\func{\Lieder{\Man{}}}
{\liebracket{\avecf{1}}{\avecf{2}}{\Man{}}}\]}{\cmp{\cf}{\asmooth{}}}\cr
&=\cmp{\Rder{\[\cmp{\Rder{\(\cmp{\cf}{\asmooth{}}\)}{\Man{}}}{\avecf{2}}\]}{\Man{}}}{\avecf{1}}-
\cmp{\Rder{\[\cmp{\Rder{\(\cmp{\cf}{\asmooth{}}\)}{\Man{}}}{\avecf{1}}\]}{\Man{}}}{\avecf{2}}.
\end{align}
According to \Ref{eqspecialchainrule1} (the special chain-rule of differentiation),
\Ref{thmliebracketofrelatedvectorfieldspeq1}, \Ref{thmliebracketofrelatedvectorfieldspeq2},
\Ref{thmliebracketofrelatedvectorfieldspeq3}, \Ref{thmliebracketofrelatedvectorfieldspeq4},
\Ref{thmliebracketofrelatedvectorfieldspeq6}, abd \Ref{thmliebracketofrelatedvectorfieldspeq7} yield,
\begin{align}\label{thmliebracketofrelatedvectorfieldspeq8}
&\begin{aligned}
\Foreach{\cf}{\mapdifclass{\infty}{\Man{1}}{\RR}}
\end{aligned}\cr
&\begin{aligned}
&\hskip 0.5\baselineskip\cmp{\(\Rder{\cf}{\Man{1}}\)}{\[\cmp{\liebracket{\avecff{1}}{\avecff{2}}{\Man{1}}}{\asmooth{}}\]}\cr
&=\cmp{\(\func{\[\func{\Lieder{\Man{1}}}{\liebracket{\avecff{1}}{\avecff{2}}{\Man{1}}}\]}{\cf}\)}{\asmooth{}}\cr
&=\cmp{\[\cmp{\Rder{\(\cmp{\Rder{\cf}{\Man{1}}}{\avecff{2}}\)}{\Man{1}}}{\avecff{1}}-
\cmp{\Rder{\(\cmp{\Rder{\cf}{\Man{1}}}{\avecff{1}}\)}{\Man{1}}}{\avecff{2}}\]}{\asmooth{}}\cr
&=\cmp{\[\cmp{\Rder{\(\cmp{\Rder{\cf}{\Man{1}}}{\avecff{2}}\)}{\Man{1}}}{\avecff{1}}\]}{\asmooth{}}-
\cmp{\[\cmp{\Rder{\(\cmp{\Rder{\cf}{\Man{1}}}{\avecff{1}}\)}{\Man{1}}}{\avecff{2}}\]}{\asmooth{}}\cr
&=\cmp{\Rder{\(\cmp{\Rder{\cf}{\Man{1}}}{\avecff{2}}\)}{\Man{1}}}{\(\cmp{\avecff{1}}{\asmooth{}}\)}-
\cmp{\Rder{\(\cmp{\Rder{\cf}{\Man{1}}}{\avecff{1}}\)}{\Man{1}}}{\(\cmp{\avecff{2}}{\asmooth{}}\)}\cr
&=\cmp{\Rder{\(\cmp{\Rder{\cf}{\Man{1}}}{\avecff{2}}\)}{\Man{1}}}{\(\cmp{\[\der{\asmooth{}}{\Man{}}{\Man{1}}\]}{\avecf{1}}\)}-
\cmp{\Rder{\(\cmp{\Rder{\cf}{\Man{1}}}{\avecff{1}}\)}{\Man{1}}}{\(\cmp{\[\der{\asmooth{}}{\Man{}}{\Man{1}}\]}{\avecf{2}}\)}\cr
&=\cmp{\(\cmp{\[\Rder{\(\cmp{\Rder{\cf}{\Man{1}}}{\avecff{2}}\)}{\Man{1}}\]}{\[\der{\asmooth{}}{\Man{}}{\Man{1}}\]}\)}{\avecf{1}}-
\cmp{\(\cmp{\[\Rder{\(\cmp{\Rder{\cf}{\Man{1}}}{\avecff{1}}\)}{\Man{1}}\]}{\[\der{\asmooth{}}{\Man{}}{\Man{1}}\]}\)}{\avecf{2}}\cr
&=\cmp{\Rder{\[\cmp{\Rder{\cf}{\Man{1}}}{\(\cmp{\avecff{2}}{\asmooth{}}\)}\]}{\Man{}}}{\avecf{1}}-
\cmp{\Rder{\[\cmp{\Rder{\cf}{\Man{1}}}{\(\cmp{\avecff{1}}{\asmooth{}}\)}\]}{\Man{}}}{\avecf{2}}\cr
&=\cmp{\Rder{\[\cmp{\Rder{\cf}{\Man{1}}}{\(\cmp{\[\der{\asmooth{}}{\Man{}}{\Man{1}}\]}{\avecf{2}}\)}\]}{\Man{}}}{\avecf{1}}-
\cmp{\Rder{\[\cmp{\Rder{\cf}{\Man{1}}}{\(\cmp{\[\der{\asmooth{}}{\Man{}}{\Man{1}}\]}{\avecf{1}}\)}\]}{\Man{}}}{\avecf{2}}\cr
&=\cmp{\Rder{\[\cmp{\Rder{\(\cmp{\cf}{\asmooth{}}\)}{\Man{}}}{\avecf{2}}\]}{\Man{}}}{\avecf{1}}-
\cmp{\Rder{\[\cmp{\Rder{\(\cmp{\cf}{\asmooth{}}\)}{\Man{}}}{\avecf{1}}\]}{\Man{}}}{\avecf{2}}\cr
&=\func{\[\func{\Lieder{\Man{}}}{\liebracket{\avecf{1}}{\avecf{2}}{\Man{}}}\]}{\cmp{\cf}{\asmooth{}}}\cr
&=\cmp{\(\Rder{\cf}{\Man{1}}\)}{\[\cmp{\(\der{\asmooth{}}{\Man{}}{\Man{1}}\)}{\liebracket{\avecf{1}}{\avecf{2}}{\Man{}}}\]}.
\end{aligned}
\end{align}
Therefore,
\begin{align}\label{thmliebracketofrelatedvectorfieldspeq9}
&\Foreach{\point}{\G{}}
\Foreach{\cf}{\mapdifclass{\infty}{\Man{1}}{\RR}}\cr
&\func{\(\Rder{\cf}{\Man{1}}\)}{\func{\[\cmp{\liebracket{\avecff{1}}{\avecff{2}}{\Man{1}}}{\asmooth{}}\]}{\point}}=
\func{\(\Rder{\cf}{\Man{1}}\)}{\func{\[\cmp{\(\der{\asmooth{}}{\Man{}}{\Man{1}}\)}{\liebracket{\avecf{1}}{\avecf{2}}{\Man{}}}\]}{\point}}.
\end{align}
Since every smooth vector-field on a manifold sends each point of that manifold into the tangent-space of the manifold at that point,
and the differential operator of a smooth map between a pair of manifolds maps the tangent-space of the source manifold at a given point of it
to the tangent-space of the target manifold at the image of that point under the considered smooth map, it is clear that,
\begin{align}
\func{\[\cmp{\liebracket{\avecff{1}}{\avecff{2}}{\Man{1}}}{\asmooth{}}\]}{\point}&\in\tanspace{\func{\asmooth{}}{\point}}{\Man{1}},
\label{thmliebracketofrelatedvectorfieldspeq10}\\
\func{\[\cmp{\(\der{\asmooth{}}{\Man{}}{\Man{1}}\)}{\liebracket{\avecf{1}}{\avecf{2}}{\Man{}}}\]}{\point}&\in\tanspace{\func{\asmooth{}}{\point}}{\Man{1}}.
\label{thmliebracketofrelatedvectorfieldspeq11}
\end{align}
Thus, by considering that $\Rder{\cf}{\Man{1}}$ operates linearly when restricted to $\tanspace{\func{\asmooth{}}{\point}}{\Man{1}}$,
\Ref{thmliebracketofrelatedvectorfieldspeq9} implies,
\begin{align}
&\Foreach{\point}{\G{}}
\Foreach{\cf}{\mapdifclass{\infty}{\Man{1}}{\RR}}\cr
&\func{\(\Rder{\cf}{\Man{1}}\)}{\func{\[\cmp{\liebracket{\avecff{1}}{\avecff{2}}{\Man{1}}}{\asmooth{}}\]}{\point}-
\func{\[\cmp{\(\der{\asmooth{}}{\Man{}}{\Man{1}}\)}{\liebracket{\avecf{1}}{\avecf{2}}{\Man{}}}\]}{\point}}=0.
\end{align}
Thus, according to \reflem{lemmanonvanishingderivative},
\begin{equation}
\Foreach{\point}{\G{}}
\func{\[\cmp{\liebracket{\avecff{1}}{\avecff{2}}{\Man{1}}}{\asmooth{}}\]}{\point}-
\func{\[\cmp{\(\der{\asmooth{}}{\Man{}}{\Man{1}}\)}{\liebracket{\avecf{1}}{\avecf{2}}{\Man{}}}\]}{\point}=
\zerovec{\tanspace{\func{\asmooth{}}{\point}}{\Man{1}}},
\end{equation}
and hence,
\begin{equation}
\cmp{\liebracket{\avecff{1}}{\avecff{2}}{\Man{1}}}{\asmooth{}}=
\cmp{\(\der{\asmooth{}}{\Man{}}{\Man{1}}\)}{\liebracket{\avecf{1}}{\avecf{2}}{\Man{}}}.
\end{equation}
\endthm
\section{The Theorem of Frobenius}
\textit{The theorems of this section, in particular the theorem of Frobenius, are accepted as classical facts of
differential geometry, ergo stated without proofs on this account.}
\fixed
$r$ is taken as a positive integer less than or equal to the dimension of $\Man{}$ ($r\leq n$).
Also $q$ is defined to be $n-r$.
\endfixed
\definition\label{defconsistentvectorfieldswithdistributions}
$\distribution{}$ is taken as an element of $\Distributions{r}{\Man{}}$.
\begin{equation}
\distvecf{\Man{}}{\distribution{}}{r}:=
\defset{\avecf{}}{\vecf{\Man{}}{\infty}}
{\big[\Foreach{\point}{\M{}}\func{\avecf{}}{\point}\in\func{\distribution{}}{\point}\big]}.
\end{equation}
Each element of $\distvecf{\Man{}}{\distribution{}}{r}$ is called a $\quotl$smooth vector-field on $\Man{}$
consistent with the smooth distribution $\distribution{}$ of $\Man{}$$\quotr$.
\endef
\definition\label{definvolutivedistribution}
\begin{itemize}
\item[$\centerdot$]
$\distribution{}$ is taken as an element of $\Distributions{r}{\Man{}}$.
$\distribution{}$ is referred to as an $\quotl$$r$-dimensional involutive smooth distribution of $\Man{}$$\quotr$ iff
for every pair $\avecf{1}$ and $\avecf{2}$ of smooth vector-fields on $\Man{}$ consistent with the distribution
$\distribution{}$ of $\Man{}$, their Lie-bracket $\liebracket{\avecf{1}}{\avecf{2}}{\Man{}}$ is also consistent
with the distribution $\distribution{}$.
\item[$\centerdot$]
The set of all $r$-dimensional involutive smooth distributions of the manifold $\Man{}$ is denoted by $\involutivedist{r}{\Man{}}$.
That is,
\begin{align}
&~\involutivedist{r}{\Man{}}\cr
:=&\defset{\distribution{}}{\Distributions{r}{\Man{}}}
{\(\begin{aligned}
&\Foreach{\opair{\avecf{1}}{\avecf{2}}}{\[\distvecf{\Man{}}{\distribution{}}{r}\]^{\times 2}}\cr
&\liebracket{\avecf{1}}{\avecf{2}}{\Man{}}\in\distvecf{\Man{}}{\distribution{}}{r}
\end{aligned}\)}.
\end{align}
\end{itemize}
\endef
\theorem\label{thmFrobenius}
$\[The~Theorem~of~Frobenius\]$
Any $r$-dimensional smooth distribution of $\Man{}$ is an integrable smooth distribution of $\Man{}$ if and only if it
is an involutive smooth distribution of $\Man{}$. That is, the set of all $r$-dimensional integrable smooth distributions of
$\Man{}$ equals the set of all $r$-dimensional involutive smooth distributions of $\Man{}$.
\begin{equation}
\integrabledist{r}{\Man{}}=
\involutivedist{r}{\Man{}}.
\end{equation}
\endthm
\theorem\label{thmflowsofcommutativevectorfieldscommute}
Each $\avecf{1}$ and $\avecf{2}$ is taken as a complete vector-field on $\Man{}$.
If $\avecf{1}$ and $\avecf{2}$ commute, then for every pair of real numbers $t$ and $s$,
$\func{\vfFlow{\Man{}}{\avecf{1}}}{t}$ and $\func{\vfFlow{\Man{}}{\avecf{2}}}{s}$ also commute
in the group of $\infty$-automorphisms of $\Man{}$.
\begin{align}
&\hskip0.9\baselineskip\liebracket{\avecf{1}}{\avecf{2}}{\Man{}}=\zerovec{\Vecf{\Man{}}{\infty}}\cr
&\then\
\Foreach{\opair{t}{s}}{\R^2}
\cmp{\func{\vfFlow{\Man{}}{\avecf{1}}}{t}}{\func{\vfFlow{\Man{}}{\avecf{2}}}{s}}=
\cmp{\func{\vfFlow{\Man{}}{\avecf{2}}}{s}}{\func{\vfFlow{\Man{}}{\avecf{1}}}{t}}.
\end{align}
\endthm
\section{The Module structure of Smooth Vector-Fields}
\subsection{Basics of Modules over Rings}
\textit{The definitions of algebraic structures ring and module are subject to subtle variations in the literature.
The definitions used here is in agreement with those proposed in \cite{MacLane}. The properties related to rings and modules utilized here
are in alignment with these forms of definitions, which are asserted without proofs.
So for the sake of clarity, the intended definitions are explicitly stated.}
\definition\label{defring}
$\ring{}$ is taken  as a non-empty set, and each $\vsum{r}$ and $\vsprod{r}$ as an element of $\Func{\Cprod{\ring{}}{\ring{}}}{\ring{}}$
(a binary operations on $\ring{}$). The triple $\triple{\ring{}}{\vsum{r}}{\vsprod{r}}$ is referred to as a $\quotl$ring$\quotr$ iff
the following axioms are hold.
\begin{itemize}
\item[\myitem{Rin~1.}]
$\opair{\ring{}}{\vsum{r}}$ possesses the structure of an abelian group.
\item[\myitem{Rin~2.}]
$\opair{\ring{}}{\vsprod{r}}$ possesses the structure of a monoid.
\item[\myitem{Rin~3.}]
For every $\rr{1}$, $\rr{2}$, and $\rr{3}$ in $\ring{}$,
\begin{align}
&\rr{1}\vsprod{r}\(\rr{2}\vsum{r}\rr{3}\)=\(\rr{1}\vsprod{r}\rr{2}\)\vsum{r}\(\rr{1}\vsprod{r}\rr{3}\),\\
&\(\rr{1}\vsum{r}\rr{2}\)\vsprod{r}\rr{3}=\(\rr{1}\vsprod{r}\rr{3}\)\vsum{r}\(\rr{2}\vsprod{r}\rr{3}\).
\end{align}
\end{itemize}
When $\Ring{}:=\triple{\ring{}}{\vsum{r}}{\vsprod{r}}$ is a ring, $\ringzero{\Ring{}}$ is defined to be
the unique neutral element of the group $\opair{\ring{}}{\vsum{r}}$, and $\ringid{\Ring{}}$ is defined to be
the unique neutral element of the monoid $\opair{\ring{}}{\vsprod{r}}$. $\ringzero{\Ring{}}$ and $\ringid{\Ring{}}$
are called the $\quotl$zero (element) of the ring $\Ring{}$$\quotr$ and the $\quotl$unit (element) of the ring $\Ring{}$$\quotr$,
respectively.\\
When $\Ring{}:=\triple{\ring{}}{\vsum{r}}{\vsprod{r}}$ is a ring, $\vsum{r}$ and $\vsprod{r}$ are called the
$\quotl$addition (operation) of the ring $\Ring{}$$\quotr$ and the
$\quotl$multiplication (operation) of the ring $\Ring{}$$\quotr$, respectively.
\endef
\definition\label{defcommutativering}
$\ring{}$ is taken  as a non-empty set, and each $\vsum{r}$ and $\vsprod{r}$ as an element of $\Func{\Cprod{\ring{}}{\ring{}}}{\ring{}}$.
The triple $\triple{\ring{}}{\vsum{r}}{\vsprod{r}}$ is referred to as a $\quotl$commutative-ring$\quotr$ iff
the following properties are hold.
\begin{itemize}
\item[\myitem{CR~1.}]
$\triple{\ring{}}{\vsum{r}}{\vsprod{r}}$ is a ring.
\item[\myitem{CR~2.}]
The multiplication operation of $\Ring{}$ is commutative. That is,
\begin{equation}
\Foreach{\opair{\rr{1}}{\rr{2}}}{\Cprod{\ring{}}{\ring{}}}
\rr{1}\vsprod{r}\rr{2}=\rr{2}\vsprod{r}\rr{1}.
\end{equation}
\end{itemize}
\endef
\definition\label{defmodule}
$\Ring{}=\triple{\ring{}}{\vsum{r}}{\vsprod{r}}$ is taken as a ring, $\module{}$ a non-empty set,
$\vsum{m}$ an element of $\Func{\Cprod{\module{}}{\module{}}}{\module{}}$ (a binary operation on $\module{}$), and
$\spro{}$ as an element of $\Func{\Cprod{\ring{}}{\module{}}}{\module{}}$. The $4$-tuple $\tuple{\module{}}{\Ring{}}{\vsum{m}}{\spro{}}$
is referred to as a $\quotl$left module over the ring $\Ring{}$$\quotr$ or a $\quotl$left $\Ring{}$-module$\quotr$ iff
the following axioms are hold.
\begin{itemize}
\item[\myitem{Mo~1.}]
$\opair{\module{}}{\vsum{m}}$ is an abelian group.
\item[\myitem{Mo~2.}]
\begin{align}
\Foreach{\opair{\rr{1}}{\rr{2}}}{\Cprod{\ring{}}{\ring{}}}
\Foreach{\mm{}}{\module{}}
&\(\rr{1}\vsum{r}\rr{2}\)\spro{}\mm{}=\(\rr{1}\spro{}\mm{}\)\vsum{m}\(\rr{2}\spro{}\mm{}\),\\
\Foreach{\rr{}}{\ring{}}
\Foreach{\opair{\mm{1}}{\mm{2}}}{\Cprod{\module{}}{\module{}}}
&\rr{}\spro{}\(\mm{1}\vsum{m}\mm{2}\)=\(\rr{}\spro{}\mm{1}\)\vsum{m}\(\rr{}\spro{}\mm{2}\),\\
\Foreach{\opair{\rr{1}}{\rr{2}}}{\Cprod{\ring{}}{\ring{}}}
\Foreach{\mm{}}{\module{}}
&\rr{1}\spro{}\(\rr{2}\spro{}\mm{}\)=\(\rr{1}\vsprod{r}\rr{2}\)\spro{}\mm{}.
\end{align}
\item[\myitem{Mo~3.}]
\begin{equation}
\Foreach{\mm{}}{\module{}}
\ringid{\Ring{}}\spro{}\mm{}=\mm{}.
\end{equation}
\end{itemize}
When $\Module{}:=\tuple{\module{}}{\Ring{}}{\vsum{m}}{\spro{}}$ is a left $\Ring{}$-module, then $\modulezero{\Module{}}$ is defined to be
the unique neutral element of the group $\opair{\module{}}{\vsum{m}}$, and is called the
$\quotl$zero (element) of the $\Ring{}$-module $\Module{}$$\quotr$.\\
When $\Module{}:=\tuple{\module{}}{\Ring{}}{\vsum{m}}{\spro{}}$ is a left $\Ring{}$-module, then $\vsum{m}$ and $\spro{}$
are referred to as the $\quotl$addition (operation) of the $\Ring{}$-module $\Module{}$$\quotr$ and the
$\quotl$scalar-product (operation) of the $\Ring{}$-module $\Module{}$$\quotr$, respectively.\\
When $\Module{}:=\tuple{\module{}}{\Ring{}}{\vsum{m}}{\spro{}}$ is a left $\Ring{}$-module and $\Ring{}$ is a commutative ring, then
$\Module{}$ is called a $\quotl$module over the commutative ring $\Ring{}$$\quotr$.
\endef
\definition\label{defsubmodule}
$\Ring{}=\triple{\ring{}}{\vsum{r}}{\vsprod{r}}$ is taken as a  ring, $\Module{}=\tuple{\module{}}{\Ring{}}{\vsum{m}}{\spro{}}$
a left $\Ring{}$-module, and $\module{1}$ a non-empty subset of $\module{}$. $\module{1}$ is called a
$\quotl$submodule of the left $\Ring{}$-module $\Module{}$$\quotr$ iff $\module{1}$ when endowed with the restrictions of the addition and
scalar-multiplication operations of $\Module{}$ to $\module{}$, that is
$\tuple{\module{1}}{\Ring{}}{\func{\res{\vsum{m}}}{\Cprod{\module{1}}{\module{1}}}}{\func{\res{\spro{}}}{\Cprod{\ring{}}{\module{1}}}}$
is itself a left $\Ring{}$-module. The set of all submodules of the left $\Ring{}$-module $\Module{}$ is denoted by $\Submodules{\Module{}}$.\\
When $\module{1}$ is a submodule of $\Module{}$, the left $\Ring{}$-module
$\tuple{\module{1}}{\Ring{}}{\func{\res{\vsum{m}}}{\Cprod{\module{1}}{\module{1}}}}{\func{\res{\spro{}}}{\Cprod{\ring{}}{\module{1}}}}$
is denoted by $\subspace{\Module{}}{\module{1}}$ by definition.
\endef
\theorem\label{thmsubmodule1}
$\Ring{}=\triple{\ring{}}{\vsum{r}}{\vsprod{r}}$ is taken as a  ring, $\Module{}=\tuple{\module{}}{\Ring{}}{\vsum{m}}{\spro{}}$,
and $\module{1}$ as a subset of $\module{}$. $\module{1}$ is a submodule of $\Module{}$ if and only if
$\module{1}$ is non-empty and closed under the addition and scalar-product operations of $\Module{}$. That is,
\begin{align}
\module{1}\in\Submodules{\Module{}}\thenn
\(\begin{aligned}
&\Foreach{\opair{\mm{1}}{\mm{2}}}{\Cprod{\module{1}}{\module{1}}}
\mm{1}\vsum{m}\mm{2}\in\module{1},\cr
&\Foreach{\opair{\rr{}}{\mm{}}}{\Cprod{\ring{}}{\module{}}}
\rr{}\spro{}\mm{}\in\module{1}
\end{aligned}\).
\end{align}
\endthm
\theorem
$\Ring{}=\triple{\ring{}}{\vsum{r}}{\vsprod{r}}$ is taken as a  ring, and $\Module{}=\tuple{\module{}}{\Ring{}}{\vsum{m}}{\spro{}}$
a left $\Ring{}$-module. For every non-void collection $\collection{}$ of submodules of $\Module{}$, that is for every non-empty subset of
$\Submodules{\Module{}}$, the intersection of all elements of $\collection{}$ is a submodule of $\Module{}$.
\begin{equation}
\Foreach{\collection{}}{\[\compl{\CSs{\Submodules{\Module{}}}}{\seta{\empty}}\]}
\intersection{\collection{}}\in\Submodules{\Module{}}.
\end{equation}
\endthm
\definition\label{defmodulegenerator}
$\Ring{}=\triple{\ring{}}{\vsum{r}}{\vsprod{r}}$ is taken as a  ring, $\Module{}=\tuple{\module{}}{\Ring{}}{\vsum{m}}{\spro{}}$
a left $\Ring{}$-module, and $\SET{}$ a subset of $\module{}$. The intersection of all submodules of $\Module{}$ that include
$\SET{}$ is called the $\quotl$submodule of $\Module{}$ generated by $\SET{}$$\quotr$, and denoted by $\func{\modulegen{\Module{}}}{\SET{}}$.
\begin{equation}
\func{\modulegen{\Module{}}}{\SET{}}:=
\intersection{\defset{\module{1}}{\Submodules{\Module{}}}{\module{1}\supseteq\SET{}}}.
\end{equation}
\endef
\theorem
$\Ring{}=\triple{\ring{}}{\vsum{r}}{\vsprod{r}}$ is taken as a  ring, $\Module{}=\tuple{\module{}}{\Ring{}}{\vsum{m}}{\spro{}}$
a left $\Ring{}$-module, and $\SET{}$ a non-empty subset of $\module{}$.
\begin{equation}
\func{\modulegen{\Module{}}}{\SET{}}=
\defset{\mm{}}{\module{}}{
\[\begin{aligned}
&\Exists{n}{\Zp}
\Exists{\mtuple{\rr{1}}{\rr{n}}}{{\ring{}}^{\times n}}
\Exists{\mtuple{\mm{1}}{\mm{n}}}{{\SET{}}^{\times n}}\cr
&\mm{}=\sum_{k=1}^{n}\rr{k}\spro{}\mm{k}
\end{aligned}\]}.
\end{equation}
In particular, if $\SET{}=\seta{\suc{\mm{1}}{\mm{n}}}$ and $\Card{\SET{}}=n$ for some $n\in\Zp$, then
\begin{equation}
\func{\modulegen{\Module{}}}{\SET{}}=
\defSet{\sum_{k=1}^{n}\rr{k}\spro{}\mm{k}}
{\mtuple{\rr{1}}{\rr{n}}\in{\ring{}}^{\times n}}.
\end{equation}
\endthm
\definition\label{deflinearlyindependentsetsofmodule}
$\Ring{}=\triple{\ring{}}{\vsum{r}}{\vsprod{r}}$ is taken as a  ring, $\Module{}=\tuple{\module{}}{\Ring{}}{\vsum{m}}{\spro{}}$
a left $\Ring{}$-module, and $\SET{}$ a subset of $\module{}$. $\SET{}$ is called a $\quotl$linearly-independent set of the left
$\Ring{}$-module $\Module{}$$\quotr$ iff,
\begin{align}
&\Foreach{n}{\Zp}
\Foreach{\seta{\suc{\mm{1}}{\mm{n}}}}{\defset{\asubset}{\CSs{\SET{}}}{\Card{\asubset}=n}}
\Foreach{\mtuple{\rr{1}}{\rr{n}}}{{\ring{}}^{\times n}}\cr
&\[\(\sum_{k=1}^{n}\rr{k}\spro{}\mm{k}=\modulezero{\Module{}}\)\then\(\Foreach{k}{\seta{\suc{1}{n}}}\rr{k}=\ringzero{\Ring{}}\)\].
\end{align}
\endef
\definition\label{defbasisofmodule}
$\Ring{}=\triple{\ring{}}{\vsum{r}}{\vsprod{r}}$ is taken as a  ring, $\Module{}=\tuple{\module{}}{\Ring{}}{\vsum{m}}{\spro{}}$
a left $\Ring{}$-module, and $\SET{}$ a subset of $\module{}$. $\SET{}$ is called a
$\quotl$module-basis of the left $\Ring{}$-module $\Module{}$$\quotr$ iff the following properties are hold.
\begin{itemize}
\item[\myitem{MB~1.}]
$\SET{}$ is a linearly-independent set of the left $\Ring{}$-module $\Module{}$.
\item[\myitem{MB~2.}]
The submodule of $\Module{}$ generated by $\SET{}$ equals $\module{}$ itself. That is,
\begin{equation}
\func{\modulegen{\Module{}}}{\SET{}}=\module{}.
\end{equation}
\end{itemize}
The set of all module-bases of $\Module{}$ is denoted by $\Modulebases{\Module{}}$.
\endef
\definition\label{deffreemodule}
$\Ring{}=\triple{\ring{}}{\vsum{r}}{\vsprod{r}}$ is taken as a  ring, $\Module{}=\tuple{\module{}}{\Ring{}}{\vsum{m}}{\spro{}}$
a left $\Ring{}$-module. $\Module{}$ is referred to as a $\quotl$free left $\Ring{}$-module$\quotr$ or simply as a
$\quotl$free module$\quotr$ iff there exists at least one module-basis of $\Module{}$, that is,
\begin{equation}
\Modulebases{\Module{}}\neq\empty.
\end{equation}
\endef
\definition\label{defmodulehomomorphism}
$\Ring{}=\triple{\ring{}}{\vsum{r}}{\vsprod{r}}$ is taken as a  ring, each $\Module{}=\tuple{\module{}}{\Ring{}}{\vsum{m}}{\spro{}}$
and $\Module{1}=\tuple{\module{1}}{\Ring{}}{\vsum{m_1}}{\spro{1}}$ as a left $\Ring{}$-module.
\begin{itemize}
\item[$\centerdot$]
$\amodmor{}$ is taken as an element of
$\Func{\module{}}{\module{1}}$ (a function from $\module{}$ to $\module{1}$). $\amodmor{}$ is referred to as a
$\quotl$module-homomorphism from the module $\Module{}$ to the module $\Module{1}$$\quotr$ or an
$\quotl$$\Ring{}$-module-homomorphism from $\Module{}$ to $\Module{1}$$\quotr$ iff,
\begin{align}
&\Foreach{\opair{\mm{}}{\p{\mm{}}}}{\Cprod{\module{}}{\module{}}}
\func{\amodmor{}}{\mm{}\vsum{m}\p{\mm{}}}=\func{\amodmor{}}{\mm{}}\vsum{m_1}\func{\amodmor{}}{\p{\mm{}}},\cr
&\Foreach{\opair{\rr{}}{\mm{}}}{\Cprod{\ring{}}{\module{}}}
\func{\amodmor{}}{\rr{}\spro{}\mm{}}=\rr{}\spro{1}\[\func{\amodmor{}}{\mm{}}\].
\end{align}
\item[$\centerdot$]
The set of all $\Ring{}$-module-homomorphisms from $\Module{}$ to $\Module{1}$ is denoted by $\ModuleMors{\Module{}}{\Module{1}}$.
That is,
\begin{align}
&\ModuleMors{\Module{}}{\Module{1}}:=\cr
&~~~\defset{\amodmor{}}{\Func{\module{}}{\module{1}}}
{\begin{aligned}
&\Foreach{\opair{\mm{}}{\p{\mm{}}}}{\Cprod{\module{}}{\module{}}}
\func{\amodmor{}}{\mm{}\vsum{m}\p{\mm{}}}=\func{\amodmor{}}{\mm{}}\vsum{m_1}\func{\amodmor{}}{\p{\mm{}}},\cr
&\Foreach{\opair{\rr{}}{\mm{}}}{\Cprod{\ring{}}{\module{}}}
\func{\amodmor{}}{\rr{}\spro{}\mm{}}=\rr{}\spro{1}\[\func{\amodmor{}}{\mm{}}\]
\end{aligned}}.\cr
&{}
\end{align}
\item[$\centerdot$]
$\amodmor{}$ is taken as an element of $\Func{\module{}}{\module{1}}$.
$\amodmor{}$ is referred to as a
$\quotl$module-isomorphism from the module $\Module{}$ to the module $\Module{1}$$\quotr$ or an
$\quotl$$\Ring{}$-module-isomorphism from $\Module{}$ to $\Module{1}$$\quotr$ iff
$\amodmor{}$ is a bijection from $\module{}$ to $\module{1}$ and both $\amodmor{}$ and $\finv{\amodmor{}}$
are $\Ring{}$-module-homomorphisms from their target left $\Ring{}$-module to their source left $\Ring{}$-module, that is,
\begin{align}
\amodmor{}&\in\IF{\module{}}{\module{1}},\\
\amodmor{}&\in\ModuleMors{\Module{}}{\Module{1}},\\
\finv{\amodmor{}}&\in\ModuleMors{\Module{1}}{\Module{}}.
\end{align}
\item[$\centerdot$]
The set of all $\Ring{}$-module-isomorphisms from $\Module{}$ to $\Module{1}$ is denoted by $\ModuleIsoms{\Module{}}{\Module{1}}$.
That is,
\begin{equation}
\ModuleIsoms{\Module{}}{\Module{1}}:=
\defset{\amodmor{}}{\IF{\module{}}{\module{1}}}
{\amodmor{}\in\ModuleMors{\Module{}}{\Module{1}},
\finv{\amodmor{}}\in\ModuleMors{\Module{1}}{\Module{}}}.
\end{equation}
\end{itemize}
\endef
\theorem\label{thmmoduleisomorphismequiv1}
$\Ring{}=\triple{\ring{}}{\vsum{r}}{\vsprod{r}}$ is taken as a  ring, each $\Module{}=\tuple{\module{}}{\Ring{}}{\vsum{m}}{\spro{}}$
and $\Module{1}=\tuple{\module{1}}{\Ring{}}{\vsum{m_1}}{\spro{1}}$ as a left $\Ring{}$-module. For every $\amodmor{}$ in
$\IF{\module{}}{\module{1}}$, $\amodmor{}$ is an $\Ring{}$-module-isomorphism from $\Module{}$ t $\Module{1}$ if and only if
it is an $\Ring{}$-module-homomorphism from $\Module{}$ t $\Module{1}$. That is,
\begin{align}
\ModuleIsoms{\Module{}}{\Module{1}}&=
\defset{\amodmor{}}{\IF{\module{}}{\module{1}}}
{\amodmor{}\in\ModuleMors{\Module{}}{\Module{1}}}\cr
&=\defset{\amodmor{}}{\IF{\module{}}{\module{1}}}
{\finv{\amodmor{}}\in\ModuleMors{\Module{1}}{\Module{}}}.
\end{align}
\endthm
\subsection{Smooth Vector-Fields and Derivations as Modules}
\definition\label{defsrealvaluedmoothmapsproduct}
The mappings $\smoothfsum{\Man{}}$, $\smoothfprod{\Man{}}$, $\smoothvfsum{\Man{}}$,
$\smoothvfprod{\Man{}}$, $\smoothdersum{\Man{}}$, and $\smoothderprod{\Man{}}$ are defined as the following.
\begin{align}
&\smoothfsum{\Man{}}\indef\Func{\Cprod{\mapdifclass{\infty}{\Man{}}{\RR}}{\mapdifclass{\infty}{\Man{}}{\RR}}}
{\mapdifclass{\infty}{\Man{}}{\RR}},\cr
&\Foreach{\opair{\cf}{\cg}}{\Cprod{\mapdifclass{\infty}{\Man{}}{\RR}}{\mapdifclass{\infty}{\Man{}}{\RR}}}
\[\Foreach{\point}{\M{}}\func{\(\cf\smoothfsum{\Man{}}\cg\)}{\point}\eqdef
\func{\cf}{\point}+\func{\cg}{\point}\].
\end{align}
\begin{align}
&\smoothfprod{\Man{}}\indef\Func{\Cprod{\mapdifclass{\infty}{\Man{}}{\RR}}{\mapdifclass{\infty}{\Man{}}{\RR}}}
{\mapdifclass{\infty}{\Man{}}{\RR}},\cr
&\Foreach{\opair{\cf}{\cg}}{\Cprod{\mapdifclass{\infty}{\Man{}}{\RR}}{\mapdifclass{\infty}{\Man{}}{\RR}}}
\[\Foreach{\point}{\M{}}\func{\(\cf\smoothfprod{\Man{}}\cg\)}{\point}\eqdef
\func{\cf}{\point}\func{\cg}{\point}\].
\end{align}
\begin{align}
&\smoothvfsum{\Man{}}\indef\Func{\Cprod{\vecf{\Man{}}{\infty}}{\vecf{\Man{}}{\infty}}}
{\vecf{\Man{}}{\infty}},\cr
&\Foreach{\opair{\avecf{1}}{\avecf{2}}}{\Cprod{\vecf{\Man{}}{\infty}}{\vecf{\Man{}}{\infty}}}
\[\Foreach{\point}{\M{}}\func{\(\avecf{1}\smoothvfsum{\Man{}}\avecf{2}\)}{\point}\eqdef
\func{\avecf{1}}{\point}\vsum{\opair{\Man{}}{\point}}\func{\avecf{2}}{\point}\].\cr
&{}
\end{align}
\begin{align}
&\smoothvfprod{\Man{}}\indef\Func{\Cprod{\mapdifclass{\infty}{\Man{}}{\RR}}{\vecf{\Man{}}{\infty}}}{\vecf{\Man{}}{\infty}},\cr
&\Foreach{\opair{\cf}{\avecf{}}}{\Cprod{\mapdifclass{\infty}{\Man{}}{\RR}}{\vecf{\Man{}}{\infty}}}
\[\Foreach{\point}{\M{}}\func{\(\cf\smoothvfprod{\Man{}}\avecf{}\)}{\point}\eqdef
\func{\cf}{\point}\vsprod{\opair{\Man{}}{\point}}\func{\avecf{}}{\point}\].\cr
&{}
\end{align}
\begin{align}
&\smoothdersum{\Man{}}\indef\Func{\Cprod{\Derivation{\Man{}}{\infty}}{\Derivation{\Man{}}{\infty}}}
{\Derivation{\Man{}}{\infty}},\cr
&\Foreach{\opair{\aderivation{1}}{\aderivation{2}}}{\Cprod{\Derivation{\Man{}}{\infty}}{\Derivation{\Man{}}{\infty}}}\cr
&\hskip6\baselineskip\[\begin{aligned}
\Foreach{\cf}{\mapdifclass{\infty}{\Man{}}}
\func{\(\aderivation{1}\smoothdersum{\Man{}}\aderivation{2}\)}{\cf}\eqdef
\func{\aderivation{1}}{\cf}\smoothfsum{\Man{}}\func{\aderivation{2}}{\cf}
\end{aligned}\].\cr
&{}
\end{align}
\begin{align}
&\smoothderprod{\Man{}}\indef\Func{\Cprod{\mapdifclass{\infty}{\Man{}}{\RR}}{\Derivation{\Man{}}{\infty}}}{\Derivation{\Man{}}{\infty}},\cr
&\Foreach{\opair{\cf}{\aderivation{}}}{\Cprod{\mapdifclass{\infty}{\Man{}}{\RR}}{\Derivation{\Man{}}{\infty}}}\cr
&\hskip6\baselineskip\[\begin{aligned}
\Foreach{\cg}{\mapdifclass{\infty}{\Man{}}{\RR}}
\func{\(\cf\smoothderprod{\Man{}}\aderivation{}\)}{\cg}\eqdef
\cf\smoothfprod{\Man{}}\func{\aderivation{}}{\cg}
\end{aligned}\].\cr
&{}
\end{align}
\caution
The well-definedness of each of these mappings is a standard fact of differential geometry.\\
When the context is clear enough, both $\smoothfsum{\Man{}}$ and $\smoothvfsum{\Man{}}$ can be replaced simply by
$\vsum{}$, and $\cf\smoothfprod{\Man{}}\cg$ and $\cf\smoothvfprod{\Man{}}\avecf{}$ can be denoted by
$\cf\cg$ and $\cf\avecf{}$ respectively.
\endef
\theorem
$\triple{\mapdifclass{\infty}{\Man{}}{\RR}}{\smoothfsum{\Man{}}}{\smoothfprod{\Man{}}}$ is a commutative ring.
\proof
The verification of axioms of commutative ring is trivial.
\endthm
\definition
$\smoothring{\infty}{\Man{}}$ denotes the commutative ring
$\triple{\mapdifclass{\infty}{\Man{}}{\RR}}{\smoothfsum{\Man{}}}{\smoothfprod{\Man{}}}$ by definition,
and is called the $\quotl$commutative ring of real-valued smooth maps on the manifold $\Man{}$$\quotr$.
Also, $\smoothfsub{\Man{}}$ denotes the subtraction operation of this ring.
\endef
\theorem\label{thmvecfmodule}
$\tuple{\vecf{\Man{}}{\infty}}{\smoothring{\infty}{\Man{}}}{\smoothvfsum{\Man{}}}{\smoothvfprod{\Man{}}}$ is
a $\smoothring{\infty}{\Man{}}$-module (a module over the commutative ring $\smoothring{\infty}{\Man{}}$).
\proof
The verification of module axioms in this case is trivial.
\endthm
\definition\label{defvecfmodule}
$\smoothvfmodule{\infty}{\Man{}}$ denotes the $\smoothring{\infty}{\Man{}}$-module
$\tuple{\vecf{\Man{}}{\infty}}{\smoothring{\infty}{\Man{}}}{\smoothvfsum{\Man{}}}{\smoothvfprod{\Man{}}}$ by definition,
and is called the $\quotl$$\smoothring{\infty}{\Man{}}$-module (or module) of smooth vector-fields on the manifold
$\Man{}$$\quotr$ or simply as the $\quotl$vector-fields module of $\Man{}$$\quotr$.
Also, $\smoothvfsub{\Man{}}$ denotes the subtraction operation of this module.
\endef
\theorem\label{thmdermodule}
$\tuple{\Derivation{\Man{}}{\infty}}{\smoothring{\infty}{\Man{}}}{\smoothdersum{\Man{}}}{\smoothderprod{\Man{}}}$ is
a $\smoothring{\infty}{\Man{}}$-module (a module over the commutative ring $\smoothring{\infty}{\Man{}}$).
\proof
The verification of module axioms in this case is also trivial.
\endthm
\definition\label{defdermodule}
$\smoothdermodule{\infty}{\Man{}}$ denotes the $\smoothring{\infty}{\Man{}}$-module
$\tuple{\Derivation{\Man{}}{\infty}}{\smoothring{\infty}{\Man{}}}{\smoothdersum{\Man{}}}{\smoothderprod{\Man{}}}$ by definition,
and is called the $\quotl$$\smoothring{\infty}{\Man{}}$-module (or module) of $\infty$-derivations on the manifold
$\Man{}$$\quotr$ or simply as the $\quotl$derivation module of $\Man{}$$\quotr$.
Also, $\smoothdersub{\Man{}}$ denotes the subtraction operation of this module.
\endef
\theorem\label{thmliederivativeisamodulehomomorphism}
The Lie-derivative operator on $\Man{}$, $\Lieder{\Man{}}$,
is an $\smoothring{\infty}{\Man{}}$-module-isomorphism from $\smoothvfmodule{\infty}{\Man{}}$ to
$\smoothdermodule{\infty}{\Man{}}$.
\begin{equation}
\Lieder{\Man{}}\in\ModuleIsoms{\smoothvfmodule{\infty}{\Man{}}}
{\smoothdermodule{\infty}{\Man{}}}.
\end{equation}
\proof
Since $\smoothvfsum{\Man{}}$ and $\smoothdersum{\Man{}}$ are identically the addition operations of
the canonical linear structures of smooth vector-fields and smooth derivations on $\Man{}$, that is
the addition operations of the vector-spaces $\Vecf{\Man{}}{\infty}$ and $\LDerivation{\Man{}}{\infty}$, respectively,
and considering that $\Lieder{\Man{}}$ is a linear map from $\Vecf{\Man{}}{\infty}$ to $\LDerivation{\Man{}}{\infty}$
(based on \reflem{lemmaliederivativeislinear}), it is evident that,
\begin{equation}\label{thmliederivativeisamodulehomomorphismpeq1}
\Foreach{\opair{\avecf{1}}{\avecf{2}}}{\Cprod{\vecf{\Man{}}{\infty}}{\vecf{\Man{}}{\infty}}}
\func{\Lieder{\Man{}}}{\avecf{1}\smoothvfsum{\Man{}}\avecf{2}}=
\func{\Lieder{\Man{}}}{\avecf{1}}\smoothdersum{\Man{}}\func{\Lieder{\Man{}}}{\avecf{2}}.
\end{equation}
\begin{itemize}
\item[$\pr{1}$]
$\cf$ is taken as an arbitrary element of $\mapdifclass{\infty}{\Man{}}{\RR}$, and
$\avecf{}$ as an arbitrary element of $\vecf{\Man{}}{\infty}$. According to \refdef{defliederivative},
\refdef{defliederivative00}, and \refdef{defsrealvaluedmoothmapsproduct}, and considering that the derivative
of a real-valued smooth map on $\Man{}$ operates linearly when restricted to the tangent-space of any point of $\Man{}$,
\begin{align}
\Foreach{\cg}{\mapdifclass{\infty}{\Man{}}{\RR}}
\Foreach{\point}{\M{}}
\func{\(\func{\[\func{\Lieder{\Man{}}}{\cf\smoothvfprod{\Man{}}\avecf{}}\]}{\cg}\)}{\point}&=
\func{\(\cmp{\Rder{\cg}{\Man{}}}{\cf\avecf{}}\)}{\point}\cr
&=\func{\Rder{\cg}{\Man{}}}{\func{\[\cf\avecf{}\]}{\point}}\cr
&=\func{\Rder{\cg}{\Man{}}}{\func{\cf}{\point}\func{\avecf{}}{\point}}\cr
&=\func{\cf}{\point}\[\func{\Rder{\cg}{\Man{}}}{\func{\avecf{}}{\point}}\]\cr
&=\func{\cf}{\point}\[\func{\(\cmp{\Rder{\cg}{\Man{}}}{\avecf{}}\)}{\point}\]\cr
&=\func{\cf}{\point}\[\func{\(\func{\[\func{\Lieder{\Man{}}}{\avecf{}}\]}{\cg}\)}{\point}\]\cr
&=\func{\[\func{\(\cf\smoothderprod{\Man{}}\[\func{\Lieder{\Man{}}}{\avecf{}}\]\)}{\cg}\]}{\point},\cr
&{}
\end{align}
and hence,
\begin{equation}
\func{\Lieder{\Man{}}}{\cf\smoothvfprod{\Man{}}\avecf{}}=
\cf\smoothderprod{\Man{}}\[\func{\Lieder{\Man{}}}{\avecf{}}\].
\end{equation}
\endp
\end{itemize}
Therefore,
\begin{align}\label{thmliederivativeisamodulehomomorphismpeq2}
\Foreach{\opair{\cf}{\avecf{}}}{\Cprod{\mapdifclass{\infty}{\Man{}}{\RR}}{\vecf{\Man{}}{\infty}}}
\func{\Lieder{\Man{}}}{\cf\smoothvfprod{\Man{}}\avecf{}}=
\cf\smoothderprod{\Man{}}\[\func{\Lieder{\Man{}}}{\avecf{}}\].
\end{align}
Thus according to \Ref{thmliederivativeisamodulehomomorphismpeq1} and \Ref{thmliederivativeisamodulehomomorphismpeq2}, and
based on \refdef{defmodulehomomorphism} and also \refthm{thmvecfmodule}, \refdef{defvecfmodule},
\refthm{thmdermodule}, and \refdef{defdermodule}, $\Lieder{\Man{}}$ is an
$\smoothring{\infty}{\Man{}}$-module-homomorphism from $\smoothvfmodule{\infty}{\Man{}}$ to
$\smoothdermodule{\infty}{\Man{}}$.
\begin{equation}\label{thmliederivativeisamodulehomomorphismpeq3}
\Lieder{\Man{}}\in\ModuleMors{\smoothvfmodule{\infty}{\Man{}}}
{\smoothdermodule{\infty}{\Man{}}}.
\end{equation}
Moreover, according to \refcor{corliederivativeisalinearisomorphism}, $\Lieder{\Man{}}$ is a bijection
from $\vecf{\Man{}}{\infty}$ to $\Derivation{\Man{}}{\infty}$.
\begin{equation}\label{thmliederivativeisamodulehomomorphismpeq4}
\Lieder{\Man{}}\in\IF{\vecf{\Man{}}{\infty}}{\Derivation{\Man{}}{\infty}}.
\end{equation}
Therefore based on \refthm{thmmoduleisomorphismequiv1}, \Ref{thmliederivativeisamodulehomomorphismpeq3}
and \Ref{thmliederivativeisamodulehomomorphismpeq4} imply that $\Lieder{\Man{}}$ is an
$\smoothring{\infty}{\Man{}}$-module-isomorphism from $\smoothvfmodule{\infty}{\Man{}}$ to
$\smoothdermodule{\infty}{\Man{}}$.
\begin{equation}
\Lieder{\Man{}}\in\ModuleIsoms{\smoothvfmodule{\infty}{\Man{}}}
{\smoothdermodule{\infty}{\Man{}}}.
\end{equation}
\endthm
\lemma\label{lemliebracketonthemoduleofvectorfields}
\begin{align}
&\Foreach{\cf}{\mapdifclass{\infty}{\Man{}}{\RR}}
\Foreach{\opair{\avecf{1}}{\avecf{2}}}{\Cprod{\vecf{\Man{}}{\infty}}{\vecf{\Man{}}{\infty}}}\cr
&\liebracket{\avecf{1}}{\cf\smoothvfprod{\Man{}}\avecf{2}}{\Man{}}=
\[\cf\smoothvfprod{\Man{}}\liebracket{\avecf{1}}{\avecf{2}}{\Man{}}\]\smoothvfsum{\Man{}}
\[\bigg(\func{\[\func{\Lieder{\Man{}}}{\avecf{1}}\]}{\cf}\bigg)\smoothvfprod{\Man{}}\avecf{2}\].
\end{align}
\proof
Since for every smooth vector-field $\avecf{}$ on $\Man{}$, $\func{\Lieder{\Man{}}}{\avecf{}}$ is an $\infty$-derivation on
$\Man{}$, based on the definition of $\infty$-derivations on a manifold,
\begin{align}\label{lemliebracketonthemoduleofvectorfieldspeq1}
\Foreach{\avecf{}}{\vecf{\Man{}}{\infty}}
&\Foreach{\opair{\cf}{\cg}}{\Cprod{\mapdifclass{\infty}{\Man{}}{\RR}}{\mapdifclass{\infty}{\Man{}}{\RR}}}\cr
&\func{\[\func{\Lieder{\Man{}}}{\avecf{}}\]}{\cf\smoothfprod{\Man{}}\cg}=
\bigg(\func{\[\func{\Lieder{\Man{}}}{\avecf{}}\]}{\cf}\bigg)\smoothfprod{\Man{}}\cg\smoothfsum{\Man{}}
\cf\smoothfprod{\Man{}}\bigg(\func{\[\func{\Lieder{\Man{}}}{\avecf{}}\]}{\cg}\bigg).\cr
&{}
\end{align}
Each $\avecf{1}$ and $\avecf{2}$ is taken as an arbitrary element of $\vecf{\Man{}}{\infty}$, and
$\cf$ as an arbitray element of $\mapdifclass{\infty}{\Man{}}{\RR}$.
According to \refdef{defliederivative} and \refdef{defliederivative00}, \refthm{thmliederivativeisamodulehomomorphism},
\refdef{defmodulehomomorphism}, \refdef{defsrealvaluedmoothmapsproduct}, and \Ref{lemliebracketonthemoduleofvectorfieldspeq1},
\begin{align}\label{lemliebracketonthemoduleofvectorfieldspeq2}
&\Foreach{\cg}{\mapdifclass{\infty}{\Man{}}{\RR}}\cr
&\begin{aligned}
&\hskip0.5\baselineskip\func{\[\cmp{\func{\Lieder{\Man{}}}{\avecf{1}}}{\func{\Lieder{\Man{}}}{\cf\smoothvfprod{\Man{}}\avecf{2}}}\]}{\cg}\cr
&=\func{\[\func{\Lieder{\Man{}}}{\avecf{1}}\]}{\func{\[\func{\Lieder{\Man{}}}{\cf\smoothvfprod{\Man{}}\avecf{2}}\]}{\cg}}\cr
&=\func{\[\func{\Lieder{\Man{}}}{\avecf{1}}\]}{\func{\[\cf\smoothderprod{\Man{}}\func{\Lieder{\Man{}}}{\avecf{2}}\]}{\cg}}\cr
&=\func{\[\func{\Lieder{\Man{}}}{\avecf{1}}\]}{\cf\smoothfprod{\Man{}}\(\func{\[\func{\Lieder{\Man{}}}{\avecf{2}}\]}{\cg}\)}\cr
&=\cf\smoothfprod{\Man{}}\bigg(\func{\[\func{\Lieder{\Man{}}}{\avecf{1}}\]}{\func{\[\func{\Lieder{\Man{}}}{\avecf{2}}\]}{\cg}}\bigg)\smoothfsum{\Man{}}
\bigg(\func{\[\func{\Lieder{\Man{}}}{\avecf{1}}\]}{\cf}\bigg)
\smoothfprod{\Man{}}\bigg(\func{\[\func{\Lieder{\Man{}}}{\avecf{2}}\]}{\cg}\bigg),
\end{aligned}
\end{align}
and,
\begin{align}\label{lemliebracketonthemoduleofvectorfieldspeq3}
&\Foreach{\cg}{\mapdifclass{\infty}{\Man{}}{\RR}}\cr
&\begin{aligned}
&\hskip0.5\baselineskip\func{\[\cmp{\func{\Lieder{\Man{}}}{\cf\smoothvfprod{\Man{}}\avecf{2}}}{\func{\Lieder{\Man{}}}{\avecf{1}}}\]}{\cg}\cr
&=\func{\[\cf\smoothderprod{\Man{}}\func{\Lieder{\Man{}}}{\avecf{2}}\]}{\func{\[\func{\Lieder{\Man{}}}{\avecf{1}}\]}{\cg}}\cr
&=\cf\smoothfprod{\Man{}}\bigg(\func{\[\func{\Lieder{\Man{}}}{\avecf{2}}\]}{\func{\[\func{\Lieder{\Man{}}}{\avecf{1}}\]}{\cg}}\bigg),
\end{aligned}
\end{align}
and therefore using the distributive properties of operations of a ring (\refdef{defring}),
\begin{align}\label{lemliebracketonthemoduleofvectorfieldspeq4}
&\Foreach{\cg}{\mapdifclass{\infty}{\Man{}}{\RR}}\cr
&\begin{aligned}
&\hskip0.5\baselineskip\func{\[\cmp{\func{\Lieder{\Man{}}}{\avecf{1}}}{\func{\Lieder{\Man{}}}{\cf\smoothvfprod{\Man{}}\avecf{2}}}\smoothdersub{\Man{}}
\cmp{\func{\Lieder{\Man{}}}{\cf\smoothvfprod{\Man{}}\avecf{2}}}{\func{\Lieder{\Man{}}}{\avecf{1}}}\]}{\cg}\cr
&=~\cf\smoothfprod{\Man{}}\bigg(\func{\[\func{\Lieder{\Man{}}}{\avecf{1}}\]}{\func{\[\func{\Lieder{\Man{}}}{\avecf{2}}\]}{\cg}}\smoothfsub{\Man{}}
\func{\[\func{\Lieder{\Man{}}}{\avecf{2}}\]}{\func{\[\func{\Lieder{\Man{}}}{\avecf{1}}\]}{\cg}}\bigg)\cr
&~~\smoothfsum{\Man{}}
\bigg(\func{\[\func{\Lieder{\Man{}}}{\avecf{1}}\]}{\cf}\bigg)
\smoothfprod{\Man{}}\bigg(\func{\[\func{\Lieder{\Man{}}}{\avecf{2}}\]}{\cg}\bigg)\cr
&=~\func{\bigg(\cf\smoothderprod{\Man{}}\[\cmp{\func{\Lieder{\Man{}}}{\avecf{1}}}{\func{\Lieder{\Man{}}}{\avecf{2}}}\smoothdersub{\Man{}}
\cmp{\func{\Lieder{\Man{}}}{\avecf{2}}}{\func{\Lieder{\Man{}}}{\avecf{1}}}\]\bigg)}{\cg}\cr
&~~\smoothfsum{\Man{}}
\func{\bigg[\bigg(\func{\[\func{\Lieder{\Man{}}}{\avecf{1}}\]}{\cf}\bigg)
\smoothderprod{\Man{}}\[\func{\Lieder{\Man{}}}{\avecf{2}}\]\bigg]}{\cg}\cr
&=~
\func{\left\{\bigg(\cf\smoothderprod{\Man{}}\[\cmp{\func{\Lieder{\Man{}}}{\avecf{1}}}{\func{\Lieder{\Man{}}}{\avecf{2}}}\smoothdersub{\Man{}}
\cmp{\func{\Lieder{\Man{}}}{\avecf{2}}}{\func{\Lieder{\Man{}}}{\avecf{1}}}\]\bigg)\right.\cr
&~~\left.\smoothdersum{\Man{}}
\bigg[\bigg(\func{\[\func{\Lieder{\Man{}}}{\avecf{1}}\]}{\cf}\bigg)
\smoothderprod{\Man{}}\[\func{\Lieder{\Man{}}}{\avecf{2}}\]\bigg]\right\}}{\cg}.
\end{aligned}
\end{align}
Therefore,
\begin{align}\label{lemliebracketonthemoduleofvectorfieldspeq5}
&\hskip0.5\baselineskip\[\cmp{\func{\Lieder{\Man{}}}{\avecf{1}}}{\func{\Lieder{\Man{}}}{\cf\smoothvfprod{\Man{}}\avecf{2}}}\smoothdersub{\Man{}}
\cmp{\func{\Lieder{\Man{}}}{\cf\smoothvfprod{\Man{}}\avecf{2}}}{\func{\Lieder{\Man{}}}{\avecf{1}}}\]\cr
&=~\bigg(\cf\smoothderprod{\Man{}}\[\cmp{\func{\Lieder{\Man{}}}{\avecf{1}}}{\func{\Lieder{\Man{}}}{\avecf{2}}}\smoothdersub{\Man{}}
\cmp{\func{\Lieder{\Man{}}}{\avecf{2}}}{\func{\Lieder{\Man{}}}{\avecf{1}}}\]\bigg)\cr
&~~\smoothdersum{\Man{}}
\bigg[\bigg(\func{\[\func{\Lieder{\Man{}}}{\avecf{1}}\]}{\cf}\bigg)
\smoothderprod{\Man{}}\[\func{\Lieder{\Man{}}}{\avecf{2}}\]\bigg].
\end{align}
According to \refthm{thmliederivativeisamodulehomomorphism} and \refdef{defmodulehomomorphism},
$\finv{\(\Lieder{\Man{}}\)}$ is an $\smoothring{\infty}{\Man{}}$-module-isomorphism from $\smoothdermodule{\infty}{\Man{}}$ to
$\smoothvfmodule{\infty}{\Man{}}$, and hence according to \Ref{lemliebracketonthemoduleofvectorfieldspeq5} and by invoking
the definition of Lie-bracket of smooth vector-fields stated in \refdef{defliebracket},
\begin{align}
\liebracket{\avecf{1}}{\cf\smoothvfprod{\Man{}}\avecf{2}}{\Man{}}&=
\func{\finv{\Lieder{\Man{}}}}{\cmp{\func{\Lieder{\Man{}}}{\avecf{1}}}{\func{\Lieder{\Man{}}}
{\cf\smoothvfprod{\Man{}}\avecf{2}}}\smoothdersub{\Man{}}
\cmp{\func{\Lieder{\Man{}}}{\cf\smoothvfprod{\Man{}}\avecf{2}}}{\func{\Lieder{\Man{}}}{\avecf{1}}}}\cr
&=~\func{\finv{\Lieder{\Man{}}}}{\cf\smoothderprod{\Man{}}\[\cmp{\func{\Lieder{\Man{}}}{\avecf{1}}}
{\func{\Lieder{\Man{}}}{\avecf{2}}}\smoothdersub{\Man{}}
\cmp{\func{\Lieder{\Man{}}}{\avecf{2}}}{\func{\Lieder{\Man{}}}{\avecf{1}}}\]}\cr
&~~\smoothvfsum{\Man{}}
\func{\finv{\Lieder{\Man{}}}}{\bigg(\func{\[\func{\Lieder{\Man{}}}{\avecf{1}}\]}{\cf}\bigg)
\smoothderprod{\Man{}}\[\func{\Lieder{\Man{}}}{\avecf{2}}\]}\cr
&=~\cf\smoothvfprod{\Man{}}\[\func{\finv{\Lieder{\Man{}}}}{\[\cmp{\func{\Lieder{\Man{}}}{\avecf{1}}}
{\func{\Lieder{\Man{}}}{\avecf{2}}}\smoothdersub{\Man{}}
\cmp{\func{\Lieder{\Man{}}}{\avecf{2}}}{\func{\Lieder{\Man{}}}{\avecf{1}}}\]}\]\cr
&~~\smoothvfsum{\Man{}}
\bigg(\func{\[\func{\Lieder{\Man{}}}{\avecf{1}}\]}{\cf}\bigg)\smoothvfprod{\Man{}}
\[\func{\finv{\Lieder{\Man{}}}}{\[\func{\Lieder{\Man{}}}{\avecf{2}}\]}\]\cr
&=\cf\smoothvfprod{\Man{}}\liebracket{\avecf{1}}{\avecf{2}}{\Man{}}\smoothvfsum{\Man{}}
\bigg(\func{\[\func{\Lieder{\Man{}}}{\avecf{1}}\]}{\cf}\bigg)\smoothvfprod{\Man{}}\avecf{2}.
\end{align}
\endthm
\theorem\label{thmliebracketonthemoduleofvectorfields1}
\begin{align}
&\Foreach{\opair{\cf_1}{\cf_2}}{\Cprod{\mapdifclass{\infty}{\Man{}}{\RR}}{\mapdifclass{\infty}{\Man{}}{\RR}}}
\Foreach{\opair{\avecf{1}}{\avecf{2}}}{\Cprod{\vecf{\Man{}}{\infty}}{\vecf{\Man{}}{\infty}}}\cr
&\begin{aligned}
\liebracket{\cf_1\smoothvfprod{\Man{}}\avecf{1}}{\cf_2\smoothvfprod{\Man{}}\avecf{2}}{\Man{}}
&=\hskip0.5\baselineskip\(\cf_1\smoothfprod{\Man{}}\cf_2\)\smoothvfprod{\Man{}}\liebracket{\avecf{1}}{\avecf{2}}{\Man{}}\cr
&\hskip0.5\baselineskip\smoothvfsum{\Man{}}
\[\cf_1\smoothfprod{\Man{}}\bigg(\func{\[\func{\Lieder{\Man{}}}{\avecf{1}}\]}{\cf_2}\bigg)\]\smoothvfprod{\Man{}}\avecf{2}\cr
&\hskip0.5\baselineskip\smoothvfsub{\Man{}}
\[\cf_2\smoothfprod{\Man{}}\bigg(\func{\[\func{\Lieder{\Man{}}}{\avecf{2}}\]}{\cf_1}\bigg)\]\smoothvfprod{\Man{}}\avecf{1}.
\end{aligned}
\end{align}
\proof
According to \reflem{lemliebracketonthemoduleofvectorfields},
\refthm{thmliealgebraofvectorfields0},
\refdef{defliealgebra}, and \refdef{defsrealvaluedmoothmapsproduct}, it is straightforward.
\endthm
\subsection{Submodules of Vector-Fields Module Induced by Distributions}
\theorem\label{thmsubmoduleinducedbydistribution}
$\distribution{}$ is taken as an element of $\Distributions{r}{\Man{}}$ (an $r$-dimensional smooth distribution of $\Man{}$),
for some $r$ in $\Zp$.
The set of all smooth vector-fields on $\Man{}$ consistent with the smooth distribution $\distribution{}$ of $\Man{}$ is a submodule
of the $\smoothring{\infty}{\Man{}}$-module of smooth vector-fields on $\Man{}$. That is,
\begin{equation}
\distvecf{\Man{}}{\distribution{}}{r}\in
\Submodules{\smoothvfmodule{\infty}{\Man{}}}.
\end{equation}
\proof
According to \refdef{defconsistentvectorfieldswithdistributions},
\begin{equation}\label{thmsubmoduleinducedbydistributionpeq1}
\distvecf{\Man{}}{\distribution{}}{r}=
\defset{\avecf{}}{\vecf{\Man{}}{\infty}}
{\big[\Foreach{\point}{\M{}}\func{\avecf{}}{\point}\in\func{\distribution{}}{\point}\big]}.
\end{equation}
Thus according to \refdef{defdistribution},
it is evident that the neutral element of the vector-space $\Vecf{\Man{}}{\infty}$ is in $\distvecf{\Man{}}{\distribution{}}{r}$. That is,
\begin{equation}
\zerovec{\Vecf{\Man{}}{\infty}}\in
\distvecf{\Man{}}{\distribution{}}{r},
\end{equation}
and hence $\distvecf{\Man{}}{\distribution{}}{r}$ is non-empty.
\begin{itemize}
\item[$\pr{1}$]
Each $\avecf{}$ abd $\avecf{1}$ is taken as an arbitrary element of $\distvecf{\Man{}}{\distribution{}}{r}$, and $\cf$
as an arbitrary element of $\mapdifclass{\infty}{\Man{}}{\RR}$.
According to \Ref{thmsubmoduleinducedbydistributionpeq1}, and considering that for every point $\point$ of $\Man{}$
the value of the $\infty$-distribution at $\point$ is a vector-subspace of $\tanspace{\point}{\Man{}}$ (\refdef{defdistribution}),
it is clear that,
\begin{align}
\Foreach{\point}{\M{}}
\left\{\begin{aligned}
&\func{\avecf{}}{\point}+\func{\avecf{1}}{\point}\in\func{\distribution{}}{\point},\cr
&\func{\cf}{\point}\func{\avecf{}}{\point}\in\func{\distribution{}}{\point},
\end{aligned}\right.
\end{align}
and consequently, according to \refdef{defsrealvaluedmoothmapsproduct} and \Ref{thmsubmoduleinducedbydistributionpeq1},
\begin{align}
\avecf{}+\avecf{1}&\in\distvecf{\Man{}}{\distribution{}}{r},\\
\cf\avecf{}&\in\distvecf{\Man{}}{\distribution{}}{r}.
\end{align}
\endp
\end{itemize}
Therefore, $\distvecf{\Man{}}{\distribution{}}{r}$ is a non-empty subset of $\vecf{\Man{}}{\infty}$ that is closed under the
addition and scalar-product operations of the $\smoothring{\infty}{\Man{}}$-module $\smoothvfmodule{\infty}{\Man{}}$,
and hence according to \refthm{thmsubmodule1}, is a submodule of $\smoothvfmodule{\infty}{\Man{}}$.
\endthm
\chapteR{
Abstract Smooth Group}
\thispagestyle{fancy}
\section{Basic Structure of a Smooth Group}
\definition\label{defliegroup}
$\G{}$ is taken as a set, $\gop{}$ as an element of $\Func{\Cprod{\G{}}{\G{}}}{\G{}}$ (a binary operation on $\G{}$),
and $\maxatlas{}$ as an element of $\maxatlases{\infty}{\G{}}{\R^n}$ (a maximal-atlas of differentiablity class $\difclass{\infty}$
on $\G{}$ constructed upon the Banach-space $\R^n$) for some positive integer $n$.
The triple $\triple{\G{}}{\gop{}}{\maxatlas{}}$ is referred to as a $\quotl$smooth group (of dimension $n$)$\quotr$ or
a $\quotl$real Lie-Group (of dimension $n$)$\quotr$ iff these properties are satisfied.
\begin{itemize}
\item[\myitem{LG~1.}]
$\Group{}:=\opair{\G{}}{\gop{}}$ is a group.
\item[\myitem{LG~2.}]
$\opair{\G{}}{\maxatlas{}}$ is a manifold. This means that $\G{}$ endowd with the topology on $\G{}$ induced by the maximal-atlas
$\maxatlas{}$, that is $\mantop{\opair{\G{}}{\maxatlas{}}}$, is a Hausdorff and second-countable topological-space.
\item[\myitem{LG~3.}]
$\gop{}$ is a smooth map from the manifold-product
$\manprod{\opair{\G{}}{\maxatlas{}}}{\opair{\G{}}{\maxatlas{}}}$ to the manifold $\opair{\G{}}{\maxatlas{}}$. That is,
\begin{equation}
\gop{}\in\mapdifclass{\infty}{\manprod{\opair{\G{}}{\maxatlas{}}}{\opair{\G{}}{\maxatlas{}}}}{\opair{\G{}}{\maxatlas{}}}.
\end{equation}
\item[\myitem{LG~4.}]
$\ginv{\Group{}}$ is a smooth map from $\opair{\G{}}{\maxatlas{}}$ to $\opair{\G{}}{\maxatlas{}}$. That is,
\begin{equation}
\ginv{\Group{}}\in\mapdifclass{\infty}{\opair{\G{}}{\maxatlas{}}}{\opair{\G{}}{\maxatlas{}}}.
\end{equation}
\end{itemize}
\endef
\theorem\label{thmliegroupequiv0}
$\G{}$ is taken as a set, $\gop{}$ as an element of $\Func{\Cprod{\G{}}{\G{}}}{\G{}}$,
and $\maxatlas{}$ as an element of $\maxatlases{\infty}{\G{}}{\R^n}$ for some positive integer $n$, such that
$\opair{\G{}}{\gop{}}$ is a group and $\opair{\G{}}{\maxatlas{}}$ is a manifold.
$\triple{\G{}}{\gop{}}{\maxatlas{}}$ is a smooth group if and only if $\gopr{\gop{}}{\Group{}}$
is a smooth map, that is,
\begin{equation}
\gopr{\gop{}}{\Group{}}\in\mapdifclass{\infty}{\manprod{\opair{\G{}}{\maxatlas{}}}{\opair{\G{}}{\maxatlas{}}}}{\opair{\G{}}{\maxatlas{}}}.
\end{equation}
\proof
Let $\IG{}$ denote the identity element of the group $\opair{\G{}}{\gop{}}$.
It can be easily verified that,
\begin{align}
\gopr{\gop{}}{\Group{}}&=\cmp{\gop{}}{\(\identity{\G{}}\times\ginv{\Group{}}\)},
\label{thmliegroupequiv0peq1}\\
\ginv{\Group{}}&=\cmp{\cmp{\gopr{\gop{}}{\Group{}}}{\(\constmap{\G{}}{\IG{}}\times\identity{\G{}}\)}}{\diagmap{\G{}}},
\label{thmliegroupequiv0peq2}\\
\gop{}&=\cmp{\gopr{\gop{}}{\Group{}}}{\(\identity{\G{}}\times\ginv{\Group{}}\)},
\label{thmliegroupequiv0peq3}
\end{align}
where $\function{\constmap{\G{}}{\IG{}}}{\G{}}{\G{}}$ denotes the constnant map on $\G{}$ with single value $\IG{}$,
$\function{\diagmap{\G{}}}{\G{}}{\Cprod{\G{}}{\G{}}}$ denotes the diagonal map on $\G{}$ sending
each $\g{}$ in $\G{}$ to $\opair{\g{}}{\g{}}$, and $\function{\identity{\G{}}}{\G{}}{\G{}}$ denotes the
identity map on $\G{}$ sending each $\g{}$ in $\G{}$ to $\g{}$. All mappings $\constmap{\G{}}{\IG{}}$,
$\identity{\G{}}$, and $\diagmap{\G{}}$ are smooth. That is,
\begin{align}
\constmap{\G{}}{\IG{}}\in\mapdifclass{\infty}{\opair{\G{}}{\maxatlas{}}}{\opair{\G{}}{\maxatlas{}}},
\label{thmliegroupequiv0peq4}\\
\identity{\G{}}\in\mapdifclass{\infty}{\opair{\G{}}{\maxatlas{}}}{\opair{\G{}}{\maxatlas{}}},
\label{thmliegroupequiv0peq5}\\
\diagmap{\G{}}\in\mapdifclass{\infty}{\opair{\G{}}{\maxatlas{}}}{\manprod{\opair{\G{}}{\maxatlas{}}}{\opair{\G{}}{\maxatlas{}}}}.
\label{thmliegroupequiv0peq6}
\end{align}
\begin{itemize}
\item[$\pr{1}$]
It is assumed that $\triple{\G{}}{\gop{}}{\maxatlas{}}$ is a real Lie-group. So according to \refdef{defliegroup},
\begin{align}
\gop{}&\in\mapdifclass{\infty}{\manprod{\opair{\G{}}{\maxatlas{}}}{\opair{\G{}}{\maxatlas{}}}}{\opair{\G{}}{\maxatlas{}}},\\
\ginv{\Group{}}&\in\mapdifclass{\infty}{\opair{\G{}}{\maxatlas{}}}{\opair{\G{}}{\maxatlas{}}}.
\end{align}
So since $\identity{\G{}}\in\mapdifclass{\infty}{\opair{\G{}}{\maxatlas{}}}{\opair{\G{}}{\maxatlas{}}}$,
obviously,
\begin{equation}
\(\identity{\G{}}\times\ginv{\Group{}}\)\in
\mapdifclass{\infty}{\manprod{\opair{\G{}}{\maxatlas{}}}{\opair{\G{}}{\maxatlas{}}}}
{\manprod{\opair{\G{}}{\maxatlas{}}}{\opair{\G{}}{\maxatlas{}}}},
\end{equation}
and hence,
\begin{equation}
\cmp{\gop{}}{\(\identity{\G{}}\times\ginv{\Group{}}\)}
\in\mapdifclass{\infty}{\manprod{\opair{\G{}}{\maxatlas{}}}{\opair{\G{}}{\maxatlas{}}}}{\opair{\G{}}{\maxatlas{}}},
\end{equation}
which according to \Ref{thmliegroupequiv0peq1} means,
\begin{equation}
\gopr{\gop{}}{\Group{}}
\in\mapdifclass{\infty}{\manprod{\opair{\G{}}{\maxatlas{}}}{\opair{\G{}}{\maxatlas{}}}}{\opair{\G{}}{\maxatlas{}}}.
\end{equation}
\endp
\end{itemize}
\begin{itemize}
\item[$\pr{2}$]
It is now assumed that $\gopr{\gop{}}{\Group{}}
\in\mapdifclass{\infty}{\manprod{\opair{\G{}}{\maxatlas{}}}{\opair{\G{}}{\maxatlas{}}}}{\opair{\G{}}{\maxatlas{}}}$.
Then according to \Ref{thmliegroupequiv0peq2}, \Ref{thmliegroupequiv0peq4}, \Ref{thmliegroupequiv0peq5}, and
\Ref{thmliegroupequiv0peq6},
\begin{equation}
\ginv{\Group{}}\in\mapdifclass{\infty}{\opair{\G{}}{\maxatlas{}}}{\opair{\G{}}{\maxatlas{}}},
\end{equation}
and accordingly, \Ref{thmliegroupequiv0peq3} and \Ref{thmliegroupequiv0peq5} imply,
\begin{equation}
\gop{}\in\mapdifclass{\infty}{\manprod{\opair{\G{}}{\maxatlas{}}}{\opair{\G{}}{\maxatlas{}}}}{\opair{\G{}}{\maxatlas{}}}.
\end{equation}
Therefore $\triple{\G{}}{\gop{}}{\maxatlas{}}$ is a real Lie-group.
\endp
\end{itemize}
\endthm
\definition\label{deftopologyofliegroup}
$\Liegroup{}=\triple{\G{}}{\gop{}}{\maxatlas{}}$ is taken to be a smooth group.
\begin{itemize}
\item[$\centerdot$]
The intrinsic (or underlying) group structure of $\Liegroup{}$ is denoted by $\LieG{\Liegroup{}}$. That is,
$\LieG{\Liegroup{}}:=\opair{\G{}}{\gop{}}$.
\item[$\centerdot$]
The underlying manifold of $\Liegroup{}$ is denoted by $\Lieman{\Liegroup{}}$. That is,
$\Lieman{\Liegroup{}}:=\opair{\G{}}{\maxatlas{}}$.
\item[$\centerdot$]
The topology on $\G{}$ induced by the maximal-atlas $\maxatlas{}$ is denoted by $\lietop{\Liegroup{}}$
which is referred to as the $\quotl$intrinsic (or underlying) topology of the smooth group $\Liegroup{}$$\quotr$.
The topological-space corresponded with the topology $\lietop{\Liegroup{}}$ is denoted by $\lietops{\Liegroup{}}$, that is,
$\lietops{\Liegroup{}}:=\opair{\G{}}{\lietop{\Liegroup{}}}$, which is called the $\quotl$intrinsic topological-structure
of the smooth group $\Liegroup{}$$\quotr$.
\item[$\centerdot$]
Each element of $\lietop{\Liegroup{}}$, that is each open set of the topological-space $\lietops{\Liegroup{}}$,
is called an $\quotl$open set of the smooth group $\Liegroup{}$$\quotr$.
\item[$\centerdot$]
Each closed set of the topological-space $\lietops{\Liegroup{}}$,
is called a $\quotl$closed set of the smooth group $\Liegroup{}$$\quotr$.
\item[$\centerdot$]
Each element of $\G{}$ is called a $\quotl$point of the smooth group $\Liegroup{}$$\quotr$. The set of all points of the
smooth group $\Liegroup{}$, that is $\G{}$, can alternatively denoted by $\Liepoints{\Liegroup{}}$. 
\end{itemize}
\endef
\fixed
\begin{itemize}
\item[$\centerdot$]
$\Liegroup{}=\triple{\G{}}{\gop{}}{\maxatlas{}}$ is fixed as a smooth group of dimension $n$,
for some positive integer $n$. $\IG{}$ is defined to be the identity element of the group $\LieG{\Liegroup{}}$.
\item[$\centerdot$]
$\Liegroup{1}=\triple{\G{1}}{\gop{1}}{\maxatlas{1}}$ is fixed as a smooth group of dimension $n_1$,
for some positive integer $n_1$. $\IG{1}$ is defined to be the identity element of the group $\LieG{\Liegroup{1}}$.
\item[$\centerdot$]
$\Liegroup{2}=\triple{\G{2}}{\gop{2}}{\maxatlas{2}}$ is fixed as a smooth group of dimension $n_2$,
for some positive integer $n_2$. $\IG{2}$ is defined to be the identity element of the group $\LieG{\Liegroup{2}}$.
\item[$\centerdot$]
$\Man{}=\opair{\M{}}{\maxatlas{0}}$ is fixed as an $n_0$-dimensional and $\difclass{\infty}$ manifold
modeled on the Banach-space $\R^{n_0}$.
\end{itemize}
\endfixed
\theorem\label{thminversemappingisdiffeomorphism}
The inverse-mapping of the group $\LieG{\Liegroup{}}$ is an
$\infty$-automorphisms of the underlying manifold $\Lieman{\Liegroup{}}$ of $\Liegroup{}$. That is,
\begin{equation}
\ginv{\LieG{\Liegroup{}}}\in\Diff{\infty}{\Lieman{\Liegroup{}}}.
\end{equation}
\proof
According to \refdef{defliegroup}, $\ginv{\LieG{\Liegroup{}}}$ is a smooth from $\Lieman{\Liegroup{}}$
to $\Lieman{\Liegroup{}}$.
\begin{equation}\label{thminversemappingisadiffeomorphismpeq1}
\ginv{\LieG{\Liegroup{}}}\in\mapdifclass{\infty}{\Lieman{\Liegroup{}}}{\Lieman{\Liegroup{}}}.
\end{equation}
Moreover, it is well-known that $\ginv{\LieG{\Liegroup{}}}$ is a bijective mapping
from $\G{}$ to $\G{}$, and coincides with its own inverse, that is,
\begin{equation}
\finv{\ginv{\LieG{\Liegroup{}}}}=\ginv{\LieG{\Liegroup{}}}.
\end{equation}
Thus, in addition to \Ref{thminversemappingisadiffeomorphismpeq1},
\begin{equation}
\ginv{\LieG{\Liegroup{}}}\in\mapdifclass{\infty}{\Lieman{\Liegroup{}}}{\Lieman{\Liegroup{}}},
\end{equation}
and consequently it becomes evident that $\ginv{\LieG{\Liegroup{}}}$ is an $\infty$-automorphism of the manifold
$\Lieman{\Liegroup{}}$.
\endthm
\theorem\label{thmtranslationsarediffeomorphisms}
For every point $\g{}$ of $\Liegroup{}$, the left-translation and right-translation of $\LieG{\Liegroup{}}$ by $\g{}$, and
the $\g{}$-conjugation of $\LieG{\Liegroup{}}$
are $\infty$-automorphisms of the underlying manifold $\Lieman{\Liegroup{}}$ of the smooth group $\Liegroup{}$. That is,
\begin{align}
\Foreach{\g{}}{\G{}}
\begin{cases}
\gltrans{\LieG{\Liegroup{}}}{\g{}}\in\Diff{\infty}{\Lieman{\Liegroup{}}},\cr
\grtrans{\LieG{\Liegroup{}}}{\g{}}\in\Diff{\infty}{\Lieman{\Liegroup{}}},\cr
\gconj{\LieG{\Liegroup{}}}{\g{}}\in\Diff{\infty}{\Lieman{\Liegroup{}}}.
\end{cases}
\end{align}
\proof
According to \refdef{defgrouptranslations} and considering the smoothness of the mappings $\constmap{\G{}}{\g{}}$ for all $\g{}$ in $\G{}$,
$\identity{\G{}}$, and $\diagmap{\G{}}$, it is clear that,
\begin{equation}\label{thmtranslationsarediffeomorphismspeq1}
\Foreach{\g{}}{\G{}}
\gltrans{\LieG{\Liegroup{}}}{\g{}}\in\mapdifclass{\infty}
{\Lieman{\Liegroup{}}}{\Lieman{\Liegroup{}}}.
\end{equation}
In addition, it is known that for every $\g{}$ in $\G{}$, $\gltrans{\LieG{\Liegroup{}}}{\g{}}$ is a bijective mapping
and,
\begin{equation}\label{thmtranslationsarediffeomorphismspeq2}
\Foreach{\g{}}{\G{}}
\finv{\(\gltrans{\LieG{\Liegroup{}}}{\g{}}\)}=\gltrans{\LieG{\Liegroup{}}}{\invg{\g{}}{}}.
\end{equation}
\Ref{thmtranslationsarediffeomorphismspeq1} and \Ref{thmtranslationsarediffeomorphismspeq2} imply that,
\begin{align}
\Foreach{\g{}}{\G{}}
\gltrans{\LieG{\Liegroup{}}}{\g{}}\in\mapdifclass{\infty}
{\Lieman{\Liegroup{}}}{\Lieman{\Liegroup{}}},~
\finv{\(\gltrans{\LieG{\Liegroup{}}}{\g{}}\)}\in\mapdifclass{\infty}
{\Lieman{\Liegroup{}}}{\Lieman{\Liegroup{}}}.
\end{align}
This means $\gltrans{\LieG{\Liegroup{}}}{\g{}}$ is an $\infty$-automorphism of the manifold $\Lieman{\Liegroup{}}$,
for all $\g{}$ in $\G{}$. That is,
\begin{equation}
\Foreach{\g{}}{\G{}}
\gltrans{\LieG{\Liegroup{}}}{\g{}}\in\Diff{\infty}{\Lieman{\Liegroup{}}}.
\end{equation}
In an obviously similar way it can be seen that,
\begin{equation}
\Foreach{\g{}}{\G{}}
\grtrans{\LieG{\Liegroup{}}}{\g{}}\in\Diff{\infty}{\Lieman{\Liegroup{}}}.
\end{equation}
Furthermore, since $\gltrans{\LieG{\Liegroup{}}}{\g{}}$ and $\grtrans{\LieG{\Liegroup{}}}{\g{}}$
are $\infty$-automorphisms of $\Lieman{\Liegroup{}}$ and
$\gconj{\LieG{\Liegroup{}}}{\g{}}=\cmp{\gltrans{\LieG{\Liegroup{}}}{\g{}}}{\grtrans{\Group{}}{\g{}}}$, for all
$\g{}$ is in $\G{}$, and considering that the composition of any pair of $\infty$-automorphisms of $\Lieman{\Liegroup{}}$
is again an $\infty$-automorphism of $\Lieman{\Liegroup{}}$, it is evident that,
\begin{equation}
\Foreach{\g{}}{\G{}}
\gconj{\LieG{\Liegroup{}}}{\g{}}\in\Diff{\infty}{\Lieman{\Liegroup{}}}.
\end{equation}
\endthm
\theorem\label{thmliegroupisatopologicalgroup}
The intrinsic group structure of the smooth group $\Liegroup{}$ endowed with its intrinsic topology
is a topological group. That is, the triple $\triple{\G{}}{\gop{}}{\lietop{\Liegroup{}}}$ is a topological group.
\proof
The smoothness of the group operation of $\LieG{\Liegroup{}}$ from the product manifold
$\manprod{\Lieman{\Liegroup{}}}{\Lieman{\Liegroup{}}}$ to the manifold $\Lieman{\Liegroup{}}$ trivially implies
its continuity from the underlying topological-space of $\manprod{\Lieman{\Liegroup{}}}{\Lieman{\Liegroup{}}}$
to the underlying topological-space of $\Lieman{\Liegroup{}}$, and since the underlying topological-space of
$\manprod{\Lieman{\Liegroup{}}}{\Lieman{\Liegroup{}}}$ is $\topprod{\lietops{\Liegroup{}}}{\lietops{\Liegroup{}}}$,
\begin{equation}
\gop\in\CF{\topprod{\lietops{\Liegroup{}}}{\lietops{\Liegroup{}}}}{\lietops{\Liegroup{}}}.
\end{equation}
Moreover, since the inverse mapping of $\LieG{\Liegroup{}}$ is a smooth map from $\Lieman{\Liegroup{}}$ to $\Lieman{\Liegroup{}}$,
obviously it is a continuous map from the topological-space $\lietops{\Liegroup{}}$ to itself. That is,
\begin{equation}
\ginv{\LieG{\Liegroup{}}}\in\CF{\lietops{\Liegroup{}}}{\lietops{\Liegroup{}}}.
\end{equation}
Therefore, since $\lietops{\Liegroup{}}=\opair{\G{}}{\lietop{\Liegroup{}}}$, according to \refdef{deftopologicalgroup},
the triple $\triple{\G{}}{\gop{}}{\lietop{\Liegroup{}}}$ is a topological group.
\endthm
\definition\label{deftopologicalgroupstructureofliegroup}
The topological group $\triple{\G{}}{\gop{}}{\lietop{\Liegroup{}}}$ is denoted by $\Lietopg{\Liegroup{}}$, and
is called the $\quotl$underlying (intrinsic) topological group structure of the smooth group $\Liegroup{}$$\quotr$.
\endef
\section{Lie-Algebra of a Smooth Group}
\definition\label{defleftinvariantvectorfields}
$\Leftinvvf{\Liegroup{}}$ is defined as the set of all mappings $\function{\avecf{}}{\G{}}{\tanbun{\Lieman{\Liegroup{}}}}$
such that $\avecf{}$ sends every point of $\Lieman{\Liegroup{}}$ to the tangent-space of $\Lieman{\Liegroup{}}$ at that point,
and for every $\g{}$ in $\G{}$, the differential of the left-translation of $\LieG{\Liegroup{}}$ by $\g{}$ maps $\func{\avecf{}}{\point}$
to $\func{\avecf{}}{\g{}\gop{}\point}$ for each point $\point$ of $\Liegroup{}$. That is,
\begin{align}
\Leftinvvf{\Liegroup{}}:=\defset{\avecf{}}{\Func{\G{}}{\tanbun{\Lieman{\Liegroup{}}}}}
{\(\cmp{\basep{\Lieman{\Liegroup{}}}}{\avecf{}}=\identity{\G{}},~
\Foreach{\g{}}{\G{}}\cmp{\(\der{\gltrans{\LieG{\Liegroup{}}}{\g{}}}{\Lieman{\Liegroup{}}}
{\Lieman{\Liegroup{}}}\)}{\avecf{}}=\cmp{\avecf{}}{\gltrans{\LieG{\Liegroup{}}}{\g{}}}\)}.
\end{align}
Each element of $\Leftinvvf{\Liegroup{}}$ is referred to as a
$\quotl$left-invariant vector-field on the Lie-group $\Liegroup{}$$\quotr$.\\
Illustrating with commutative diagrams, a mapping $\function{\avecf{}}{\G{}}{\tanbun{\Lieman{\Liegroup{}}}}$
is defined to be a left-invariant vector-field on $\Liegroup{}$ iff for every $\g{}$ in $\G{}$, the following
diagrams commute.
\begin{center}
\vskip0.5\baselineskip
\hskip-2\baselineskip
\begin{tikzcd}[row sep=6em, column sep=6em]
& \G{}
\arrow[r,"\avecf{}" description]
\arrow[swap,d,"\gltrans{\LieG{\Liegroup{}}}{\g{}}" description]
& \tanbun{\Lieman{\Liegroup{}}}
\arrow[d,"\derr{\gltrans{\LieG{\Liegroup{}}}{\g{}}}" description] \\
& \G{}
\arrow[swap,r,"\avecf{}" description] & \tanbun{\Lieman{\Liegroup{}}}
\end{tikzcd}
\hskip\baselineskip
\begin{tikzcd}[row sep=6em, column sep=6em]
\G{}
\arrow[r,"\avecf{}" description]
\arrow[swap,rd,"\identity{\G{}}" description]
& \tanbun{\Lieman{\Liegroup{}}}
\arrow[d,"\basep{\Lieman{\Liegroup{}}}" description]\\
& \G{}
\end{tikzcd}
\end{center}
\endef
\lemma\label{lemleftinvariantvectorfieldsequiv0}
For every $\avecf{}$ in $\Func{\G{}}{\tanbun{\Lieman{\Liegroup{}}}}$,
$\avecf{}$ is a left-invariant vector-field on $\Liegroup{}$
if and only if $\func{\avecf{}}{\IG{}}$ belongs to the tangent-space of $\Lieman{\Liegroup{}}$
at $\IG{}$ ($\tanspace{\IG{}}{\Lieman{\Liegroup{}}}$) and
for every point $\g{}$ of $\Liegroup{}$, the differential of $\gltrans{\LieG{\Liegroup{}}}{\g{}}$
sends $\func{\avecf{}}{\IG{}}$ to $\func{\avecf{}}{\g{}}$.
That is,
\begin{align}
\Leftinvvf{\Liegroup{}}=\defset{\avecf{}}{\Func{\G{}}{\tanbun{\Lieman{\Liegroup{}}}}}
{\(\func{\avecf{}}{\IG{}}\in\tanspace{\IG{}}{\Lieman{\Liegroup{}}},~
\Foreach{\g{}}{\G{}}\func{\(\der{\gltrans{\LieG{\Liegroup{}}}{\g{}}}{\Lieman{\Liegroup{}}}
{\Lieman{\Liegroup{}}}\)}{\func{\avecf{}}{\IG{}}}=\func{\avecf{}}{\g{}}\)}.
\end{align}
\proof
It is trivial that for every $\avecf{}$ in $\Func{\G{}}{\tanbun{\Lieman{\Liegroup{}}}}$,
\begin{equation}
\cmp{\basep{\Lieman{\Liegroup{}}}}{\avecf{}}=\identity{\G{}}
\then
\func{\avecf{}}{\IG{}}\in\tanspace{\IG{}}{\Lieman{\Liegroup{}}},
\end{equation}
and,
\begin{equation}
\[\Foreach{\g{}}{\G{}}\cmp{\(\der{\gltrans{\LieG{\Liegroup{}}}{\g{}}}{\Lieman{\Liegroup{}}}
{\Lieman{\Liegroup{}}}\)}{\avecf{}}=\cmp{\avecf{}}{\gltrans{\LieG{\Liegroup{}}}{\g{}}}\]
\then
\[\Foreach{\g{}}{\G{}}\func{\(\der{\gltrans{\LieG{\Liegroup{}}}{\g{}}}{\Lieman{\Liegroup{}}}
{\Lieman{\Liegroup{}}}\)}{\func{\avecf{}}{\IG{}}}=\func{\avecf{}}{\g{}}\].
\end{equation}
Thus,
\begin{equation}
\Leftinvvf{\Liegroup{}}\subseteq\defset{\avecf{}}{\Func{\G{}}{\tanbun{\Lieman{\Liegroup{}}}}}
{\(\func{\avecf{}}{\IG{}}\in\tanspace{\IG{}}{\Lieman{\Liegroup{}}},~
\Foreach{\g{}}{\G{}}\func{\(\der{\gltrans{\LieG{\Liegroup{}}}{\g{}}}{\Lieman{\Liegroup{}}}
{\Lieman{\Liegroup{}}}\)}{\func{\avecf{}}{\IG{}}}=\func{\avecf{}}{\g{}}\)}.
\end{equation}
\begin{itemize}
\item[$\pr{1}$]
$\avecf{}$ is taken as an element of $\Func{\G{}}{\tanbun{\Lieman{\Liegroup{}}}}$ such that,
\begin{equation}\label{lemleftinvariantvectorfieldsequiv0p1eq1}
\func{\avecf{}}{\IG{}}\in\tanspace{\IG{}}{\Lieman{\Liegroup{}}},
\end{equation}
and,
\begin{equation}\label{lemleftinvariantvectorfieldsequiv0p1eq2}
\Foreach{\g{}}{\G{}}\func{\(\der{\gltrans{\LieG{\Liegroup{}}}{\g{}}}{\Lieman{\Liegroup{}}}
{\Lieman{\Liegroup{}}}\)}{\func{\avecf{}}{\IG{}}}=\func{\avecf{}}{\g{}}.
\end{equation}
Exploiting the chain-rule of differentiation and the composition rule of left-translations of a group,
\Ref{lemleftinvariantvectorfieldsequiv0p1eq2} implies,
\begin{align}
\Foreach{\g{}}{\G{}}\Foreach{\point}{\G{}}
\func{\[\cmp{\(\der{\gltrans{\LieG{\Liegroup{}}}{\g{}}}{\Lieman{\Liegroup{}}}
{\Lieman{\Liegroup{}}}\)}{\avecf{}}\]}{\point}&=
\func{\(\der{\gltrans{\LieG{\Liegroup{}}}{\g{}}}{\Lieman{\Liegroup{}}}
{\Lieman{\Liegroup{}}}\)}{\func{\avecf{}}{\point}}\cr
&=\func{\(\der{\gltrans{\LieG{\Liegroup{}}}{\g{}}}{\Lieman{\Liegroup{}}}
{\Lieman{\Liegroup{}}}\)}{\func{\(\der{\gltrans{\LieG{\Liegroup{}}}{\point}}{\Lieman{\Liegroup{}}}
{\Lieman{\Liegroup{}}}\)}{\func{\avecf{}}{\IG{}}}}\cr
&=\func{\[\cmp{\(\der{\gltrans{\LieG{\Liegroup{}}}{\g{}}}{\Lieman{\Liegroup{}}}
{\Lieman{\Liegroup{}}}\)}{\(\der{\gltrans{\LieG{\Liegroup{}}}{\point}}{\Lieman{\Liegroup{}}}
{\Lieman{\Liegroup{}}}\)}\]}{\func{\avecf{}}{\IG{}}}\cr
&=\func{\[\der{\(\cmp{\gltrans{\LieG{\Liegroup{}}}{\g{}}}{\gltrans{\LieG{\Liegroup{}}}{\point}}\)}{\Lieman{\Liegroup{}}}
{\Lieman{\Liegroup{}}}\]}{\func{\avecf{}}{\IG{}}}\cr
&=\func{\[\der{\gltrans{\LieG{\Liegroup{}}}{\g{}\gop{}\point}}{\Lieman{\Liegroup{}}}
{\Lieman{\Liegroup{}}}\]}{\func{\avecf{}}{\IG{}}}\cr
&=\func{\avecf{}}{\g{}\gop{}\point}\cr
&=\func{\[\cmp{\avecf{}}{\gltrans{\LieG{\Liegroup{}}}{\g{}}}\]}{\point}.
\end{align}
Thus,
\begin{equation}
\Foreach{\g{}}{\G{}}\cmp{\(\der{\gltrans{\LieG{\Liegroup{}}}{\g{}}}{\Lieman{\Liegroup{}}}
{\Lieman{\Liegroup{}}}\)}{\avecf{}}=\cmp{\avecf{}}{\gltrans{\LieG{\Liegroup{}}}{\g{}}}.
\end{equation}
In addition, for every $\g{}$ in $\G{}$, since $\gltrans{\LieG{\Liegroup{}}}{\g{}}$ is an
$\infty$-diffeomorphism of $\Lieman{\Liegroup{}}$ and sends
$\IG{}$ to $\g{}$, clealy its differential $\der{\gltrans{\LieG{\Liegroup{}}}{\g{}}}{\Lieman{\Liegroup{}}}
{\Lieman{\Liegroup{}}}$ maps $\tanspace{\IG{}}{\Lieman{\Liegroup{}}}$ injectively onto
$\tanspace{\g{}}{\Lieman{\Liegroup{}}}$. That is,
\begin{equation}
\Foreach{\g{}}{\G{}}
\func{\image{\[\der{\gltrans{\LieG{\Liegroup{}}}{\g{}}}{\Lieman{\Liegroup{}}}
{\Lieman{\Liegroup{}}}\]}}{\tanspace{\IG{}}{\Lieman{\Liegroup{}}}}=
\tanspace{\g{}}{\Lieman{\Liegroup{}}},
\end{equation}
and hence according to \Ref{lemleftinvariantvectorfieldsequiv0p1eq1} and
\Ref{lemleftinvariantvectorfieldsequiv0p1eq2},
\begin{equation}
\Foreach{\g{}}{\G{}}
\func{\avecf{}}{\g{}}\in\tanspace{\g{}}{\Lieman{\Liegroup{}}},
\end{equation}
which means,
\begin{equation}
\cmp{\basep{\Lieman{\Liegroup{}}}}{\avecf{}}=\identity{\G{}}.
\end{equation}
\endp
\end{itemize}
\endlem
\theorem\label{thmleftinvariantvectorfieldsaresmooth}
Every left-invariant vector-field on $\Liegroup{}$ is a smooth map from $\Lieman{\Liegroup{}}$ to $\Tanbun{\Lieman{\Liegroup{}}}$, that is,
\begin{equation}
\Foreach{\avecf{}}{\Leftinvvf{\Liegroup{}}}
\avecf{}\in\mapdifclass{\infty}{\Lieman{\Liegroup{}}}{\Tanbun{\Lieman{\Liegroup{}}}}.
\end{equation}
This means every left-invariant vector-field on $\Liegroup{}$ is a smooth vector-field on the manifold $\Lieman{\Liegroup{}}$, that is,
\begin{equation}
\Leftinvvf{\Liegroup{}}\subseteq\vecf{\Lieman{\Liegroup{}}}{\infty}.
\end{equation}
\proof
$\avecf{}$ is taken as an arbitrary element of $\Leftinvvf{\Liegroup{}}$. So according to
\refdef{defleftinvariantvectorfields} and \reflem{lemleftinvariantvectorfieldsequiv0},
$\avecf{}$ is an element of $\Func{\G{}}{\tanbun{\Lieman{\Liegroup{}}}}$ such that,
\begin{align}
&\cmp{\basep{\Lieman{\Liegroup{}}}}{\avecf{}}=\identity{\G{}},
\label{thmleftinvariantvectorfieldsaresmoothpeq1}\\
&\Foreach{\g{}}{\G{}}\func{\(\der{\gltrans{\LieG{\Liegroup{}}}{\g{}}}{\Lieman{\Liegroup{}}}
{\Lieman{\Liegroup{}}}\)}{\func{\avecf{}}{\IG{}}}=\func{\avecf{}}{\g{}}.
\label{thmleftinvariantvectorfieldsaresmoothpeq2}
\end{align}
\begin{itemize}
\item[$\pr{1}$]
$\g{}$ is taken as an arbitrary element of $\G{}$. $\phi$ is taken as a chart of $\Lieman{\Liegroup{}}$ centered at $\g{}$,
that is an element of $\defset{\p{\phi}}{\maxatlas{}}{\g{}\in\domain{\p{\phi}},~\func{\p{\phi}}{\g{}}=\zerovec{}}$.
According to \Ref{eqtangentbundlemaps},
\begin{align}\label{thmleftinvariantvectorfieldsaresmoothp1eq1}
&\tanchart{\Lieman{\Liegroup{}}}{\phi}\in
\Func{\func{\pimage{\basep{\Lieman{\Liegroup{}}}}}{\domain{\phi}}}{\Cprod{\funcimage{\phi}}{\R^n}},\cr
&\Foreach{\avec{}}{\func{\pimage{\basep{\Lieman{\Liegroup{}}}}}{\domain{\phi}}}
\func{\tanchart{\Lieman{\Liegroup{}}}{\phi}}{\avec{}}\eqdef
\opair{\func{\(\cmp{\phi}{\basep{\Lieman{\Liegroup{}}}}\)}{\avec{}}}
{\func{\(\tanspaceiso{\func{\basep{\Lieman{\Liegroup{}}}}{\avec{}}}{\Lieman{\Liegroup{}}}{\phi}\)}{\avec{}}},
\end{align}
where $\tanchart{\Lieman{\Liegroup{}}}{\phi}$ is the tangent-bundle chart of $\Lieman{\Liegroup{}}$ associated with $\phi$.
$\tanchart{\Lieman{\Liegroup{}}}{\phi}$ is an element of $\tanatlas{\Lieman{\Liegroup{}}}$
(the maximal-atlas of the tangent-bundle of $\Lieman{\Liegroup{}}$),
and hence a chart of the tangent-bundle of $\Lieman{\Liegroup{}}$.
\Ref{thmleftinvariantvectorfieldsaresmoothpeq1} and \Ref{thmleftinvariantvectorfieldsaresmoothp1eq1} imply that,
\begin{equation}\label{thmleftinvariantvectorfieldsaresmoothp1eq2}
\func{\image{\avecf{}}}{\domain{\phi}}\subseteq
\func{\pimage{\basep{\Lieman{\Liegroup{}}}}}{\domain{\phi}}=
\domain{\tanchart{\Lieman{\Liegroup{}}}{\phi}}.
\end{equation}
In addition, it is clear that $\gltrans{\LieG{\Liegroup{}}}{\g{}}$ transfers the chart $\phi$ of
$\Lieman{\Liegroup{}}$ centered at $\g{}$ to the chart $\cmp{\phi}{\gltrans{\LieG{\Liegroup{}}}{\g{}}}$
of $\Lieman{\Liegroup{}}$ centered at $\IG{}$. That is,
\begin{equation}\label{thmleftinvariantvectorfieldsaresmoothp1eq3}
\bar{\phi}:=\cmp{\phi}{\gltrans{\LieG{\Liegroup{}}}{\g{}}}\in
\defset{\p{\phi}}{\maxatlas{}}{\IG{}\in\domain{\p{\phi}},~\func{\p{\phi}}{\IG{}}=\zerovec{}}.
\end{equation}
So the domain of $\bar{\phi}$
is an open set of $\lietops{\Liegroup{}}$ containing $\IG{}$, and hence
according to \refdef{defnucleioftopologicalgroup}, is a nucleus of $\Lietopg{\Liegroup{}}$. That is,
\begin{equation}\label{thmleftinvariantvectorfieldsaresmoothp1eq4}
\domain{\bar{\phi}}\in\nuclei{\Lietopg{\Liegroup{}}}.
\end{equation}
Thus, according to \reflem{lemeachnucleusincludesarestrictednucleus},
\begin{equation}\label{thmleftinvariantvectorfieldsaresmoothp1eq5}
\Existsis{\V}{\nuclei{\Lietopg{\Liegroup{}}}}
\bigg(\V\subseteq\domain{\bar{\phi}},~
\func{\image{\gop{}}}{\Cprod{\V}{\V}}\subseteq\domain{\bar{\phi}}\bigg).
\end{equation}
It is clear that,
\begin{equation}\label{thmleftinvariantvectorfieldsaresmoothp1eq6}
\domain{\bar{\phi}}=
\func{\pimage{\[\gltrans{\LieG{\Liegroup{}}}{\g{}}\]}}{\domain{\phi}},
\end{equation}
and therefore considering that $\V\subseteq\domain{\bar{\phi}}$, evidently,
\begin{equation}\label{thmleftinvariantvectorfieldsaresmoothp1eq7}
\WW{}:=\func{\image{\[\gltrans{\LieG{\Liegroup{}}}{\g{}}\]}}{\V}\subseteq
\domain{\phi}.
\end{equation}
Moreover, since $\IG{}\in\V$, clearly,
\begin{equation}\label{thmleftinvariantvectorfieldsaresmoothp1eq8}
\g{}\in\WW{},
\end{equation}
and additionally, since $\V$ is an open set of $\lietops{\Liegroup{}}$ and $\gltrans{\LieG{\Liegroup{}}}{\g{}}$
is a homeomorphism from $\lietops{\Liegroup{}}$ to itself, $\WW{}$ is also an open set of $\lietops{\Liegroup{}}$.
\begin{equation}\label{thmleftinvariantvectorfieldsaresmoothp1eq9}
\WW{}\in\lietop{\Liegroup{}}.
\end{equation}
Therefore the restriction of $\bar{\phi}$ to $\V$ and the restriction of $\phi$ to $\WW{}$
yield charts of $\Lieman{\Liegroup{}}$ centered at $\IG{}$ and $\g{}$, respectively. That is,
\begin{align}
\psi:=&\func{\res{\phi}}{\WW{}}\in\defset{\p{\phi}}{\maxatlas{}}{\g{}\in\domain{\p{\phi}},~\func{\p{\phi}}{\g{}}=\zerovec{}},
\label{thmleftinvariantvectorfieldsaresmoothp1eq10}\\
\bar{\psi}:=&\func{\res{\bar{\phi}}}{\V}\in\defset{\p{\phi}}{\maxatlas{}}{\IG{}\in\domain{\p{\phi}},~\func{\p{\phi}}{\IG{}}=\zerovec{}}.
\label{thmleftinvariantvectorfieldsaresmoothp1eq11}
\end{align}
So since $\WW{}=\domain{\psi}$, \Ref{thmleftinvariantvectorfieldsaresmoothp1eq2}
and \Ref{thmleftinvariantvectorfieldsaresmoothp1eq7} imply,
\begin{equation}\label{thmleftinvariantvectorfieldsaresmoothp1eq12}
\func{\image{\avecf{}}}{\domain{\psi}}\subseteq
\domain{\tanchart{\Lieman{\Liegroup{}}}{\phi}}.
\end{equation}
Therefore according to \Ref{thmleftinvariantvectorfieldsaresmoothpeq1},
\Ref{thmleftinvariantvectorfieldsaresmoothp1eq1} and
\Ref{thmleftinvariantvectorfieldsaresmoothp1eq10},
\begin{align}\label{thmleftinvariantvectorfieldsaresmoothp1eq13}
\cmp{\tanchart{\Lieman{\Liegroup{}}}{\phi}}{\cmp{\avecf{}}{\finv{\psi}}}\in
\Func{\func{\image{\phi}}{\WW{}}}{\Cprod{\funcimage{\phi}}{\R^n}},
\end{align}
and,
\begin{align}\label{thmleftinvariantvectorfieldsaresmoothp1eq14}
&\begin{aligned}
\Foreach{\x}{\func{\image{\phi}}{\WW{}}}\cr
\end{aligned}\cr
&\begin{aligned}
\func{\(\cmp{\tanchart{\Lieman{\Liegroup{}}}{\phi}}{\cmp{\avecf{}}{\finv{\psi}}}\)}{\x}&=
\opair{\func{\[\cmp{\phi}{\basep{\Lieman{\Liegroup{}}}}\]}{\func{\avecf{}}{\func{\finv{\phi}}{\x}}}}
{\func{\[\tanspaceiso{\func{\basep{\Lieman{\Liegroup{}}}}{\func{\avecf{}}{\func{\finv{\phi}}{\x}}}}
{\Lieman{\Liegroup{}}}{\phi}\]}{\func{\avecf{}}{\func{\finv{\phi}}{\x}}}}\cr
&=\opair{\x}{\func{\[\tanspaceiso{\func{\finv{\phi}}{\x}}
{\Lieman{\Liegroup{}}}{\phi}\]}{\func{\avecf{}}{\func{\finv{\phi}}{\x}}}}\cr
&=\opair{\func{\Injection{\func{\image{\phi}}{\WW{}}}{\funcimage{\phi}}}{\x}}{\func{\alpha}{\x}},
\end{aligned}
\end{align}
where $\alpha$ is defined to be the mapping,
\begin{align}\label{thmleftinvariantvectorfieldsaresmoothp1eq15}
&\alpha\in\Func{\func{\image{\phi}}{\WW{}}}{\R^n},\cr
&\Foreach{\x}{\func{\image{\phi}}{\WW{}}}\func{\alpha}{\x}\eqdef
\func{\[\tanspaceiso{\func{\finv{\phi}}{\x}}
{\Lieman{\Liegroup{}}}{\phi}\]}{\func{\avecf{}}{\func{\finv{\phi}}{\x}}}.
\end{align}
Thus according to \Ref{thmleftinvariantvectorfieldsaresmoothpeq2},
\begin{equation}\label{thmleftinvariantvectorfieldsaresmoothp1eq16}
\Foreach{\x}{\func{\image{\phi}}{\WW{}}}
\func{\alpha}{\x}=\func{\[\tanspaceiso{\func{\finv{\phi}}{\x}}
{\Lieman{\Liegroup{}}}{\phi}\]}{\func{\[\der{\gltrans{\LieG{\Liegroup{}}}{\func{\finv{\phi}}{\x}}}{\Lieman{\Liegroup{}}}
{\Lieman{\Liegroup{}}}\]}{\func{\avecf{}}{\IG{}}}}.
\end{equation}
Considering that $\func{\basep{\Lieman{\Liegroup{}}}}{\func{\avecf{}}{\IG{}}}=\IG{}$,
$\bar{\psi}$ is a chart of $\Lieman{\Liegroup{}}$ around $\IG{}$,
and $\phi$ is a chart of $\Lieman{\Liegroup{}}$ around
$\func{\finv{\phi}}{\x}=\func{\gltrans{\LieG{\Liegroup{}}}{\func{\finv{\phi}}{\x}}}{\IG{}}$ for every $\x$ in $\func{\image{\phi}}{\WW{}}$,
and according to \Ref{eqdefinitionofdifferentialofamap} and \Ref{thmleftinvariantvectorfieldsaresmoothp1eq11},
\begin{align}\label{thmleftinvariantvectorfieldsaresmoothp1eq17}
\Foreach{\x}{\func{\image{\phi}}{\WW{}}}
&\hskip 0.5\baselineskip
\func{\[\der{\gltrans{\LieG{\Liegroup{}}}{\func{\finv{\phi}}{\x}}}{\Lieman{\Liegroup{}}}
{\Lieman{\Liegroup{}}}\]}{\func{\avecf{}}{\IG{}}}\cr
&=\func{\bigg[\cmp{\finv{\(\tanspaceiso{\func{\finv{\phi}}{\x}}{\Lieman{\Liegroup{}}}{\phi}\)}}
{\cmp{\bigg(\func{\[\banachder{\(\cmp{\phi}
{\cmp{\gltrans{\LieG{\Liegroup{}}}{\func{\finv{\phi}}{\x}}}{\finv{\bar{\psi}}}}\)}{\R^n}{\R^n}\]}
{\func{\bar{\phi}}{\IG{}}}\bigg)}{\tanspaceiso{\IG{}}{\Lieman{\Liegroup{}}{}}{\bar{\psi}}}}\bigg]}{\avec{0}},\cr
&{}
\end{align}
where,
\begin{equation}\label{thmleftinvariantvectorfieldsaresmoothp1eq18}
\avec{0}:=\func{\avecf{}}{\IG{}}.
\end{equation}
Combining \Ref{thmleftinvariantvectorfieldsaresmoothp1eq16} and \Ref{thmleftinvariantvectorfieldsaresmoothp1eq17}
yields,
\begin{align}\label{thmleftinvariantvectorfieldsaresmoothp1eq19}
\Foreach{\x}{\func{\image{\phi}}{\WW{}}}
\func{\alpha}{\x}
=\func{\bigg[\cmp{\bigg(\func{\[\banachder{\(\cmp{\phi}
{\cmp{\gltrans{\LieG{\Liegroup{}}}{\func{\finv{\phi}}{\x}}}{\finv{\bar{\psi}}}}\)}{\R^n}{\R^n}\]}
{\func{\bar{\phi}}{\IG{}}}\bigg)}{\tanspaceiso{\IG{}}{\Lieman{\Liegroup{}}{}}{\bar{\psi}}}\bigg]}{\avec{0}}.
\end{align}
According to \Ref{thmleftinvariantvectorfieldsaresmoothp1eq3}, \Ref{thmleftinvariantvectorfieldsaresmoothp1eq7},
\Ref{thmleftinvariantvectorfieldsaresmoothp1eq10}, and \Ref{thmleftinvariantvectorfieldsaresmoothp1eq11},
it is evident that,
\begin{equation}\label{thmleftinvariantvectorfieldsaresmoothp1eq20}
\bar{\psi}=\cmp{\psi}{\gltrans{\LieG{\Liegroup{}}}{\g{}}},
\end{equation}
and hence considering that $\gltrans{\LieG{\Liegroup{}}}{\invg{\g{}}{}}$ is the inverse of
$\gltrans{\LieG{\Liegroup{}}}{\g{}}$,
\begin{align}\label{thmleftinvariantvectorfieldsaresmoothp1eq21}
\finv{\bar{\psi}}&=\cmp{\finv{\[\gltrans{\LieG{\Liegroup{}}}{\g{}}\]}}{\finv{\psi}}\cr
&=\cmp{\gltrans{\LieG{\Liegroup{}}}{\invg{\g{}}{}}}{\finv{\psi}}.
\end{align}
The mapping $\eta$ is defined as,
\begin{equation}\label{thmleftinvariantvectorfieldsaresmoothp1eq22}
\eta=\cmp{\cmp{\phi}{\[\cmp{\gltrans{\LieG{\Liegroup{}}}{\g{}}}
{\cmp{\gop{}}{\(\funcprod{\gltrans{\LieG{\Liegroup{}}}{\invg{\g{}}{}}}{\gltrans{\LieG{\Liegroup{}}}{\invg{\g{}}{}}}\)}}\]}}{\finv{\(\funcprod{\psi}{\psi}\)}}.
\end{equation}
Clearly,
\begin{equation}\label{thmleftinvariantvectorfieldsaresmoothp1eq22a}
\domain{\funcprod{\psi}{\psi}}=\Cprod{\WW{}}{\WW{}},
\end{equation}
and,
\begin{align}\label{thmleftinvariantvectorfieldsaresmoothp1eq22b}
\Foreach{\opair{\point_1}{\point_2}}{\Cprod{\WW{}}{\WW{}}}
\func{\[\cmp{\gltrans{\LieG{\Liegroup{}}}{\g{}}}
{\cmp{\gop{}}{\(\funcprod{\gltrans{\LieG{\Liegroup{}}}{\invg{\g{}}{}}}{\gltrans{\LieG{\Liegroup{}}}{\invg{\g{}}{}}}\)}}\]}
{\binary{\point_1}{\point_2}}=\g{}\gop{}\(\invg{\g{}}{}\gop{}\point_1\)\gop{}\(\invg{\g{}}{}\gop{}\point_2\).
\end{align}
In addition, according to \Ref{thmleftinvariantvectorfieldsaresmoothp1eq7},
\begin{align}\label{thmleftinvariantvectorfieldsaresmoothp1eq22c}
\Foreach{\point}{\WW{}}
\(\invg{\g{}}{}\gop{}\point\)\in\V,
\end{align}
and hence according to \Ref{thmleftinvariantvectorfieldsaresmoothp1eq5},
\begin{equation}\label{thmleftinvariantvectorfieldsaresmoothp1eq22d}
\Foreach{\opair{\point_1}{\point_2}}{\Cprod{\WW{}}{\WW{}}}
\(\invg{\g{}}{}\gop{}\point_1\)\gop{}\(\invg{\g{}}{}\gop{}\point_2\)\in\domain{\bar{\phi}},
\end{equation}
and thus according to \Ref{thmleftinvariantvectorfieldsaresmoothp1eq3},
\begin{equation}\label{thmleftinvariantvectorfieldsaresmoothp1eq22e}
\Foreach{\opair{\point_1}{\point_2}}{\Cprod{\WW{}}{\WW{}}}
\g{}\gop{}\(\invg{\g{}}{}\gop{}\point_1\)\gop{}\(\invg{\g{}}{}\gop{}\point_2\)\in\domain{\phi}
\end{equation}
Therefore, according to \Ref{thmleftinvariantvectorfieldsaresmoothp1eq22a} and \Ref{thmleftinvariantvectorfieldsaresmoothp1eq22b},
and using \Ref{thmleftinvariantvectorfieldsaresmoothp1eq10},
\begin{align}\label{thmleftinvariantvectorfieldsaresmoothp1eq23}
&\eta\in\Func{\Cprod{\func{\image{\phi}}{\WW{}}}{\func{\image{\phi}}{\WW{}}}}{\funcimage{\phi}},\cr
&\Foreach{\opair{\x}{\y}}{\Cprod{\func{\image{\phi}}{\WW{}}}{\func{\image{\phi}}{\WW{}}}}
\func{\eta}{\binary{\x}{\y}}=\func{\phi}{\func{\finv{\phi}}{\x}\gop{}\invg{\g{}}{}\gop{}\func{\finv{\phi}}{\y}},
\end{align}
and thus according to \Ref{thmleftinvariantvectorfieldsaresmoothp1eq10} and \Ref{thmleftinvariantvectorfieldsaresmoothp1eq21},
\begin{align}\label{thmleftinvariantvectorfieldsaresmoothp1eq24}
\Foreach{\opair{\x}{\y}}{\Cprod{\func{\image{\phi}}{\WW{}}}{\func{\image{\phi}}{\WW{}}}}
\func{\eta}{\binary{\x}{\y}}&=\func{\(\cmp{\phi}
{\cmp{\gltrans{\LieG{\Liegroup{}}}{\func{\finv{\phi}}{\x}}}{\cmp{\gltrans{\LieG{\Liegroup{}}}{\invg{\g{}}{}}}{\finv{\psi}}}}\)}{\y}\cr
&=\func{\(\cmp{\phi}
{\cmp{\gltrans{\LieG{\Liegroup{}}}{\func{\finv{\phi}}{\x}}}{\finv{\bar{\psi}}}}\)}{\y}.
\end{align}
By definition,
\begin{align}\label{thmleftinvariantvectorfieldsaresmoothp1eq25}
&\lambda\indef\Func{\R^n}{\Func{\R^n}{\Cprod{\R^n}{\R^n}}},\cr
&\Foreach{\x}{\R^n}\Foreach{\y}{\R^n}
\func{\[\func{\lambda}{\x}\]}{\y}\eqdef\opair{\x}{\y}.
\end{align}
Using \Ref{thmleftinvariantvectorfieldsaresmoothp1eq25} in \Ref{thmleftinvariantvectorfieldsaresmoothp1eq24} yields,
\begin{equation}\label{thmleftinvariantvectorfieldsaresmoothp1eq26}
\Foreach{\x}{\func{\image{\phi}}{\WW{}}}
\Foreach{\y}{\func{\image{\phi}}{\WW{}}}
\func{\(\cmp{\phi}
{\cmp{\gltrans{\LieG{\Liegroup{}}}{\func{\finv{\phi}}{\x}}}{\finv{\bar{\psi}}}}\)}{\y}=
\func{\[\cmp{\eta}{\func{\lambda}{\x}}\]}{\y}.
\end{equation}
According to \Ref{thmleftinvariantvectorfieldsaresmoothp1eq23} and \Ref{thmleftinvariantvectorfieldsaresmoothp1eq25},
it is clear that,
\begin{equation}\label{thmleftinvariantvectorfieldsaresmoothp1eq27}
\Foreach{\x}{\func{\image{\phi}}{\WW{}}}
\cmp{\eta}{\func{\lambda}{\x}}\in\Func{\func{\image{\phi}}{\WW{}}}{\funcimage{\phi}}.
\end{equation}
Also, as a trivial consequence of \Ref{thmleftinvariantvectorfieldsaresmoothp1eq23}
\begin{equation}\label{thmleftinvariantvectorfieldsaresmoothp1eq28}
\Foreach{\x}{\func{\image{\phi}}{\WW{}}}
\(\cmp{\phi}
{\cmp{\gltrans{\LieG{\Liegroup{}}}{\func{\finv{\phi}}{\x}}}{\finv{\bar{\psi}}}}\)
\in\Func{\func{\image{\phi}}{\WW{}}}{\funcimage{\phi}}.
\end{equation}
\Ref{thmleftinvariantvectorfieldsaresmoothp1eq26}, \Ref{thmleftinvariantvectorfieldsaresmoothp1eq27},
and \Ref{thmleftinvariantvectorfieldsaresmoothp1eq28} imply that,
\begin{equation}\label{thmleftinvariantvectorfieldsaresmoothp1eq29}
\Foreach{\x}{\func{\image{\phi}}{\WW{}}}
\(\cmp{\phi}
{\cmp{\gltrans{\LieG{\Liegroup{}}}{\func{\finv{\phi}}{\x}}}{\finv{\bar{\psi}}}}\)=
\cmp{\eta}{\func{\lambda}{\x}}.
\end{equation}
Considering that $\gop{}$, $\gltrans{\LieG{\Liegroup{}}}{\g{}{}}$, and
$\gltrans{\LieG{\Liegroup{}}}{\invg{\g{}}{}}\times\gltrans{\LieG{\Liegroup{}}}{\invg{\g{}}{}}$ are smooth maps, that is,
\begin{align}\label{thmleftinvariantvectorfieldsaresmoothp1eq30}
\funcprod{\gltrans{\LieG{\Liegroup{}}}{\invg{\g{}}{}}}{\gltrans{\LieG{\Liegroup{}}}{\invg{\g{}}{}}}&\in
\mapdifclass{\infty}{\manprod{\Lieman{\Liegroup{}}}{\Lieman{\Liegroup{}}}}
{\manprod{\Lieman{\Liegroup{}}}{\Lieman{\Liegroup{}}}},\cr
\gop{}&\in\mapdifclass{\infty}{\manprod{\Lieman{\Liegroup{}}}{\Lieman{\Liegroup{}}}}{\Lieman{\Liegroup{}}},\cr
\gltrans{\LieG{\Liegroup{}}}{\g{}{}}&\in
\mapdifclass{\infty}{\Lieman{\Liegroup{}}}{\Lieman{\Liegroup{}}},
\end{align}
based on the composition rule of smooth maps it is evident that,
\begin{equation}\label{thmleftinvariantvectorfieldsaresmoothp1eq31}
\[\cmp{\gltrans{\LieG{\Liegroup{}}}{\g{}}}
{\cmp{\gop{}}{\(\funcprod{\gltrans{\LieG{\Liegroup{}}}{\invg{\g{}}{}}}{\gltrans{\LieG{\Liegroup{}}}{\invg{\g{}}{}}}\)}}\]\in
\mapdifclass{\infty}{\manprod{\Lieman{\Liegroup{}}}{\Lieman{\Liegroup{}}}}{\Lieman{\Liegroup{}}}.
\end{equation}
$\phi$ is a chart of $\Lieman{\Liegroup{}}$ and
since $\psi$ is a chart of $\Lieman{\Liegroup{}}$, $\psi\times\psi$ is a chart of
$\manprod{\Lieman{\Liegroup{}}}{\Lieman{\Liegroup{}}}$. So \Ref{thmleftinvariantvectorfieldsaresmoothp1eq31}
together with \Ref{thmleftinvariantvectorfieldsaresmoothp1eq22} and
\Ref{thmleftinvariantvectorfieldsaresmoothp1eq23} implies that $\eta$ is a smooth map from
$\Cprod{\func{\image{\phi}}{\WW{}}}{\func{\image{\phi}}{\WW{}}}$
to $\funcimage{\phi}$, that is,
\begin{equation}\label{thmleftinvariantvectorfieldsaresmoothp1eq32}
\eta\in\banachmapdifclass{\infty}{\Cprod{\R^{n}}{\R^{n}}}{\R^n}{\Cprod{\func{\image{\phi}}{\WW{}}}{\func{\image{\phi}}{\WW{}}}}{\funcimage{\phi}},
\end{equation}
and hence the derived map of $\eta$, that is 
$\banachder{\eta}{\Cprod{\R^n}{\R^n}}{\R^n}$, is a smooth map from $\Cprod{\func{\image{\phi}}{\WW{}}}{\func{\image{\phi}}{\WW{}}}$ to
$\Lin{\Cprod{\R^n}{\R^n}}{\R^n}$. That is,
\begin{equation}\label{thmleftinvariantvectorfieldsaresmoothp1eq33}
\banachder{\eta}{\Cprod{\R^n}{\R^n}}{\R^n}\in
\banachmapdifclass{\infty}{\Cprod{\R^n}{\R^n}}
{\NVLin{\Cprod{\R^n}{\R^n}}{\R^n}}{\Cprod{\func{\image{\phi}}{\WW{}}}
{\func{\image{\phi}}{\WW{}}}}{\Lin{\Cprod{\R^n}{\R^n}}{\R^n}}.
\end{equation}
Equivalent to \Ref{thmleftinvariantvectorfieldsaresmoothp1eq25},
\begin{equation}\label{thmleftinvariantvectorfieldsaresmoothp1eq34}
\Foreach{\x}{\R^n}
\func{\lambda}{\x}=\cmp{\[\funcprod{\constmap{\R^n}{\x}}{\identity{\R^n}}\]}{\diagmap{\R^n}}.
\end{equation}
So since both $\identity{\R^n}$ and $\constmap{\R^n}{\x}$
are smooth maps from $\R^n$ to $\R^n$,
and $\diagmap{\R^n}$ is a smooth map from $\R^n$ to
$\Cprod{\R^n}{\R^n}$,
it is evident that,
\begin{equation}\label{thmleftinvariantvectorfieldsaresmoothp1eq35}
\Foreach{\x}{\func{\image{\phi}}{\WW{}}}
\func{\lambda}{\x}\in\banachmapdifclass{\infty}{\R^n}{\Cprod{\R^{n}}{\R^n}}
{\R^n}{\Cprod{\R^n}{\R^n}}.
\end{equation}
Thus, based on the composition rule of smooth maps, \Ref{thmleftinvariantvectorfieldsaresmoothp1eq32}
and \Ref{thmleftinvariantvectorfieldsaresmoothp1eq35} imply that,
\begin{equation}\label{thmleftinvariantvectorfieldsaresmoothp1eq36}
\Foreach{\x}{\func{\image{\phi}}{\WW{}}}
\cmp{\eta}{\func{\lambda}{\x}}\in\banachmapdifclass{\infty}{\R^n}{\R^n}{\func{\image{\phi}}{\WW{}}}{\funcimage{\phi}}.
\end{equation}
Therefore,
\begin{equation}\label{thmleftinvariantvectorfieldsaresmoothp1eq37}
\Foreach{\x}{\func{\image{\phi}}{\WW{}}}
\banachder{\(\cmp{\eta}{\func{\lambda}{\x}}\)}{\R^n}{\R^n}\in
\banachmapdifclass{\infty}{\R^n}{\NVLin{\R^n}{\R^n}}{\func{\image{\phi}}{\WW{}}}{\Lin{\R^n}{\R^n}}.
\end{equation}
In addition, according to \Ref{thmleftinvariantvectorfieldsaresmoothp1eq19} and \Ref{thmleftinvariantvectorfieldsaresmoothp1eq29},
\begin{equation}\label{thmleftinvariantvectorfieldsaresmoothp1eq38}
\Foreach{\x}{\func{\image{\phi}}{\WW{}}}
\func{\alpha}{\x}=\func{\bigg(\func{\[\banachder{\(\cmp{\eta}{\func{\lambda}{\x}}\)}{\R^n}{\R^n}\]}{\func{\bar{\phi}}{\IG{}}}\bigg)}
{\func{\tanspaceiso{\IG{}}{\Lieman{\Liegroup{}}{}}{\bar{\psi}}}{\avec{0}}}.
\end{equation}
According to the chain-rule of differentiation of differentiable maps between Banach-spaces,
\begin{align}\label{thmleftinvariantvectorfieldsaresmoothp1eq39}
\Foreach{\x}{\func{\image{\phi}}{\WW{}}}\
&\hskip 0.5\baselineskip\func{\[\banachder{\(\cmp{\eta}{\func{\lambda}{\x}}\)}{\R^n}{\R^n}\]}{\func{\bar{\phi}}{\IG{}}}\cr
&=\cmp{\bigg(\func{\[\cmp{\(\banachder{\eta}{\Cprod{\R^n}{\R^n}}{\R^n}\)}{\func{\lambda}{\x}}\]}{\func{\bar{\phi}}{\IG{}}}\bigg)}
{\[\func{\(\banachder{\func{\lambda}{\x}}{\R^n}{\Cprod{\R^n}{\R^n}}\)}{\func{\bar{\phi}}{\IG{}}}\]}.\cr
&{}
\end{align}
According to \Ref{thmleftinvariantvectorfieldsaresmoothp1eq34}, it can be easily verified that,
\begin{equation}\label{thmleftinvariantvectorfieldsaresmoothp1eq40}
\Foreach{\x}{\func{\image{\phi}}{\WW{}}}
\func{\(\banachder{\func{\lambda}{\x}}{\R^n}{\Cprod{\R^n}{\R^n}}\)}{\func{\bar{\phi}}{\IG{}}}=\u,
\end{equation}
where,
\begin{equation}\label{thmleftinvariantvectorfieldsaresmoothp1eq41}
\u:=\cmp{\[\funcprod{\constmap{\R^n}{\zerovec{}}}
{\identity{\R^n}}\]}{\diagmap{\R^n}},
\end{equation}
or equivalently,
\begin{align}\label{thmleftinvariantvectorfieldsaresmoothp1eq42}
&\u\indef\Func{\R^n}{\Cprod{\R^n}{\R^n}}\cr
&\Foreach{z}{\R^n}\func{\u}{z}\eqdef\opair{\zerovec{}}{z},
\end{align}
where $\zerovec{}$ denotes the neutral element of addition in $\R^n$.
\Ref{thmleftinvariantvectorfieldsaresmoothp1eq38}-\Ref{thmleftinvariantvectorfieldsaresmoothp1eq42}
imply that,
\begin{equation}\label{thmleftinvariantvectorfieldsaresmoothp1eq43}
\Foreach{\x}{\func{\image{\phi}}{\WW{}}}
\func{\alpha}{\x}=\func{\[\func{\(\banachder{\eta}{\Cprod{\R^n}{\R^n}}{\R^n}\)}{\func{\[\func{\lambda}{\x}\]}{\func{\bar{\phi}}{\IG{}}}}\]}
{\func{\u}{\func{\tanspaceiso{\IG{}}{\Lieman{\Liegroup{}}{}}{\bar{\psi}}}{\avec{0}}}}.
\end{equation}
By defining the mappings,
\begin{align}\label{thmleftinvariantvectorfieldsaresmoothp1eq44}
&\mu\indef\Func{\R^n}{\Cprod{\R^n}{\R^n}},\cr
&\Foreach{\x}{\R^n}\func{\mu}{\x}\eqdef\func{\[\func{\lambda}{\x}\]}{\func{\bar{\phi}}{\IG{}}}=
\opair{\x}{\zerovec{}},
\end{align}
and,
\begin{align}\label{thmleftinvariantvectorfieldsaresmoothp1eq45}
&\nu\indef\Func{\Lin{\Cprod{\R^n}{\R^n}}{\R^n}}{\R^n},\cr
&\Foreach{T}{\Lin{\Cprod{\R^n}{\R^n}}{\R^n}}\func{\nu}{T}\eqdef
\func{T}{\func{\u}{\func{\tanspaceiso{\IG{}}{\Lieman{\Liegroup{}}{}}{\bar{\psi}}}{\avec{0}}}},
\end{align}
\Ref{thmleftinvariantvectorfieldsaresmoothp1eq15}, \Ref{thmleftinvariantvectorfieldsaresmoothp1eq33}, and
\Ref{thmleftinvariantvectorfieldsaresmoothp1eq43} imply that,
\begin{equation}\label{thmleftinvariantvectorfieldsaresmoothp1eq46}
\alpha=
\cmp{\cmp{\nu}{\(\banachder{\eta}{\Cprod{\R^n}{\R^n}}{\R^n}\)}}{\mu}.
\end{equation}
It is trivial that the mappings $\mu$ and $\nu$ are both linear maps, and since their domain is
finite-dimensional Banach-spaces, they are linear-continuous maps and hence infinitely-differentiable. That is,
\begin{align}
\mu&\in\banachmapdifclass{\infty}{\R^n}{\Cprod{\R^n}{\R^n}}{\R^n}{\Cprod{\R^n}{\R^n}},
\label{thmleftinvariantvectorfieldsaresmoothp1eq47}\\
\nu&\in\banachmapdifclass{\infty}{\NVLin{\Cprod{\R^n}{\R^n}}{\R^n}}{\R^n}{\Lin{\Cprod{\R^n}{\R^n}}{\R^n}}{\R^n}.
\label{thmleftinvariantvectorfieldsaresmoothp1eq48}
\end{align}
\Ref{thmleftinvariantvectorfieldsaresmoothp1eq33}, \Ref{thmleftinvariantvectorfieldsaresmoothp1eq46},
\Ref{thmleftinvariantvectorfieldsaresmoothp1eq47}, and \Ref{thmleftinvariantvectorfieldsaresmoothp1eq48} imply that
$\alpha$ is a smooth map from $\func{\image{\phi}}{\WW{}}$ to $\R^n$, that is,
\begin{equation}\label{thmleftinvariantvectorfieldsaresmoothp1eq49}
\alpha\in\banachmapdifclass{\infty}{\R^n}{\R^n}{\func{\image{\phi}}{\WW{}}}{\R^n}.
\end{equation}
\Ref{thmleftinvariantvectorfieldsaresmoothp1eq13}, \Ref{thmleftinvariantvectorfieldsaresmoothp1eq14},
and \Ref{thmleftinvariantvectorfieldsaresmoothp1eq15} imply that,
\begin{equation}\label{thmleftinvariantvectorfieldsaresmoothp1eq50}
\(\cmp{\tanchart{\Lieman{\Liegroup{}}}{\phi}}{\cmp{\avecf{}}{\finv{\psi}}}\)=
\cmp{\[\funcprod{\Injection{\func{\image{\phi}}{\WW{}}}{\funcimage{\phi}}}{\alpha}\]}
{\diagmap{\func{\image{\phi}}{\WW{}}}},
\end{equation}
and additionally, it is trivial that,
\begin{align}
&\Injection{\func{\image{\phi}}{\WW{}}}{\funcimage{\phi}}\in
\banachmapdifclass{\infty}{\R^n}{\R^n}{\func{\image{\phi}}{\WW{}}}{\funcimage{\phi}},
\label{thmleftinvariantvectorfieldsaresmoothp1eq51}\\
&\diagmap{\func{\image{\phi}}{\WW{}}}\in\banachmapdifclass{\infty}{\R^n}{\Cprod{\R^n}{\R^n}}
{\func{\image{\phi}}{\WW{}}}{\Cprod{\func{\image{\phi}}{\WW{}}}{\func{\image{\phi}}{\WW{}}}}.
\label{thmleftinvariantvectorfieldsaresmoothp1eq52}
\end{align}
\Ref{thmleftinvariantvectorfieldsaresmoothp1eq49}, \Ref{thmleftinvariantvectorfieldsaresmoothp1eq50},
\Ref{thmleftinvariantvectorfieldsaresmoothp1eq51}, and \Ref{thmleftinvariantvectorfieldsaresmoothp1eq52} imply that
the local mapping $\cmp{\tanchart{\Lieman{\Liegroup{}}}{\phi}}{\cmp{\avecf{}}{\finv{\psi}}}$ is a smooth map, that is,
\begin{equation}\label{thmleftinvariantvectorfieldsaresmoothp1eq53}
\cmp{\tanchart{\Lieman{\Liegroup{}}}{\phi}}{\cmp{\avecf{}}{\finv{\psi}}}\in
\banachmapdifclass{\infty}{\R^n}{\Cprod{\R^n}{\R^n}}{\func{\image{\phi}}{\WW{}}}{\Cprod{\funcimage{\phi}}{\R^n}}.
\end{equation}
\endp
\end{itemize}
Therefore, for every point $\g{}$ of $\Liegroup{}$, there exists a chart $\psi$ of the manifold $\Lieman{\Liegroup{}}$
around $\g{}$ and a chart $\tanchart{\Lieman{\Liegroup{}}}{\phi}$ of the manifold $\Tanbun{\Lieman{\Liegroup{}}}$
(the tangent-bundle of the underlying manifold of the smooth group $\Liegroup{}$)
around $\func{\avecf{}}{\g{}}$ such that $\avecf{}$ maps the domain of $\psi$ into the domain of $\tanchart{\Lieman{\Liegroup{}}}{\phi}$,
and the composite mapping $\cmp{\tanchart{\Lieman{\Liegroup{}}}{\phi}}{\cmp{\avecf{}}{\finv{\psi}}}$ is smooth. That is.
\begin{align}\label{thmleftinvariantvectorfieldsaresmoothpeq3}
&\Foreach{\g{}}{\G{}}
\Exists{\psi}{\defset{\p{\phi}}{\maxatlas{}}{\g{}\in\domain{\p{\phi}}}}
\Exists{\tanchart{\Lieman{\Liegroup{}}}{\phi}}{\defset{\p{\phi}}{\tanatlas{\Lieman{\Liegroup{}}}}
{\func{\avecf{}}{\g{}}\in\domain{\p{\phi}}}}\cr
&\bigg(\func{\image{\avecf{}}}{\domain{\psi}}\subseteq\domain{\tanchart{\Lieman{\Liegroup{}}}{\phi}},~
\cmp{\tanchart{\Lieman{\Liegroup{}}}{\phi}}{\cmp{\avecf{}}{\finv{\psi}}}\in
\banachmapdifclass{\infty}{\R^n}{\Cprod{\R^n}{\R^n}}{\func{\image{\phi}}{\WW{}}}{\Cprod{\funcimage{\phi}}{\R^n}}\bigg).
\end{align}
This means $\avecf{}$ is a smooth map from the manifold $\Lieman{\Liegroup{}}$ to the manifold
$\Tanbun{\Lieman{\Liegroup{}}}$, that is,
\begin{equation}\label{thmleftinvariantvectorfieldsaresmoothpeq4}
\avecf{}\in\mapdifclass{\infty}{\Lieman{\Liegroup{}}}{\Tanbun{\Lieman{\Liegroup{}}}}.
\end{equation}
\endthm
\theorem\label{thmvectorspaceofleftinvariantvectorfields}
The set of all left-invariant vector-fields on $\Liegroup{}$ is a vector-subspace of
$\Vecf{\Lieman{\Liegroup{}}}{\infty}$ (the vector-space of all
smooth vector-fields on $\Lieman{\Liegroup{}}$). That is, by invoking
the addition and scalar-multiplication operations of the vector-space $\Vecf{\Lieman{\Liegroup{}}}{\infty}$,
\begin{align}
&\empty\neq\Leftinvvf{\Liegroup{}}\subseteq\vecf{\Lieman{\Liegroup{}}}{\infty},\\
&\Foreach{\c}{\R}\Foreach{\opair{\avecf{1}}{\avecf{2}}}{{\Leftinvvf{\Liegroup{}}}^{\times 2}}
\(\c\avecf{1}+\avecf{2}\)\in\Leftinvvf{\Liegroup{}}.
\end{align}
\proof
It can be easily verified that the trivial smooth vector-field on $\Lieman{\Liegroup{}}$ that assigns to each point of $\Liegroup{}$
the neutral element of the tangent-space of $\Lieman{\Liegroup{}}$ at that point, that is $\zerovec{\Vecf{\Lieman{\Liegroup{}}}{\infty}}$
(the neutral element (zero-vector) of the vector-space $\Vecf{\Lieman{\Liegroup{}}}{\infty}$), is a left-invariant vector-field on
$\Liegroup{}$.
\begin{equation}
\zerovec{\Vecf{\Lieman{\Liegroup{}}}{\infty}}\in\Leftinvvf{\Liegroup{}}.
\end{equation}
Hence there exists at least one left-invariant vector-field on $\Liegroup{}$ and thus,
\begin{equation}
\Leftinvvf{\Liegroup{}}\neq\empty.
\end{equation}
Also according to \refthm{thmleftinvariantvectorfieldsaresmooth},
$\Leftinvvf{\Liegroup{}}\subseteq\vecf{\Lieman{\Liegroup{}}}{\infty}$.
\begin{itemize}
\item[$\pr{1}$]
Each $\avecf{1}$ and $\avecf{2}$ is taken as an element of $\Leftinvvf{\Liegroup{}}$, and $\c$ as an arbitrary real number.
Since $\(\c\avecf{1}+\avecf{2}\)\in\vecf{\Lieman{\Liegroup{}}}{\infty}$,
\begin{equation}\label{thmvectorspaceofleftinvariantvectorfieldsp1eq1}
\cmp{\basep{\Lieman{\Liegroup{}}}}{\(\c\avecf{1}+\avecf{2}\)}=\identity{\G{}}.
\end{equation}
Additionally, according to \refdef{defleftinvariantvectorfields} and considering that $\der{\gltrans{\LieG{\Liegroup{}}}{\g{}}}{\Lieman{\Liegroup{}}}
{\Lieman{\Liegroup{}}}$ acts linearly when restricted to  the tangent-space of $\Lieman{\Liegroup{}}$ at any point of $\Liegroup{}$,
\begin{align}\label{thmvectorspaceofleftinvariantvectorfieldsp1eq2}
\Foreach{\g{}}{\G{}}
\Foreach{\point}{\G{}}
&~~~\func{\[\cmp{\(\der{\gltrans{\LieG{\Liegroup{}}}{\g{}}}{\Lieman{\Liegroup{}}}
{\Lieman{\Liegroup{}}}\)}{\(\c\avecf{1}+\avecf{2}\)}\]}{\point}\cr
&=\func{\(\der{\gltrans{\LieG{\Liegroup{}}}{\g{}}}{\Lieman{\Liegroup{}}}
{\Lieman{\Liegroup{}}}\)}{\c\func{\avecf{1}}{\point}+\func{\avecf{2}}{\point}}\cr
&=\c\[\func{\(\der{\gltrans{\LieG{\Liegroup{}}}{\g{}}}{\Lieman{\Liegroup{}}}
{\Lieman{\Liegroup{}}}\)}{\func{\avecf{1}}{\point}}\]+
\func{\(\der{\gltrans{\LieG{\Liegroup{}}}{\g{}}}{\Lieman{\Liegroup{}}}
{\Lieman{\Liegroup{}}}\)}{\func{\avecf{2}}{\point}}\cr
&=\c\func{\avecf{1}}{\g{}\gop{}\point}+\func{\avecf{2}}{\g{}\gop{}\point}\cr
&=\func{\(\c\avecf{1}+\avecf{2}\)}{\g{}\gop{}\point}\cr
&=\func{\[\cmp{\(\c\avecf{1}+\avecf{2}\)}{\gltrans{\LieG{\Liegroup{}}}{\g{}}}\]}{\point},
\end{align}
and thus,
\begin{equation}\label{thmvectorspaceofleftinvariantvectorfieldsp1eq3}
\Foreach{\g{}}{\G{}}
\cmp{\(\der{\gltrans{\LieG{\Liegroup{}}}{\g{}}}{\Lieman{\Liegroup{}}}
{\Lieman{\Liegroup{}}}\)}{\(\c\avecf{1}+\avecf{2}\)}=
\cmp{\(\c\avecf{1}+\avecf{2}\)}{\gltrans{\LieG{\Liegroup{}}}{\g{}}}.
\end{equation}
Based on \refdef{defleftinvariantvectorfields}, \Ref{thmvectorspaceofleftinvariantvectorfieldsp1eq1} and
\Ref{thmvectorspaceofleftinvariantvectorfieldsp1eq3} imply that,
\begin{equation}
\(\c\avecf{1}+\avecf{2}\)\in\Leftinvvf{\Liegroup{}}.
\end{equation}
\endp
\end{itemize}
\endthm
\definition\label{defvectorspaceleftinvariantvectorfields0}
The vector-space obtained by endowing $\Leftinvvf{\Liegroup{}}$ (the set of all left-invariant vector-fields on $\Liegroup{}$)
with the linear-structure inherited from
that of $\Vecf{\Lieman{\Liegroup{}}}{\infty}$ is denoted by $\VLeftinvvf{\Liegroup{}}$, and it is referred to as the
$\quotl$vector-space of left-invariant vector-fields on $\Liegroup{}$$\quotr$.
\endef
\theorem\label{thmleftinvariantvectorfieldsareclosedunderliebracket}
The set of all left-invariant vector-fields is closed under the Lie-bracket operation on $\vecf{\Lieman{\Liegroup{}}}{\infty}$.
This means the Lie-bracket (with respect to $\Lieman{\Liegroup{}}$) of any pair of left-invariant vector-fields on $\Liegroup{}$
is a left-invariant vector-field on $\Liegroup{}$. That is,
\begin{equation}
\Foreach{\opair{\avecf{1}}{\avecf{2}}}
{{\Leftinvvf{\Liegroup{}}}^{\times 2}}
\func{\Liebracket{\Lieman{\Liegroup{}}}}{\binary{\avecf{1}}{\avecf{2}}}\in\Leftinvvf{\Liegroup{}},
\end{equation}
or equivalently,
\begin{equation}
\func{\image{\Liebracket{\Lieman{\Liegroup{}}}}}{\Cprod{\Leftinvvf{\Liegroup{}}}{\Leftinvvf{\Liegroup{}}}}\subseteq
\Leftinvvf{\Liegroup{}}.
\end{equation}
\proof
Each $\avecf{1}$ and $\avecf{2}$ is taken as an arbitrary element of $\Leftinvvf{\Liegroup{}}$.
So according to \refdef{defleftinvariantvectorfields},
\begin{align}
\Foreach{\g{}}{\G{}}\cmp{\(\der{\gltrans{\LieG{\Liegroup{}}}{\g{}}}{\Lieman{\Liegroup{}}}
{\Lieman{\Liegroup{}}}\)}{\avecf{1}}&=\cmp{\avecf{1}}{\gltrans{\LieG{\Liegroup{}}}{\g{}}},
\label{thmleftinvariantvectorfieldsareclosedunderliebracketpeq00}\\
\Foreach{\g{}}{\G{}}\cmp{\(\der{\gltrans{\LieG{\Liegroup{}}}{\g{}}}{\Lieman{\Liegroup{}}}
{\Lieman{\Liegroup{}}}\)}{\avecf{2}}&=\cmp{\avecf{2}}{\gltrans{\LieG{\Liegroup{}}}{\g{}}}.
\label{thmleftinvariantvectorfieldsareclosedunderliebracketpeq0}
\end{align}
According to \refdef{definvliederivative}, \refdef{defliederinv}, \refcor{corliederivativeisalinearisomorphism},
and \refdef{defliebracket},
\begin{align}\label{thmleftinvariantvectorfieldsareclosedunderliebracketpeq1}
\func{\[\func{\Liebracket{\Lieman{\Liegroup{}}}}{\binary{\avecf{1}}{\avecf{2}}}\]}{\IG{}}&=
\func{\[\func{\finv{\(\Lieder{\Lieman{\Liegroup{}}}\)}}
{\aderivation{}}\]}{\IG{}}\cr
&=
\sum_{k=1}^{n}\bigg(\func{\[\func{\aderivation{}}{\cf_k}\]}{\IG{}}\bigg)\bu{k},~~~~
\end{align}
where, $\aderivation{}$ is the $\infty$-derivation on $\Lieman{\Liegroup{}}$ defined by,
\begin{equation}\label{thmleftinvariantvectorfieldsareclosedunderliebracketpeq2}
\aderivation{}:=\cmp{\func{\Lieder{\Lieman{\Liegroup{}}}}{\avecf{1}}}{\func{\Lieder{\Lieman{\Liegroup{}}}}{\avecf{2}}}-
\cmp{\func{\Lieder{\Lieman{\Liegroup{}}}}{\avecf{2}}}{\func{\Lieder{\Lieman{\Liegroup{}}}}{\avecf{1}}},
\end{equation}
and,
\begin{equation}\label{thmleftinvariantvectorfieldsareclosedunderliebracketpeq3}
\Foreach{k}{\seta{\suc{1}{n}}}
\bu{k}:=\func{\finv{\[\tanspaceiso{\IG{}}{\Man{}}{\phi}\]}}{\Eucbase{n}{k}},
\end{equation}
and $\phi$ is an element of $\defset{\p{\phi}}{\maxatlas{}}{\IG{}\in\domain{\p{\phi}},~\func{\p{\phi}}{\IG{}}=\zerovec{}}$,
and $\suc{\cf_1}{\cf_n}$ are elements of
$\suc{\fextension{\Man{}}{\IG{}}{\phi}{\projection{n}{1}}}{\fextension{\Man{}}{\IG{}}{\phi}{\projection{n}{n}}}$, respectively.\\
According to \refdef{defliebracket}, $\func{\Liebracket{\Lieman{\Liegroup{}}}}{\binary{\avecf{1}}{\avecf{2}}}$ is
a smooth vector-field on $\Lieman{\Liegroup{}}$ and hence an element of $\Func{\G{}}{\tanbun{\Lieman{\Liegroup{}}}}$.
So trivially $\func{\[\func{\Liebracket{\Lieman{\Liegroup{}}}}{\binary{\avecf{1}}{\avecf{2}}}\]}{\IG{}}$ belongs to
the tangent-space of $\Lieman{\Liegroup{}}$ at $\IG{}$.
\begin{align}
\func{\Liebracket{\Lieman{\Liegroup{}}}}{\binary{\avecf{1}}{\avecf{2}}}&\in
\Func{\G{}}{\tanbun{\Lieman{\Liegroup{}}}},
\label{thmleftinvariantvectorfieldsareclosedunderliebracketpeq4}\\
\func{\[\func{\Liebracket{\Lieman{\Liegroup{}}}}{\binary{\avecf{1}}{\avecf{2}}}\]}{\IG{}}&\in
\tanspace{\IG{}}{\Lieman{\Liegroup{}}}.
\label{thmleftinvariantvectorfieldsareclosedunderliebracketpeq5}
\end{align}
According to \refdef{defliederivative}, and by invoking the inherent linear-structures of the vector-spaces
$\Lmapdifclass{\infty}{\Lieman{\Liegroup{}}}{\RR}$ and
$\VLin{\Lmapdifclass{\infty}{\Lieman{\Liegroup{}}}{\RR}}{\Lmapdifclass{\infty}{\Lieman{\Liegroup{}}}{\RR}}$,
\begin{align}\label{thmleftinvariantvectorfieldsareclosedunderliebracketpeq6}
&\begin{aligned}
\Foreach{\cf}{\mapdifclass{\infty}{\Lieman{\Liegroup{}}}{\RR}}
\Foreach{\point}{\G{}}
\end{aligned}\cr
&\begin{aligned}
\func{\[\func{\aderivation{}}{\cf}\]}{\point}&=
\func{\[\func{\(\cmp{\func{\Lieder{\Lieman{\Liegroup{}}}}{\avecf{1}}}{\func{\Lieder{\Lieman{\Liegroup{}}}}{\avecf{2}}}-
\cmp{\func{\Lieder{\Lieman{\Liegroup{}}}}{\avecf{2}}}{\func{\Lieder{\Lieman{\Liegroup{}}}}{\avecf{1}}}\)}{\cf}\]}{\point}\cr
&=\func{\bigg[\func{\func{\Lieder{\Lieman{\Liegroup{}}}}{\avecf{1}}}{\func{\[\func{\Lieder{\Lieman{\Liegroup{}}}}{\avecf{2}}\]}{\cf}}-
\func{\func{\Lieder{\Lieman{\Liegroup{}}}}{\avecf{2}}}{\func{\[\func{\Lieder{\Lieman{\Liegroup{}}}}{\avecf{1}}\]}{\cf}}\bigg]}{\point}\cr
&=\func{\bigg(\func{\func{\Lieder{\Lieman{\Liegroup{}}}}{\avecf{1}}}
{\func{\[\func{\Lieder{\Lieman{\Liegroup{}}}}{\avecf{2}}\]}{\cf}}\bigg)}{\point}-
\func{\bigg(\func{\func{\Lieder{\Lieman{\Liegroup{}}}}{\avecf{2}}}
{\func{\[\func{\Lieder{\Lieman{\Liegroup{}}}}{\avecf{1}}\]}{\cf}}\bigg)}{\point}\cr
&=\func{\(\cmp{\Rder{\[\cmp{\Rder{\cf}{\Lieman{\Liegroup{}}}}{\avecf{2}}\]}{\Lieman{\Liegroup{}}}}{\avecf{1}}\)}{\point}-
\func{\(\cmp{\Rder{\[\cmp{\Rder{\cf}{\Lieman{\Liegroup{}}}}{\avecf{1}}\]}{\Lieman{\Liegroup{}}}}{\avecf{2}}\)}{\point}\cr
&=\func{\(\Rder{\[\cmp{\Rder{\cf}{\Lieman{\Liegroup{}}}}{\avecf{2}}\]}{\Lieman{\Liegroup{}}}\)}{\func{\avecf{1}}{\point}}-
\func{\(\Rder{\[\cmp{\Rder{\cf}{\Lieman{\Liegroup{}}}}{\avecf{1}}\]}{\Lieman{\Liegroup{}}}\)}{\func{\avecf{2}}{\point}}.
\end{aligned}
\end{align}
Particularly,
\begin{align}\label{thmleftinvariantvectorfieldsareclosedunderliebracketpeq7}
&\Foreach{k}{\seta{\suc{1}{n}}}\cr
&\func{\[\func{\aderivation{}}{\cf_k}\]}{\IG{}}=
\func{\(\Rder{\[\cmp{\Rder{\cf_{k}}{\Lieman{\Liegroup{}}}}{\avecf{2}}\]}{\Lieman{\Liegroup{}}}\)}{\func{\avecf{1}}{\IG{}}}-
\func{\(\Rder{\[\cmp{\Rder{\cf_{k}}{\Lieman{\Liegroup{}}}}{\avecf{1}}\]}{\Lieman{\Liegroup{}}}\)}{\func{\avecf{2}}{\IG{}}}.
\end{align}
\begin{itemize}
\item[$\pr{1}$]
$\g{}$ is taken as an arbitrary element of $\G{}$.
Since $\gltrans{\LieG{\Liegroup{}}}{\g{}}$ is an $\infty$-diffeomorphism of $\Lieman{\Liegroup{}}$ and
$\gltrans{\LieG{\Liegroup{}}}{\invg{\g{}}{}}=\finv{\(\gltrans{\LieG{\Liegroup{}}}{\g{}}\)}$,
according to \Ref{eqdiffeomorphismdifrule},
\begin{align}\label{thmleftinvariantvectorfieldsareclosedunderliebracketp1eq1}
\func{\avecf{1}}{\IG{}}&=\func{\identity{\tanbun{\Lieman{\Liegroup{}}}}}{\func{\avecf{1}}{\IG{}}}\cr
&=\func{\[\cmp{\(\der{\gltrans{\LieG{\Liegroup{}}}{\invg{\g{}}{}}}{\Lieman{\Liegroup{}}}
{\Lieman{\Liegroup{}}}\)}{\(\der{\gltrans{\LieG{\Liegroup{}}}{\g{}}}{\Lieman{\Liegroup{}}}
{\Lieman{\Liegroup{}}}\)}\]}{\func{\avecf{1}}{\IG{}}}\cr
&=\func{\(\der{\gltrans{\LieG{\Liegroup{}}}{\ginv{\g{}}{}}}{\Lieman{\Liegroup{}}}
{\Lieman{\Liegroup{}}}\)}{\func{\[\der{\gltrans{\LieG{\Liegroup{}}}{\g{}}}{\Lieman{\Liegroup{}}}
{\Lieman{\Liegroup{}}}\]}{\func{\avecf{1}}{\IG{}}}}.
\end{align}
Moreover, since $\avecf{1}$ is a left-invariant vector-field on $\Liegroup{}$, according to \reflem{lemleftinvariantvectorfieldsequiv0},
\begin{equation}\label{thmleftinvariantvectorfieldsareclosedunderliebracketp1eq2}
\func{\(\der{\gltrans{\LieG{\Liegroup{}}}{\g{}}}{\Lieman{\Liegroup{}}}
{\Lieman{\Liegroup{}}}\)}{\func{\avecf{1}}{\IG{}}}=\func{\avecf{1}}{\g{}}.
\end{equation}
Therefore,
\begin{equation}\label{thmleftinvariantvectorfieldsareclosedunderliebracketp1eq3}
\func{\avecf{1}}{\IG{}}=
\func{\(\der{\gltrans{\LieG{\Liegroup{}}}{\invg{\g{}}{}}}{\Lieman{\Liegroup{}}}
{\Lieman{\Liegroup{}}}\)}{\func{\avecf{1}}{\g{}}}.
\end{equation}
Similarly, since $\avecf{2}$ is also a left-invariant vector-field on $\Liegroup{}$,
\begin{equation}\label{thmleftinvariantvectorfieldsareclosedunderliebracketp1eq4}
\func{\avecf{2}}{\IG{}}=
\func{\(\der{\gltrans{\LieG{\Liegroup{}}}{\invg{\g{}}{}}}{\Lieman{\Liegroup{}}}
{\Lieman{\Liegroup{}}}\)}{\func{\avecf{2}}{\g{}}}.
\end{equation}
According to \Ref{thmleftinvariantvectorfieldsareclosedunderliebracketpeq1} and
\Ref{thmleftinvariantvectorfieldsareclosedunderliebracketpeq7}, and considering that the differential of
$\gltrans{\LieG{\Liegroup{}}}{\g{}}$ operates linearly when restricted to $\tanspace{\IG{}}{\Lieman{\Liegroup{}}}$,
\begin{align}\label{thmleftinvariantvectorfieldsareclosedunderliebracketp1eq5}
&~~~\func{\(\der{\gltrans{\LieG{\Liegroup{}}}{\g{}}}{\Lieman{\Liegroup{}}}
{\Lieman{\Liegroup{}}}\)}{\func{\[\func{\Liebracket{\Lieman{\Liegroup{}}}}
{\binary{\avecf{1}}{\avecf{2}}}\]}{\IG{}}}\cr
&=\sum_{k=1}^{n}\bigg(\func{\[\func{\aderivation{}}{\cf_k}\]}{\IG{}}\bigg)
\[\func{\(\der{\gltrans{\LieG{\Liegroup{}}}{\g{}}}{\Lieman{\Liegroup{}}}
{\Lieman{\Liegroup{}}}\)}{\bu{k}}\]\cr
&=\sum_{k=1}^{n}\bigg[\func{\(\Rder{\[\cmp{\Rder{\cf_{k}}{\Lieman{\Liegroup{}}}}{\avecf{2}}\]}{\Lieman{\Liegroup{}}}\)}{\func{\avecf{1}}{\IG{}}}-
\func{\(\Rder{\[\cmp{\Rder{\cf_{k}}{\Lieman{\Liegroup{}}}}{\avecf{1}}\]}{\Lieman{\Liegroup{}}}\)}{\func{\avecf{2}}{\IG{}}}\bigg]
\[\func{\(\der{\gltrans{\LieG{\Liegroup{}}}{\g{}}}{\Lieman{\Liegroup{}}}
{\Lieman{\Liegroup{}}}\)}{\bu{k}}\].\cr
&{}
\end{align}
According to \Ref{thmleftinvariantvectorfieldsareclosedunderliebracketp1eq3} and
\Ref{eqspecialchainrule1} (the special chain-rule of differentiation),
\begin{align}\label{thmleftinvariantvectorfieldsareclosedunderliebracketp1eq6}
\func{\(\Rder{\[\cmp{\Rder{\cf_{k}}{\Lieman{\Liegroup{}}}}{\avecf{2}}\]}
{\Lieman{\Liegroup{}}}\)}{\func{\avecf{1}}{\IG{}}}&=
\func{\(\Rder{\[\cmp{\Rder{\cf_{k}}{\Lieman{\Liegroup{}}}}{\avecf{2}}\]}
{\Lieman{\Liegroup{}}}\)}{\func{\[\der{\gltrans{\LieG{\Liegroup{}}}{\invg{\g{}}{}}}{\Lieman{\Liegroup{}}}
{\Lieman{\Liegroup{}}}\]}{\func{\avecf{1}}{\g{}}}}\cr
&=\func{\(\cmp{\(\Rder{\[\cmp{\Rder{\cf_{k}}{\Lieman{\Liegroup{}}}}{\avecf{2}}\]}
{\Lieman{\Liegroup{}}}\)}{\[\der{\gltrans{\LieG{\Liegroup{}}}{\invg{\g{}}{}}}{\Lieman{\Liegroup{}}}
{\Lieman{\Liegroup{}}}\]}\)}{\func{\avecf{1}}{\g{}}}\cr
&=\func{\(\Rder{\[\cmp{\cmp{\Rder{\cf_{k}}{\Lieman{\Liegroup{}}}}{\avecf{2}}}
{\gltrans{\LieG{\Liegroup{}}}{\invg{\g{}}{}}}\]}{\Lieman{\Liegroup{}}}\)}{\func{\avecf{1}}{\g{}}}.
\end{align}
Additionally, according to \Ref{thmleftinvariantvectorfieldsareclosedunderliebracketpeq00},
and again \Ref{eqspecialchainrule1} (the special chain-rule of differentiation),
\begin{align}\label{thmleftinvariantvectorfieldsareclosedunderliebracketp1eq7}
\cmp{\Rder{\cf_{k}}{\Lieman{\Liegroup{}}}}{\(\cmp{\avecf{2}}
{\gltrans{\LieG{\Liegroup{}}}{\invg{\g{}}{}}}\)}&=
\cmp{\Rder{\cf_{k}}{\Lieman{\Liegroup{}}}}{\[\cmp{\(\der{\gltrans{\LieG{\Liegroup{}}}{\invg{\g{}}{}}}{\Lieman{\Liegroup{}}}
{\Lieman{\Liegroup{}}}\)}{\avecf{2}}\]}\cr
&=\cmp{\[\cmp{\Rder{\cf_{k}}{\Lieman{\Liegroup{}}}}{\(\der{\gltrans{\LieG{\Liegroup{}}}{\invg{\g{}}{}}}{\Lieman{\Liegroup{}}}
{\Lieman{\Liegroup{}}}\)}\]}{\avecf{2}}\cr
&=\cmp{\[\Rder{\(\cmp{\cf_{k}}{\gltrans{\LieG{\Liegroup{}}}{\invg{\g{}}{}}}\)}{\Lieman{\Liegroup{}}}\]}{\avecf{2}}.
\end{align}
\endp
Combining \Ref{thmleftinvariantvectorfieldsareclosedunderliebracketp1eq6} and
\Ref{thmleftinvariantvectorfieldsareclosedunderliebracketp1eq7} yields,
\begin{equation}\label{thmleftinvariantvectorfieldsareclosedunderliebracketp1eq8}
\func{\(\Rder{\[\cmp{\Rder{\cf_{k}}{\Lieman{\Liegroup{}}}}{\avecf{2}}\]}
{\Lieman{\Liegroup{}}}\)}{\func{\avecf{1}}{\IG{}}}=
\func{\(\Rder{\[\cmp{\Rder{\(\cmp{\cf_{k}}{\gltrans{\LieG{\Liegroup{}}}{\invg{\g{}}{}}}\)}
{\Lieman{\Liegroup{}}}}{\avecf{2}}\]}{\Lieman{\Liegroup{}}}\)}{\func{\avecf{1}}{\g{}}}.
\end{equation}
In completely a similar manner,
\begin{equation}\label{thmleftinvariantvectorfieldsareclosedunderliebracketp1eq9}
\func{\(\Rder{\[\cmp{\Rder{\cf_{k}}{\Lieman{\Liegroup{}}}}{\avecf{1}}\]}
{\Lieman{\Liegroup{}}}\)}{\func{\avecf{2}}{\IG{}}}=
\func{\(\Rder{\[\cmp{\Rder{\(\cmp{\cf_{k}}{\gltrans{\LieG{\Liegroup{}}}{\invg{\g{}}{}}}\)}
{\Lieman{\Liegroup{}}}}{\avecf{1}}\]}{\Lieman{\Liegroup{}}}\)}{\func{\avecf{2}}{\g{}}}.
\end{equation}	
Using \Ref{thmleftinvariantvectorfieldsareclosedunderliebracketp1eq8} and
\Ref{thmleftinvariantvectorfieldsareclosedunderliebracketp1eq9} in
\Ref{thmleftinvariantvectorfieldsareclosedunderliebracketp1eq5} yields,
\begin{align}\label{thmleftinvariantvectorfieldsareclosedunderliebracketp1eq10}
&~~~\func{\(\der{\gltrans{\LieG{\Liegroup{}}}{\g{}}}{\Lieman{\Liegroup{}}}
{\Lieman{\Liegroup{}}}\)}{\func{\[\func{\Liebracket{\Lieman{\Liegroup{}}}}
{\binary{\avecf{1}}{\avecf{2}}}\]}{\IG{}}}\cr
&=\sum_{k=1}^{n}\bigg[\func{\(\Rder{\[\cmp{\Rder{\(\cmp{\cf_{k}}{\gltrans{\LieG{\Liegroup{}}}{\invg{\g{}}{}}}\)}
{\Lieman{\Liegroup{}}}}{\avecf{2}}\]}{\Lieman{\Liegroup{}}}\)}{\func{\avecf{1}}{\g{}}}-
\func{\(\Rder{\[\cmp{\Rder{\(\cmp{\cf_{k}}{\gltrans{\LieG{\Liegroup{}}}{\invg{\g{}}{}}}\)}
{\Lieman{\Liegroup{}}}}{\avecf{1}}\]}{\Lieman{\Liegroup{}}}\)}{\func{\avecf{2}}{\g{}}}\bigg]\bw{k},\cr
&{}
\end{align}
where,
\begin{equation}\label{thmleftinvariantvectorfieldsareclosedunderliebracketp1eq11}
\Foreach{k}{\seta{\suc{1}{n}}}
\bw{k}:=\[\func{\(\der{\gltrans{\LieG{\Liegroup{}}}{\g{}}}{\Lieman{\Liegroup{}}}
{\Lieman{\Liegroup{}}}\)}{\bu{k}}\].
\end{equation}
According to \Ref{eqcharttransfer}, $\cmp{\phi}{\gltrans{\LieG{\Liegroup{}}}{\invg{\g{}}{}}}$ is a chart of $\Man{}$
around $\func{\finv{\[\gltrans{\LieG{\Liegroup{}}}{\invg{\g{}}{}}\]}}{\IG{}}$, and since
$\finv{\[\gltrans{\LieG{\Liegroup{}}}{\invg{\g{}}{}}\]}=\gltrans{\LieG{\Liegroup{}}}{\g{}}$,
$\cmp{\phi}{\gltrans{\LieG{\Liegroup{}}}{\invg{\g{}}{}}}$ is a chart of $\Man{}$
around $\g{}$. Also since $\func{\phi}{\IG{}}=\zerovec{}$, clearly this chart is centered at $\g{}$, that is,
$\func{\[\cmp{\phi}{\gltrans{\LieG{\Liegroup{}}}{\invg{\g{}}{}}}\]}{\g{}}=\zerovec{}$. So briefly,
\begin{equation}\label{thmleftinvariantvectorfieldsareclosedunderliebracketp1eq12}
\psi:=\cmp{\phi}{\gltrans{\LieG{\Liegroup{}}}{\invg{\g{}}{}}}\in
\defset{\p{\phi}}{\maxatlas{}}{\g{}\in\domain{\p{\phi}},~\func{\p{\phi}}{\g{}}=\zerovec{}}.
\end{equation}
In addition, since for every $k\in\seta{\suc{1}{n}}$, $\cf_k\in\fextension{\Man{}}{\IG{}}{\phi}{\projection{n}{k}}$,
according to \refdef{defsmoothextensionoflocalmaps} it is trivial that each
$\cmp{\cf_{k}}{\gltrans{\LieG{\Liegroup{}}}{\invg{\g{}}{}}}$ is a
smooth extension of $\cmp{\projection{n}{k}}{\psi}$ on $\Man{}$ fixed at $\g{}$.
That is,
\begin{align}\label{thmleftinvariantvectorfieldsareclosedunderliebracketp1eq13}
\Foreach{k}{\seta{\suc{1}{n}}}
\hf_{k}\in
\fextension{\Man{}}{\g{}}{\psi}{\projection{n}{k}},
\end{align}
where,
\begin{equation}\label{thmleftinvariantvectorfieldsareclosedunderliebracketp1eq14}
\Foreach{k}{\seta{\suc{1}{n}}}
\hf_{k}:=\cmp{\cf_{k}}{\gltrans{\LieG{\Liegroup{}}}{\invg{\g{}}{}}}.
\end{equation}
Furthermore, since $\phi$ is a chart of $\Man{}$ around $\IG{}$ and $\psi$ is a chart of $\Man{}$
around $\g{}=\func{\gltrans{\LieG{\Liegroup{}}}{\g{}}}{\IG{}}$,
according to \Ref{eqdefinitionofdifferentialofamap} (definition of the differential operator),
\begin{align}\label{thmleftinvariantvectorfieldsareclosedunderliebracketp1eq15}
&\begin{aligned}
\Foreach{k}{\seta{\suc{1}{n}}}
\end{aligned}\cr
&\begin{aligned}
\bw{k}&=\func{\(\der{\gltrans{\LieG{\Liegroup{}}}{\g{}}}{\Lieman{\Liegroup{}}}
{\Lieman{\Liegroup{}}}\)}{\func{\finv{\[\tanspaceiso{\IG{}}{\Man{}}{\phi}\]}}{\Eucbase{n}{k}}}\cr
&=
\func{\bigg[\cmp{\finv{\(\tanspaceiso{\g{}}{\Lieman{\Liegroup{}}}{\psi}\)}}
{\cmp{\bigg(\func{\[\banachder{\(\cmp{\psi}{\cmp{\gltrans{\LieG{\Liegroup{}}}{\g{}}}{\finv{\phi}}}\)}{\R^n}{\R^n}\]}
{\func{\phi}{\IG{}}}\bigg)}{\tanspaceiso{\IG{}}{\Man{}}{\phi}}}\bigg]}
{\func{\finv{\[\tanspaceiso{\IG{}}{\Man{}}{\phi}\]}}{\Eucbase{n}{k}}}\cr
&=\func{\bigg[\cmp{\finv{\(\tanspaceiso{\g{}}{\Lieman{\Liegroup{}}}{\psi}\)}}
{\bigg(\func{\[\banachder{\(\cmp{\psi}{\cmp{\gltrans{\LieG{\Liegroup{}}}{\g{}}}{\finv{\phi}}}\)}{\R^n}{\R^n}\]}
{\func{\phi}{\IG{}}}\bigg)}\bigg]}{\Eucbase{n}{k}}.
\end{aligned}\cr
&{}
\end{align}
In addition, \Ref{thmleftinvariantvectorfieldsareclosedunderliebracketp1eq12} implies that,
$\cmp{\psi}{\cmp{\gltrans{\LieG{\Liegroup{}}}{\g{}}}{\finv{\phi}}}$ is the restriction of
the identity-map on $\R^n$ to the open set $\funcimage{\phi}$ of $\R^n$. That is,
\begin{equation}\label{thmleftinvariantvectorfieldsareclosedunderliebracketp1eq16}
\cmp{\psi}{\cmp{\gltrans{\LieG{\Liegroup{}}}{\g{}}}{\finv{\phi}}}=
\func{\res{\identity{\R^n}}}{\funcimage{\phi}}.
\end{equation}
Therefore, the derived map of
$\cmp{\psi}{\cmp{\gltrans{\LieG{\Liegroup{}}}{\g{}}}{\finv{\phi}}}$
is clearly the constant map on $\funcimage{\phi}$ assigning the identy-map on $\R^n$ to each point of it. That is,
\begin{equation}\label{thmleftinvariantvectorfieldsareclosedunderliebracketp1eq17}
\Foreach{\x}{\funcimage{\phi}}
\func{\[\banachder{\(\cmp{\psi}{\cmp{\gltrans{\LieG{\Liegroup{}}}{\g{}}}{\finv{\phi}}}\)}{\R^n}{\R^n}\]}{\x}=
\identity{\R^n}.
\end{equation}
\Ref{thmleftinvariantvectorfieldsareclosedunderliebracketp1eq15} and
\Ref{thmleftinvariantvectorfieldsareclosedunderliebracketp1eq17} imply that,
\begin{align}\label{thmleftinvariantvectorfieldsareclosedunderliebracketp1eq18}
\Foreach{k}{\seta{\suc{1}{n}}}
\bw{k}&=\func{\finv{\(\tanspaceiso{\g{}}{\Lieman{\Liegroup{}}}{\psi}\)}}{\Eucbase{n}{k}}.
\end{align}
\Ref{thmleftinvariantvectorfieldsareclosedunderliebracketpeq6},
\Ref{thmleftinvariantvectorfieldsareclosedunderliebracketp1eq10},
\Ref{thmleftinvariantvectorfieldsareclosedunderliebracketp1eq14}, and
\Ref{thmleftinvariantvectorfieldsareclosedunderliebracketp1eq18} imply,
\begin{align}\label{thmleftinvariantvectorfieldsareclosedunderliebracketp1eq19}
&~~~\func{\(\der{\gltrans{\LieG{\Liegroup{}}}{\g{}}}{\Lieman{\Liegroup{}}}
{\Lieman{\Liegroup{}}}\)}{\func{\[\func{\Liebracket{\Lieman{\Liegroup{}}}}
{\binary{\avecf{1}}{\avecf{2}}}\]}{\IG{}}}\cr
&=\sum_{k=1}^{n}\bigg[\func{\(\Rder{\[\cmp{\Rder{\hf_{k}}
{\Lieman{\Liegroup{}}}}{\avecf{2}}\]}{\Lieman{\Liegroup{}}}\)}{\func{\avecf{1}}{\g{}}}-
\func{\(\Rder{\[\cmp{\Rder{\hf_{k}}
{\Lieman{\Liegroup{}}}}{\avecf{1}}\]}{\Lieman{\Liegroup{}}}\)}{\func{\avecf{2}}{\g{}}}\bigg]
\[\func{\finv{\(\tanspaceiso{\g{}}{\Lieman{\Liegroup{}}}{\psi}\)}}{\Eucbase{n}{k}}\]\cr
&=\sum_{k=1}^{n}\func{\[\func{\aderivation{}}{\hf_{k}}\]}{\g{}}
\[\func{\finv{\(\tanspaceiso{\g{}}{\Lieman{\Liegroup{}}}{\psi}\)}}{\Eucbase{n}{k}}\].\cr
&{}
\end{align}
So since $\psi$ is a chart of $\Man{}$ centered at $\g{}$ (\Ref{thmleftinvariantvectorfieldsareclosedunderliebracketp1eq12}),
and each $\hf_{k}$ is a smooth extension of $\cmp{\projection{n}{k}}{\psi}$ on $\Man{}$ fixed at $\g{}$
(\Ref{thmleftinvariantvectorfieldsareclosedunderliebracketp1eq13}), based on
\refdef{definvliederivative}, \refdef{defliederinv}, \refcor{corliederivativeisalinearisomorphism},
and \refdef{defliebracket},
\Ref{thmleftinvariantvectorfieldsareclosedunderliebracketpeq2} and
\Ref{thmleftinvariantvectorfieldsareclosedunderliebracketp1eq19} yield,
\begin{align}\label{thmleftinvariantvectorfieldsareclosedunderliebracketp1eq20}
\func{\(\der{\gltrans{\LieG{\Liegroup{}}}{\g{}}}{\Lieman{\Liegroup{}}}
{\Lieman{\Liegroup{}}}\)}{\func{\[\func{\Liebracket{\Lieman{\Liegroup{}}}}
{\binary{\avecf{1}}{\avecf{2}}}\]}{\IG{}}}&=
\func{\[\func{\finv{\(\Lieder{\Lieman{\Liegroup{}}}\)}}
{\aderivation{}}\]}{\g{}}\cr
&=\func{\[\func{\Liebracket{\Lieman{\Liegroup{}}}}{\binary{\avecf{1}}{\avecf{2}}}\]}{\g{}}.
\end{align}
\endp
\end{itemize}
Therefore,
\begin{equation}\label{thmleftinvariantvectorfieldsareclosedunderliebracketpeq8}
\Foreach{\g{}}{\G{}}
\func{\(\der{\gltrans{\LieG{\Liegroup{}}}{\g{}}}{\Lieman{\Liegroup{}}}
{\Lieman{\Liegroup{}}}\)}{\func{\[\func{\Liebracket{\Lieman{\Liegroup{}}}}
{\binary{\avecf{1}}{\avecf{2}}}\]}{\IG{}}}=
\func{\[\func{\Liebracket{\Lieman{\Liegroup{}}}}{\binary{\avecf{1}}{\avecf{2}}}\]}{\g{}}.
\end{equation}
According to \reflem{lemleftinvariantvectorfieldsequiv0},
\Ref{thmleftinvariantvectorfieldsareclosedunderliebracketpeq4},
\Ref{thmleftinvariantvectorfieldsareclosedunderliebracketpeq5}, and
\Ref{thmleftinvariantvectorfieldsareclosedunderliebracketpeq8}
imply that $\func{\Liebracket{\Lieman{\Liegroup{}}}}{\binary{\avecf{1}}{\avecf{2}}}$
is a left-invariant vector-field on $\Liegroup{}$.
\begin{equation}\label{thmleftinvariantvectorfieldsareclosedunderliebracketpeq9}
\func{\Liebracket{\Lieman{\Liegroup{}}}}{\binary{\avecf{1}}{\avecf{2}}}\in
\Leftinvvf{\Liegroup{}}.
\end{equation}
\endthm
\corollary
The set of all left-invariant vector-fields on $\Liegroup{}$ is a Lie-subalgebra of the
Lie-algebra of smooth vector-fields on $\Lieman{\Liegroup{}}$. That is,
\begin{equation}
\Leftinvvf{\Liegroup{}}\in\sublie{\Lievecf{\Lieman{\Liegroup{}}}{\infty}}.
\end{equation}
Consequently, $\subspace{\Lievecf{\Lieman{\Liegroup{}}}{\infty}}{\Leftinvvf{\Liegroup{}}}$
is a Lie-algebra (over the field $\R$). Furthermore, the Lie-operation of
$\subspace{\Lievecf{\Lieman{\Liegroup{}}}{\infty}}{\Leftinvvf{\Liegroup{}}}$
clearly coincides with the Lie-bracket operation on $\vecf{\Lieman{\Liegroup{}}}{\infty}$
when its domain of operation is restricted to the set of all left-invariant vector-fields on $\Liegroup{}$.
\proof
It is an immediate consequence of
\refdef{defliesubalgebra}, \refcor{corliesubalgebra}, \refdef{defliealgebraofvectorfields},
\refthm{thmvectorspaceofleftinvariantvectorfields}, and
\refthm{thmleftinvariantvectorfieldsareclosedunderliebracket}
\endcor
\definition\label{defliealgebraofliegroup}
The Lie-algebra $\subspace{\Lievecf{\Lieman{\Liegroup{}}}{\infty}}{\Leftinvvf{\Liegroup{}}}$
is referred to as the $\quotl$canonical Lie-algebra of the smooth group $\Liegroup{}$ (via the left-invariant vector-fields)$\quotr$,
and is denoted by $\LiegroupLiealgebra{\Liegroup{}}$. The Lie-operation of the Lie-algebra
$\LiegroupLiealgebra{\Liegroup{}}$ is denoted by $\Lliebracket{\Liegroup{}}$,that is,
\begin{equation}
\Lliebracket{\Liegroup{}}:=\func{\res{\Liebracket{\Lieman{\Liegroup{}}}}}{\Cprod{\Leftinvvf{\Liegroup{}}}{\Leftinvvf{\Liegroup{}}}}.
\end{equation}
which is simply
referred to as the $\quotl$Lie-bracket operation on the set of all left-invariant vector-fields on $\Liegroup{}$$\quotr$.
Also, for every $\avecf{1}$ and $\avecf{2}$ in $\Leftinvvf{\Liegroup{}}$,
$\func{\Lliebracket{\Liegroup{}}}{\binary{\avecf{1}}{\avecf{2}}}$ can be denoted alternatively by
$\liebracket{\avecf{1}}{\avecf{2}}{\Liegroup{}}$.
\endef
\definition\label{deflinvvftanspacecorrespondence}
The mapping $\function{\liegvftan{\Liegroup{}}}{\Leftinvvf{\Liegroup{}}}{\tanspace{\IG{}}{\Lieman{\Liegroup{}}}}$
is defined as,
\begin{equation}
\Foreach{\avecf{}}{\Leftinvvf{\Liegroup{}}}
\func{\liegvftan{\Liegroup{}}}{\avecf{}}\eqdef\func{\avecf{}}{\IG{}}.
\end{equation}
It is called the $\quotl$canonical correspondence between the set of all left-invariant vector-fields on $\Liegroup{}$
and the tangent-space of $\Lieman{\Liegroup{}}$ at the identity element of $\LieG{\Liegroup{}}$$\quotr$.
\endef 
\corollary\label{corlinvvftanspalinearcecorrespondence}
$\liegvftan{\Liegroup{}}$ is a linear-isomorphism from $\VLeftinvvf{\Liegroup{}}$ to $\Tanspace{\IG{}}{\Lieman{\Liegroup{}}}$. That is,
\begin{align}
&\Foreach{\opair{\avecf{1}}{\avecf{2}}}{\Cprod{\Leftinvvf{\Liegroup{}}}{\Leftinvvf{\Liegroup{}}}}
\Foreach{\c}{\R}\cr
&\func{\liegvftan{\Liegroup{}}}{\c\avecf{1}+\avecf{2}}=\c\func{\liegvftan{\Liegroup{}}}{\avecf{1}}+
\func{\liegvftan{\Liegroup{}}}{\avecf{2}}.
\end{align}
\proof
The bijectivity of $\liegvftan{\Liegroup{}}$ is an immediate result of \reflem{lemleftinvariantvectorfieldsequiv0},
and its linear behavior is evident when the linear-structures of
$\Leftinvvf{\Liegroup{}}$ and $\tanspace{\IG{}}{\Lieman{\Liegroup{}}}$ are invoked.
\endcor
\definition\label{deftangentspaceliebracketofliegroup}
The mapping $\function{\tanliebracket{\Liegroup{}}}{\Cprod{\tanspace{\IG{}}{\Lieman{\Liegroup{}}}}{\tanspace{\IG{}}{\Lieman{\Liegroup{}}}}}
{\tanspace{\IG{}}{\Lieman{\Liegroup{}}}}$ is defined as,
\begin{align}
\Foreach{\opair{\avec{1}}{\avec{2}}}{\Cprod{\tanspace{\IG{}}{\Lieman{\Liegroup{}}}}{\tanspace{\IG{}}{\Lieman{\Liegroup{}}}}}
\func{\tanliebracket{\Liegroup{}}}{\binary{\avec{1}}{\avec{2}}}&\eqdef
\func{\Lliebracket{\Liegroup{}}}{\binary{\func{\finv{\liegvftan{\Liegroup{}}}}{\avec{1}}}{\func{\finv{\liegvftan{\Liegroup{}}}}{\avec{2}}}}\cr
&=\liebracket{\func{\finv{\liegvftan{\Liegroup{}}}}{\avec{1}}}
{\func{\finv{\liegvftan{\Liegroup{}}}}{\avec{2}}}{\Lieman{\Liegroup{}}}.
\end{align}
It is called the $\quotl$Lie-bracket operation on the tangent-space of $\Lieman{\Liegroup{}}$ at the identity
element of $\LieG{\Liegroup{}}$$\quotr$.
\endef
\corollary\label{corliealgebraoftangentspaceatidentity}
$\opair{\Tanspace{\IG{}}{\Lieman{\Liegroup{}}}}{\tanliebracket{\Liegroup{}}}$ is a Lie-algebra.
\proof
According to \refdef{defliealgebraofliegroup}, \refcor{corlinvvftanspalinearcecorrespondence},
\refcor{corlinvvftanspalinearcecorrespondence} and
\refdef{deftangentspaceliebracketofliegroup}, it is clear.
Simply, the linear-isomorphism $\liegvftan{\Liegroup{}}$ transfers the Lie-algebra strucure of
the set of all left-invariant vector-fields on $\Liegroup{}$ to the tangent-space of $\Lieman{\Liegroup{}}$
at the identity element of $\LieG{\Liegroup{}}$.
\endcor
\definition\label{deftangentspaceliealgebraofliegroup}
The Lie-algebra $\opair{\Tanspace{\IG{}}{\Lieman{\Liegroup{}}}}{\tanliebracket{\Liegroup{}}}$
is referred to as the $\quotl$canonical Lie-algebra of the smooth group $\Liegroup{}$ via the tangent-space at identity$\quotr$,
and is denoted by $\LiegroupLiealgebratan{\Liegroup{}}$.
\endef
\corollary
The dimension of the canonical Lie-algebra of the smooth group $\Liegroup{}$
(and also the dimension of canonical Lie-algebra of the smooth group $\Liegroup{}$ via the tangent-space at identity)
is the same as the dimension of the underlying manifold of the smooth group $\Liegroup{}$. That is,
\begin{equation}
\liedim{\LiegroupLiealgebra{\Liegroup{}}}=
\liedim{\LiegroupLiealgebratan{\Liegroup{}}}=\mandim{\Lieman{\Liegroup{}}}=n.
\end{equation}
\refthm{thmliederivativeisinjective},
\refthm{thmliegroupequiv0}
\endcor
\section{Morphisms of Smooth Groups}
\definition\label{defsmoothgroupmorphism}
$\LieMor{\Liegroup{}}{\Liegroup{1}}$ is defined to be the set of all mappings $\aliemor{}$ from $\G{}$ to $\G{1}$
such that $\aliemor{}$ is simultaneously a group-homomorphism from the underlying group structure of the smooth group $\Liegroup{}$
(that is $\LieG{\Liegroup{}}$)
to the underlying group structure of the smooth group $\Liegroup{1}$ (that is $\LieG{\Liegroup{1}}$),
and a smooth map from the underlying manifold
of the smooth group $\Liegroup{}$ (that is $\Lieman{\Liegroup{}}$) to the underlying manifold of the smooth group $\Liegroup{1}$
(that is $\Lieman{\Liegroup{1}}$). That is,
\begin{equation}
\LieMor{\Liegroup{}}{\Liegroup{1}}:=\GHom{\LieG{\Liegroup{}}}{\LieG{\Liegroup{1}}}\cap
\mapdifclass{\infty}{\Lieman{\Liegroup{}}}{\Lieman{\Liegroup{1}}}.
\end{equation}
Each element of $\LieMor{\Liegroup{}}{\Liegroup{1}}$ is referred to as a
$\quotl$smooth-group-morphism from the smooth group $\Liegroup{}$ to the smooth group $\Liegroup{1}$$\quotr$,
or a $\quotl$Lie-morphism from $\Liegroup{}$ to $\Liegroup{1}$$\quotr$.
\endef
\theorem\label{thmcompositionofsmoothgroupmorphisms}
The composition of every smooth-group-morphism from $\Liegroup{}$ to $\Liegroup{1}$
with each smooth-group-morphism from $\Liegroup{1}$ to $\Liegroup{2}$ is a
smooth-group-morphism from $\Liegroup{}$ to $\Liegroup{2}$. That is,
\begin{equation}
\Foreach{\aliemor{}}{\LieMor{\Liegroup{}}{\Liegroup{1}}}
\Foreach{\aliemor{1}}{\LieMor{\Liegroup{1}}{\Liegroup{2}}}
\cmp{\aliemor{1}}{\aliemor{}\in\LieMor{\Liegroup{}}{\Liegroup{2}}}.
\end{equation}
\proof
According to \refdef{defsmoothgroupmorphism},
it is obvious by considering that the composition of a group-homomorphism from $\LieG{\Liegroup{}}$ to $\LieG{\Liegroup{1}}$
with a group-homomorphism from $\LieG{\Liegroup{1}}$ to $\LieG{\Liegroup{2}}$ is a group-homomorphism from
$\LieG{\Liegroup{}}$ to $\LieG{\Liegroup{2}}$, and
the composition of a smooth map from $\Lieman{\Liegroup{}}$ to $\Lieman{\Liegroup{1}}$
with a smooth map from $\Lieman{\Liegroup{1}}$ to $\Lieman{\Liegroup{2}}$ is a smooth map from
$\Lieman{\Liegroup{}}$ to $\Lieman{\Liegroup{2}}$.
\endthm
\theorem\label{thmtheidentitymapisasmoothgroupmorphism}
The identity map of $\G{}$ is a smooth-group-morphism from $\Liegroup{}$ to itself. That is,
\begin{equation}
\identity{\G{}}\in\LieMor{\Liegroup{}}{\Liegroup{}}.
\end{equation}
\proof
Clearly the identity map of $\G{}$ is both a group-homomorphism from $\LieG{\Liegroup{}}$ to itself,
and a smooth map from $\Lieman{\Liegroup{}}$ to itself. Thus according to \refdef{defsmoothgroupmorphism},
it is a smooth-group-morphism from $\Liegroup{}$ to itself.
\endthm
\definition\label{defsmoothgroupisomorphism}
$\LieIsom{\Liegroup{}}{\Liegroup{1}}$ is defined to be the set of all bijective mappings $\aliemor{}$ from $\G{}$ to $\G{1}$
such that $\aliemor{}$ is a smooth-group-morphism from $\Liegroup{}$ to $\Liegroup{1}$ and $\finv{\aliemor{}}$ is
a smooth-group-morphism from $\Liegroup{1}$ to $\Liegroup{}$. That is,
\begin{equation}
\LieIsom{\Liegroup{}}{\Liegroup{1}}:=\defset{\aliemor{}}{\IF{\G{}}{\G{1}}}{\aliemor{}\in\LieMor{\Liegroup{}}{\Liegroup{1}},~
\finv{\aliemor{}}\in\LieMor{\Liegroup{1}}{\Liegroup{}}}.
\end{equation}
Each element of $\LieIsom{\Liegroup{}}{\Liegroup{1}}$ is referred to as a
$\quotl$smooth-group-isomorphism from the smooth group $\Liegroup{}$ to the smooth group $\Liegroup{1}$$\quotr$,
or a $\quotl$Lie-isomorphism from $\Liegroup{}$ to $\Liegroup{1}$$\quotr$.\\
By definition,
\begin{equation}
\lieisomorphic{\Liegroup{}}{\Liegroup{1}}
:\thenn
\LieIsom{\Liegroup{}}{\Liegroup{1}}\neq\empty.
\end{equation}
It is said that $\quotl$$\Liegroup{}$ is Lie-isomorphic to $\Liegroup{1}$$\quotr$ iff
$\lieisomorphic{\Liegroup{}}{\Liegroup{1}}$.
\endef
\corollary\label{corsmoothgroupisomorphismequiv0}
$\LieIsom{\Liegroup{}}{\Liegroup{1}}$ equals the set of all mappings $\aliemor{}$ from $\G{}$ to $\G{1}$
such that $\aliemor{}$ is simultaneously a group-isomorphism from the underlying group structure of the smooth group $\Liegroup{}$
(that is $\LieG{\Liegroup{}}$)
to the underlying group structure of the smooth group $\Liegroup{1}$ (that is $\LieG{\Liegroup{1}}$),
and an $\infty$-diffeomorphism from the underlying manifold
of the smooth group $\Liegroup{}$ (that is $\Lieman{\Liegroup{}}$) to the underlying manifold of the smooth group $\Liegroup{1}$
(that is $\Lieman{\Liegroup{1}}$). That is,
\begin{equation}
\LieIsom{\Liegroup{}}{\Liegroup{1}}=\GIsom{\LieG{\Liegroup{}}}{\LieG{\Liegroup{1}}}\cap
\Diffeo{\infty}{\Lieman{\Liegroup{}}}{\Lieman{\Liegroup{1}}}.
\end{equation}
\endcor
\theorem\label{thmcompositionofsmoothgroupisomorphisms}
The composition of every smooth-group-isomorphism from $\Liegroup{}$ to $\Liegroup{1}$
with any smooth-group-isomorphism from $\Liegroup{1}$ to $\Liegroup{2}$ is a
smooth-group-isomorphism from $\Liegroup{}$ to $\Liegroup{2}$. That is,
\begin{equation}
\Foreach{\aliemor{}}{\LieIsom{\Liegroup{}}{\Liegroup{1}}}
\Foreach{\aliemor{1}}{\LieIsom{\Liegroup{1}}{\Liegroup{2}}}
\cmp{\aliemor{1}}{\aliemor{}\in\LieIsom{\Liegroup{}}{\Liegroup{2}}}.
\end{equation}
\proof
According to \refcor{corsmoothgroupisomorphismequiv0},
and considering the fact that the composition of a group-isomorphism from the group $\LieG{\Liegroup{}}$ to $\LieG{\Liegroup{1}}$
with a group-isomorphism from $\LieG{\Liegroup{1}}$ to $\LieG{\Liegroup{2}}$ is a group-isomorphism from
$\LieG{\Liegroup{}}$ to $\LieG{\Liegroup{2}}$, and
the composition of a $\infty$-diffeomorphism from $\Lieman{\Liegroup{}}$ to $\Lieman{\Liegroup{1}}$
with an $\infty$-diffeomorphism from $\Lieman{\Liegroup{1}}$ to $\Lieman{\Liegroup{2}}$ is an $\infty$-diffeomorphism from
the manifold $\Lieman{\Liegroup{}}$ to $\Lieman{\Liegroup{2}}$, the truth of the regarded assertion is clear.
\endthm
\theorem\label{thmlieisomorphismisliemorphism}
Every Lie-isomorphism from $\Liegroup{}$ to $\Liegroup{1}$ is a bijective Lie-morphism from $\Liegroup{}$
to $\Liegroup{1}$. That is,
\begin{equation}
\LieIsom{\Liegroup{}}{\Liegroup{1}}\subseteq\LieMor{\Liegroup{}}{\Liegroup{1}}\cap
\IF{\G{}}{\G{1}}.
\end{equation}
\proof
It is clear according to \refdef{defsmoothgroupmorphism} and \refcor{corsmoothgroupisomorphismequiv0},
and considering that every group-isomorphism between a pair of groups is a  bijective group-homomorphism between them,
and every $\infty$-diffeomorphism between a pair of manifolds is a smooth map between them.
\endthm
\theorem\label{thmsmoothgroupisomorphismequiv0}
Every $\infty$-diffeomorphism from the manifold $\Lieman{\Liegroup{}}$ to the manifold $\Lieman{\Liegroup{1}}$ that is simulaneously a
group-homomorphism from $\LieG{\Liegroup{}}$ to $\LieG{\Liegroup{1}}$ is a Lie-isomorphism from $\Liegroup{}$
to $\Liegroup{1}$. That is,
\begin{equation}
\LieIsom{\Liegroup{}}{\Liegroup{1}}:=\GHom{\LieG{\Liegroup{}}}{\LieG{\Liegroup{1}}}\cap
\Diffeo{\infty}{\Lieman{\Liegroup{}}}{\Lieman{\Liegroup{1}}}.
\end{equation}
\proof
It is clear since every bijective group-homomorphism between a pair of groups is necessarily a group-isomorphism
between them.
\endthm
\theorem\label{inverseofalieisomorphismisalieisomorphism}
The inverse mapping of every smooth-group-isomorphism from $\Liegroup{}$ to $\Liegroup{1}$ is
a smooth-group-isomorphism from $\Liegroup{1}$ to $\Liegroup{}$. That is,
\begin{equation}
\Foreach{\aliemor{}}{\LieMor{\Liegroup{}}{\Liegroup{1}}}
\finv{\aliemor{}}\in\LieMor{\Liegroup{1}}{\Liegroup{}}.
\end{equation}
\proof
It is clear according to \refcor{corsmoothgroupisomorphismequiv0}, and considering that the inverse of an $\infty$-diffeomorphism
between a pair of manifolds is an $\infty$-diffeomorphism between them in reverse, and the inverse of a group-isomorphism
between a pair of groups is necessarily a group-isomorphism in reverse.
\endthm
\corollary\label{corgroupoflieautomorphisms}
The set of all smooth-group-isomorphisms from $\Liegroup{}$ to $\Liegroup{}$ endowed with the binary operation of
function composition on it, that is $\opair{\LieIsom{\Liegroup{}}{\Liegroup{}}}{\cmp{}{}}$, possesses the structure of group,
having $\identity{\G{}}$ as the identity element the group-structure. It is a subgroup of the group of all
$\infty$-automorphisms of $\Lieman{\Liegroup{}}$, that is the group $\GDiff{\infty}{\Lieman{\Liegroup{}}}$.
\endcor
\definition\label{deflieautomorphism}
The set of all smooth-group-isomorphisms from $\Liegroup{}$ to $\Liegroup{}$ is denoted by $\LieAut{\Liegroup{}}$.
In addition, the group-structure obtained by equipping $\LieAut{\Liegroup{}}$ with the binary operation of function composition on
it is denoted by $\GLieAut{\Liegroup{}}$.
\begin{align}
\LieAut{\Liegroup{}}:=&\LieIsom{\Liegroup{}}{\Liegroup{}},\cr
\GLieAut{\Liegroup{}}:=&\opair{\LieIsom{\Liegroup{}}{\Liegroup{}}}{\cmp{}{}}.
\end{align}
Each element of $\LieAut{\Liegroup{}}$ is called an $\quotl$Lie-automorphism of the smooth group $\Liegroup{}$, and
$\GLieAut{\Liegroup{}}$ is referred to as the $\quotl$group of Lie-automorphisms of the smooth group $\Liegroup{}$$\quotr$
or the $\quotl$Lie-automorphism group of $\Liegroup{}$$\quotr$.
\endef
\theorem\label{thmconjugationisautomorphism}
For every point $\g{}$ of the smooth group $\Liegroup{}$, the $\g{}$-conjugation of $\LieG{\Liegroup{}}$ is
a Lie-automorphism of $\Liegroup{}$.
\begin{equation}
\Foreach{\g{}}{\G{}}
\gconj{\g{}}{\LieG{\Liegroup{}}}\in\LieAut{\Liegroup{}}.
\end{equation}
\proof
For every $\g{}$ in $\G{}$, it is known that $\gconj{\g{}}{\LieG{{\Liegroup{}}}}$ is a group-automorphism of the manifold
$\LieG{\Liegroup{}}$ and additionally according to \refthm{thmtranslationsarediffeomorphisms} it is an $\infty$-automorphism of
$\Lieman{\Liegroup{}}$. So according to \refdef{thmlieisomorphismisliemorphism}, the result is obtained.
\endthm
\theorem\label{thmliemorphismlocality}
$\aliemor{}$ is taken as an element of $\GHom{\LieG{\Liegroup{}}}{\LieG{\Liegroup{1}}}$ (a group-homomorphism from
$\LieG{\Liegroup{}}$ to $\LieG{\Liegroup{1}}$). $\aliemor{}$ is a smooth-group-morphism from $\Liegroup{}$ to
$\Liegroup{1}$ if and only if $\aliemor{}$ is infinitely differentiable at the identity element of the goup
$\LieG{\Liegroup{}}$, that is $\IG{}$.
\begin{align}
&\hskip\baselineskip\aliemor{}\in\LieMor{\Liegroup{}}{\Liegroup{1}}\cr
&\thenn
\(\begin{aligned}
&\Exists{\phi}{\defset{\p{\phi}}{\maxatlas{}}{\IG{}\in\domain{\p{\phi}}}}\cr
&\Exists{\psi}{\defset{\p{\psi}}{\maxatlas{1}}{\IG{1}=\func{\aliemor{}}{\IG{}}\in\domain{\p{\psi}}}}\cr
&\func{\image{\aliemor{}}}{\domain{\phi}}\subseteq\domain{\psi},~
\cmp{\psi}{\cmp{\aliemor{}}{\finv{\phi}}}\in\banachmapdifclass{\infty}{\R^n}{\R^{n_1}}{\funcimage{\phi}}{\funcimage{\psi}}
\end{aligned}\).
\end{align}
\proof
It is clear that if $\aliemor{}$ is a smooth-group-morphism from $\Liegroup{}$ to $\Liegroup{1}$, then it is
infinitely differentiable at every point of $\Liegroup{}$ and hence at $\IG{}$ in particular.\\
Since $\aliemor{}$ is a group-homomorphism from $\LieG{\Liegroup{}}$ to $\LieG{\Liegroup{1}}$, clearly,
\begin{align}
&\Foreach{\opair{\g{1}}{\g{2}}}{\Cprod{\G{}}{\G{}}}
\func{\aliemor{}}{\g{1}\gop{}\g{2}}=\func{\aliemor{}}{\g{1}}\gop{1}\func{\aliemor{}}{\g{2}},
\label{thmliemorphismlocalitypeq1}\\
&\Foreach{\g{}}{\G{}}
\invg{\[\func{\aliemor{}}{\g{}}\]}{\LieG{\Liegroup{1}}}=\func{\aliemor{}}{\invg{\g{}}{\LieG{\Liegroup{}}}},
\label{thmliemorphismlocalitypeq2}
\end{align}
and accordingly by considering that $\finv{\gltrans{\LieG{\Liegroup{1}}}{\g{1}}}=\gltrans{\LieG{\Liegroup{1}}}{\invg{\g{1}}{}}$
for every $\g{1}$ in $\G{1}$,
\begin{equation}\label{thmliemorphismlocalitypeq3}
\finv{\[\gltrans{\LieG{\Liegroup{1}}}{\func{\aliemor{}}{\g{}}}\]}=
\gltrans{\LieG{\Liegroup{1}}}{\func{\aliemor{}}{\invg{\g{}}{}}}
\end{equation}
\begin{itemize}
\item[$\pr{1}$]
It is assumed that $\aliemor{}$ is infinitely differentiable at $\IG{}$, that is,
\begin{align}\label{thmliemorphismlocalityp1eq1}
&\Existsis{\phi}{\defset{\p{\phi}}{\maxatlas{}}{\IG{}\in\domain{\p{\phi}}}}\cr
&\Existsis{\psi}{\defset{\p{\psi}}{\maxatlas{1}}{\IG{1}=\func{\aliemor{}}{\IG{}}\in\domain{\p{\psi}}}}\cr
&\func{\image{\aliemor{}}}{\domain{\phi}}\subseteq\domain{\psi},~
\cmp{\psi}{\cmp{\aliemor{}}{\finv{\phi}}}\in\banachmapdifclass{\infty}{\R^n}{\R^{n_1}}{\funcimage{\phi}}{\funcimage{\psi}}
\end{align}
\begin{itemize}
\item[$\pr{1-1}$]
$\g{}$ is taken as an arbitrary element of $\G{}$.
According to \refthm{thmtranslationsarediffeomorphisms}, the left-translation of the underlying group
of $\Liegroup{}$ by $\invg{\g{}}{}$ and the left-translation of the underlying group of
$\Liegroup{1}$ by $\func{\aliemor{}}{\invg{\g{}}{}}$ are $\infty$-automorphisms of the
underlying manifolds of $\Liegroup{}$ and $\Liegroup{1}$, respectively. That is,
\begin{align}
\gltrans{\LieG{\Liegroup{}}}{\invg{\g{}}{}}&\in\Diff{\infty}{\Lieman{\Liegroup{}}},
\label{thmliemorphismlocalityp11eq1}\\
\gltrans{\LieG{\Liegroup{1}}}{\func{\aliemor{}}{\invg{\g{}}{}}}&\in\Diff{\infty}{\Lieman{\Liegroup{1}}}.
\label{thmliemorphismlocalityp11eq2}
\end{align}
Thus since $\phi$ is a chart of $\Lieman{\Liegroup{}}$ around $\IG{}$, $\gltrans{\LieG{\Liegroup{}}}{\invg{\g{}}{}}$
transfers $\phi$ to a chart of $\Lieman{\Liegroup{}}$ around
$\func{\finv{\[\gltrans{\LieG{\Liegroup{}}}{\invg{\g{}}{}}\]}}{\IG{}}=\g{}$, and similarly since
$\psi$ is a chart of $\Lieman{\Liegroup{1}}$ around $\IG{1}$, $\gltrans{\LieG{\Liegroup{1}}}{\func{\aliemor{}}{\invg{\g{}}{}}}$
transfers $\psi$ to a chart of $\Lieman{\Liegroup{}}$ around
$\func{\finv{\[\gltrans{\LieG{\Liegroup{1}}}{\func{\aliemor{}}{\invg{\g{}}{}}}\]}}{\IG{1}}=\func{\aliemor{}}{\g{}}$
(\Ref{thmliemorphismlocalitypeq2} implies this equality). That is,
\begin{align}
\bar{\phi}:=&\cmp{\phi}{\gltrans{\LieG{\Liegroup{}}}{\invg{\g{}}{}}}\in
\defset{\p{\phi}}{\maxatlas{}}{\g{}\in\domain{\p{\phi}}},
\label{thmliemorphismlocalityp11eq3}\\
\bar{\psi}:=&\cmp{\psi}{\gltrans{\LieG{\Liegroup{1}}}{\func{\aliemor{}}{\invg{\g{}}{}}}}\in
\defset{\p{\psi}}{\maxatlas{1}}{\func{\aliemor{}}{\g{}}\in\domain{\p{\psi}}}.
\label{thmliemorphismlocalityp11eq4}
\end{align}
It is clear that,
\begin{align}
\domain{\bar{\phi}}&=\func{\image{\[\gltrans{\LieG{\Liegroup{}}}{\g{}}\]}}{\domain{\phi}},
\label{thmliemorphismlocalityp11eq5}\\
\domain{\bar{\psi}}&=\func{\image{\[\gltrans{\LieG{\Liegroup{1}}}{\func{\aliemor{}}{\g{}}}\]}}{\domain{\psi}}.
\label{thmliemorphismlocalityp11eq6}
\end{align}
\Ref{thmliemorphismlocalitypeq1}, \Ref{thmliemorphismlocalityp1eq1}, \Ref{thmliemorphismlocalityp11eq5}, and
\Ref{thmliemorphismlocalityp11eq6} imply that,
\begin{align}\label{thmliemorphismlocalityp11eq7}
\func{\image{\aliemor{}}}{\domain{\bar{\phi}}}&=
\func{\image{\aliemor{}}}{\func{\image{\[\gltrans{\LieG{\Liegroup{}}}{\g{}}\]}}{\domain{\phi}}}\cr
&=\func{\image{\[\cmp{\aliemor{}}{\gltrans{\LieG{\Liegroup{}}}{\g{}}}\]}}{\domain{\phi}}\cr
&=\func{\image{\[\cmp{\gltrans{\LieG{\Liegroup{}}}{\func{\aliemor{}}{\g{}}}}{\aliemor{}}\]}}{\domain{\phi}}\cr
&=\func{\image{\[\gltrans{\LieG{\Liegroup{}}}{\func{\aliemor{}}{\g{}}}\]}}{\func{\image{\aliemor{}}}{\domain{\phi}}}\cr
&\subseteq\func{\image{\[\gltrans{\LieG{\Liegroup{}}}{\func{\aliemor{}}{\g{}}}\]}}{\domain{\psi}}\cr
&=\domain{\bar{\psi}}.
\end{align}
According to \Ref{thmliemorphismlocalityp11eq3},
\begin{align}\label{thmliemorphismlocalityp11eq8}
\finv{\bar{\phi}}&=\finv{\[\cmp{\phi}{\gltrans{\LieG{\Liegroup{}}}{\invg{\g{}}{}}}\]}\cr
&=\cmp{\finv{\[\gltrans{\LieG{\Liegroup{}}}{\invg{\g{}}{}}\]}}{\finv{\phi}}\cr
&=\cmp{\gltrans{\LieG{\Liegroup{}}}{\g{}}}{\finv{\phi}}.
\end{align}
\Ref{thmliemorphismlocalityp11eq4} and \Ref{thmliemorphismlocalityp11eq8} yield,
\begin{equation}\label{thmliemorphismlocalityp11eq9}
\cmp{\bar{\psi}}{\cmp{\aliemor{}}{\finv{\bar{\phi}}}}=
\cmp{\psi}{\cmp{\(\cmp{\cmp{\gltrans{\LieG{\Liegroup{1}}}{\func{\aliemor{}}{\invg{\g{}}{}}}}{\aliemor{}}}
{\gltrans{\LieG{\Liegroup{}}}{\g{}}}\)}{\finv{\phi}}}.
\end{equation}
According to \Ref{thmliemorphismlocalitypeq1}, it can be easily verified that,
\begin{equation}\label{thmliemorphismlocalityp11eq10}
\cmp{\cmp{\gltrans{\LieG{\Liegroup{1}}}{\func{\aliemor{}}{\invg{\g{}}{}}}}{\aliemor{}}}
{\gltrans{\LieG{\Liegroup{}}}{\g{}}}=\aliemor{},
\end{equation}
and hence \Ref{thmliemorphismlocalityp11eq9} becomes,
\begin{equation}\label{thmliemorphismlocalityp11eq11}
\cmp{\bar{\psi}}{\cmp{\aliemor{}}{\finv{\bar{\phi}}}}=
\cmp{\psi}{\cmp{\aliemor{}}{\finv{\phi}}},
\end{equation}
and thus according to \Ref{thmliemorphismlocalityp1eq1} and \Ref{thmliemorphismlocalityp11eq8},
\begin{equation}\label{thmliemorphismlocalityp11eq12}
\cmp{\bar{\psi}}{\cmp{\aliemor{}}{\finv{\bar{\phi}}}}\in
\banachmapdifclass{\infty}{\R^n}{\R^{n_1}}{\funcimage{\phi}}{\funcimage{\psi}}.
\end{equation}
According to \Ref{thmliemorphismlocalityp11eq3} and \Ref{thmliemorphismlocalityp11eq4},
it is obvious that,
\begin{align}
\funcimage{\bar{\phi}}&=\funcimage{\phi},
\label{thmliemorphismlocalityp11eq13}\\
\funcimage{\bar{\psi}}&=\funcimage{\psi},
\label{thmliemorphismlocalityp11eq14}
\end{align}
and hence \Ref{thmliemorphismlocalityp11eq12} becomes,
\begin{equation}\label{thmliemorphismlocalityp11eq15}
\cmp{\bar{\psi}}{\cmp{\aliemor{}}{\finv{\bar{\phi}}}}\in
\banachmapdifclass{\infty}{\R^n}{\R^{n_1}}{\funcimage{\bar{\phi}}}{\funcimage{\bar{\psi}}}.
\end{equation}
\endp
Therefore,
\begin{align}\label{thmliemorphismlocalityp1eq2}
\Foreach{\g{}}{\G{}}&\cr
&\Exists{\bar{\phi}}{\defset{\p{\phi}}{\maxatlas{}}{\g{}\in\domain{\p{\phi}}}}
\Exists{\bar{\psi}}{\defset{\p{\psi}}{\maxatlas{1}}{\func{\aliemor{}}{\g{}}\in\domain{\p{\psi}}}}\cr
&\func{\image{\aliemor{}}}{\domain{\bar{\phi}}}\subseteq\domain{\bar{\psi}},~
\cmp{\bar{\psi}}{\cmp{\aliemor{}}{\finv{\bar{\phi}}}}\in
\banachmapdifclass{\infty}{\R^n}{\R^{n_1}}{\funcimage{\bar{\phi}}}{\funcimage{\bar{\psi}}},\cr
&{}
\end{align}
which means $\aliemor{}$ is a smooth map from the underlying manifold of $\Liegroup{}$ to
the underlying manifold of $\Liegroup{1}$, that is,
\begin{equation}
\aliemor{}\in\mapdifclass{\infty}{\Lieman{\Liegroup{}}}{\Lieman{\Liegroup{1}}},
\end{equation}
and hence according to \refdef{defsmoothgroupmorphism} and considering that $\aliemor{}$ is already assumed as an element of
$\GHom{\LieG{\Liegroup{}}}{\LieG{\Liegroup{1}}}$, it becomes evident that,
\begin{equation}
\aliemor{}\in\LieMor{\Liegroup{}}{\Liegroup{1}}.
\end{equation}
\end{itemize}
\endp
\end{itemize}
\endthm
\definition\label{definducedliealgebramorphismfromliemorphism}
The mapping $\indliemor{\Liegroup{}}{\Liegroup{1}}$ is defined as,
\begin{align}
&\indliemor{\Liegroup{}}{\Liegroup{1}}\indef\Func{\LieMor{\Liegroup{}}{\Liegroup{1}}}
{\Func{\Leftinvvf{\Liegroup{}}}{\Leftinvvf{\Liegroup{1}}}},\cr
&\Foreach{\aliemor{}}{\LieMor{\Liegroup{}}{\Liegroup{1}}}
\func{\indliemor{\Liegroup{}}{\Liegroup{1}}}{\aliemor{}}\eqdef
\cmp{\finv{\liegvftan{\Liegroup{1}}}}{\cmp{\[\der{\aliemor{}}{\Lieman{\Liegroup{}}}{\Lieman{\Liegroup{1}}}\]}{\liegvftan{\Liegroup{}}}},\cr
&{}
\end{align}
or equivalently,
\begin{equation}
\Foreach{\aliemor{}}{\LieMor{\Liegroup{}}{\Liegroup{1}}}
\Foreach{\avecf{}}{\Leftinvvf{\Liegroup{}}}
\func{\[\func{\indliemor{\Liegroup{}}{\Liegroup{1}}}{\aliemor{}}\]}{\avecf{}}\eqdef
\func{\finv{\liegvftan{\Liegroup{1}}}}{\func{\[\der{\aliemor{}}{\Lieman{\Liegroup{}}}{\Lieman{\Liegroup{1}}}\]}{\func{\avecf{}}{\IG{}}}}.
\end{equation}
In other words, given a smooth-group-morphism $\aliemor{}$ from $\Liegroup{}$ to $\Liegroup{1}$
and a left-invariant vector-field $\avecf{}$ on $\Liegroup{}$,
$\func{\[\func{\indliemor{\Liegroup{}}{\Liegroup{1}}}{\aliemor{}}\]}{\avecf{}}$ is defined to be the unique left-invariant
vector-field $\avecff{}$ on $\Liegroup{1}$ such that,
\begin{equation}
\func{\avecff{}}{\IG{1}}\eqdef
\func{\[\der{\aliemor{}}{\Lieman{\Liegroup{}}}{\Lieman{\Liegroup{1}}}\]}{\func{\avecf{}}{\IG{}}}.
\end{equation}
\endef
\theorem\label{thmdefinducedliealgebramorphismfromliemorphismislinear}
For every smooth-group-morphism $\aliemor{}$ from $\Liegroup{}$ to $\Liegroup{1}$,
$\func{\indliemor{\Liegroup{}}{\Liegroup{1}}}{\aliemor{}}$ is a linear map from
the vector-space of left-invariant vector-fields on $\Liegroup{}$ to the vector-space of
left-invariant vector-fields on $\Liegroup{1}$. That is,
\begin{equation}
\Foreach{\aliemor{}}{\LieMor{\Liegroup{}}{\Liegroup{1}}}
\func{\indliemor{\Liegroup{}}{\Liegroup{1}}}{\aliemor{}}\in
\Lin{\VLeftinvvf{\Liegroup{}}}{\VLeftinvvf{\Liegroup{1}}},
\end{equation}
which means,
\begin{align}
&\Foreach{\aliemor{}}{\LieMor{\Liegroup{}}{\Liegroup{1}}}\cr
&\Foreach{\opair{\avecf{1}}{\avecf{2}}}{\Cprod{\Leftinvvf{\Liegroup{}}}{\Leftinvvf{\Liegroup{}}}}
\Foreach{\c}{\R}\cr
&\func{\[\func{\indliemor{\Liegroup{}}{\Liegroup{1}}}{\aliemor{}}\]}{\c\avecf{1}+\avecf{2}}=
\c\func{\[\func{\indliemor{\Liegroup{}}{\Liegroup{1}}}{\aliemor{}}\]}{\avecf{1}}+
\func{\[\func{\indliemor{\Liegroup{}}{\Liegroup{1}}}{\aliemor{}}\]}{\avecf{2}}.
\end{align}
\proof
Each $\avecf{1}$ and $\avecf{2}$ is taken as an arbitrary element of $\Leftinvvf{\Liegroup{}}$, and $\c$
as an arbitrary real number. According to \refcor{corlinvvftanspalinearcecorrespondence} and
\refdef{definducedliealgebramorphismfromliemorphism}, and by invoking the canonical linear-structure of
$\Leftinvvf{\Liegroup{}}$, and considering that $\der{\aliemor{}}{\Lieman{\Liegroup{}}}{\Lieman{\Liegroup{1}}}$ operates
linearly when restricted to the tangent-space of $\Lieman{\Liegroup{}}$ at $\IG{}$,
\begin{align}
\func{\[\func{\indliemor{\Liegroup{}}{\Liegroup{1}}}{\aliemor{}}\]}{\c\avecf{1}+\avecf{2}}&=
\func{\finv{\liegvftan{\Liegroup{1}}}}{\func{\[\der{\aliemor{}}{\Lieman{\Liegroup{}}}{\Lieman{\Liegroup{1}}}\]}{\func{\[\c\avecf{1}+\avecf{2}\]}{\IG{}}}}\cr
&=\func{\finv{\liegvftan{\Liegroup{1}}}}{\func{\[\der{\aliemor{}}{\Lieman{\Liegroup{}}}{\Lieman{\Liegroup{1}}}\]}{\c\func{\avecf{1}}{\IG{}}+\func{\avecf{2}}{\IG{}}}}\cr
&=\func{\finv{\liegvftan{\Liegroup{1}}}}{\c\func{\[\der{\aliemor{}}{\Lieman{\Liegroup{}}}{\Lieman{\Liegroup{1}}}\]}{\func{\avecf{1}}{\IG{}}}+
\func{\[\der{\aliemor{}}{\Lieman{\Liegroup{}}}{\Lieman{\Liegroup{1}}}\]}{\func{\avecf{2}}{\IG{}}}}\cr
&=\c\[\func{\finv{\liegvftan{\Liegroup{1}}}}{\func{\[\der{\aliemor{}}{\Lieman{\Liegroup{}}}{\Lieman{\Liegroup{1}}}\]}{\func{\avecf{1}}{\IG{}}}}\]+
\func{\finv{\liegvftan{\Liegroup{1}}}}{\func{\[\der{\aliemor{}}{\Lieman{\Liegroup{}}}{\Lieman{\Liegroup{1}}}\]}{\func{\avecf{2}}{\IG{}}}}\cr
&=\c\func{\[\func{\indliemor{\Liegroup{}}{\Liegroup{1}}}{\aliemor{}}\]}{\avecf{1}}+
\func{\[\func{\indliemor{\Liegroup{}}{\Liegroup{1}}}{\aliemor{}}\]}{\avecf{2}}.
\end{align}
\endthm
\lemma\label{leminducedliealgebramorphismfromliemorphism0}
For every smooth-group-morphism from $\Liegroup{}$ to $\Liegroup{1}$, and every
left-invariant vector-field $\avecf{}$ on $\Liegroup{}$, the following diagram commutes.
\begin{center}
\vskip0.5\baselineskip
\begin{tikzcd}[row sep=6em, column sep=6em]
& \G{}
\arrow{r}{\aliemor{}}
\arrow[swap]{d}{\avecf{}}
& \G{1}
\arrow{d}{\func{\[\func{\indliemor{\Liegroup{}}{\Liegroup{1}}}{\aliemor{}}\]}{\avecf{}}} \\
& \tanbun{\Lieman{\Liegroup{}}}
\arrow[swap]{r}{\der{\aliemor{}}{\Lieman{\Liegroup{}}}{\Lieman{\Liegroup{1}}}} & \tanbun{\Lieman{\Liegroup{1}}}
\end{tikzcd}
\end{center}
That is,
\begin{align}
&\Foreach{\aliemor{}}{\LieMor{\Liegroup{}}{\Liegroup{1}}}
\Foreach{\avecf{}}{\Leftinvvf{\Liegroup{}}}\cr
&\cmp{\bigg(\func{\[\func{\indliemor{\Liegroup{}}{\Liegroup{1}}}{\aliemor{}}\]}{\avecf{}}\bigg)}{\aliemor{}}=
\cmp{\[\der{\aliemor{}}{\Lieman{\Liegroup{}}}{\Lieman{\Liegroup{1}}}\]}{\avecf{}}.
\end{align}
\proof
$\aliemor{}$ is taken as an arbitrary element of $\LieMor{\Liegroup{}}{\Liegroup{1}}$,
$\avecf{}$ as an arbitrary element of $\Leftinvvf{\Liegroup{}}$, and $\g{}$ as an arbitrary element of $\G{}$. Also,
\begin{equation}\label{leminducedliealgebramorphismfromliemorphism0peq1}
\avecff{}:=\func{\[\func{\indliemor{\Liegroup{}}{\Liegroup{1}}}{\aliemor{}}\]}{\avecf{}}.
\end{equation}
Since $\func{\[\func{\indliemor{\Liegroup{}}{\Liegroup{1}}}{\aliemor{}}\]}{\avecf{}}$ is a left-invariant vector-field on $\Liegroup{1}$,
according to \reflem{lemleftinvariantvectorfieldsequiv0}, \refdef{definducedliealgebramorphismfromliemorphism},
and the chain rule of differentiation,
\begin{align}\label{leminducedliealgebramorphismfromliemorphism0peq2}
\func{\avecff{}}{\func{\aliemor{}}{\g{}}}&=
\func{\[\der{\gltrans{\LieG{\Liegroup{1}}}{\func{\aliemor{}}{\g{}}}}{\Lieman{\Liegroup{1}}}{\Lieman{\Liegroup{1}}}\]}
{\func{\avecff{}}{\IG{1}}}\cr
&=\func{\[\der{\gltrans{\LieG{\Liegroup{1}}}{\func{\aliemor{}}{\g{}}}}{\Lieman{\Liegroup{1}}}{\Lieman{\Liegroup{1}}}\]}
{\func{\[\der{\aliemor{}}{\Lieman{\Liegroup{}}}{\Lieman{\Liegroup{1}}}\]}{\func{\avecf{}}{\IG{}}}}\cr
&=\func{\(\cmp{\[\der{\gltrans{\LieG{\Liegroup{1}}}{\func{\aliemor{}}{\g{}}}}{\Lieman{\Liegroup{1}}}{\Lieman{\Liegroup{1}}}\]}
{\[\der{\aliemor{}}{\Lieman{\Liegroup{}}}{\Lieman{\Liegroup{1}}}\]}\)}{\func{\avecf{}}{\IG{}}}\cr
&=\func{\[\der{\(\cmp{\gltrans{\LieG{\Liegroup{1}}}{\func{\aliemor{}}{\g{}}}}{\aliemor{}}\)}{\Lieman{\Liegroup{}}}{\Lieman{\Liegroup{1}}}\]}
{\func{\avecf{}}{\IG{}}}.
\end{align}
Furthermore, according to \refdef{defgrouptranslations} and considering that $\aliemor{}$ is a group-homomorphism
from $\LieG{\Liegroup{}}$ to $\LieG{\Liegroup{1}}$,
\begin{align}\label{leminducedliealgebramorphismfromliemorphism0peq3}
\Foreach{\point}{\G{}}
\func{\(\cmp{\gltrans{\LieG{\Liegroup{1}}}{\func{\aliemor{}}{\g{}}}}{\aliemor{}}\)}{\point}&=
\func{\aliemor{}}{\g{}}\gop{1}\func{\aliemor{}}{\point}\cr
&=\func{\aliemor{}}{\g{}\gop{}\point}\cr
&=\func{\(\cmp{\aliemor{}}{\gltrans{\LieG{\Liegroup{}}}{\g{}}}\)}{\point},
\end{align}
and hence,
\begin{equation}\label{leminducedliealgebramorphismfromliemorphism0peq4}
\cmp{\gltrans{\LieG{\Liegroup{1}}}{\func{\aliemor{}}{\g{}}}}{\aliemor{}}=
\cmp{\aliemor{}}{\gltrans{\LieG{\Liegroup{}}}{\g{}}}.
\end{equation}
Since $\avecf{}$ is a left-invariant vector-field on $\Liegroup{}$,
based on \reflem{lemleftinvariantvectorfieldsequiv0} and the chain rule of differentiation,
\Ref{leminducedliealgebramorphismfromliemorphism0peq2} and \Ref{leminducedliealgebramorphismfromliemorphism0peq4} imply,
\begin{align}
\func{\avecff{}}{\func{\aliemor{}}{\g{}}}&=
\func{\[\der{\(\cmp{\aliemor{}}{\gltrans{\LieG{\Liegroup{}}}{\g{}}}\)}{\Lieman{\Liegroup{}}}{\Lieman{\Liegroup{1}}}\]}
{\func{\avecf{}}{\IG{}}}\cr
&=\func{\[\cmp{\der{\aliemor{}}{\Lieman{\Liegroup{}}}{\Lieman{\Liegroup{1}}}}
{\der{\gltrans{\LieG{\Liegroup{}}}{\g{}}}{\Lieman{\Liegroup{}}}{\Lieman{\Liegroup{1}}}}\]}{\func{\avecf{}}{\IG{}}}\cr
&=\func{\[\der{\aliemor{}}{\Lieman{\Liegroup{}}}{\Lieman{\Liegroup{1}}}\]}{
\func{\[\der{\gltrans{\LieG{\Liegroup{}}}{\g{}}}{\Lieman{\Liegroup{}}}{\Lieman{\Liegroup{1}}}\]}{\func{\avecf{}}{\IG{}}}}\cr
&=\func{\[\der{\aliemor{}}{\Lieman{\Liegroup{}}}{\Lieman{\Liegroup{1}}}\]}{\func{\avecf{}}{\g{}}}.
\end{align}
\endlem
\theorem\label{thmdefinducedliealgebramorphismfromliemorphismpreservesliebracket}
$\aliemor{}$ is taken as an element of $\LieMor{\Liegroup{}}{\Liegroup{1}}$.
\begin{align}
&\Foreach{\opair{\avecf{1}}{\avecf{2}}}{\Cprod{\Leftinvvf{\Liegroup{}}}{\Leftinvvf{\Liegroup{}}}}\cr
&\func{\[\func{\indliemor{\Liegroup{}}{\Liegroup{1}}}{\aliemor{}}\]}{\liebracket{\avecf{1}}{\avecf{2}}{\Liegroup{}}}=
\liebracket{\func{\[\func{\indliemor{\Liegroup{}}{\Liegroup{1}}}{\aliemor{}}\]}{\avecf{1}}}
{\func{\[\func{\indliemor{\Liegroup{}}{\Liegroup{1}}}{\aliemor{}}\]}{\avecf{2}}}{\Liegroup{1}}.
\end{align}
\proof
Each $\avecf{1}$ and $\avecf{2}$ is taken as an arbitrary element of $\Leftinvvf{\Liegroup{}}$, and
\begin{align}
\avecff{1}:=&\func{\[\func{\indliemor{\Liegroup{}}{\Liegroup{1}}}{\aliemor{}}\]}{\avecf{1}},\\
\avecff{2}:=&\func{\[\func{\indliemor{\Liegroup{}}{\Liegroup{1}}}{\aliemor{}}\]}{\avecf{2}}.
\end{align}
According to \reflem{leminducedliealgebramorphismfromliemorphism0},
\begin{align}
\cmp{\avecff{1}}{\aliemor{}}&=
\cmp{\[\der{\aliemor{}}{\Lieman{\Liegroup{}}}{\Lieman{\Liegroup{1}}}\]}{\avecf{1}}\\
\cmp{\avecff{2}}{\aliemor{}}&=
\cmp{\[\der{\aliemor{}}{\Lieman{\Liegroup{}}}{\Lieman{\Liegroup{1}}}\]}{\avecf{2}}.
\end{align}
Thus according to \refdef{definducedliealgebramorphismfromliemorphism} and
\refthm{thmliebracketofrelatedvectorfields},
\begin{equation}
\cmp{\liebracket{\avecff{1}}{\avecff{2}}{\Liegroup{1}}}{\aliemor{}}=
\cmp{\[\der{\aliemor{}}{\Lieman{\Liegroup{}}}{\Lieman{\Liegroup{1}}}\]}{\liebracket{\avecf{1}}{\avecf{2}}{\Liegroup{}}},
\end{equation}
and hence, considering that $\func{\aliemor{}}{\IG{}}=\IG{1}$,
\begin{equation}
\func{\liebracket{\avecff{1}}{\avecff{2}}{\Liegroup{1}}}{\IG{1}}=
\func{\[\der{\aliemor{}}{\Lieman{\Liegroup{}}}{\Lieman{\Liegroup{1}}}\]}{\func{\liebracket{\avecf{1}}{\avecf{2}}{\Liegroup{}}}{\IG{}}},
\end{equation}
which according to \refdef{definducedliealgebramorphismfromliemorphism} means,
\begin{equation}
\liebracket{\avecff{1}}{\avecff{2}}{\Liegroup{1}}=
\func{\[\func{\indliemor{\Liegroup{}}{\Liegroup{1}}}{\aliemor{}}\]}{\liebracket{\avecf{1}}{\avecf{2}}{\Liegroup{}}}.
\end{equation}
\endthm
\theorem\label{thminducedliealgebramorphismfromliemorphism1}
For every smooth-group-morphism from the smooth group $\Liegroup{}$ to the smooth group
$\Liegroup{1}$, the mapping $\func{\indliemor{\Liegroup{}}{\Liegroup{1}}}{\aliemor{}}$ is a
Lie-algebra-morphism from the canonical Lie-algebra of $\Liegroup{}$ to the canonical Lie-algebra of
$\Liegroup{1}$. That is,
\begin{equation}
\Foreach{\aliemor{}}{\LieMor{\Liegroup{}}{\Liegroup{1}}}
\func{\indliemor{\Liegroup{}}{\Liegroup{1}}}{\aliemor{}}\in
\LiealgMor{\LiegroupLiealgebra{\Liegroup{}}}{\LiegroupLiealgebra{\Liegroup{1}}}.
\end{equation}
\proof
According to \refdef{defliealgebramorphism}, \refdef{defliealgebraofliegroup},
\refdef{definducedliealgebramorphismfromliemorphism},
\refthm{thmdefinducedliealgebramorphismfromliemorphismislinear},
and \refthm{thmdefinducedliealgebramorphismfromliemorphismpreservesliebracket}, it is clear.
\endthm
\lemma\label{lemthminducedliealgebraisomorphismfromlieisomorphism0}
For every smooth-group-isomorphism from $\Liegroup{}$ to $\Liegroup{1}$, $\func{\indliemor{\Liegroup{}}{\Liegroup{1}}}{\aliemor{}}$ is a
bijective mapping, and
$\func{\indliemor{\Liegroup{}}{\Liegroup{1}}}{\finv{\aliemor{}}}$
is the inverse map of $\func{\indliemor{\Liegroup{}}{\Liegroup{1}}}{\aliemor{}}$. That is,
\begin{align}
\begin{aligned}
\Foreach{\aliemor{}}{\LieIsom{\Liegroup{}}{\Liegroup{1}}}
\end{aligned}
\left\{\begin{aligned}
&\func{\indliemor{\Liegroup{}}{\Liegroup{1}}}{\aliemor{}}\in
\IF{\Leftinvvf{\Liegroup{}}}{\Leftinvvf{\Liegroup{1}}},\\
&\func{\indliemor{\Liegroup{}}{\Liegroup{1}}}{\finv{\aliemor{}}}=
\finv{\[\func{\indliemor{\Liegroup{}}{\Liegroup{1}}}{\aliemor{}}\]}.
\end{aligned}\right.
\end{align}
\proof
$\aliemor{}$ is taken as an element of $\LieIsom{\Liegroup{}}{\Liegroup{1}}$. Thus according to \refdef{defsmoothgroupisomorphism},
\begin{align}
\aliemor{}&\in\LieMor{\Liegroup{}}{\Liegroup{1}},\\
\finv{\aliemor{}}&\in\LieMor{\Liegroup{1}}{\Liegroup{}}.
\end{align}
Also it is known that,
\begin{equation}
\IF{\tanbun{\Lieman{\Liegroup{1}}}}{\tanbun{\Lieman{\Liegroup{}}}}
\ni
\der{\finv{\aliemor{}}}{\Lieman{\Liegroup{1}}}{\Lieman{\Liegroup{}}}=
\finv{\[\der{\aliemor{}}{\Lieman{\Liegroup{}}}{\Lieman{\Liegroup{1}}}\]}.
\end{equation}
Thus according to \refdef{definducedliealgebramorphismfromliemorphism} and \refcor{corlinvvftanspalinearcecorrespondence},
\begin{align}
\func{\indliemor{\Liegroup{1}}{\Liegroup{}}}{\finv{\aliemor{}}}&=
\cmp{\finv{\liegvftan{\Liegroup{}}}}{\cmp{\[\der{\finv{\aliemor{}}}{\Lieman{\Liegroup{1}}}{\Lieman{\Liegroup{}}}\]}{\liegvftan{\Liegroup{1}}}}\cr
&=\cmp{\finv{\liegvftan{\Liegroup{}}}}{\cmp{\finv{\[\der{\aliemor{}}{\Lieman{\Liegroup{}}}{\Lieman{\Liegroup{1}}}\]}}{\liegvftan{\Liegroup{1}}}}\in
\IF{\Leftinvvf{\Liegroup{}}}{\Leftinvvf{\Liegroup{1}}},
\end{align}
and accordingly,
\begin{align}
\func{\indliemor{\Liegroup{1}}{\Liegroup{}}}{\finv{\aliemor{}}}
&=\finv{\bigg(\cmp{\finv{\liegvftan{\Liegroup{1}}}}{\cmp{\[\der{\aliemor{}}
{\Lieman{\Liegroup{}}}{\Lieman{\Liegroup{1}}}\]}{\liegvftan{\Liegroup{}}}}\bigg)}\cr
&=\finv{\[\func{\indliemor{\Liegroup{}}{\Liegroup{1}}}{\aliemor{}}\]}.
\end{align}
\endlem
\theorem\label{thmthminducedliealgebraisomorphismfromlieisomorphism1}
For every smooth-group-isomorphism from the smooth group $\Liegroup{}$ to the smooth group
$\Liegroup{1}$, the mapping $\func{\indliemor{\Liegroup{}}{\Liegroup{1}}}{\aliemor{}}$ is a
Lie-algebra-isomorphism from the canonical Lie-algebra of $\Liegroup{}$ to the canonical Lie-algebra of
$\Liegroup{1}$. That is,
\begin{equation}
\Foreach{\aliemor{}}{\LieIsom{\Liegroup{}}{\Liegroup{1}}}
\func{\indliemor{\Liegroup{}}{\Liegroup{1}}}{\aliemor{}}\in
\LiealgIsom{\LiegroupLiealgebra{\Liegroup{}}}{\LiegroupLiealgebra{\Liegroup{1}}}.
\end{equation}
\proof
$\aliemor{}$ is taken as an element of $\LieIsom{\Liegroup{}}{\Liegroup{1}}$. Thus
according to \refdef{defsmoothgroupisomorphism},
\begin{align}
\aliemor{}&\in\LieMor{\Liegroup{}}{\Liegroup{1}},\\
\finv{\aliemor{}}&\in\LieMor{\Liegroup{1}}{\Liegroup{}},
\end{align}
and hence according to \refthm{thminducedliealgebramorphismfromliemorphism1},
\begin{align}
\func{\indliemor{\Liegroup{}}{\Liegroup{1}}}{\aliemor{}}&\in
\LiealgMor{\LiegroupLiealgebra{\Liegroup{}}}{\LiegroupLiealgebra{\Liegroup{1}}},\\
\func{\indliemor{\Liegroup{1}}{\Liegroup{}}}{\finv{\aliemor{}}}&\in
\LiealgMor{\LiegroupLiealgebra{\Liegroup{1}}}{\LiegroupLiealgebra{\Liegroup{}}}.
\end{align}
Therefore, since $\func{\indliemor{\Liegroup{}}{\Liegroup{1}}}{\aliemor{}}$ is a bijective mapping and
$\func{\indliemor{\Liegroup{}}{\Liegroup{1}}}{\finv{\aliemor{}}}=\finv{\[\func{\indliemor{\Liegroup{}}{\Liegroup{1}}}{\aliemor{}}\]}$
(as asserted in \reflem{lemthminducedliealgebraisomorphismfromlieisomorphism0}), according to \refdef{defliealgebramorphism},
\begin{equation}
\func{\indliemor{\Liegroup{}}{\Liegroup{1}}}{\aliemor{}}\in
\LiealgIsom{\LiegroupLiealgebra{\Liegroup{}}}{\LiegroupLiealgebra{\Liegroup{1}}}.
\end{equation}
\endthm
\corollary\label{corthminducedliealgebraisomorphismfromlieisomorphism2}
If $\Liegroup{}$ is Lie-isomorphic to $\Liegroup{1}$, then the canonical Lie-algebra of $\Liegroup{}$ is Lie-algebraically-isomorphic
to the canonical Lie-algebra of $\Liegroup{1}$. That is,
\begin{equation}
\(\lieisomorphic{\Liegroup{}}{\Liegroup{1}}\)\then
\(\liealgisomorphic{\aliealgebra{}}{\aliealgebra{1}}\).
\end{equation}
\endcor
\section{Smooth Subgroups of a Smooth Group}
\definition\label{defimmersedliesubgroup}
\begin{itemize}
\item[$\centerdot$]
$\Liegroup{1}$ is said to be an $\quotl$($n_1$-dimensional) immersed smooth (Lie) subgroup of the smooth group $\Liegroup{}$$\quotr$ iff
$\Lieman{\Liegroup{1}}$ is an ($n_1$-dimensional) $\infty$-immersed submanifold of $\Lieman{\Liegroup{}}$, and $\LieG{\Liegroup{1}}$
is a subgroup of $\LieG{\Liegroup{}}$. The set of all immersed smooth subgroups of $\Liegroup{}$ is denoted by $\imsubgroup{\Liegroup{}}$.
Moreover, the set of all immersed smooth subgroups of $\Liegroup{}$ whose underlying topological-space is connected is denoted by
$\connectedimsubgroup{\Liegroup{}}$.
\item[$\centerdot$]
$\Liegroup{1}$ is said to be an $\quotl$($n_1$-dimensional) embedded smooth (Lie) subgroup of the smooth group $\Liegroup{}$$\quotr$ iff
$\Lieman{\Liegroup{1}}$ is an ($n_1$-dimensional) $\infty$-embedded submanifold of $\Lieman{\Liegroup{}}$, and $\LieG{\Liegroup{1}}$
is a subgroup of $\LieG{\Liegroup{}}$. The set of all embedded smooth subgroups of $\Liegroup{}$ is denoted by $\emsubgroup{\Liegroup{}}$.
\end{itemize}
\endef
\lemma\label{leminjectionofsmoothsubgroupisamorphism}
If $\Liegroup{1}$ is an immersed smooth subgroup of $\Liegroup{}$, then the injection of $\G{1}$ into $\G{}$
is a smooth-group-morphism from $\Liegroup{1}$ to $\Liegroup{}$.
\begin{equation}
\Liegroup{1}\in\imsubgroup{\Liegroup{}}\then
\Injection{\G{1}}{\G{}}\in\LieMor{\Liegroup{1}}{\Liegroup{}}.
\end{equation}
\proof
It is assumed that $\Liegroup{1}$ is an immersed smooth subgroup of $\Liegroup{}$.
So $\Injection{\G{1}}{\G{}}$ is an $\infty$-immersion and hence a smooth map from $\Lieman{\Liegroup{1}}$ to $\Lieman{\Liegroup{}}$.
Additionally, $\Injection{\G{1}}{\G{}}$ is obviously a group-homomorphism from $\LieG{\Liegroup{1}}$ to $\LieG{\Liegroup{}}$.
Therefore, according to \refdef{defsmoothgroupmorphism}, $\Injection{\G{1}}{\G{}}$ is a smooth-group-morphism from
$\Liegroup{1}$ to $\Liegroup{}$.
\endlem
\theorem\label{thmsmoothsubgroupliealgebra}
If $\Liegroup{1}$ is an immersed smooth subgroup of $\Liegroup{}$, then the canonical Lie-algebra of $\Liegroup{1}$
is isomorphic to a Lie-subalgebra of the canonical Lie-algebra of $\Liegroup{}$.
\begin{equation}
\Liegroup{1}\in\imsubgroup{\Liegroup{}}\then
\[\Exists{\WW{}}{\sublie{\LiegroupLiealgebra{\Liegroup{}}}}
\LiealgIsom{\LiegroupLiealgebra{\Liegroup{1}}}{\subspace{\LiegroupLiealgebra{\Liegroup{}}}{\WW{}}}\neq\empty\].
\end{equation}
Actually, the image of $\func{\indliemor{\Liegroup{1}}{\Liegroup{}}}{\Injection{\G{1}}{\G{}}}$, that is
$\funcimage{\func{\indliemor{\Liegroup{1}}{\Liegroup{}}}{\Injection{\G{1}}{\G{}}}}$ is an instance of such $\WW{}$,
and in this case $\func{\rescd{\func{\indliemor{\Liegroup{1}}{\Liegroup{}}}{\Injection{\G{1}}{\G{}}}}}{\WW{}}\in
\LiealgIsom{\LiegroupLiealgebra{\Liegroup{1}}}{\subspace{\LiegroupLiealgebra{\Liegroup{}}}{\WW{}}}$.
\proof
It is assumed that $\Liegroup{1}$ is an immersed smooth subgroup of $\Liegroup{}$.
Then according to \reflem{leminjectionofsmoothsubgroupisamorphism}, $\Injection{\G{1}}{\G{}}$ is a smooth-group-morpism
from $\Liegroup{1}$ to $\Liegroup{}$. Thus according to \refthm{thminducedliealgebramorphismfromliemorphism1},
$\func{\indliemor{\Liegroup{1}}{\Liegroup{}}}{\Injection{\G{1}}{\G{}}}$ is a Lie-algebra-morphism from
$\LiegroupLiealgebra{\Liegroup{1}}$ to $\LiegroupLiealgebra{\Liegroup{}}$.\\
In addition, according to \refdef{defimmersedliesubgroup}, $\Lieman{\Liegroup{1}}$ is
an immersed submanifold of $\Lieman{\Liegroup{}}$ and hence $\der{\Injection{\G{1}}{\G{}}}{\Lieman{\Liegroup{1}}}{\Lieman{\Liegroup{}}}$
becomes a linear monomorphism from $\Tanspace{\IG{1}}{\Lieman{\Liegroup{1}}}$ to $\Tanspace{\IG{}}{\Lieman{\Liegroup{}}}$ when restricted to
$\tanspace{\IG{1}}{\Lieman{\Liegroup{1}}}$.
Thus according to \refdef{definducedliealgebramorphismfromliemorphism} and \refcor{corlinvvftanspalinearcecorrespondence},
$\func{\indliemor{\Liegroup{1}}{\Liegroup{}}}{\Injection{\G{1}}{\G{}}}$ is a Lie-algebra-monomorphism, more specifically.
Therefore according to \refthm{thmliealgebramorphismimageandkernel} and \refcor{corimageofliealgebramonomorphism},
by setting $\WW{}$ as the image of the mapping $\func{\indliemor{\Liegroup{1}}{\Liegroup{}}}{\Injection{\G{1}}{\G{}}}$,
and $\aliealgmor{}$ as the mapping $\func{\rescd{\func{\indliemor{\Liegroup{1}}{\Liegroup{}}}{\Injection{\G{1}}{\G{}}}}}{\WW{}}$ (the codomain-restriction
of the mapping $\func{\indliemor{\Liegroup{1}}{\Liegroup{}}}{\Injection{\G{1}}{\G{}}}$ to $\WW{}$),
$\subspace{\LiegroupLiealgebra{\Liegroup{}}}{\WW{}}$ is a Lie-subalgebra of $\LiegroupLiealgebra{\Liegroup{}}$
and $\aliealgmor{}$ is a Lie-algebra-isomorphism from $\LiegroupLiealgebra{\Liegroup{1}}$ to
$\subspace{\LiegroupLiealgebra{\Liegroup{}}}{\WW{}}$.
\endthm
\lemma\label{lemimmersedsubmanifoldsubgroup}
If $\Man{}$ is an $\infty$-immersed submanifold of $\Lieman{\Liegroup{}}$, and the set of all points $\M{}$ of $\Man{}$
is a subgroup of $\LieG{\Liegroup{}}$, then $\triple{\M{}}{\func{\res{\gop{}}}{\Cprod{\M{}}{\M{}}}}{\maxatlas{0}}$ is a smooth group,
and hence an immersed smooth subgroup of $\Liegroup{}$. Equivalently, if
$\Man{}$ is an $\infty$-immersed submanifold of $\Lieman{\Liegroup{}}$, and the set of all points $\M{}$ of $\Man{}$
is a subgroup of $\LieG{\Liegroup{}}$, then by defining the binary operation $\gop{\M{}}$ on $\M{}$ as the restriction of
$\gop{}$ to $\Cprod{\M{}}{\M{}}$ (that is, $\gop{\M{}}:=\func{\res{\gop{}}}{\Cprod{\M{}}{\M{}}}$),
\begin{equation}
\gop{\M{}}\in\mapdifclass{\infty}{\topprod{\Man{}}{\Man{}}}{\Man{}}.
\end{equation}
\proof
It is assumed that $\Man{}$ is an $\infty$-immersed submanifold of $\Lieman{\Liegroup{}}$, and the set of all points $\M{}$ of $\Man{}$
is a subgroup of $\LieG{\Liegroup{}}$, and ${\gop{\M{}}}:=\func{\res{\gop{}}}{\Cprod{\M{}}{\M{}}}$. Also,
$\Injec{}:=\Injection{\M{}}{\G{}}$.\\
Since $\Man{}$ is an $\infty$-immersed submanifold of $\Lieman{\Liegroup{}}$, $\Injec{}$ is an $\infty$-immersion from
$\Man{}$ to $\Lieman{\Liegroup{}}$, and therefore according to \cite[page~88,~2.6.10.~Proposition]{Berger},
for every $\point$ in $\M{}$ there exists a neighborhood $\U_{\point}$ of $\point$ in the underlying topological-space
$\mantops{\Man{}}$ of $\Man{}$ such that the image of $\U_{\point}$ under $\Injec{}$, and hence the $\U_{\point}$ itself, is an element of
$\Emsubman{\Lieman{\Liegroup{}}}$ (an embedded set of $\Lieman{\Liegroup{}}$), and $\func{\res{\Injec{}}}{\U_{\point}}$ is an
$\infty$-diffeomorphism from $\emsubman{\Man{}}{\U_{\point}}$ to
$\emsubman{\Lieman{\Liegroup{}}}{\U_{\point}}$. It is worthwhile to remark that every open set of a manifold is an embedded set of it,
and hence $\emsubman{\Man{}}{\U_{\point}}$ is naturally well-defined. So briefly,
\begin{align}\label{lemimmersedsubmanifoldsubgrouppeq1}
&\Foreach{\point}{\M{}}
\Existsis{\U_{\point}}{\mantop{\Man{}}}
\point\in\U_{\point},\cr
&\U_{\point}\in\emsubman{\Lieman{\Liegroup{}}}{\U_{\point}},~
\func{\res{\Injec{}}}{\U_{\point}}\in\Diffeo{\infty}{\emsubman{\Man{}}{\U_{\point}}}{\emsubman{\Lieman{\Liegroup{}}}{\U_{\point}}}.
\end{align}
\begin{itemize}
\item[$\pr{1}$]
Each $\point_1$ and $\point_2$ is taken as an arbitrary element of $\M{}$, and $\point_0:=\point_1\gop{}\point_2$ and
$\Injec{0}:=\func{\res{\Injec{}}}{\U_{\point_0}}$.
Since $\M{}$ is a subgroup of $\G{}$, clearly $\point_0\in\M{}$. $\U_{\point_0}$ is chosen in accordance with
\Ref{lemimmersedsubmanifoldsubgrouppeq1}.
According to \cite[page~86,~2.6.3.~Lemma]{Berger}, there exists a neighborhood $\U$ of $\point_0$ in the underlying topological-space of
$\Lieman{\Liegroup{}}$ such that $\p{\U}:=\M{}\cap\U$ is a subset of $\U_{\point_0}$, which is also inevitably an open set of the
manifold $\emsubman{\Lieman{\Liegroup{}}}{\U_{\point_0}}$. That is,
\begin{equation}\label{lemimmersedsubmanifoldsubgroupp1eq1}
\Existsis{\U}{\mantop{\Lieman{\Liegroup{}}}}
\point_0\in\p{\U}:=\M{}\cap\U\subseteq\U_{\point_0},~
\p{\U}\in\mantop{\emsubman{\Lieman{\Liegroup{}}}{\U_{\point_0}}}.
\end{equation}
Since $\p{\U}$ is an open set of the manifold $\emsubman{\Lieman{\Liegroup{}}}{\U_{\point}}$, it is also an
embedded set of this manifold and therefore according to the transitivity property of embedded submanifolds, $\p{\U}$
is ultimately an embedded set of the manifold $\Lieman{\Liegroup{}}$.
\begin{equation}\label{lemimmersedsubmanifoldsubgroupp1eq2}
\p{\U}\in\Emsubman{\Lieman{\Liegroup{}}}.
\end{equation}
By definition,
\begin{equation}\label{lemimmersedsubmanifoldsubgroupp1eq3}
\widehat{\gop{}}:=\cmp{\gop{}}{\(\funcprod{\Injec{}}{\Injec{}}\)}.
\end{equation}
Since $\Injec{}$ is a smooth map from $\Man{}$ to $\Lieman{\Liegroup{}}$, clearly
$\funcprod{\Injec{}}{\Injec{}}$ is a smooth map from the product manifold $\manprod{\Man{}}{\Man{}}$ to the
product manifold $\manprod{\Lieman{\Liegroup{}}}{\Lieman{\Liegroup{}}}$. That is,
\begin{equation}\label{lemimmersedsubmanifoldsubgroupp1eq4}
\funcprod{\Injec{}}{\Injec{}}\in\mapdifclass{\infty}{\manprod{\Man{}}{\Man{}}}
{\manprod{\Lieman{\Liegroup{}}}{\Lieman{\Liegroup{}}}}.
\end{equation}
Accordingly, since $\gop{}\in\mapdifclass{\infty}{\manprod{\Lieman{\Liegroup{}}}{\Lieman{\Liegroup{}}}}{\Lieman{\Liegroup{}}}$,
it follows from the chain-rule of differentiation that $\widehat{\gop{}}$ is a smooth map from
$\manprod{\Man{}}{\Man{}}$ to $\Lieman{\Liegroup{}}$. That is,
\begin{equation}\label{lemimmersedsubmanifoldsubgroupp1eq5}
\widehat{\gop{}}\in\mapdifclass{\infty}{\manprod{\Man{}}{\Man{}}}{\Lieman{\Liegroup{}}}.
\end{equation}
Therefore, $\widehat{\gop{}}$ is a continuous map from the underlying topological-space of the manifold $\manprod{\Man{}}{\Man{}}$
to the underlying topological-space of the manifold $\Lieman{\Liegroup{}}$ (or the smooth group $\Liegroup{}$). That is,
\begin{equation}\label{lemimmersedsubmanifoldsubgroupp1eq6}
\widehat{\gop{}}\in
\CF{\topprod{\mantops{\Man{}}}{\mantops{\Man{}}}}{\lietops{\Liegroup{}}}.
\end{equation}
Therefore, since $\U$ is an open set of the topological-space $\lietops{\Liegroup{}}=\opair{\G{}}{\mantop{\Lieman{\Liegroup{}}}}$,
$\func{\pimage{\widehat{\gop{}}}}{\U}$ is an open set of the topological-space $\topprod{\mantops{\Man{}}}{\mantops{\Man{}}}$. That is,
\begin{equation}\label{lemimmersedsubmanifoldsubgroupp1eq7}
\V:=\func{\pimage{\widehat{\gop{}}}}{\U}\in\topologyofspace{\topprod{\mantops{\Man{}}}{\mantops{\Man{}}}},
\end{equation}
and thus $\V$ is naturally an embedded set of the product manifold $\manprod{\Man{}}{\Man{}}$, that is,
\begin{equation}\label{lemimmersedsubmanifoldsubgroupp1eq8}
\V\in\Emsubman{\manprod{\Man{}}{\Man{}}},
\end{equation}
So, since $\widehat{\gop{}}$ is an smooth map from $\manprod{\Man{}}{\Man{}}$ to $\Lieman{\Liegroup{}}$,
the domain-restriction of $\widehat{\gop{}}$ to $\V$ is a smooth map from $\emsubman{\manprod{\Man{}}{\Man{}}}{\V}$
to $\Lieman{\Liegroup{}}$. That is,
\begin{equation}\label{lemimmersedsubmanifoldsubgroupp1eq9}
\ddot{\bullet}:=\func{\resd{\widehat{\gop{}}}}{\V}\in\mapdifclass{\infty}{\emsubman{\manprod{\Man{}}{\Man{}}}{\V}}{\Lieman{\Liegroup{}}}.
\end{equation}
According to \Ref{lemimmersedsubmanifoldsubgroupp1eq1}, \Ref{lemimmersedsubmanifoldsubgroupp1eq3}, and
\Ref{lemimmersedsubmanifoldsubgroupp1eq7}, and considering that $\M{}$ is a subgroup of $\LieG{\Liegroup{}}$, it is evident that
the image of the mapping $\func{\resd{\widehat{\gop{}}}}{\V}$ lies in $\p{\U}$. That is,
\begin{equation}\label{lemimmersedsubmanifoldsubgroupp1eq10}
\funcimage{\ddot{\bullet}}=\func{\image{\[\func{\resd{\widehat{\gop{}}}}{\V}\]}}{\V}\subseteq\p{\U}.
\end{equation}
Thus according to \Ref{lemimmersedsubmanifoldsubgroupp1eq2} and \Ref{lemimmersedsubmanifoldsubgroupp1eq9},
and considering that the codomain-restriction of a smooth map to any embedded submanifold of the target manifold of the smooth map
including the image of that smooth map is again a smooth map, it is clear that,
\begin{equation}\label{lemimmersedsubmanifoldsubgroupp1eq11}
\dddot{\bullet}:=\func{\rescd{\ddot{\bullet}}}{\p{\U}}\in
\mapdifclass{\infty}{\emsubman{\manprod{\Man{}}{\Man{}}}{\V}}{\emsubman{\Lieman{\Liegroup{}}}{\p{\U}}}.
\end{equation}
According to \Ref{lemimmersedsubmanifoldsubgrouppeq1},
\begin{equation}\label{lemimmersedsubmanifoldsubgroupp1eq12}
\Injec{0}\in\Diffeo{\infty}{\emsubman{\Man{}}{\U_{\point_0}}}{\emsubman{\Lieman{\Liegroup{}}}{\U_{\point_0}}}.
\end{equation}
Additionally, since $\p{\U}$ is an open set of $\emsubman{\Man{}}{\U_{\point_0}}$, $\p{\U}=\func{\image{\Injec{0}}}{\p{\U}}$
is an open set of $\emsubman{\Lieman{\Liegroup{}}}{\U_{\point_0}}$. Therefore, $\p{\U}$ is an embedded set of both of the manifolds
$\emsubman{\Man{}}{\U_{\point_0}}$ and $\emsubman{\Lieman{\Liegroup{}}}{\U_{\point_0}}$, and hence an embedded set of both $\Man{}$
and $\Lieman{\Liegroup{}}$ based on the transitivity property of embedded submanifolds. That is,
\begin{align}
&\p{\U}\in\Emsubman{\Man{}},
\label{lemimmersedsubmanifoldsubgroupp1eq13}\\
&\p{\U}\in\Emsubman{\Lieman{\Liegroup{}}}.
\label{lemimmersedsubmanifoldsubgroupp1eq14}
\end{align}
Additionally, according to \Ref{lemimmersedsubmanifoldsubgroupp1eq12},
the restriction of $\Injec{0}$ to $\p{\U}$ is a diffeomorphism from
$\emsubman{\emsubman{\Man{}}{\U_{\point_0}}}{\p{\U}}$ to
$\emsubman{\emsubman{\Lieman{\Liegroup{}}}{\U_{\point_0}}}{\p{\U}}$, and hence since,
\begin{align}
\emsubman{\emsubman{\Man{}}{\U_{\point_0}}}{\p{\U}}&=\emsubman{\Man{}}{\p{\U}},
\label{lemimmersedsubmanifoldsubgroupp1eq15}\\
\emsubman{\emsubman{\Lieman{\Liegroup{}}}{\U_{\point_0}}}{\p{\U}}&=\emsubman{\Lieman{\Liegroup{}}}{\p{\U}},
\label{lemimmersedsubmanifoldsubgroupp1eq16}
\end{align}
clearly,
\begin{equation}\label{lemimmersedsubmanifoldsubgroupp1eq17}
\func{\res{\Injec{0}}}{\p{\U}}\in
\Diffeo{\infty}{\emsubman{\Man{}}{\p{\U}}}{\emsubman{\Lieman{\Liegroup{}}}{\p{\U}}},
\end{equation}
which yields ultimately,
\begin{equation}\label{lemimmersedsubmanifoldsubgroupp1eq18}
\Injec{1}:=\finv{\[\func{\res{\Injec{0}}}{\p{\U}}\]}\in
\mapdifclass{\infty}{\emsubman{\Lieman{\Liegroup{}}}{\p{\U}}}{\emsubman{\Man{}}{\p{\U}}}.
\end{equation}
It can be easily verified that,
\begin{equation}\label{lemimmersedsubmanifoldsubgroupp1eq19}
\func{\rescd{\[\func{\resd{\gop{\M{}}}}{\V}\]}}{\p{\U}}=
\cmp{\Injec{1}}{\dddot{\bullet}}.
\end{equation}
Moreover, according to the chain-rule of differentiation,
\Ref{lemimmersedsubmanifoldsubgroupp1eq11}, \Ref{lemimmersedsubmanifoldsubgroupp1eq18},
and \Ref{lemimmersedsubmanifoldsubgroupp1eq19} yield,
\begin{equation}
\func{\rescd{\[\func{\resd{\gop{\M{}}}}{\V}\]}}{\p{\U}}\in
\mapdifclass{\infty}{\emsubman{\manprod{\Man{}}{\Man{}}}{\V}}{\emsubman{\Man{}}{\p{\U}}}.
\end{equation}
As it is known, extending the codomain of $\func{\rescd{\[\func{\resd{\gop{\M{}}}}{\V}\]}}{\p{\U}}$
to the whole $\M{}$ yields again a smooth map from $\emsubman{\manprod{\Man{}}{\Man{}}}{\V}$
to $\Man{}$. But this codomain extension trivially gives $\func{\resd{\gop{\M{}}}}{\V}$. Thus,
\begin{equation}
\func{\resd{\gop{\M{}}}}{\V}\in\mapdifclass{\infty}{\emsubman{\manprod{\Man{}}{\Man{}}}{\V}}{\Man{}}.
\end{equation}
\endp
\end{itemize}
Therefore, for every $\opair{\point_1}{\point_2}$ in $\Cprod{\M{}}{\M{}}$, there exists an open embedded submanifold $\V$
of $\manprod{\Man{}}{\Man{}}$ such that the domain-restriction of $\gop{\M{}}$ to $\V$ is smooth. That is,
\begin{align}
&\Foreach{\opair{\point_1}{\point_2}}{\Cprod{\M{}}{\M{}}}
\Exists{\V}{\mantop{\Man{}}}\cr
&\func{\resd{\gop{\M{}}}}{\V}\in\mapdifclass{\infty}{\emsubman{\manprod{\Man{}}{\Man{}}}{\V}}{\Man{}},
\end{align}
and hence it is trivially inferred that $\gop{\M{}}$ is smooth on the whole manifold $\Man{}$. That is,
$\gop{\M{}}\in\mapdifclass{\infty}{\topprod{\Man{}}{\Man{}}}{\Man{}}$.
\endlem
\lemma\label{lemdistributioninducedbyliesubalgebraofliegroupliealgebra0}
$\LS{}$ is taken as an $r$-dimensional Lie-subalgebra of $\LiegroupLiealgebra{\Liegroup{}}$
(the canonical Lie-algebra of $\Liegroup{}$) for some $r$ in $\Zp$, and $\function{\avecf{}}{\seta{\suc{1}{r}}}{\LS{}}$
as an ordered-basis of the vector-space $\subspace{\VLeftinvvf{\Liegroup{}}}{\LS{}}$.
For every point $\g{}$ of $\Liegroup{}$, $\defSet{\func{\avecf{}}{\g{}}}{\avecf{}\in\LS{}}$ is an $r$-dimensional vector-subspace
of $\Tanspace{\g{}}{\Lieman{\Liegroup{}}}$ (and cosequently an element of $\subtanbun{\Lieman{\Liegroup{}}}{r}$), and
$\mtuple{\func{\avecf{1}}{\g{}}}{\func{\avecf{r}}{\g{}}}$ is an ordered-basis for
this vector-subspace of the tangent-space of $\Lieman{\Liegroup{}}$ at $\g{}$. That is,
\begin{align}
\Foreach{\g{}}{\G{}}
\left\{\begin{aligned}
&\defSet{\func{\avecf{}}{\g{}}}{\avecf{}\in\LS{}}\in\subvec{\Tanspace{\g{}}{\Lieman{\Liegroup{}}}}{r}\subseteq\subtanbun{\Lieman{\Liegroup{}}}{r},\cr
&\mtuple{\func{\avecf{1}}{\g{}}}{\func{\avecf{r}}{\g{}}}\in
\ovecbasis{\subspace{\Tanspace{\g{}}{\Lieman{\Liegroup{}}}}{\defSet{\func{\avecf{}}{\g{}}}{\avecf{}\in\LS{}}}}.
\end{aligned}\right.
\end{align}
\proof
\begin{itemize}
\item[$\pr{1}$]
$\g{}$ is taken as an arbitrary element of $\G{}$. Since $\gltrans{\LieG{\Liegroup{}}}{\g{}}\in\Diff{\infty}{\Lieman{\Liegroup{}}}$,
$\func{\gltrans{\LieG{\Liegroup{}}}{\g{}}}{\IG{}}=\g{}$, clearly the restrction of
$\der{\gltrans{\LieG{\Liegroup{}}}{\g{}}}{\Lieman{\Liegroup{}}}{\Lieman{\Liegroup{}}}$ to the tangent-space of $\Lieman{\Liegroup{}}$ at $\IG{}$
is a linear-isomorphism from $\Tanspace{\IG{}}{\Lieman{\Liegroup{}}}$ to $\Tanspace{\g{}}{\Lieman{\Liegroup{}}}$.
\begin{equation}
\[\func{\res{\(\der{\gltrans{\LieG{\Liegroup{}}}{\g{}}}{\Lieman{\Liegroup{}}}{\Lieman{\Liegroup{}}}\)}}
{\Tanspace{\IG{}}{\Lieman{\Liegroup{}}}}\]\in
\Linisom{\Tanspace{\IG{}}{\Lieman{\Liegroup{}}}}{\Tanspace{\g{}}{\Lieman{\Liegroup{}}}}.
\end{equation}
Moreover, according to \refcor{corlinvvftanspalinearcecorrespondence}, $\liegvftan{\Liegroup{}}$ is a linear-isomorphism
from $\VLeftinvvf{\Liegroup{}}$ to $\Tanspace{\IG{}}{\Lieman{\Liegroup{}}}$, that is,
\begin{equation}
\liegvftan{\Liegroup{}}\in\Linisom{\VLeftinvvf{\Liegroup{}}}{\Tanspace{\IG{}}{\Lieman{\Liegroup{}}}}.
\end{equation}
Therefore,
\begin{equation}
\eta:=\cmp{\[\func{\res{\(\der{\gltrans{\LieG{\Liegroup{}}}{\g{}}}{\Lieman{\Liegroup{}}}{\Lieman{\Liegroup{}}}\)}}
{\Tanspace{\IG{}}{\Lieman{\Liegroup{}}}}\]}{\liegvftan{\Liegroup{}}}\in
\Linisom{\VLeftinvvf{\Liegroup{}}}{\Tanspace{\g{}}{\Lieman{\Liegroup{}}}}.
\end{equation}
Moreover, according to \reflem{lemleftinvariantvectorfieldsequiv0} and \refdef{deflinvvftanspacecorrespondence},
\begin{align}
\Foreach{\avecf{}}{\Leftinvvf{\Liegroup{}}}
\func{\eta}{\avecf{}}&=\func{\(\der{\gltrans{\LieG{\Liegroup{}}}{\g{}}}{\Lieman{\Liegroup{}}}{\Lieman{\Liegroup{}}}\)}
{\func{\liegvftan{\Liegroup{}}}{\avecf{}}}\cr
&=\func{\(\der{\gltrans{\LieG{\Liegroup{}}}{\g{}}}{\Lieman{\Liegroup{}}}{\Lieman{\Liegroup{}}}\)}{\func{\avecf{}}{\IG{}}}\cr
&=\func{\avecf{}}{\g{}},
\end{align}
and thus,
\begin{equation}
\func{\image{\eta}}{\LS{}}=
\defSet{\func{\avecf{}}{\g{}}}{\avecf{}\in\LS{}}.
\end{equation}
Therefore, considering that $\LS{}$ is a vector-subspace of $\VLeftinvvf{\Liegroup{}}$ and
a linear-isomorphism maps a vector-subspace injectively onto a vector-subspace of the same dimension, $\func{\image{\eta}}{\LS{}}$
is an $r$-dimensional vector-subspace of $\Tanspace{\g{}}{\Lieman{\Liegroup{}}}$.
\begin{equation}
\defSet{\func{\avecf{}}{\g{}}}{\avecf{}\in\LS{}}\in\subvec{\Tanspace{\g{}}{\Lieman{\Liegroup{}}}}{r}.
\end{equation}
Also, since a linear-isomorphism sends an ordered-basis of a vector-subspace of the source space
to an ordered-basis of the corresponded vector-subspace of the target space,
\begin{align}
\mtuple{\func{\avecf{1}}{\g{}}}{\func{\avecf{r}}{\g{}}}=
\mtuple{\func{\eta}{\avecf{1}}}{\func{\eta}{\avecf{r}}}\in
\ovecbasis{\subspace{\Tanspace{\g{}}{\Lieman{\Liegroup{}}}}{\defSet{\func{\avecf{}}{\g{}}}{\avecf{}\in\LS{}}}}.
\end{align}
\end{itemize}
\endlem
\definition\label{defdistributioninducedbyliesubalgebraofliegroupliealgebra}
$r$ is taken as an element of $\Zp$ (a positive integer). The mapping $\subliealgdist{\Liegroup{}}{r}$ is defined as,
\begin{align}
&\subliealgdist{\Liegroup{}}{r}\indef\Func{\subliedim{r}{\LiegroupLiealgebra{\Liegroup{}}}}
{\Func{\M{}}{\subtanbun{\Lieman{\Liegroup{}}}{r}}},\cr
&\Foreach{\LS{}}{\subliedim{r}{\LiegroupLiealgebra{\Liegroup{}}}}
\Foreach{\g{}}{\G{}}
\func{\[\func{\subliealgdist{\Liegroup{}}{r}}{\LS{}}\]}{\g{}}\eqdef
\defSet{\func{\avecf{}}{\g{}}}{\avecf{}\in\LS{}}.
\end{align}
\endef
\theorem\label{thmdistributioninducedbyliesubalgebraofliegroupliealgebraissmooth}
$r$ is taken as an element of $\Zp$. For every $\LS{}$ in
$\subliedim{r}{\LiegroupLiealgebra{\Liegroup{}}}$ (every $r$-dimensional Lie-subalgebra of
the canonical Lie-algebra of $\Liegroup{}$, $\LiegroupLiealgebra{\Liegroup{}}$),
$\func{\subliealgdist{\Liegroup{}}{r}}{\LS{}}$ is an $r$-dimensional $\infty$-distribution of the manifold $\Lieman{\Liegroup{}}$. That is,
\begin{equation}
\func{\subliealgdist{\Liegroup{}}{r}}{\LS{}}\in\Distributions{r}{\Lieman{\Liegroup{}}}.
\end{equation}
\proof
The verification is straightforward according to \reflem{lemdistributioninducedbyliesubalgebraofliegroupliealgebra0},
\refdef{defdistribution}, and \refdef{defdistributioninducedbyliesubalgebraofliegroupliealgebra}.
\endthm
\theorem\label{thmconsistentvectorfieldsofsubalgebradistributionisafreemodule0}
$r$ is taken as an element of $\Zp$, $\LS{}$ as an element of $\subliedim{r}{\LiegroupLiealgebra{\Liegroup{}}}$
(an $r$-dimensional Lie-subalgebra of
the canonical Lie-algebra of the smooth group $\Liegroup{}$, that is $\LiegroupLiealgebra{\Liegroup{}}$), and
$\function{\avecf{}}{\seta{\suc{1}{r}}}{\Leftinvvf{\Liegroup{}}}$ as an ordered-basis of the Lie-algebra
$\subspace{\LiegroupLiealgebra{\Liegroup{}}}{\LS{}}$ (or equivalently an ordered-basis of the vector-space
$\subspace{\VLeftinvvf{\Liegroup{}}}{\LS{}}$).
$\seta{\suc{\avecf{1}}{\avecf{r}}}$ is also a module-basis of the $\smoothring{\infty}{\Lieman{\Liegroup{}}}$-module
$\subspace{\smoothvfmodule{\infty}{\Lieman{\Liegroup{}}}}{\distvecf{\Lieman{\Liegroup{}}}{\func{\subliealgdist{\Liegroup{}}{r}}{\LS{}}}{r}}$.
\begin{equation}\label{thmconsistentvectorfieldsofsubalgebradistributionisafreemodule0eq1}
\seta{\suc{\avecf{1}}{\avecf{r}}}\in\Modulebases{
\subspace{\smoothvfmodule{\infty}{\Lieman{\Liegroup{}}}}{\distvecf{\Lieman{\Liegroup{}}}{\func{\subliealgdist{\Liegroup{}}{r}}{\LS{}}}{r}}}.
\end{equation}
Equivalently, the set $\seta{\suc{\avecf{1}}{\avecf{r}}}$ of cardinality $r$ is a linearly-independent set of the submodule
$\distvecf{\Lieman{\Liegroup{}}}{\func{\subliealgdist{\Liegroup{}}{r}}{\LS{}}}{r}$
of the $\smoothring{\infty}{\Lieman{\Liegroup{}}}$-module $\smoothvfmodule{\infty}{\Lieman{\Liegroup{}}}$,
and also generates this submodule, that is,
\begin{align}\label{thmconsistentvectorfieldsofsubalgebradistributionisafreemodule0eq2}
\distvecf{\Lieman{\Liegroup{}}}{\func{\subliealgdist{\Liegroup{}}{r}}{\LS{}}}{r}&=
\func{\modulegen{\[\subspace{\smoothvfmodule{\infty}{\Lieman{\Liegroup{}}}}{\distvecf{\Lieman{\Liegroup{}}}
{\func{\subliealgdist{\Liegroup{}}{r}}{\LS{}}}{r}}\]}}{\seta{\suc{\avecf{1}}{\avecf{r}}}}\cr
&=\defSet{\sum_{k=1}^{r}\cf_k\smoothvfprod{\Lieman{\Liegroup{}}}\avecf{k}}
{\mtuple{\cf_1}{\cf_r}\in{\mapdifclass{\infty}{\Lieman{\Liegroup{}}}{\RR}}^{\times r}}.
\end{align}
\caution\refthm{thmsubmoduleinducedbydistribution} guaranties that the set of all smooth vector-fields on $\Lieman{\Liegroup{}}$
consistent with the $\infty$-distribution $\func{\subliealgdist{\Liegroup{}}{r}}{\LS{}}$ of $\Lieman{\Liegroup{}}$, that is
$\distvecf{\Lieman{\Liegroup{}}}{\func{\subliealgdist{\Liegroup{}}{r}}{\LS{}}}{r}$, is a
submodule of the $\smoothring{\infty}{\Lieman{\Liegroup{}}}$-module $\smoothvfmodule{\infty}{\Lieman{\Liegroup{}}}$
consisting of all smooth vector-fields on the underlying manifold $\Lieman{\Liegroup{}}$ of the smooth group $\Liegroup{}$.
\proof
It is immediate to verify that,
\begin{equation}\label{thmconsistentvectorfieldsofsubalgebradistributionisafreemodule0peq1}
\defSet{\sum_{k=1}^{r}\cf_k\smoothvfprod{\Lieman{\Liegroup{}}}\avecf{k}}
{\mtuple{\cf_1}{\cf_r}\in{\mapdifclass{\infty}{\Lieman{\Liegroup{}}}{\RR}}^{\times r}}\subseteq
\distvecf{\Lieman{\Liegroup{}}}{\func{\subliealgdist{\Liegroup{}}{r}}{\LS{}}}{r}.
\end{equation}
\begin{itemize}
\item[$\pr{1}$]
$\avecff{}$ is taken as an arbitrary element of $\distvecf{\Lieman{\Liegroup{}}}{\func{\subliealgdist{\Liegroup{}}{r}}{\LS{}}}{r}$.
Then according to \refdef{defconsistentvectorfieldswithdistributions}, $\avecff{}$ is a smooth vector-field of $\Lieman{\Liegroup{}}$
whose value at every point $\g{}$ of $\Liegroup{}$ lies in the vector-subspace
$\func{\[\func{\subliealgdist{\Liegroup{}}{r}}{\LS{}}\]}{\g{}}$ of the tangent-space of $\Lieman{\Liegroup{}}$ at $\g{}$. That is,
according to \refdef{defdistributioninducedbyliesubalgebraofliegroupliealgebra},
\begin{align}
&\avecff{}\in\vecf{\Lieman{\Liegroup{}}}{\infty},
\label{thmconsistentvectorfieldsofsubalgebradistributionisafreemodule0p1eq1}\\
&\Foreach{\g{}}{\G{}}\func{\avecff{}}{\g{}}\in
\func{\[\func{\subliealgdist{\Liegroup{}}{r}}{\LS{}}\]}{\g{}}=
\defSet{\func{\avecf{}}{\g{}}}{\avecf{}\in\LS{}}.
\label{thmconsistentvectorfieldsofsubalgebradistributionisafreemodule0p1eq2}
\end{align}
According to \reflem{lemdistributioninducedbyliesubalgebraofliegroupliealgebra0},
for every point $\g{}$ of $\Liegroup{}$, $\mtuple{\func{\avecf{1}}{\g{}}}{\func{\avecf{r}}{\g{}}}$ is an ordered-basis
of the vector-subspace $\func{\[\func{\subliealgdist{\Liegroup{}}{r}}{\LS{}}\]}{\g{}}$ of $\Tanspace{\g{}}{\Lieman{\Liegroup{}}}$.
Therefore, since $\func{\avecff{}}{\g{}}$ lies in this vector-space for every point $\g{}$ of $\Liegroup{}$,
$\func{\avecff{}}{\g{}}$ can be represented by a unique sequence of $r$ scalars $\suc{\zeta_1}{\zeta_r}$
as the coefficients of the expansion of $\func{\avecff{}}{\g{}}$
with respect to the ordered-basis $\mtuple{\func{\avecf{1}}{\g{}}}{\func{\avecf{r}}{\g{}}}$ for every $\g{}$ in $\G{}$. Precisely,
there is a unique sequence of well-defined functions $\suc{\zeta_1}{\zeta_r}$ such that,
\begin{align}\label{thmconsistentvectorfieldsofsubalgebradistributionisafreemodule0p1eq3}
&\Foreach{k}{\seta{\suc{1}{r}}}\zeta_k\indef\Func{\G{}}{\R},\cr
&\Foreach{\g{}}{\G{}}\func{\avecff{}}{\g{}}\eqdef\sum_{k=1}^{n}\func{\zeta_k}{\g{}}\func{\avecf{k}}{\g{}}.
\end{align}
\begin{itemize}
\item[$\pr{1-1}$]
$\point$ is taken as an arbitrary element of $\G{}$ and
$\phi$ is chosen as a chart of the manifold $\Lieman{\Liegroup{}}$ around $\point$, that is an element of
$\defset{\p{\phi}}{\maxatlas{}}{\point\in\domain{\phi}}$. Also
$\U{}:=\domain{\phi}$ and $z:=\func{\phi}{\point}$. It is known that $\tanspaceiso{\g{}}{\Lieman{\Liegroup{}}}{\phi}$ is a linear-isomorphism
from the vector-space $\Tanspace{\g{}}{\Lieman{\Liegroup{}}}$ to the vector-space $\R^n$ for every $\g{}$ in $\U{}$.
\begin{equation}\label{thmconsistentvectorfieldsofsubalgebradistributionisafreemodule0p11eq1}
\Foreach{\g{}}{\U{}}
\tanspaceiso{\g{}}{\Lieman{\Liegroup{}}}{\phi}\in\Linisom{\Tanspace{\g{}}{\Lieman{\Liegroup{}}}}{\R^n}.
\end{equation}
Since $\avecff{}$, and $\suc{\avecf{1}}{\avecf{r}}$ are vector-fields on $\Lieman{\Liegroup{}}$, they assign to each point
of $\Lieman{\Liegroup{}}$ a vector of the tangent-space of $\Lieman{\Liegroup{}}$ at that point, and hence,
\begin{align}\label{thmconsistentvectorfieldsofsubalgebradistributionisafreemodule0p11eq2}
\Foreach{\g{}}{\U{}}
\left\{
\begin{aligned}
&\func{\avecff{}}{\g{}}\in\tanspace{\g{}}{\Lieman{\Liegroup{}}},\cr
&\Foreach{k}{\seta{\suc{1}{r}}}\func{\avecf{k}}{\g{}}\in\tanspace{\g{}}{\Lieman{\Liegroup{}}}.
\end{aligned}\right.
\end{align}
\Ref{thmconsistentvectorfieldsofsubalgebradistributionisafreemodule0p1eq3},
\Ref{thmconsistentvectorfieldsofsubalgebradistributionisafreemodule0p11eq1}, and
\Ref{thmconsistentvectorfieldsofsubalgebradistributionisafreemodule0p11eq2} yield,
\begin{equation}\label{thmconsistentvectorfieldsofsubalgebradistributionisafreemodule0p11eq3}
\Foreach{\g{}}{\U{}}
\func{\[\cmp{\tanspaceiso{\g{}}{\Lieman{\Liegroup{}}}{\phi}}{\avecff{}}\]}{\g{}}=
\sum_{k=1}^{n}\func{\zeta_k}{\g{}}\(\func{\[\cmp{\tanspaceiso{\g{}}{\Lieman{\Liegroup{}}}{\phi}}{\avecf{k}}\]}{\g{}}\),
\end{equation}
which immediately implies,
\begin{align}\label{thmconsistentvectorfieldsofsubalgebradistributionisafreemodule0p11eq4}
&\Foreach{\x}{\func{\image{\phi}}{\U{}}}\cr
&\func{\[\cmp{\cmp{\tanspaceiso{\func{\finv{\phi}}{\x}}{\Lieman{\Liegroup{}}}{\phi}}{\avecff{}}}{\finv{\phi}}\]}{\x}=
\sum_{k=1}^{n}\bigg(\func{\(\cmp{\zeta_k}{\finv{\phi}}\)}{\x}\bigg)
\(\func{\[\cmp{\cmp{\tanspaceiso{\func{\finv{\phi}}{\x}}{\Lieman{\Liegroup{}}}{\phi}}{\avecf{k}}}{\finv{\phi}}\]}{\x}\).
\end{align}
Since $\avecff{}$, and $\suc{\avecf{1}}{\avecf{r}}$ are smooth vector-fields on $\Lieman{\Liegroup{}}$,
\begin{align}\label{thmconsistentvectorfieldsofsubalgebradistributionisafreemodule0p11eq5}
&\[\cmp{\cmp{\tanchart{\Lieman{\Liegroup{}}}{\phi}}{\avecff{}}}{\finv{\phi}}\]\in
\banachmapdifclass{\infty}{\R^n}{\Cprod{\R^n}{\R^n}}{\func{\image{\phi}}{\U{}}}{\Cprod{\func{\image{\phi}}{\U{}}}{\R^n}},\cr
&\Foreach{k}{\seta{\suc{1}{r}}}
\[\cmp{\cmp{\tanchart{\Lieman{\Liegroup{}}}{\phi}}{\avecf{k}}}{\finv{\phi}}\]\in
\banachmapdifclass{\infty}{\R^n}{\Cprod{\R^n}{\R^n}}{\func{\image{\phi}}{\U{}}}{\Cprod{\func{\image{\phi}}{\U{}}}{\R^n}}.\cr
&{}
\end{align}
Therefore by defining the mappings $\beta$, and $\suc{\alpha_1}{\alpha_r}$ as,
\begin{align}\label{thmconsistentvectorfieldsofsubalgebradistributionisafreemodule0p11eq6}
&\beta\indef\Func{\func{\image{\phi}}{\U{}}}{\R^n},\cr
&\Foreach{\x}{\func{\image{\phi}}{\U{}}}\func{\beta}{\x}\eqdef
\func{\[\cmp{\cmp{\tanspaceiso{\func{\finv{\phi}}{\x}}{\Lieman{\Liegroup{}}}{\phi}}{\avecff{}}}{\finv{\phi}}\]}{\x},
\end{align}
and,
\begin{align}\label{thmconsistentvectorfieldsofsubalgebradistributionisafreemodule0p11eq7}
\Foreach{k}{\seta{\suc{1}{r}}}
\left\{\begin{aligned}
&\alpha_k\indef\Func{\func{\image{\phi}}{\U{}}}{\R^n},\cr
&\Foreach{\x}{\func{\image{\phi}}{\U{}}}\func{\alpha_k}{\x}\eqdef
\func{\[\cmp{\cmp{\tanspaceiso{\func{\finv{\phi}}{\x}}{\Lieman{\Liegroup{}}}{\phi}}{\avecf{k}}}{\finv{\phi}}\]}{\x},
\end{aligned}\right.
\end{align}
these mappings are smooth maps from $\func{\image{\phi}}{\U{}}$ to $\R^n$, that is,
\begin{align}
\beta&\in\banachmapdifclass{\infty}{\R^n}{\R^n}{\func{\image{\phi}}{\U{}}}{\R^n},
\label{thmconsistentvectorfieldsofsubalgebradistributionisafreemodule0p11eq8}\\
\Foreach{k}{\seta{\suc{1}{r}}}
\alpha_k&\in\banachmapdifclass{\infty}{\R^n}{\R^n}{\func{\image{\phi}}{\U{}}}{\R^n},
\label{thmconsistentvectorfieldsofsubalgebradistributionisafreemodule0p11eq9}
\end{align}
because,
\begin{align}\label{thmconsistentvectorfieldsofsubalgebradistributionisafreemodule0p11eq10}
\[\cmp{\cmp{\tanchart{\Lieman{\Liegroup{}}}{\phi}}{\avecff{}}}{\finv{\phi}}\]&=
\funcprod{\identity{\func{\image{\phi}}{\U{}}}}{\beta},\cr
\Foreach{k}{\seta{\suc{1}{r}}}
\[\cmp{\cmp{\tanchart{\Lieman{\Liegroup{}}}{\phi}}{\avecf{k}}}{\finv{\phi}}\]&=
\funcprod{\identity{\func{\image{\phi}}{\U{}}}}{\alpha_k}.
\end{align}
By defining the mappings,
\begin{align}
\Foreach{j}{\seta{\suc{1}{n}}}&\beta^{j}:=\cmp{\projection{n}{j}}{\beta},
\label{thmconsistentvectorfieldsofsubalgebradistributionisafreemodule0p11eq11}\\
\Foreach{j}{\seta{\suc{1}{n}}}\Foreach{k}{\seta{\suc{1}{r}}}
&\alpha_k^{j}:=\cmp{\projection{n}{j}}{\alpha_k}.
\label{thmconsistentvectorfieldsofsubalgebradistributionisafreemodule0p11eq12}
\end{align}
it is a trivial consequence of \Ref{thmconsistentvectorfieldsofsubalgebradistributionisafreemodule0p11eq8}
and \Ref{thmconsistentvectorfieldsofsubalgebradistributionisafreemodule0p11eq9} that,
\begin{align}
\Foreach{j}{\seta{\suc{1}{n}}}&\beta^{j}\in\banachmapdifclass{\infty}{\R^n}{\R}{\func{\image{\phi}}{\U{}}}{\R},
\label{thmconsistentvectorfieldsofsubalgebradistributionisafreemodule0p11eq13}\\
\Foreach{j}{\seta{\suc{1}{n}}}\Foreach{k}{\seta{\suc{1}{r}}}
&\alpha_k^{j}\in\banachmapdifclass{\infty}{\R^n}{\R}{\func{\image{\phi}}{\U{}}}{\R}.
\label{thmconsistentvectorfieldsofsubalgebradistributionisafreemodule0p11eq14}
\end{align}
According to \Ref{thmconsistentvectorfieldsofsubalgebradistributionisafreemodule0p11eq4},
\Ref{thmconsistentvectorfieldsofsubalgebradistributionisafreemodule0p11eq6},
\Ref{thmconsistentvectorfieldsofsubalgebradistributionisafreemodule0p11eq7}, and
\Ref{thmconsistentvectorfieldsofsubalgebradistributionisafreemodule0p11eq11}, it is clear that,
\begin{equation}\label{thmconsistentvectorfieldsofsubalgebradistributionisafreemodule0p11eq15}
\Foreach{j}{\seta{\suc{1}{n}}}
\Foreach{\x}{\func{\image{\phi}}{\U{}}}
\func{\beta^{j}}{\x}=\sum_{k=1}^{n}\func{\alpha_k^{j}}{\x}
\bigg(\func{\(\cmp{\zeta_k}{\finv{\phi}}\)}{\x}\bigg).
\end{equation}
By defining the $\R$-matrix-valued mappings $\tilde{\alpha}$, $\tilde{\beta}$, and $\tilde{\gamma}$ as,
\begin{align}\label{thmconsistentvectorfieldsofsubalgebradistributionisafreemodule0p11eq16}
&\tilde{\alpha}\indef\Func{\func{\image{\phi}}{\U{}}}{\Mat{\R}{n}{r}},\cr
&\Foreach{\x}{\func{\image{\phi}}{\U{}}}
\Foreach{j}{\seta{\suc{1}{n}}}\Foreach{k}{\seta{\suc{1}{r}}}
\matelement{\func{\tilde{\alpha}}{\x}}{j}{k}\eqdef\func{\alpha_k^{j}}{\x},
\end{align}
and,
\begin{align}\label{thmconsistentvectorfieldsofsubalgebradistributionisafreemodule0p11eq17}
&\tilde{\beta}\indef\Func{\func{\image{\phi}}{\U{}}}{\Mat{\R}{n}{1}},\cr
&\Foreach{\x}{\func{\image{\phi}}{\U{}}}
\Foreach{j}{\seta{\suc{1}{n}}}
\matelement{\func{\tilde{\beta}}{\x}}{j}{}\eqdef\func{\beta^{j}}{\x},
\end{align}
and,
\begin{align}\label{thmconsistentvectorfieldsofsubalgebradistributionisafreemodule0p11eq18}
&\tilde{\gamma}\indef\Func{\func{\image{\phi}}{\U{}}}{\Mat{\R}{n}{1}},\cr
&\Foreach{\x}{\func{\image{\phi}}{\U{}}}
\Foreach{k}{\seta{\suc{1}{r}}}
\matelement{\func{\tilde{\gamma}}{\x}}{k}{}\eqdef\func{\(\cmp{\zeta_k}{\finv{\phi}}\)}{\x},
\end{align}
\Ref{thmconsistentvectorfieldsofsubalgebradistributionisafreemodule0p11eq15} becomes,
\begin{equation}\label{thmconsistentvectorfieldsofsubalgebradistributionisafreemodule0p11eq19}
\Foreach{\x}{\func{\image{\phi}}{\U{}}}
\func{\tilde{\beta}}{\x}=\func{\tilde{\alpha}}{\x}\cdot\func{\tilde{\gamma}}{\x},
\end{equation}
where $\cdot$ denotes the matrix-multiplication operation. By denoting the matrices by their entries,
\begin{equation}\label{thmconsistentvectorfieldsofsubalgebradistributionisafreemodule0p11eq20}
\Foreach{\x}{\func{\image{\phi}}{\U{}}}
\begin{pmatrix}
\func{\beta_1}{\x}\\
\vdots\\
\func{\beta_n}{\x}
\end{pmatrix}=
\begin{pmatrix}
\func{\alpha^{1}_{1}}{\x} & \cdots & \func{\alpha^{1}_r}{\x}\\
\vdots & \ddots & \vdots\\
\func{\alpha^{n}_{1}}{\x} & \cdots & \func{\alpha^{n}_r}{\x}
\end{pmatrix}
\begin{pmatrix}
\func{\(\cmp{\zeta_1}{\finv{\phi}}\)}{\x}\\
\vdots\\
\func{\(\cmp{\zeta_r}{\finv{\phi}}\)}{\x}
\end{pmatrix}.
\end{equation}
Since
$\mtuple{\func{\avecf{1}}{\func{\finv{\phi}}{z}}}{\func{\avecf{r}}{\func{\finv{\phi}}{z}}}$
is a linearly independet
sequence of $\Tanspace{\point}{\Lieman{\Liegroup{}}}$ and $\tanspaceiso{\func{\finv{\phi}}{z}}{\Lieman{\Liegroup{}}}{\phi}$
is a linear-isomorphism from $\Tanspace{\point}{\Lieman{\Liegroup{}}}$ to $\R^n$,
trivially according to \Ref{thmconsistentvectorfieldsofsubalgebradistributionisafreemodule0p11eq7}, the $n$-tuple
$\mtuple{\func{\alpha_1}{z}}{\func{\alpha_r}{z}}$ is a liearly-independent sequence of $\R^n$.
Thus according to \Ref{thmconsistentvectorfieldsofsubalgebradistributionisafreemodule0p11eq12}
and \Ref{thmconsistentvectorfieldsofsubalgebradistributionisafreemodule0p11eq16},
all $r$ columns of $\func{\tilde{\alpha}}{z}$ form
a linearly independent sequence of the space $\BMat{\R}{n}{1}$, and hence the rank of the matrix $\func{\tilde{\alpha}}{z}$
is $r$. Therefore, the row-rank of the matrix $\func{\tilde{\alpha}}{z}$ is again $r$ and thus there exists
a selection $\function{s}{\seta{\suc{1}{r}}}{\seta{\suc{1}{n}}}$ (which is an injection) of rows of $\func{\tilde{\alpha}}{z}$
associated with a non-singular square sub-matrix of $\func{\tilde{\alpha}}{z}$ of degree $r$. That is,
\begin{align}\label{thmconsistentvectorfieldsofsubalgebradistributionisafreemodule0p11eq21}
\Existsis{s}{\Func{\seta{\suc{1}{r}}}{\seta{\suc{1}{n}}}}
\Det{r}
\begin{pmatrix}
\func{\alpha^{s_1}_{1}}{z} & \cdots & \func{\alpha^{s_1}_r}{z}\\
\vdots & \ddots & \vdots\\
\func{\alpha^{s_n}_{1}}{z} & \cdots & \func{\alpha^{s_n}_r}{z}
\end{pmatrix}
\neq 0.
\end{align}
Therefore by defining the $\R$-matrix-valued mapping ${\tilde{\alpha}}_0$ as,
\begin{align}\label{thmconsistentvectorfieldsofsubalgebradistributionisafreemodule0p11eq22}
&{\tilde{\alpha}}_0\indef\Func{\func{\image{\phi}}{\U{}}}{\Mat{\R}{r}{r}},\cr
&\Foreach{\x}{\func{\image{\phi}}{\U{}}}
\Foreach{\opair{j}{k}}{\seta{\suc{1}{r}}^{\times 2}}
\func{{\tilde{\alpha}}_0}{\x}\eqdef\matelement{\tilde{\alpha}}{s_j}{k}=\alpha^{s_j}_{k},
\end{align}
it is evident that,
\begin{equation}\label{thmconsistentvectorfieldsofsubalgebradistributionisafreemodule0p11eq23}
\func{\Det{r}}{\func{\tilde{\alpha}_0}{z}}\neq 0.
\end{equation}
It is also clear that by defining the $\R$-matrix-valued mapping ${\tilde{\beta}}_0$ as,
\begin{align}\label{thmconsistentvectorfieldsofsubalgebradistributionisafreemodule0p11eq24}
&{\tilde{\beta}}_0\indef\Func{\func{\image{\phi}}{\U{}}}{\Mat{\R}{r}{1}},\cr
&\Foreach{\x}{\func{\image{\phi}}{\U{}}}
\Foreach{j}{\seta{\suc{1}{r}}}
\func{{\tilde{\beta}}_0}{\x}\eqdef\matelement{\tilde{\beta}}{s_j}{}=\beta^{s_j},
\end{align}
\Ref{thmconsistentvectorfieldsofsubalgebradistributionisafreemodule0p11eq19} implies,
\begin{equation}\label{thmconsistentvectorfieldsofsubalgebradistributionisafreemodule0p11eq25}
\Foreach{\x}{\func{\image{\phi}}{\U{}}}
\func{\tilde{\beta}_0}{\x}=\func{\tilde{\alpha}_0}{\x}\cdot\func{\tilde{\gamma}}{\x},
\end{equation}
or by using the standard matrix notation via entries,
\begin{equation}\label{thmconsistentvectorfieldsofsubalgebradistributionisafreemodule0p11eq26}
\Foreach{\x}{\func{\image{\phi}}{\U{}}}
\begin{pmatrix}
\func{\beta_{s_1}}{\x}\\
\vdots\\
\func{\beta_{s_r}}{\x}
\end{pmatrix}=
\begin{pmatrix}
\func{\alpha^{s_1}_{1}}{\x} & \cdots & \func{\alpha^{s_1}_r}{\x}\\
\vdots & \ddots & \vdots\\
\func{\alpha^{s_n}_{1}}{\x} & \cdots & \func{\alpha^{s_n}_r}{\x}
\end{pmatrix}
\begin{pmatrix}
\func{\(\cmp{\zeta_1}{\finv{\phi}}\)}{\x}\\
\vdots\\
\func{\(\cmp{\zeta_r}{\finv{\phi}}\)}{\x}
\end{pmatrix}.
\end{equation}
According to \Ref{thmconsistentvectorfieldsofsubalgebradistributionisafreemodule0p11eq13},
\Ref{thmconsistentvectorfieldsofsubalgebradistributionisafreemodule0p11eq14},
\Ref{thmconsistentvectorfieldsofsubalgebradistributionisafreemodule0p11eq16},
and \Ref{thmconsistentvectorfieldsofsubalgebradistributionisafreemodule0p11eq17},
all of the entries of the $\R$-matrix-valued mappings $\tilde{\beta}_0$ and $\tilde{\alpha}_0$ are smooth
real-valued mappings and hence the $\R$-matrix-valued mappings $\tilde{\alpha}_0$ and $\tilde{\beta}_0$
are smooth when the space of matrices are endowed with their canonical norm. That is,
\begin{align}
\tilde{\alpha}_0&\in\banachmapdifclass{\infty}{\R^n}{\BMat{\R}{r}{r}}{\func{\image{\phi}}{\U{}}}{\Mat{\R}{r}{r}},
\label{thmconsistentvectorfieldsofsubalgebradistributionisafreemodule0p11eq27}\\
\tilde{\beta}_0&\in\banachmapdifclass{\infty}{\R^n}{\BMat{\R}{r}{1}}{\func{\image{\phi}}{\U{}}}{\Mat{\R}{r}{1}},
\label{thmconsistentvectorfieldsofsubalgebradistributionisafreemodule0p11eq28}
\end{align}
and as an immediate consequence, $\tilde{\alpha}_0$ is also continuous, that is,
\begin{align}\label{thmconsistentvectorfieldsofsubalgebradistributionisafreemodule0p11eq29}
\tilde{\alpha}_0&\in\banachmapdifclass{0}{\R^n}{\BMat{\R}{r}{r}}{\func{\image{\phi}}{\U{}}}{\Mat{\R}{r}{r}}.
\end{align}
Since $\Det{r}$ is a continuous map from the Banach-space of the square $\R$-matrices of degree $r$ to $\R$,
\Ref{thmconsistentvectorfieldsofsubalgebradistributionisafreemodule0p11eq29} implies that the composition of $\Det{r}$ and
the $\R$-matrix-valued mapping $\tilde{\alpha}_0$ is also a continuous map from $\func{\image{\phi}}{\U{}}$ to $\R$. That is,
\begin{equation}\label{thmconsistentvectorfieldsofsubalgebradistributionisafreemodule0p11eq30}
\cmp{\Det{r}}{\tilde{\alpha}_0}\in
\banachmapdifclass{0}{\R^n}{\R}{\func{\image{\phi}}{\U{}}}{\R}.
\end{equation}
Therefore since $\cmp{\Det{r}}{\tilde{\alpha}_0}$ is non-zero at $z$
(\Ref{thmconsistentvectorfieldsofsubalgebradistributionisafreemodule0p11eq23}), there exists an open neighborhood $\V{}$ of
$z$ included in $\func{\image{\phi}}{\U{}}$ such that the restriction of $\cmp{\Det{r}}{\tilde{\alpha}_0}$ to $\V{}$
is non-vanishing. That is,
\begin{align}\label{thmconsistentvectorfieldsofsubalgebradistributionisafreemodule0p11eq31}
&\Existsis{\V{}}{\topologyofspace{\R^n}}z\in\V{}\subseteq\U{},\cr
&\Foreach{\x}{\V{}}\func{\Det{r}}{\func{\tilde{\alpha}_0}{\x}}\neq 0.
\end{align}
Thus, for every $\x$ in $\V{}$, the matrix $\func{\tilde{\alpha}_0}{\x}$ is invertible.
Since $\tilde{\alpha}_0$ is smooth (\Ref{thmconsistentvectorfieldsofsubalgebradistributionisafreemodule0p11eq27})
and $\V{}$ is an open subset of its domain $\func{\image{\phi}}{\U{}}$, clearly the restriction of $\tilde{\alpha}_0$
to $\V{}$ is also smooth. That is,
\begin{equation}\label{thmconsistentvectorfieldsofsubalgebradistributionisafreemodule0p11eq32}
\func{\res{\tilde{\alpha}_0}}{\V{}}\in\banachmapdifclass{\infty}{\R^n}{\BMat{\R}{n}{r}}{\V{}}{\Mat{\R}{r}{r}}.
\end{equation}
Moreover, according to \Ref{thmconsistentvectorfieldsofsubalgebradistributionisafreemodule0p11eq31}, the value
the $\R$-matrix-valued mapping $\func{\res{\tilde{\alpha}_0}}{\V{}}$ is an invertible matrix all over its domain,
and therefore it is permissible to define the mapping $\omega$ as,
\begin{align}\label{thmconsistentvectorfieldsofsubalgebradistributionisafreemodule0p11eq33}
&\tilde{\omega}\indef\Func{\V{}}{\Mat{\R}{r}{r}},\cr
&\Foreach{\x}{\V{}}\func{\tilde{\omega}}{\x}\eqdef\(\func{\tilde{\alpha}_0}{\x}\)^{-1}.
\end{align}
It is a trivial consequence of \Ref{thmconsistentvectorfieldsofsubalgebradistributionisafreemodule0p11eq32} and
\Ref{thmconsistentvectorfieldsofsubalgebradistributionisafreemodule0p11eq33} that $\omega$ is a smooth map from $\V{}$
to the Banach-space of all square $\R$-matrices of degree $r$. That is,
\begin{equation}\label{thmconsistentvectorfieldsofsubalgebradistributionisafreemodule0p11eq34}
\tilde{\omega}\in\banachmapdifclass{\infty}{\R^n}{\BMat{\R}{r}{r}}{\V{}}{\Mat{\R}{r}{r}}.
\end{equation}
In addition, \Ref{thmconsistentvectorfieldsofsubalgebradistributionisafreemodule0p11eq25} and
\Ref{thmconsistentvectorfieldsofsubalgebradistributionisafreemodule0p11eq33} imply that,
\begin{align}\label{thmconsistentvectorfieldsofsubalgebradistributionisafreemodule0p11eq35}
\Foreach{\x}{\V{}}
\func{\tilde{\gamma}}{\x}&=\(\func{\tilde{\alpha}_0}{\x}\)^{-1}\cdot\func{\tilde{\beta}_0}{\x}\cr
&=\func{\tilde{\omega}}{\x}\cdot\func{\tilde{\beta}_0}{\x}.
\end{align}
Thus by defining the $\R$-matrix-valued mappings $\tilde{\gamma}_0$ and $\tilde{\beta}_1$ as the restrictions of
$\tilde{\gamma}$ and $\tilde{\beta}_0$ to $\V{}$, that is,
\begin{align}
\tilde{\gamma}_0&:=\func{\res{\tilde{\gamma}}}{\V{}},
\label{thmconsistentvectorfieldsofsubalgebradistributionisafreemodule0p11eq36}\\
\tilde{\beta}_1&:=\func{\res{\tilde{\beta}_0}}{\V{}},
\label{thmconsistentvectorfieldsofsubalgebradistributionisafreemodule0p11eq37}
\end{align}
clearly,
\begin{equation}\label{thmconsistentvectorfieldsofsubalgebradistributionisafreemodule0p11eq38}
\Foreach{\x}{\V{}}
\func{\tilde{\gamma}_0}{\x}=\func{\tilde{\omega}}{\x}\cdot\func{\tilde{\beta}_1}{\x}.
\end{equation}
Furthermore, since $\tilde{\beta}_0$ is smooth )\Ref{thmconsistentvectorfieldsofsubalgebradistributionisafreemodule0p11eq25}),
and $\V{}$ is an open subset of its domain $\func{\image{\phi}}{\U{}}$, its restriction to $\V{}$, that is $\tilde{\beta}_1$,
is also smooth.
\begin{equation}\label{thmconsistentvectorfieldsofsubalgebradistributionisafreemodule0p11eq39}
\tilde{\beta}_1\in\banachmapdifclass{\infty}{\R^n}{\BMat{\R}{r}{1}}{\V{}}{\Mat{\R}{r}{1}}.
\end{equation}
Since both $\tilde{\omega}$ and $\tilde{\beta}_1$ are smooth, \Ref{thmconsistentvectorfieldsofsubalgebradistributionisafreemodule0p11eq37}
implies that $\tilde{\gamma}_0$ is also smooth.
\begin{equation}\label{thmconsistentvectorfieldsofsubalgebradistributionisafreemodule0p11eq40}
\tilde{\gamma}_0\in\banachmapdifclass{\infty}{\R^n}{\BMat{\R}{r}{1}}{\V{}}{\Mat{\R}{r}{1}}.
\end{equation}
Therefore the entry mappings of $\tilde{\gamma}_0$ are also smooth, which are obtained by restricting each $\cmp{\zeta_k}{\finv{\phi}}$
to $\V{}$. Precisely, by defining,
\begin{equation}\label{thmconsistentvectorfieldsofsubalgebradistributionisafreemodule0p11eq41}
\psi:=\func{\res{\phi}}{\func{\pimage{\phi}}{\V{}}},
\end{equation}
the smoothness of $\tilde{\gamma}_0$ implies the smoothness of the mappings $\cmp{\zeta_k}{\finv{\psi}}$ as its entry mappings.
That is,
\begin{equation}\label{thmconsistentvectorfieldsofsubalgebradistributionisafreemodule0p11eq42}
\Foreach{k}{\seta{\suc{1}{r}}}
\cmp{\zeta_k}{\finv{\psi}}\in\banachmapdifclass{\infty}{\R^n}{\R}{\V{}}{\R}.
\end{equation}
Since $\phi$ is a chart of $\Lieman{\Liegroup{}}$, \Ref{thmconsistentvectorfieldsofsubalgebradistributionisafreemodule0p11eq31}
and \Ref{thmconsistentvectorfieldsofsubalgebradistributionisafreemodule0p11eq41} imply that,
$\psi$ is also a chart of $\Lieman{\Liegroup{}}$ around $\point$. That is,
\begin{equation}\label{thmconsistentvectorfieldsofsubalgebradistributionisafreemodule0p11eq43}
\psi\in\defset{\p{\phi}}{\maxatlas{}}{\point\in\domain{\p{\phi}}}.
\end{equation}
\endp
\end{itemize}
Therefore, for every point $\point$ of $\Liegroup{}$ there exists a chart of $\Lieman{\Liegroup{}}$ around $\point$ such that
$\cmp{\zeta_k}{\finv{\psi}}$ is smooth for each $k$ in $\seta{\suc{1}{r}}$.
\begin{align}
&\Foreach{k}{\seta{\suc{1}{r}}}\cr
&\Foreach{\point}{\G{}}
\Exists{\psi}{\defset{\p{\phi}}{\maxatlas{}}{\point\in\domain{\p{\phi}}}}
\cmp{\zeta_k}{\finv{\psi}}\in\banachmapdifclass{\infty}{\R^n}{\R}{\V{}}{\R}.
\end{align}
Trivially, this is equivalent to the fact that each $\zeta_k$ is a smooth map from the manifold $\Lieman{\Liegroup{}}$
to the manifold $\RR$. That is,
\begin{equation}
\Foreach{k}{\seta{\suc{1}{r}}}\zeta_k\in\mapdifclass{\infty}{\Lieman{\Liegroup{}}}{\RR},
\end{equation}
and hence according to \Ref{thmconsistentvectorfieldsofsubalgebradistributionisafreemodule0p1eq3}, and
\refdef{defsrealvaluedmoothmapsproduct},
\begin{equation}
\avecff{}\in\defSet{\sum_{k=1}^{r}\cf_k\smoothfprod{\Lieman{\Liegroup{}}}\avecf{k}}
{\mtuple{\cf_1}{\cf_r}\in{\mapdifclass{\infty}{\Lieman{\Liegroup{}}}{\RR}}^{\times r}}.
\end{equation}
\endp
\end{itemize}
Therefore,
\begin{equation}\label{thmconsistentvectorfieldsofsubalgebradistributionisafreemodule0peq2}
\distvecf{\Lieman{\Liegroup{}}}{\func{\subliealgdist{\Liegroup{}}{r}}{\LS{}}}{r}\subseteq
\defSet{\sum_{k=1}^{r}\cf_k\smoothvfprod{\Lieman{\Liegroup{}}}\avecf{k}}
{\mtuple{\cf_1}{\cf_r}\in{\mapdifclass{\infty}{\Lieman{\Liegroup{}}}{\RR}}^{\times r}}.
\end{equation}
\Ref{thmconsistentvectorfieldsofsubalgebradistributionisafreemodule0peq1} and
\Ref{thmconsistentvectorfieldsofsubalgebradistributionisafreemodule0peq2} imply
\Ref{thmconsistentvectorfieldsofsubalgebradistributionisafreemodule0eq2}.\\
The fact that $\seta{\suc{\avecf{1}}{\avecf{r}}}$ forms a linearly-independent set (in the sense of modules)
of the $\smoothring{\infty}{\Lieman{\Liegroup{}}}$-module
$\subspace{\smoothvfmodule{\infty}{\Lieman{\Liegroup{}}}}{\distvecf{\Lieman{\Liegroup{}}}
{\func{\subliealgdist{\Liegroup{}}{r}}{\LS{}}}{r}}$, can be verified trivially by a point-wise analysis of the
problem based on \refdef{defsrealvaluedmoothmapsproduct} and considering that
$\seta{\suc{\func{\avecf{1}}{\point}}{\func{\avecf{r}}{\point}}}$ is
a linearly-independent set of the vector-space $\Tanspace{\point}{\Lieman{\Liegroup{}}}$ at every point
$\point$ of $\Liegroup{}$.
\endthm
\textcolor{Blue}{
\corollary
The $\smoothring{\infty}{\Lieman{\Liegroup{}}}$-module $\smoothvfmodule{\infty}{\Lieman{\Liegroup{}}}$
(the module of smooth vector-fields of the underlying manifold of the smooth group $\Liegroup{}$) is a free module.
}
\proof
Fix $\function{\avecf{}}{\seta{\suc{1}{n}}}{\Leftinvvf{\Liegroup{}}}$ as an an ordered-basis of the canonical Lie-algebra
of the smooth group $\Liegroup{}$, $\LiegroupLiealgebra{\Liegroup{}}$.\\
According to \refdef{defdistributioninducedbyliesubalgebraofliegroupliealgebra},
it is straightforward to verify that the smooth distribution $\func{\subliealgdist{\Liegroup{}}{n}}{\LiegroupLiealgebra{\Liegroup{}}}$
assigns to every point of $\Liegroup{}$ the full set of vectors of its corresponded tangent-space. That is,
\begin{equation}
\Foreach{\g{}}{\G{}}
\func{\[\func{\subliealgdist{\Liegroup{}}{n}}{\LiegroupLiealgebra{\Liegroup{}}}\]}{\g{}}=
\tanspace{\g{}}{\Lieman{\Liegroup{}}}.
\end{equation}
Therefore the set of all smooth vector-fields of $\Lieman{\Liegroup{}}$ consistent with this smooth distribution trivially
coincides with the set of all smooth vector-fields of $\Lieman{\Liegroup{}}$, according to
\refdef{defconsistentvectorfieldswithdistributions}. That is,
\begin{equation}
\distvecf{\Lieman{\Liegroup{}}}{\func{\subliealgdist{\Liegroup{}}{n}}{\LiegroupLiealgebra{\Liegroup{}}}}{n}=
\vecf{\Lieman{\Liegroup{}}}{\infty}.
\end{equation}
Therefore according to \refthm{thmconsistentvectorfieldsofsubalgebradistributionisafreemodule0},
$\seta{\suc{\avecf{1}}{\avecf{n}}}$ is a
module-basis of the $\smoothring{\infty}{\Lieman{\Liegroup{}}}$-module
$\smoothvfmodule{\infty}{\Lieman{\Liegroup{}}}$.
\begin{equation}
\seta{\suc{\avecf{1}}{\avecf{n}}}\in\Modulebases{
\smoothvfmodule{\infty}{\Lieman{\Liegroup{}}}}.
\end{equation}
Therefore, $\smoothvfmodule{\infty}{\Lieman{\Liegroup{}}}$ is a free module.
\endcor
\theorem\label{thmdisubalgebradistributionisinvolutive}
$r$ is taken as an element of $\Zp$, $\LS{}$ as an element of $\subliedim{r}{\LiegroupLiealgebra{\Liegroup{}}}$
(an $r$-dimensional Lie-subalgebra of
the canonical Lie-algebra of the smooth group $\Liegroup{}$, that is $\LiegroupLiealgebra{\Liegroup{}}$).
$\func{\subliealgdist{\Liegroup{}}{r}}{\LS{}}$ is an involutive smooth distribution of $\Lieman{\Liegroup{}}$.
\begin{equation}
\func{\subliealgdist{\Liegroup{}}{r}}{\LS{}}\in\involutivedist{r}{\Lieman{\Liegroup{}}}.
\end{equation}
\proof
$\function{\avecf{}}{\seta{\suc{1}{r}}}{\Leftinvvf{\Liegroup{}}}$ is chosen as an ordered-basis of the Lie-subalgebra
$\LS{}$ of $\LiegroupLiealgebra{\Liegroup{}}$
(precisely, an ordered-basis of the vector-space $\subspace{\VLeftinvvf{\Liegroup{}}}{\LS{}}$).
\begin{itemize}
\item[$\pr{1}$]
Each $\avecff{1}$ and $\avecff{2}$ is taken as an arbitrary element of
$\distvecf{\Lieman{\Liegroup{}}}{\func{\subliealgdist{\Liegroup{}}{r}}{\LS{}}}{r}$, that is a smooth vector-field on
$\Lieman{\Liegroup{}}$ consistent with the smooth distribution $\func{\subliealgdist{\Liegroup{}}{r}}{\LS{}}$ of $\Lieman{\Liegroup{}}$.
Based on \refthm{thmconsistentvectorfieldsofsubalgebradistributionisafreemodule0},
\begin{align}
&\Existsis{\mtuple{\cf_1}{\cf_r}}{{\mapdifclass{\infty}{\Lieman{\Liegroup{}}}{\RR}}^{\times 2}}
\avecff{1}=\sum_{k=1}^{r}\cf_{k}\smoothvfprod{\Lieman{\Liegroup{}}}\avecf{k},
\label{thmdisubalgebradistributionisinvolutivep1eq1}\\
&\Existsis{\mtuple{\cg_1}{\cg_r}}{{\mapdifclass{\infty}{\Lieman{\Liegroup{}}}{\RR}}^{\times 2}}
\avecff{2}=\sum_{k=1}^{r}\cg_{k}\smoothvfprod{\Lieman{\Liegroup{}}}\avecf{k}.
\label{thmdisubalgebradistributionisinvolutivep1eq2}
\end{align}
According to \refthm{thmliealgebraofvectorfields0}, \refthm{thmliebracketonthemoduleofvectorfields1},
and considering the bilinearity of the Lie-operation of a Lie-algebra, it is evident that,
\begin{align}\label{thmdisubalgebradistributionisinvolutivep1eq3}
\liebracket{\avecff{1}}{\avecff{2}}{\Lieman{\Liegroup{}}}
&=\liebracket{\sum_{k=1}^{r}\cf_{k}\smoothvfprod{\Lieman{\Liegroup{}}}\avecf{k}}
{\sum_{k=1}^{r}\cg_{k}\smoothvfprod{\Lieman{\Liegroup{}}}\avecf{k}}{\Lieman{\Liegroup{}}}
=\sum_{j=1}^{r}\sum_{k=1}^{r}\liebracket{\cf_{j}\smoothvfprod{\Lieman{\Liegroup{}}}\avecf{j}}
{\cg_{k}\smoothvfprod{\Lieman{\Liegroup{}}}\avecf{k}}{\Lieman{\Liegroup{}}}\cr
&=\hskip0.5\baselineskip\sum_{j=1}^{r}\sum_{k=1}^{r}\left\{\(\cf_j
\smoothfprod{\Lieman{\Liegroup{}}}\cg_k\)\smoothvfprod{\Lieman{\Liegroup{}}}\liebracket{\avecf{j}}{\avecf{k}}{\Lieman{\Liegroup{}}}
\smoothvfsum{\Lieman{\Liegroup{}}}
\[\cf_j\smoothfprod{\Lieman{\Liegroup{}}}
\bigg(\func{\[\func{\Lieder{\Lieman{\Liegroup{}}}}{\avecf{j}}\]}{\cg_k}\bigg)\]\smoothvfprod{\Lieman{\Liegroup{}}}\avecf{k}\right.\cr
&\hskip0.5\baselineskip\left.\smoothvfsub{\Lieman{\Liegroup{}}}
\[\cg_k\smoothfprod{\Lieman{\Liegroup{}}}
\bigg(\func{\[\func{\Lieder{\Lieman{\Liegroup{}}}}{\avecf{k}}\]}{\cf_j}\bigg)\]\smoothvfprod{\Lieman{\Liegroup{}}}\avecf{j}\right\}.
\end{align}
According to \refdef{defconsistentvectorfieldswithdistributions} and
\refdef{defdistributioninducedbyliesubalgebraofliegroupliealgebra}, it is evident that,
\begin{equation}\label{thmdisubalgebradistributionisinvolutivep1eq4}
\LS{}\subseteq\distvecf{\Lieman{\Liegroup{}}}{\func{\subliealgdist{\Liegroup{}}{r}}{\LS{}}}{r},
\end{equation}
and since $\LS{}$ is a Lie-subalgebra of $\LiegroupLiealgebra{\Liegroup{}}$, $\LS{}$ is closed under the Lie-bracket operation
on $\vecf{\Lieman{\Liegroup{}}}{\infty}$, and thus,
\begin{equation}\label{thmdisubalgebradistributionisinvolutivep1eq5}
\Foreach{\opair{\avecf{}}{\p{\avecf{}}}}{\Cprod{\LS{}}{\LS{}}}
\liebracket{\avecf{}}{\p{\avecf{}}}{\Lieman{\Liegroup{}}}\in
\distvecf{\Lieman{\Liegroup{}}}{\func{\subliealgdist{\Liegroup{}}{r}}{\LS{}}}{r}.
\end{equation}
Therefore, since $\seta{\suc{\avecf{1}}{\avecf{r}}}\subseteq\LS{}$, it is obvious that,
\begin{align}
&\seta{\suc{\avecf{1}}{\avecf{r}}}\subseteq
\distvecf{\Lieman{\Liegroup{}}}{\func{\subliealgdist{\Liegroup{}}{r}}{\LS{}}}{r},
\label{thmdisubalgebradistributionisinvolutivep1eq6}\\
&\Foreach{\opair{j}{k}}{{\seta{\suc{1}{r}}}^{\times 2}}
\liebracket{\avecf{j}}{\avecf{k}}{\Lieman{\Liegroup{}}}\in
\distvecf{\Lieman{\Liegroup{}}}{\func{\subliealgdist{\Liegroup{}}{r}}{\LS{}}}{r}.
\label{thmdisubalgebradistributionisinvolutivep1eq7}
\end{align}
Therefore, since according to \refthm{thmsubmoduleinducedbydistribution},
$\distvecf{\Lieman{\Liegroup{}}}{\func{\subliealgdist{\Liegroup{}}{r}}{\LS{}}}{r}$ is a submodule of the
$\smoothring{\infty}{\Lieman{\Liegroup{}}}$-module $\smoothvfmodule{\infty}{\Lieman{\Liegroup{}}}$,
based on \refdef{defvecfmodule}, \Ref{thmdisubalgebradistributionisinvolutivep1eq3} implies that,
\begin{equation}
\liebracket{\avecff{1}}{\avecff{2}}{\Lieman{\Liegroup{}}}\in
\distvecf{\Lieman{\Liegroup{}}}{\func{\subliealgdist{\Liegroup{}}{r}}{\LS{}}}{r}.
\end{equation}
\endp
\end{itemize}
Therefore,
\begin{equation}
\Foreach{\opair{\avecff{1}}{\avecff{2}}}
{{\distvecf{\Lieman{\Liegroup{}}}{\func{\subliealgdist{\Liegroup{}}{r}}{\LS{}}}{r}}^{\times 2}}
\liebracket{\avecff{1}}{\avecff{2}}{\Lieman{\Liegroup{}}}\in
\distvecf{\Lieman{\Liegroup{}}}{\func{\subliealgdist{\Liegroup{}}{r}}{\LS{}}}{r},
\end{equation}
which, according to \refdef{definvolutivedistribution} means that
$\distvecf{\Lieman{\Liegroup{}}}{\func{\subliealgdist{\Liegroup{}}{r}}{\LS{}}}{r}$ is an $r$-dimensional
involutive smooth distribution of $\Lieman{\Liegroup{}}$.
\endthm
\corollary\label{thmdistributioninducedbyliesubalgebraofliegroupliealgebraisintegrable}
$r$ is taken as an element of $\Zp$, $\LS{}$ as an element of $\subliedim{r}{\LiegroupLiealgebra{\Liegroup{}}}$
(an $r$-dimensional Lie-subalgebra of
the canonical Lie-algebra of the smooth group $\Liegroup{}$, that is $\LiegroupLiealgebra{\Liegroup{}}$).
$\func{\subliealgdist{\Liegroup{}}{r}}{\LS{}}$ is an integrable smooth distribution of $\Lieman{\Liegroup{}}$.
\begin{equation}
\func{\subliealgdist{\Liegroup{}}{r}}{\LS{}}\in\integrabledist{r}{\Lieman{\Liegroup{}}}.
\end{equation}
\proof
According to \refthm{thmdisubalgebradistributionisinvolutive}, and the theorem of Frobenius
(\refthm{thmFrobenius}), it is obvious.
\endcor
\theorem\label{thmlefttranslationspreservesubalgebradistributions}
$r$ is taken as an element of $\Zp$, $\LS{}$ as an element of $\subliedim{r}{\LiegroupLiealgebra{\Liegroup{}}}$
(an $r$-dimensional Lie-subalgebra of
the canonical Lie-algebra of the smooth group $\Liegroup{}$, that is $\LiegroupLiealgebra{\Liegroup{}}$).
For every $\g{}$ in $\G{}$, the left-translation of $\LieG{\Liegroup{}}$ by $\g{}$, that is
$\gltrans{\LieG{\Liegroup{}}}{\g{}}$, is a distribution-preserving $\infty$-diffeomorphism from the
$\difclass{\infty}$ distributed-space $\opair{\Lieman{\Liegroup{}}}{\func{\subliealgdist{\Liegroup{}}{r}}{\LS{}}}$
to itself. That is,
\begin{align}
\Foreach{\g{}}{\G{}}
\Foreach{\point}{\G{}}
\func{\image{\[\der{\gltrans{\LieG{\Liegroup{}}}{\g{}}}{\Lieman{\Liegroup{}}}
{\Lieman{\Liegroup{}}}\]}}{\func{\[\func{\subliealgdist{\Liegroup{}}{r}}{\LS{}}\]}{\point}}&=
\func{\[\func{\subliealgdist{\Liegroup{}}{r}}{\LS{}}\]}{\func{\gltrans{\LieG{\Liegroup{}}}{\g{}}}{\point}}\cr
&=\func{\[\func{\subliealgdist{\Liegroup{}}{r}}{\LS{}}\]}{\g{}\gop{}\point}.
\end{align}
\proof
Each $\g{}$ and $\point$ is taken as an arbitrary element of $\G{}$.
According to \refthm{thmtranslationsarediffeomorphisms}, $\gltrans{\LieG{\Liegroup{}}}{\g{}}$ is an $\infty$-diffeomorphism from
$\Lieman{\Liegroup{}}$ to itself. Moreover,
According to \refdef{defleftinvariantvectorfields} and \refdef{defdistributioninducedbyliesubalgebraofliegroupliealgebra},
\begin{align}
\func{\image{\[\der{\gltrans{\LieG{\Liegroup{}}}{\g{}}}{\Lieman{\Liegroup{}}}
{\Lieman{\Liegroup{}}}\]}}{\func{\[\func{\subliealgdist{\Liegroup{}}{r}}{\LS{}}\]}{\point}}&=
\func{\image{\[\der{\gltrans{\LieG{\Liegroup{}}}{\g{}}}{\Lieman{\Liegroup{}}}
{\Lieman{\Liegroup{}}}\]}}{\defSet{\func{\avecf{}}{\point}}{\avecf{}\in\LS{}}}\cr
&=\defSet{\func{\(\cmp{\[\der{\gltrans{\LieG{\Liegroup{}}}{\g{}}}{\Lieman{\Liegroup{}}}
{\Lieman{\Liegroup{}}}\]}{\avecf{}}\)}{\point}}{\avecf{}\in\LS{}}\cr
&=\defSet{\func{\(\cmp{\avecf{}}{\gltrans{\LieG{\Liegroup{}}}{\g{}}}\)}{\point}}{\avecf{}\in\LS{}}\cr
&=\func{\[\func{\subliealgdist{\Liegroup{}}{r}}{\LS{}}\]}{\func{\gltrans{\LieG{\Liegroup{}}}{\g{}}}{\point}}.
\end{align}
Therefore according to \refdef{defdistributionpreservingmap}, $\gltrans{\LieG{\Liegroup{}}}{\g{}}$ is
a distribution-preserving $\infty$-diffeomorphism from the $\difclass{\infty}$ distributed-space
$\opair{\Lieman{\Liegroup{}}}{\func{\subliealgdist{\Liegroup{}}{r}}{\LS{}}}$
to itself.
\endthm
\theorem\label{thmlefttranslationspermutetheleavesofthesolutionofsubalgebradistribution}
$r$ is taken as an element of $\Zp$, $\LS{}$ as an element of $\subliedim{r}{\LiegroupLiealgebra{\Liegroup{}}}$
(an $r$-dimensional Lie-subalgebra of
the canonical Lie-algebra of the smooth group $\Liegroup{}$, that is $\LiegroupLiealgebra{\Liegroup{}}$).
For every $\g{}$ in $\G{}$, the left-translation of $\LieG{\Liegroup{}}$ by $\g{}$ permutes the
leaves of the smooth foliation $\intdistsolution{\Lieman{\Liegroup{}}}{r}\func{\subliealgdist{\Liegroup{}}{r}}{\LS{}}$,
that is the unique solution of the smooth distribution $\func{\subliealgdist{\Liegroup{}}{r}}{\LS{}}$, of
$\Lieman{\Liegroup{}}$.
That is,
\begin{align}
&\Foreach{\g{}}{\G{}}\cr
&\EqClass{\G{}}{\(\intdistsolution{\Lieman{\Liegroup{}}}{r}\func{\subliealgdist{\Liegroup{}}{r}}{\LS{}}\)}=
\defSet{\func{\image{\[\gltrans{\LieG{\Liegroup{}}}{\g{}}\]}}{\leaf{}}}{\leaf{}\in\EqClass{\G{}}
{\(\intdistsolution{\Lieman{\Liegroup{}}}{r}\func{\subliealgdist{\Liegroup{}}{r}}{\LS{}}\)}}.
\end{align}
\proof
According to
\refthm{thmlefttranslationspreservesubalgebradistributions} and
\refthm{thmdistpreservingdiffinterchangestheleavesofsolutions}, it is trivial.
\endthm
\theorem\label{thmidentityleafofthesolutionofsubalgebradistributionisanimmersedsubgroup}
$r$ is taken as an element of $\Zp$, $\LS{}$ as an element of $\subliedim{r}{\LiegroupLiealgebra{\Liegroup{}}}$
(an $r$-dimensional Lie-subalgebra of
the canonical Lie-algebra of the smooth group $\Liegroup{}$, that is $\LiegroupLiealgebra{\Liegroup{}}$).
The leaf of the smooth foliation $\intdistsolution{\Lieman{\Liegroup{}}}{r}\func{\subliealgdist{\Liegroup{}}{r}}{\LS{}}$
containing the identity element of the group $\LieG{\Liegroup{}}$, that is $\IG{}$, (when endowed with its canonical maximal-atlas induced by
the foliation $\intdistsolution{\Lieman{\Liegroup{}}}{r}\func{\subliealgdist{\Liegroup{}}{r}}{\LS{}}$)
is an $r$-dimensional immersed smooth subgroup of the Lie group $\Liegroup{}$. Precisely, the manifold
$\func{\leafman{\Lieman{\Liegroup{}}}{\intdistsolution{\Lieman{\Liegroup{}}}{r}\func{\subliealgdist{\Liegroup{}}{r}}{\LS{}}}{r}}
{\PEqclass{\intdistsolution{\Lieman{\Liegroup{}}}{r}\func{\subliealgdist{\Liegroup{}}{r}}{\LS{}}}{\IG{}}}$, as defined in
\refdef{defleaves}, is an $r$-dimensional immersed smooth subgroup of $\Liegroup{}$, when endowed with the group operation $\gop{}$
of $\Liegroup{}$ restricted to $\PEqclass{\intdistsolution{\Lieman{\Liegroup{}}}{r}\func{\subliealgdist{\Liegroup{}}{r}}{\LS{}}}{\IG{}}$.
By a subtle abuse of notation,
\begin{equation}
\opair{\func{\leafman{\Lieman{\Liegroup{}}}{\intdistsolution{\Lieman{\Liegroup{}}}{r}\func{\subliealgdist{\Liegroup{}}{r}}{\LS{}}}{r}}
{\PEqclass{\intdistsolution{\Lieman{\Liegroup{}}}{r}\func{\subliealgdist{\Liegroup{}}{r}}{\LS{}}}{\IG{}}}}{\gop{}}
\in\imsubgroup{\Liegroup{}}.
\end{equation}
\proof
To ease the notation,
\begin{align}
&\pEqclass{\IG{}}{}:=\PEqclass{\intdistsolution{\Lieman{\Liegroup{}}}{r}\func{\subliealgdist{\Liegroup{}}{r}}{\LS{}}}{\IG{}},\\
&\Leafman{\pEqclass{\IG{}}{}}:=
\func{\leafman{\Lieman{\Liegroup{}}}{\intdistsolution{\Lieman{\Liegroup{}}}{r}\func{\subliealgdist{\Liegroup{}}{r}}{\LS{}}}{r}}
{\PEqclass{\intdistsolution{\Lieman{\Liegroup{}}}{r}\func{\subliealgdist{\Liegroup{}}{r}}{\LS{}}}{\IG{}}}.
\end{align}
According to \refthm{thmleavesareimmersedsubmanifolds},
$\Leafman{\pEqclass{\IG{}}{}}$
is an $r$-dimensional immersed submanifold of $\Lieman{\Liegroup{}}$.
\begin{itemize}
\item[$\pr{1}$]
$\g{}$ is taken as an arbitrary element of $\pEqclass{\IG{}}{}$. Clearly,
\begin{align}
\func{\[\gltrans{\LieG{\Liegroup{}}}{\invG{\g{}}{}}\]}{\g{}}&=
\invG{\g{}}{}\gop{}\g{}\cr
&=\IG{}\in\pEqclass{\IG{}}{},
\end{align}
and since based on \refthm{thmlefttranslationspermutetheleavesofthesolutionofsubalgebradistribution},
$\gltrans{\LieG{\Liegroup{}}}{\invG{\g{}}{}}$ permutes the
leaves of the smooth foliation $\intdistsolution{\Lieman{\Liegroup{}}}{r}\func{\subliealgdist{\Liegroup{}}{r}}{\LS{}}$, evidently,
\begin{equation}
\func{\image{\[\gltrans{\LieG{\Liegroup{}}}{\invG{\g{}}{}}\]}}{\pEqclass{\IG{}}{}}=\pEqclass{\IG{}}{}.
\end{equation}
This trivially implies that,
\begin{equation}
\Foreach{\hh}{\pEqclass{\IG{}}{}}
\invG{\g{}}{}\gop{}\hh\in\pEqclass{\IG{}}{},
\end{equation}
\endp
\end{itemize}
Therefore,
\begin{equation}
\Foreach{\opair{\g{}}{\hh}}{\Cprod{\pEqclass{\IG{}}{}}{\pEqclass{\IG{}}{}}}
\invG{\g{}}{}\gop{}\hh\in\pEqclass{\IG{}}{},
\end{equation}
and hence considering that $\pEqclass{\IG{}}{}$ is non-empty (because it contains $\IG{}$),
$\pEqclass{\IG{}}{}$ is a subgroup of the group $\LieG{\Liegroup{}}$.\\
Therefore, according to \reflem{lemimmersedsubmanifoldsubgroup},
since $\Leafman{\pEqclass{\IG{}}{}}$ is an $\infty$-immersed submanifold of $\Lieman{\Liegroup{}}$,
and the set of all points $\pEqclass{\IG{}}{}$ of this submanifold
is a subgroup of $\LieG{\Liegroup{}}$, $\Leafman{\pEqclass{\IG{}}{}}$ is an immersed smooth subgroup of $\Liegroup{}$.
The dimension of this submanifold is trivially $r$, since it is a leaf of an $r$-dimensional smooth foliation of
$\Lieman{\Liegroup{}}$.
\endthm
\definition\label{defidentityleafofthefoliationoftheinduceddistributionofasubalgebra}
$r$ is taken as an element of $\Zp$, $\LS{}$ as an element of $\subliedim{r}{\LiegroupLiealgebra{\Liegroup{}}}$
(an $r$-dimensional Lie-subalgebra of
the canonical Lie-algebra of the smooth group $\Liegroup{}$, that is $\LiegroupLiealgebra{\Liegroup{}}$).
The immersed smooth subgroup $\func{\liesublie{\Liegroup{}}{r}}{\LS{}}$ of $\Liegroup{}$ is defined as,
\begin{equation}
\func{\liesublie{\Liegroup{}}{r}}{\LS{}}:=
\opair{\func{\leafman{\Lieman{\Liegroup{}}}{\intdistsolution{\Lieman{\Liegroup{}}}{r}\func{\subliealgdist{\Liegroup{}}{r}}{\LS{}}}{r}}
{\PEqclass{\intdistsolution{\Lieman{\Liegroup{}}}{r}\func{\subliealgdist{\Liegroup{}}{r}}{\LS{}}}{\IG{}}}}{\gop{}}.
\end{equation}
Also,
\begin{equation}
\func{\liesublieset{\Liegroup{}}{r}}{\LS{}}:=
\PEqclass{\intdistsolution{\Lieman{\Liegroup{}}}{r}\func{\subliealgdist{\Liegroup{}}{r}}{\LS{}}}{\IG{}}.
\end{equation}
\endef
\theorem\label{thmliesubalgebraofliealgebraofliegroupcorrespondstoaliesubgroup}
$r$ is taken as an element of $\Zp$, $\LS{}$ as an element of $\subliedim{r}{\LiegroupLiealgebra{\Liegroup{}}}$
(an $r$-dimensional Lie-subalgebra of
the canonical Lie-algebra of the smooth group $\Liegroup{}$, that is $\LiegroupLiealgebra{\Liegroup{}}$).
The canonical Lie-algebra of the immersed subgroup
$\func{\liesublie{\Liegroup{}}{r}}{\LS{}}$ of
$\Liegroup{}$ is Lie-algebraically isomorphic to $\subspace{\LiegroupLiealgebra{\Liegroup{}}}{\LS{}}$. That is,
\begin{equation}
\liealgisomorphic{\LiegroupLiealgebra{\func{\liesublie{\Liegroup{}}{r}}{\LS{}}}}
{\subspace{\LiegroupLiealgebra{\Liegroup{}}}{\LS{}}},
\end{equation}
where,
\begin{equation}
\funcimage{\func{\indliemor{\func{\liesublie{\Liegroup{}}{r}}{\LS{}}}
{\Liegroup{}}}{\Injection{\func{\liesublieset{\Liegroup{}}{r}}{\LS{}}}{\G{}}}}=\LS{}.
\end{equation}
Moreover, the underlying topological-space of the smooth group $\func{\liesublie{\Liegroup{}}{r}}{\LS{}}$ is connected.
\proof
According to \refthm{thmidentityleafofthesolutionofsubalgebradistributionisanimmersedsubgroup} and
\refthm{thmsmoothsubgroupliealgebra},
the canonical Lie-algebra of the smooth group $\func{\liesublie{\Liegroup{}}{r}}{\LS{}}$, that is
$\LiegroupLiealgebra{\func{\liesublie{\Liegroup{}}{r}}{\LS{}}}$, is
Lie-algebraically isomorphic to the Lie-subalgebra
$\funcimage{\func{\indliemor{\func{\liesublie{\Liegroup{}}{r}}{\LS{}}}
{\Liegroup{}}}{\Injection{\func{\liesublieset{\Liegroup{}}{r}}{\LS{}}}{\G{}}}}$ of the canonical Lie-algebra of $\Liegroup{}$
(as the image of the Lie-algebra-morphism $\func{\indliemor{\func{\liesublie{\Liegroup{}}{r}}{\LS{}}}
{\Liegroup{}}}{\Injection{\func{\liesublieset{\Liegroup{}}{r}}{\LS{}}}{\G{}}}$). That is,
\begin{equation}
\liealgisomorphic{\LiegroupLiealgebra{\func{\liesublie{\Liegroup{}}{r}}{\LS{}}}}
{\subspace{\LiegroupLiealgebra{\Liegroup{}}}{\WW{}}},
\end{equation}
where,
\begin{equation}
\WW{}:=\funcimage{\func{\indliemor{\func{\liesublie{\Liegroup{}}{r}}{\LS{}}}
{\Liegroup{}}}{\Injection{\func{\liesublieset{\Liegroup{}}{r}}{\LS{}}}{\G{}}}}.
\end{equation}
Thus according to \refdef{definducedliealgebramorphismfromliemorphism} and
\refdef{deflinvvftanspacecorrespondence}, \refdef{deftangentbundleoffoliation},
\refdef{defsolutionofanintegrabledistributionasafoliation}, and
\refdef{defdistributioninducedbyliesubalgebraofliegroupliealgebra}, considering that
$\func{\liesublie{\Liegroup{}}{r}}{\LS{}}$ is the solution (as a foliation) of the smooth distribution
$\func{\subliealgdist{\Liegroup{}}{r}}{\LS{}}$, and considering that $\liegvftan{\Liegroup{}}$
is a linear-isomorphism,
\begin{align}
&\hskip0.5\baselineskip\funcimage{\func{\indliemor{\func{\liesublie{\Liegroup{}}{r}}{\LS{}}}
{\Liegroup{}}}{\Injection{\func{\liesublieset{\Liegroup{}}{r}}{\LS{}}}{\G{}}}}\cr
&=\func{\image{\(\cmp{\finv{\liegvftan{\Liegroup{}}}}
{\cmp{\[\der{\Injection{\func{\liesublieset{\Liegroup{}}{r}}{\LS{}}}{\G{}}}
{\Lieman{\func{\liesublie{\Liegroup{}}{r}}{\LS{}}}}{\Lieman{\Liegroup{}}}\]}
{\liegvftan{\func{\liesublie{\Liegroup{}}{r}}{\LS{}}}}}\)}}{\Leftinvvf{\func{\liesublie{\Liegroup{}}{r}}{\LS{}}}}\cr
&=\func{\pimage{\liegvftan{\Liegroup{}}}}{\func{\image{\[\der{\Injection{\func{\liesublieset{\Liegroup{}}{r}}{\LS{}}}{\G{}}}
{\Lieman{\func{\liesublie{\Liegroup{}}{r}}{\LS{}}}}{\Lieman{\Liegroup{}}}\]}}{\tanspace{\IG{}}{\func{\liesublie{\Liegroup{}}{r}}{\LS{}}}}}\cr
&=\func{\pimage{\liegvftan{\Liegroup{}}}}{\func{\[\func{\subliealgdist{\Liegroup{}}{r}}{\LS{}}\]}{\IG{}}}\cr
&=\func{\pimage{\liegvftan{\Liegroup{}}}}{\defSet{\func{\avecf{}}{\IG{}}}{\avecf{}\in\LS{}}}\cr
&=\func{\pimage{\liegvftan{\Liegroup{}}}}{\defSet{\func{\liegvftan{\Liegroup{}}}{\avecf{}}}{\avecf{}\in\LS{}}}\cr
&=\LS{}.
\end{align}
The underlying topological-space of the smooth group $\func{\liesublie{\Liegroup{}}{r}}{\LS{}}$ is connected, because its underlying manifold
$\Lieman{\func{\liesublie{\Liegroup{}}{r}}{\LS{}}}$ is a leaf of a foliation and leaves of any foliation are connected with respect
to their natural topology (which is not necessarily the same as the topology of the base manifold).
\endthm
\theorem\label{thmleavesofthefoliationoftheinduceddistributionofasubalgebra}
$r$ is taken as an element of $\Zp$, $\LS{}$ as an element of $\subliedim{r}{\LiegroupLiealgebra{\Liegroup{}}}$
(an $r$-dimensional Lie-subalgebra of
the canonical Lie-algebra of the smooth group $\Liegroup{}$, that is $\LiegroupLiealgebra{\Liegroup{}}$).
The set of all leaves of the foliation
$\intdistsolution{\Lieman{\Liegroup{}}}{r}\func{\subliealgdist{\Liegroup{}}{r}}{\LS{}}$,
(actually the set of all equivalence classes of the foliation
$\intdistsolution{\Lieman{\Liegroup{}}}{r}\func{\subliealgdist{\Liegroup{}}{r}}{\LS{}}$)
coincides with the set of all
left-cosets (or right-cosets) of the subgroup $\func{\liesublieset{\Liegroup{}}{r}}{\LS{}}$ of $\LieG{\Liegroup{}}$.
\begin{equation}
\EqClass{\G{}}{\[\intdistsolution{\Lieman{\Liegroup{}}}{r}\func{\subliealgdist{\Liegroup{}}{r}}{\LS{}}\]}=
\func{\LCoset{\LieG{\Liegroup{}}}}{\func{\liesublieset{\Liegroup{}}{r}}{\LS{}}}.
\end{equation}
\proof
Each left-coset of $\func{\liesublieset{\Liegroup{}}{r}}{\LS{}}$ is the image of $\func{\liesublieset{\Liegroup{}}{r}}{\LS{}}$
under the left-translation of $\LieG{\Liegroup{}}$ by an element $\g{}$ of $\G{}$. Actually,
\begin{equation}
\func{\LCoset{\LieG{\Liegroup{}}}}{\func{\liesublieset{\Liegroup{}}{r}}{\LS{}}}=
\defSet{\func{\image{\[\gltrans{\LieG{\Liegroup{}}}{\g{}}\]}}{\func{\liesublieset{\Liegroup{}}{r}}{\LS{}}}}
{\g{}\in\G{}}.
\end{equation}
Moreover, since $\func{\liesublieset{\Liegroup{}}{r}}{\LS{}}$ is an equivalence class of the foliation
$\intdistsolution{\Lieman{\Liegroup{}}}{r}\func{\subliealgdist{\Liegroup{}}{r}}{\LS{}}$,
according to \refthm{thmlefttranslationspermutetheleavesofthesolutionofsubalgebradistribution},
\begin{equation}
\defSet{\func{\image{\[\gltrans{\LieG{\Liegroup{}}}{\g{}}\]}}{\func{\liesublieset{\Liegroup{}}{r}}{\LS{}}}}{\g{}\in\G{}}\subseteq
\EqClass{\G{}}{\[\intdistsolution{\Lieman{\Liegroup{}}}{r}\func{\subliealgdist{\Liegroup{}}{r}}{\LS{}}\]},
\end{equation}
and hence,
\begin{equation}
\func{\LCoset{\LieG{\Liegroup{}}}}{\func{\liesublieset{\Liegroup{}}{r}}{\LS{}}}\subseteq
\EqClass{\G{}}{\[\intdistsolution{\Lieman{\Liegroup{}}}{r}\func{\subliealgdist{\Liegroup{}}{r}}{\LS{}}\]}.
\end{equation}
Furthermore, considering that both the set of all left-cosets of $\func{\liesublieset{\Liegroup{}}{r}}{\LS{}}$
and the set of all equivalence classes of the foliation
$\intdistsolution{\Lieman{\Liegroup{}}}{r}\func{\subliealgdist{\Liegroup{}}{r}}{\LS{}}$ are partitions of
$\G{}$, they must be equal.
\endthm
\definition\label{defimmersedsubgroupsandliesubalgebracorrespondence}
The mapping $\liesubgalcor{\Liegroup{}}$ from the set of all connected immersed smooth subgroups of $\Liegroup{}$
to the set of all non-trivial Lie-subalgebras of the canonical Lie-algebra of $\Liegroup{}$ (those with positive dimension),
is defined as,
\begin{align}
&\liesubgalcor{\Liegroup{}}\indef\Func{\connectedimsubgroup{\Liegroup{}}}
{\compl{\sublie{\LiegroupLiealgebra{\Liegroup{}}}}{\subliedim{0}{\LiegroupLiealgebra{\Liegroup{}}}}},\cr
&\Foreach{\aliegroup{}}{\imsubgroup{\Liegroup{}}}
\func{\liesubgalcor{\Liegroup{}}}{\aliegroup{}}\eqdef
\funcimage{\func{\indliemor{\aliegroup{}}
{\Liegroup{}}}{\Injection{\Liepoints{\aliegroup{}}}{\G{}}}}.
\end{align}
\caution
This mapping is well-defined according to \refthm{thminducedliealgebramorphismfromliemorphism1}, and considering that
the injection of an immersed smooth subgroup of a smooth group into that smooth group is a smooth group morphism
(\reflem{leminjectionofsmoothsubgroupisamorphism}), and that
the image of a Lie-algebra-morphism is a Lie-subalgebra of the target Lie-algebra.
\endef
\theorem\label{thmimmersedsubgroupsandliesubalgebracorrespondence}
The mapping $\liesubgalcor{\Liegroup{}}$ is a bijection from the set of all connected immersed smooth subgroups of $\Liegroup{}$
to the set of all non-trivial Lie-subalgebras of the canonical Lie-algebra of $\Liegroup{}$.
\begin{equation}
\liesubgalcor{\Liegroup{}}\in\IF{\connectedimsubgroup{\Liegroup{}}}
{\compl{\sublie{\LiegroupLiealgebra{\Liegroup{}}}}{\subliedim{0}{\LiegroupLiealgebra{\Liegroup{}}}}}.
\end{equation}
\proof
According to \refthm{thmliesubalgebraofliealgebraofliegroupcorrespondstoaliesubgroup}, and
\refdef{defimmersedsubgroupsandliesubalgebracorrespondence}, the mapping
$\liesubgalcor{\Liegroup{}}$ is obviously surjective.
\begin{itemize}
\item[$\pr{1}$]
$\aliegroup{}$ is taken as an arbitrary connected immersed smooth subgroup of $\Liegroup{}$, and,
\begin{align}
\LS{}:&=\func{\liesubgalcor{\Liegroup{}}}{\aliegroup{}},
\label{thmimmersedsubgroupsandliesubalgebracorrespondencep1eq1}\\
r:&=\liedim{\subspace{\LiegroupLiealgebra{\Liegroup{}}}{\LS{}}}.
\label{thmimmersedsubgroupsandliesubalgebracorrespondencep1eq2}
\end{align}
Thus $\LS{}$ is a non-trivial Lie-subalgebra of the canonical Lie-algebra $\LiegroupLiealgebra{\Liegroup{}}$ of $\Liegroup{}$,
and therefore a vector-subspace of the vector-space $\Leftinvvf{\Liegroup{}}$ of the left-invariant vector-fields of $\Liegroup{}$.
According to \refdef{defimmersedsubgroupsandliesubalgebracorrespondence},
\refdef{deflinvvftanspacecorrespondence}, and
\refdef{definducedliealgebramorphismfromliemorphism},
\begin{align}\label{thmimmersedsubgroupsandliesubalgebracorrespondencep1eq2}
\LS{}&=\funcimage{\func{\indliemor{\aliegroup{}}
{\Liegroup{}}}{\Injection{\Liepoints{\aliegroup{}}}{\G{}}}}\cr
&=\funcimage{\cmp{\finv{\liegvftan{\Liegroup{}}}}
{\cmp{\[\der{\Injection{\Liepoints{\aliegroup{}}}{\G{}}}{\Lieman{\aliegroup{}}}{\Lieman{\Liegroup{}}}\]}{\liegvftan{\aliegroup{}}}}},\cr
&=\func{\pimage{\liegvftan{\Liegroup{}}}}
{\func{\image{\[\der{\Injection{\Liepoints{\aliegroup{}}}{\G{}}}{\Lieman{\aliegroup{}}}{\Lieman{\Liegroup{}}}\]}}
{\func{\image{\liegvftan{\aliegroup{}}}}{\LiegroupLiealgebra{\aliegroup{}}}}}.
\end{align}
Therefore,
\begin{align}\label{thmimmersedsubgroupsandliesubalgebracorrespondencep1eq3}
\defSet{\func{\avecff{}}{\IG{}}}{\avecff{}\in\LS{}}=
\func{\image{\liegvftan{\Liegroup{}}}}{\LS{}}&=
\func{\image{\[\der{\Injection{\Liepoints{\aliegroup{}}}{\G{}}}{\Lieman{\aliegroup{}}}
{\Lieman{\Liegroup{}}}\]}}{\func{\image{\liegvftan{\aliegroup{}}}}{\LiegroupLiealgebra{\aliegroup{}}}}\cr
&=\defSet{\func{\[\der{\Injection{\Liepoints{\aliegroup{}}}{\G{}}}{\Lieman{\aliegroup{}}}
{\Lieman{\Liegroup{}}}\]}{\func{\avecf{}}{\IG{}}}}{\avecf{}\in\LiegroupLiealgebra{\aliegroup{}}}.
\end{align}
\begin{itemize}
\item[$\pr{1-1}$]
$\g{}$ is taken as an arbitrary element of $\Liepoints{\aliegroup{}}$.
According to \reflem{lemleftinvariantvectorfieldsequiv0}, and considering that $\LS{}$ is a subset of $\Leftinvvf{\Liegroup{}}$,
\Ref{thmimmersedsubgroupsandliesubalgebracorrespondencep1eq3} together with the chain rule of differentiation implies,
\begin{align}\label{thmimmersedsubgroupsandliesubalgebracorrespondencep11eq1}
\defSet{\func{\avecff{}}{\g{}}}{\avecff{}\in\LS{}}&=
\defSet{\func{\(\der{\gltrans{\LieG{\Liegroup{}}}{\g{}}}
{\Lieman{\Liegroup{}}}{\Lieman{\Liegroup{}}}\)}{\func{\avecff{}}{\IG{}}}}{\avecff{}\in\LS{}}\cr
&=\defSet{\func{\[\cmp{\(\der{\gltrans{\LieG{\Liegroup{}}}{\g{}}}
{\Lieman{\Liegroup{}}}{\Lieman{\Liegroup{}}}\)}{\(\der{\Injection{\Liepoints{\aliegroup{}}}{\G{}}}{\Lieman{\aliegroup{}}}
{\Lieman{\Liegroup{}}}\)}\]}{\func{\avecf{}}{\IG{}}}}{\avecf{}\in\LiegroupLiealgebra{\aliegroup{}}}\cr
&=\defSet{\func{\[\der{\(\cmp{\gltrans{\LieG{\Liegroup{}}}{\g{}}}{\Injection{\Liepoints{\aliegroup{}}}{\G{}}}\)}
{\Lieman{\Liegroup{}}}{\Lieman{\Liegroup{}}}\]}{\func{\avecf{}}{\IG{}}}}{\avecf{}\in\LiegroupLiealgebra{\aliegroup{}}}.
\end{align}
since the intrinsic group structure of the smooth group $\aliegroup{}$ is a subgroup of the intrinsic group of
$\Liegroup{}$, it is straightforward to verify that,
\begin{equation}\label{thmimmersedsubgroupsandliesubalgebracorrespondencep11eq2}
\cmp{\gltrans{\LieG{\Liegroup{}}}{\g{}}}{\Injection{\Liepoints{\aliegroup{}}}{\G{}}}=
\cmp{\Injection{\Liepoints{\aliegroup{}}}{\G{}}}{\gltrans{\LieG{\aliegroup{}}}{\g{}}}.
\end{equation}
Therefore, using the chain rule of differentiation,
\begin{align}\label{thmimmersedsubgroupsandliesubalgebracorrespondencep11eq3}
\der{\(\cmp{\gltrans{\LieG{\Liegroup{}}}{\g{}}}{\Injection{\Liepoints{\aliegroup{}}}{\G{}}}\)}
{\Lieman{\Liegroup{}}}{\Lieman{\Liegroup{}}}&=
\der{\(\cmp{\Injection{\Liepoints{\aliegroup{}}}{\G{}}}{\gltrans{\LieG{\aliegroup{}}}{\g{}}}\)}
{\Lieman{\Liegroup{}}}{\Lieman{\Liegroup{}}}\cr
&=\cmp{\(\der{\Injection{\Liepoints{\aliegroup{}}}{\G{}}}{\Lieman{\aliegroup{}}}{\Lieman{\Liegroup{}}}\)}
{\(\der{\gltrans{\LieG{\aliegroup{}}}{\g{}}}{\Lieman{\aliegroup{}}}{\Lieman{\aliegroup{}}}\)},
\end{align}
and hence according to \reflem{lemleftinvariantvectorfieldsequiv0},
and \Ref{thmimmersedsubgroupsandliesubalgebracorrespondencep11eq1},
\begin{align}\label{thmimmersedsubgroupsandliesubalgebracorrespondencep11eq4}
\defSet{\func{\avecff{}}{\g{}}}{\avecff{}\in\LS{}}&=
\defSet{\func{\[\cmp{\(\der{\Injection{\Liepoints{\aliegroup{}}}{\G{}}}{\Lieman{\aliegroup{}}}{\Lieman{\Liegroup{}}}\)}
{\(\der{\gltrans{\LieG{\aliegroup{}}}{\g{}}}{\Lieman{\aliegroup{}}}{\Lieman{\aliegroup{}}}\)}\]}
{\func{\avecf{}}{\IG{}}}}{\avecf{}\in\LiegroupLiealgebra{\aliegroup{}}}\cr
&=\defSet{\func{\(\der{\Injection{\Liepoints{\aliegroup{}}}{\G{}}}{\Lieman{\aliegroup{}}}{\Lieman{\Liegroup{}}}\)}
{\func{\[\der{\gltrans{\LieG{\aliegroup{}}}{\g{}}}{\Lieman{\aliegroup{}}}{\Lieman{\aliegroup{}}}\]}{\func{\avecf{}}{\IG{}}}}}
{\avecf{}\in\LiegroupLiealgebra{\aliegroup{}}}\cr
&=\defSet{\func{\(\der{\Injection{\Liepoints{\aliegroup{}}}{\G{}}}{\Lieman{\aliegroup{}}}{\Lieman{\Liegroup{}}}\)}
{\func{\avecf{}}{\g{}}}}
{\avecf{}\in\LiegroupLiealgebra{\aliegroup{}}}\cr
&=\func{\image{\(\der{\Injection{\Liepoints{\aliegroup{}}}{\G{}}}{\Lieman{\aliegroup{}}}{\Lieman{\Liegroup{}}}\)}}
{\defSet{\func{\avecf{}}{\g{}}}{\avecf{}\in\LiegroupLiealgebra{\aliegroup{}}}}\cr
&=\func{\image{\(\der{\Injection{\Liepoints{\aliegroup{}}}{\G{}}}{\Lieman{\aliegroup{}}}{\Lieman{\Liegroup{}}}\)}}
{\tanspace{\g{}}{\Lieman{\aliegroup{}}}}.
\end{align}
Moreover, according to \refdef{defdistributioninducedbyliesubalgebraofliegroupliealgebra}
\begin{equation}
\func{\[\func{\subliealgdist{\Liegroup{}}{r}}{\LS{}}\]}{\g{}}=
\defSet{\func{\avecff{}}{\g{}}}{\avecff{}\in\LS{}},
\end{equation}
where, $\func{\subliealgdist{\Liegroup{}}{r}}{\LS{}}$ is an $r$-dimensional integrable smooth distribution of the underlying
manifold $\Lieman{\Liegroup{}}$ of the smooth group $\Liegroup{}$, according to
\refthm{thmdistributioninducedbyliesubalgebraofliegroupliealgebraisintegrable}.
\endp
\end{itemize}
Therefore,
\begin{equation}\label{thmimmersedsubgroupsandliesubalgebracorrespondencep1eq4}
\Foreach{\g{}}{\Liepoints{\aliegroup{}}}
\func{\image{\(\der{\Injection{\Liepoints{\aliegroup{}}}{\G{}}}{\Lieman{\aliegroup{}}}{\Lieman{\Liegroup{}}}\)}}
{\tanspace{\g{}}{\Lieman{\aliegroup{}}}}=
\func{\[\func{\subliealgdist{\Liegroup{}}{r}}{\LS{}}\]}{\g{}}.
\end{equation}
Thus, based on \refthm{thmintegralmanifoldofadistributionisanopensetofaleafogthesolutionfoliation},
$\Lieman{\aliegroup{}}$ is a connected integral manifold of the integrable smooth distribution
$\func{\subliealgdist{\Liegroup{}}{r}}{\LS{}}$, and hence the set of its points is
an open set of the intrinsic topological-space of the leaf of the foliation
$\intdistsolution{\Lieman{\Liegroup{}}}{r}\func{\subliealgdist{\Liegroup{}}{r}}{\LS{}}$
(the unique solution of the integrable smooth distribution $\func{\subliealgdist{\Liegroup{}}{r}}{\LS{}}$ of $\Lieman{\Liegroup{}}$).
That is, according to \refdef{defidentityleafofthefoliationoftheinduceddistributionofasubalgebra},
$\Liepoints{\aliegroup{}}$ is an open set of the smooth group $\func{\liesublie{\Liegroup{}}{r}}{\LS{}}$.
\begin{equation}
\Liepoints{\aliegroup{}}\in\lietop{\func{\liesublie{\Liegroup{}}{r}}{\LS{}}}.
\end{equation}
Therefore, since $\Liepoints{\aliegroup{}}$ also contains $\IG{}$ which is the identity element of the underlying topological group
$\Lietopg{\func{\liesublie{\Liegroup{}}{r}}{\LS{}}}$, clearly $\Liepoints{\aliegroup{}}$ is a nucleus of this topological group,
and thus according to \refthm{thmeverynucleusgeneratestheconnectedtopologicalgroup},
and considering that the underlying topological-space of $\func{\liesublie{\Liegroup{}}{r}}{\LS{}}$ is connected
(and hence so is that of its underlying topological group structure),
$\Liepoints{\aliegroup{}}$ is a generator
of the intrinsic group of $\func{\liesublie{\Liegroup{}}{r}}{\LS{}}$. That is,
\begin{equation}
\Gen{\Liepoints{\aliegroup{}}}{\LieG{\func{\liesublie{\Liegroup{}}{r}}{\LS{}}}}=\func{\liesublieset{\Liegroup{}}{r}}{\LS{}}.
\end{equation}
In addition, since $\Liepoints{\aliegroup{}}$ is a subgroup of the intrinsic group of $\func{\liesublie{\Liegroup{}}{r}}{\LS{}}$,
clearly it is also a generator of itself, that is,
\begin{equation}
\Gen{\Liepoints{\aliegroup{}}}{\LieG{\func{\liesublie{\Liegroup{}}{r}}{\LS{}}}}=\Liepoints{\aliegroup{}}.
\end{equation}
Therefore,
\begin{equation}
\Liepoints{\aliegroup{}}=
\func{\liesublieset{\Liegroup{}}{r}}{\LS{}},
\end{equation}
and since based on \refthm{thmintegralmanifoldofadistributionisanopensetofaleafogthesolutionfoliation},
$\Lieman{\aliegroup{}}$ is an embedded submanifold of $\Lieman{\func{\liesublie{\Liegroup{}}{r}}{\LS{}}}$, clearly,
\begin{equation}
\Lieman{\aliegroup{}}=\Lieman{\func{\liesublie{\Liegroup{}}{r}}{\LS{}}}.
\end{equation}
Furthermore, considering that both $\aliegroup{}$ and $\func{\liesublie{\Liegroup{}}{r}}{\LS{}}$ are
smooth subgroups of $\Liegroup{}$, their group operations also coincide, and ultimately it is inferred that
both coincide as smooth groups.
\begin{equation}
\aliegroup{}=\func{\liesublie{\Liegroup{}}{r}}{\LS{}}.
\end{equation}
\endp
\end{itemize}
Therefore it becomes evident that the mapping $\liesubgalcor{\Liegroup{}}$ is injective.
\endthm
\corollary\label{corconnectedimmersedsubgroupsandliesubalgebrascorrespondence}
The cardinality of the set of all connected immersed smooth subgroups of the smooth group $\Liegroup{}$
equals the cardinality of the set of all non-trivial Lie-subalgebras of the canonical Lie-algebra of
the $\Liegroup{}$.
\begin{equation}
\Card{\connectedimsubgroup{\Liegroup{}}}=
\Card{\compl{\sublie{\LiegroupLiealgebra{\Liegroup{}}}}{\subliedim{0}{\LiegroupLiealgebra{\Liegroup{}}}}}.
\end{equation}
\endcor
\theorem
If $\Liegroup{1}$ is an embedded smooth subgroup of $\Liegroup{}$, then the set of its points $\G{1}$ is a closed set of
the underlying topological-space $\lietops{\Liegroup{}}$ of the smooth group $\Liegroup{}$.
\begin{equation}
\Liegroup{1}\in\emsubgroup{\Liegroup{}}\then
\G{1}\in\closedsets{\lietops{\Liegroup{}}}.
\end{equation}
\proof
It is assumed that $\Liegroup{1}$ is an embedded smooth subgroup of $\Liegroup{}$.
Then according to \refdef{defimmersedliesubgroup}, $\Lieman{\Liegroup{1}}$ is an embedded submanifold of
$\Lieman{\Liegroup{}}$.
It is also a well-known fact in differential geometry that the set of points of any embedded submanifold of a manifold
is a locally-closed set of the underlying topological-space of that manifold \cite[p.~23,~chap.~2,~sec.~2]{Lang}.
Therefore, $\G{1}$ is a locally-closed
set of the topological-space $\lietops{\Liegroup{}}$. In addition, according to \refdef{defimmersedliesubgroup},
$\G{1}$ is also a subgroup of the intrinsic group $\opair{\G{}}{\gop{}}$ of the smooth group $\Liegroup{}$.
Moreover, based on \refthm{thmliegroupisatopologicalgroup}, $\triple{\G{}}{\gop{}}{\lietops{\Liegroup{}}}$ is
a topological group. Therefore $\G{1}$ is a locally-closed subgroup of the topological group
$\triple{\G{}}{\gop{}}{\lietops{\Liegroup{}}}$, and hence according to \refthm{thmlocallyclosedtopologicalsubgroupisclosed},
it is a closed set of the topological-space $\lietops{\Liegroup{}}$.
\endthm
\section{The Exponential Mapping of a Smooth Group}
\theorem\label{thmleftinvariantvectorfieldsarecomplete}
Every left-invariant vector-field on the smooth group $\Liegroup{}$ is a complete vector-field on
the underlying manifold $\Lieman{\Liegroup{}}$ of $\Liegroup{}$. That is, for every $\avecf{}$ in $\Leftinvvf{\Liegroup{}}$,
the domain of the maximal integral-curve of $\avecf{}$ with the initial-condition $\g{}$ equals $\R$, for every point $\g{}$
of $\Liegroup{}$.
\begin{equation}
\Foreach{\avecf{}}{\Leftinvvf{\Liegroup{}}}
\Foreach{\g{}}{\G{}}
\func{\maxinterval{\Lieman{\Liegroup{}}}}{\binary{\avecf{}}{\g{}}}=\R.
\end{equation}
\proof
$\avecf{}$ is taken as an arbitrary element of $\Leftinvvf{\Liegroup{}}$, and $\g{}$ as an arbitrary element of $\G{}$.
According to \refcor{corsubcurvesofmaximalintegralcurve}, there exists a positive real number $\varepsilon$ and an integral-curve
$\curve{0}$ of the smooth vector-field $\avecf{}$ on $\Lieman{\Liegroup{}}$
with the initial condition $\g{}$ such that the domain of $\curve{0}$ is $\Ointerval{-\varepsilon}{\varepsilon}$.
\begin{equation}\label{thmleftinvariantvectorfieldsarecompletepeq1}
\Existsis{\varepsilon}{\Rp}
\Existsis{\curve{0}}{\func{\integralcurves{\Lieman{\Liegroup{}}}}{\binary{\avecf{}}{\g{}}}}
\domain{\curve{0}}=\Ointerval{-\varepsilon}{\varepsilon}.
\end{equation}
Therefore, by defining $\aninterval{}:=\Ointerval{-\varepsilon}{\varepsilon}$, according to \refdef{defintegralcurvesofasmoothvectorfield},
$\curve{}\in\Func{\aninterval{}}{\G{}}$, and,
\begin{align}
&\func{\curve{0}}{0}=\g{},
\label{thmleftinvariantvectorfieldsarecompletepeq2}\\
&\Foreach{t}{\Ointerval{-\varepsilon}{\varepsilon}}
\func{\[\der{\curve{0}}{\Ropenman{\aninterval{}}{}}{\Lieman{\Liegroup{}}}\]}{\Rtanidentity{\aninterval{}}{t}}=\func{\(\cmp{\avecf{}}{\curve{0}}\)}{t}.
\label{thmleftinvariantvectorfieldsarecompletepeq3}
\end{align}
Define,
\begin{align}
\delta:=&\frac{\varepsilon}{2},
\label{thmleftinvariantvectorfieldsarecompletepeq4}\\
\hh:=&\func{\curve{0}}{\delta},
\label{thmleftinvariantvectorfieldsarecompletepeq5}\\
\p{\hh}:=&\func{\curve{0}}{-\delta},
\label{thmleftinvariantvectorfieldsarecompletepeq6}
\end{align}
and the mappings $\function{\beta}{\Ointerval{0}{\varepsilon+\delta}}{\G{}}$ and
$\function{\p{\beta}}{\Ointerval{-\varepsilon-\delta}{0}}{\G{}}$ as,
\begin{align}
&\Foreach{t}{\Ointerval{0}{\varepsilon+\delta}}
\func{\beta}{t}\eqdef\(\hh\gop{}\invG{\g{}}{}\)\gop{}\func{\curve{0}}{t-\delta},
\label{thmleftinvariantvectorfieldsarecompletepeq7}\\
&\Foreach{t}{\Ointerval{-\varepsilon-\delta}{0}}
\func{\p{\beta}}{t}\eqdef\(\p{\hh}\gop{}\invG{\g{}}{}\)\gop{}\func{\curve{0}}{t+\delta}.
\label{thmleftinvariantvectorfieldsarecompletepeq8}
\end{align}
Clearly, by defining the function $\function{T}{\Ointerval{0}{\varepsilon+\delta}}{\Ointerval{-\varepsilon}{\varepsilon}}$
as,
\begin{equation}\label{thmleftinvariantvectorfieldsarecompletepeq9}
\Foreach{t}{\Ointerval{\delta}{\varepsilon+\delta}}
\func{T}{t}\eqdef t-\delta,
\end{equation}
it is evident that,
\begin{equation}\label{thmleftinvariantvectorfieldsarecompletepeq10}
\beta=\cmp{\cmp{\gltrans{\LieG{\Liegroup{}}}{\hh}}{\alpha_0}}{T}.
\end{equation}
By defining $J:=\Ointerval{0}{\varepsilon+\delta}$, it is trivial that,
\begin{equation}\label{thmleftinvariantvectorfieldsarecompletepeq11}
\Foreach{t}{J}
\func{\[\der{T}{\Ropenman{J}{}}{\Ropenman{\aninterval{}}{}}\]}{\Rtanidentity{J}{t}}=
\Rtanidentity{\aninterval{}}{t-\delta}.
\end{equation}
Therefore, according to \Ref{thmleftinvariantvectorfieldsarecompletepeq10} and using the chain rule of differentiation,
and according to \Ref{thmleftinvariantvectorfieldsarecompletepeq3} and
\refdef{defleftinvariantvectorfields},
\begin{align}\label{thmleftinvariantvectorfieldsarecompletepeq12}
\Foreach{t}{J}
\func{\[\der{\beta}{\Ropenman{J}{}}{\Lieman{\Liegroup{}}}\]}{\Rtanidentity{J}{t}}&=
\func{\(\cmp{\cmp
{\[\der{\gltrans{\LieG{\Liegroup{}}}{\hh}}{\Lieman{\Liegroup{}}}{\Lieman{\Liegroup{}}}\]}
{\[\der{\curve{0}}{\Ropenman{\aninterval{}}{}}{\Lieman{\Liegroup{}}}\]}}
{\[\der{T}{\Ropenman{J}{}}{\Ropenman{\aninterval{}}{}}\]}\)}{\Rtanidentity{J}{t}}\cr
&=\func{\[\der{\gltrans{\LieG{\Liegroup{}}}{\hh}}{\Lieman{\Liegroup{}}}{\Lieman{\Liegroup{}}}\]}
{\func{\[\der{\curve{0}}{\Ropenman{\aninterval{}}{}}{\Lieman{\Liegroup{}}}\]}{\Rtanidentity{\aninterval{}}{t-\delta}}}\cr
&=\func{\(\cmp{\[\der{\gltrans{\LieG{\Liegroup{}}}{\hh}}{\Lieman{\Liegroup{}}}{\Lieman{\Liegroup{}}}\]}{\avecf{}}\)}
{\func{\curve{0}}{t-\delta}}\cr
&=\func{\[\cmp{\avecf{}}{\gltrans{\LieG{\Liegroup{}}}{\hh}}\]}{\func{\curve{0}}{t-\delta}}\cr
&=\func{\avecf{}}{\hh\gop{}\func{\curve{0}}{t-\delta}}\cr
&=\func{\(\cmp{\avecf{}}{\beta}\)}{t}.
\end{align}
In a completely similar manner, by defining $\p{J}:=\Ointerval{-\delta-\varepsilon}{0}$, it can be seen that,
\begin{equation}\label{thmleftinvariantvectorfieldsarecompletepeq13}
\Foreach{t}{\p{J}}
\func{\[\der{\p{\beta}}{\Ropenman{\p{J}}{}}{\Man{}}\]}{\Rtanidentity{\p{J}}{t}}=
\func{\(\cmp{\avecf{}}{\p{\beta}}\)}{t}.
\end{equation}
Moreover,
\begin{equation}\label{thmleftinvariantvectorfieldsarecompletepeq14}
\func{\beta}{\delta}=\(\hh\gop{}\invG{\g{}}{}\)\gop{}\func{\curve{0}}{0}
=\(\hh\gop{}\invG{\g{}}{}\)\gop{}\g{}=\hh=\func{\curve{0}}{\delta}.
\end{equation}
Thus by defining the mappings $\function{\overline{\curve{0}}}{\Ointerval{-\varepsilon-\delta}{\varepsilon-\delta}}{\G{}}$,
and $\function{\overline{\beta}}{\Ointerval{-\delta}{\varepsilon}}{\G{}}$ as,
\begin{align}
&\Foreach{t}{\Ointerval{-\varepsilon-\delta}{\varepsilon-\delta}}
\func{\overline{\curve{0}}}{t}\eqdef\func{\curve{0}}{t+\delta},
\label{thmleftinvariantvectorfieldsarecompletepeq15}\\
&\Foreach{t}{\Ointerval{-\delta}{\varepsilon}}
\func{\overline{\beta}}{t}\eqdef\func{\beta}{t+\delta},
\label{thmleftinvariantvectorfieldsarecompletepeq16}
\end{align}
it is clear that,
\begin{equation}\label{thmleftinvariantvectorfieldsarecompletepeq17}
\func{\overline{\curve{0}}}{0}=\func{\overline{\beta}}{0}=\hh.
\end{equation}
Therefore, according to \Ref{thmleftinvariantvectorfieldsarecompletepeq2} and
\Ref{thmleftinvariantvectorfieldsarecompletepeq12}, and based on \refdef{defintegralcurvesofasmoothvectorfield},
it is straightforward to verify that
$\overline{\curve{0}}$ and $\overline{\beta}$ are integral-curves of $\avecf{}$ with the initial condition $\hh$.
Thus, according to \refthm{defintegralcurvesofasmoothvectorfield}, $\overline{\curve{0}}$ and $\overline{\beta}$
must coincide within the intersection of their domains. That is,
\begin{equation}\label{thmleftinvariantvectorfieldsarecompletepeq18}
\Foreach{t}{\Ointerval{-\delta}{\varepsilon-\delta}}
\func{\overline{\beta}}{t}=
\func{\overline{\curve{0}}}{t},
\end{equation}
which trivially implies that $\curve{0}$ and $\beta$ coincide when restricted to the intersection of their domains. That is,
\begin{equation}\label{thmleftinvariantvectorfieldsarecompletepeq19}
\Foreach{t}{\Ointerval{0}{\varepsilon}}
\func{\beta}{t}=
\func{\curve{0}}{t}.
\end{equation}
It can be similarly inferred that $\p{\beta}$ and $\curve{0}$ coincide when restricted to the intersection of their domains.
\begin{equation}\label{thmleftinvariantvectorfieldsarecompletepeq20}
\Foreach{t}{\Ointerval{-\varepsilon}{0}}
\func{\p{\beta}}{t}=
\func{\curve{0}}{t}.
\end{equation}
Therefore a mapping can be obtained by gluing the mappings $\curve{0}$, $\beta$, and $\p{\beta}$ together.
Hence the mapping $\curve{1}$ is well-defined, which is defined as,
\begin{align}\label{thmleftinvariantvectorfieldsarecompletepeq21}
&\curve{1}\indef{\Ointerval{\delta-\varepsilon}{\delta+\varepsilon}},\cr
&\begin{aligned}
\func{\curve{1}}{t}\eqdef
\begin{cases}
\func{\p{\beta}}{t}=\p{\hh}\gop{}\func{\curve{0}}{t+\delta}, &t\in\OCinterval{-\delta-\varepsilon}{-\delta},\cr
\func{\curve{0}}{t}, &t\in\Ointerval{-\delta}{\delta},\cr
\func{\beta}{t}=\hh\gop{}\func{\curve{0}}{t-\delta}, &t\in\COinterval{\delta}{\varepsilon+\delta}.
\end{cases}
\end{aligned}
\end{align}
Thus, by defining ${\aninterval{}}^{\(1\)}:=\Ointerval{\delta-\varepsilon}{\delta+\varepsilon}$,
according to \Ref{thmleftinvariantvectorfieldsarecompletepeq3}, \Ref{thmleftinvariantvectorfieldsarecompletepeq12},
and \Ref{thmleftinvariantvectorfieldsarecompletepeq13}, it is straightforward to verify that,
\begin{align}
&\func{\curve{1}}{0}=\g{},
\label{thmleftinvariantvectorfieldsarecompletepeq15}\\
&\Foreach{t}{\Ointerval{-\delta-\varepsilon}{\varepsilon+\delta}}
\func{\[\der{\curve{1}}{\Ropenman{{\aninterval{}}^{\(1\)}}{}}{\Lieman{\Liegroup{}}}\]}{\Rtanidentity{{\aninterval{}}^{\(1\)}}{t}}=
\func{\(\cmp{\avecf{}}{\curve{1}}\)}{t}.
\label{thmleftinvariantvectorfieldsarecompletepeq16}
\end{align}
Hence, according to \refdef{defintegralcurvesofasmoothvectorfield}, $\curve{1}$ is also an integral-curve of $\avecf{}$
with the initial condition $\g{}$. That is,
\begin{equation}
\curve{1}\in
\func{\integralcurves{\Lieman{\Liegroup{}}}}{\binary{\avecf{}}{\g{}}}.
\end{equation}
Also, it is clear that,
\begin{equation}
\domain{\curve{1}}=\Ointerval{-\frac{3}{2}\varepsilon}{\frac{3}{2}\varepsilon}.
\end{equation}
By iteration of the above process inductively, it becomes evident that for every $n$ in $\Zp$,
there exists an integral-curve of $\avecf{}$ with the initial condition $\g{}$ whose domain is
$\Ointerval{-{\(\frac{3}{2}\)}^{n}\varepsilon}{{\(\frac{3}{2}\)}^{n}\varepsilon}$. That is,
\begin{equation}
\Foreach{n}{\Zp}
\Existsis{\curve{n}}
{\func{\integralcurves{\Lieman{\Liegroup{}}}}{\binary{\avecf{}}{\g{}}}}
\domain{\curve{n}}=\Ointerval{-{\(\frac{3}{2}\)}^{n}\varepsilon}{{\(\frac{3}{2}\)}^{n}\varepsilon}.
\end{equation}
Therefore, according to \refdef{deftimeintervalofmaximalintegralcurve},
\begin{align}
\Unionn{n}{1}{\infty}{\Ointerval{-{\(\frac{3}{2}\)}^{n}\varepsilon}{{\(\frac{3}{2}\)}^{n}\varepsilon}}&\subseteq
\union{\defSet{\domain{\curve{}}}{\curve{}\in\func{\integralcurves{\Lieman{\Liegroup{}}}}{\binary{\avecf{}}{\g{}}}}}\cr
&=\func{\maxinterval{\Lieman{\Liegroup{}}}}{\binary{\avecf{}}{\g{}}},
\end{align}
and hence trivially,
\begin{equation}
\func{\maxinterval{\Lieman{\Liegroup{}}}}{\binary{\avecf{}}{\g{}}}=\R.
\end{equation}
\endthm
\definition\label{defexponentialmappingofsmoothgroup}
The mapping $\expLie{\Liegroup{}}$ is defined as,
\begin{align}
&\expLie{\Liegroup{}}\indef\Func{\Leftinvvf{\Liegroup{}}}{\G{}},\cr
&\Foreach{\avecf{}}{\Leftinvvf{\Liegroup{}}}
\func{\expLie{\Liegroup{}}}{\avecf{}}\eqdef
\func{\[\func{\maxintcurve{\Lieman{\Liegroup{}}}}{\binary{\avecf{}}{\IG{}}}\]}{1},
\end{align}
and is referred to as the $\quotl$exponential-mapping of the smooth group $\Liegroup{}$$\quotr$.\\
Simply, the exponential mapping of $\Liegroup{}$ assigns to each left-invariant vector-field $\avecf{}$ of $\Liegroup{}$,
the point of $\Liegroup{}$ that is reached by the maximal integral-curve of $\avecf{}$ on $\Lieman{\Liegroup{}}$
with the initial condition $\IG{}$
(the identity element of the group $\LieG{\Liegroup{}}$) at the time $t=1$.
\endef
\lemma\label{lemidentitycurveisaoneparametersubgroup}
$\avecf{}$ is taken as an element of $\Leftinvvf{\Liegroup{}}$. The maximal integral-curve of
$\avecf{}$ on $\Lieman{\Liegroup{}}$ with the initial condition $\IG{}$ is a group-homomorphism from the additive group
of $\R$, that is $\opair{\R}{+}$, to the underlying group structure $\LieG{\Liegroup{}}$ of the smooth group $\Liegroup{}$. That is,
\begin{equation}
\func{\maxintcurve{\Lieman{\Liegroup{}}}}{\binary{\avecf{}}{\IG{}}}\in
\GHom{\opair{\R}{+}}{\LieG{\Liegroup{}}},
\end{equation}
or equivalently,
\begin{equation}
\Foreach{\opair{t}{s}}{\R^2}
\func{\[\func{\maxintcurve{\Lieman{\Liegroup{}}}}{\binary{\avecf{}}{\IG{}}}\]}{s+t}
=\func{\[\func{\maxintcurve{\Lieman{\Liegroup{}}}}{\binary{\avecf{}}{\IG{}}}\]}{s}
\gop{}
\func{\[\func{\maxintcurve{\Lieman{\Liegroup{}}}}{\binary{\avecf{}}{\IG{}}}\]}{t}.
\end{equation}
Consequently, $\func{\maxintcurve{\Lieman{\Liegroup{}}}}{\binary{\avecf{}}{\IG{}}}$ is a smooth-group-morphism
from the smooth group $\opair{\RR}{+}$ to $\Liegroup{}$. That is,
\begin{equation}
\func{\maxintcurve{\Lieman{\Liegroup{}}}}{\binary{\avecf{}}{\IG{}}}\in
\LieMor{\opair{\RR}{+}}{\Liegroup{}}.
\end{equation}
\proof
$s$ is taken as an arbitrary real number. The mappings
$\function{\curve{}}{\R}{\G{}}$ and $\function{\beta}{\R}{\G{}}$ are defined as,
\begin{align}
&\Foreach{t}{\R}\func{\curve{}}{t}\eqdef
\func{\[\func{\maxintcurve{\Lieman{\Liegroup{}}}}{\binary{\avecf{}}{\IG{}}}\]}{s+t},\\
&\Foreach{t}{\R}\func{\beta}{t}=
\func{\[\func{\maxintcurve{\Lieman{\Liegroup{}}}}{\binary{\avecf{}}{\IG{}}}\]}{s}
\gop{}
\func{\[\func{\maxintcurve{\Lieman{\Liegroup{}}}}{\binary{\avecf{}}{\IG{}}}\]}{t}.
\end{align}
Since $\func{\maxintcurve{\Lieman{\Liegroup{}}}}{\binary{\avecf{}}{\IG{}}}$ is the maximal integral-curve of
$\avecf{}$ with the initial condition $\IG{}$, it is evident that $\curve{0}$ must be the maximal integral-curve
of $\avecf{}$ with the initial condition $\hh$, where by definition,
\begin{equation}
\hh:=\func{\[\func{\maxintcurve{\Lieman{\Liegroup{}}}}{\binary{\avecf{}}{\IG{}}}\]}{s}.
\end{equation}
That is,
\begin{equation}
\curve{}=\func{\maxintcurve{\Lieman{\Liegroup{}}}}{\binary{\avecf{}}{\hh}}.
\end{equation}
Considering the definition of $\beta$, it is clear that,
\begin{equation}
\beta=
\cmp{\gltrans{\LieG{\Liegroup{}}}{\hh}}{\func{\maxintcurve{\Lieman{\Liegroup{}}}}{\binary{\avecf{}}{\IG{}}}},
\end{equation}
and hence using the chain rule of differentiation,
considering that $\func{\maxintcurve{\Lieman{\Liegroup{}}}}{\binary{\avecf{}}{\IG{}}}$ is an integral-curve of $\avecf{}$
with the initial condition $\IG{}$,
and according to \refdef{defleftinvariantvectorfields},
\begin{align}
\Foreach{t}{\R}
\func{\[\der{\beta}{\RR}{\Lieman{\Liegroup{}}}\]}{\Rtanidentity{\R}{t}}&=
\func{\(\cmp
{\[\der{\gltrans{\LieG{\Liegroup{}}}{\hh}}{\Lieman{\Liegroup{}}}{\Lieman{\Liegroup{}}}\]}
{\[\der{\[\func{\maxintcurve{\Lieman{\Liegroup{}}}}{\binary{\avecf{}}{\IG{}}}\]}{\RR}{\Lieman{\Liegroup{}}}\]}\)}{\Rtanidentity{\R}{t}}\cr
&=\func{\[\der{\gltrans{\LieG{\Liegroup{}}}{\hh}}{\Lieman{\Liegroup{}}}{\Lieman{\Liegroup{}}}\]}
{\func{\[\der{\[\func{\maxintcurve{\Lieman{\Liegroup{}}}}{\binary{\avecf{}}{\IG{}}}\]}{\RR}{\Lieman{\Liegroup{}}}\]}{\Rtanidentity{\R}{t}}}\cr
&=\func{\(\cmp{\[\der{\gltrans{\LieG{\Liegroup{}}}{\hh}}{\Lieman{\Liegroup{}}}{\Lieman{\Liegroup{}}}\]}{\avecf{}}\)}
{\func{\[\func{\maxintcurve{\Lieman{\Liegroup{}}}}{\binary{\avecf{}}{\IG{}}}\]}{t}}\cr
&=\func{\[\cmp{\avecf{}}{\gltrans{\LieG{\Liegroup{}}}{\hh}}\]}{\func{\[\func{\maxintcurve{\Lieman{\Liegroup{}}}}{\binary{\avecf{}}{\IG{}}}\]}{t}}\cr
&=\func{\avecf{}}{\hh\gop{}\func{\[\func{\maxintcurve{\Lieman{\Liegroup{}}}}{\binary{\avecf{}}{\IG{}}}\]}{t}}\cr
&=\func{\(\cmp{\avecf{}}{\beta}\)}{t}.
\end{align}
In addition,
\begin{equation}
\func{\beta}{0}=\hh\gop{}\func{\[\func{\maxintcurve{\Lieman{\Liegroup{}}}}{\binary{\avecf{}}{\IG{}}}\]}{t}=
\hh\gop{}\IG{}=\hh.
\end{equation}
Therefore, it becomes evident that $\beta$ is the maximal integral-curve of $\avecf{}$ with the initial condition $\hh$. That is,
\begin{equation}
\beta=\func{\maxintcurve{\Lieman{\Liegroup{}}}}{\binary{\avecf{}}{\hh}}.
\end{equation}
So ultimately,
\begin{equation}
\beta=\curve{},
\end{equation}
and thus,
\begin{align}
\Foreach{t}{\R}
\func{\[\func{\maxintcurve{\Lieman{\Liegroup{}}}}{\binary{\avecf{}}{\IG{}}}\]}{s+t}=
\func{\[\func{\maxintcurve{\Lieman{\Liegroup{}}}}{\binary{\avecf{}}{\IG{}}}\]}{s}
\gop{}
\func{\[\func{\maxintcurve{\Lieman{\Liegroup{}}}}{\binary{\avecf{}}{\IG{}}}\]}{t}.
\end{align}
Furthermore, considering that $\func{\maxintcurve{\Lieman{\Liegroup{}}}}{\binary{\avecf{}}{\IG{}}}$ is a smooth map from
the manifold $\RR$ to the maifold $\Lieman{\Liegroup{}}$, evidently it is also a smooth-froup-morphism from
$\opair{\RR}{+}$ to $\Liegroup{}$.
\endlem
\lemma\label{lemoneparametersubgroupisidentitycurve}
$\avecf{}$ is taken as an arbitrary element of $\Leftinvvf{\Liegroup{}}$, and $\curve{}$ as a map from $\R$ to $\G{}$.
If $\curve{}$ is a smooth-group-morphism from $\opair{\RR}{+}$ to $\G{}$,
and the value of the differential of $\curve{}$ at $\Rtanidentity{\R}{0}$ coincides with $\func{\avecf{}}{\IG{}}$,
then $\curve{}$ is the maximal integral-curve of $\avecf{}$ on $\Lieman{\Liegroup{}}$ with the initial condition $\IG{}$.
That is,
\begin{equation}
\[\curve{}\in\bigg(\GHom{\opair{\R}{+}}{\LieG{\Liegroup{}}}\cap\mapdifclass{\infty}{\RR}{\Lieman{\Liegroup{}}}\bigg),~
\func{\[\der{\curve{}}{\RR}{\Lieman{\Liegroup{}}}\]}{\Rtanidentity{\R}{0}}=\func{\avecf{}}{\IG{}}\]\then
\curve{}=\func{\maxintcurve{\Lieman{\Liegroup{}}}}{\binary{\avecf{}}{\IG{}}}.
\end{equation}
\proof
$\curve{}$ is taken as an arbitrary element of
$\bigg(\GHom{\opair{\R}{+}}{\LieG{\Liegroup{}}}\cap\mapdifclass{\infty}{\RR}{\Lieman{\Liegroup{}}}\bigg)$
such that
\begin{equation}\label{lemoneparametersubgroupisidentitycurvepeq0}
\func{\[\der{\curve{}}{\RR}{\Lieman{\Liegroup{}}}\]}{\Rtanidentity{\R}{0}}=\func{\avecf{}}{\IG{}}.
\end{equation}
Then,
\begin{equation}\label{lemoneparametersubgroupisidentitycurvepeq1}
\Foreach{\opair{t}{s}}{\R^2}\func{\curve{}}{s+t}=\func{\curve{}}{s}\gop{}\func{\curve{}}{t},
\end{equation}
and consequently,
\begin{equation}\label{lemoneparametersubgroupisidentitycurvepeq2}
\func{\curve{}}{0}=\IG{}.
\end{equation}
\begin{itemize}
\item[$\pr{1}$]
$s$ is taken as an arbitrary element of $\R$, and the function $\function{T_s}{\R}{\R}$ is defined as,
\begin{equation}\label{lemoneparametersubgroupisidentitycurvep1eq1}
\Foreach{t}{\R}\func{T_s}{t}\eqdef t+s.
\end{equation}
It is trivial that,
\begin{equation}\label{lemoneparametersubgroupisidentitycurvep1eq2}
\func{\der{T_s}{\RR}{\RR}}{\Rtanidentity{\R}{0}}=\Rtanidentity{\R}{s},
\end{equation}
and hence, using the chain rule of differentiation,
\begin{align}\label{lemoneparametersubgroupisidentitycurvep1eq3}
\func{\[\der{\curve{}}{\RR}{\Lieman{\Liegroup{}}}\]}{\Rtanidentity{\R}{s}}&=
\func{\[\der{\curve{}}{\RR}{\Lieman{\Liegroup{}}}\]}{\func{\der{T_s}{\RR}{\RR}}{\Rtanidentity{\R}{0}}}\cr
&=\func{\cmp{\[\der{\curve{}}{\RR}{\Lieman{\Liegroup{}}}\]}{\[\der{T_s}{\RR}{\RR}\]}}{\Rtanidentity{\R}{0}}\cr
&=\func{\[\der{\(\cmp{\curve{}}{T_s}\)}{\RR}{\Lieman{\Liegroup{}}}\]}{\Rtanidentity{\R}{0}}.
\end{align}
In addition, according to \Ref{lemoneparametersubgroupisidentitycurvepeq1} and \Ref{lemoneparametersubgroupisidentitycurvep1eq1},
\begin{equation}\label{lemoneparametersubgroupisidentitycurvep1eq4}
\cmp{\curve{}}{T_s}=\cmp{\gltrans{\LieG{\Liegroup{}}}{\func{\curve{}}{s}}}{\curve{}},
\end{equation}
and therefore since $\avecf{}$ is a left-invariant vector-field of $\Liegroup{}$, according to
\Ref{lemoneparametersubgroupisidentitycurvepeq0} and
\reflem{lemleftinvariantvectorfieldsequiv0},
\begin{align}\label{lemoneparametersubgroupisidentitycurvep1eq5}
\func{\[\der{\curve{}}{\RR}{\Lieman{\Liegroup{}}}\]}{\Rtanidentity{\R}{s}}&=
\func{\[\der{\(\cmp{\gltrans{\LieG{\Liegroup{}}}{\func{\curve{}}{s}}}{\curve{}}\)}{\RR}{\Lieman{\Liegroup{}}}\]}{\Rtanidentity{\R}{0}}\cr
&=\func{\(\cmp{\[\der{\gltrans{\LieG{\Liegroup{}}}{\func{\curve{}}{s}}}{\Lieman{\Liegroup{}}}{\Lieman{\Liegroup{}}}\]}
{\[\der{\curve{}}{\RR}{\Lieman{\Liegroup{}}}\]}\)}{\Rtanidentity{\R}{0}}\cr
&=\func{\[\der{\gltrans{\LieG{\Liegroup{}}}{\func{\curve{}}{s}}}{\Lieman{\Liegroup{}}}{\Lieman{\Liegroup{}}}\]}{\func{\avecf{}}{\IG{}}}\cr
&=\func{\(\cmp{\avecf{}}{\gltrans{\LieG{\Liegroup{}}}{\func{\curve{}}{s}}}\)}{\IG{}}\cr
&=\func{\avecf{}}{\func{\curve{}}{s}}.
\end{align}
\endp
\end{itemize}
Therefore,
\begin{equation}\label{lemoneparametersubgroupisidentitycurvepeq3}
\Foreach{s}{\R}
\func{\[\der{\curve{}}{\RR}{\Lieman{\Liegroup{}}}\]}{\Rtanidentity{\R}{s}}=
\func{\avecf{}}{\func{\curve{}}{s}}.
\end{equation}
Based on \refdef{defintegralcurvesofasmoothvectorfield}, \Ref{lemoneparametersubgroupisidentitycurvepeq2} and
\Ref{lemoneparametersubgroupisidentitycurvepeq2} imply that $\curve{}$ is an integral-curve of $\avecf{}$
with the initial condition $\IG{}$, and since domain of $\curve{}$ is $\R$, clearly $\curve{}$ must be the maximal integral-curve
of $\avecf{}$ with the initial condition $\IG{}$. That is,
\begin{equation}
\curve{}=\func{\maxintcurve{\Lieman{\Liegroup{}}}}{\binary{\avecf{}}{\IG{}}}.
\end{equation}
\endlem
\corollary\label{coridentitycurveistheonlyoneparametersubgroup}
There exists only one smooth-group-morphism from $\opair{\RR}{+}$ to $\G{}$,
such that the value of the differential of $\curve{}$ at $\Rtanidentity{\R}{0}$ coincides with $\func{\avecf{}}{\IG{}}$,
which is the maximal integral-curve of $\avecf{}$ on $\Lieman{\Liegroup{}}$ with the initial condition $\IG{}$. That is,
\begin{align}
\defset{\curve{}}{\GHom{\opair{\R}{+}}{\LieG{\Liegroup{}}}\cap\mapdifclass{\infty}{\RR}{\Lieman{\Liegroup{}}}}{
\func{\[\der{\curve{}}{\RR}{\Lieman{\Liegroup{}}}\]}{\Rtanidentity{\R}{0}}=\func{\avecf{}}{\IG{}}}=
\seta{\func{\maxintcurve{\Lieman{\Liegroup{}}}}{\binary{\avecf{}}{\IG{}}}}.
\end{align}
\endcor
\theorem\label{thmnaturalityofexponentialmapping}
$\aliemor{}$ is taken as an element of $\LieMor{\Liegroup{}}{\Liegroup{1}}$ (a Lie-morphism from the smooth group $\Liegroup{}$
to the smooth group $\Liegroup{1}$).
\begin{equation}
\cmp{\expLie{\Liegroup{1}}}{\[\func{\indliemor{\Liegroup{}}{\Liegroup{1}}}{\aliemor{}}\]}=
\cmp{\aliemor{}}{\expLie{\Liegroup{}}}.
\end{equation}
Equivalently, the following diagram commutes.
\begin{center}
\vskip0.5\baselineskip
\hskip-3\baselineskip
\begin{tikzcd}[row sep=7em, column sep=6em]
& \Leftinvvf{\Liegroup{}}
\arrow{r}{\func{\indliemor{\Liegroup{}}{\Liegroup{1}}}{\aliemor{}}}
\arrow[swap,d,"\expLie{\Liegroup{}}" description]
& \Leftinvvf{\Liegroup{1}}
\arrow[d,"\expLie{\Liegroup{1}}" description] \\
& \G{}
\arrow[swap,r,"\aliemor{}" description] & \G{1}
\end{tikzcd}
\end{center}
\hskip\baselineskip
\proof
$\avecf{}$ is taken as an arbitrary element of $\Leftinvvf{\Liegroup{}}$, that is a left-invariant vector-field of $\Liegroup{}$.
According to \refdef{defintegralcurvesofasmoothvectorfield} and \refthm{thmmaximalintegralcurveisanintegralcurve},
the maximal integral-curve of $\avecf{}$ on $\Lieman{\Liegroup{}}$
with the initial condition $\IG{}$, that is $\func{\maxintcurve{\Lieman{\Liegroup{}}}}{\binary{\avecf{}}{\IG{}}}$,
satisfies the following properties.
\begin{align}
&\func{\[\func{\maxintcurve{\Lieman{\Liegroup{}}}}{\binary{\avecf{}}{\IG{}}}\]}{0}=\IG{},
\label{thmnaturalityofexponentialmappingpeq1}\\
&\Foreach{t}{\R}
\func{\[\der{\[\func{\maxintcurve{\Lieman{\Liegroup{}}}}{\binary{\avecf{}}{\IG{}}}\]}{\RR}{\Lieman{\Liegroup{}}}\]}{\Rtanidentity{\R}{t}}=
\func{\[\cmp{\avecf{}}{\func{\maxintcurve{\Lieman{\Liegroup{}}}}{\binary{\avecf{}}{\IG{}}}}\]}{t}.
\label{thmnaturalityofexponentialmappingpeq2}
\end{align}
According to \reflem{lemidentitycurveisaoneparametersubgroup}, $\func{\maxintcurve{\Lieman{\Liegroup{}}}}{\binary{\avecf{}}{\IG{}}}$ is
a smooth-group-morphism from $\opair{\RR}{+}$ to $\Liegroup{}$. Therefore, since $\aliemor{}$ is a smooth-group-morphism from
$\Liegroup{}$ to $\Liegroup{1}$, \refthm{thmcompositionofsmoothgroupmorphisms} implies that their composition, that is
$\cmp{\aliemor{}}{\[\func{\maxintcurve{\Lieman{\Liegroup{}}}}{\binary{\avecf{}}{\IG{}}}\]}$ must be a smooth-group-morphism from
$\opair{\RR}{+}$ to $\Liegroup{1}$.
\begin{equation}\label{thmnaturalityofexponentialmappingpeq3}
\cmp{\aliemor{}}{\[\func{\maxintcurve{\Lieman{\Liegroup{}}}}{\binary{\avecf{}}{\IG{}}}\]}\in
\LieMor{\opair{\RR}{+}}{\Liegroup{1}}.
\end{equation}
Moreover, using the chain rule of differentiation, and according to
\Ref{thmnaturalityofexponentialmappingpeq1} and \Ref{thmnaturalityofexponentialmappingpeq2},
\begin{align}\label{thmnaturalityofexponentialmappingpeq4}
\func{\[\der{\(\cmp{\aliemor{}}
{\[\func{\maxintcurve{\Lieman{\Liegroup{}}}}{\binary{\avecf{}}{\IG{}}}}\]\)}{\RR}{\Lieman{\Liegroup{1}}}\]}{\Rtanidentity{\R}{0}}&=
\func{\[\der{\aliemor{}}{\Liegroup{}}{\Liegroup{1}}\]}{
\func{\[\der{\func{\maxintcurve{\Lieman{\Liegroup{}}}}{\binary{\avecf{}}{\IG{}}}}{\RR}{\Lieman{\Liegroup{}}}\]}{\Rtanidentity{\R}{0}}
}\cr
&=\func{\[\der{\aliemor{}}{\Liegroup{}}{\Liegroup{1}}\]}{\func{\avecf{}}{\IG{}}}.
\end{align}
In addition, by defining,
\begin{equation}\label{thmnaturalityofexponentialmappingpeq5}
\avecff{}:=\func{\[\func{\indliemor{\Liegroup{}}{\Liegroup{1}}}{\aliemor{}}\]}{\avecf{}},
\end{equation}
according to \refdef{definducedliealgebramorphismfromliemorphism} and
\refdef{deflinvvftanspacecorrespondence},
$\avecff{}$ is a left-invariant vector-field of $\Liegroup{1}$, that is,
\begin{equation}\label{thmnaturalityofexponentialmappingpeq6}
\avecff{}\in\Leftinvvf{\Liegroup{1}},
\end{equation}
and,
\begin{equation}\label{thmnaturalityofexponentialmappingpeq7}
\func{\avecff{}}{\IG{1}}=
\func{\[\der{\aliemor{}}{\Liegroup{}}{\Liegroup{1}}\]}{\func{\avecf{}}{\IG{}}}.
\end{equation}
\Ref{thmnaturalityofexponentialmappingpeq5} and \Ref{thmnaturalityofexponentialmappingpeq7} yield,
\begin{equation}\label{thmnaturalityofexponentialmappingpeq8}
\func{\[\der{\(\cmp{\aliemor{}}
{\[\func{\maxintcurve{\Lieman{\Liegroup{}}}}{\binary{\avecf{}}{\IG{}}}}\]\)}{\RR}{\Lieman{\Liegroup{1}}}\]}{\Rtanidentity{\R}{0}}=
\func{\avecff{}}{\IG{1}}.
\end{equation}\label{thmnaturalityofexponentialmappingpeq9}
Therefore, considering that $\cmp{\aliemor{}}
{\func{\maxintcurve{\Lieman{\Liegroup{}}}}{\binary{\avecf{}}{\IG{}}}}$ is a
smooth-group-morphism from $\opair{\RR}{+}$ to $\Liegroup{1}$, \reflem{lemoneparametersubgroupisidentitycurve} implies that
$\cmp{\aliemor{}}{\[\func{\maxintcurve{\Lieman{\Liegroup{}}}}{\binary{\avecf{}}{\IG{}}}\]}$ is the maximal integral-curve of
$\avecff{}$ on $\Lieman{\Liegroup{1}}$ with the initial condition $\IG{1}$. That is,
\begin{equation}\label{thmnaturalityofexponentialmappingpeq10}
\cmp{\aliemor{}}
{\[\func{\maxintcurve{\Lieman{\Liegroup{}}}}{\binary{\avecf{}}{\IG{}}}\]}=
\func{\maxintcurve{\Lieman{\Liegroup{1}}}}{\binary{\avecff{}}{\IG{1}}}.
\end{equation}
Therefore, according to \Ref{thmnaturalityofexponentialmappingpeq6} and
\refdef{defexponentialmappingofsmoothgroup},
\begin{align}\label{thmnaturalityofexponentialmappingpeq11}
\func{\(\cmp{\expLie{\Liegroup{1}}}{\[\func{\indliemor{\Liegroup{}}{\Liegroup{1}}}{\aliemor{}}\]}\)}{\avecf{}}&=
\func{\expLie{\Liegroup{1}}}{\avecff{}}\cr
&=\func{\[\func{\maxintcurve{\Lieman{\Liegroup{1}}}}{\binary{\avecff{}}{\IG{1}}}\]}{1}\cr
&=\func{\(\cmp{\aliemor{}}{\[\func{\maxintcurve{\Lieman{\Liegroup{}}}}{\binary{\avecf{}}{\IG{}}}\]}\)}{1}\cr
&=\func{\aliemor{}}{\func{\[\func{\maxintcurve{\Lieman{\Liegroup{}}}}{\binary{\avecf{}}{\IG{}}}\]}{1}}\cr
&=\func{\aliemor{}}{\func{\expLie{\Liegroup{}}}{\avecf{}}}\cr
&=\func{\[\cmp{\aliemor{}}{\expLie{\Liegroup{}}}\]}{\avecf{}}.
\end{align}
\endthm
\corollary\label{corexponetialmappingandliesubgroups}
$\Liegroup{1}$ is supposed to be an immersed smooth subgroup of $\Liegroup{}$.
\begin{equation}
\cmp{\expLie{\Liegroup{}}}{\[\func{\indliemor{\Liegroup{1}}{\Liegroup{}}}{\Injection{\G{1}}{\G{}}}\]}=
\cmp{\Injection{\G{1}}{\G{}}}{\expLie{\Liegroup{1}}}.
\end{equation}
Equivalently, the following diagram commutes.
\begin{center}
\vskip0.5\baselineskip
\hskip-3\baselineskip
\begin{tikzcd}[row sep=6em, column sep=8em]
& \Leftinvvf{\Liegroup{1}}
\arrow{r}{\func{\indliemor{\Liegroup{1}}{\Liegroup{}}}{\Injection{\G{1}}{\G{}}}}
\arrow[swap,d,"\expLie{\Liegroup{1}}" description]
& \Leftinvvf{\Liegroup{}}
\arrow[d,"\expLie{\Liegroup{}}" description] \\
& \G{1}
\arrow[swap,r,"\Injection{\G{1}}{\G{}}{}" description] & \G{}
\end{tikzcd}
\end{center}
An immediate consequence of this commutative diagram is that the image of $\expLie{\Liegroup{1}}$ is a subset of
the image of $\expLie{\Liegroup{}}$.
\begin{equation}
\funcimage{\expLie{\Liegroup{1}}}\subseteq\funcimage{\expLie{\Liegroup{}}}.
\end{equation}
\caution
Since every embedded smooth subgroup of $\Liegroup{}$ is an immersed smooth subgroup of it,
the above properties are also satisfied for embedded smooth subgroups of $\Liegroup{}$ in particular.
\proof
According to \refthm{leminjectionofsmoothsubgroupisamorphism}, $\Injection{\G{1}}{\G{}}$ is a smooth-group-morphism from
$\Liegroup{1}$ to $\Liegroup{}$. So the commutativity of the diagram follows trivially from
\refthm{thmnaturalityofexponentialmapping}.
\endcor
\corollary\label{corexponentialmappingandsmoothgroupisomorphisms}
It is supposed the smooth groups $\Liegroup{}$ and $\Liegroup{1}$ are Lie-isomorphic, and
$\aliemor{}$ is taken as a smooth-group-isomorphism from $\Liegroup{}$ to $\Liegroup{1}$. The following diagram commutes.
\begin{center}
\vskip0.5\baselineskip
\hskip-3\baselineskip
\begin{tikzcd}[row sep=7em, column sep=6em]
& \Leftinvvf{\Liegroup{}}
\arrow[r, shift left, "\func{\indliemor{\Liegroup{}}{\Liegroup{1}}}{\aliemor{}}"]
\arrow[swap,d,"\expLie{\Liegroup{}}" description]
& \Leftinvvf{\Liegroup{1}}
\arrow[l, shift left, "\func{\indliemor{\Liegroup{1}}{\Liegroup{}}}{\finv{\aliemor{}}}"]
\arrow[d,"\expLie{\Liegroup{1}}" description] \\
& \G{}
\arrow[swap, r, shift right, "\aliemor{}"] & \G{1}
\arrow[swap, l, shift right, "\finv{\aliemor{}}"]
\end{tikzcd}
\end{center}
\hskip\baselineskip\\
Therefore, if one of $\expLie{\Liegroup{}}$ and $\expLie{\Liegroup{1}}$ is injective or surjective (or bijective), so must be the other one.
\endcor
\lemma\label{lemleftinvvfintegralcurveswitharbitraryIC}
$\avecf{}$ is taken as an element of $\Leftinvvf{\Liegroup{}}$, and $\g{}$ as a point of $\Liegroup{}$.
The maximal integral-curve of $\avecf{}$ on $\Lieman{\Liegroup{}}$ with the initial condition $\g{}$ is obtained by
composing the left-translation of $\LieG{\Liegroup{}}$ by $\g{}$ and the maximal integral-curve of $\avecf{}$ with
the initial condition $\IG{}$. That is,
\begin{equation}
\func{\maxintcurve{\Lieman{\Liegroup{}}}}{\binary{\avecf{}}{\g{}}}=
\cmp{\gltrans{\LieG{\Liegroup{}}}{\g{}}}{\[\func{\maxintcurve{\Lieman{\Liegroup{}}}}{\binary{\avecf{}}{\IG{}}}\]}.
\end{equation}
Equivalently,
\begin{equation}
\Foreach{t}{\R}
\func{\[\func{\maxintcurve{\Lieman{\Liegroup{}}}}{\binary{\avecf{}}{\g{}}}\]}{t}=
\g{}\gop{}\func{\[\func{\maxintcurve{\Lieman{\Liegroup{}}}}{\binary{\avecf{}}{\IG{}}}\]}{t}.
\end{equation}
\proof
Clearly,
\begin{align}\label{lemleftinvvfintegralcurveswitharbitraryICpeq1}
\func{\(\cmp{\gltrans{\LieG{\Liegroup{}}}{\g{}}}{\[\func{\maxintcurve{\Lieman{\Liegroup{}}}}{\binary{\avecf{}}{\IG{}}}\]}\)}{0}&=
\g{}\gop{}\func{\[\func{\maxintcurve{\Lieman{\Liegroup{}}}}{\binary{\avecf{}}{\IG{}}}\]}{0}\cr
&=\g{}\gop{}\IG{}=\g{}.
\end{align}
Since $\gltrans{\LieG{\Liegroup{}}}{\g{}}$ and $\func{\maxintcurve{\Lieman{\Liegroup{}}}}{\binary{\avecf{}}{\IG{}}}$ both are
smooth maps, trivially so is their composition, that is,
\begin{equation}\label{lemleftinvvfintegralcurveswitharbitraryICpeq2}
\cmp{\gltrans{\LieG{\Liegroup{}}}{\g{}}}{\[\func{\maxintcurve{\Lieman{\Liegroup{}}}}{\binary{\avecf{}}{\IG{}}}\]}\in
\mapdifclass{\infty}{\RR}{\Lieman{\Liegroup{}}}.
\end{equation}
In addition, considering that $\avecf{}$ is a left-invariant vector-field on $\Liegroup{}$,
according to \refdef{defleftinvariantvectorfields}, and \refdef{defintegralcurvesofasmoothvectorfield},
and using the chain rule of differentiation,
\begin{align}\label{lemleftinvvfintegralcurveswitharbitraryICpeq3}
\Foreach{t}{\R}
\func{\[\der{\(\cmp{\gltrans{\LieG{\Liegroup{}}}{\g{}}}{\[\func{\maxintcurve{\Lieman{\Liegroup{}}}}{\binary{\avecf{}}{\IG{}}}\]}\)}
{\RR}{\Lieman{\Liegroup{}}}\]}
{\Rtanidentity{\R}{t}}&=
\func{\(\cmp{\der{\gltrans{\LieG{\Liegroup{}}}{\g{}}}{\Lieman{\Liegroup{}}}{\Lieman{\Liegroup{}}}}
{\der{\[\func{\maxintcurve{\Lieman{\Liegroup{}}}}{\binary{\avecf{}}{\IG{}}}\]}{\RR}{\Lieman{\Liegroup{}}}}\)}{\Rtanidentity{\R}{t}}\cr
&=\func{\(\cmp{\cmp{\der{\gltrans{\LieG{\Liegroup{}}}{\g{}}}{\Lieman{\Liegroup{}}}{\Lieman{\Liegroup{}}}}{\avecf{}}}
{\func{\maxintcurve{\Lieman{\Liegroup{}}}}{\binary{\avecf{}}{\IG{}}}}\)}{t}\cr
&=\func{\(\cmp{\cmp{\avecf{}}{\gltrans{\LieG{\Liegroup{}}}{\g{}}}}{\func{\maxintcurve{\Lieman{\Liegroup{}}}}{\binary{\avecf{}}{\IG{}}}}\)}{t}\cr
&=\func{\avecf{}}{\func{\[\cmp{\gltrans{\LieG{\Liegroup{}}}{\g{}}}{\func{\maxintcurve{\Lieman{\Liegroup{}}}}{\binary{\avecf{}}{\IG{}}}}\]}{t}}.
\end{align}
Based on \refdef{defintegralcurvesofasmoothvectorfield},
\Ref{lemleftinvvfintegralcurveswitharbitraryICpeq1} and \Ref{lemleftinvvfintegralcurveswitharbitraryICpeq3}
imply that $\cmp{\gltrans{\LieG{\Liegroup{}}}{\g{}}}{\func{\maxintcurve{\Lieman{\Liegroup{}}}}{\binary{\avecf{}}{\IG{}}}}$ is
an integral-curve of $\avecf{}$ with the initial condition $\g{}$. In addition, considering that the domain of
$\cmp{\gltrans{\LieG{\Liegroup{}}}{\g{}}}{\func{\maxintcurve{\Lieman{\Liegroup{}}}}{\binary{\avecf{}}{\IG{}}}}$ is $\R$,
it must clearly be the maximal integral-curve of $\avecf{}$ with the initial condition $\g{}$.
\endlem
\lemma\label{lemtheidentitycurvesofcommutativeleftinvvectorfieldscommute}
Each $\avecf{1}$ and $\avecf{2}$ is taken as an arbitrary element of $\Leftinvvf{\Liegroup{}}$. If $\avecf{1}$ and $\avecf{2}$
commute, then $\func{\[\func{\maxintcurve{\Lieman{\Liegroup{}}}}{\binary{\avecf{1}}{\IG{}}}\]}{t}$ and
$\func{\[\func{\maxintcurve{\Lieman{\Liegroup{}}}}{\binary{\avecf{2}}{\IG{}}}\]}{s}$ commute in the group $\LieG{\Liegroup{}}$,
for every $t$ and $s$ in $\R$. That is,
\begin{align}
&\hskip0.8\baselineskip\liebracket{\avecf{}}{\avecf{2}}{\Lieman{\Liegroup{}}}=\zerovec{\Vecf{\Lieman{\Liegroup{}}}{\infty}}\cr
&\then
\Foreach{\opair{t}{s}}{\R^2}
\func{\[\func{\maxintcurve{\Lieman{\Liegroup{}}}}{\binary{\avecf{1}}{\IG{}}}\]}{t}\gop{}
\func{\[\func{\maxintcurve{\Lieman{\Liegroup{}}}}{\binary{\avecf{2}}{\IG{}}}\]}{s}=
\func{\[\func{\maxintcurve{\Lieman{\Liegroup{}}}}{\binary{\avecf{2}}{\IG{}}}\]}{s}\gop{}
\func{\[\func{\maxintcurve{\Lieman{\Liegroup{}}}}{\binary{\avecf{1}}{\IG{}}}\]}{t}.\cr
&{}
\end{align}
\proof
It is assumed that,
\begin{equation}\label{lemtheidentitycurvesofcommutativeleftinvvectorfieldscommutepeq0}
\liebracket{\avecf{}}{\avecf{2}}{\Lieman{\Liegroup{}}}=\zerovec{\Vecf{\Lieman{\Liegroup{}}}{\infty}}.
\end{equation}
Since left-invariant vector-fields on $\Liegroup{}$ are complete vector-fields on $\Lieman{\Liegroup{}}$
(\refthm{thmleftinvariantvectorfieldsarecomplete}), according to \refthm{thmflowsofcommutativevectorfieldscommute},
for every $t$ and $s$ in $\R$,
$\func{\vfFlow{\Lieman{\Liegroup{}}}{\avecf{1}}}{t}$ and $\func{\vfFlow{\Lieman{\Liegroup{}}}{\avecf{2}}}{s}$ commute
in the group of $\infty$-automorphisms of $\Lieman{\Liegroup{}}$.
\begin{align}\label{lemtheidentitycurvesofcommutativeleftinvvectorfieldscommutepeq1}
\Foreach{\opair{t}{s}}{\R^2}
\cmp{\func{\vfFlow{\Lieman{\Liegroup{}}}{\avecf{1}}}{t}}{\func{\vfFlow{\Lieman{\Liegroup{}}}{\avecf{2}}}{s}}=
\cmp{\func{\vfFlow{\Lieman{\Liegroup{}}}{\avecf{2}}}{s}}{\func{\vfFlow{\Lieman{\Liegroup{}}}{\avecf{1}}}{t}}.
\end{align}
According to \refdef{defcompletevectorfieldFlow}, \refdef{defsmoothvectorfieldflow},
and \reflem{lemleftinvvfintegralcurveswitharbitraryIC},
\begin{align}\label{lemtheidentitycurvesofcommutativeleftinvvectorfieldscommutepeq2}
\Foreach{\opair{t}{s}}{\R^2}
\func{\[\func{\vfFlow{\Lieman{\Liegroup{}}}{\avecf{1}}}{t}\]}{\func{\[\func{\vfFlow{\Lieman{\Liegroup{}}}{\avecf{2}}}{s}\]}{\IG{}}}
&=\func{\vfflow{\Lieman{\Liegroup{}}}{\avecf{1}}}{\binary{t}{\func{\vfflow{\Lieman{\Liegroup{}}}{\avecf{2}}}{\binary{s}{\IG{}}}}}\cr
&=\func{\[\func{\maxintcurve{\Lieman{\Liegroup{}}}}{\binary{\avecf{1}}
{\func{\[\func{\maxintcurve{\Lieman{\Liegroup{}}}}{\binary{\avecf{2}}{\IG{}}}\]}{s}}}\]}{t}\cr
&=\func{\[\func{\maxintcurve{\Lieman{\Liegroup{}}}}{\binary{\avecf{2}}{\IG{}}}\]}{s}\gop{}
\func{\[\func{\maxintcurve{\Lieman{\Liegroup{}}}}{\binary{\avecf{1}}{\IG{}}}\]}{t}.\cr
&{}
\end{align}
In a completely similar manner, it is inferred that,
\begin{equation}\label{lemtheidentitycurvesofcommutativeleftinvvectorfieldscommutepeq3}
\Foreach{\opair{t}{s}}{\R^2}
\func{\[\func{\vfFlow{\Lieman{\Liegroup{}}}{\avecf{2}}}{s}\]}{\func{\[\func{\vfFlow{\Lieman{\Liegroup{}}}{\avecf{1}}}{t}\]}{\IG{}}}=
\func{\[\func{\maxintcurve{\Lieman{\Liegroup{}}}}{\binary{\avecf{1}}{\IG{}}}\]}{t}\gop{}
\func{\[\func{\maxintcurve{\Lieman{\Liegroup{}}}}{\binary{\avecf{2}}{\IG{}}}\]}{s}.
\end{equation}
Therefore,
\begin{equation}
\Foreach{\opair{t}{s}}{\R^2}
\func{\[\func{\maxintcurve{\Lieman{\Liegroup{}}}}{\binary{\avecf{1}}{\IG{}}}\]}{t}\gop{}
\func{\[\func{\maxintcurve{\Lieman{\Liegroup{}}}}{\binary{\avecf{2}}{\IG{}}}\]}{s}=
\func{\[\func{\maxintcurve{\Lieman{\Liegroup{}}}}{\binary{\avecf{2}}{\IG{}}}\]}{s}\gop{}
\func{\[\func{\maxintcurve{\Lieman{\Liegroup{}}}}{\binary{\avecf{1}}{\IG{}}}\]}{t}.
\end{equation}
\endlem
\lemma\label{lemproductofidentitycurvesofapairofleftinvfsisaLiemorphism}
Each $\avecf{1}$ and $\avecf{2}$ is taken as an arbitrary element of $\Leftinvvf{\Liegroup{}}$. If $\avecf{1}$ and $\avecf{2}$
commute, then the mapping $\function{\eta}{\R}{\G{}}$ defined as,
\begin{equation}
\Foreach{t}{\R}
\func{\eta}{t}\eqdef
\func{\[\func{\maxintcurve{\Lieman{\Liegroup{}}}}{\binary{\avecf{1}}{\IG{}}}\]}{t}\gop{}
\func{\[\func{\maxintcurve{\Lieman{\Liegroup{}}}}{\binary{\avecf{2}}{\IG{}}}\]}{t},
\end{equation}
is the maximal integral-curve of $\avecf{1}+\avecf{2}$ on $\Lieman{\Liegroup{}}$
with the initial condition $\IG{}$.
\begin{equation}
\eta=\func{\maxintcurve{\Lieman{\Liegroup{}}}}{\binary{\avecf{1}+\avecf{2}}{\IG{}}}.
\end{equation}
\proof
It is assumed that,
\begin{equation}\label{lemproductofidentitycurvesofapairofleftinvfsisaLiemorphismpeq1}
\liebracket{\avecf{}}{\avecf{2}}{\Lieman{\Liegroup{}}}=\zerovec{\Vecf{\Lieman{\Liegroup{}}}{\infty}}.
\end{equation}
Then, according to \reflem{lemtheidentitycurvesofcommutativeleftinvvectorfieldscommute},
$\func{\[\func{\maxintcurve{\Lieman{\Liegroup{}}}}{\binary{\avecf{1}}{\IG{}}}\]}{s}$ and
$\func{\[\func{\maxintcurve{\Lieman{\Liegroup{}}}}{\binary{\avecf{2}}{\IG{}}}\]}{t}$ commute
for every $s$ and $t$ in $\R$, and hence considering that each $\func{\maxintcurve{\Lieman{\Liegroup{}}}}{\binary{\avecf{1}}{\IG{}}}$
and $\func{\maxintcurve{\Lieman{\Liegroup{}}}}{\binary{\avecf{2}}{\IG{}}}$ is a group-morphism from $\opair{\R}{+}$ to $\LieG{\Liegroup{}}$
(\reflem{lemidentitycurveisaoneparametersubgroup}),
\begin{align}\label{lemproductofidentitycurvesofapairofleftinvfsisaLiemorphismpeq2}
&\Foreach{\opair{t}{s}}{\R^2}\cr
&\begin{aligned}
\func{\eta}{t+s}&=\func{\[\func{\maxintcurve{\Lieman{\Liegroup{}}}}{\binary{\avecf{1}}{\IG{}}}\]}{t+s}\gop{}
\func{\[\func{\maxintcurve{\Lieman{\Liegroup{}}}}{\binary{\avecf{2}}{\IG{}}}\]}{t+s}\cr
&=\func{\[\func{\maxintcurve{\Lieman{\Liegroup{}}}}{\binary{\avecf{1}}{\IG{}}}\]}{t}\gop{}
\func{\[\func{\maxintcurve{\Lieman{\Liegroup{}}}}{\binary{\avecf{1}}{\IG{}}}\]}{s}\gop{}
\func{\[\func{\maxintcurve{\Lieman{\Liegroup{}}}}{\binary{\avecf{2}}{\IG{}}}\]}{t}\gop{}
\func{\[\func{\maxintcurve{\Lieman{\Liegroup{}}}}{\binary{\avecf{2}}{\IG{}}}\]}{s}\cr
&=\bigg(\func{\[\func{\maxintcurve{\Lieman{\Liegroup{}}}}{\binary{\avecf{1}}{\IG{}}}\]}{t}\gop{}
\func{\[\func{\maxintcurve{\Lieman{\Liegroup{}}}}{\binary{\avecf{2}}{\IG{}}}\]}{t}\bigg)\gop{}
\bigg(\func{\[\func{\maxintcurve{\Lieman{\Liegroup{}}}}{\binary{\avecf{1}}{\IG{}}}\]}{s}\gop{}
\func{\[\func{\maxintcurve{\Lieman{\Liegroup{}}}}{\binary{\avecf{2}}{\IG{}}}\]}{s}\bigg)\cr
&=\func{\eta}{t}\gop{}\func{\eta}{s},
\end{aligned}\cr
&{}
\end{align}
which means $\eta$ is a group-homomorphism from $\opair{\R}{+}$ to $\opair{\G{}}{\gop{}}$.
\begin{equation}\label{lemproductofidentitycurvesofapairofleftinvfsisaLiemorphismpeq3}
\eta\in\GHom{\opair{\R}{+}}{\LieG{\Liegroup{}}}.
\end{equation}
Moreover, it is evident that,
\begin{equation}\label{lemproductofidentitycurvesofapairofleftinvfsisaLiemorphismpeq4}
\eta=\cmp{\cmp{\gop{}}{\(\funcprod{\[\func{\maxintcurve{\Lieman{\Liegroup{}}}}{\binary{\avecf{1}}{\IG{}}}\]}
{\[\func{\maxintcurve{\Lieman{\Liegroup{}}}}{\binary{\avecf{2}}{\IG{}}}\]}\)}}{\diagmap{\R}}.
\end{equation}
Since each $\func{\maxintcurve{\Lieman{\Liegroup{}}}}{\binary{\avecf{1}}{\IG{}}}$ and
$\func{\maxintcurve{\Lieman{\Liegroup{}}}}{\binary{\avecf{2}}{\IG{}}}$ is a smooth map from $\RR$ to $\Lieman{\Liegroup{}}$,
evidently,
\begin{equation}\label{lemproductofidentitycurvesofapairofleftinvfsisaLiemorphismpeq5}
\funcprod{\[\func{\maxintcurve{\Lieman{\Liegroup{}}}}{\binary{\avecf{1}}{\IG{}}}\]}
{\[\func{\maxintcurve{\Lieman{\Liegroup{}}}}{\binary{\avecf{2}}{\IG{}}}\]}\in\mapdifclass{\infty}
{\RR^2}{\manprod{\Lieman{\Liegroup{}}}{\Lieman{\Liegroup{}}}}.
\end{equation}
Additionally, it is known that,
\begin{align}
&\gop{}\in\mapdifclass{\infty}{\manprod{\Lieman{\Liegroup{}}}{\Lieman{\Liegroup{}}}}{\Lieman{\Liegroup{}}},
\label{lemproductofidentitycurvesofapairofleftinvfsisaLiemorphismpeq6}\\
&\diagmap{\R}\in\mapdifclass{\infty}{\RR}{\RR^2}.
\label{lemproductofidentitycurvesofapairofleftinvfsisaLiemorphismpeq7}
\end{align}
Therefore, $\eta$ is a smooth map from $\RR$ to $\Lieman{\Liegroup{}}$.
\begin{equation}\label{lemproductofidentitycurvesofapairofleftinvfsisaLiemorphismpeq8}
\eta\in\mapdifclass{\infty}{\RR}{\Lieman{\Liegroup{}}}.
\end{equation}
From a different point of view, by defining the mapping $\function{\zeta}{\R^2}{\G{}}$ as,
\begin{align}\label{lemproductofidentitycurvesofapairofleftinvfsisaLiemorphismpeq9}
\Foreach{\opair{t}{s}}{\R^2}
\func{\zeta}{\binary{t}{s}}&\eqdef
\func{\[\func{\maxintcurve{\Lieman{\Liegroup{}}}}{\binary{\avecf{1}}{\IG{}}}\]}{t}\gop{}
\func{\[\func{\maxintcurve{\Lieman{\Liegroup{}}}}{\binary{\avecf{2}}{\IG{}}}\]}{s}\cr
&=\func{\[\func{\maxintcurve{\Lieman{\Liegroup{}}}}{\binary{\avecf{2}}{\IG{}}}\]}{s}\gop{}
\func{\[\func{\maxintcurve{\Lieman{\Liegroup{}}}}{\binary{\avecf{1}}{\IG{}}}\]}{t},
\end{align}
clearly,
\begin{equation}\label{lemproductofidentitycurvesofapairofleftinvfsisaLiemorphismpeq10}
\eta=\cmp{\zeta}{\diagmap{\R}}.
\end{equation}
According to \Ref{EQpartialsLemma},
\begin{align}\label{lemproductofidentitycurvesofapairofleftinvfsisaLiemorphismpeq11}
\func{\[\der{\zeta}{\manprod{\RR}{\RR}}{\Lieman{\Liegroup{}}}\]}{\func{\prodmantan{\RR}{\RR}{t}{s}}{\binary{\Rtanidentity{\R}{0}}{\Rtanidentity{\R}{0}}}}=
&\hskip0.5\baselineskip\func{\[\der{\(\cmp{\zeta}{\leftparinj{\R}{\R}{0}}\)}{\RR}{\Lieman{\Liegroup{}}}\]}{\Rtanidentity{\R}{0}}\cr
&+\func{\[\der{\(\cmp{\zeta}{\rightparinj{\R}{\R}{0}}\)}{\RR}{\Lieman{\Liegroup{}}}\]}{\Rtanidentity{\R}{0}},
\end{align}
It is straightforward to check that,
\begin{align}
&\cmp{\zeta}{\leftparinj{\R}{\R}{0}}=\cmp{\gltrans{\LieG{\Liegroup{}}}{\func{\[\func{\maxintcurve{\Lieman{\Liegroup{}}}}{\binary{\avecf{2}}{\IG{}}}\]}{0}}}
{\[\func{\maxintcurve{\Lieman{\Liegroup{}}}}{\binary{\avecf{1}}{\IG{}}}\]},
\label{lemproductofidentitycurvesofapairofleftinvfsisaLiemorphismpeq12}\\
&\cmp{\zeta}{\rightparinj{\R}{\R}{0}}=\cmp{\gltrans{\LieG{\Liegroup{}}}{\func{\[\func{\maxintcurve{\Lieman{\Liegroup{}}}}{\binary{\avecf{1}}{\IG{}}}\]}{0}}}
{\[\func{\maxintcurve{\Lieman{\Liegroup{}}}}{\binary{\avecf{2}}{\IG{}}}\]},
\label{lemproductofidentitycurvesofapairofleftinvfsisaLiemorphismpeq13}
\end{align}
and hence according to \reflem{lemleftinvvfintegralcurveswitharbitraryIC},
\begin{align}
&\cmp{\zeta}{\leftparinj{\R}{\R}{0}}=\func{\maxintcurve{\Lieman{\Liegroup{}}}}{\binary{\avecf{1}}
{\func{\[\func{\maxintcurve{\Lieman{\Liegroup{}}}}{\binary{\avecf{2}}{\IG{}}}\]}{0}}}=
\func{\maxintcurve{\Lieman{\Liegroup{}}}}{\binary{\avecf{1}}{\IG{}}},
\label{lemproductofidentitycurvesofapairofleftinvfsisaLiemorphismpeq14}\\
&\cmp{\zeta}{\rightparinj{\R}{\R}{0}}=\func{\maxintcurve{\Lieman{\Liegroup{}}}}{\binary{\avecf{2}}
{\func{\[\func{\maxintcurve{\Lieman{\Liegroup{}}}}{\binary{\avecf{1}}{\IG{}}}\]}{0}}}=
\func{\maxintcurve{\Lieman{\Liegroup{}}}}{\binary{\avecf{2}}{\IG{}}}.
\label{lemproductofidentitycurvesofapairofleftinvfsisaLiemorphismpeq15}
\end{align}
That is,
$\cmp{\zeta}{\leftparinj{\R}{\R}{0}}$ is the maximal integral-curve of $\avecf{1}$ on $\Lieman{\Liegroup{}}$ with the initial condition
$\func{\[\func{\maxintcurve{\Lieman{\Liegroup{}}}}{\binary{\avecf{2}}{\IG{}}}\]}{0}=\IG{}$,
and $\cmp{\zeta}{\rightparinj{\R}{\R}{0}}$ is the maximal integral-curve of $\avecf{2}$ with the initial condition
$\func{\[\func{\maxintcurve{\Lieman{\Liegroup{}}}}{\binary{\avecf{1}}{\IG{}}}\]}{0}=\IG{}$.
Therefore, according to \refdef{defintegralcurvesofasmoothvectorfield},
\begin{align}
&\func{\[\der{\(\cmp{\zeta}{\leftparinj{\R}{\R}{0}}\)}{\RR}{\Lieman{\Liegroup{}}}\]}{\Rtanidentity{\R}{0}}=
\func{\avecf{1}}{\IG{}},
\label{lemproductofidentitycurvesofapairofleftinvfsisaLiemorphismpeq16}\\
&\func{\[\der{\(\cmp{\zeta}{\rightparinj{\R}{\R}{0}}\)}{\RR}{\Lieman{\Liegroup{}}}\]}{\Rtanidentity{\R}{0}}=
\func{\avecf{2}}{\IG{}},
\label{lemproductofidentitycurvesofapairofleftinvfsisaLiemorphismpeq17}
\end{align}
and thus, \Ref{lemproductofidentitycurvesofapairofleftinvfsisaLiemorphismpeq11},
\Ref{lemproductofidentitycurvesofapairofleftinvfsisaLiemorphismpeq16},
\Ref{lemproductofidentitycurvesofapairofleftinvfsisaLiemorphismpeq17} yield,
\begin{align}\label{lemproductofidentitycurvesofapairofleftinvfsisaLiemorphismpeq18}
\func{\[\der{\zeta}{\manprod{\RR}{\RR}}{\Lieman{\Liegroup{}}}\]}{\func{\prodmantan{\RR}{\RR}{t}{s}}{\binary{\Rtanidentity{\R}{0}}{\Rtanidentity{\R}{0}}}}&=
\func{\avecf{1}}{\IG{}}+\func{\avecf{2}}{\IG{}}\cr
&=\func{\(\avecf{1}+\avecf{2}\)}{\IG{}}.
\end{align}
Furthermore, it is known that,
\begin{equation}\label{lemproductofidentitycurvesofapairofleftinvfsisaLiemorphismpeq19}
\func{\der{\diagmap{\R}}{\RR}{\RR^2}}{\Rtanidentity{\R}{0}}=
\func{\prodmantan{\RR}{\RR}{t}{s}}{\binary{\Rtanidentity{\R}{0}}{\Rtanidentity{\R}{0}}}.
\end{equation}
According to \Ref{lemproductofidentitycurvesofapairofleftinvfsisaLiemorphismpeq10},
\Ref{lemproductofidentitycurvesofapairofleftinvfsisaLiemorphismpeq18},
\Ref{lemproductofidentitycurvesofapairofleftinvfsisaLiemorphismpeq19},
and using the chain rule of differentiation,
\begin{align}\label{lemproductofidentitycurvesofapairofleftinvfsisaLiemorphismpeq20}
\func{\[\der{\eta}{\RR}{\Lieman{\Liegroup{}}}\]}{\Rtanidentity{\R}{0}}&=
\func{\[\cmp{\(\der{\zeta}{\manprod{\RR}{\RR}}{\Lieman{\Liegroup{}}}\)}{\(\der{\diagmap{\R}}{\RR}{\RR^2}\)}\]}
{\Rtanidentity{\R}{0}}\cr
&=\func{\(\der{\zeta}{\manprod{\RR}{\RR}}{\Lieman{\Liegroup{}}}\)}
{\func{\[\der{\diagmap{\R}}{\RR}{\RR^2}\]}{\Rtanidentity{\R}{0}}}\cr
&=\func{\[\der{\zeta}{\manprod{\RR}{\RR}}{\Lieman{\Liegroup{}}}\]}{\func{\prodmantan{\RR}{\RR}{t}{s}}{\binary{\Rtanidentity{\R}{0}}{\Rtanidentity{\R}{0}}}}\cr
&=\func{\(\avecf{1}+\avecf{2}\)}{\IG{}}.
\end{align}
Therefore, according to \Ref{lemproductofidentitycurvesofapairofleftinvfsisaLiemorphismpeq3},
\Ref{lemproductofidentitycurvesofapairofleftinvfsisaLiemorphismpeq8}, and
\Ref{lemproductofidentitycurvesofapairofleftinvfsisaLiemorphismpeq20},
$\eta$ is a smooth-group-morphism from $\opair{\RR}{+}$ to $\Liegroup{}$, such that its differential at $\Rtanidentity{\R}{0}$
coincides with the value of the smooth vector-field $\avecf{1}+\avecf{2}$ at $\IG{}$. Also, since each $\avecf{1}$ and $\avecf{2}$ is
a left-invariant vector-field on $\Liegroup{}$, so is their sum. Thus, bease on
\reflem{lemoneparametersubgroupisidentitycurve}, it is inferred that $\eta$ is the maximal integral-curve of $\avecf{1}+\avecf{2}$
with the initial condition $\IG{}$. That is,
\begin{equation}
\eta=\func{\maxintcurve{\Lieman{\Liegroup{}}}}{\binary{\avecf{1}+\avecf{2}}{\IG{}}}.
\end{equation}
\endlem
\theorem\label{thmexpofsumoftwocommutativeleftinvvfs}
Each $\avecf{1}$ and $\avecf{2}$ is taken as an arbitrary element of $\Leftinvvf{\Liegroup{}}$.
\begin{align}
\liebracket{\avecf{1}}{\avecf{2}}{\Lieman{\Liegroup{}}}=\zerovec{\Vecf{\Lieman{\Liegroup{}}}{\infty}}\then
\func{\expLie{\Liegroup{}}}{\avecf{1}+\avecf{2}}=\func{\expLie{\Liegroup{}}}{\avecf{1}}\gop{}
\func{\expLie{\Liegroup{}}}{\avecf{2}}.
\end{align}
\proof
Assume that $\liebracket{\avecf{1}}{\avecf{2}}{\Lieman{\Liegroup{}}}=\zerovec{\Vecf{\Lieman{\Liegroup{}}}{\infty}}$.
Then according to \reflem{lemproductofidentitycurvesofapairofleftinvfsisaLiemorphism}, and taking $\eta$
as the mapping defined in it, and also according to the definition of the exponential-mapping of a smooth group
(\refdef{defexponentialmappingofsmoothgroup}),
\begin{align}
\func{\expLie{\Liegroup{}}}{\avecf{1}+\avecf{2}}&=
\func{\[\func{\maxintcurve{\Lieman{\Liegroup{}}}}{\binary{\avecf{1}+\avecf{2}}{\IG{}}}\]}{1}\cr
&=\func{\eta}{1}\cr
&=\func{\[\func{\maxintcurve{\Lieman{\Liegroup{}}}}{\binary{\avecf{1}}{\IG{}}}\]}{1}\gop{}
\func{\[\func{\maxintcurve{\Lieman{\Liegroup{}}}}{\binary{\avecf{2}}{\IG{}}}\]}{1}\cr
&=\func{\expLie{\Liegroup{}}}{\avecf{1}}\gop{}\func{\expLie{\Liegroup{}}}{\avecf{2}}.
\end{align}
\endthm
\corollary
Each $\avecf{1}$ and $\avecf{2}$ is taken as an arbitrary element of $\Leftinvvf{\Liegroup{}}$. If $\avecf{1}$ and $\avecf{2}$
commute, then $\func{\expLie{\Liegroup{}}}{\avecf{1}}$ and $\func{\expLie{\Liegroup{}}}{\avecf{2}}$ also commute in the group
$\LieG{\Liegroup{}}$.
\begin{equation}
\liebracket{\avecf{1}}{\avecf{2}}{\Lieman{\Liegroup{}}}=\zerovec{\Vecf{\Lieman{\Liegroup{}}}{\infty}}\then
\func{\expLie{\Liegroup{}}}{\avecf{1}}\gop{}
\func{\expLie{\Liegroup{}}}{\avecf{2}}=
\func{\expLie{\Liegroup{}}}{\avecf{2}}\gop{}
\func{\expLie{\Liegroup{}}}{\avecf{21}}
\end{equation}
\proof
It is an immediate consequence of \refthm{thmexpofsumoftwocommutativeleftinvvfs}.
\endcor
\theorem
$\avecf{}$ is taken as an element of $\Leftinvvf{\Liegroup{}}$.
\begin{equation}
\Foreach{\opair{t}{s}}{\R^2}
\func{\expLie{\Liegroup{}}}{\(t+s\)\avecf{}}=
\func{\expLie{\Liegroup{}}}{t\avecf{}+s\avecf{}}=
\func{\expLie{\Liegroup{}}}{t\avecf{}}\gop{}
\func{\expLie{\Liegroup{}}}{s\avecf{}}.
\end{equation}
\proof
Considering the Lie-algebraic structure of left-invariant vecor-fields on $\Liegroup{}$,
it is clear that $t\avecf{}$ and $s\avecf{}$ commute for every real numbers $t$ and $s$.
\begin{equation}
\Foreach{\opair{t}{s}}{\R^2}
\liebracket{t\avecf{}}{s\avecf{}}{\Lieman{\Liegroup{}}}=\zerovec{\Vecf{\Lieman{\Liegroup{}}}{\infty}}.
\end{equation}
Thus the intended result trivially follows from \refthm{thmexpofsumoftwocommutativeleftinvvfs}.
\endthm
\lemma\label{lemexponentialmappingrescaled}
$\avecf{}$ is taken as an element of $\Leftinvvf{\Liegroup{}}$.
\begin{equation}
\Foreach{t}{\R}
\func{\[\func{\maxintcurve{\Lieman{\Liegroup{}}}}{\binary{\avecf{}}{\IG{}}}\]}{t}=
\func{\expLie{\Liegroup{}}}{t\avecf{}}.
\end{equation}
\proof
According to \refthm{thmintegralcurvesofvectorfieldswithvariablescale} and \refdef{defexponentialmappingofsmoothgroup},
\begin{align}
\Foreach{t}{\compl{\R}{0}}
\func{\expLie{\Liegroup{}}}{t\avecf{}}&=
\func{\[\func{\maxintcurve{\Lieman{\Liegroup{}}}}{\binary{t\avecf{}}{\IG{}}}\]}{1}\cr
&=\func{\[\func{\maxintcurve{\Lieman{\Liegroup{}}}}{\binary{\avecf{}}{\IG{}}}\]}{t}.
\end{align}
The case for $t=0$ is trivial, considering that $0\avecf{}$ is the zero vector-field on $\Lieman{\Liegroup{}}$,
which is a complete vector-field whose maximal integral-curves are the stationary points, and hence,
\begin{align}
\func{\expLie{\Liegroup{}}}{0\avecf{}}&=
\func{\[\func{\maxintcurve{\Lieman{\Liegroup{}}}}{\binary{\zerovec{}}{\IG{}}}\]}{1}\cr
&=\IG{}\cr
&=\func{\[\func{\maxintcurve{\Lieman{\Liegroup{}}}}{\binary{\avecf{}}{\IG{}}}\]}{0}.
\end{align}
\endlem
\corollary\label{corinverseofexpofaleftinvvf}
For every left-invarian vector-field $\avecf{}$ on $\Liegroup{}$,
$\func{\expLie{\Liegroup{}}}{-\avecf{}}$ equals the inverse of $\func{\expLie{\Liegroup{}}}{\avecf{}}$
in the group $\LieG{\Liegroup{}}$.
\begin{equation}
\Foreach{\avecf{}}{\Leftinvvf{\Liegroup{}}}
\func{\expLie{\Liegroup{}}}{-\avecf{}}=\invG{\[\func{\expLie{\Liegroup{}}}{\avecf{}}\]}{\LieG{\Liegroup{}}}.
\end{equation}
\proof
$\avecf{}$ is taken as an arbitrary left-invariant vector-field on $\Liegroup{}$.
According to \refdef{defexponentialmappingofsmoothgroup} and
\reflem{lemexponentialmappingrescaled}, invoking the canonical linear structure of $\Leftinvvf{\Liegroup{}}$,
and considering that $\func{\maxintcurve{\Lieman{\Liegroup{}}}}{\binary{\avecf{}}{\IG{}}}$ is a group-homomorphism
from $\opair{\R}{+}$ to $\LieG{\Liegroup{}}$ based on \reflem{lemidentitycurveisaoneparametersubgroup},
\begin{align}
\func{\expLie{\Liegroup{}}}{-\avecf{}}&=
\func{\expLie{\Liegroup{}}}{-1\avecf{}}\cr
&=\func{\[\func{\maxintcurve{\Lieman{\Liegroup{}}}}{\binary{\avecf{}}{\IG{}}}\]}{-1}\cr
&=\invG{\[\func{\expLie{\Liegroup{}}}{\avecf{}}\]}{\LieG{\Liegroup{}}}.
\end{align}
\endcor
\theorem\label{thmderivativeofexpatidentity}
The differential of the inverse mapping of the group $\LieG{\Liegroup{}}$ when restricted to the tangent-space
of the identity element of $\LieG{\Liegroup{}}$, equals $-\identity{\tanspace{\IG{}}{\Lieman{\Liegroup{}}}}$. Equivalently,
\begin{equation}
\Foreach{\vv{}}{\tanspace{\IG{}}{\Lieman{\Liegroup{}}}}
\func{\[\der{\ginv{\LieG{\Liegroup{}}}}{\Lieman{\Liegroup{}}}{\Lieman{\Liegroup{}}}\]}{\vv{}}=-\vv{}.
\end{equation}
\proof
$\vv{}$ is taken as an arbitrary element of the tangent-space of $\Lieman{\Liegroup{}}$ at $\IG{}$. Since
the set of left-invariant vector-fields on $\Liegroup{}$ is in a one-to-one correpondence with $\tanspace{\IG{}}{\Lieman{\Liegroup{}}}$
via the mapping $\liegvftan{\Liegroup{}}$ (\refcor{corlinvvftanspalinearcecorrespondence}), according to
\refdef{deflinvvftanspacecorrespondence} it is clear that there exists a left-invariant vector-field $\avecf{}$ on $\Liegroup{}$
such that $\vv{}=\func{\avecf{}}{\IG{}}$. Moreover, considering that
$\func{\maxintcurve{\Lieman{\Liegroup{}}}}{\binary{\avecf{}}{\IG{}}}$ is a group-homomorphism
from $\opair{\R}{+}$ to $\LieG{\Liegroup{}}$ based on \reflem{lemidentitycurveisaoneparametersubgroup}, and
according to \refthm{thmintegralcurvesofvectorfieldswithvariablescale},
\begin{align}
\Foreach{t}{\R}
\func{\[\cmp{\ginv{\LieG{\Liegroup{}}}}{\func{\maxintcurve{\Lieman{\Liegroup{}}}}{\binary{\avecf{}}{\IG{}}}}\]}{t}&=
\invG{\bigg(\func{\[\func{\maxintcurve{\Lieman{\Liegroup{}}}}{\binary{\avecf{}}{\IG{}}}\]}{t}\bigg)}{\LieG{\Liegroup{}}}\cr
&=\func{\[\func{\maxintcurve{\Lieman{\Liegroup{}}}}{\binary{\avecf{}}{\IG{}}}\]}{-t}\cr
&=\func{\[\func{\maxintcurve{\Lieman{\Liegroup{}}}}{\binary{-\avecf{}}{\IG{}}}\]}{t},
\end{align}
which means,
\begin{equation}
\cmp{\ginv{\LieG{\Liegroup{}}}}{\func{\maxintcurve{\Lieman{\Liegroup{}}}}{\binary{\avecf{}}{\IG{}}}}=
\func{\maxintcurve{\Lieman{\Liegroup{}}}}{\binary{-\avecf{}}{\IG{}}}.
\end{equation}
Therefore, according to \refdef{defintegralcurvesofasmoothvectorfield} and
\refthm{thmmaximalintegralcurveisanintegralcurve}, applying the chain rule of differentiation,
\begin{align}
\func{\[\der{\ginv{\LieG{\Liegroup{}}}}{\Lieman{\Liegroup{}}}{\Lieman{\Liegroup{}}}\]}{\vv{}}&=
\func{\[\der{\ginv{\LieG{\Liegroup{}}}}{\Lieman{\Liegroup{}}}{\Lieman{\Liegroup{}}}\]}{\func{\avecf{}}{\IG{}}}\cr
&=\func{\[\der{\ginv{\LieG{\Liegroup{}}}}{\Lieman{\Liegroup{}}}{\Lieman{\Liegroup{}}}\]}
{\func{\avecf{}}{\func{\[\func{\maxintcurve{\Lieman{\Liegroup{}}}}{\binary{\avecf{}}{\IG{}}}\]}{0}}}\cr
&=\func{\(\cmp{\[\der{\ginv{\LieG{\Liegroup{}}}}{\Lieman{\Liegroup{}}}{\Lieman{\Liegroup{}}}\]}
{\[\der{\func{\maxintcurve{\Lieman{\Liegroup{}}}}{\binary{\avecf{}}{\IG{}}}}{\RR}{\Lieman{\Liegroup{}}}\]}\)}{\Rtanidentity{\R}{0}}\cr
&=\func{\(\der{\[\cmp{\ginv{\LieG{\Liegroup{}}}}{\func{\maxintcurve{\Lieman{\Liegroup{}}}}{\binary{\avecf{}}{\IG{}}}}\]}{\RR}{\Lieman{\Liegroup{}}}\)}
{\Rtanidentity{\R}{0}}\cr
&=\func{\bigg(\der{\func{\maxintcurve{\Lieman{\Liegroup{}}}}{\binary{-\avecf{}}{\IG{}}}}{\RR}{\Lieman{\Liegroup{}}}\bigg)}{\Rtanidentity{\R}{0}}\cr
&=\func{\(-\avecf{}\)}{\func{\[\func{\maxintcurve{\Lieman{\Liegroup{}}}}{\binary{-\avecf{}}{\IG{}}}\]}{0}}\cr
&=\func{\(-\avecf{}\)}{\IG{}}\cr
&=-\func{\avecf{}}{\IG{}}=-\vv{}.
\end{align}
\endthm
\theorem\label{thmliealgebraofabelianliegroupistrivial}
If $\Liegroup{}$ is an abelian smooth group (that is, its intrinsic group structure $\LieG{\Liegroup{}}$ is commutative),
then the canonical Lie-algebra of $\Liegroup{}$ is also abelian. That is,
\begin{align}
\bigg(\Foreach{\opair{\g{1}}{\g{2}}}{\Cprod{\G{}}{\G{}}}\g{1}\gop{}\g{2}=\g{2}\gop{}\g{1}\bigg)
\then
\bigg(\Foreach{\opair{\avecf{1}}{\avecf{2}}}{\Leftinvvf{\Liegroup{}}^{\times 2}}
\liebracket{\avecf{1}}{\avecf{2}}{\Liegroup{}}=\zerovec{\Vecf{\Lieman{\Liegroup{}}}{\infty}}\bigg).
\end{align}
\proof
It is assumed that $\Liegroup{}$ is abelian, that is,
\begin{equation}
\Foreach{\opair{\g{1}}{\g{2}}}{\Cprod{\G{}}{\G{}}}\g{1}\gop{}\g{2}=\g{2}\gop{}\g{1}.
\end{equation}
Then clearly,
\begin{align}
\Foreach{\opair{\g{1}}{\g{2}}}{\Cprod{\G{}}{\G{}}}
\func{\ginv{\LieG{\Liegroup{}}}}{\g{1}\gop{}\g{2}}&=
\invG{\(\g{1}\gop{}\g{2}\)}{}\cr
&=\invG{\g{2}}{}\gop{}\invG{\g{1}}{}\cr
&=\invG{\g{1}}{}\gop{}\invG{\g{2}}{}\cr
&=\func{\ginv{\LieG{\Liegroup{}}}}{\g{1}}\gop{}
\func{\ginv{\LieG{\Liegroup{}}}}{\g{2}},
\end{align}
which means the inverse mapping of the group $\LieG{\Liegroup{}}$ is a group-homomorphism (actually a group-isomorphism).
\begin{equation}
\ginv{\LieG{\Liegroup{}}}\in\GHom{\LieG{\Liegroup{}}}{\LieG{\Liegroup{}}}.
\end{equation}
In addition, $\ginv{\LieG{\Liegroup{}}}$ is also a smooth map from $\Lieman{\Liegroup{}}$ to $\Lieman{\Liegroup{}}$ according to
\refthm{thminversemappingisdiffeomorphism}.
Therefore $\ginv{\LieG{\Liegroup{}}}$ is a smooth group-morphism (actually a smooth-group-isomorphism) from $\Liegroup{}$ to itself.
Thus $\func{\indliemor{\Liegroup{}}{\Liegroup{1}}}{\ginv{\LieG{\Liegroup{}}}}$ is a Lie-algebra-morphism from
the canonical Lie-algebra of $\Liegroup{}$ to itself.
That is (according to \refthm{thmdefinducedliealgebramorphismfromliemorphismpreservesliebracket}),
\begin{align}
&\Foreach{\opair{\avecf{1}}{\avecf{2}}}{\Cprod{\Leftinvvf{\Liegroup{}}}{\Leftinvvf{\Liegroup{}}}}\cr
&\func{\[\func{\indliemor{\Liegroup{}}{\Liegroup{}}}{\ginv{\LieG{\Liegroup{}}}}\]}{\liebracket{\avecf{1}}{\avecf{2}}{\Liegroup{}}}=
\liebracket{\func{\[\func{\indliemor{\Liegroup{}}{\Liegroup{}}}{\ginv{\LieG{\Liegroup{}}}}\]}{\avecf{1}}}
{\func{\[\func{\indliemor{\Liegroup{}}{\Liegroup{}}}{\ginv{\LieG{\Liegroup{}}}}\]}{\avecf{2}}}{\Liegroup{}}.
\end{align}
Additionally, according to \refdef{definducedliealgebramorphismfromliemorphism},
\refdef{deflinvvftanspacecorrespondence}, considering that $\liegvftan{\Liegroup{}}$ is a bijection
from $\Leftinvvf{\Liegroup{}}$ to $\tanspace{\IG{}}{\Lieman{\Liegroup{}}}$, \refthm{thmderivativeofexpatidentity} yields,
\begin{align}
\Foreach{\avecf{}}{\Leftinvvf{\Liegroup{}}}
\func{\[\func{\indliemor{\Liegroup{}}{\Liegroup{}}}{\ginv{\LieG{\Liegroup{}}}}\]}{\avecf{}}&=
\func{\(\cmp{\finv{\liegvftan{\Liegroup{}}}}{\cmp{\[\der{\ginv{\LieG{\Liegroup{}}}}
{\Lieman{\Liegroup{}}}{\Lieman{\Liegroup{}}}\]}{\liegvftan{\Liegroup{}}}}\)}{\avecf{}}\cr
&=\func{\finv{\liegvftan{\Liegroup{}}}}
{\func{\[\func{\indliemor{\Liegroup{}}{\Liegroup{}}}{\ginv{\LieG{\Liegroup{}}}}\]}{\func{\avecf{}}{\IG{}}}}\cr
&=\func{\finv{\liegvftan{\Liegroup{}}}}{-\func{\avecf{}}{\IG{}}}\cr
&=\func{\finv{\liegvftan{\Liegroup{}}}}{\func{\(-\avecf{}\)}{\IG{}}}\cr
&=-\avecf{}.
\end{align}
Thus, considering that $\Leftinvvf{\Liegroup{}}$ endowed with the Lie-bracket operation on it is a Lie-algebra
(the canonical Lie-algebra of $\Liegroup{}$), according to axioms of the Lie-algebra structure, it is clear that,
\begin{align}
&\Foreach{\opair{\avecf{1}}{\avecf{2}}}{\Cprod{\Leftinvvf{\Liegroup{}}}{\Leftinvvf{\Liegroup{}}}}\cr
&\begin{aligned}
-\liebracket{\avecf{1}}{\avecf{2}}{\Liegroup{}}&=
\func{\[\func{\indliemor{\Liegroup{}}{\Liegroup{}}}{\ginv{\LieG{\Liegroup{}}}}\]}{\liebracket{\avecf{1}}{\avecf{2}}{\Liegroup{}}}\cr
&=\liebracket{\func{\[\func{\indliemor{\Liegroup{}}{\Liegroup{}}}{\ginv{\LieG{\Liegroup{}}}}\]}{\avecf{1}}}
{\func{\[\func{\indliemor{\Liegroup{}}{\Liegroup{}}}{\ginv{\LieG{\Liegroup{}}}}\]}{\avecf{2}}}{\Liegroup{}}\cr
&=\liebracket{-\avecf{1}}{-\avecf{2}}{\Liegroup{}}\cr
&=\liebracket{\avecf{1}}{\avecf{2}}{\Liegroup{}}.
\end{aligned}
\end{align}
So since the underlying field of the vector-space $\VLeftinvvf{\Liegroup{}}$ (or equivalenty that of the Lie-algebra
$\LiegroupLiealgebra{\Liegroup{}}$) is $\R$ (a field of charactersitic $0$), it then becomes evident that,
\begin{equation}
\Foreach{\opair{\avecf{1}}{\avecf{2}}}{\Cprod{\Leftinvvf{\Liegroup{}}}{\Leftinvvf{\Liegroup{}}}}
\liebracket{\avecf{1}}{\avecf{2}}{\Liegroup{}}=\zerovec{\Vecf{\Lieman{\Liegroup{}}}{\infty}}.
\end{equation}
\endthm
\theorem
The image of the exponential-mapping of $\Liegroup{}$ is included in the connected-component of the topological-space
$\lietops{\Liegroup{}}$ that contains $\IG{}$.
\proof
According to \refdef{defintegralcurvesofasmoothvectorfield} and
\refdef{defexponentialmappingofsmoothgroup},
for every left-invariant vector-field $\avecf{}$ on $\Liegroup{}$, the points $\IG{}$ and
$\func{\expLie{\Liegroup{}}}{\avecf{}}$ of the topological-space $\lietops{\Liegroup{}}$ are connected via the curve
$\func{\maxintcurve{\Lieman{\Liegroup{}}}}{\binary{\avecf{}}{\IG{}}}$. In addition, since
$\func{\maxintcurve{\Lieman{\Liegroup{}}}}{\binary{\avecf{}}{\IG{}}}$ is a smooth map from $\RR$ to\
$\Lieman{\Liegroup{}}$, it is clearly a continuous map from $\topR{\R}$ to the topological-spac $\lietops{\Liegroup{}}$, and hence
so is its restriction to the closed interval $\cinterval{0}{1}$ of $\R$, which is then a path of $\lietops{\Liegroup{}}$. Therefore,
there exists a path in $\lietops{\Liegroup{}}$ (a continuous map from $\topR{\R}$ to $\lietops{\Liegroup{}}$) that joins
$\IG{}$ and $\func{\expLie{\Liegroup{}}}{\avecf{}}$. Hence, these two points lie in the same path-connected component of
the underlying topological-space of the manifold $\Lieman{\Liegroup{}}$ which is obviously the topological-space $\lietops{\Liegroup{}}$.
Moreover, considering that
path-connectedness and connectedness are equivalent properties of the underlying topological-space of a manifold
(the set of all connected sets equals the set of all path-connected sets of the topological-space), it is ultimately
inferred that $\IG{}$ and $\func{\expLie{\Liegroup{}}}{\avecf{}}$ belong to the same connected-component of $\lietops{\Liegroup{}}$.
\endthm
\newpage
\Bibliography{}
\renewcommand{\addcontentsline}[3]{}

\let\addcontentsline\oldaddcontentsline

\end{document}